\documentclass[12pt]{article}
\usepackage{amssymb}
\newcommand{\N}{\mbox{I$\!$N} }
\newcommand{\R}{\mbox{I$\!$R}}

\newcommand{\qed}{{\hfill {$\rlap{$\sqcap$}\sqcup$}}\\[0.2in]\hspace*{0.5in}}
\newcommand{\qedwh}{{\hfill {$\rlap{$\sqcap$}\sqcup$}}\\[0.2in]}
\newcommand{\bk}{\\[0.1in] \hspace*{0.5in} }
\newcommand{\btd}{\bigtriangledown}
\newcommand{\mfor}{\ \ \ \ {\mbox{for}} \ \ }

\newcommand{\parallelsum}{\mathbin{\!/\mkern-5mu/\!}}
\begin{document}

\begin{center}
{\LARGE {\bf Construction of Blow\,-\,up Sequences for  the}} \medskip \medskip \smallskip \smallskip\\
{\LARGE {\bf Prescribed Scalar Curvature Equation }} \medskip \medskip \smallskip \smallskip  \\
{\LARGE {\bf on $S^n$. IV.  Clustered Blow\,-\,ups  }}
\end{center}

\vspace{0.8in}

\centerline{\Large  {Man Chun  {\LARGE L}EUNG${\,}^{*}$  }}

\vspace*{0.2in}

\centerline{\large {National University of Singapore }}

\vspace{0.73in}


\begin{abstract}
\vspace{0.33in}
\noindent For the prescribed scalar curvature equation on $\,S^n\,$ ($\,n\,\ge\ 6\,$)\,,\, we consider the situation where  the  number of bubbles tends to infinity in the
Lyapunov\,-\,Schmidt (\,finite dimension\,) reduction method\,.\, In an outstanding paper by Wei and Yan \cite{Wei-Yan}\,,\, the special case where the bubbles are arranged ``\,evenly\,"\, (\,close to
a great circle\,) is considered\,.\, Here we are concerned with the generic  scenario, where the
 bubbles are ``\,planted\," (\,arranged\,) in general position\,,\, in adjacent  to the critical points of the prescribed function\,.\,  The main interest of this article is to   extract key information in the finite dimensional reduced functional,
as well as in  its first partial derivatives\,.\, Overall\,,\, there are two main contributions in the \,``expansions\,"\,,\,  namely, interaction between bubbles and curvature involvement.  As each bubble becomes infinitesimally  close to other bubbles in the neighborhood,   the errors in the  expansions are carefully estimated in terms of the geometric parameters. This lays the foundation for our next step on constructing blow\,-\,up sequences of solutions with dense blow\,-\,up points (clustered blow\,-\,up)\,,\, where the (\,non\,-\,radial\,-\,symmetric\,) prescribed scalar curvature function  has non\,-\,isolated critical points. 
\end{abstract}

\vspace*{0.5in}

{\bf Key Words}\,:  Scalar Curvature Equation\,;\, Blow\,-\,up\,;\, Critical Points\,;\, Lyapunov\,-\,Schmidt \\[0.075in]
\hspace*{1in} (\,Finite Dimension\,) Reduction Method. \\[0.1in]
{\bf 2020 AMS MS Classification}\,:\ \ Primary 35\,J\,60\,;\ \ Secondary 53\,C\,21. \\[0.1in]
$\overline{\ \ \ \ \ \ \ \ \ \ \ \ \ \ \ \ \ \ \ \ \  \ \ \ \ \ \ \ \ \ \ \ \ \ \ \ \ \ \ \ \ \ \ \ }$\\[0.05in]
\noindent{$\!\!{\!}^{*}$ \,{\tt{matlmc\,@\,nus.edu.sg}}

\newpage


{\bf \large {\bf  1. \ \ Introduction.}}\\[0.2in]
With vast varieties in geometry\,,\,   one inclines to look for a representative which is ``\,as simple as possible\,"\,.\, Yamabe problem is one good example along  this line. On the other hand, given a simple object, one can also ask how many ``varieties\," are within. Minkowski problem is a classical example. Another prominent example is the prescribed scalar curvature problem on $\,S^n\,$ (\,also known as the  Kazdan\,-\,Warner\,/\,Nirenberg problem \cite{Kazdan} \cite{Tristan}\,)\,,\, which is governed by the semilinear partial differential equation
$$
\ \ \ \   \Delta_{1}\, U\ -\ [\ {\tilde c}_n\,n\,(\,n\,-\,1\,)\ ]\,U\ +\ [\ {\tilde c}_n\,{\cal K}\,]\,U^{\,{{n\,+\,2}\over {\,n\,-\ 2\,}}}\ \, = \ 0\ \ \ \ \ \   {\mbox{in}}\ \ \ \ \ S^n\ \ \ \ \ \ \ (\,U \ > \ 0\,)\,.\leqno (1.1)
$$
Here $\,{\cal K}\,$  is a given $\,C^2\,$ function on $\,S^n\,$ (\,potentially the scalar curvature of the conformal metric $\,U^{4\over {n\,-\,2}}\,g_1\,$)\,.\, In (\,1.1\,)\,,\, $\,{\tilde c}_n \,=\,(n\,-\,2)/\,[\, 4\,(\,n\,-\,1\,)\,]\ .\,$ See \,{\bf \S 1\,f}\, for the rather standard notations we use\,.\, \bk
Close to half a century (\,cf. an early work in 1972   by Dimitri Koutroufiotis \cite{1972}, whose thesis adviser is   Louis Nirenberg\,)\,,\,  the apparently simple equation (1.1) serves as a vehicle for sophisticated techniques in nonlinear partial differential equations to be deployed and developed\,.\, The equation can also be branched out to complete  open manifolds, CR manifolds, Q\,-\,curvature\,.\, It is also related to the mean field equations. See  some recent works \, [\,1\,]\,,\, [\,3\,]\,--\,[\,14\,]\,,\,  [\,16\,]\,--\,[\,17\,]\,,\, [\,21\,]\,--\,[\,22\,]\,,\, [\,24\,]\,--\,[\,33\,]\,,\, [\,35\,]\,--\,[\,36\,]\,,\,  and the references therein\,.\, \bk
Via a stereographic projection $\,\dot{\cal P}\,$ onto $\,\R^n\,$ [\,see \,(\,4.42\,) }\,]\,,\,  equation   (1.1) can be transformed to
\begin{eqnarray*}
(\,1.2\,) \ \ \ \ \ \ \ \ \ \ \ \ \ \ \ \ \ \ \ \ \  \Delta\,u \ + \,[\ {\tilde c}_n\,K\,]\,u^{{n\,+\,2}\over {\,n\,-\ 2\,}} & = & 0\ \ \ \ \ \   {\mbox{in}}\ \ \ \ \ \R^n,\\[0.2in]
u \ > \ 0\,; \ \ K\,( \,y\,) & = & {\cal K}\,(\,x\,)\,, \ \ \ \ \ {\mbox{where}} \ \ y \ = \ \dot{\cal P}\,(\,x\,)\,. \ \ \ \  \ \ \ \ \ \ \ \ \ \ \ \ \ \ \ \ \ \ \
\end{eqnarray*}
A  stand\,-\,out feature of equation (\,1.2\,) is the exponent  $\,{{n\,+\,2}\over {\,n\,-\ 2\,}}\,.\,$ Indeed\,,\,   it is related to the  critical Sobolev embedding\,:\,   the injection $\,H^{\,1\,,\,2} \,(\,S^n\,)\, \hookrightarrow\, L^{{2\,n}\over {\,n\,-\ 2\,}}\, (\,S^n\,)\,$ is not compact, leading to the concentration phenomenon, or the blow\,-\,up process\,. \bk
The blow\,-\,up phenomenon has a rich structure in the case of $\,S^n\,.\,$
Intuitively, the building block is a simple bubble, which can be by itself (simple blow up), or in juxtaposition with another bubble (aggregated blow\,-\,up), or being stacked up with other bubble(s) (towering blow\,-\,up), or many  bubbles gathered together (clustered blow\,-\,up\,,\, in which the number of bubbles tends to infinity\,)\,.\, See \cite{Leung-Supported} for more detail.\bk
Each type of (\,collpased\,) blow\,-\,up is constructed with the help of Lyapunov\,-\,Schmidt (\,finite dimension\,)  reduction method (\,also known as the finite dimension reduction method\,)\,.\, See \cite{I} \cite{II} \cite{III} \cite {Wei-Yan}\,.\, In particular, for clustered blow\,-\,up,  Wei and Yan   consider bubbles arranged symmetrically along a great circle (\,number of bubbles tends to infinity\,)\,.\, They also put forward the following conjecture \cite {Wei-Yan}\,.

\newpage

{\it Assume that the maximal set} $\ \displaystyle{ {\cal M}  \, := \,  {\Large\{} \ y \, \in \, \R^n \ | \  K\,(\,y\,)\, = \ \max \,K \,\Large\}}\,$ {\it is an $\,m\,$-dimensional smooth manifold without boundary, where $\,1\,\le\,m\,\le \,n\,-\,1\,.$\, Then equation} (1.2) {\it admits infinitely many positive solutions.}\bk
In the interest of this article, we may also add  that the infinite number of solutions predicted in the conjecture together forms a clustered blow\,-\,up sequence of solutions\,.\,\bk
The case in \cite{Wei-Yan} (\,by Wei and Yan\,) is when $\,{\cal M}\,=\,S^{n\,-\,1}\ ,\ $ and  $\,K \,\in\,C^2\,(\,\R^n\,)\,$ is radially symmetric, that is, $\,K\,(\,y\,) \,=\, K_{\,\small{\mbox{r}}}\,(\,r\,)\,$,\, where $\,r\,=\,\Vert\,y\,\Vert\,,\,$ which fulfills  the conditions that there is a number $\,\ell \,\in\,[\ 2\,, \ n\,-\,2\,)\,$ such that
$$ \left(\,{\tilde c}_n \cdot  K_{\,\small{\mbox{r}}}\ \right)\,(\,r)\ = \  n\,(\,n\,-\,2\,) \,- \,c_o\,|\,r\,-\,1\,|^{\,\ell} \ +\ O\,\left(\,|\,r\,-\,1\,|^{\,\ell\ +\ \theta}\,\right) \leqno (\,1.3\,)
$$
for $\,r \,\in\,(\,1\,-\,\delta\,, \ 1\,+\,\delta\,)\,.\,$
By  arranging bubbles ``\,evenly\,"  close to a great circle of $\,S^{\,n\,-\,1}\,\subset\,\R^n$\,,\, they construct a sequence of solutions of equation (1.2) [\ which can be pulled back to $\,S^n\,$ as solutions of equation (1.1)\ ]\,,\, so that (\,effectively\,) every point along the great circle is a blow\,-\,up point (\,clustered blow\,-\,up point\,)\,.\,   Wei and Yan  remark that the spherical symmetric condition can be relaxed a bit    (\,see Remark 1.4 in \cite{Wei-Yan}\,)\,.\, Their proof makes use of this radial symmetry in an essential way\,,\, bring down  the complexity of the reduced functional.\, \{\ In (\,1.3\,)\,,\, $\,\theta\,\in\,(\,0\,, \ 1\,]\,$ is fixed\,,\,  and $\,\delta\,$ a small positive number.\, We refer to  \cite{Wei-Yan} for the full detail.\ \}

\vspace*{0.3in}

{\large{\bf \S\,1\,a.}}\ \  {\it Many bubbles case \,--\, Lyapunov\,-\,Schmidt} (\,{\it finite dimension}\,) {\it  reduction method  for equation} (1.2)\,.  \ \ Given a positive integer $\,\flat \ (\  \ge\,2\,;\ \,\to \,\infty\,)\,,\,$ consider (\,small\,) positive numbers
$$
 \!\!\!\!\!\!\!\!\!\!\!\!\!\!\!\lambda_{\,1}\,, \ \,\cdot\, \cdot\, \cdot\,, \ \,\lambda_{\,\,\flat}\ ,
$$
together with  distinct points (\,can be close to each other\,) in $\,\R^n\,$:
$$\ \ \ \ \ \ \ \ \ \ \ \ \ \ \ \ \ \ \ \ \ \ \ \ \ \ \ \
\xi_{\,1}\,, \ \cdot \cdot \cdot\,, \ \xi_{\,\flat} \ \ \ \ \ \ \ (\,\xi_{\,k} \ \not= \ \xi_{\,l} \ \ \ \ {\mbox{for}} \ \ j \ \not= \ k\,)\
$$
(\,they are called {\it bubble parameters}\,)\,.\,
 We restrict ourselves to the situation in which all the $\,\xi_{\,k}\,$ are inside a fixed compact set. That is\,,\, there exists a fixed compact set $\,{\cal C}_{\,[]} \, \subset \, \R^n\,$   such that
$$
 \{\  \xi_{\,1}\,, \ \cdot \cdot \cdot\,, \ \xi_{\,\flat}  \ \} \ \subset \ {\cal C}_{\,[]} \ \ \ \ \ [ \ {\mbox{independent \   \ on}} \ \  \flat \ ] \ . \leqno (\,1.4\,)
$$
Set
$$
V_{\,l}\,(\,y\,) \,=\,V_{\lambda_{\ \!l}\,,\ \xi_{\,l}}  \,=\,\left( {{\lambda_{\ \!l}}\over {\ \lambda_{\ \!l}^2 \ + \ \Vert\,y\,-\,\xi_{\,l}\,\Vert^{\,2}\ }} \right)^{\!\!{{\,n\,-\,2\,}\over 2}  } \ \ \mfor \ \ y \ \in \ \R^n \ \ \  (\ l\,=\,1\,,\ \cdot \cdot \cdot\,, \ \flat\,)\,. \leqno (1.5)
$$
This is referred to as a   standard\ \ bubble\,,\, and $\,\xi_{\,l} \,\in\,\R^n\,$ as the center of the bubble\,:
$$
\lim_{\lambda_{\ \!l} \,\to\,0^+ } \left( {{\lambda_{\ \!l}}\over {\ \lambda_{\ \!l}^2 \ + \ \Vert\,y\,-\,\xi_{\,l}\,\Vert^{\,2}\ }} \right)^{\!\!{{\,n\,-\,2\,}\over 2}  }_{\,\bigg\vert_{\ y \ = \ \xi_{\,l} }} \!\! = \  \infty\,.
$$
$V_{\,l}\,\,$ satisfies the equation
 $$
\Delta \,V_{\,l}\ + \ n\,(\,n\,-\,2)\,V_{\,l}^{{n\,+\,2}\over {\ n\,-\,2\ }} = \ 0 \ \ \ \ \ \ \ {\mbox{in}} \ \ \ \R^n\,. \leqno (1.6)
 $$
{\it With this simplified notation\,, we recognize that  the numbers  $\,\{\, \lambda_{\ \!l}\, \}_{\,l\,=\,1}^{\,\flat} \,$ and the points  $\,\{\,\xi_{\,l}\, \}_{\,l\,=\,1}^{\,\flat} \,$ may vary with\,} $\flat\,.\,$ That is\,,\, for a different  positive integer $\,\flat\,$,\ \,  $\,\{\, \lambda_{\ \!l}\, \}_{\,l\,=\,1}^{\,\flat}  \,$ and  $\,\{\, \xi_{\,l}\, \}_{\,l\,=\,1}^{\,\flat} \,$  may have different values and different positions\,,\, respectively\,.\,  The main object of this study is the ``\,plantation\," of the bubbles given by
$$
W_{\,\flat} \,:=\,\sum_{l\,=\,1}^\flat\, V_{\,l}\,. \leqno (1.7)
$$
To start with, we assume the following uniformity condition\,:
$$
{1\over { { {\bar C}_{\,1}}^{\ } }} \cdot {\tilde \lambda}_{\,\flat} \ \le \ \lambda_{\ \!l} \ \le \ {\bar C}_{\,1} \cdot   {\tilde \lambda}_{\,\flat}  \ \ \ \ \mfor \ \ l\,=\,1\,,\ \cdot \cdot \cdot\,, \ \flat\,, \ \ \ \ {\mbox{where}} \ \   {\tilde \lambda}_{\,\flat}  \ = \  \sqrt[^\flat\,]{\,\,\lambda_{\,1}\, \cdot \, \cdot\, \cdot\, \lambda_{\,\,\flat}\,\, } \   .\leqno (1.8)
$$
Here $\,{\bar C}_{\,1}\,$ is a fixed positive constant\,.\,
We also require that\,,\, for $\,1 \ \le \ j\,, \ k  \ \le \ \flat\,$,\,
$$
{\bf d}_{\,j\,,\,k} \ := \ {{\ \Vert\, \,\xi_{\,j}\ - \ \xi_{\,k}\,\Vert \  }\over
{\sqrt{\,\lambda_{\ \!j} \cdot \,\lambda_k\,}}}  \ \to \ \infty
 \ \ \ \ \ \ {\mbox{as}} \ \ \ {\tilde \lambda}_{\,\flat}  \ \to \ 0^+
 \ \ \ \ \ (\,j\ \not= \ k  \,)  \leqno (1.9)
$$
\{\ ${\bf d}_{\,j\,,\,k}\,$ is known as the   quasi\,-\,hyperbolic distance (\,see \cite{III}\,) between any two points $\,\xi_{\,j}\,\,$ and $\ \xi_{\,k}$\,\}\,,\,
together with the condition that
$$
\xi_{\,l} \ \to \ \ {\mbox{``\ \,a \ \ critical \ \ point\," \ \ of \ \ }} K  \ \ \ \ \ \   {\mbox{as}} \ \ \ {\tilde \lambda}_{\,\flat}  \ \to \ 0^+\,. \leqno (1.10)
$$
The functional corresponding to equation (1.2) is given by
\begin{eqnarray*}
(1.11) \ \ \ \ \ {\bf I}\,(\,f\,)& = & {1\over 2}\,\int_{\R^n}\,\langle\,\btd\,f\,,\,\btd\,f\,\rangle\, -\ \left(\,{{n\,-\,2}\over {2n}}\,\right) \cdot \!\int_{\R^n}\,(\,{\tilde c}_n\!\cdot K\,)\,f_+^{\,{{2n}\over {\,n\,-\ 2\,}}}\ \ \ \mfor f\,\in\,{\cal D}^{\,1,\,\,2}\ , \ \ \ \ \
\end{eqnarray*}
where $\,{\cal D}^{\,1,\,\,2}\,$ is the Hilbert space defined in (\,1.45\,)\,.\, The next expression  emphasises on  the separation of the constant curvature (\,unperturbed\,)  part and the curvature changing part.
\begin{eqnarray*}
(1.11)\,' \ \ \ \ \ \ \ {\bf I}\,(\,f\,)& = &  {1\over 2} \int_{\, \R^n} \left[\,  \langle\, \btd \,f\,,\,\, \btd\,f \  \rangle\  - \,  (\,n\, -\, 2\,)^{\,2}   \,f_+^{\,{{2n}\over {n - 2}}} \,\right] \ + \ \ \ \ \ \ \ \ \ \ \ \ \ \ \ \ \ \ \ \ \ \ \ \ \ \ \ \ \ \ \ \ \ \ \ \ \\[0.1in]
& \ & \ \ \ \  \leftarrow \ \ \ \  {\mbox{unperturbed \ \ \ \  part }} \ \ \rightarrow \\[0.1in]
& \ & \ \ \ \ \ \ \ \ \ \  \ \ \ \ \ \ \ \ \ \ \ \ \ \ \ \  +\   \left( \,{{n-2}\over {2n}}\ \right) \cdot  \int_{\R^n}\ [ \ n\,(\,n\ - \ 2\,) \ - \  (\,{\tilde c}_n\,K\,)\ ]  \,f_+^{\,{{2n}\over {n - 2}}} \ . \ \ \ \ \ \ \  \\[0.1in]
& \ &\ \ \ \ \ \ \ \ \ \  \ \ \ \ \ \ \ \ \ \ \ \ \ \ \ \  \ \ \ \ \ \ \ \ \ \ \ \ \ \ \ \ \ \ \    \leftarrow \ \    {\mbox{curvature \ \ \ \  portion }} \ \rightarrow
\end{eqnarray*}
 Indeed, the solution we look for has the form
$$
W_{\,\flat} \ + \ \phi_{\,\,\flat} \ = \ \left(\ \sum_{l\,=\,1}^\flat\, V_{\,l} \right) \ + \ \phi_{\,\,\flat}\,,
$$
where the first team comes from the unperturbed part\,,\, and the ``\,small\," contribution $\,\phi_{\,\,\flat} \,$ is used to adjust the curvature changing part.

\newpage

\hspace*{0.5in}We find it more natural, and the expressions become cleaner, by working with the  coupled quasi\,-\,hyperbolic gradient (\,introduced in \cite{III}\,)\,:  $$\displaystyle{\,(\,\lambda \cdot \btd\,) \,=\,\left( \  \lambda \cdot \partial_\lambda\,,\  \, \lambda \cdot \partial_{\xi_{\,|_{\,1}}}\, \,, \ \ \cdot \cdot \cdot\,, \  \lambda \cdot \partial_{\xi_{\,|_{\,n}}}\,   \right)\,,}$$

\vspace*{-0.2in}

where
 $$
 (\,\lambda \cdot \partial_\lambda\,)\,V_{\lambda\,,\ \xi} \ = \ \lambda \cdot {{\partial \,V_{\lambda\,,\ \xi}}\over {\partial \lambda}} \,\bigg\vert_{\,(\,\lambda\,,\ \xi\,)}\ \ \ \ {\mbox{and}} \ \ \ \  (\,\lambda \cdot \partial_{\xi_{\,|_{\,1}}}\,)\,V_{\lambda\,,\ \xi}  \ = \ \lambda \cdot {{\partial \,V_{\lambda\,,\ \xi}}\over {\partial \xi_{\,|_1}}} \,\bigg\vert_{\,(\,\lambda\,,\ \xi\,)}\,,\  \cdot \cdot \cdot\,. \leqno (1.12)
 $$
Here $\,(\,\lambda\,,\ \xi \,)\,=\,(\,\lambda\,;\ \xi_{\,|_{\,1}}\,, \ \cdot \cdot \cdot\,,  \ \xi_{\,|_{\,n}}\,)\,\in\,\R^+\times\,\R^n\,$ is treated as a point in the upper half space.
With this, the ``\,tangent space\," [\,of dimension $\,(\,n\,+\,1\,) \cdot \flat\,$]\, is given by [\,refer to {\bf \S\,A\,2} in the \cite{Notes}  concerning the linear independence of the collection\,,\, and the conditions to be applied\,]\\[0.1in]
(1.13)
$$
{\bf T}_{\,\flat} \,:=\,{\mbox{Span}} \ \ \bigg\{ \ (\,\lambda_{\ \!l} \cdot \partial_{\,\lambda_{\ \!l}}\,)\,V_{\,l}\,, \ \ \ (\,\lambda_{\ \!l} \cdot \partial_{\,\xi_{\,l_{|_{\,j}}}}\,)\ V_{\,l} \,, \   \ \ \  j \ = \ 1\,,\ \cdot \cdot \cdot\,, \ n   \ \ \ \ {\mbox{and}} \ \ \ l \ = \ 1\,,\ \cdot \cdot \cdot\,, \ \flat \  \bigg\}  \ \subset \ {\cal D}^{\,1,\,\,2}
$$
[\,see \,{\bf \S\,1\,g}\, for the definition of the Sobolev's space $\,{\cal D}^{\,1,\,\,2}\,$,\, and the inner product $\,\langle \ \,, \ \rangle_{\btd}\,$\,]\,.\,  The orthogonal complement is given by\\[0.1in]
(1.14)
\begin{eqnarray*}
{\cal D}^{\,1,\,\,2}_{\,\flat_\perp}&=&\bigg\{ \ \phi \,\in\, {\cal D}^{\,1,\,\,2}  \ \ \bigg\vert \ \  \langle\,\phi\,,\ \,[\ (\,\lambda_{\ \!l} \cdot \partial_{\lambda_{\ \!l}}\,)\,V_{\,l}\ ]\ \rangle_{\,\btd} \ = \ 0\ = \  \langle\,\phi\,,\ \,[\ (\,\lambda_{\ \!l} \cdot \partial_{\xi_{\,l_{|_{\,j}}}}\,)\,V_{\,l}\,]\ \rangle_{\,\btd}  \\[0.15in]
& \ & \ \ \ \ \ \ \ \ \ \ \  \ \ \ \ \  \ \ \ \  \ \ \ \ \   \ \ \ \  \ \  \ \ \ \ \ \ \ \ \ \ {\mbox{for \ \ all}} \ \ \   j \ = \ 1\,,\ \cdot \cdot \cdot\,, \ n   \ \ \ \ {\mbox{and}} \ \ \ l \ = \ 1\,,\ \cdot \cdot \cdot\,, \ \flat \  \bigg\}\ .
\end{eqnarray*}
It follows that
$$
{\cal D}^{\,1,\,\,2} \ = \ {\bf T}_{\,\flat} \ \oplus \ {\cal D}^{\,1,\,\,2}_{\,\flat_\perp}\ \,. \leqno (1.15)
$$

\vspace*{-0.05in}

Let

\vspace*{-0.25in}

$$
{\cal P}_{\,\flat} \ : \ {\cal D}^{\,1,\,\,2} \ \to \ {\cal D}^{\,1,\,\,2}_{\,\flat_\perp} \leqno (1.16)
$$
be the orthogonal projection. Boardly Speaking, the first step toward the finite dimension reduction is to solve the equation in the perpendicular direction, that is, to find a solution $\,\phi \,\in\,{\cal D}^{\,1,\,\,2}_{\,\flat_\perp}\,$ of the ``\,equation\,"
$$
{\cal P}_{\,\flat}  \left( \ {}^{``}{\,\Delta \,(\,W_{\,\flat} \,+\,\phi\,) \ + \  (\,{\tilde c}_n\, K)\,(\,W_{\,\flat} \,+\,\phi\,)_+^{{n\,+\,2}\over {\,n\,-\ 2\,}}\,}^{\,"}\ \right) \ = \ 0\ . \leqno (1.17)
$$
where the term inside \,$\,``\ \cdot \cdot \cdot\cdot\cdot \,"$\ \,is to be interpreted by means of the Riesz\,  Representation\,  Theory \,[\ cf. {\bf \S\,A\,1\,.\,b}\, of the \cite{Notes}\,]\,.\, The precise statement is found in Proposition 2.7\,.\, The cases for one and two bubbles are handled in \cite{I} \cite{II} \cite{III}\,.\, When the number of bubbles tends to infinity, this is done succinctly in \cite{Wei-Yan}\,,\,  under  the condition suitable for that paper, which is modified for the general situation in this present article. \{\,See also \cite{Bahri} for an earlier treatment\,,\, as well as \cite{Pino}\,.\.\}

\vspace*{0.3in}

{\large{\bf \S\,1\,b\,.}}\ \  {\it The reduced functional.  } \ \ Given $\,W_{\,\flat\,}\,$ as in (\,1.7\,)\,,\, let $\,\phi_{\,\,\flat}  \ \in \  {\cal D}^{1\,,\  2}_\perp\,$  be the \,``\,small\,"\, solution, specified by   equation (1.17) [\,to be made precise in {\bf \S} 2\ ]\,.\, Introduce the reduced functional
$$
 {\bf I}_{\,\cal R} \ = \ {\bf I}\ (\,W_{\,\flat\,} \ + \ \phi_{\,\,\flat}\,)\,. \leqno (1.18)
$$
Note that $\,{\bf I}_{\,\cal R}\,$
 depends only on the $(\,n\,+\,1\,) \cdot \flat$   bubble parameters\,:
 $$
 {\bf P}_{\,(\,\flat)\,} \ = \ (\,\lambda_{\,1}\,, \ \cdot \cdot \cdot\,, \ \lambda_{\,\,\flat}\,; \ \xi_{\,1}\,, \ \cdot \cdot \cdot\,, \ \xi_{\,\flat}\,)\,, \leqno (1.19)
 $$
 $$
\lambda_{\,1}\,, \ \cdot \cdot \cdot\,, \ \lambda_{\,\,\flat} \ \in \ \R^+\,, \  \ \ \ \ {\mbox{and}} \ \ \ \
\xi_{\,1}\,, \ \cdot \cdot \cdot\,, \ \xi_{\,\flat}\ \ \in \ \R^n\,. \leqno {\mbox{where}}
$$
{\it Provided that $\,\lambda_{\,\,\flat}$ is small enough}\, [\,see Lemma 3.1 for detail\ ]\,,\, {\it if}\, $\,W_{\,\flat\,}\,+\,\phi_{\,\,\flat}\,$ {\it is a critical point of the reduced functional $\,{\bf I}_{\,\cal R}\,$,\,  that is}
$$
 {{\ \partial \,{\bf I}_{\,\cal R}\ }\over {\partial\, \lambda_{\ \!l} }}\,\bigg\vert_{\  {\bf P}_{\,(\,\flat)\,} } \ = \ {{\ \partial \,{\bf I}_{\,\cal R}\  }\over {\partial \,\xi_{{\,l \,|_{\,j}}} }} \,\bigg\vert_{\  {\bf P}_{\,(\,\flat)\,} } \ = \ 0 \ \ \ \ \ {\it{for \ \ all }}  \ \ \ l\  = \ 1\,,\, \cdot \cdot \cdot\,, \   \flat\,, \ \ \ \ j \ = \ 1\,, \ 2\,, \ \cdot \cdot \cdot \,, \ n\,, \leqno (1.20)
$$
{\it then indeed\,}  $\,W_{\,\flat\,} \ + \ \phi_{\,\,\flat}\,$ {\it is  a critical point of}\, (\,{\it the full functional\,}\,) $\,{\bf I} \,$,\, {\it that is\,,\,}
$$
{\bf I}\,' \, (\ W_{\,\flat\,} \ + \ \phi_{\,\,\flat}) \ = \ 0\,. \leqno (1.21)
$$
In view of finding a set of solution  $\,{\bf P}_{\,(\,\flat\,)}\,$ [\,cf. (1.19)\,] to the finite number of equations in (\,1.20\,)\,,\, the goal of this paper is to expression the key information in
$$
 {\bf I}_{\,\cal R}\,, \ \ \ \ \ \ \
{{\ \partial \,{\bf I}_{\,\cal R}\ }\over {\partial\, \lambda_{\ \!l} }}  \ \ \ \ \ \ {\mbox{and}} \ \ \ \ \ \ \  {{\ \partial \,{\bf I}_{\,\cal R} \ }\over {\,\partial \,\xi_{{\,l \,|_{\,j}}}\, }}
$$
for $\,\,l\  = \ 1\,,\, \cdot \cdot \cdot\,, \   \flat\,, \ \ \ \ j \ = \ 1\,, \ 2\,, \ \cdot \cdot \cdot \,, \ n\,.
$

\vspace{0.2in}

{\large{\bf \S\,1\,c\,.}}\ \   {\it Key geometric parameters $\,\gamma\,,\,$ $\,\kappa\,$,\, $\nu\,$ and $\,\,\ell\,.$} \ \   We first assume that the arrangement of the bubbles at $\ \xi_{\,1}\,, \ \cdot \cdot \cdot\,, \ \xi_{\,\flat}$\, satisfies
 $$
\max_{\,l} \ \left\{  \ \sum_{k\ \not= \ l} \  \left( {1\over { \ {\bf d}_{\ k\,,\ l} \    }} \right)  \ \right\}  \ \le  \  {\bar\lambda}_{\,\,\flat}^{\,\gamma}   \,, \ \ \ \ \  \ {\mbox{where}} \  \ \gamma\ \in \ \left(\ {1\over 2}\,,\ \,1\,\right) \ \ {\mbox{is \ \ fixed}}\ . \leqno (1.22)
 $$
 Here $\,1 \ \le \ l \ \le \ \flat\,.\,$
Intuitively, it means that the points  $\,\xi_{\,1}\,, \ \cdot \cdot \cdot\,,\ \xi_{\,\flat\,-\,1}\,$ and $\,\xi_{\,\flat}$\, are well\,-\,separated in the quasi\,-\,hyperbolic distance [\ see  (1.8) and (\,1.9\,) for the definition of
 $\, {\bar\lambda}_{\,\,\flat}\,$ and $\,
 {\bf d}_{\,k\,,\ l}  $\ ,\, respectively\,]\,.\\[0.15in]
\noindent {\it Remark}\, 1.23\,.\, \ \   As our main results hold for a small perturbation of the index $\,\gamma\,,\,$ we take the approach on not including a positive constant $\,``\ C\,"\,$ in front of $\,{\bar\lambda}_{\,\,\flat}^{\,\gamma}\,$ in (1.22)\,.\, This is specially clear because  $\,{\bar\lambda}_{\,\,\flat} \ = \ o\,(\,1\,)\,,\,$  as $\, C \cdot  {\bar\lambda}_{\,\,\flat}^{\,\gamma} \ = \  [ \ C \cdot  {\bar\lambda}_{\,\,\flat}^{\,\epsilon}\, ] \cdot [ \  {\bar\lambda}_{\,\,\flat}^{\,\gamma \ - \ \epsilon} \, ] \ \le \ {\bar\lambda}_{\,\,\flat}^{\,\gamma \ - \ \epsilon}\,$\  when $\,{\bar\lambda}_{\,\,\flat}\,$ is small enough\,.\, Here $\,\epsilon\,$ is a small positive number\,.\, (\,1.28\,)\,,\, (\,1.29\,) and (\,1.32\,) below are treated likewise.

\newpage

\hspace*{0.5in}\emph{}Although portions of our discussion can be fit in general geometric situation, for clarity sake we first approach the situation where the (\,critical\,) set\,:
$$
{\cal {C}}{\it rt} \ := \ \{\  y \,\in\,\R^n\ \ \big\vert \ \ y \ \,{\mbox{is \ a \ critical \ point \ of }} \ K\,, \ \ {\mbox{that \ \ is}}\,, \ \ \btd\,K\,(\,y\,) \ = \ 0\ \}
$$
is non\,-\,isolated, and indeed it contains a hypersurface $\,{\cal H}\,$ (\,of dimension equal to $\,n\,-\,1\,$)\,.\, Moreover, we assume that in a neighborhood of $\,{\cal H}\,$,\, there is no other critical point of $\,K\,$ except those in $\,{\cal H}\,$.\, That is\,,\, there exists an open set $\,\cal O\,$  such that
$$
{\cal H} \,\subset\,{\cal {C}}{\it rt}\ \ \ \ {\mbox{and}} \ \ \ {\cal H} \,\subset\, {\cal O}\ , \ \ \ \ {\mbox{with}} \ \ \  \  [\ \cal O\,\setminus \,{\cal H} \ ]\ \cap \ {\cal {C}}{\it rt} \ = \ \emptyset\,. \leqno (1.24)
$$
We are restricted to the case that $\,{\cal O}\,$ can be chosen to be ``\,staying close\," to $\,{\cal H}\,$ so that the ``\,projection\," of $\,y\,$ unto $\,{\cal H}\,$,\, denoted by $\,{\bf p}_{\,y}\,\in\,{\cal H}\,,\,$ is uniquely defined\,:
$$
\exists\,! \ \,{\bf p}_{\,y} \ \in  \,{\cal H} \ \ \ {\mbox{so \ \  that \ \ Dist}}\  (\,y\,,\ {\cal H}\,) \ = \ \Vert  \,y \ - \ {\bf p}_{\,y}\ \Vert \ \ \ \ (\,y \,\in\, \cal O\,)\,.\, \leqno (1.25)
$$
This is clearly the case if ${\,\cal H}\,\subset \,\R^n$ is compact\,.\, Now we assume that \\[0.1in]
(\,1.26\,)
$$
(\,{\tilde c}_n\,K\,)\,(\,y\,) \,=\,n\,(\,n\,-\,2\,) \ - \ C\,(\,{\bf p}_{\,y}\,)\cdot \Vert\,y \,-\,{\bf p}_{\,y}\,\Vert^{\,\ell} \ + \  {\bf R}_{\,\,\ell \  + \ 1}(\,y\,) \ \  \mfor \ \ y \ \in \ {\cal O}\,.
$$
Here $\, C\,(\,{\bf p}_{\,y}\,)\,$ is a positive number depending on $\,\,{\bf p}_{\,y}\,$,\, and
$$
 {\bf R}_{\,\, \ell \  + \ 1}(\,y\,)  \ \le \ {\bar C}_{R} \,\Vert\,y \,-\,{\bf p}_{\,y}\,\Vert^{\,\ell\ +\,1}\ \ \  \mfor \ \ y \ \in \ {\cal O}\,.
$$
The integer $\,\ell\,$ satisfies
$$\,\ell \ \in \ [\ 2\,, \ n \ - \ 2\ )\,. \leqno (1.27)$$
[\ One can also replace the index $\,\ell \ + \ 1\,$ in the remainder by  $\,\ell \ + \ \theta\,$,\, as in (\,1.3\,)\,.\ ]
 \bk
With reference to (\,1.5\,) and (\,1.10\,)\,,\, the center of each bubble\,,\, namely\,,\, $\xi_{\,l}$\,,\, is assumed to be close to $\,{\cal H}\,$,\, but not exactly on $\,{\cal H}\,$.\, Precisely,
$$
0 \ < \   \eta_{\,\,l} \ := \   {\mbox{dist}}\,(\,\xi_{\,l}\,,\ {\cal H}\,) \ \le  \   {\bar\lambda}_{\,\,\flat}^{\,1\ + \ \kappa}\ , \ \ \ \ \ \ \ \ \   {\mbox{where}} \ \ \kappa \, \in  \, (\,0\,,\   1\,) \ \ {\mbox{is \ \ fixed }}\,.\leqno (1.28)
$$
Often when we perform various integrations, we cut\,-\,off via the ball $\,B_{\,\xi_{\,l}}\,(\,\rho_\nu\,)\,$ with radius satisfying
$$
\rho_\nu \ = \ {\bar\lambda}_{\,\,\flat}^{\,\nu} \ \ \ \ {\mbox{where}} \ \ \ \ \rho_\nu \ = \ o_{\,+}\,(\,1\,) \cdot \min \ \{  \ \Vert \,\xi_{\,l} \ - \ \xi_{\,k}\, \Vert \ \ | \ \ 1 \ \le \ l\ \not= \ k\ \le \ \flat \ \}\ . \leqno (1.29)
$$
We always assume that $\,{\bar\lambda}_{\,\,\flat}\,$  is small enough so that $\,B_{\,\xi_{\,l}}\,(\,\rho_\nu\,)\,\subset \, \cal O$\, for all $\, l \ = \ 1\,,\, \ \cdot \cdot \cdot\,,\, \ \flat\,.$\,
Once $\,\gamma \ (\ >  \ {1\over 2}\ )\ $ is fixed\,,\, we adjust $\,\nu\,$  so that
$$
(\,1 \ > \ ) \ \gamma \ >  \ \nu\ >  \ {1\over 2}  \ . \leqno (1.30)
$$
It follows that
$$
\gamma \ > \ {1\over 2} \ \ \ \ {\mbox{and}} \ \ \ \ \nu \ > \ {1\over 2}  \ \ \Longrightarrow \ \ \gamma \ + \ \nu \ >  \ 1 \ \ \Longrightarrow \ \ \gamma \ > \ 1 \ - \ \nu\,.\leqno (1.31)
$$
Finally, the number of bubbles $\,\flat\,$ is bounded from above by
$$
2 \ \le \ \flat \ \le\  {1\over { { {\bar\lambda}_{\,\,\flat}   }^{\,\sigma} }} \ \ ( \ \to \ \infty \ \ {\mbox{as}} \ \ {\bar\lambda}_{\,\,\flat} \ \to \ 0^+\,)\, . \leqno (1.32)
$$
%


\newpage


{\large{\bf \S\,1\,d.}}\ \ \    {\it Main information within the reduced functional and its first derivatives.  }\\[0.15in]
{\bf Main\,  Theorem\, 1.33\,.\,} \ \ {\it For  $\,n \ \ge \ 6$\,,\, under the conditions in}\, (\,1.4\,)\,,  (\,1.8\,)\,, (\,1.22\,)\,, (\,1.24\,)\,--\,(\,1.28\,)\,,\, (\,1.30\,)\, {\it and} \,(\,1.32\,)\,,\, {\it assume also that the geometric parameters $\,\gamma\,,\,$  $\nu\,$ and $\,\,\ell\,$ satisfy}
$$
\mbox{Min}  \, \left\{ \     {{n}\over 2} \,\cdot\,(\,1\, - \ \nu\,)\,, \  \ \ \  \ell    \cdot  \nu \  \right\}  \ > \ \left( {{\,n\ + \ 10}\over 8}\  \right) \cdot \sigma\ \ ( \ \ge \ 2\cdot \sigma\,)\ , \leqno (\,1.34\,)
$$
{\it and}

\vspace*{-0.3in}

$$
\gamma \ > \ {\sigma\over {\,n \ - \ 2\,}} \ . \leqno (\,1.35\,)
$$

{\it We can find a small positive number  $\,{\underline\lambda}_{\ \epsilon}\,$  so that if $$\ {\bar\lambda}_{\,\,\flat} \ \le \ {\underline\lambda}_{\ \epsilon}\,,\,$$ then we have the  following}\, ``\,{\it expansions}\," (\,1.36\,)\,\,--\,(\,1.38\,)\,.
\begin{eqnarray*}
(\,1.36\,)& \ &  {\bf I}_{\,{\cal R}} \,(\,\lambda_{\,1}\,, \ \cdot \cdot\,,\ \lambda_{\ \!l}\,,\  \cdot\,, \ \lambda_{\,\,\flat}\,;\, \ \xi_{\,1}\,, \ \cdot \cdot\,, \ \xi_{\,l}\,, \ \cdot \cdot\,, \ \xi_{\,\flat}\,)\\[0.15in]
 & = &{1\over 2} \,n\, (\,n\,-\,2\,)\cdot V\,(\,n\,) \cdot \flat  \ \ + \\[0.2in]
  & \ &  \!\!\!\!\!\!\!\!\!\!\!\!\!+ \ \,  {\hat C}_{\,0\,,\ 1}  \sum_{l \ =\,1}^\flat C\,(\,{\bf p}_{\,\,l}\,) \, \lambda_{\,\,l}^\ell \ - \    {\hat C}_{\,0\,,\ 2}   \cdot \left\{ \   \sum_{l\,=\,1}^\flat \ \sum_{k \,\not=\,l} \,  \left( \  {1\over { {\bf d}_{\,\,l\,,\ k} }} \  \right)^{\!n\,-\ 2} \ \right\}   \cdot\bigg[\ 1  + \, O\,\left(\ {\bar\lambda}_{\,\,\flat}^{\, 2 \,(\,1 \ - \ \nu\,) }\ \right) \ \bigg] \, + \\[0.05in]
 & \ &\hspace*{-0.85in}(\,{\mbox{curvature \ \ contribution}} \ \ \uparrow \ )\hspace*{1.7in} ( \ \ \uparrow \ \ {\mbox{bubble \ \ interactions}} \ )\\[0.1in]
  & \ &   \ \ \ \  \ \ \ \ \ \ \ \ \ \ \ \ \ \ \ \ \ \ \ \ \ +   \ \ {\hat C}_{\,0\,,\ 3}  \cdot  \sum_{l\,=\,1}^\flat  \ C\,(\,{\bf p}_{\,\,l}\,) \cdot \left[\ \lambda_{\,\,l}^{\ell\ -\ 2} \cdot \, \eta^{\,2}_{\,\,l} \ \right]   \ \  + \ \ {\bf E}_{\,o} \ . \\[0.11in]
(\,1.37\,)
& \ & \left( \lambda_{\ \!l}\cdot {\partial\over {\partial\, \lambda_{\ \!l}}}\right)   \,(\,\lambda_{\,1}\,, \ \cdot \cdot\,,\ \lambda_{\ \!l}\,,\  \cdot\,, \ \lambda_{\,\,\flat}\,; \ \, \xi_{\,1}\,, \ \cdot \cdot\,, \ \,\xi_{\,l}\,, \ \cdot \cdot\,, \ \xi_{\,\flat}\,)  \ \ \ \ \ \ \ \ \ \ \ \ \  [ \ 1 \ \le \ l \ \le \ \flat\ ] \\[0.1in]& = &    {\hat C}_{\,1\,,\ 1}  \cdot C\,(\,{\bf p}_{\,l}\,) \cdot \lambda^{\ell}_{\,l}  \ -\ {\hat C}_{\,1\,,\ 2}  \left\{ \  \sum_{k\ \not= \ l}   \left( \  {1\over { {\bf d}_{\,\,l\,,\ k} }} \  \right)^{\!n\,-\ 2} \ \right\}
   \cdot\left[\,1 \ + \ O\,\left(\ {\bar\lambda}_{\,\,\flat}^{\, 2 \,(\,1 \ - \ \nu\,) }\ \right) \ \right]  \ + \ \\[0.2in]
  & \ & \ \ \ \ \ \ \ \ \ \ \ \ \ \ \  \ \ \ \ \ \ \  \ \  \ \     + \ \  {\hat C}_{\,1\,,\ 3}  \cdot C\,(\ {\bf p}_{\,l}\ ) \cdot \lambda^{\,\ell\ -\ 2}_{\,l} \cdot  \eta^2_{\,l}    \ \ \  + \ \ {\bf E}_{\,\lambda_{\ \!l}}\ .  \\[0.2in]
(\,1.38\,)
& \ & (\,  \lambda_{\ \!l} \cdot  \btd_{\xi_{\,l} }\,)\ {\bf I}_{\,\cal R}\,(\,\lambda_{\,1}\,, \ \cdot \cdot\,,\ \lambda_{\ \!l}\,,\  \cdot\,, \ \lambda_{\,\,\flat}\,; \ \,\xi_{\,1}\,, \ \cdot \cdot\,, \ \xi_{\,l}\,, \ \cdot \cdot\,, \ \xi_{\,\flat}\,) \ \ \ \ \ \ \ \ \ \ \ \ \  [ \ 1 \ \le \ l \ \le \ \flat\ ] \\[0.1in]
& = &     {\hat C}_{\,2\,,\ 2} \left\{  \   \sum_{k\ \not= \ l}   \left( \  {1\over { {\bf d}_{\,\,l\,,\ k} }} \  \right)^{\!\!n} \!\cdot \left( {{\,\xi_{{\,l} }\ - \ \xi_{k}\,}\over {\lambda_k}} \right)  \  \right\}    \cdot\left[\,1\ + \ O\,\left(\ {\bar\lambda}_{\,\,\flat}^{\, 2 \,(\,1 \ - \ \nu\,) }\ \right) \ \right]  \\[0.2in]
& \ &   \ \ \  \ \ \ \ \ \ \ \ \ \ \ \ \ \ \ \ \ \ \ \ \ \ \ \ \ \ \ \ + \  \left[\ {\hat C}_{\,2\,,\ 3}\, \cdot   C\, (\ {\bf p}_{\,l}\ ) \cdot \lambda_{\ \!l}^\ell \cdot\left( {{  \eta_{\,l}  }\over {  \lambda_{\ \!l} }} \right) \ {\bf n}_{\,l} \  \right] \  \ + \ {\vec{\,\bf E}}_{\,\,\xi_{\,l}}  \ .
 \end{eqnarray*}

  \newpage

{\it Here} ${\mbox{ \ Dist}}\  (\ \xi_{\,l}\,,\ {\cal H}\,) \ = \ \Vert  \,\xi_{\,l} \ - \  {\bf p}_{\,l}\ \Vert \ > \ 0$\,\, {\it{and}}\  $\displaystyle{ \ \
  {\bf n}_{\,\,l} \ = \  {{ \xi_{\,l} \,-\, {\bf p}_{\,l} }\over { \Vert \  \xi_{\,l} \,-\, {\bf p}_{\,l} \,\Vert }}}$ \ ,\    {\it{where}} $\ \ {\bf p}_{\,l}\,\in\,{\cal H}\,$ {\it for}\\[0.1in]
   \noindent$\, 1 \ \le \ l \ \le  \ \flat\,.$\, {\it See} \,{\bf \S\,1\,e}\, {\it regarding the constants, and}  \,{\bf \S\,1\,f\ }\, {\it for the lower order} (\,{\it{error}}\,) {\it terms.}

 \vspace*{0.5in}

{\large{\bf \S\,1\,e.}}\ \ \    {\it The constants.} \ \  In (\,1.36\,)\,,\, (\,1.37\,) and (\,1.38\,)\,,\,
\begin{eqnarray*}
  V(\,n\,) & = & \int_{\R^n} \left(\ {1\over {1\ + \ \Vert\,Y\,\Vert^2 }}\ \right)^{\!n}\, d\,Y\ ,\\[0.2in]
  {\hat C}_{\,0\,, \ 1} & = & \left(\ {{\,n\,-\ 2\,}\over {2\,n}}\,\right)\,\cdot\,\int_{\R^n} \left(\ {1\over { 1 \ + \ \Vert\ Y\,\Vert^2 }}\ \right)^{\!\!n} \,d\,Y\ = \  \left(\ {{\,n\,-\ 2\,}\over {2\,n}}\,\right)\,\cdot\, V\,(\,n\,)\ ,\\[0.2in]
{\hat C}_{\,1\,,\ 1}  & = &  \left(\ {{\,n\,-\ 2\,}\over {2}}\,\right)\,\cdot\,\int_{\R^n} |\,Y_n\,|^{\,\ell} \cdot  \left(\ {1\over { 1 \ + \ \Vert\ Y\,\Vert^2 }}\ \right)^{\!\!n} \!\cdot \left(\, {{\ \Vert\ Y\,\Vert^2 \ - \ 1 }\over { \ \  \Vert\ Y\,\Vert^2 \  + \ 1  }}\, \right)\, d\,Y \\[0.2in]  (\!\!& = & \ell \cdot {\hat C}_{\,0\,,\ 1}\,, \ \ {\mbox{see \ \ {\bf \S\,6\,c}}}  \ ) \\[0.2in]
    {\hat C}_{\,0\,, \ 2} & = &  \left( \, {{\,n \ - \ 2\,}\over 2} \, \right) \cdot \omega_n \ \ \ \ \ \ \ \ \ (\ \omega_n \ \ {\mbox{is \ \ the \ \ volume \ \ of  \ \ the \ \ unit \ \ sphere \ \ in}} \ \ \R^n \ )\  , \\[0.3in]
      {\hat C}_{\,1\,, \ 2} & = &    {{\,(\,n \ - \ 2\,)^{\,2}\,}\over 2}  \cdot \omega_n \ \left(\ = \    2 \cdot \left(\,{{\,n\,-\,2\,}\over 2}\,\right) \cdot {\hat C}_{\,0\,,\ 2} \ \right)\ , \\[0.2in]
    {\hat C}_{\,2\,,\ 2}  & = & (\,n\,-\,2\,)^{\,2} \cdot  \omega_n \ \left( \  = \  2 \cdot (\,n\,-\,2\,) \cdot {\hat C}_{\,0\,,\ 2} \ \right)\,,\\[0.3in]
  {\hat C}_{\,0\,, \ 3} & = & \left(\ {{n\,-\,2}\over {2n}}\,\right)\,\cdot\, {{\ell\,\cdot\, (\,\ell\ - \ 1\,)}\over 2}\,\cdot\,\int_{\R^n} |\,Y_n\,|^{\,\ell\ - \ 2}\,\cdot\, \left(\ {1\over { 1 \ + \ \Vert\ Y\,\Vert^2 }}\ \right)^{\!\!n} \, d\,Y\ ,\\[0.2in]
  {\hat C}_{\,1\,, \ 3} & = &   \left({{n\,-\,2}\over {2}}\,\right)   \cdot {{\ell \cdot (\,\ell\, - \, 1\,)}\over 2} \cdot \int_{\R^n} |\,Y_n\,|^{\,\,\ell\, - \, 2} \cdot \left(\, {1\over { 1 \, + \, \Vert\,Y\,\Vert^2 }}\right)^{\!\!n}  \!\cdot \left(\, {{\Vert\ Y\,\Vert^2 \, - \, 1 }\over { \ \  \Vert\ Y\,\Vert^2 \, + \, 1 \ }}\, \right) \, d\,{\,Y}\ , \\[0.2in]
  (\!\!& = & (\,\ell\,-\,2\,) \cdot {\hat C}_{\,0\,,\ 3}\,, \ \ {\mbox{see \ \ {\bf \S\,6\,c}}}  \ )  \\[0.2in]
 {\hat C}_{\,2\,, \ 3}  & = &  \left(\,n\,-\,2\,\right)   \cdot \ell \cdot \int_{\R^n} |\,Y_n\,|^{\,\,\ell } \cdot \left( {1\over { 1 \ + \ \Vert\,Y\,\Vert^{\,2} }}\right)^{\!\!n\,+\,1}  d\,Y \ ( \,= \ 2 \cdot {\hat C}_{\,0\,,\ 3}\,, \ \ {\mbox{see \ \ {\bf \S\,6\,c}}}  \ ) \ .\\
\end{eqnarray*}

\newpage

{\large{\bf \S\,1\,f.}}\ \
 {\it The lower order terms.} \ \ In (\,1.38\,)\,,\,
 $$
 {\vec{\,\bf E}}_{\,\,\xi_{\,\,l}}  \ = \  \left( \ E_{{\,\xi_{\,l}}_{\,|_1}}\,, \ E_{{\,\xi_{\,l}}_{\,|_1}}\,, \ \cdot \cdot \cdot\, \ E_{{\,\xi_{\,n}}_{\,|_1}}\,\,\right) \ . \leqno (\,1.39\,)
  $$
The error terms in (\,1.36\,)\,,\, (\,1.37\,) and (\,1.38\,)\,, namely\,,\,
 $$
 {\bf E}_{\,o} \,, \ \ \  {\bf E}_{\,\lambda_{\ \!l}}\,, \ \ \ E_{{\,\xi_{\,l}}_{\,|_1}}\,, \ \ E_{{\,\xi_{\,l}}_{\,|_2}}\,, \ \cdot \cdot \cdot\, \ \ {\mbox{and}} \ \ \  E_{{\,\xi_{\,l}}_{\,|_n}}  \ \ \ \ \ ( \ l \ = \ 1 \,, \ \cdot \cdot \cdot\,, \ \flat \ )\,,
 $$
 are of the following order
$$
  O \,
 \left( \ {\bar\lambda}_{\,\,\flat}^{\,\ell \ + \ \mu_{\,\vert_{\cal K}} } \right)  \    + \  O\,\left(\, {\bar\lambda}_{\,\,\flat}^{\,{\varpi} \ + \  \mu_{\,\vert_{\,\Phi}}  }\, \right)  \ + \ O\left( \  {\bar\lambda}_{\,\,\flat}^{ \,n \,\gamma \, - \ \,\sigma}\ \right) \ + \     O\left(\   {\bar\lambda}_{\,\,\flat}^{ \  n\cdot\,(\,1\ - \ \nu\,) \ - \ o_{\,+}\,(\,1\,)} \ \right)   \ . \leqno  (\,1.40\,)
   $$
Here
 \begin{eqnarray*}
 (\,1.41\,) \ \  \mu_{\,\vert_{\cal K}} \!\!& = & \!\! \mbox{Min} \ \ \bigg\{ \ (\,n\,-\ \ell\,)\, \cdot\, (\,1 \ - \ \nu)\,, \ \ 1 \ -\ o_{\,{\bar\lambda}_{\ \flat}}\,(\,1\,)\,, \ \ (\,n\ -\ 2\,)\,\cdot\, [\ (\,\gamma \ + \ \nu \,)\ - \  1\,]\,, \\[0.2in]
   & \ &  \ \ \ \ \ \ \ \ \ \ \ \ \ \ \ \ \ \ \ \ \ \ \ \ \ \ \ \  \ \ \ \ \ \ \ \ \ \ \  \ \ \  \ \ \ \ \ \  \ \ \ \ \ \ \ \ \ \ \ \ (\,\ell \ + \ 1\,)\cdot  (\,2\,\nu \ - \ 1\,) \ ,  \ \  \ \ 3\,\kappa\ \bigg\}\ , \ \ \
   \end{eqnarray*}
\begin{eqnarray*}
& \ & \\[-0.3in]
(\,1.42\,) \ \ \ \ \ \  \  \ \ \ \ \ \ \ \   \mu_{\,\vert_{\,\Phi}}  & = & \mbox{Min} \ \ \left\{ \ {{4\over {n\,-\,2}} \,\cdot\ \varpi}  \ , \  \ \ \ 4\,\gamma \ -\ o_{\,+}\,(\,1\,) \ , \ \ \ \  \varpi \ - \ \sigma\,, \ \ \ \ \ell \ \  \right\}\ , \ \ \ \ \ \  \ \ \ \ \ \ \ \ \  \ \  \ \ \ \ \ \ \ \ \ \ \ \\[0.2in]
(\,1.43\,)  \ \ \ \ \ \ \ \ \ \ \ \ \ \ \ \ \  \varpi &  = & \mbox{Min}  \, \left\{ \     {{n}\over 2} \,\cdot\,(\,1\, - \ \nu\,)\,, \  \ \ \  \ell    \cdot  \nu \  \right\}\ - \ {1\over 2} \cdot \sigma \ \,-\ o_{\,+}\,(\,1\,)\ .
\end{eqnarray*}
Refer to\, {\bf \S\,1\,h} (\,$\bullet_{\,2}$\,) for the way we use
  $\,o_{\,+}\,(\,1\,)$\,  and $\,\,o_{\,{\bar\lambda}_{\,\,\flat}}\,(\,1\,)\,$.

%
%

\vspace*{0.3in}

{\large{\bf \S\,1\,g\,.}}\ \  {\it Layout of the article.}\ \ In {\bf \S\,2}\,,\, we discuss how to find the small solution to the equation in the perpendicular direction (\,1.17\,)\,.\, {\bf \S\,3} is devoted to the finite dimnesion reduced functional  and its critical points. We present the arguments towards the expansion (\,1.37\,) in {\bf \S\,4}\,,\, (\,1.38\,) in {\bf \S\,5}\,,\,  and (\,last but not least\,)\, (\,1.36\,) in {\bf \S\,6}\,.


\vspace*{0.3in}

{\large{\bf \S\,1\,h\,.}}\ \  {\it General conditions, assumptions and conventions.}\ \ Throughout this work,\,
$$
\ \ S^n\ =\ \bigg\{ \,x\,=\,(\,x_1\,,\ \cdots\ ,\ x_{n +1}\ )\ \in\ \R^{n+1}\ \ \bigg|\ \ x_1^2\,+\,\cdots\,+\,x_{n+1}^2\ =\ 1\,\bigg\} \ \ \ \ \ (\,n\,\ge\,3)\,, \leqno (\,1.44\,)
$$
with the (\,starndard\,) induced metric $\,g_{\,1}\,$\,. $\,\Delta_{\,1}\,$ is the Laplace\,-\,Beltrami operator associated with $\,g_1\,$ on $\,S^n\,$. Likewise, $\,\Delta\,$ is the Laplace\,-\,Beltrami operator associated with Euclidean metric $\,g_o\,$ on $\,\R^n\,$,\, with coordinates $\,y \ = \  (y_1\,,\ \cdots\ ,\, y_{n})\,\in \,\R^n\,.\,$ Denote by $\,\omega_n\,$ the volume of $\,S^{\,n\,-\,1} \subset \ \R^n\,.$\, Moreover, the norm $\,\Vert\ \ \Vert\,$ and the inner product $\,\langle\ \,\,,\ \ \, \rangle\,$ are defined via Euclidean metric $\,g_o\,$ on $\,\R^n\,$. \,\\[0.1in]

\newpage

$\bullet_{\ 1}$\ \, In this article, the constant $\,{\tilde c}_n\,=\,{{n\,-\,2}\over{\ 4\,(\,n\,-\,1\,)\ }}$\ \,.\, We  observe the practice on  using `$\,C\,$'\,,\, possibly with sub\,-\,indices, to denote various positive constants\,,\, which may be rendered {\it differently\,} from line to line according to contents.  These constants are independent on geometric parameters
$$\,\gamma\,, \ \ \sigma\,, \ \ \nu \ \ \ {\mbox{and}} \ \ \  \sigma\,,$$ satisfying the conditions (\,1.4\,)\,,  (\,1.8\,)\,, (\,1.22\,)\,, (\,1.24\,)\,--\,(\,1.28\,)\,,  \,(\,1.30\,)\, and  \,(\,1.32\,)\, \,.\, Some of these constants may depend on  the compact set $\,{\cal C}_{\,[]}\,$ in (\,1.4\,)\,,\, the constant $\,{\bar C}_{\,1}\, $ in (\,1.8\,)\,,\, and\,/\,or the constant $\,{\bar C}_{\,R}$\, involved in the remainder in (\,1.26\,) concerning $\,K\,$.\,   $\,{\cal C}_{\,[]}\,$,\,     $\,{\bar C}_{\,1}\, $ and $\,{\bar C}_{\,R}$\ are supposed to be given and fixed\,.\, {\it Whilst we use `$\ {\bar C}\,$' \,or\, `$\ {\hat C}\,$'\,,\, possibly with sub\,-\,indices, to denote\,} fixed\, {\it positive constants which always keep the same values as they are first defined\,}.\,\\[0.1in]
$\bullet_{\ 2}$\ \,  \ \ We record that an expression is of $\,O\,(\,1\,)\,,\,$ if the absolute value of the expression is bounded from above by a positive constant\,.\, $o_{\,+}\,(\,1\,)$\, denote a fixed positive number that can be chosen to be arbitrarily  small\,,\, whereas the positive number $\,o_{\,{\bar\lambda}_{\,\,\flat}}(\,1\,)\ \to \ 0^{\,+}$\  \,as $\ {\bar\lambda}_{\,\,\flat} \ \to \ 0^+\,$.\,  \\[0.1in]
$\bullet_{\ 3}$\ \, Denote by $\,B_{\,y}\,(\,r\,)\,$ the open ball in $\,(\,\R^n\,,\ g_{\,o})\,$ with center at $\,y\,$ and radius $\,r\,>\,0\,$, $\,\partial B_{\,y}\,(\,r\,)\,$ its boundary\,.\, Whenever there is no risk of misunderstanding, we suppress `$\,d\,y\,$'\, from the integration expressions on domains in $\,\R^n\,$.\\[0.1in]
$\bullet_{\ 4}$\ \, $\,f_{\,+}\,$ denotes the positive part of $\,f\,$.\,  In this article, we work on the Hilbert space given by
$$
{\cal D}^{\,1,\,\,2}\,=\,{\cal D}^{\,1,\,\,2}\,(\,\R^n)\,:=\,\left\{\,f\in L^{\,{{2n}\over {\,n\,-\,2\,}}}\,(\,\R^n)\:\bigcap\: W^{\,1,\, 2}_{\small\mbox{loc}}\, (\,\R^n)\ \,\bigg\vert\,\ \int_{\R^n}\langle\ \btd\,f\,,\,\btd\,f\ \rangle\,<\,\infty\,\right\}\,,\leqno (\,1.47\,)
$$
with the inner product   defined by
$$
\langle\ f\,,\ \psi\,\rangle_{\,\btd}\ :=\ \int_{\R^n}\langle\ \btd\,f\,,\,\btd\,\psi\ \rangle\ \ \ \ \mfor\ f\,,\ \psi\ \in\ {\cal D}^{\,1,\,\,2}  \ .  \leqno (\,1.48\,)
$$

\vspace*{0.2in}

{\large{\bf \S\,1\,i\,.}}\ \   \textbf{e}\,\,-{\it Appendix.}\ \ \ Many  of the preparatory lemmas and  estimates of errors are situational modifications of well\,-\,established arguments, mainly found in \cite{Wei-Yan}\,.\, Necessarily some of these arguments have to be modified to fit into the present generic situation. These steps should be familiar to people working in the area. In order to keep this article within manageable size\,,\, and to present these  lengthy and tedious calculations clearly and precisely,  we choose to put those details in {\bf e}\,-\,Appendix\,,\, which is   available at
arXiv.org $\rightarrow$ math.AP $\rightarrow$ Analysis of PDEs (\,{\tt{https://arxiv.org/pdf/????.??????}}\,)\,.\, In this way\,,\,  interested readers can refer to it conveniently.

\vspace*{0.2in}

{\large{\bf \S\,1\,j\,.}}\ \ {\it Acknowledgement.} \ \ The author is grateful to F. Pacard his insight into  the conjecture by J. Wei and S. Yan.

 \newpage

{\bf \large {\bf \S\,2\,. Solving the  equation in the $\,\perp\,$-\,direction\,-\, the case of multiple}}\\[0.1in]
\hspace*{0.4in} {\bf \large {\bf    bubbles\,.\,}}  \\[0.2in]
We now describe equation (\,1.17\,) properly\,.\, We seek $\,\phi_{\,\,\flat} \,\in\,{\cal D}^{\,1,\,\,2}_{\,\flat_\perp}\ \,$ so that [\ in a weak sense first\,,\, see (\,1.14\,)\ ]\\[0.1in]
(2.1)

\vspace*{-0.3in}

$$
\Delta \,(\,W_{\,\flat} \,+\,\phi_{\,\,\flat} ) \ + \  (\,{\tilde c}_n\, K)\,(W_{\,\flat} \,+\,\phi_{\,\,\flat} )_+^{{n\,+\,2}\over {\,n\,-\ 2\,}} = \ \sum_{i\,=\,1}^\flat a_{\,l} \,[\,(\,\lambda_{\ \!l} \cdot \partial_{\lambda_{\ \!l}}\,)\,V_{\,l}\,] \ \ + \ \sum_{l\,=\,1}^\flat \left( \ \sum_{j\,=\,1}^n b_{\,l\,,\,j} \ [\,(\,\lambda_{\ \!l} \cdot \partial_{\xi_{\,l_{|_{\,j}}}}\,)\,V_{\,l} \ ]\, \right) \!.
$$

Before we proceed, let us take note that [\ refer to (\,1.16\,)\ ] \,,\,
\begin{eqnarray*}
& \ & {\cal P}_{\,\flat}\, \left( \   {\mbox{right \ \ hand \ \ side \ \ of \ \ (\,2.1\,)}} \ \right) \\[0.2in]& = &
{\cal P}_{\,\flat}\, \left( \  \sum_{i\,=\,1}^\flat a_{\,l} \,[\,(\,\lambda_{\ \!l} \cdot \partial_{\lambda_{\ \!l}}\,)\,V_{\,l}\,] \ + \ \sum_{i\,=\,1}^\flat \left\{ \ \sum_{j\,=\,1}^n b_{\,l\,,\,j} \ [\,(\,\lambda_{\ \!l} \cdot \partial_{\xi_{\,l_{|_{\,j}}}}\,)\,V_{\,l}\ ]\ \right\}\ \  \right)
  \ \ =\ 0\,.\,
\end{eqnarray*}
With respect to the terms in the right hand side of (\,2.1\,)\,,\, for $\,\psi \ \in \ {\cal D}^{\,1,\,\,2}\,,\,$\\[0.1in]
(2.2)
$$
\langle\,\psi\,,\ \,[\,(\,\lambda_{\ \!l} \cdot \partial_{\lambda_{\ \!l}}\,)\,V_{\,l}\,]\ \rangle_{\,\btd} \ = \ 0 \ \ \ \Longleftrightarrow \ \   {1\over { n\,(\,n\ + \ 1\,)}} \cdot \bigg\langle  \ \psi\,, \  \left[ \    V_{\,l}^{{4}\over {n\,-\,2}} \cdot (\,\lambda_{\ \!l} \cdot \partial_{\lambda_{\ \!l}}\,)\,  \,V_{\,l} \ \right]  \  \bigg\rangle_{\int} \ = \ 0\,,
$$
where
$$
\bigg\langle  \ \psi\,, \  \left[ \    V_{\,l}^{{4}\over {n\,-\,2}} \cdot (\,\lambda_{\ \!l} \cdot \partial_{\lambda_{\ \!l}}\,)\,  \,V_{\,l} \ \right]  \  \bigg\rangle_{\int} \ = \ \int_{\R^n} \psi \cdot \left[ \    V_{\,l}^{{4}\over {n\,-\,2}} \cdot (\,\lambda_{\ \!l} \cdot \partial_{\lambda_{\ \!l}}\,)\,  \,V_{\,l} \ \right] \ .
$$
This can be shown by  using  equation (\,1.6\,)
and integration by parts formula (\,see for example\,,\, \cite{I}\ )\,.\, Back to (2.1)\,,\,
 $\,a_{\,l}\,$ and $\,b_{\,l\,,\,\,j}\,$ (\,depending on $\,W_{\,\flat}\,$) are coefficients  to be determined (\,see Proposition \,{\bf 2.7}\, concerning uniqueness\,.\,)\,

\vspace*{0.2in}

{\large{\bf \S\,2\,a.}}\ \ \  {\it The weighted sup\,-\,norms\,.\,}\ \
Juxtaposition of multiple bubbles is initially discussed in \cite{Bahri}\,.\, As shown in \cite{Wei-Yan}\,,\, by introducing selected weighted norms specific to the arrangement \{\,$\,\xi_{\,1}\,,\ \cdot \cdot \cdot\,, \ \xi_{\,\flat}\,$\}\, of the  bubbles, the discussion can be brought to its clarity\,.\, We borrow from \cite{Wei-Yan} the following definitions for our use here\,.\,
\begin{eqnarray*}(\,2.3\,) \ \ \ \
\Vert \,\tt{f}\,\Vert_*\!&  = &  \sup_{y\, \in\, \R^n} \left( \ \sum_{l\,=\,1}^\flat\  {1\over { \   (  \,1\ + \ \Vert\,y \,-\,\xi_{\,l}\,\Vert\,)^{ {{n\,-\,2}\over 2 } \,+\,\tau_{\,>\,1} }\  }}\
\right)^{\!\!-\,1} \cdot |\,{\tt{f}}\,(\,y)\,|\,, \ \ \ \ \\[0.2in]
{\mbox{and}} \ \  \Vert \, \tt{g} \,\Vert_{**} \!&  = &  \sup_{y\,\in\,\R^n} \left( \ \sum_{l\,=\,1}^\flat\  {1\over {  \ (\,  1\ + \ \Vert\,y \,-\,\xi_{\,l}\,\Vert\,)^{ {{n\,-\,2}\over 2 } \,+\,2 \,+\,\tau_{\,>\,1} }\  }}\
\right)^{\!\!-\,1}\!\! \cdot |\,{\tt{g}}\,(\,y)\,|\,, \ \  {\tt{f}}\,, \ \  {{\tt{g}}} \,\in \,C^o\,(\,\R^n)\,.
\end{eqnarray*}
Here $\,\tau_{\,>\,1} \,$ is a fixed number slightly {\it bigger}\, than one\,,\, that is\,,\,
$$\,\tau_{\,>\,1} \ = \ 1  \ + \ \varepsilon_\tau \ \ \ \ \  \ (\,\varepsilon_\tau \ > \ 0\,)\ . \leqno (\,2.4\,) $$
Following \cite{Wei-Yan} (\,see also \cite{Pino}\,)\,,\, we record the following. \\[0.1in]
(\,2.5\,)
\begin{eqnarray*}
W_*^\flat\!\!  & = &\!\!\left\{ \, f \,\in\, C^o\,(\,\R^n\,) \ \,\bigg\vert \ \ \Vert \,\tt{f}\,\Vert_* \ < \ \infty \ \right\}\  \ \ \ {\mbox{and}} \ \ \ \ W_{**}^\flat  \ = \ \left\{ \ g \,\in\, C^o\,(\,\R^n\,) \ \,\bigg\vert \ \ \Vert \,\tt{g}\,\Vert_{**} \ < \ \infty \ \right\}\ ,
\end{eqnarray*}

 With reference to (\,2.2\,)\,,\, we also define \\[0.1in]
(\,2.6\,)
 \begin{eqnarray*}
W_{*_\perp}^\flat\!\!  & = &\!\!\bigg\{ \, \phi \,\in\, W_*^\flat \ \,\bigg\vert \ \   \bigg\langle  \ \phi\,, \  \left[ \    V_{\,l}^{{4}\over {n\,-\,2}} \cdot (\,\lambda_{\ \!l} \cdot \partial_{\lambda_{\ \!l}}\,)\,  \,V_{\,l} \ \right]  \  \bigg\rangle_{\int}  \! =  \bigg\langle  \ \phi\,, \  \left[ \    V_{\,l}^{{4}\over {n\,-\,2}} \cdot (\,\lambda_{\ \!l} \cdot  \partial_{\xi_{\,l_{|_{\,j}}}}\,)\,  \,V_{\,l} \ \right]  \  \bigg\rangle_{\int}  = \, 0 \\[0.2in]
& \ & \ \ \ \ \ \ \ \ \ \ \  \ \ \ \ \  \ \ \ \  \ \ \ \ \   \ \ \ \  \ \ \ \ \   \ \ \ \ \ \ \ \ \ \ \ \ \ \ \ \ \ \  {\mbox{for}} \ \ \    j \ = \ 1\,,\ \cdot \cdot \cdot\,, \ n   \ \ \ \ {\mbox{and}} \ \ \ l \ = \ 1\,,\ \cdot \cdot \cdot\,, \ \flat \  \bigg\}\ .
\end{eqnarray*}

 \vspace*{0.15in}

{\bf Proposition 2.7.} \ \ {\it Let $\,n \ \ge \ 6\,.\,$  Under the conditions in\,} (\,1.4\,)\,,  (\,1.8\,)\,, (\,1.22\,)\,, (\,1.24\,)\,--\,(\,1.28\,)\,,\,  \,(\,1.30\,)\, {\it and}  \,(\,1.32\,)\,,\, {\it  we assume that }
$$
\vartheta\ := \ {\mbox{Min}} \ \left\{ \    {{n}\over 2}\,\cdot  (\,1\ - \ \nu\,) \,, \  \ \   \nu\cdot \ell  \  \right\} \ > \  \left(\,  {{ \ n \ + \ 10\  }\over 8}\, \right)\, \cdot \sigma \ .  \leqno (\,2.8\,)
$$
Let
$$
\varpi \ = \ \vartheta \ - \ {1\over 2}\cdot \sigma \ - \ o_{\,+}\,(\,1\,)\,,
$$
 {\it where $\,o_{\,+}\,(\,1\,)$ is a fixed small positive number\,.\, There exists a small positive number}  $\,{{\underline\lambda}}_{\  \epsilon}\,$ {\it such that for all $\,{\bar\lambda}_{\,\,\flat}$ satisfying}
 $$\, {\bar\lambda}_{\,\,\flat} \ \le \  {{\underline\lambda}}_{\  \epsilon}\,,\,\leqno (\,2.9\,)$$
 {\it equation\,} (2.1) {\it has a unique solution\,}  $\,\phi_{\,\,\flat}  \, \in \, W_{*_\perp}^\flat\,\cap\,{\cal D}^{1\,,\,2}\,$ {\it with the following properties.}\\[0.1in]
 (\,2.10\,)

 \vspace*{-0.35in}

 $$
 \Vert \, \phi_{\,\,\flat}  \,\Vert_\btd  \ \le \   {\bar\lambda}_{\,\,\flat}^{\,\varpi}\,. \leqno{\large{\bf *}_1}
 $$
  $$
\Vert \, \phi_{\,\,\flat}  \,\Vert_* \ \le \  {\bar\lambda}_{\,\,\flat}^{\,\varpi\, -\ \tau_{\ >\,1} }  \ .  \leqno{\large{\bf *}_2}
 $$
  $$ |\ a_{\,l} \,| \ + \
|\ b_{\,l_{\,j}}\,| \ \le \  \lambda_{\,\,\flat}^{\,\varpi}    \leqno{\large{\bf *}_3}
 $$
 {\it for}\, $\,j \ = \ 1\,,\ \cdot \cdot \cdot\,\ n\,,$\,   {\it{and}} \ $l \ = \ 1\,,\ \cdot \cdot \cdot\,, \ \flat\,.\,$
{\it In addition\,,\, $\phi_{\,\,\flat} \,$ depends on the bubbles parameters $\,{\bf B}_{\,\flat}$ \,in a $\,C^1$\,-\,manner\, }\,
[\ {\it fulflling\,} (\,1.4\,)\,,  (\,1.8\,)\,, (\,1.22\,)\,, (\,1.24\,)\,--\,(\,1.28\,)\,,\,  \,(\,1.30\,)\, {\it and}  \,(\,1.32\,)\, ]\,.\, {

\vspace*{0.15in}
\hspace*{0.5in}The proof of Propositon 2.7 proceeds in a similar fashion as in the proof of Proposition 2.2 in \cite{Wei-Yan}\,,\, with modification based on the generic setting here, and additional resource from \cite{Pino}\,.\, We present the detials in \,{\bf \S\,A\,3}\,,\,  \,{\bf \S\,A\,4}\,  and \,{\bf \S\,A\,5}\, in the  {\bf e}\,-\,Appendix. Refer to Remark 1.23 for the lack of constants in the right hands sides of $\,{\bf *}_{\,1}\,,\,$ $\,{\bf *}_{\,2}\,,\,$ and $\,{\bf *}_{\,3}\,$ [\,\,in Proposition 2.7 \,]\,.

\newpage

{\bf \large {\bf \S\,3. \ (\,Finite dimension\,) reduced functional and its critical points\,.}}\\[0.2in]
%
{\bf Lemma 3.1.} \ \   {\it Suppose that the bubble parameters}
$$
 \lambda_{\,1}\,, \ \cdot \cdot \cdot\,, \ \lambda_{\,\,\flat}\,, \ \  \xi_{\,1}\,, \ \cdot \cdot \cdot\,, \ \xi_{\,\flat}
$$
{\it satisfy the conditions in Proposition} {\bf 2.7}} (\,{\it with a smaller choice of $\ {\underline{\lambda}}_{\ \epsilon}\,$})\,,\, {\it  assume also that}
$$
   \gamma \ > \ {{\sigma}\over {\,n\ - \ 2\,  }} \ .
   \leqno (\,3.2\,)
$$
{\it Let $\,\phi_{\,\,\flat}\,$ be the  small solution to equation} (\,2.1\,) {\it associated with}\, $\,W_{\,\flat}$\,,\, {\it as described in Proposition}\, {\bf 2.7}\,.\, {\it Here}
 $$
  W_{\,\flat} \ = \ \sum_{l\,=\,1}^\flat\,V_{\,l}\ = \ \sum_{l\,=\,1}^\flat\,V_{\,\lambda_{\,\,l}\,,\ \xi_{\,\,l} }\ . \leqno (\,3.3\,)
 $$
  {\it Define}\\[0.1in]
(\,3.4\,)
\begin{eqnarray*}
& \ & \\[-0.4in]
& \ &  {\bf I}_{\,{\cal R}} \,(\,\lambda_{\,1}\,, \ \cdot \cdot \cdot\,, \ \lambda_{\,\,\flat}\,; \ \, \xi_{\,1}\,, \ \cdot \cdot \cdot\,, \ \xi_{\,\flat}\,)\ = \ {\bf I}\,(\ W_{\,\flat}\ + \ \phi_{\,\,\flat}\ )  \ \ \ \ \ \  \ \ \\[0.3in]
& = &  {1\over 2}\,\int_{\R^n}\,\langle\,\btd\,(\ W_{\,\flat} \ + \ \phi_{\,\,\flat}\,)\,,\ \btd\,(\ W_{\,\flat} \ + \ \phi_{\,\,\flat}\,)\,\rangle\,-\, \left(\,{{n\,-\,2}\over {2n}}\,\right) \cdot \int_{\R^n}\,(\,{\tilde c}_n\!\cdot K\,)\,(\ W_{\,\flat} \ + \ \phi_{\,\,\flat}\,)_+^{\,{{2n}\over {\,n\,-\ 2\,}}}\ \ .
\end{eqnarray*}

{\it If }
$$
 {\bf P}_{\,(\,\flat)\,} \ = \ (\,\lambda_{\,1}\,, \ \cdot \cdot \cdot\,, \ \lambda_{\,\,\flat}\,; \ \xi_{\,1}\,, \ \cdot \cdot \cdot\,, \ \xi_{\,\flat}\,)\,, \leqno (\,3.5\,)
 $$
{\it  is a critical point of the reduced functional $\, {\bf I}_{\,\cal R}$\ ,\,} {\it that is}\,,\\[0.01in]

\vspace*{-0.3in}

$$
 {{\partial \,{\bf I}_{\,\cal R}}\over {\partial\, \lambda_{\ \!l} }}\,\bigg\vert_{\  {\bf P}_{\,(\,\flat\,)} } \ = \ {{\partial \,{\bf I}_{\,\cal R} }\over {\partial \,\xi_{{\,l \,|_{\,j}}} }} \,\bigg\vert_{\   {\bf P} _{\,(\,\flat\,)} } = \ 0 \ \ \ \ \ {\it{for}}  \ \ \ l\  = \ 1\,,\, \cdot \cdot \cdot\,, \   \flat\,, \ \ \ \ j \ = \ 1\,, \ 2\,, \ \cdot \cdot \cdot \,, \ n\,, \leqno  (\,3.6\,)
$$

\vspace*{0.1in}

{\it then\,}  $\,W_{\,\flat\,} \, + \ \phi_{\,\,\flat}\,\,$ [\ {\it corresponding to} $\  {\bf P} _{\,(\, \flat\,)}\,$ {\it via} (\,3.3\,) and (\,3.5\,) \,]\, {\it is also a critical point of}\, (\,{\it the full functional}\ )\, $\,{\bf I} \,$,\, {\it that is\,,\,}
$$
{\bf I}\,' \, (\ W_{\,\flat\,} \ + \ \phi_{\,\,\flat}\ ) \ = \ 0\,.
$$

\vspace*{0.15in}

\hspace*{.5in}The proof of Lemma follows a similar scheme as in  the proof of Theorem 2.8 in   \cite{III} (\ see also  Theorem 2.12 in \cite{Progress-Book}\ )\,.\, More care has to be taken here as the number of bubble $\,\flat\,\to\,\infty\,.\,$ We present the details in {\bf \S\,A\,6}\, in the {\bf e}\,-\,Appendix.

\vspace*{0.2in}

\newpage

{\bf \large {\bf \S\,4. \ \  Extracting the key information in the reduced functional\ - }}\\[0.1in]
\hspace*{0.48in}{\bf \large \ \  {\bf  derivative with respect to  $\,\lambda$\,.}}\\[0.2in]
{\bf \large {\bf \S\,4 a\,.}} \  \ {\it Separation of $\,W_{\flat}\,$ and $\,\phi_{\,\,\flat}\,$ in the reduced functional and its derivatives\,.}\\[0.15in]
{\bf Lemma 4.1.} \ \ {\it Under the conditions in Theorem} \,{\bf 1.33}\,,\, {\it we have}
\begin{eqnarray*}(4.2) \ \ \ \ \ \ \ \ \ \ \ \ \ \ \ \ \ \
{\bf I}_{\cal R} \ (\, W_{\,\flat} \ + \ \phi_{\,\,\flat}\,) \,|_{\,{\bf P}_{\,(\,\flat\,)} }& = & {\bf I}\,(\,W_{\,\flat\,}\,) \,|_{\,{\bf P}_{\,(\,\flat\,)} }\ + \ {\cal E}_{\,o}\\[0.2in]
(4.3) \ \ \ \ \ \ \ \ \ \ \ \ \ \ \ \ \,   \left(\ \lambda_{\ \!l}\cdot {{\partial }\over  {\partial\, \lambda_{\ \!l}}} \, \right)\,  {\bf I}_{\,R}\ \Bigg\vert_{\,{\bf P}_{\,(\,\flat\,)} }& = &  \left(\ \lambda_{\ \!l}\cdot {{\partial }\over  {\partial\, \lambda_{\ \!l} }} \, \right)\,  {\bf I}\,(\,W_{\,\flat\,}\,) \ \Bigg\vert_{\,{\bf P}_{\,(\,\flat\,)} }\ + \   {\cal E}_{\,l}\\[0.2in]
(4.4)  \ \ \ \ \ \ \ \ \ \ \ \ \ \ \,\left(\ \lambda_{\ \!l}\cdot {{\partial }\over {\partial\, \xi_{\,l_{\,|_{\,j}}}}}  \, \right)\,  {\bf I}_R\ \Bigg\vert_{\,{\bf P}_{\,(\,\flat\,)} } & = &  \left(\ \lambda_{\ \!l}\cdot {{\partial }\over  {\partial\, \xi_{\,l_{\,|_{\,j}}}}}  \, \right)\,  {\bf I}\,(\,W_{\,\flat\,}\,) \ \Bigg\vert_{\,{\bf P}_{\,(\,\flat\,)} } \ + \    {\cal E}_{\,1_{|\,j}}  \ .\ \ \ \ \ \ \ \ \ \ \ \ \ \ \ \ \ \ \ \ \ \ \ \ \ \ \ \ \ \ \ \ \ \ \ \
 \end{eqnarray*}
 {\it Here}\\[0.1in]
 (4.5)

 \vspace*{-0.3in}

 $$
 \max_{1\,\le \,l\,\le \,\flat\,;\ 1\,\le\,j\,\le\,n} \ \{ \ {\cal E}_o\,, \  {\cal E}_{\,l}\,, \ {\cal E}_{\,1_{|\,j}}  \ \} \ = \ O\,\left(\, {\bar\lambda}_{\,\,\flat}^{{\varpi}  }\, \right)\,\cdot\,\left[ \   {\bar\lambda}_{\,\,\flat}^{{{4\over {n\,-\,2}} \,\cdot\ \varpi}  }  \ \ + \    {\bar\lambda}_{\,\,\flat}^{{ \, \,4\,\gamma \ -\ o\,(\,1\,)\,}  }  \ \ + \ {\bar\lambda}_{\,\,\flat}^{ \varpi \ - \ \sigma} \  \right]\ .
 $$
{\it In the above\,,\, $o\,(\,1\,)$\, is a positive number that can be chosen to be small\,.\,}

 \vspace*{0.15in}


\hspace*{0.5in}For the detail in proof of Lemma {\bf 4.1}\,,\, we refer to {\bf \S\,A\,7}\, in the {\bf e}\,-\,Appendix\,,\, which essentially follows the argument in \cite{Wei-Yan}\,.\, We pick it up from (\,4.3\,) here [\,with bubble parameters $\,{\bf P}_{\,(\,\flat\,)}\,$,\, integration over the variable $\,y\,$ \,] \,.\, \\[0.1in]
(4.6)
\begin{eqnarray*}
\hspace*{-0.5in}& \  &  \left(\ \lambda_{\ \!l}\cdot {{\partial }\over  {\partial\, \lambda_{\ \!l} }} \, \right)\,  {\bf I}\,(\,W_{\,\flat\,}\,)
\ = \  - \int_{\R^n} \left[ \  \Delta\,W_{\,\flat\,}  \ + \  (\,{\tilde c}_n\!\cdot K\,)\,(\,W_{\,\flat\,}\, )^{{n\,+\,2}\over {\,n\,-\ 2\,}}  \ \right]\,\cdot\,\left[ \, \left(\ \lambda_{\ \!l}\cdot {\partial\over {\partial\, \lambda_{\ \!l}}}\ \right) V_{\,l}   \, \right]  \ \ \ \ \ \  \\[0.2in]
  & = & \ -\int_{\R^n} \left[ \  \Delta\,W_{\,\flat\,} \ + \  n\,(\,n\,-\,2)\cdot\,(\,W_{\,\flat\,}\, )^{{n\,+\,2}\over {\,n\,-\ 2\,}} \ \right]\,\cdot\,\left[ \, \left(\ \lambda_{\ \!l}\cdot {\partial\over {\partial\, \lambda_{\ \!l}}}\ \right) V_{\,l}   \, \right]    \ \  \,\cdot\,\cdot\,\cdot\,\cdot\,\cdot\,\,\cdot\,\cdot\,\cdot\,\,\cdot\,\cdot\,\cdot \ \   {\bf (A)}_{\,(\,4.6\,)}\\[0.2in]
 & \ &  \ \ \ \ \ \ \ \ \ \ \  \ \  - \    \int_{\R^n}\, \bigg[\,(\,{\tilde c}_n\!\cdot K\,) \ - \ n\,(\,n\,-\,2)\,\bigg]  \cdot\,(\,W_{\,\flat\,}\, )^{{n\,+\,2}\over {\,n\,-\ 2\,}} \cdot\left[ \, \left(\ \lambda_{\ \!l}\cdot {\partial\over {\partial\, \lambda_{\ \!l}}}\ \right) V_{\,l}   \, \right] \ .   \cdot\,\cdot \,\cdot\,\cdot \ \ {\bf (B)}_{\,(\,4.6\,)}\\
 \end{eqnarray*}

\newpage

{\bf \large {\bf \S\,4.\,b.}} \  \ {\it Extracting the key information in}\, ${\bf (\,A\,)}_{\,(\,4.6\,)}$\,--\ {\it bubble interactions}\,.\\[0.1in]
 Through equation (\,1.6\,) and expression (1.7)\,,\,  we obtain
  \begin{eqnarray*}
  (4.7) \ \ \ \ \ \   (\,A\,)}_{\,(\,4.6\,) & = &  n\,(\,n\,-\,2)\,\cdot\,\int_{\R^n} \Bigg\{  \,\left[\  V_{\,l}^{{n\,+\,2}\over {\,n\,-\ 2\,}} \ + \  V_2^{{n\,+\,2}\over {\,n\,-\ 2\,}} \ +\,\cdot\,\cdot\,\cdot\, + \ V_\flat^{{n\,+\,2}\over {\,n\,-\ 2\,}} \ \right] \\[0.2in]
    & \ & \ \ \ \ \ \  \ \ \ \ \ \
    - \ \left[\  V_{\,l}\ + \  V_2  \ + \,\cdot\,\cdot\,\cdot\, + \ V_\flat  \ \right]^{{n\,+\,2}\over {\,n\,-\ 2\,}} \  \, \Bigg\} \cdot\left[ \, \left(\ \lambda_{\ \!l}\cdot {\partial\over {\partial\, \lambda_{\ \!l}}}\ \right) V_{\,l}   \, \right] \ .\ \ \ \ \ \  \ \ \ \ \ \ \ \ \ \ \ \  \ \ \ \ \ \
\end{eqnarray*}
Inside $\,B_{\,\xi_{\,1}}(\,\rho_{\,\nu})\, $ [ \ cf. (\,1.29\,) for the definition of $\,\rho_\nu$\ ]\,,
 \begin{eqnarray*}
V_1\,(\,y\,) &  = & \left(\ {{\lambda_{\,1}}\over { \lambda_{\,1}^2 \ + \ \Vert\,y\ - \ \xi_{\,1}\,\Vert^2}}\
\right)^{\!\! {{\,n\,-\,2\,}\over 2}} \ \ge \  \left(\ {{\lambda_{\,1}}\over { \lambda_{\,1}^2 \ + \ \rho_{\,\nu}^2}}\
\right)^{\!\! {{\,n\,-\,2\,}\over 2}} \\[0.2in]
& \ge & C_1\,\cdot\, \left(\ {{1 }\over {  {\bar\lambda}_{\,\,\flat}}}
\right)^{\!\! {{\,n\,-\,2\,}\over 2}} \,\cdot\,  \  {\bar\lambda}_{\,\,\flat}^{(\,n\,-\,2\,)\,\cdot\, (\,1\,-\ \nu\,)} \ \ \ \  \ \ \ \ \ \ \ \ \ \ \ \ \ \ \ \ \ \ \ \ [ \ {\mbox{for}}\ \ \  y\,\in\,B_{\,\xi_{\,1}}(\,\rho_{\,\nu})\ ]\,,\\[0.2in]
{\mbox{and}} \ \ \ \ & \ & V_2\,(\,y\,) \ + \ \cdot\,\cdot\,\cdot\,\ + \ V_\flat\,(\,y\,) \\[0.15in]
&  = &  \sum_{k\,=\,2}^\flat  \left(\ {{\lambda_k}\over { \lambda_k^2 \ + \ \Vert\,y\ - \ \xi_k\,\Vert^2}}
\right)^{\!\! {{\,n\,-\,2\,}\over 2}}
\  = \ \sum_{k\,=\,2}^\flat \left(\ {{1 }\over { {\bar\lambda}_{\,k}}}
\right)^{\!\! {{\,n\,-\,2\,}\over 2}} \cdot \left(\ {{\lambda_k}\over { 1\ + \ \bigg\Vert\,{y\over {\lambda_k}} \ - \ {{\xi_k}\over {\lambda_k}} \,\bigg\Vert^2}}
\right)^{\!\! {{\,n\,-\,2\,}\over 2}}\\[0.2in] & \le & C_2\cdot
\left(\ {{1 }\over { {\bar\lambda}_{\,\,\flat}}}
\right)^{\!\! {{\,n\,-\,2\,}\over 2}}\!\! \cdot \left[ \ \sum_{k\,=\,2}^\flat  \left(\ {{1}\over { 1\ + \ \Vert\,\Xi_1\ - \ \Xi_{\,k}\,\Vert^2}}
\right)^{\!\! {{\,n\,-\,2\,}\over 2}}   \ \right]
\ \le \  C_3\cdot \left(\ {{1 }\over {{\bar\lambda}_{\,\,\flat}}}
\right)^{\!\! {{\,n\,-\,2\,}\over 2}}  \!\! \cdot  {\bar\lambda}_{\,\,\flat}^{(\,n\,-\,2\,)\,\cdot\, \gamma}\\[0.1in]
& \ & \hspace*{2.5in} \ \ \ \  \ \ \ \ \ \ \ \ \ \ \ \ \ \ \ \ \ [ \ {\mbox{for}}\ \ \  y\,\in\,B_{\,\xi_{\,1}}(\,\rho_{\,\nu})\ ]\,.
\end{eqnarray*}
Together we obtain
$$
{{\ V_2\,(\,y\,) \ + \ \cdot\,\cdot\,\cdot\,\ + \ V_\flat\,(\,y\,) \  }\over {V_1\,(\,y\,)  }} \ \le \ C_4 \cdot {\bar\lambda}_{\,\,\flat}^{\,`(\,n\,-\,2\,)\,\cdot\, [ \ (\,\gamma\,+\ \nu\,) \ - \ 1\ ] } \ \  \ \ \mfor \ \ \  y\,\in\,B_{\,\xi_{\,1}}(\,\rho_{\,\nu}\,)\,.
$$
Recall from (\,1.31\,) that $\,\gamma \ + \ \nu \ > \ 1\,.$
On account of the above\,,\, in $\,B_{\,\xi_{\,1}}(\,\rho_{\,\nu})\,,\,$  we make use of the expansion
\begin{eqnarray*}
(4.8) & \ &   \left[\  V_1\ + \  V_2  \ + \ \cdot\,\cdot\,\cdot\,\ + \ V_\flat  \ \right]^{{n\,+\,2}\over {\,n\,-\ 2\,}}  \ = \ V_1^{{n\,+\,2}\over {\,n\,-\ 2\,}}\,\cdot\,\left[\  1 \ + \ \bigg( \ {{  V_2  \ + \ \cdot\,\cdot\,\cdot\,\ + \ V_\flat}\over {V_1}} \ \bigg)  \, \right]^{{n\,+\,2}\over {\,n\,-\ 2\,}}\\[0.15in]
& = & V_1^{{n\,+\,2}\over {\,n\,-\ 2\,}}\,\cdot\,\left[\  1 \ + \ {{n\,+\,2}\over {\,n\,-\ 2\,}}\,\cdot\,\left(\ \ {{  V_2  \ + \ \cdot\,\cdot\,\cdot\,\ + \ V_\flat}\over {V_1}}  \ \right)\ + \ O \left(\ \  {{  V_2  \ + \ \cdot\,\,\cdot\,\cdot \ + \ V_\flat}\over {V_1}}\  \right)^{{n\,+\,2}\over {\,n\,-\ 2\,}}\  \right] \ \ \ \ \  \ \\[0.15in]
& = & V_1^{{n\,+\,2}\over {\,n\,-\ 2\,}} \ + \  {{n\,+\,2}\over {\,n\,-\ 2\,}} \,\cdot\,V_1^{4\over {n\,-\,2}}\,\cdot\, \left[\   V_2  \ + \ \cdot\,\cdot\,\cdot\,\ + \ V_\flat  \ \right] \ + \ O \left(\,  \left[ \, V_2  \ + \ \cdot\,\cdot\,\cdot\,\ + \ V_\flat \, \right]^{  {{n\,+\,2}\over {\,n\,-\ 2\,}}} \right) \ .
\end{eqnarray*}

\newpage

It follows from (4.8)  that\\[0.1in]
(4.9)
\begin{eqnarray*}
& \ & \\[-0.4in]
& \ &  n\,(\,n\,-\,2)\,\cdot\,\int_{B_{\,\xi_{\,1}} (\,{\rho_{\,\nu}}\,)}     \left\{  \,\left[\  V_1^{{n\,+\,2}\over {\,n\,-\ 2\,}} \ + \  V_2^{{n\,+\,2}\over {\,n\,-\ 2\,}} \ +\,\cdot\,\cdot\,\cdot\, + \ V_\flat^{{n\,+\,2}\over {\,n\,-\ 2\,}} \ \right] \right. \\[0.2in]
    & \ & \left. \ \ \ \ \ \  \ \ \ \ \ \  \ \ \ \ \ \  \ \ \ \ \ \  \ \ \ \ \ \  \ \ \ \ \ \  \ \ \ \ \ \
    - \ \left[\  V_1\ + \  V_2  \ + \,\cdot\,\cdot\,\cdot\, + \ V_\flat  \ \right]^{{n\,+\,2}\over {\,n\,-\ 2\,}} \  \right\}\,\cdot\,\left[ \, \left(\ \lambda_{\,1}\cdot {\partial\over {\partial\, \lambda_{\,1}}}\ \right) V_1   \, \right] \\[0.2in]
& = & -\,n\,(\,n\,-\,2)\,\cdot\,\int_{B_{\,\xi_{\,1}} (\,{\rho_{\,\nu}}\,)} \left\{ \ \left(\, {{n\,+\,2}\over {\,n\,-\ 2\,}} \ \right)\,\cdot\,V_1^{4\over {n\,-\,2}}  \cdot  \left[\   V_2  \ + \ \cdot\,\cdot\,\cdot\,\ + \ V_\flat  \ \right] \   \right\}\cdot \left(\ \lambda_{\,1}\cdot {{\partial\,V_1}\over {\partial\, \lambda_{\,1}}}\ \right) \\ [0.2in]
& \ & \ \ \ \ \ \ \ \ \ + \ {\bf{error}}_{\,(\,4.9\,)}\,,
\end{eqnarray*}
where
\begin{eqnarray*}
 |\,{\bf{error}}_{\,(\,4.9\,)} \,| & \le &   O\,\left( \  {\bar\lambda}_{\,\,\flat}^{\,n\,\cdot\, \gamma}\ \right)\,\cdot\, O\,\left( \  {\bar\lambda}_{\,\,\flat}^{\ 2\ [\ (\,\gamma \ + \ \nu \,) \ - \ 1\ ]}\ \right)\ \ \ \ \ \ \ \ \   \ \ \ \ \ (\ \gamma \ + \ \nu \ >  \ 1\ )\,,
\end{eqnarray*}
which is estimated in  the {\bf \S\,A\,8\,.b} in the \,{\bf e}\,-\,Appendix\,.\, \bk Let us focus on the ``\,interaction\," term between the bubbles $\,V_1\,$ and (\,for example\,)   $\,V_2\,$:
$$
   V_1^{4\over {n\,-\,2}}  \cdot V_2 \cdot \left(\ \lambda_{\,1}\cdot {{\partial\,V_1}\over {\partial\, \lambda_{\,1}}}\ \right) \ .
$$

\vspace*{-0.13in}

We begin with [\,inside $\,B_{\,\xi_{\,1}}(\,\rho_{\,\nu})\ $]\,,\, \\[0.1in]
(4.10)
 \begin{eqnarray*}
 V_2 (\,y\,) & = & \left(\ {{\lambda_2}\over { \lambda^2_2 \ + \ \Vert\,y\ - \ \xi_2\,\Vert^2  }}\ \right)^{\!\!{{\,n\,-\,2\,}\over 2}  }\\[0.125in]
& = &  \left(\ {{{1\over\lambda_{\,1}}}\over { \left(\ {{\lambda_2}\over {\lambda_{\,1}}}  \right)\ + \ {{\Vert\,y\ - \ \xi_2\,\Vert^2}\over
{\lambda_{\,1}\,\cdot\,\lambda_2}}  }}\ \right)^{\!\!{{\,n\,-\,2\,}\over 2}  } \ = \  {1\over {\lambda_{\,1}^{{\,n\,-\,2\,}\over 2} }}\,\cdot\,\left(\  {{{1}}
\over { \left(\ {{\lambda_2}\over {\lambda_{\,1}}}  \right)\ + \ {{\Vert\,(\,y\ - \ \xi_{\,1})\ + \ (\,\xi_{\,1}\ - \ \xi_2)\,\Vert^2}\over
{\lambda_{\,1}\,\cdot\,\,\lambda_2}}  }}\ \right)^{\!\!{{\,n\,-\,2\,}\over 2}  } \\[0.125in]
& = &  {1\over {\lambda_{\,1}^{{\,n\,-\,2\,}\over 2} }}\,\cdot\,\left(\ {{{1}}
\over { \left(\ {{\lambda_2}\over {\lambda_{\,1}}}  \right) \ + \ {{\Vert\, \,\xi_{\,1}\ - \ \xi_2\,\Vert^2}\over
{\lambda_{\,1}\,\cdot\,\,\lambda_2}}   \ + \ {{\Vert\,\,y\ - \ \xi_{\,1}\,\Vert^2}\over
{\lambda_{\,1}\,\cdot\,\,\lambda_2}}   \ + \ {{\,2\,(\,y\ - \ \xi_{\,1})\cdot(\,\xi_{\,1}\ - \ \xi_2)\, }\over
{\lambda_{\,1}\,\cdot\,\,\lambda_2}}}}\ \right)^{\!\!{{\,n\,-\,2\,}\over 2}  }\\[0.125in]
&  \ & \hspace*{-0.85in} = \  {1\over {\lambda_{\,1}^{{\,n\,-\,2\,}\over 2} }}\cdot{1\over  { {\bf d}_{1\,,\,2}^{\,n\,-\,2}  }}\cdot\! \left(\,{{{1}}
\over { \, 1  + \,\left\{\ {1\over { {\bf d}_{1\,,\,2}^{\,2}  }}\,\cdot\,\left(\ {{\lambda_2}\over {\lambda_{\,1}}}  \right) \ + \ {1\over  { {\bf d}_{1\,,\,2}^{\,2}  }}\,\cdot\,{{\Vert\,\,y\ - \ \xi_{\,1}\,\Vert^2}\over
{\lambda_{\,1}\,\cdot\,\,\lambda_2}}   \ + \ {1\over { {\bf d}_{1\,,\,2}^{\,2}  }}\cdot {{\,2\,(\,y\ - \ \xi_{\,1})\,\cdot\,(\,\xi_{\,1}\ - \ \xi_2)\, }\over
{\lambda_{\,1}\,\cdot\,\,\lambda_2}}\, \right\} \ }}\,\right)^{\!\!\!{{n\,-\,2\,}\over 2}  }\\[0.125in]
& \ & \ \ \ \ \ \ \ \ \ \ \ \ \ \ \ \ \ \ \  \leftarrow \ \ \ \ \ \ \ \ \ \ \ \  \ \ \ \ \ \ \ \ \  \ \ \ \ \ \ {\bf T}\   \ \ (\,\leftarrow \ \ {\mbox{small}}\,) \ \ \ \ \  \ \ \ \ \ \ \ \ \ \ \ \ \ \ \    \ \ \ \    \rightarrow \\[0.125in]
&  &\hspace*{-0.85in} = \ {1\over {\lambda_{\,1}^{{\,n\,-\,2\,}\over 2} }}\,\cdot\,{1\over { {\bf d}_{1\,,\,2}^{n\,-\,2} }} \times \   \left\{ \  1 \ - \  {1\over 2} \cdot (\,n\ - \ 2\,)\,\cdot\,{\bf T}  \ + \  O\,(\,|\, {\bf T} \,|^{\,\,2}\,) \ \right\} \ \ \  \mfor \ \ y\,\in\,{B_{\,\xi_{\,1}} (\,{\rho_{\,\nu}}\,)} \ .\\[0.1in]
& \ & \hspace*{0.7in}[\  \uparrow \ \ {\mbox{Leading\ \  order\ \  team}}\,. \  ]
\end{eqnarray*}
Recall that from (\,1.9\,) that
$$
{\bf d}_{\,1\,,\ 2}  \  :=   \  {{\ \Vert\,\xi_{\,1}\ -\ \xi_2\,\Vert\ }\over{\,\sqrt{\,\lambda_{\,1} \,\cdot\,\lambda_2\,}\,}} \ \ \ \  (\ \to \ \infty \  \ {\mbox{as}} \ \ {\bar\lambda}_{\,\,\flat} \ \to \ 0^+\,) \ .
$$
Also\,,\, the term $\,{\bf T}\,$ is defined via
\begin{eqnarray*}
(4.11) \ \ \ {\bf T} & := &  {1\over { {\bf d}_{1\,,\,2}^{\,2}  }}\,\cdot\,\left(\ {{\lambda_2}\over {\lambda_{\,1}}}  \right) \ + \ {1\over  { \ {\bf d}_{1\,,\,2}^{\,2}\   }}\,\cdot\,{{\Vert\,\,y\ - \ \xi_{\,1}\,\Vert^2}\over
{\lambda_{\,1}\,\cdot\,\,\lambda_2}}   \ + \ {1\over { \ {\bf d}_{1\,,\,2}^{\,2}\  }}\cdot {{\,2\,(\,y\ - \ \xi_{\,1})\,\cdot\,(\,\xi_{\,1}\ - \ \xi_2)\, }\over
{\lambda_{\,1}\,\cdot\,\,\lambda_2}} \ \ \ \ \ \ \ \ \ \ \ \ \ \\[0.2in]
& = &  O\,\left(\  {\bar\lambda}_{\,\,\flat}^{\ (\,\gamma\,+\,\nu\,) \, - \ 1 } \ \right)\ .
\end{eqnarray*}
 See {\bf \S\,A\,8\,.\,c\,} in the \,{\bf e}\,-\,Appendix\,.\,  As a record,
 \begin{eqnarray*}
  \lambda_{\,1} \cdot {{\partial V_{\,1}}\over {\partial \lambda_{\,1}}} =   \lambda_{\,1} \cdot {{\partial}\over {\partial \lambda_{\,1}}} \left[\,\left({\lambda_{\,1}\over {\lambda_{\,1}^2 + |\, y - \xi_{\,1}\,|^{\,2}}} \right)^{{n - 2}\over 2}\  \right] \ = \ -\left( \,{{n - 2}\over 2}\,\right)  \cdot \lambda_{\,1}^{{n - 2}\over 2} \cdot {{ (\ \lambda^2_1 - |\, y \,- \,\xi_{\,1}\,|^{\,2}\ )}\over {\ (\ \lambda^2_1 + |\, y \,-\, \xi_{\,1}\,|^{\,2}\,)^{{n }\over 2}\ }}\ .
\end{eqnarray*}
We collect the crucial term in the following\,.
 $$
\bigg[ \ \cdot\,\cdot \cdot\ \bigg]_{\,(\,4.13\,)}\ \Bigg\vert_{\,y}  \, \ = \  V_1^{4\over {n\,-\,2}}
\cdot \left(\ {1\over {  \lambda_{\,1}^{{n\,-\,2}\over 2 }   }}\,\cdot\,{1\over {\large{\bf d}_{1\,,\,2}^{\,n\,-\,2}    }}\   \right)\,\cdot\,\left\{  \  \lambda_{\,1}^{{n \,-\, 2}\over 2}
\cdot {{ (\,\lambda_{\,1}^2\, -\, \Vert\, y \,-\, \xi_{\,1}\,\Vert^{\,2})}\over {\ (\,\lambda_{\,1}^2 \,+\, \Vert\, y \,-\, \xi_{\,1}\,\Vert^{\,2})^{{n }\over 2}\  }}\ \right\} \ . \leqno (4.12)
 $$
 Returning  to (4.9)\,,\, for the part which involves $V_1\,$ and $\,V_2\,,\,$ we continue with \\[0.1in]
(4.13)
 \begin{eqnarray*}
 & \ & \\[-0.34in]
 & \ & -\,n\,(\,n\,-\,2)\,\cdot\,\int_{B_{\,\xi_{\,1}} (\,{\rho_{\,\nu}}\,)} \left\{ \, \left(\ \,{{n\,+\,2}\over {\,n\,-\ 2\,}} \,\right)\,\cdot\,V_1^{4\over {n\,-\,2}}\,\cdot\,\left[\   V_2   \ \right]   \right\}\,\cdot\, \left(\,\lambda_{\,1}\cdot { { \partial \,V_1} \over {\partial\, \lambda_{\,1}}}\ \right)\\[0.2in]
& =  &   +\,n\,(\,n\,-\,2) \cdot\left( \,{{n\,-\,2}\over { 2}} \,\right)\cdot  \left(\ \,{{n\,+\,2}\over {\,n\,-\ 2\,}} \,\right)\,\cdot\,\int_{B_{\,\xi_{\,1}}\,(\,{\rho_{\,\nu}})}  {\bf[} \ \cdot\,\cdot \cdot\ {\bf ]}_{\,(\,4.13\,)} \times\\[0.1in]
& \ & \hspace*{2.5in}\ \ \ \  \times \  \left\{ \  1 \ - \  \left(\,{{\,n\,-\,2\,}\over 2} \right)\,\cdot\,{\bf T}  \ + \  O\,(\,|\, {\bf T} \,|^{\,\,2}\,) \ \right\}\\[0.2in]
& =  &   +\,n\,(\,n\,-\,2) \cdot\left( \,{{n\,-\,2}\over { 2}} \,\right)\cdot  \left(\ \,{{n\,+\,2}\over {\,n\,-\ 2\,}} \,\right)\,\cdot\,\int_{B_{\,\xi_{\,1}}\,(\,{\rho_{\,\nu}})}  {\bf[} \ \cdot\,\cdot \cdot\ {\bf ]}_{\,(\,4.13\,)} \ \ + \ {\bf E}_{\,(\,4.13\,)} \\[0.2in]
&   = &  {\tilde C}_1\,(\,n)\,\cdot\,\left\{ \ {{1}\over {\large{\bf d}_{\,1\,,\ 2}^{\,n\,-\ 2}  }}   \  \right\}\,\cdot\,\left\{ \
{\displaystyle{\int}_0^{\,\rho_{\,\nu}}
\left[  \   \left(\ {{\lambda_{\,1}}\over {\lambda^2_1 \,+\, {\bar r}^2}} \right)^{\!\!2}\!\!
\cdot {{ (\,\lambda_{\,1}^2 \ -\  {\bar r}^{\,2})}\over {\ (\,\lambda_{\,1}^2 \ + \  {\bar r}^{\,2})^{{n }\over 2}\  }}\right]\,\cdot\, {\bar r}^{n\,-\,1}\,\cdot\,d\,{\bar r}  } \ \right\}  \ \ + \ {\bf E}_{\,(\,4.13\,)}  \ .\\
    \end{eqnarray*}
We estimate  \,${\bf E}_{\,(\,4.13\,)}$\, in the {\bf \S\,A\,8\,.\,d} and  {\bf \S\,A\,8\,.\,e}\, of the \,{\bf e}\,-Appendix, and
$$
{\bar y} \ =\ y \ - \ \xi_{\,1} \,,\ \ \ \ \ \ \  {\bar r} \ = \ \Vert\,\bar y\,\Vert \ ,
$$

\newpage

together with
$$
{\tilde C}_1\,(\,n) \ =  \ \omega_n\cdot  n\,(\,n\,-\,2)\,\cdot\,{{n \,- \,2}\over 2}\,\cdot\,  {{n\,+\,2}\over {\,n\,-\ 2\,}}  \ = \ {\omega_n \over 2}\,\cdot\, n\,(\,n\,-\,2\,)\,\cdot\,(\,n \,+ \,2\,)\ . \leqno (\,4.14\,)
$$

    Via the change of variables $\,\bar r \ = \ \lambda_{\,1}\,\cdot\,\tan\,\theta \,$,\, we continue with
    \begin{eqnarray*}
 (4.15)\ \    & \ & \\[-0.25in]
   & \ & \int_0^{\rho_{\,\nu}}\!
\left[  \   \left(\ {{\lambda_{\,1}}\over {\lambda^2_1 \,+\, {\bar r}^2}} \right)^{\!\!2}\!\!
\cdot {{ (\,\lambda_{\,1}^2 \ -\  {\bar r}^{\,2})}\over {(\,\lambda_{\,1}^2 \ + \  {\bar r}^{\,2})^{{n }\over 2} }}\ \right]\,\cdot\, {\bar r}^{n\,-\,1}\,\cdot\,d\,{\bar r} \ \ \  \left( \ {{ \rho_\nu}\over { {\bar{\lambda}}_{\,\flat} }} \ = \ {1\over { \  {\bar{\lambda}^{1\,-\,\,\nu}}_{\,\flat} \  }} \ \to \ \infty \, \right) \ \ \ \  \ \ \ \ \ \ \ \ \ \ \\[0.15in]
& = & \int_0^\infty
\left[  \   \left(\ {{\lambda_{\,1}}\over {\lambda^2_1 \,+\, {\bar r}^2}} \right)^{\!\!2}\!\!
\cdot {{ (\,\lambda_{\,1}^2 \ -\  {\bar r}^{\,2})}\over {(\,\lambda_{\,1}^2 \ + \  {\bar r}^{\,2})^{{n }\over 2} }}\ \right]\,\cdot\, {\bar r}^{n\,-\,1}\,\cdot\,d\,{\bar r}  \ + \ {\bf E}_{\,(\,4.15\,)}\\[0.2in]
& = &
\int_0^{ \,{\pi\over 2}}\   \left\{\  [\ \sin\,\theta\,]^{\,n\,-\,1} \,\cdot\,[\ \cos\,\theta\,]^{\,3} \ - \
[\ \sin\,\theta\,]^{\,n\,+\,1} \,\cdot\,[\ \cos\,\theta\,]\  \right\} \,d\theta  \ + \ {\bf E}_{\,(\,4.15\,)}\\[0.2in]
& = &
\int_0^{ \,{\pi\over 2}}\   \left\{\  [\ \sin\,\theta\,]^{\,n\,-\,1} \,\cdot\,[\ \cos\,\theta\,]  \ - \
2\,[\ \sin\,\theta\,]^{\,n\,+\,1} \,\cdot\,[\ \cos\,\theta\,]\  \right\} \,d\theta \ + \ {\bf E}_{\,(\,4.15\,)}\\[0.2in]
& = &
    \int_0^{ \,{\pi\over 2}} \left\{\, [\ \sin\,\theta\,]^{\,n\,-\,1}     \ - \
2\,[\ \sin\,\theta\,]^{\,n\,+\,1}  \,\right\} \,d\,[\,\sin \,\theta\,] \ + \ {\bf E}_{\,(\,4.15\,)}\\[0.2in]\emph{}
& = &
\left[ \ {1\over n} \ - \ {2\over {n\,+\,2}}\,\right]  \ + \ {\bf E}_{\,(\,4.15\,)}\ = \  -\, {{n\,-\,2}\over {n\,(\,n\,+\,2)  }} \ + \ {\bf E}_{\,(\,4.15\,)}\ .
\end{eqnarray*}
Here
$${\bf E}_{\,(\,4.15\,)} \ = \  O\,\left(\ {\bar{\lambda}}_{\,\flat}^{\,2\, (\ 1\,-\ \nu\,)}\  \right) \ , $$ which is estimated in {\bf \S\,A\,8\,.\,j}\, of the \,{\bf e}\,-\,Appendix\,.\,  We arrive at
 \begin{eqnarray*}
& \ &  -\,n\,(\,n\,-\,2)\,\cdot\,\int_{B_{\,\xi_{\,1}} (\,{\rho_{\,\nu}}\,)} \left\{ \, \left(\ \,{{n\,+\,2}\over {\,n\,-\ 2\,}} \,\right)\,\cdot\,V_1^{4\over {n\,-\,2}}\,\cdot\,\left[\   V_2   \ \right]  \  \right\}\,\cdot\, \left(\,\lambda_{\,1}\cdot { { \partial \,V_1} \over {\partial\, \lambda_{\,1}}}\ \right) \\[0.2in]
&  = &  -\,{\hat C}_{\,1\,,\,2} \,\cdot\,{{\lambda_{\,1}^{{n\,-\,2 }\over 2}\!\cdot\lambda_2^{{n\,-\ 2}\over 2}}\over {\ \Vert\,\xi_{\,1}\,-\,\xi_{\,2}\,\Vert^{\,n\,-\,2}\  }}      \,\cdot\,   \, \left[\ 1  \,+\,  O\,\left(\ {\bar{\lambda}}_{\,\flat}^{\,2\, (\ 1\,-\ \nu\,)}\  \right) \ \right]\ .
 \end{eqnarray*}
In the above [\,combining the numbers in (\,4.14\,) and  (\,4.15\,)\ ]\,,\,
$$
{\hat C}_{\,1\,,\,2} \ =    {{n\,-\,2}\over {n\,(\,n\,+\,2)  }} \,\cdot\,\left[ \ {\omega_n \over 2}\,\cdot\, n\,(\,n\,-\,2)\,\cdot\,(\,n \,+ \,2) \ \right] \ = \ {{(\,n\,-\,2\,)^2}\over 2}\,\cdot\,\, \omega_n  \ . \leqno (4.16)
$$
Similarly we estimate the interactions of $\,V_{\,1}\,$ with other bubbles $\,V_3\,, \ \cdot \cdot \cdot\,, \ V_{\,\flat}\,,\,$  arriving at\\[0.1in]
 (\,4.17\,)
 \begin{eqnarray*}
  {\bf (\,A\,)}_{(\ 4.6\,)}  & = &  -\,{\hat C}_{\,1\,,\,2}  \,\cdot\,\left\{ \ \sum_{k\,\not=\,1}  {{\lambda_{\,1}^{{n\,-\,2 }\over 2}\!\cdot\lambda_k^{{n\,-\ 2}\over 2}}\over {\ \Vert\,\xi_{\,1}\,-\,\xi_k\,\Vert^{\,n\,-\,2}\  }}    \ \right\}  \,\cdot\,   \, \left[\ 1  \,+\,  O\,\left(\ {\bar{\lambda}}_{\,\flat}^{\,2\, (\ 1\,-\ \nu\,)}\  \right) \ \right] \ +   \ \ {\cal E}_{\,(4.17\,)}   \, \ \ \ \
 \end{eqnarray*}
where
$$
{\cal E}_{\,(4.17\,)}  \ = \  O\left(\   {\bar\lambda}_{\ \flat}^{ \   {{n\,+\,2}\over 2}\,\,\cdot\,\,\gamma \ + \ {{\,n\,-\,2\,}\over 2}\,\cdot\,(\,1\ - \ \nu\,) \ - \ o_{\,+\,}\,(\,1\,) } \ \right) \ + \ O\left( \  {\bar\lambda}_{\ \flat}^{ \,n \,\gamma \ - \ \sigma}\ \right)\ .
$$
See \,{\bf \S\,A\,8\,.\,k}\, of the {\bf e}\,-\,Appendix. Here the bubble parameters $\,\left( \ (\,\lambda_{\,1}\,, \ \xi_{\,1}\,)\,, \ \cdot \cdot \cdot\,, \ (\,\lambda_{\,\flat}\,, \ \xi_{\ \flat} \,) \ \right)\,$ satisfy the conditions in Theorem 1.33\,.
Note that
$$\,\gamma \ + \ \nu \ > \ 1 \ \ \Longrightarrow \ \ \gamma \ >  \ 1 \ - \ \nu \ \ \Longrightarrow \ \ {\bar{\lambda}}_{\,\flat}^{\,2\,\gamma}  \ \le \  {\bar{\lambda}}_{\,\flat}^{\ 2\, (\ 1\,-\ \nu\,)}\  \ \ \ \ \ {\mbox{for}} \ \ {\bar{\lambda}}_{\,\flat} \ \le \ 1\ .$$

\vspace*{0.3in}

{\bf \large {\bf \S\,4\,c\,.}} \  \ {\it Extracting the key information in}\, ${\bf (\,B\,)}_{\,(\,4.6\,)}\,-\,$ {\it curvature contribution.}\\[0.1in]
As in the previous sub\,-\,section (\,{\bf \S\,4\,b}\,)\,,\, we begin by focusing on what is inside $\,B_{\,\xi_{\,1}}\,(\,\rho_{\,\nu})\,.\, $
\begin{eqnarray*}
(4.18) \ \ \ \  {\bf{(\,B)} }_{\,(\,4.6\,)} & = &    - \    \int_{\R^n}\, [\,(\,{\tilde c}_n\!\cdot K\,) \ - \ n\,(\,n\,-\,2)\,]  \cdot\,(\,W_{\,\flat\,}\, )^{{n\,+\,2}\over {\,n\,-\ 2\,}} \cdot\left[ \, \left(\ \lambda_{\ \!l}\cdot {\partial\over {\partial\, \lambda_{\ \!l}}}\ \right) V_{\,l}   \, \right] \\[0.2in]
 &  \ & \hspace*{-1.1in} =\  - \    \int_{B_{\,\xi_{\,1}}(\,\rho_{\,\nu})}\, [\,(\,{\tilde c}_n\!\cdot K\,) \ - \ n\,(\,n\,-\,2)\,]  \cdot\,(\,W_{\,\flat\,}\, )^{{n\,+\,2}\over {\,n\,-\ 2\,}} \cdot\left[ \, \left(\ \lambda_{\ \!l}\cdot {\partial\over {\partial\, \lambda_{\ \!l}}}\ \right) V_{\,l}   \, \right] \ + \ {\bf E}_{\,(\,4.18\,)}\ , \ \ \ \ \ \ \ \ \ \ \ \ \
 \end{eqnarray*}
where
$$ \ \ \ \ \ \ \
{\bf E}_{\,(\,4.18\,)} \ = \  O\left( \ \,{\bar\lambda}_{\,\,\flat}^{\,n\,(\,1\,-\ \nu\,) \ - \ o_{\,+}\,(\,1\,)} \ \right) \ \ \ \ \ \ \ \ \ \ [\ \gamma \ >  \ (\ 1\ - \ \nu\ ) \ ] \leqno (4.19)
$$
is estimated in {\bf \S\,A\,9} of the \,{\bf e}\,-\,Appendix [\ specifically\,,\, see (\,A.9.3\,)\,]\,.\, Let us find in (\,4.18\,)\,.
$$
 V_1^{{n\,+\,2}\over {\,n\,-\ 2\,}}\,\cdot\,\left(\ \lambda_{\,1}\,\cdot\,{ { \partial  V_1 } \over {\partial\, \lambda_{\,1}}}\ \right).
$$
Similar to (4.8)\,,\, we obtain
\begin{eqnarray*}
(4.20) \ \ &  \ &   -\,\int_{B_{\,\xi_{\,1}} (\,{\rho_{\,\nu}}\,)}\,[\,(\,{ c}_n\!\cdot K\,) \ - \ n\,(\,n\,-\,2)\,] \cdot\,(\,W_{\,\flat\,}\, )^{{n\,+\,2}\over {\,n\,-\ 2\,}}\,\cdot\,\left(\ \lambda_{\,1} \cdot { { \partial \,V_1 } \over {\partial\, \lambda_{\,1}}}\ \right)  \\[0.2in]   &  = &  \left\{ \ -\, \int_{B_{\,\xi_{\,1}} (\,{\rho_{\,\nu}}\,)}\, [\,(\,{ c}_n\!\cdot K\,) \ - \ n\,(\,n\,-\,2)\,] \cdot\,V_1^{{n\,+\,2}\over {\,n\,-\ 2\,}}\,\cdot\,\left(\ \lambda_{\,1} \cdot { { \partial \,V_1 } \over {\partial\, \lambda_{\,1}}}\ \right)  \  \right\}\,\times \,[ \ 1 \ + \ {\bf E}_{\,(\,4.20\,)}\ ]\ \ \ \ \ \ \ \ \ \\[0.2in]
& = &  -\,{{n\ -\ 2}\over 2}\,\cdot\,\left\{ \  \int_{B_{\,\xi_{\,1}} (\,{\rho_{\,\nu}}\,)} [\ \cdot\,\cdot\,\cdot\,\ ]_{\,(\,4.20\,)} \ d\,y \ \right\} \,\times \,[ \ 1 \ + \ {\bf E}_{\,(\,4.20\,)}\ ]\ .
\end{eqnarray*}
Here
$$
  [\ \cdot\,\cdot\,\cdot\,\ ]_{\,(\,4.20\,)}\ \Bigg\vert_{\ y}   \ = \   \,[\,n\,(\,n\,-\,2) \,-\, (\,{ c}_n\!\cdot K\,)\,(\,y\,)  \ ] \,\cdot\, \left(\ {{\lambda^2_1 \ - \ \Vert\ y\,-\,\xi_{\,1}\,\Vert^2}\over {\lambda^2_1 \ + \ \Vert\ y\,-\,\xi_{\,1}\,\Vert^2}} \right)\,\cdot\,\left(\ {\lambda_{\,1}\over {\lambda^2_1 \ + \ \Vert\ y\,-\,\xi_{\,1}\,\Vert^2}} \right)^{\!\!n }  \ , \leqno (4.21) $$
  and
$$
{\bf E}_{\,(\,4.20\,)}  \ = \   O\,\left( \ {\bar\lambda}_{\,\,\flat}^{\ (\,n\ -\,2\,)\,\cdot\, [ \ (\,\gamma \ + \ \nu \,)\ - \  1\ ] }\ \,\right)\,, \leqno (4.22)
$$
which is estimated in {\bf \S\,A\,9} of the \,{\bf e}\,-\,Appendix [\,in particular\,,\, (\,A.9.2\,)\ ]\,.\, \bk
Recall that
there exists a fixed compact set $\,{\cal C}_{\,[]} \, \subset \, \R^n\,$   such that
$$
 \{\  \xi_{\,1}\,, \ \cdot \cdot \cdot\,, \ \xi_{\ \flat}  \ \} \ \subset \ {\cal C}_{\,[]} \ \ \ \ \ [ \ {\mbox{independent \   \ on}} \ \  \flat \ ] \ .
$$

{\it Geometric set\,-\,up\,.} \ \  Fix $\,\xi_{\,1}\,\not\in\,{\cal H}\,,\,$ which satisfies condition (\,1.28\,)\,,\,\emph{} let $\,{\bf p}_{\,\xi_{\,1}}\,\in\,{\cal H}\,$ be the projection of $\,\xi_{\,1}\,$ to $\,{\cal H}\,,$ according to (\,1.25\,)\,.\, That is\,,\,
$$
(\ 0 \ <  \ ) \ \  \Vert\ \xi_{\,1} \ - \ {\bf p}_{\,\xi_{\,1}} \,\Vert \ = \ {\mbox{dist}}\,(\,\xi_{\,1}\,,\,{\cal H}\,) \ = \ O\, \left(\ {\bar\lambda}_{\,\,\flat}^{\ 1\,+\,\kappa}\,\right)\,.
$$
As the last integral in (\,4.20\,) is invariant under   translations and rotations, we may arrange things so that\\[0.1in]
(\,4.23)
$$
{\bf p}_{\,\xi_{\,1}} \ = \ 0 \ \ \ \ {\mbox{and \   the \   hyperplane   \ defined \   by }} \ \ \{ \,y_{\,n} \ = \ 0\,\} \ \ {\mbox{is \   the \   tangent \    space \   of}} \ \ \ {\cal H} \ \ {\mbox{at \ \ }} 0\,.
$$
Based on this arrangement\,,\, we have
$$
\xi_{\,1}\ = \ (\,0\, \ \cdot\,\cdot \cdot\,, \ 0\,, \ \xi_{1_{\,|_{\,n}}}\,)\ ,  \ \ \ \ \  {\mbox{and \ \ we \ \ choose \ \ the \ \ orientation \ \ so \ \ that}} \ \ \xi_{1_{\,|_{\,n}}} \ > \ 0\,. \leqno (4.24)
$$
Note that in this position
$$
\xi_{1_{\,|_{\,n}}}  \ = \ {\mbox{dist}}\ (\ \xi_{\,1}\,,\,{\cal H}\,) \,.
$$
In $\, B_{\,\xi_{\,1}} \,(\,2\,\rho_{\,\nu\,}\,)\,,\,$ for $\,\rho_{\,\nu\,}\,\approx\,0\,$ (\,that is\,,\, when $\,{\bar\lambda}_{\,\,\flat}\,$ is small\,)\,,\,  $\,{\cal H}\,$ is realized as the level set of a smooth function $\,F\,$ [\,defined (\,locally\,)  inside $\, B_{\,\xi_{\,1}} \,(\,2\rho_{\,\nu}\,)\,$\,]\,:
$$
{\cal H}\,\cap\,B_{\,\xi_{\,1}} \,(\,2\,\rho_{\,\nu\,}\,) \ = \ \bigg\{ \ (\,y_1\,,\ y_2\,, \ \cdot\,\cdot \cdot\,, \ y_{\,n}\,)  \,\in\,B_{\,\xi_{\,1}} \,(\,2\,\rho_{\,\nu\,}\,)  \ \ | \ \ F\,(\,y_1\,,\ y_2\,, \ \cdot\,\cdot \cdot\,, \ y_{\,n}\,) \ = \ 0 \ \bigg\}\,. \leqno (4.25) \ \
$$
This can be done ``\, uniformly\," in teams of size of the ball $\,B_{\,\xi_j} \,(\,2\,\rho_{\,\nu\,}\,)\,,\,$ $\,1 \ \le \ j \ \le \ \flat\,,\,$ as we work under the compactness assumption in (\,1.4\,)\,.\,
Moreover, there is loss of generality   to select $\,F\,$ so that
$$
\btd\,F\,(\,0\,) \ \ {\mbox{is \ \ along \ \ the \ \ positive }}\  y_{\,n}{\mbox{\,\,-\,axis}} \ \ {\mbox{direction\,.}} \leqno (4.26)
$$
By performing a rescaling (\, $F \ \to \ c \cdot F\,$)\,,\, we may also assume that
$$
\Vert\, \btd\,F\,(\,0\,)  \,\Vert \ = \ \Vert\, (\,0\,,\ \cdot \cdot \cdot\,, \ 0\,,\ 1\,)\,\Vert \ = \ 1\,.\leqno (4.27)
$$
Furthermore, via continuity and the compactness condition (\,1.4\,)\,,\, we can take it that  the $n$\,-\,th component of $\,\btd\,F\,(\,{\bf p}\,)\,$ satisfies
$$
\btd\,F\,(\,{\bf p}\,)\,|_{\,n} \ \not= \ 0 \ \ \ \ \mfor \ \ {\mbox{all}} \ \ \ \ \emph{} {\bf p} \ \in \ B_o\,(\,2\,\rho_\nu\,)\ . \leqno (4.28)
$$
For a fixed $\,y\,\in\,  B_o\,(\,\rho_{\,\nu}\,)$\,,\, let $\,{\bf p}_{\,y}\,$ be the projection of $\,y\,$ unto  $\,{\cal H}\,$ [\,cf. (\,1.25\,)\,]\,.\, By (\,4.27\,) and choosing $\rho_{\,\nu}$ to be small [ \ uniformly so because of (\,1.4\,)\ ]\,,\, we have
$$
{\bf p}_{\,y} \,\in B_o\,(\,2\,\rho_{\,\nu}\,) \ \ \  \ \mfor \ \ y \,\in\,B_o\,(\,\rho_{\,\nu}\,)\ .
$$

 In terms of  $\ {\bf p}_{\,y} \ = (\,{\bf p}_{\,y_{\,|\,1}} \,,\ {\bf p}_{\,y_{\,|\,2}} \,, \ \cdot\,\cdot \cdot\,, \ {\bf p}_{\,y_{\,|\,n}} \,)  \,  \in \,  B_o\,(\,2\,\rho_{\,\nu}\,)\,,$\,
we apply the multi\,-\,variable Taylor expansion at the origin [\, making use of $\,F\,(\,0\,)\,=\,0\,$,\,  (\,4.26\,) and (\,4.27\,)\ ]
$$
F\,(\,{\bf p}_{\,y_{\,|\,1}} \,,\ {\bf p}_{\,y_{\,|\,2}} \,, \ \cdot\,\cdot \cdot\,, \ {\bf p}_{\,y_{\,|\,n}} \,)   \ = \ (\,+\,1\,) \,\cdot {\bf p}_{\,y_{\,|\,n}} \ + \ {1\over 2}\,\times\!\!\sum_{1\ \le \ i\,,\ j\ \le \ n} \!\!a_{i\,j}\   {\bf p}_{\,y_{\,|\,i}} \, {\bf p}_{\,y_{\,|\,j}} \  + \ {\vec{\,O}}\,(\,\Vert\ {\bf p}_{\,y}\,\Vert^{\,3}\,)\,.\leqno (4.29) $$
Here $\ \displaystyle{a_{\,i\,j} \ = \ {{\partial^2\,F}\over {\partial p_{\,i}\,\partial p_{\,j}} }\ \bigg\vert_{\,0\,}}\  $ is the entry in the Hessian matrix for $\,F\,$ at the origin\,,\, and $\,{\vec{\,O}}\,(\,\Vert\ {\bf p}_{\,y}\,\Vert^{\,3}\,)\,$\\[0.1in]
\noindent is a point in $\,\R^n\,$ with the absolute value of each component of the order $\,{{\,O}}\,\left(\,\Vert\ {\bf p}_{\,y}\,\Vert^{\,3}\,\right)$\,.\bk
%
%
%
%
For a fixed $\,y\,\in\,  B_o\,(\,\rho_{\,\nu}\,)$\,,\, via  the Lagrange multiplier method and (4.29), we obtain
\begin{eqnarray*}
(4.30) \ \ y \ - \ {\bf p}_{\,y} & = & \chi\,\cdot\,\btd_{\bf p}\, F\,|_{\,{\bf p}_{\,y}} \ = \ \chi \,\left( {{\partial F}\over {\partial \,p_1}}\,, \ {{\partial F}\over {\partial \,p_2}}\,, \ \,\cdot\,\cdot\,\cdot\,\,, \ {{\partial F}\over {\partial \,p_{\,n}}} \ \right)\bigg\vert_{\ {\bf p}_{\,y}} \\[0.2in]
&  =&   \chi \,\cdot\, \left(\ \sum_{1\,\le \,j\,\le \,n}\,a_{1\,j}\   {\bf p}_{\,y_{\,|\,j}} \,, \ \ \  \sum_{1\,\le \,j\,\le \,n}\,a_{2\,j}\   {\bf p}_{\,y_{\,|\,j}} \,, \  \ \cdot\,\cdot \cdot\ , \ \ \  1 \, +   \sum_{1\,\le \,j\,\le \,n}\,a_{n\,j}\  {\bf p}_{\,y_{\,|\,j}} \ \right) \\[0.2in]
& \ & \ \ \ \ \ \ \ \ \  + \ |\,\chi\,|\,\cdot\,{\vec{\,O}}\,\left(\ \Vert\ {\bf p}_{\,y}\,\Vert^2\,\right)\ .
\end{eqnarray*}
In the above, $\,\chi\,$ denotes the multiplier (\,in general\,,\, it depends on $\,y\,$)\,.\, The constraint is given by the equation [\ for $\,{\bf p}_{\,y} \ = \ (\,{\bf p}_{\,y_{\,|\,1}} \,,\ {\bf p}_{\,y_{\,|\,2}} \,, \ \cdot\,\cdot \cdot\,, \ {\bf p}_{\,y_{\,|\,n}} \,) \ \in \ {\cal H} \,\cap\,B_o\,(\,2\,\rho_{\,\nu}\,) $\ ] \\[0.1in]
(4.31)
$$
F\,(\,{\bf p}_{\,y_{\,|\,1}} \,,\ {\bf p}_{\,y_{\,|\,2}} \,, \ \cdot\,\cdot \cdot\,, \ {\bf p}_{\,y_{\,|\,n}} \,)   \ = \  0 \ \ \Longrightarrow \ \ {\bf p}_{\,y_{\,|\,n}} \ + \  \ {1\over 2}\,\cdot \!\!\!\sum_{1\ \le \ i,\ j\ \le \ n} \!\!a_{i\,j}\   {\bf p}_{\,y_{\,|\,i}} \, {\bf p}_{\,y_{\,|\,j}} \ \, + \ {\vec{\,O}}\,(\ \Vert\ {\bf p}_{\,y}\,\Vert^{\,3}\,)  \ = \ 0\ .
$$
Based on (4.30)\,,\, we obtain
\begin{eqnarray*}
y_1 \, - \, {\bf p}_{\,y_{\,|\,1}} \!\!& = & \chi\,\cdot\,\left(\ \, a_{\,11}\,{\bf p}_{\,y_{\,|\,1}}  \, + \,  a_{\,12}\,{\bf p}_{\,y_{\,|\,2}}  \, + \,\,\cdot\,\cdot\,\cdot\, \ + \ a_{\,1n}\,{\bf p}_{\,y_{\,|\,n}}  \, \right) \ + \ |\,\chi\,| \cdot O\,\left(\  \Vert\,{\bf p}_{\,y}\, \Vert^2 \,\right) \\[0.2in]
& \ &  \hspace*{-1.5in}   \Longrightarrow \  \ \ \ \  \left(\ 1 \ + \ \chi
\cdot \,a_{\,11}    \ \right)\,\cdot\,{\bf p}_{\,y_{\,|\,1}} \, + \,  [\,\chi
\cdot a_{\,12}\,]\cdot {\bf p}_{\,y_{\,|\,2}} \, + \,\,\cdot\,\cdot\,\cdot\,\ + \ [\,\chi
\cdot a_{\,1n}\,]\,\cdot\,{\bf p}_{\,y_{\,|\,n}}
\\[0.2in]&   \ & \hspace*{-1.2in} \ \  = \   y_1  \ + \ |\,\chi\,| \cdot O\,\left(\  \Vert\,{\bf p}_{\,y}\, \Vert^2 \,\right)\,,\\[0.01in] &:&\\[0.01in]&:&\\[0.01in]
y_{n\,-\,1} \, - \, {\bf p}_{\,y_{\,|\,n\,-\,1}}  & = & \chi \cdot\left( \, a_{\,(n\,-\,1)\,1}\,{\bf p}_{\,y_{\,|\,1}}  \, + \,  a_{\,(n\,-\,1)\,2}\,{\bf p}_{\,y_{\,|\,2}} \, + \,\,\cdot\, \cdot\,\cdot\,\ + \ a_{\,(n\,-\,1)\,n}\,{\bf p}_{\,y_{\,|\,n}}   \, \right)\\[0.1in]
& \ &  \hspace*{-1.5in}\Longrightarrow  \  \  \ \ \ \  [\,\chi\cdot a_{\,(n\,-\,1)\,1}\,]\,\cdot\,{\bf p}_{\,y_{\,|\,1}}  \ + \ \left(\,1 \ + \ \chi\,\cdot\,\,a_{\,(n\,-\,1)\,2}\,\right) \ {\bf p}_{\,y_{\,|\,2}}  \, + \, [\,\chi\cdot a_{\,(n\,-\,1)\,3}\,]\,\cdot\,{\bf p}_{\,y_{\,|\,3}}  \, +  \,\cdot\,\cdot\,\cdot\,\,+ \,  \\[0.2in]&  \  & \hspace*{-0.75in}\ + \  [\,\chi\,\cdot\,a_{\,(n\,-\,1)\,n}\,]\,\cdot\,{\bf p}_{\,y_{\,|\,n}} \\[0.2in]&   \ & \hspace*{-1.2in} \ \  = \     y_{n\,-\,1}  \ + \ |\,\chi\,| \cdot O\,\left(\  \Vert\,{\bf p}_{\,y}\, \Vert^2 \,\right)\,,
\end{eqnarray*}
\begin{eqnarray*}
y_{\,n} \, - \, {\bf p}_{\,y_{\,|\,n}} & = & \chi \cdot\left( \, 1 \ + \    a_{\,n\,1}\,{\bf p}_{\,y_{\,|\,1}}  \ + \    a_{\,n\,2}\,{\bf p}_{\,y_{\,|\,2}} \, + \,\,\cdot\, \cdot\,\cdot\,\ + \ a_{\,n\,n}\,{\bf p}_{\,y_{\,|\,n}}   \, \right)  \ + \ |\,\chi\,| \cdot O\,\left(\  \Vert\,{\bf p}_{\,y}\, \Vert^2 \,\right)\\[0.1in]
& \ &  \hspace*{-1.5in} \Longrightarrow \  \ \ \  [\, \chi\,
\cdot \,a_{\,n\,1}\,]\,\cdot\,{\bf p}_{\,y_{\,|\,1}}  \, + \,  [\,\chi
\cdot a_{\,n\,2}\,]\cdot {\bf p}_{\,y_{\,|\,n}} \, + \,\,\cdot\,\cdot\,\cdot\,\ + \  [\, \chi\,
\cdot \,a_{\,n\,(n\,-\,1)}\,]\,\cdot\,{\bf p}_{\,y_{\,|\,(n\,-\,1)}}  \\[0.2in]&  \  & \hspace*{-0.75in}\ + \  [\,1 \ + \ \chi
\cdot a_{\,n\,n}\,]\,\cdot\,{\bf p}_{\,y_{\,|\,n}}    \  = \   (\,y_{\,n} \ - \ \chi\,)  \ + \ |\,\chi\,| \cdot O\,\left(\  \Vert\,{\bf p}_{\,y}\, \Vert^2 \,\right)\ .\\
\end{eqnarray*}
Cramer's rule shows that the inverse of a matrix of the form \\[0.1in]
{\scriptsize{\[  {\large{ }}     \ \  \ \ \  \ \ \  \ \ \  \ \ \  \ \  \  \ \ \  \ \ \  \ \ \  \left(\ \begin{array}{cccc}
1 \ + \ O\,(\ |\,\chi\,|\ ) &    O\,(\ |\,\chi\,|\ )  & {\,\cdot\,\ \,\cdot\,\cdot\,\cdot\,} & {\cdot  \ \,\cdot\,\cdot \ \  O\,(\ |\,\chi\,|\ )  }\\[0.15in]
O\,(\ |\,\chi\,|\ )   &  1 \ + \ O\,(\ |\,\chi\,|\ ) & {   \,\cdot\,\cdot\,\cdot\,} & {  \ \,\cdot\,\cdot\,\cdot\,\ \  O\,(\ |\,\chi\,|\ )  }\\[0.15in]
\cdot & \\[0.15in]
\cdot & \\[0.15in]
O\,(\ |\,\chi\,|\ )   &  O\,(\ |\,\chi\,|\ )  &  {\,\cdot\,\cdot\,\cdot\,} &  \,\cdot\,\cdot\,\cdot\,\ \ 1 \ + \ O\,(\ |\,\chi\,|\ )   \\[0.15in]
\end{array} \right)\ \ \  \ \ \  \ \ \  \ \ \  \ \ \  \ \ \  \ \ \  \ \ \  \ \ \  \ \ \  \ \ \ \ \ \  \ \ \  \ \ \  \ \ \  \ \ \  \ \ \  \ \ \  \ \ \  \ \ \  \ \ \  \ \ \   \]}}

\vspace*{-1in}

\hspace*{5in}[ \ $ O\,(\,|\,\chi\,|\ )  \ \approx \ 0$ \ ]

\vspace*{.75in}

has similar expression [\,via (\,4.23\,) and (\,4.31\,)\,,\, \,,\, see also (\,4.33\,)\ ]\,.\, It follows that
 \begin{eqnarray*}(\,4.32\,) \ \ \ \ \ \ \ \ \ \ \ \ \ \  \ \ \ \ \ \ \ \ \ \ \ \ \ \
{\bf p}_{\,y_{\,|\,1}} & = & y_1   \ + \ O\,(\ |\,\chi\,|\ )\,\cdot\,O\,(\,\Vert\ y \,\Vert\ )\ ,   \\[0.01in]&:&\\[-0.1in]&:&\\[0.01in]
{\bf p}_{\,y_{\,|\,n\,-\,1}}  & = & y_{n\,-\,1}   \ + \ O\,(\ |\,\chi\,|\ )\,\cdot\,O\,(\,\Vert\ y \,\Vert\ )\ ,  \\[0.2in]
{\bf p}_{\,y_{\,|\,n}} & = & (\, y_{\,n} \ - \ \chi\,) \ + \  O\,(\ |\,\chi\,|\ )\,\cdot\,O\,(\,\Vert\ y \,\Vert\ ) \ .\ \ \ \ \ \ \ \ \ \ \ \ \ \ \ \ \ \ \ \ \   \ \ \ \ \ \ \ \ \ \ \ \ \ \
\end{eqnarray*}
From (\,4.32\,)\,,\, we obtain
$$
 y_{\,n}  \ - \ {\bf p}_{\,y_{\,|\,n}} \ = \ \chi \ + \ O\,(\ |\,\chi\,|\ )\,\cdot\,O\,(\,\Vert\ y \,\Vert\ ) \ , \leqno (\,4.33\,)
$$
and
 \begin{eqnarray*}
(\,4.34\,) \ \ \ \ \ \ \  & \ & {\bf p}_{\,y_{\,|\,1}}^2 \ + \ \cdot\,\cdot\,\cdot\,\ + \ {\bf p}_{\,y_{\,|\,n\,-\,1}}^2 = \ y_1^2 \ + \ \cdot\,\cdot\,\cdot\,\ + \ y_{n\,-\,1}^2 \ +\ O\,(\ |\,\chi\,|\ )\,\cdot\,O\,(\,\Vert\ y \,\Vert^2\,) \ \ \ \ \ \ \ \ \ \ \ \\[0.2in]
& \ &  \Longrightarrow \ \  \ \
\Vert\,\sigma_{\,{\bf p}}\,\Vert^2 \ = \ \Vert\,\sigma_{\,y} \,\Vert^2 \ + \ O\,(\ |\,\chi\,|\ )\,\cdot\,O\,(\,\Vert\ y \,\Vert^2\,)\ ,\\[0.2in]
& \ & \hspace*{-1.18in} {\mbox{where}} \ \ \ \ \sigma_{\,{\bf p}}\ =\ \left(\ {\bf p}_{\,y_{\,|\,1}}\,, \  \cdot\,\cdot\,\cdot\,, \ {\bf p}_{\,y_{\,|\,n\,-\,1}}\ \right)   \ \ \ \ \ \ {\mbox{and}} \ \ \ \ \ \    \ \ \sigma_{\,y}\ =\ \left(\ y_1\,, \ \cdot\,\cdot\,\cdot\,,\ y_{n\,-\,1}\ \right)\ . \ \ \ \ \ \ \ \ \ \ \
\end{eqnarray*}
Moreover, the constraint equation (\,4.32\,)  shows that
\begin{eqnarray*}
& \ & {\bf p}_{\,y_{\,|\,n}} \ + \ {1\over 2}\,\cdot\,\sum_{1\,\le\, i\,,\,j\,\le\, n}\!\!a_{i\,j}\  {\bf p}_{\,y_{\,|\,i}}\ {\bf p}_{\,y_{\,|\,j}} \  = \ O\,(\,\Vert\ {\bf p}\,\Vert^3\,)\\[0.15in]
\Longrightarrow &  \ &  {\bf p}_{\,y_{\,|\,n}}\,\cdot\,\left(\ \ 1 \ + \ \sum_{1\,\le\, i\,\le\, n}  \!\!a_{n\,i}\   {\bf p}_{\,y_{\,|\,i}}   \ \right) \ =   \ - \, \sum_{i\,,\,j \ \not= \ n} \!\!a_{i\,j}\ {\bf p}_{\,y_{\,|\,i}}\ {\bf p}_{\,y_{\,|\,j}}  \ + \  O\,(\,\Vert\ {\bf p}\,\Vert^3\,)\\[0.15in]
\Longrightarrow &  \ &  {\bf p}_{\,y_{\,|\,n}} \ = \  \ {{ \ -  {\displaystyle{\sum_{i\,,\,j \ \not= \ n} }} \,a_{i\,j}\ {\bf p}_{\,y_{\,|\,i}}\ {\bf p}_{\,y_{\,|\,j}} \  }\over  { 1 \ + \  {\displaystyle{\sum_{1\,\le\, i\,\le\, n} }} a_{n\,i}\   {\bf p}_{\,y_{\,|\,i}}     }} \ \  + \  O\,(\,\Vert\ {\bf p}\,\Vert^3\,)\\[0.15in]
\Longrightarrow &  \ &  {\bf p}_{\,y_{\,|\,n}} \ = \   - \, \sum_{i\,,\,j \ \not= \ n} \,a_{i\,j}\ {\bf p}_{\,y_{\,|\,i}}\ {\bf p}_{\,y_{\,|\,j}}    \ + \  O\,(\,\Vert\ {\bf p}\,\Vert^3\,)\\[0.1in]
& \ &  \ \ \ \ \ \ \ \ \ \ \ \ \  \ \ \ \ \ \ \ \ \ \ \ \ \  \ \ \ \ \    \left\{ \ \ {\mbox{as}} \ \ \ \ \  {{ 1 }\over  { \ 1 \ + \  {\displaystyle{\sum_{1\,\le\, i\,\le\, n} }} a_{n\,i}\   {\bf p}_{\,y_{\,|\,i}}    \  }} \ = \ 1 \ +\ O\,(\,\Vert\ {\bf p}\,\Vert \,) \ \right\} \\[0.1in]
\Longrightarrow &  \ &  {\bf p}_{\,y_{\,|\,n}} \ = \ O\,(\,\sigma_{\,{\bf p}}^2 \,)   \ + \  O\,(\,\Vert\ {\bf p}\,\Vert^3\,)\\[0.2in]
(\,4.35\,)\ \cdot \cdot \  \Longrightarrow &  \ &  {\bf p}_{\,y_{\,|\,n}} \ = \   O\,\left(\,\Vert\ y\,\Vert^2\,\right) \ \ \ \ \ \ \  \ \ \ \ \ {\mbox{for}} \ \ \ \ y \,\in\,B_o\,(\,\rho_{\,\nu}\,) \ \ \ \ \ \ \ \ \ \ [ \ {\mbox{via}} \ \ (\,4.34\,)\ ] \\[0.2in]
& \ & \hspace*{-1.4in}  \bigg\{  \ {\mbox{note \ \ that \ \ }} \ \Vert\ {\bf p}\,\Vert^3 \ =\ \Vert\ {\bf p}\,\Vert\,\cdot\,(\,\sigma_{\,{\bf p}}^2 \ + \ {\bf p}_{\,y_{\,|\,n}}^2\,) \ = \ \Vert\ {\bf p}\,\Vert\,\cdot\, [\,\sigma_{\,y}^2 \ + \ o\,(\,\Vert\ y\,\Vert^2\,)\,] \ + \ o\,(\ {\bf p}_{\,y_{\,|\,n}}^2)\  \\[0.2in]
& \ & \ \ \ \ = \ \Vert\ {\bf p}\,\Vert\,\cdot\, [\ O\,(\,\Vert\ y\,\Vert^2\,)\,] \ + \ o\,(\ {\bf p}_{\,y_{\,|\,n}}^2)\  \bigg\}\ . 
\end{eqnarray*}



\hspace*{0.3in}

{\bf \large {\bf \S\,4\,d\,.}} \  \
\noindent{\it Expressing $\ \Vert\, y \ - \ {\bf p}_{\,y}\,\Vert\,$ in terms of $\ y_{\,n}\,.$}\ \  From (\,4.28)\,,\, for $\,y \,\in\,B_{o}(\,\rho_\nu\,)\,\setminus\,{\cal H}\,,\,$ we have
$$
y_{\,n} \ - \ {\bf p}_{{\,y\,|_n}} \ \not= \ 0\,.
$$
For if $\,y_{\,n} \ = \ {\bf p}_{{\,y\,|_n}}\,$,\, then the vector $\, y \, - \,  {\bf p}_{{\,y}}\,$ is perpendicular to the $\,y_{\,n}\,$ axial direction\,.\, But $\, y \, - \,  {\bf p}_{{\,y}}\,$  is also along the direction of the normal to $\,{\cal H}\,$ at $\,{\bf p}_{{\,y}}\,$,\, that is\,,\, $\, y \, - \,  {\bf p}_{{\,y}}\,=\,c \cdot \btd\, F\,(\,{\bf p}_{\,y}\,)\,$ for a non\,-\,zero number $\,c$\,,\,  implying $\,\btd\, F\,(\,{\bf p}_{\,y}\,)\ = \ 0\,,\,$  contradicting (\,4.28\,) when $\,\rho_{\,\nu}\,$ is made small enough [\ uniformly so chosen\,,\,  in view of the compactness condition (\,1.4\,)\ ]\,.\, Recall that
\begin{eqnarray*}
y & = & (\,y_1\,,\, \cdot \cdot \cdot\,,\ y_{n\,-\,1}\,, \ y_{\,n}\,)\ , \ \ \ \ \ \sigma_{\,y}\ =\ \left(\ y_1\,, \ \cdot\,\cdot\,\cdot\,,\ y_{n\,-\,1}\ \right)\ , \ \ \ \  \\[0.1in]
{\mbox{and}} \ \ \ \ \ \ \ \ \ \   \ \ \ \ \ \ \ \ \ \ \ \ \ {\bf p}_{{\,y}} & = & (\, {\bf p}_{{\,y\,|_1}}\,,\, \cdot \cdot \cdot\,,\  {\bf p}_{{\,y\,|_{n\,-\,1}}}\,, \  {\bf p}_{{\,y\,|_n}}\,)\ ,  \ \ \  \sigma_{\,{\bf p}}\ =\ \left(\ {\bf p}_{\,y_{\,|\,1}}\,, \  \cdot\,\cdot\,\cdot\,, \ {\bf p}_{\,y_{\,|\,n\,-\,1}}\ \right)  \ . \ \ \ \ \ \
\end{eqnarray*}
Via  Pythagoras' theorem, for $\,y \,\in\,B_{o}(\,\rho_\nu\,)\,$  and $\,y \,\not\in\, {\cal H}\,,\,$
\begin{eqnarray*}
\Vert \,y \ - \ {\bf p}_{\,y}\,\Vert^{\,2} & = & (\,y_{\,n} \ - \ {\bf p}_{\,y\,|_n}\,)^{\,2} \ + \ \Vert \ \sigma_{\,y} \ - \ \sigma_{\,{\bf p}_{\,y}} \Vert^{\,2}\ . \\[0.2in]
(\,4.36\,) \ \cdot \cdot \cdot \cdot \ \Vert \,y \ - \ {\bf p}_{\,y}\,\Vert & = & |\,y_{\,n} \ - \ {\bf p}_{{\,y\,|_n}}\,|\, \cdot \,\sqrt {\ 1 \ + \ {{\Vert \,\sigma_{\,y} \ - \ \sigma_{\,{\bf p}_{\,y}} \Vert^{\,2} }\over {(\ y_{\,n} \ - \ {\bf p}_{{\,y\,|_n}}\ )^{\,2}  }} \ \  }  \\[0.15in]
& = & |\,y_{\,n} \ - \ {\bf p}_{{\,y\,|_n}}\,| \cdot \sqrt {\ 1 \ + \ {{\  O\,(\,|\,\chi\,|^{\,\,2}\  ) \cdot O\,(\,\Vert\,y\,\Vert^{\,2}\,)\ }\over {(\,y_{\,n} \ - \ {\bf p}_{{\,y\,|_n}}\,)^{\,2}  }} \ }  \ \ \ \ \  \ \  [ \ {\mbox{via}} \ \ (\,4.32\,) \ ] \ \ \ \ \ \ \\[0.15in]
&   = & |\,y_{\,n} \ - \ {\bf p}_{{\,y\,|_n}}\,| \cdot  \sqrt{ \  1 \ + \ O\,(\,\Vert\,y\,\Vert^{\,2}\,) \ }  \ \ \ \ \  \ \ \ \ \ \ \ \ \ \ \ \ \ \ \ \ \ \   [ \ {\mbox{from}} \ \ (\,4.33\,) \ ] \ \ \ \ \ \ \\[0.15in]
&   = & |\,y_{\,n} \ - \ {\bf p}_{{\,y\,|_n}}\,|  \ + \ |\,y_{\,n} \ - \ {\bf p}_{{\,y\,|_n}}\,|\,\cdot  O\,(\,\Vert\,y\,\Vert^{\,2}\, ) \ .
\end{eqnarray*}
From (4.29)\,,\, we obtain

 \vspace*{-0.15in}

$$
|\, {\bf p}_{{\,y\,|_n}}\,| \ = \ O \,( \ \Vert \,y\,\Vert^2 \ )  \,\le\   {\bar C}_1 \cdot \rho_{\,\nu}^2\ \ \ \ \ \ \ \ \ \ \ {\mbox{for}} \ \ \ \ y \,\in\,B_o(\,\rho_{\,\nu}\,)\ . \leqno (\,4.37\,)
$$
Here we observe the convention to delete the positive constant in front of $\,\rho_{\,\nu}^2\,,\,$ as explained in Remark 1.23\,.\,
Let
$$
 \Delta \ = \   {\bar C}_1 \cdot  \rho_{\,\nu}^2  \ + \  \xi_{1_{|_n}}  \ \ \ \ \
 ( \ {\mbox{recall \ \ that}} \ \ \xi_{1_{|_n}} \ > \ 0\ ) \ . \leqno(\,4.38\,)
 $$
Note that\,,\, via (\,4.37\,)\,,\,
$$
y_{\,n} \ > \ \Delta \ \ \Longrightarrow \ \  |\,y_{\,n} \ - \ {\bf p}_{{\,y\,|_n}}\,|   \ = \  y_{\,n} \ - \ {\bf p}_{{\,y\,|_n}} \ >  \ 0 \ \ \  \ {\mbox{and}} \ \ \ \ \  y_{\,n} \ - \ \xi_{\,1\,|_n} \ > \ 0\,.
$$
Following (\,4.36\,)\,,\, for $\,y \,\in\,B_{o}(\,\rho_\nu\,)\ \setminus \,{\cal H}\,$ with $\, y_{\,n} \ > \ \Delta\,,\,$ we continue with
 \begin{eqnarray*}
 (\,4.39\,) \ \ \ \  \ \ \ \  \Vert \,y \ - \ {\bf p}_{\,y}\,\Vert & = & (\,y_{\,n} \ - \ {\bf p}_{{\,y\,|_n}}\,) \ + \ |\,y_{\,n} \ - \ {\bf p}_{{\,y\,|_n}}\,|\,\cdot  O\,(\,\Vert\,y\,\Vert^{\,2}\, ) \\[0.15in]
  & = &  (\ y_{\,n} \ - \ \xi_{\,1\,|_n}\ ) \ + \  (  \  \xi_{\,1\,|_n} \ - \  {\bf p}_{{\,y\,|_n}}\ )  \ + \ |\,y_{\,n} \ - \ {\bf p}_{{\,y\,|_n}}\,|\,\cdot  O\,(\,\Vert\,y\,\Vert^2\,)\\[0.2in]
  & = & (\,y_{\,n} \ - \ \xi_{1_{\,|_{\,n}}}\,) \ + \ \xi_{1_{\,|_{\,n}}} \ + \ \left[ \ O\, \left(\,\xi_{1_{\,|_{\,n}}}^2 \,\right) \ + \ O\,\left(\, \Vert \, y \ - \ \xi_{\,1}\,\Vert^2\,\right)\ \right]\ \ \ \ \ \ \ \ \ \\
 \end{eqnarray*}
 \hspace*{1.9in}\ \ \ \ $\leftarrow \ \, {\mbox{+\,ve}} \ \, \rightarrow $ \ \ \ \ \ \ \ \ \  [\ for \ \ $\,y \,\in\,B_{o}(\,\rho_\nu\,)\ \setminus \,{\cal H}\,$ with $\, y_{\,n} \ > \ \Delta\,$]\,.

\vspace*{0.2in}

 Here we make use of (\,4.37\,)\,:
 \begin{eqnarray*}
& \ & |\, {\bf p}_{\,y_{\,|\,n}}\,| \ = \   O\,(\,\Vert\ y\,\Vert^2\,)\\[0.2in]
\Longrightarrow & \ &  |\, {\bf p}_{\,y_{\,|\,n}}\,|  \ +  \ |\,y_{\,n} \ - \ {\bf p}_{{\,y\,|_n}}\,|\,\cdot  O\,(\,\Vert\,y\,\Vert^2\, ) \\[0.2in]
& = &  O\,(\,\Vert\ y\,\Vert^2\,) \ = \  O\, \left(\,\Vert \,\xi_{\,1}\,\Vert^2 \,\right) \ + \ O\,\left(\, \Vert \, y \ - \ \xi_{\,1}\,\Vert^2\,\right) \\[0.2in]
&  = & O\, \left(\,\xi_{1_{\,|_{\,n}}}^2 \,\right) \ + \ O\,\left(\, \Vert \, y \ - \ \xi_{\,1}\,\Vert^2\,\right)\
 \\[0.1in]
& \ &\hspace*{-0.5in}  [ \ {\mbox{as}} \ \ \ \Vert\,y\,\Vert \ = \ \Vert\,(\,y\,-\,\xi_{\,1}\,) \ + \ \xi_{\,1} \,\Vert  \ \le \ \Vert \, y \ - \ \xi_{\,1}\,\Vert \ + \ \Vert\,\xi_{\,1} \,\Vert \ = \ \Vert \, y \ - \ \xi_{\,1}\,\Vert \ + \ \Vert\,\xi_{1_{\,|_{\,n}}}\,\Vert \ ]\  .
\end{eqnarray*}
Recall that\,,\, in the present setting,
$$
\xi_{\,1}\ = \ \left(\ 0\, \ \cdot\,\cdot \cdot\,, \ 0\,, \ \xi_{1_{\,|_{\,n}}}\,\right) \ \ \Longrightarrow \  \ \Vert \,\xi_{\,1}\,\Vert^2 \ = \ \xi_{1_{\,|_{\,n}}}^2\ .
$$
 Likewise\,,  for $\,y \,\in\,B_{o}(\,\rho_\nu\,)\ \setminus \,{\cal H}\,$ with $\, y_{\,n} \ < \ -\,\Delta\,,\,$
\begin{eqnarray*}
(\,4.40\,) \ \ \ \ \  \   \Vert \,y \ - \ {\bf p}_{\,y}\,\Vert & = & (\,-\,y_{\,n} \ + \ {\bf p}_{{\,y\,|_n}}\,) \ + \ |\,y_{\,n} \ - \ {\bf p}_{{\,y\,|_n}}\,|\,\cdot  O\,(\,\Vert\,y\,\Vert^{\,2}\, ) \\[0.15in]
  & = & -\, (\ y_{\,n} \ - \ \xi_{\,1\,|_n}\ ) \ + \  ( \   {\bf p}_{{\,y\,|_n}} \ - \ \xi_{\,1\,|_n} \ )  \ + \ |\,y_{\,n} \ - \ {\bf p}_{{\,y\,|_n}}\,|\,\cdot  O\,(\,\Vert\,y\,\Vert^2\, ) \\[0.2in]
  & = & -\,(\,y_{\,n} \ - \ \xi_{1_{\,|_{\,n}}}\,) \ -\, \xi_{1_{\,|_{\,n}}} \ + \ \left[\ O\, \left(\,\xi_{1_{\,|_{\,n}}}^2 \,\right) \ + \ O\,\left(\, \Vert \, y \ - \ \xi_{\,1}\,\Vert^2\,\right)\ \right] \ \ \ \ \  \ \ \ \ \   \ \ \ \ \  \ \ \ \ \   \ \ \ \ \  \ \ \ \
 \end{eqnarray*}
 \hspace*{1.9in}\ \ \ \ $\leftarrow \ \, {\mbox{+\,ve}} \ \, \rightarrow $ \ \ \ \ \ \ \ \ \  [\ for \ \ $\,y \,\in\,B_{o}(\,\rho_\nu\,)\ \setminus \,{\cal H}\,$ with $\, y_{\,n} \ < \ -\, \Delta\ $]\,.\\[0.1in]
\noindent It follows that (\,when $\,\ell \ = \ 2\,$)\,,\, for $\,y \,\in\,B_{o}(\,\rho_\nu\,)\ \setminus \,{\cal H}\,$ with $\, y_{\,n} \ > \ \Delta\ , $\,
\begin{eqnarray*}
(\,4.41\,) \ \ \ \  \ \ \ \ \ \   \Vert \,y  \ - \ {\bf p}_{\,y}\,\Vert^2 & = &  (\,y_{\,n} \ - \ \xi_{1_{\,|_{\,n}}}\,)^{\,2} \ \ \ \ \left(\ \leftarrow \ \ = \ |\, y_{\,n} \ - \ \xi_{1_{\,|_{\,n}}}\,|^{\,\,2} \ \right) \\[0.2in]
  & \ &  + \ 2\,\cdot\,(\,y_{\,n}\,-\,\xi_{1_{\,|_{\,n}}}\,)\cdot \xi_{1_{{\,n}}}   \ + \ \xi_{1_{\,|_{\,n}}}^2  \\[0.2in]
  & \ & \ \ \ \ \ \ \ \  + \ O\,(\ \Vert\,y \ - \ \xi_{\,1}\,\Vert^3\ )  \ + \ O\,(\ |\,\xi_{1_{\,|_{\,n}}}\,|^3\ ) \ . \  \ \ \ \ \ \ \ \ \ \ \ \ \ \ \ \  \ \ \ \ \ \ \ \ \ \ \ \ \ \ \ \
\end{eqnarray*}

Here we apply the inequality
$$
A \cdot B \ \le \ {{A^p}\over p} \ + \ {{B^q}\over q} \ \ \ \ \ \ {\mbox{for}} \ \ A\,, \ B\,, \ p \ \,\&\, \ q \ > \ 0 \ \ \ \ {\mbox{with}} \ \ {1\over p} \ + \ {1 \over q} \ = \ 1\,. \leqno (\,4.42\,)
$$


Likewise,  for $\,y\,\in\,B_o\,(\,\rho_\nu\,)\,$ with $\,y_{\,n} \ <  \  -\,\Delta\,,\,$ and $\,\ell \ = \ 2\,,\,$ we have
\begin{eqnarray*}
(\,4.43\,) \ \ \ \  \ \ \ \ \  \Vert \,y  \ - \ {\bf p}_{\,y}\,\Vert^2 & = &  (\,-\,y_{\,n} \ + \ \xi_{1_{\,|_{\,n}}}\,)^2
  \ \ \ \ \left(\ \leftarrow \ \ = \ |\ y_{\,n} \ - \ \xi_{1_{\,|_{\,n}}}\,|^{\,\,2} \ \right) \\[0.2in]
  & \ &  + \ 2\,\cdot\,(\,-\,y_{\,n} \ + \ \xi_{1_{\,|_{\,n}}}\,)^{\ell\,-\,1}\cdot (\,-\,\xi_{1_{{\,n}}}  \,)  \ + \ \xi_{1_{\,|_{\,n}}}^2   \\[0.15in]
  & \ & \ \ \ \ \ \ \ \
   + \ O\,(\ \Vert\,y \ - \ \xi_{\,1}\,\Vert^3\ )   \ + \ O\,(\,|\,\xi_{1_{\,|_{\,n}}}\,|^3\,)  \ .\ \ \ \  \ \ \ \ \ \ \ \ \ \ \ \ \ \ \   \ \ \  \ \ \ \ \ \ \ \ \ \ \ \ \ \ \
\end{eqnarray*}
When $\,\ell \ \ge\ 3\,,\,$  we first observe that
\begin{eqnarray*}
(\,A \ + \ B\,)^\ell & = &  A^\ell \cdot \left(\, 1 \ + \ {B\over A}\,\right)^\ell\\[0.2in]
& = &  A^\ell \cdot \left[\, 1 \ + \ O\left(\   {B\over A}\, \right)\,\right] \ \ \ \ \  \ \ \ {\mbox{for}} \ \  \ A\,, \ B \ >  \ 0 \ \ {\mbox{with}} \ \ \ \left(\, {B\over A}\,\right) \ \le \ 1\\[0.2in]
& = & A^\ell \ + \ O\,(\,A^{\ell \,-\,1}\cdot B\ )\,.
\end{eqnarray*}

\vspace*{-0.15in}
We take
$$
A \ = \ (\,y_{\,n} \ - \ \xi_{1_{\,|_{\,n}}}\,)  \ + \ \xi_{1_{\,|_{\,n}}} \ = \ y_{\,n}\ \ \ \ {\mbox{and}} \ \ \ \ B\ = \ \ O\, \left(\,\xi_{1_{\,|_{\,n}}}^2 \,\right) \ + \ O\,\left(\, \Vert \, y \ - \ \xi_{\,1}\,\Vert^2\,\right)\ \le \  {\bar C}_2 \cdot \Delta\,.
$$
Thus for $\,y\,\in\,B_o\,(\,\rho_\nu\,) \ \setminus \,{\cal H}\,$ with  $\,y_{\,n} \ \ge \ {\bar C}_2 \cdot \Delta\,,\,$ $\,\ell \ \ge \ 3\,$ (\,simpler expression for $\,\ell \ = \ 2\,$)\,,\,we obtain\\[0.1in]
(\,4.44\,)
\begin{eqnarray*}
 \Vert \,y  \ - \ {\bf p}_{\,y}\,\Vert^\ell & = & [ \ (\,y_{\,n} \ - \ \xi_{1_{\,|_{\,n}}}\,)  \ + \ \xi_{1_{\,|_{\,n}}} \ ]^{\,\ell} \ + \ O\,(\,|\,y_{\,n}\,|^{\,\ell\,-\,1}\,) \,\cdot\, [\  O\, (\,\xi_{1_{\,|_{\,n}}}^2 \,) \ + \ O\,(\, \Vert \, y \ - \ \xi_{\,1}\,\Vert^2\,)\ ]\\[0.2in]
  & = &  (\,y_{\,n} \ - \ \xi_{1_{\,|_{\,n}}}\,)^{\,\ell} \ \ \ \ \left(\ \leftarrow \ \ = \ |\, y_{\,n} \ - \ \xi_{1_{\,|_{\,n}}}\,|^{\,\,\ell} \ \right) \\[0.2in]
  & \ &  + \ \ell\,\cdot\,(\,y_{\,n}\,-\,\xi_{1_{\,|_{\,n}}}\,)^{\,\ell\,-\ 1}\cdot \xi_{1_{{\,n}}}   \ \ \longleftarrow \ \ \ \ \cdot \cdot \cdot \cdot \ (\,*\,)  \\[0.2in]
   &  & \ \ \ \ \  + \  {{ \ell\,\cdot\,(\,\ell\,-\,1\,) }\over {2}}\,\cdot\,(\,y_{\,n} \ - \ \xi_{1_{\,|_{\,n}}}\,)^{\,\ell\,-\ 2}\,\cdot\, \xi_{1_{\,|_{\,n}}}^2  \ \ \longleftarrow \ \ \ \ \cdot \cdot \cdot \cdot \ (\,**\,) \\[0.1in]
   & \ & \hspace*{3in}(\ {\mbox{error \ \ grouped \ \ below}} \ \downarrow \ ) \\[0.1in]
   & \ &\hspace*{-1in} + \ \left[ \  O\,(\ \Vert \, y \ - \ \xi_{\,1}\,\Vert^{\,\ell\,-\ 3}\ ) \,\cdot\,  O\, (\,\xi_{1_{\,|_{\,n}}}^3 \,) \ + \   O\,(\, \Vert \, y \ - \ \xi_{\,1}\,\Vert^{\,\ell\,-\ 4}\ ) \,\cdot\,  O\, (\,\xi_{1_{\,|_{\,n}}}^4 \,) \ + \ \cdot \cdot \cdot \ +  \ \xi_{1_{\,|_{\,n}}}^\ell \ \right] \  + \\[0.15in]
    & \ & \ \ \ \ \ \ + \ \, O\,(\,|\,y_{\,n}\,|^{\,\ell\,-\,1} \,)\,\cdot\, \left[\  O\, (\,\xi_{1_{\,|_{\,n}}}^2 \,) \ + \ O\,(\, \Vert \, y \ - \ \xi_{\,1}\,\Vert^2\,)\ \right]\ .\\
\end{eqnarray*}


Likewise\,,\, for $\,y\,\in\,B_o\,(\,\rho_\nu\,)\ \setminus \, {\cal H}\,$ with $\,y_{\,n} \ \le \ -\, {\bar C}_2 \cdot \Delta\,,\,$ and $\,\ell \ \ge \ 3\,,\,$ we have\\[0.1in]
(\,4.45\,)
\begin{eqnarray*}
 \Vert \,y  \ - \ {\bf p}_{\,y}\,\Vert^\ell & = &  (\,-\,y_{\,n} \ + \ \xi_{1_{\,|_{\,n}}}\,)^{\,\ell}  \ \ \ \ \left(\ \leftarrow \ \ = \ |\ y_{\,n} \ - \ \xi_{1_{\,|_{\,n}}}\,|^{\,\ell} \ \right) \\[0.2in]
  & \ &  + \ \ell\,\cdot\,(\,-\,y_{\,n} \ + \ \xi_{1_{\,|_{\,n}}}\,)^{\ell\,-\,1}\cdot (\,-\,\xi_{1_{{\,n}}}  \,)     \ \ \longleftarrow \ \ \ \ \cdot \cdot \cdot \cdot \ (\,*\,)  \\[0.2in]
  & \ & \ \ \ \ \  + \  {{ \ell\,\cdot\,(\,\ell\,-\,1\,) }\over {2}}\,\cdot\,(\,-\,y_{\,n} \ + \ \xi_{1_{\,|_{\,n}}}\,)^{\ell\,-\ 2}\,\cdot\, \xi_{1_{\,|_{\,n}}}^2  \ \ \longleftarrow \ \ \ \ \cdot \cdot \cdot \cdot \ (\,**\,) \\[0.1in]
   & \ & \hspace*{3in}(\ {\mbox{error \ \ grouped \ \ below}} \ \downarrow \ ) \\[0.1in]
    & \ &\hspace*{-1in} + \ \left[ \  O\,(\ \Vert \, y \ - \ \xi_{\,1}\,\Vert^{\,\ell\,-\ 3}\ ) \,\cdot\,  O\, (\,\xi_{1_{\,|_{\,n}}}^3 \,) \ + \   O\,(\, \Vert \, y \ - \ \xi_{\,1}\,\Vert^{\,\ell\,-\ 4}\ ) \,\cdot\,  O\, (\,\xi_{1_{\,|_{\,n}}}^4 \,) \ + \ \cdot \cdot \cdot \ +  \ \xi_{1_{\,|_{\,n}}}^\ell \ \right] \  + \\[0.2in]
    & \ & \ \ \ \ \ \ + \ \, O\,(\,|\,y_{\,n}\,|^{\,\ell\,-\,1} \,)\,\cdot\, \left[\  O\, (\,\xi_{1_{\,|_{\,n}}}^2 \,) \ + \ O\,(\, \Vert \, y \ - \ \xi_{\,1}\,\Vert^2\,)\ \right]\ .
\end{eqnarray*}

\newpage

Now we are ready to discern the key contribution in (\,4.20\,).\\[0.1in]
(\,4.46\,)
\begin{eqnarray*}
&  \ &     \int_{B_{\,\xi_{\,1}} (\,{\rho_{\,\nu}}\,)}
\,[\,n\,(\,n\,-\,2) \,-\, (\,{ c}_n\!\cdot K\,)\,(\,y\,)  \ ] \,\cdot\, \left(\ {{\lambda^2_1 \ - \ \Vert\ y\,-\,\xi_{\,1}\,\Vert^2}\over {\lambda^2_1 \ + \ \Vert\ y\,-\,\xi_{\,1}\,\Vert^2}} \right)\,\cdot\,\left(\ {\lambda_{\,1}\over {\lambda^2_1 \ + \ \Vert\ y\,-\,\xi_{\,1}\,\Vert^2}} \right)^{\!\!n } dy \\[0.2in]
& = &  \int_{B_{\,\xi_{\,1}}\,(\,\rho_{\,\nu}\,)   }  C\,(\,{\bf p}_{\,y}\,) \,\cdot \Vert \,y  \ - \ {\bf p}_{\,y}\,\Vert^{\,\ell} \,\cdot\, \left(\ {{\lambda^2_1 \ - \ \Vert\ y\,-\,\xi_{\,1}\,\Vert^2}\over {\lambda^2_1 \ + \ \Vert\ y\,-\,\xi_{\,1}\,\Vert^2}} \right)\,\cdot\,\left(\ {\lambda_{\,1}\over {\lambda^2_1 \ + \ \Vert\ y\,-\,\xi_{\,1}\,\Vert^2}} \right)^{\!\!n } dy \ + \\[0.1in]
& \ &  \hspace*{5.5in} \ + \ \ {\bf E}_{\,(\,4.46\,)} \\[0.1in]
& = &   C\,(\,{\bf p}_{\,1}\,) \,\cdot   \int_{B_{o}\,(\,\rho_{\,\nu}\,)}\ |\ y_{\,n} \ - \ \xi_{1_{\,|_{\,n}}}\,|^{\,\ell}   \,\cdot\, \left(\ {{\lambda^2_1 \ - \ \Vert\ y\,-\,\xi_{\,1}\,\Vert^2}\over {\lambda^2_1 \ + \ \Vert\ y\,-\,\xi_{\,1}\,\Vert^2}} \right)\,\cdot\,\left(\ {\lambda_{\,1}\over {\lambda^2_1 \ + \ \Vert\ y\,-\,\xi_{\,1}\,\Vert^2}} \right)^{\!\!n } dy\ + \ \\[0.1in]
& \ &  \hspace*{5.5in} \ + \ \ {\bf E}'_{\,(\,4.46\,)} \\[0.1in]
& = & \lambda^{\,\ell}_1\cdot  \int_{B_o  \,\left(\,\lambda^{-1}\,\cdot\,{\rho_{\,\nu}}\ \right)} |\,Y_n\,|^{\,\,\ell}\,\cdot\,\left(\ {1\over {1 \ + \ \Vert\ Y\,\Vert^2}} \right)^{\!\!n }  \,\cdot\, \left(\ {{1 \ - \ \Vert\ Y\,\Vert^2}\over {1 \ + \ \Vert\ Y\,\Vert^2}} \right)\ d\,Y    \ + \ {\bf E}'_{\,(\,4.46\,)}  \ \ \ \ \\[0.1in]
& \ &  \ \ \ \ \ \ \ \ \ \ \ \ \ \ \ \ \ \ \ \ \hspace*{3.5in} [ \ Y \ = \ \lambda^{-1}_1\,(\,  y\,-\,\xi_{\,1}\,)   \ ]\\[0.1in]
& = & \lambda^{\,\ell}_1\cdot  \left[ \ -\,\int_{\R^n} |\,Y_n\,|^{\,\,\ell}\,\cdot\,\left(\ {1\over {1 \ + \ \Vert\ Y\,\Vert^2}} \right)^{\!\!n }  \,\cdot\, \left(\ {{\Vert\ Y\,\Vert^2 \ - \ 1   }\over {\Vert\ Y\,\Vert^2 \ + \ 1  }} \right)\ d\,Y  \ \right]  \ + \ {\bf E}''_{\,(\,4.46\,)} \\[0.1in]
& \ & \ \ \ \ \ \ \ \ \uparrow \ \ {\mbox{negative}}
\end{eqnarray*}
In the above\,,\, the errors $ {\bf E}_{\,(\,4.46\,)} \,,$ $\  {\bf E}'_{\,(\,4.46\,)} \,$  and  $\, {\bf E}''_{\,(\,4.46\,)} \,$ are estimated in {\bf \S\,A\,9} \,of the \,{\bf e}\,-\,Appendix\,.\,\\[0.1in]
\hspace*{0.5in}Next, we show that
$$
{\hat C}_{\,1\,,\ 1}\ := \ \left( {{\ n \ - \ 2\ }\over 2}\right) \cdot  \int_{\R^n} |\,Y_n\,|^{\,\,\ell}\,\cdot\,\left(\ {1\over {1 \ + \ \Vert\ Y\,\Vert^2}} \right)^{\!\!n }  \,\cdot\, \left(\ {{\Vert\ Y\,\Vert^2 \ - \ 1   }\over {\Vert\ Y\,\Vert^2 \ + \ 1  }} \right)\ dY\ > \ 0\ . \leqno (\,4.47\,)
$$
To this end, we consider the stereographic projection
\begin{eqnarray*}
(\,4.48\,)   \  \ \ \ \ \dot{\cal P} : S^n \setminus \{ \,{\bf N}\, \}\, &\to &  {\R}^n\,\\[0.2in]
   X & \mapsto & Y \ = \ \dot{\cal P}\,(X)\,,\ \ \ \
{\mbox{where}} \ \   \ \  Y_i \ = \ {{X_{\,i}}\over{\ 1 \,-\, X_{\,n \,+\, 1}\ }}\ ,
\ \ \ \ \ \ \ 1 \ \le\  i \ \le \ n\,.\ \ \ \ \   \ \ \ \ \ \ \ \ \ \   \ \ \ \ \
\end{eqnarray*}
Here $\,X \,=\,(\,X_{\,1}\,,\,\,\cdot\,\cdot \cdot, \ X_{\,n \,+ \,1}\ )\
\in \,S^n\, \setminus \{\, {\bf N}\, \}\,\subset\,\R^{n\,+\,1}\,,\,$ and $\,{\bf N} \,=\, (\,0, \ \cdot\,\cdot \cdot, \ 0\,, \ 1\,)\,$ is the north pole. Conversely,\\[0.1in]
(\,4.49\,)
\begin{eqnarray*} X_{\,i} \ = \  {{2\,Y_{\,i}}\over{\ 1 \,+\, R^{\,2}\ }}\,, \ \ \ \ \  1 \ \le\  i \ \le \ n\,,\ \ \ \ \
{\mbox{and}} \ \ \ \ \ \ \ \ \
X_{\,n \,+ \,1} \ = \  {{R^{\,2} \,-\,1 }\over{\ R^{\,2} \,+\, 1\ }}\,,  \ \ \ \ \ \ \
{\mbox{where}} \ \  R  \ = \ \Vert\ Y\Vert\ .\,
\end{eqnarray*}
It is known that $\ \dot{\cal P}\ $ is a conformal map between $\,(\,S^n \setminus \{\, {\bf N}\, \}\,, \ g_1\,)\,$ and $\,(\R^n,\ g_o\,)\,.\,$ The conformal factor is given by
$$
{g_1}_{\,|_X}  \  =\  \left[\  {{4}\over {(\,1 \,+\, R^{\,2}\,)^2}} \ \right]\cdot {g_o}_{\,|_{\,Y}} \ \ \mfor \ \ Y \  = \ \dot{\cal P}\,(X)\,.
$$
With this, we continue with\\[0.1in]
(\,4.50\,)
 \begin{eqnarray*}
  & \ &    \int_{\R^n}   |\ Y_n\,|^{\,\,\ell}\,\cdot\,   \left(\ {1\over {1\,+\, R^{\,2}}} \right)^{\!n}\,\cdot\, \left[\  {{R^{\,2} \ - \ 1}\over { \ R^{\,2}\ + \ 1\  }}\,\right] \ d Y\ \ \ \ \ \ \ \ \ ( \ \ell \ \le \ n\,-\,2\ )\\[0.2in]
 & = &  {1\over {2^n}}\,\cdot\,\int_{S^n } \left[\ {{ |\ X_{\,n}\,|^{\,\,\ell}}\over {(1\,-\,X_{\,n\,+\,1})^{\,\ell} }}\   \right] \,\,\cdot\,X_{\,n\,+\,1} \ dS_{g_1}\\[0.2in]
 & = & {1\over {2^n}}\,\cdot\,\left(\  \int_{S^n_+ } \left[\ {{ |\ X_{\,n}\,|^{\,\,\ell}}\over {(1\,-\,X_{\,n\,+\,1})^{\,\ell} }}\ \right] \,\,\cdot\,X_{\,n\,+\,1} \ d\,S_{\,g_1} \ + \   \int_{S^n_- } \left[\ {{ |\ X_{\,n}\,|^{\,\,\ell}}\over {(1\,-\,X_{\,n\,+\,1})^{\,\ell} }}\  \right] \,\,\cdot\,X_{\,n\,+\,1} \ d\,S_{\,g_1}\,\right)\\[0.2in] & = &  {1\over {2^n}} \cdot\left(  \ \int_{S^n_+ } \left[\  {{ |\ X_{\,n}\,|^{\,\,\ell}}\over {(1\,-\,X_{\,n\,+\,1})^{\,\ell} }}\  \right] \,\,\cdot\,X_{\,n\,+\,1} \ d\,S_{\,g_1} \ + \  \right.\\[0.1in]
  & \ & \hspace*{2in}\left. \ \ \ \ \ + \  \int_{S^n_+ } \left[\ {{ |\ -\,X_{\,n}\,|^{\,\,\ell}}\over {(1\,-\,[\ -\,X_{\,n\,+\,1}\,] \ )^{\,\ell} }}\  \right] \,\,\cdot\,[\ - X_{\,n\,+\,1}\,] \ d\,S_{\,g_1} \right)\\[0.1in]
 & \ & \hspace*{3.5in}(\ \uparrow \ \ {\mbox{reflection\ :}} \ \  X_{\,n\,+\,1} \ \to \ -\, X_{\,n\,+\,1}\,)\\[0.1in]
 & = & {1\over {2^n}} \cdot\left(  \ \int_{S^n_+ } \left[\, {{ |\ X_{\,n}\,|^{\,\,\ell}}\over {(1\,-\,X_{\,n\,+\,1})^{\,\ell} }} \right] \,\,\cdot\,X_{\,n\,+\,1} \ d\,S_{\,g_1} \ -  \   \int_{S^n_+ } \left[\, {{ |\ X_{\,n}\,|^{\,\,\ell}}\over {(1\,-\,[\ -\,X_{\,n\,+\,1}\,] \ )^{\,\ell} }} \right] \,\,\cdot\, X_{\,n\,+\,1}  \ d\,S_{\,g_1}\  \right) \ .\\
 \end{eqnarray*}
In the above, $\,S^n_+\,$ is the upper hemisphere\,,\, and $\,S^n_-\,$   the lower\,.\, Consider a fixed point
\begin{eqnarray*}
  (\,X_{\,1}\,,\,\,\cdot\,\cdot \cdot, \ X_{\,n}\,,\,\ \,X_{\,n\,+\,1}\,)\! &  \in &  \!\!S^n_+ \\[0.12in]
{\mbox{and \ \ its \ \ ``reflection\," \ :}} \  \ \  \ \ \  (\,X_{\,1}\,,\,\,\cdot\,\cdot \cdot, \ X_{\,n}\,,\,\  [\,-\,X_{\,n\,+\,1}\,]\,)& \! \in & \!\! S^n_-\ \,. \ \ \ \ \ \  \ \ \  \  \ \ \ \ \ \ \ \ \ \ \ \ \ \ \ \ \ \ \ \ \ \ \
\end{eqnarray*}
For $\, 0 \ <  \ X_{\,n\,+\,1} \ <  \ 1\,,\,$ we have
$$
\left[\, {{|\ X_{\,n}\,|^{\,\,\ell}}\over {(\,1\,-\,X_{\,n\,+\,1})^2 }} \ \right] \,\ge \ \left[\, {{|\ X_{\,n}\,|^{\,\,\ell}}\over {(\,1\,+\,X_{\,n\,+\,1} )^2 }}\ \right]  \ = \ \left[\, {{|\ X_{\,n}\,|^{\,\,\ell}}\over {(\,1\,-\,[\ -\,X_{\,n\,+\,1}\ ]\,)^2 }}\ \right] \ . \leqno (4.51)
$$
With (\,4.50\,)\,,\, (\,4.51\,) yields
 \begin{eqnarray*}
{\hat C}_{\,1\,,\ 1}  & > &  0  \ .
 \end{eqnarray*}

 \vspace*{-0.25in}

That is, we verify (\,4.47\,) [\,under the conditions in Theorem 1.33\,]\,.
\newpage

\hspace*{0.5in}Likewise, as for the terms in (\,4.44\,) and (\,4.45\,) marked by (\,**\,)\,,\, to obtain the other estimate, we may replace
 $$
 \ell \ \to \ (\,\ell \ - \ 2\,)\,.
 $$
 Denote
$$
{\hat C}_{\,1\,, \ 3} \ = \ \left(\, {{\ n\,-\,2\ }\over 2}\, \right){{\ell\cdot (\,\ell\,-\,1\,)}\over 2}\,\cdot\,\int_{\R^n} |\,Y_n\,|^{\,\ell\,-\,2}\,\cdot\,\left(\ {1\over {1 \ + \ \Vert\ Y\,\Vert^2}} \right)^{\!\!n }  \,\cdot\, \left(\ {{\Vert\ Y\,\Vert^2 \ - \ 1   }\over {\Vert\ Y\,\Vert^2 \ + \ 1  }} \right)\ dY\ .
$$
Similar consideration as in  (\,4.50\,)\, and \,(\,4.51\,)\, yields $\,{\hat C}_{\,1\,, \ 3} \ >  \ 0\,.\,$ For the terms in (\,4.44\,) and (\,4.45\,) marked by (\,*\,)\,,\, we observe the cancellation in {\bf \S\,A\,9\,.\,f}\, of the \,{\bf e}\,-\,Appendix\,.\,
 After estimating the remaining error terms (\,see {\bf \S\,A\,9\,.\,b\,-- e\,, \ g\,--\,j}\, of the \,{\bf e}\,-\,Appendix\,)\,,\, we arrive at the first and thrid terms in L.H.S. of (\,1.37\,)\,.\bk
 For  the spherical symmetric case in which the   bubbles are arranged ``\,evenly\,"  close to a great circle of $\,S^{\,n\,-\,1}$\,,\, we refer to Proposition 3.2 in the work by Wei and Yan  \cite{Wei-Yan}\,.

\vspace*{0.4in}


{\bf \large {\bf \S\,5.}} \ \ {\bf \large {\bf   Extracting the key information in the reduced functional\ -}}\\[0.1in]
\hspace*{0.34in}{\bf \large \ \  {\bf  derivative with respect to  $\,\xi$\,.}}\\[0.2in]
Via (\,4.4\,)\,,\, as in (\,4.6\,)\,,\, we obtain\\[0.1in]
(\,5.1\,)
 \begin{eqnarray*}
& \ & \left(\,\lambda_{\,1}\cdot {\partial\over {\partial\, \xi_{1_{\,|_{\,j}}}   }} \right) {\bf I}_{\,\cal R}\ \bigg\vert_{ \, {\bf P} _{\,(\,\flat\,)} } \ \ \ \ \ \ \ \ \ \ \ \ \ \ \ \ \ \ \ \ \left[ \ {\bf P}_{\,(\,\flat)\,} \ = \ (\,\lambda_{\,1}\,, \ \cdot \cdot \cdot\,, \ \lambda_{\,\,\flat}\,; \ \xi_{\,1}\,, \ \cdot \cdot \cdot\,, \ \xi_{\,\flat}\,) \ \ \ \right] \\[0.2in]
& = &  \!-\int_{\R^n} (\,\Delta\,W_{\,\flat\,})\cdot\left[  \, \left(\,\lambda_{\,1}\cdot {\partial\over {\partial\, \xi_{1_{\,|_{\,j}}}}} \right)W_{\,\flat\,}  \,   \right]  \, -\!   \int_{\R^n}\,(\,{\tilde c}_n\!\cdot K\,)\,(\,W_{\,\flat\,}\, )^{{n\,+\,2}\over {\,n\,-\ 2\,}}\!\cdot \left[   \left(\lambda_{\,1}\cdot {\partial\over {\partial\, \xi_{1_{\,|_{\,j}}}}} \right)W_{\,\flat\,}    \right] + \, {\cal E}_{1_{\,|_{\,j}}} \\[0.2in]
  & = & \left\{ \ -\int_{\R^n} (\,\Delta\,W_{\,\flat\,})\,\cdot\,\left[ \,  \left(\,\lambda_{\,1}\cdot {\partial\over {\partial\, \xi_{1_{\,|_{\,j}}}}} \right)W_{\,\flat\,}   \   \right]  \right.\\[0.1in]
     & \ & \ \ \ \ \ \ \ \ \ \ \ \ \ \ \ \ \  \ \ \ \ \ \ \ \ \ \ \ \ \ \ \ \ \  \left. \ -  \    \int_{\R^n}\, \,n\,(\,n\,-\,2)\cdot\,(\,W_{\,\flat\,}\, )^{{n\,+\,2}\over {\,n\,-\ 2\,}}\,\cdot\, \left[ \,  \left(\,\lambda_{\,1}\cdot {\partial\over {\partial\, \xi_{1_{\,|_{\,j}}}}} \right)W_{\,\flat\,}   \   \right]  \ \right\} \ \\[0.1in]
  & \ & \hspace*{4.5in} \ \ \ \ \ \ \  \uparrow \ \ \cdot\,\cdot\,\cdot\,\,\cdot\, {\bf (\,C\,)}_{\,(\,5.1\,)}\\[0.1in]
 & \ &  \ \ \ \ \ \ \ \ \ \ \   - \   \left\{ \ \int_{\R^n}\, [\ (\,{\tilde c}_n\!\cdot K\,) \ - \ n\,(\,n\,-\,2)\,]  \cdot\,(\,W_{\,\flat\,}\, )^{{n\,+\,2}\over {\,n\,-\ 2\,}}\,\cdot\, \left[ \,  \left(\,\lambda_{\,1}\cdot {\partial\over {\partial\, \xi_{1_{\,|_{\,j}}}}} \right)W_{\,\flat\,}  \    \right]   \ \right\}  \\[0.1in]
  & \ & \hspace*{4.5in} \ \ \ \ \ \ \  \uparrow \ \ \cdot\,\cdot\,\cdot\,\,\cdot\, {\bf (\,D\,)}_{\,(\,5.1\,)}\\[0.1in]
 & \ & \ \ \ \ \ \ \    \ \ \ \ \ \ \   \ \ \ \ \ \ \  \ \ \    \ + \ \  {\cal E}_{1_{\,|_{\,j}}} \ \ (\ \leftarrow \ \ {\mbox{estimated \ \ in \ \ {\bf \S\,A\,10} \ \ of \ \ the \ \ {\bf e}\,-\,Appendix}}\ )\ .
 \end{eqnarray*}

 \newpage

Recall that \\[0.1in]
(\,5.2\,)
  $$\left(\,\lambda_{\,1}\,\cdot\,{\partial\over {\partial\, \xi_{1_{\,|_{\,j}}}}} \right)\left[\left(\,{\lambda_{\,1}\over {\lambda^2_1\ + \ \Vert\  y \ -\  \xi_{\,1}\,\Vert^{\,2}}} \,\right)^{{n \,-\, 2}\over 2}\  \right] \, = \, \left( \ n\,-\,2 \ \right)\cdot \left( {{ \,\lambda_{\,1} }\over {\lambda^2_1 + \Vert\  y \,-\, \xi_{\,1}\,\Vert^{\,2} }}\right)^{\!\!{n\over 2}}\,(\,y_j\,-\,\xi_{1_{\,|_{\,j}}}\,) \ .
$$

\vspace*{0.2in}


{\bf \large {\bf \S\,5 a\,.}} \  \ {\it Extracting the key information in}\, ${\bf (\,C\,)}_{(\,5.1\,)}$\,.\\[0.1in]
The presence of $\,(\,y_j\,-\,\xi_{1_{\,|_{\,j}}}\,)\,$  in (\,5.2\,) enables us to proceed with similar argument as in {\bf \S\,4\,b}\, (\,in some sense ``\,orthogonal\," to it\,)\,,\,  making use of the symmetric cancellation\,:
$$
\int_{B_{\,\xi_{\,1}}(\,R\,)} f\,(\,r\,)\,\cdot\,(\,y_1\,-\,\xi_{1_{\,|_1}}\,)^{\,\alpha_1}\ \,\cdot\,\cdot\,\cdot\,\ (\,y_{\,n}\,-\,\xi_{1_{\,|_n}}\,)^{\,\alpha_n}\ d\,y\ = \ 0 \ , \ \ \ \ \ \   \  r \ = \ \Vert\  y \ - \ \xi_{\,1}\,\Vert\ \ , \leqno (\,5.3\,)
$$
when one (\,or more\,) of the indices  is\, {\it odd} \, (\,and positive\,;\, $\alpha_1\,,\ \cdot \cdot \cdot\,, \ \alpha_{\,n\,-\,1}\,$ and  $\alpha_{\,n}$ \,are non\,-\,negative integers\,)\,.\, Here $\,f\,$ is a continuous function defined on $\,[\,0\,, \ \infty\,)\,.\,$
With (5.2) and (5.3)\,,\, one sees in (4.11)  that the first and the second terms [\ in (4.11)\ ]  are each cancelled after the integration, leaving only the third term. We leave the details to {\bf \S\,A\,10} in the \,{\bf e}\,-\,Appendix\,.\,    Expressing ${\bf (\,C\,)}_{(\,5.1\,)}$ in gradient form, we arrive at the first term in the L.H.S. of (\,1.38\,) [\,under the conditions in Theorem 1.33\,]\,. \, For the full details, see {\bf \S\,A\,10}\, in the \,{\bf e}\,-\,Appendix\,.

\vspace*{0.25in}

{\bf \large {\bf \S\,5 b\,.}} \  \  {\it Extracting the key information in}\, ${\bf (\,D\,)}_{(\,5.1\,)}$\,.\\[0.1in]
 Next, we consider the term in \,${\bf (\,D\,)}_{(\,5.1\,)}$\,,\, expressed in vector form\,:\\[0.1in]
 (\,5.4\,)
 \begin{eqnarray*}
 {\bf{[\,\vec{\,D}\,]}} & = &  -\,\left\{ \ \int_{\R^n}\, [\,(\,{\tilde c}_n\!\cdot K\,) \ - \ n\,(\,n\,-\,2)\,]  \cdot\,(\,W_{\,\flat\,}\, )^{{n\,+\,2}\over {\,n\,-\ 2\,}}\,\cdot\,\left[\ (\,\lambda_{\,1}\cdot \btd_{\xi_{\,1}}\,) \, W_{\,\flat\,}\ \right]  \ \right\}\ d\,y\\[0.2in]
 & = &  \left\{ \ \int_{\R^n}\, [\ n\,(\,n\,-\,2\,) \ - \ (\,{\tilde c}_n\!\cdot K\,) \,]  \cdot\,(\,W_{\,\flat\,}\, )^{{n\,+\,2}\over {\,n\,-\ 2\,}}\,\cdot\,\left[\ (\,\lambda_{\,1}\cdot \btd_{\xi_{\,1}}\,) \, V_1\ \right]  \
 \right\} \ d\,y\\[0.2in]
 &  \ & \hspace*{-0.5in}= \ \left\{ \,(\,n\,-\,2\,)\,\cdot
 \int_{B_o\,(\,\xi_{\,1}\,)}  [\ n\,(\,n\,-\,2\,) \ - \ (\,{\tilde c}_n\!\cdot K\,) \,]  \cdot  \,\left(  {{\lambda_{\,1}}\over {\lambda_{\,1}^2 \ + \ \Vert\, y \,-\, \xi_{\,1}\Vert^{\,2}}} \right)^{\!\!n\,+\,1}  \cdot   (\,y \ - \ \xi_{\,1} \,) \, \right\} * \\[0.2in]
  & \ & \ *\, \left[\ 1 \ + \ {\bf E}_{\,(\,5.4\,)}\  \right]   \ + \  \tilde{\bf E}_{\,(\,5.4\,)} \ .
\end{eqnarray*}
Refer to (\,A.11.1) of the {\bf e}\,-\,Appendix for $\,{\bf E}_{\,(\,5.4\,)}\,$ and $\,\tilde{\bf E}_{\,(\,5.4\,)}\,$.\,
Via the Lagrange multiplier method as in {\bf \S\,4\,c}\,,\, and  (\,5.3\,)\,,\, we have\\[0.1in]
(\,5.5\,)
$$
  \int_{B_{\,\xi_{\,1}}\,(\,{\rho_{\,\nu}})}  |\,y_{\,n}\,-\,\xi_{1_{|_n}}\,|^{\,\,\ell}   \cdot \, \left\{  \,  \left(  {{\lambda_{\,1}}\over {\lambda_{\,1}^2 \ + \ \Vert|\, y \,-\, \xi_{\,1}\Vert|^{\,\,2}}} \right)^{\!\!n\,+\,1} \!\!\!\! \cdot  (\,y_j \ - \ \xi_{1_j} \ )  \,\right\}\ dy  \ = \ 0 \ \ \ \ \ (\ 1 \ \le \ j \ \le \ n\ )\,.
$$
Consider  the next order terms in \,(4.44)\, and \,(4.45)\,,\, expressed by \\[0.1in]
(5.6)
\begin{eqnarray*}
{\bf M} \ := \ (\,n\,-\,2\,) \cdot  \ell \cdot   C \,({\bf p}_{\,1}) \cdot  \int_{B_{\,\xi_{\,1}}(\,{\rho_{\,\nu}}\,)}  |\,y_{\,n}\,-\,\xi_{1_{|_n}}\,|^{\,\ell\,-\,1}   \cdot \xi_{1_{|_n}}  \cdot \, \left\{    \,\left(  {{\lambda_{\,1}}\over {\lambda_{\,1}^2 \ + \ |\, y \,-\, \xi_{\,1}\,|^2}} \right)^{\!\!n\,+\,1} \!\!\!\!\cdot   |\,y_{\,n} \ - \ \xi_{1_n} \, |  \,\right\} \ .
\end{eqnarray*}
Under the normalization (\,4.23\,) and (\,4,24\,)\,,\,
$\,{\bf M} \,$ can be calculated using similar technique as in \,{\bf \S\,4\,c}\,,\, leading to
$$
{\bf M} \ = \   {\hat C}_{\,2\,,\,3} \cdot  C\, (\,{\bf p}_1\,) \  \cdot {{ \lambda_{\,1}^\ell }\over {\lambda_{\,1} }} \cdot \left[ \ 1 \ + \ {\bf E}_{\,(\,5.7\,)}\ \right] \cdot \xi_{1\,|_{\,n}}\ . \leqno (5.7)
$$
See  {\bf \S\,A\,11} of the {\bf e}\,-\,Appendix for the details\,,\, and {\bf \S\,1\,d} for the positive constant $\, {\hat C}_{\,2\,,\,3}\,.\,$ Moving back to the general position\,,\, as well as general index $\,l\,=\,1\,, \ \cdot \cdot \cdot\,, \ n$\,,\,  we are led to
$$ {\bf{[\,{\vec{\,D}}\,]}} \ = \  \left[\  {\hat C}_{\,2\,,\,3}\cdot   C\, (\,{\bf p}_{\,l}) \cdot \lambda_{\ \!l}^\ell \cdot\left( {{  \eta_{\,l}  }\over {  \lambda_{\ \!l} }} \right) \ {\bf n}_{\,l} \  \right] \,\cdot\, \left[\ 1 \ + \  {\bf E}_{\,(\,5.7\,)}\  \right]    \ \    + \ \ \vec{\,{\bf E}}_{\,(\,5.8\,)}\  ,\leqno (\,5.8\,)
$$
where
$$
\eta_{\,l} \ = \  \Vert \  \xi_{\,l} \,-\, {\bf p}_{\,l} \,\Vert \ \ \  \ {\mbox{and}} \ \ \ \ \left( \ {\bf n}_{\,l} \ = \  {{ \xi_{\,l} \,-\, {\bf p}_{\,l} }\over { \Vert \, \xi_{\,l} \,-\, {\bf p}_{\,l} \,\Vert }} \ \right)\,. \leqno (\,5.9\,)
$$
Thus we obtain the second term in the L.H.S. of (\,1.38\,)\, [\,under the conditions in Theorem 1.33\,]\,.\, For the full detail, we refer the readers to  {\bf \S\,A\,11} of the {\bf e}\,-\,Appendix\,,\, in which one can also find the estimates of the  errors $\, {\bf E}_{\,(\,5.7\,)}\,$ and   $\,\vec{\,{\bf E}}_{\,(\,5.8\,)}\,$.


\vspace*{0.5in}

{\bf \large {\bf \S\,6. \ \  Extracting the key information in the reduced functional\,.}}\\[0.2in]
After the derivatives\,,\, we come to the original reduced functional\,.\, We pick it up from (\,4.2\,)\,:
$$ {\bf I}_{\,{\cal R}} \,(\,\lambda_{\,1}\,, \ \cdot \cdot \cdot\,, \ \lambda_{\,\,\flat}\,; \ \xi_{\,1}\,, \ \cdot \cdot \cdot\,, \ \xi_{\,\flat}\,)\ = \ {\bf I}\,(\ W_{\,\flat}\ + \ \phi\ ) \ = \  {\bf I}\,(\,W_{\,\flat\,}\,)\ + \ {\bf E}_{\,(\,6.1\,)}\ . \leqno (\,6.1\,)
$$
Here

\begin{eqnarray*}
(\,6.2\,) \ \ \ \  \ \ \ \ {\bf I}\,(\,W_{\,\flat\,}\,) & = &   {1\over 2}\,\int_{\R^n}\,\langle\,\btd\,W_{\,\flat\,}\,,\,\btd\,W_{\,\flat\,}\,\rangle\ -\ \left(\,{{n\,-\,2}\over {2n}}\,\right)\,\cdot\,\int_{\R^n}\,(\,{\tilde c}_n\!\cdot K\,)\,W_{\,\flat\,}^{{2n}\over {\,n\,-\ 2\,}}\ \ \ \ \ \ \ \ \ \ \  \ \ \ \ \ \ \ \  \\[0.1in]
 & \ & \ \ \ \ \ \ \ \uparrow \ \ {\bf{(\,F\,)}}_{\,(\,6.2\,)} \hspace*{1.21in} \ \ \ \ \ \ \ \ \ \ \ \uparrow \ \ {\bf{(\,G\,)}}_{\,(\,6.2\,)}
\end{eqnarray*}
Refer to the {\bf \S\,A\,12}\, of the  \,{\bf e}\,-\,Appendix for the error team $\,{\bf E}_{\,(\,6.1\,)}\,$ [\ see (\,A.12.1\,)\,]\,.


\newpage

{\bf \large {\bf \S\,6 a\,.}} \  \ {\it Expansion of}\,   ${\bf{(\,F\,)}}_{\,(\,6.2\,)}$\ . \ Consider the first term in the L.H.S. of  (\,6.2\,)\,.\, Via integration by parts formula and equation (\,1.6\,)\,,\, we proceed with
\begin{eqnarray*}
(6.3) \ \ \ \ \    \ \ \ \ \ \  \
  & \ & {1\over 2}\cdot \int_{\R^n}\,\langle\,\btd\,W_{\,\flat\,}\,,\,\btd\,W_{\,\flat\,}\,\rangle\ \\[0.2in]
  &  = &   -\,   {1\over 2}\cdot \int_{\R^n}\,(\,\Delta\,W_{\,\flat\,}\,)\,\cdot\,W_{\,\flat\,} \ \\[0.2in]
  & = &  {1\over 2}\cdot n\,(\,n\,-\,2\,)\cdot \int_{\R^n}\,\left(\,V_1^{{n\,+\,2}\over {\,n\,-\ 2\,}} \ + \ \cdot\,\,\cdot\,\cdot \ + \ V_\flat^{{n\,+\,2}\over {\,n\,-\ 2\,}}\,\right)\,\cdot\,(\,V_1 \ + \ \cdot\,\cdot\,\cdot\,\ + \ V_\flat \,)\ \\[0.2in]
  & = &  {1\over 2}\cdot n\,(\,n\,-\,2\,)\cdot \int_{\R^n}\,\left(\,V_1^{{2n}\over {\,n\,-\ 2\,}} \ + \ \cdot\,\cdot\,\cdot\,\ + \ V_\flat^{{2n}\over {\,n\,-\ 2\,}}\,\right)  \
  \\[0.1in]
    & \ & \hspace*{-1in}\left[ \   {\mbox{to \ \ be \ \ combined \ \ with  \ \ similar  \ \ term \ \ in \ \ }}  (\,6.4\,)\ , \ \  \uparrow  \ = \   {1\over 2}\,\cdot\,n\,(\,n\,-\,2\,)\cdot V\,(\,n)\,\cdot\,\flat\ ; \right]\\[0.2in]
  & \ & \ \ + \  {1\over 2}\cdot n\,(\,n\,-\,2\,)\cdot \int_{\R^n}\,V_1\,\cdot\,\left(\,V_2^{{n\,+\,2}\over {\,n\,-\ 2\,}} \ + \ \cdot\,\cdot\,\cdot\,\ + \ V_\flat^{{n\,+\,2}\over {\,n\,-\ 2\,}}\,\right)  \ \ + \ \\[0.1in]
  & \ &  \ \ \ \ \hspace*{1in}: \\[0.01in]
   & \ &  \ \ \ \ \hspace*{1in}:  \\[0.1in]
  & \ & \ \   \ \ \ \ \ \ \ + \  {1\over 2}\cdot n\,(\,n\,-\,2\,)\cdot \int_{\R^n}\,V_\flat\,\cdot\, \left(\,V_1^{{n\,+\,2}\over {\,n\,-\ 2\,}} \ + \ \cdot\,\cdot\,\cdot\,\ + \ V_{\flat\,-\,1}^{{n\,+\,2}\over {\,n\,-\ 2\,}}\,\right) \ . \emph{}
\end{eqnarray*}
Refer to {\bf {\S\,1\,d}} for the constant $\,V\,(\,n\,)\,.$\,
The interaction teams such as $\,``\ V_1 \cdot V_2^{{n\,+\,2}\over {\,n\,-\ 2\,}} \ "\,$ can be estimated as in  {\bf \S\,4\,b}\,,\,  yielding  the second term in the L.H.S. of (\,1.36\,) [\,under the conditions in Theorem 1.33\,]\,.\,
Refer to the {\bf \S\,A\,12\,.\,a\,--\,d\,} of the  \,{\bf e}\,-\,Appendix\,.

\vspace*{0.25in}

{\bf \large {\bf \S\,6 b\,.}} \  \ {\it Expansion of}\, ${\bf{(\,G\,)}}_{\,(\,6.2\,)}$\ . \ \ We have\\[0.1in]
(\,6.4\,)
\begin{eqnarray*}
& \ & \\[-0.4in]
  & \ & -\ \left(\ {{n\,-\,2}\over {2n}}\,\right)\,\cdot\,\int_{\R^n}\,(\,{\tilde c}_n\!\cdot K\,)\,W_{\,\flat\,}^{{2n}\over {\,n\,-\ 2\,}} \\[0.2in]
  & = & -\ \left(\ {{n\,-\,2}\over {2n}}\,\right)\,\cdot\,\Bigg\{ \ \int_{\R^n}\,[\ (\,{\tilde c}_n\!\cdot K\,) \ - \ n\,(\,n\,-\,2\,) \ ]\,\cdot\,W_{\,\flat\,}^{{2n}\over {\,n\,-\ 2\,}} \  \   \\[0.2in]& \ &  \hspace*{1.5in}    + \ \int_{\R^n}\,  n\,(\,n\,-\,2\,) \,\cdot\,\left[ \ W_{\,\flat\,}^{{2n}\over {\,n\,-\ 2\,}} \ - \ \left(\,V_1^{{2n}\over {\,n\,-\ 2\,}} \ + \ \cdot\,\cdot\,\cdot\,\ + \ V_\flat^{{2n}\over {\,n\,-\ 2\,}}\,\right)  \ \right] \  \\[0.2in]
& \ &   \hspace*{1.71in}  +   \ \int_{\R^n}\,  n\,(\,n\,-\,2\,) \,\cdot\,\left(\,V_1^{{2n}\over {\,n\,-\ 2\,}} \ + \ \cdot\,\cdot\,\cdot\,\ + \ V_\flat^{{2n}\over {\,n\,-\ 2\,}}\,\right) \  \  \Bigg\} \\[0.2in]
  & = & -\ \left(\,{{n\,-\,2}\over {2n}}\,\right)\,\cdot\,\Bigg\{ \ \int_{\R^n}\,[\ (\,{\tilde c}_n\!\cdot K\,) \ - \ n\,(\,n\,-\,2\,) \ ]\,\cdot\,\left(\,V_1^{{2n}\over {\,n\,-\ 2\,}} \ + \ \cdot\,\cdot\,\cdot\,\ + \ V_\flat^{{2n}\over {\,n\,-\ 2\,}}\,\right)\\[0.1in]
    & \ & \hspace*{5in} \uparrow  \  \ \cdot \cdot \cdot \ (\,{\bf H}\,)_{\,(\,6.4\,)} \\[0.1in]
  & \ &  \ \ \ \ \ \ \ \ \ \  + \     \int_{\R^n}\,  [\ (\,{\tilde c}_n\!\cdot K\,) \ - \ n\,(\,n\,-\,2\,) \ ]\,\cdot\,\left[ \ W_{\,\flat\,}^{{2n}\over {\,n\,-\ 2\,}} \ - \ \left(\,V_1^{{2n}\over {\,n\,-\ 2\,}} \ + \ \cdot\,\cdot\,\cdot\,\ + \ V_\flat^{{2n}\over {\,n\,-\ 2\,}}\,\right)  \ \right] \   \ + \\[0.2in]& \ & \ \ \ \ \ \ \ \ \ \ \ \ \ \    + \ \int_{\R^n}\,  n\,(\,n\,-\,2\,) \,\cdot\,\left[ \ W_{\,\flat\,}^{{2n}\over {\,n\,-\ 2\,}} \ - \ \left(\,V_1^{{2n}\over {\,n\,-\ 2\,}} \ + \ \cdot\,\cdot\,\cdot\,\ + \ V_\flat^{{2n}\over {\,n\,-\ 2\,}}\,\right)  \ \right] \  \\[0.2in]
    & \ &  \ \ \ \  \ \ \ \ \ \ \ \ \ \ \ \ \ \ \ \ \ \  + \ \int_{\R^n}\,  n\,(\,n\,-\,2\,) \,\cdot\,\left[ \ W_{\,\flat\,}^{{2n}\over {\,n\,-\ 2\,}} \ - \ \left(\,V_1^{{2n}\over {\,n\,-\ 2\,}} \ + \ \cdot\,\cdot\,\cdot\,\ + \ V_\flat^{{2n}\over {\,n\,-\ 2\,}}\,\right)  \ \right] \  \ + \\[0.2in]
& \ &  \ \ \ \ \ \ \ \   \ \ \ \ \ \ \ \   \ \ \ \ \ \ \ \   \ \ \ \ \ \ \ \  +   \ \int_{\R^n}\,  n\,(\,n\,-\,2\,) \,\cdot\,\left(\,V_1^{{2n}\over {\,n\,-\ 2\,}} \ + \ \cdot\,\cdot\,\cdot\,\ + \ V_\flat^{{2n}\over {\,n\,-\ 2\,}}\,\right) \  \  \Bigg\}  \\[0.1in]
    & \ & \left[ \   {\mbox{to \ \ be \ \ combined \ \ with  \ \ similar  \ \ term \ \ in \ \ }}  (\,6.3\,)\ , \ \  \uparrow  \ = \   {1\over 2}\,\cdot\,n\,(\,n\,-\,2\,)\cdot V\,(\,n)\,\cdot\,\flat\ \right] \ .
\end{eqnarray*}
For the purpose on finding a critical point to the reduced functional\,,\, the key contribution lies in $\,{\bf H}_{\,(\,6.4\,)}$\,,\, which can be estimated as in \,{\bf \S\,4\,c}\,,\, leading to the first term in the L.H.S. of (\,1.36\,) [\,under the conditions in Theorem 1.33\,]\,.\, The details can be found in  {\bf \S\,A\,12\,.\,e\,-\ i\,} of the  \,{\bf e}\,-\,Appendix\,.

\vspace*{0.2in}

{\bf \large {\bf \S\,6 c\,.}} \  \  {\it The constants in } {\bf \S\,1\,d}\,. \ \
 With regards to the constants and their relations shown in {\bf \S\,1\,d}\,,\,we place it in {\bf \S\,A\,13\,} of the {\bf e}\,-\,Appendix the formal justification\,,\, as well as an  intuitive link \{\,via differentiation of the main teams in right hand side of (\,1.36\,) [\,that is\,,\, the reduced functional\,]\,\}\,.

\newpage


\centerline{\bf \LARGE {\bf e\,-Appendix for the Paper  }}

\vspace*{0.1in}

\begin{center}
{\LARGE {\bf ``\,Construction of Blow\,-\,up Sequences for  the}} \medskip \medskip \smallskip  \\
{\LARGE {\bf Prescribed Scalar Curvature Equation }} \medskip \medskip \smallskip  \\
{\LARGE {\bf on $S^n$. IV.  Clustered Blow\,-\,ups\,"  }}
\end{center}

\vspace{0.3in}

\centerline{\Large  {Man Chun  {\LARGE L}EUNG${\,}^{\sharp}$  }}

\vspace*{0.2in}

\centerline{\large {National University of Singapore  }}

\vspace{0.53in}

{\it Here we follow the notations, conventions, equation numbers, section numbers, lemma, proposition and theorem numbers as  used in the main text\,},\,  {\it unless otherwise is specifically mentioned} (\,{\it  for instances, those equation numbers started with}\, `A'\,)\,.\,  The reference citation numbers are  referred to the end of this {\bf e}\,-\,Appendix\,.\\[0.4in]
%
%
{\bf \large \S\,A\,1.} \ \ {\bf \large The case of one bubble.}\\[0.2in]
{\large{\bf \S\,A\,1\,.\,a\,.\,}} \   {\bf The set\,-\,up\,.\,} \ \  We begin with the equation
$$
\Delta_1 \, U \ - \ [\ {\tilde c}_n \,n\, (\,n\, - \,1\,) \ ] \,U \ + \ (\,{\tilde c}_n \,{\cal K}\,)\, U^{{n \,+\, 2}\over {\,n \,-\, 2\,}} \ = \ 0 \ \ \ \ \ \ {\mbox{in}} \ \ \ S^n  \ \ \ \ \ \ \ \ \ (\,U \ >  \ 0\,)\,.\leqno (A.1.1)
$$
Here $\,{\tilde c}_n\, = {{(\,n\,-\,2}\over { \,4\,(\,n \,-\, 1\,)\, }}\,.\,$
Via the stereographic projection
\begin{eqnarray*}
(A.1.2) \ \ \  \ \ \ \  \ \ \ \   \ \ \ \   \dot{\cal P} : S^n \setminus \{\, {\bf N}\, \}\, &\to &  {\R}^n\,\\
x & \mapsto & y \ = \ \dot{\cal P}\,(x)\,,\ \ \ \ \\[0.1in]
& \ &  \hspace*{0.7in}
{\mbox{where}} \ \    y_i \,=\, {{x_i}\over{\ 1 \,- \ x_{n + 1}\ }}\,,
\ \ \ \ 1 \,\le\, i \,\le \,n\,.
\end{eqnarray*}
Here $\,x\, =\,(\,x_1\,,\, \cdot \cdot \cdot, \ x_{n + 1}\,)
\in S^n \setminus \{\, {\bf N}\, \}\,\subset \R^{n + 1},\,\,$ and $\,{\bf N} \ = \ (\,0, \ \cdot \cdot \cdot, \ 0, \ 1\,)\,.\,$
Conversely,\,
$$\ \ \ \displaystyle{ x_i \,=\, {{2\,y_i}\over{\ 1 \,+\, r^2\ }}\,, \ \ \ \ \  1 \,\le\, i \,\le\, n\,,\ \ \ \ \ {\mbox{and}} \ \ \ \ \
x_{n + 1} \,=\ {{\ r^2 \,-\,1 \ }\over{\ r^2 \,+\, 1\ }}\ ,} \ \ \ \
{\mbox{where}} \ \  r\,=\, \Vert \,y\,\Vert\,.\, \leqno (A.1.3)
$$

\vspace*{0.35in}

$\overline{\ \ \ \ \ \ \ \ \ \ \ \ \ \ \ \ \ \ \ \ \  \ \ \ \ \ \ \ \ \ \ \ \ \ \ \ \ \ \ \ \ \ \ \ }$\\[0.03in]
\noindent{$\!\!{\!}^{\sharp}$ \,{\tt{matlmc@nus.edu.sg}}

 \newpage

It is known that $\,\dot{\cal P}\,$ is a conformal map between $\,(\,S^n \setminus \{ \,{\bf N} \,\}\,, \ g_1\,)\,$ and $\,(\,\R^n,\ g_o\,)\,.\,$ The conformal factor is given by
$$
g_1 \,(x) = \left[\ {{4}\over {(\,1 \,+ \,\Vert\,y\,\Vert^2\,)^{\,2}}}\  \right]\,g_o \,(\,y\,) \ \ \ \ \mfor \ \  y \ = \ \dot{\cal P}\,(\,x\,)\,.\leqno  (A.1.4)
$$
Note that the conformal factor is related to the standard bubble [\ cf. (\,A.1.7\,)\ ]\,.\,
Set
$$
v \, (\,y\,) \ :=  \  U  \, (\,{\dot{\cal P}}^{-1} (\,y\,)\,)\, \cdot \left( {2\over {1\,+\, \Vert \,y\,\Vert^{\,2}}} \right)^{{n \,-\, 2}\over 2}  \ \ \ {\mbox{and}} \ \ \ \  K\,(\,y\,)\  := \ {\cal K}  \,(\,{\dot{\cal P}}^{-1} (\,y\,)\,)   \mfor y \,\in\, \R^n. \ \ \ \ \
$$
Suppose that $\,U$\, is a positive $\,C^2\,$-\,solution of equation (A.1.1). Then  $\,v\,$ satisfies the equation
$$
\Delta_o  \,v \ + \ (\,{\tilde c}_n \,K\,)\,v^{{n \, +\,2}\over {n \,-\, 2}}  \ \,= \ 0 \ \ \ \ {\mbox{in}} \ \ \R^n\,. \leqno (A.1.5)
$$
Note that in this setting
$$
\Delta_o  \,V_{\lambda\,,\ \xi} \ +\  n\,(\,n\,-\,2\,)\,V_{\lambda\,,\ \xi}^{{n \,+\, 2}\over {n\, -\, 2}}  \ \,= \ 0 \ \ \ \ {\mbox{in}} \ \ \R^n\,, \leqno (A.1.6)
$$
where
$$
V_{\lambda, \ \emph{}\xi} \,(\,y\,)\ = \ \left( {\lambda\over {\lambda^2 \ +\  \Vert\,y \ -\  \xi\Vert^{\,2}}} \right)^{\!\!{{n-2}\over 2} } \ \ \ \ \ \mbox{with} \ \ \ \ (\,\lambda\,,\, \ \xi\,)\,\in\, \R^+\times \R^n\ .   \leqno (A.1.7)
 $$
Note that
$$
V_{\lambda, \ \emph{}\xi} \,(\,y\,) \ \le \ {C \over {\ \Vert\,y\,\Vert^{n\,-\,2}\, }} \ \ \ \ \mfor \ \  \Vert\,y\,\Vert \ \ge \ R \ \gg \ 1\,. \leqno (\,A.1.8\,)
$$

\vspace*{0.2in}

{\large{\bf \S\,A\,1.\,b\,.\,}} \   {\bf Solving the linear inhomogeneous equation \,--\, the case of one bubble\,.\,}  \\[0.1in]
{\it  Introducing orthogonality.} \ \
Associated with equation (\,A.1.5\,) is the linear inhomogeneous equation
$$
\Delta\,\phi \ +\ n\,(n\,+\,2) \cdot \left( V_{\,\lambda \,,\  \xi}^{4\over {n\,-\,2}}\cdot \phi\right) \ = \ \left( V_{\,\lambda \,,\  \xi}^{4\over {n\,-\,2}}\cdot h\right)\,, \ \ \ \ \ \ {\mbox{where}} \ \ \phi\,, \  \
 h \ \in \ {\cal D}^{\,1\,,\,2}\,.\leqno (A.1.9)
$$
We first observe that\,,\, by multiplying both sides of (A.1.9) with
$$\,\left( \lambda \cdot {{\partial }\over {\partial \lambda}} \right) V_{\,\lambda \,,\  \xi}\ , $$
and upon integration, we have\\[0.1in]
(\,A.1.10\,)
\begin{eqnarray*}
0 & = & \int_{\R^n} (\,\Delta\,\phi\,) \cdot [\ (\lambda\cdot\partial_\lambda\,)\, V_{\,\lambda \,,\  \xi} \ ] \ + \  n\,(n\,+\,2) \cdot \int_{\R^n} \left( V_{\,\lambda \,,\  \xi}^{4\over {n\,-\,2}}\cdot \phi\right) \cdot [\ (\lambda\cdot\partial_\lambda\,)\, V_{\,\lambda \,,\  \xi} \ ]  \\[0.2in]
& = &  \int_{\R^n} \left( V_{\,\lambda \,,\  \xi}^{4\over {n\,-\,2}}\cdot h\right) \cdot [\ (\lambda\cdot\partial_\lambda\,)\, V_{\,\lambda \,,\  \xi} \ ] \ = \ \langle\ h\,, \ [\ (\lambda\cdot\partial_\lambda\,)\, V_{\,\lambda \,,\  \xi} \ ]  \ \rangle_{\,\btd}
\end{eqnarray*}
Here we apply
$$
\Delta\,V_{\,\lambda \,,\  \xi} \ + \ n\,(\,n\,-\,2\,)\,V_{\,\lambda \,,\  \xi}^{{\,n\,+\,2\,}\over {n\,-\,2}} \ = \ 0 \ \ \Longrightarrow \ \
\Delta\,(\,\partial_\lambda \,V_{\,\lambda \,,\  \xi}) \ +\ n\,(n\,+\,2) \cdot \left[\, V_{\,\lambda \,,\  \xi}^{4\over {n\,-\,2}}\cdot (\,\partial_\lambda \,V_{\,\lambda \,,\  \xi}\,) \right] \ = \ 0\ ,
$$
and the integration by parts formula.
Similarly for other derivatives. Hence it is required of $\,h\,$ to satisfy the condition
$$
 \langle\ h\,, \ [\ (\lambda\,\cdot\,\partial_\lambda\,)\, V_{\,\lambda \,,\  \xi} \ ]  \ \rangle_{\,\btd} \ = \ \langle\ h\,, \ [\ (\lambda\,\cdot\,\partial_{\xi_{\,|_j}}\,)\, V_{\,\lambda \,,\  \xi} \ ]  \ \rangle_{\,\btd} \ = \ 0 \ \ \ \ \ \mfor \ \ j \ = \ 1\,, \ 2\,, \ \cdot \cdot \cdot\,, \ n\,.
$$
Equation (A.1.9) can be seen in weak form\,:  the problem is equivalent to  finding a $\,\phi\ \in \ {\cal D}^{\,1, \  2}\,$ such that
\begin{eqnarray*}
& \ & \int_{\R^n} \left\{\, \Delta \,\phi \ + \ n\,(\,n\,+\,2) \cdot \left( V_{\lambda\,,\,\xi}^{4\over{n\,-\,2}}\cdot \phi\right) \ - \ \left( V_{\,\lambda \,,\  \xi}^{4\over {n\,-\,2}}\cdot h\right)  \right\} \cdot \psi \ = \ 0\\[0.1in]
& \ & \hspace*{4in}  \ \ \ \ {\mbox{for \ \ all}} \ \ \ \psi \ \in \ {\cal D}^{\,1, \  2}\\[0.1in]
& \Longleftrightarrow &   \int_{\R^n} \langle\, \phi\,,\ \psi \,\rangle_\btd \,+ \left[\ - \,n\,(\,n\,+\,2)  \int_{\R^n} \left( V_{\lambda\,,\,\xi}^{4\over{n\,-\,2}}\cdot \phi\right) \cdot \psi \ \right]  \\[0.2in]
& \ & \hspace*{2in} \ = \ \left[\ -\, \int_{\R^n} \left( V_{\lambda\,,\,\xi}^{4\over{n\,-\,2}}\cdot h\right) \cdot \psi \ \right] \ \ \ \ \ \ \ \  {\mbox{for \ \ all}} \ \ \psi \ \in \ {\cal D}^{\,1, \  2}.
\end{eqnarray*}
Here we apply the   integration by parts formula.
Via the Riesz Representation Theorem, the problem  can be cast in the form of finding a solution (\,in $\,{\cal D}^{\,1, \  2}\,$)\, of the (functional) equation
$$
(\,{\bf Id} \ + \ {\bf K}\,) \,(\phi) \ = \  h_{\cal D}\,.\leqno (A.1.11)
$$
Here $\,{\bf Id}\,$ is the identity map interpreted by
\begin{eqnarray*}
  \langle\  {\bf Id} \,(\,\phi\,)\,, \ \psi\ \rangle_\btd     & =  & \langle\, \phi \,, \ \psi\,\rangle_\btd\ ,     \\[0.15in]
{\mbox{and}} \ \ \ \ \ \ \ \   \ \ \ \  \langle\, {\bf K }  \,(\,\phi\,)\,, \ \psi  \,\rangle_\btd  & = & \left[\ - \,n\,(\,n\,+\,2)  \int_{\R^n} \left( V_{\lambda\,,\,\xi}^{4\over{n\,-\,2}}\cdot \phi\right) \cdot \psi \ \right]  \ \ \ \ {\mbox{for \ \ all}} \ \ \psi \ \in \ {\cal D}^{\,1, \  2}  \,.   \ \ \ \ \ \ \ \   \ \ \ \ \ \ \ \   \ \ \ \ \ \ \ \
\end{eqnarray*}
From the decay of $\,V_{\lambda\,,\,\xi}\,$ as shown in (\,A.1.8\,)\,,\,  it can be checked that $\,{\bf K}\,: {\cal D}^{\,1, \  2} \ \to \ {\cal D}^{\,1, \  2}\,$ is a compact linear operator  (\,via the Rellich\,-\,Kondrachov Compactness Theorem \cite{Olsen-Holden}\ )\,.
[\,A compact linear operator is a linear operator $\,L\,$ from a Banach space $\,X\,$ to another Banach space $\,Y\,$, such that the image under $\,L\,$ of any bounded subset of $\,X\,$ is a relatively compact subset of $\,Y\,$ (\,i.e.\,,\, the closure is compact)\,.\,]\bk
Similarly, $\,h_{\cal D} \,\in\, {\cal D}^{\,1, \  2}\,$ is defined by
$$
 \langle\, h_{\cal D}\,, \ \psi  \,\rangle_\btd \ = \ \left[\ -\,\int_{\R^n}    \left( V_{\lambda\,,\,\xi}^{4\over{n\,-\,2}}\cdot h\right)   \cdot \psi \ \right] \ \  \ \ \ \ {\mbox{for \ \ all}} \ \ \psi \ \in \ {\cal D}^{\,1, \  2}.
$$

Note that\\[0.1in]
(\,A.1.12\,)
\begin{eqnarray*}
& \ &  \bigg\vert \ \int_{\R^n}  \left(  V_{\lambda\,,\,\xi}^{4\over{n\,-\,2}}\cdot h \right) \cdot \psi  \ \bigg\vert \\[0.15in]
& \le &  \left[\ \int_{\R^n}  \left(  V_{\lambda\,,\,\xi}^{4\over{n\,-\,2}}\cdot h \right)^{{2n}\over {n\,+\,2}}\, \right]^{ {{n\,+\,2}\over {2n}}} \cdot  \left( \ \int_{\R^n}  |\,\psi\,|^{{2n}\over {\,n\,-\ 2\,}}\  \right)^{\!\!{{n\,-\,2}\over {2n}}}\\[0.15in]
  & \le &  \left[\   \left( \ \int_{\R^n}  \left(  V_{\lambda\,,\,\xi}^{4\over{n\,-\,2}}  \right)^{ {{2n}\over {n\,+\,2}} \cdot  {{n\,+\,2}\over {4}} }\right)^{\!\!{{4}\over {n\,+\,2}}}  \cdot  \left( \ \int_{\R^n}  |\,h\,|^{ {{2n}\over {n\,+\,2}} \cdot  {{n\,+\,2}\over {\,n\,-\ 2\,}} }\right)^{\!\!{{n\,-\,2}\over {n\,+\,2}}}  \, \right]^{ {{n\,+\,2}\over {2n}}} \cdot  \left( \ \int_{\R^n}  |\,\psi\,|^{{2n}\over {\,n\,-\ 2\,}}\  \right)^{\!\!{{n\,-\,2}\over {2n}}}\\[0.15in]
  & \le & \left[\  \int_{\R^n}  \left(  V_{\lambda\,,\,\xi}^{4\over {n\,-\,2}}  \right)^{n\over 2} \right]^{\! {{2}\over {\,n }}}  \,\cdot  \,\Vert\,h\,\Vert_\btd \,\cdot  \,\Vert\,\psi\,\Vert_\btd \ \ 
  \\[0.15in]
 & \le & \left( \int_{\R^n}    V_{\lambda\,,\,\xi}^{{2n}\over{n\,-\,2}}\  \right)^{\!\!{2\over {\,n}}} \cdot  \Vert\,h\,\Vert_\btd \,\cdot  \,\Vert\,\psi\,\Vert_\btd \ = \  [\  V(\,n\,) \ ]^{\,{2\over {\,n}}} \cdot \Vert\,h\,\Vert_\btd \,\cdot  \,\Vert\,\psi\,\Vert_\btd \ .  \\
\end{eqnarray*}
Here
$$
  V(\,n\,) \ = \  \int_{\R^n} \left(\ {1\over {1\ + \ \Vert\,y\,\Vert^2 }}\ \right)^{\!n}\, d\,Y \ .
  $$

  \vspace*{0.1in}

{\it  Note.} \ \ If we replace
$$
V_{\lambda\,,\,\xi}^{4\over{n\,-\,2}}\cdot h
$$
in the above calculation by
$$
{\cal F} \,\cdot \, h\,, \ \ \ \ {\mbox{where}} \ \ {\cal F}\,\in\,L^{\,{n\over 2} } (\,\R^n\,)\,,\,
$$

We likewise obtain
$$ \bigg\vert \ \int_{\R^n}  \left(\, {\cal F} \,\cdot \, h\right) \cdot \psi  \ \bigg\vert \ \le \  \left[\  \int_{\R^n}  |\,{\cal F}\,|^{n\over 2} \right]^{\! {{2}\over {\,n }}}  \,\cdot  \,\Vert\,h\,\Vert_\btd \,\cdot  \,\Vert\,\psi\,\Vert_\btd \ .  $$

  \newpage

\hspace*{0.5in}The form of (A.1.11) allows us to use Fredholm alternative method to solve the  equation\,.\,  The kernel   $\,{\cal N} \ \subset \ {\cal D}^{\,1\,,\,2}\,$ of
$$
(\,{\bf Id} \ + \ {\bf K}\,)
$$
is identified with  the solution space $\,{\cal N}\,$ of the  linear equation
$$
 \Delta\,\phi_o \ +\ n\,(n\,+\,2) \cdot V_{\,\lambda \,,\  \xi}^{4\over {n\,-\,2}}\cdot \phi_o \ = \ 0 \ \ \ \ \ \ \ (\phi_o \ \in  {\cal D}^{\,1\,,\,2})\ . \leqno (A.1.13)
$$
This can be precisely determined \cite{Progress-Book}\,:
$$
{\cal N} \ = \ {\mbox{Span}} \   \left\{ \  (\lambda\cdot \partial_\lambda\,)  \,V_{\,\lambda \,,\  \xi}\,,\ \ (\lambda\cdot \partial_{\xi_{|_1}}\,)  V_{\,\lambda \,,\  \xi}     \,,\ \cdot \cdot \cdot \,, \ \ (\lambda\cdot \partial_{\xi_{|_n}}\,) \, V_{\,\lambda \,,\  \xi}  \, \right\}\,. \leqno (\,A.1.14\,)
$$
Moreover, these $\,(\,n\,+\,1\,)\,$ functions are orthogonal to each other with respect to the inner product $\,\langle\ \,\,,\, \ \,\rangle_\btd\,$.  See, for examples, \cite{Progress-Book} \cite{I}\,.
Write
$$
{\cal D}^{\,1\,,\,2}  \ = \ {\cal N} \ \oplus \, {\cal N}_\perp\,,
$$
where $\,{\cal N}_\perp\,$ is the orthogonal complement of $\, {\cal N}\,.\,$ Let
$$
{\cal P}_{\cal N} \ : \ {\cal D}^{\,1\,,\,2}  \ \to \ {\cal N} \leqno (A.1.15)
$$
be the projection\,.
 By the self\,-\,adjointness of $\,(\,{\bf I} \ + \ {\bf K}\,)\,$,
$$
(\,{\bf Id} \ + \ {\bf K}\,)  \ : \ {\cal N}_\perp \ \to \ {\cal N}_\perp \ \ \ \ {\mbox{is \ \ an \ \ isomorphism\,.}} \leqno (A.1.16)
$$
Applying the Fredholm alternative method, equation (A.1.11) is solvable if and only if
 $$
h_{\cal D}  \ \in \ {\cal N}_\perp\,. \leqno (A.1.17)
 $$
Recalling relation (4.19), we have\\[0.1in](A.1.18)
\begin{eqnarray*}
h_{\cal D}  \ \in \ {\cal N}_\perp & \Longleftrightarrow & \bigg\langle  h_{\cal D} \,, \  (\lambda\cdot \partial_\lambda\,)  \,V_{\,\lambda \,,\  \xi} \ \bigg\rangle_\btd \ = \ 0 \ \ \Longleftrightarrow \ \ \bigg\langle \left( V_{\,\lambda \,,\  \xi}^{4\over {n\,-\,2}}\cdot h\right)\,, \ \,{{\partial\, V_{\,\lambda \,,\  \xi}}\over {\partial \lambda}}\ \bigg\rangle_{\!{\int}} \ = \ 0\\[0.2in]
& \Longleftrightarrow & \bigg\langle    h\,, \ {{\partial \, V_{\,\lambda \,,\  \xi}}\over {\partial \lambda}} \ \bigg\rangle_\btd  \ = \ 0 \ .
\end{eqnarray*}
Similarly for the other derivatives with respect to the components of $\,\xi\,.\,$ Therefore
 $$
  h \ \in \ {\cal N}_\perp\,. \leqno (A.1.19)
 $$

In general, for $\,h\,\not\in\,{\cal N}_\perp\,,\,$ let
$$
{\cal P}_{\cal N}\,(\,h\,) \ = \ - \ \left[ \  a_o \cdot \lambda \,{{\partial\, V_{\,\lambda \,,\  \xi}}\over {\partial \lambda}} \ + \ \sum_{j\,=\,1}^n b_j \cdot \lambda\, {{\partial\, V_{\,\lambda \,,\  \xi}}\over {\partial\, \xi_{j}}} \ \right]\ , \leqno (\,A.1.20\,)
$$and
$$
h \ = \ h_\perp \ \ + \ {\cal P}_{\cal N}\,(\,h\,) \ \ \Longleftrightarrow \ \  h_\perp \ = \ h \ - \ {\cal P}_{\cal N}\,(\,h\,)  \ \ \ \ \  \ \ \ ( \ h_\perp \ \in \ {\cal N}_\perp \ )\leqno (A.1.21)
$$
$$
\Longrightarrow \ \ \ \ h_\perp \ = \  h \ + \ a_o \cdot \lambda \,{{\partial\, V_{\,\lambda \,,\  \xi}}\over {\partial \lambda}} \ + \ \sum_{j\,=\,1}^n b_j \cdot \lambda\, {{\partial\, V_{\,\lambda \,,\  \xi}}\over {\partial\, \xi_{j}}}\ .
$$
The coefficients $\,a_o\,$ and $\,b_j\,$ in (\,A.1.20\,) are uniquely determined\,.\, The above discussion, together with elliptic regularity theory \cite{Gilbarg-Trudinger} ,  show that there is a  unique solution $\,\phi\,\in\, {\cal N}_\perp \,$ of the inhomogeneous equation\ :\\[0.1in]
(\,A.1.22\,) \ \ \ \ \
\begin{eqnarray*}
 \Delta \,\phi \ +\ n\,(n\,+\,2) \cdot \left(  V_{\,\lambda \,,\  \xi}^{4\over {n\,-\,2}}\cdot \phi\right) & = &\left(  V_{\,\lambda \,,\  \xi}^{4\over {n\,-\,2}}\cdot h_\perp \ \right) \\[0.2in] \Longleftrightarrow  \ \ \ \ \ \Delta \,\phi \ +\ n\,(n\,+\,2) \cdot \left(  V_{\,\lambda \,,\  \xi}^{4\over {n\,-\,2}}\cdot \phi\right) & = &\left(  V_{\,\lambda \,,\  \xi}^{4\over {n\,-\,2}}\cdot h\right)  \\[0.2in]
& \ & \!\!\!\!\!\!\!\!\!\!\!\!\!\!\!\!\!\!\!\!\!\!\!\!\!\!\!\!\!\!\!\!\!\!\!\!+ \ \left[\ a_o \cdot V_{\,\lambda \,,\  \xi}^{4\over {n\,-\,2}} \cdot \lambda\, {{\partial\, V_{\,\lambda \,,\  \xi}}\over {\partial \lambda}} \ + \ \sum_{{\,j}\,=\,1}^n b_{\,j} \cdot V_{\,\lambda \,,\  \xi}^{4\over {n\,-\,2}} \cdot\lambda\, {{\partial\, V_{\,\lambda \,,\  \xi}}\over {\partial \,\xi_{|_{\,j}}}} \ \right]\,.
\end{eqnarray*}
 We sum up the discussion in the following proposition. See also \cite{Wei-Yan} for more general situation\,.\\[0.2in]
{\bf Proposition A.1.23.} \ \ {\it Given $\,h\,\in\,{\cal D}^{1\,,\ 2}\,$,\, there exists a  unique $\,\phi\,\in\, {\cal N}_\perp \,$  which solves equation\,} (A.1.20)\,,\, {\it where the    coefficients $\,a_o\,$ and $\,b_{\,j}\,$} (\,$1\,\le\,j\,\le\,n\,$) {\it  are uniquely determined by $\,h\,\,$ via the projection in\,} (\,A.1.20\,)\,.\,

\newpage


{\bf \large \S\,A\,2.} \ \  {\bf Notes to
 Gram\,–-\,Schmidt orthogonalization process.} \\[0.2in]
 Let
 $$
{\bf V}_{l}\, (\,Y\,) \ = \ \left( {{\Lambda_{\,l}^2}\over {\Lambda_{\,l} \ +  \ \Vert \,Y \, - \, \Sigma_{\,l}\ \Vert^2 }}\  \right)^{{\,n\,-\,2\,}\over 2} \ \ \ \ \mfor \ \ Y \,\in\,\R^n \ \ \ \ \ (\,1 \,\le\,l\,\le\,\flat\,)\  .
 $$
Here $\,\Lambda_{\,l} \,\in\,\R^+\,$ and $\,\Sigma_{\,l} \,\in\,\R^n\,$.\, Set
\begin{eqnarray*}
 {\bf V}_{1_{\,o}} & = & {{\partial\, {\bf V}_1}\over {\partial\, \Lambda_{\,1}}}\ , \ \ \ {\bf V}_{1_{\,1}} \ = \ {{\partial\, {\bf V}_1}\over {\partial\, \Xi_{1_{|_1}} }}\ , \ \cdot \cdot \cdot\,, \ {\bf V}_{1_{\,n}} \ = \ {{\partial\, {\bf V}_{1}}\over {\partial \,\Xi_{1_{|_n}} }}\ ; \ \ \ \ \ \   \ \ \ \ \ \  \ \ \ \ \ \   \ \ \ \ \ \  \\[0.2in]
{\bf V}_{2_{\,o}} & = & {{\partial\, {\bf V}_2}\over {\partial\, \Lambda_2}}\ , \ \ \ {\bf V}_{2_{\,1}} \ = \ {{\partial\, {\bf V}_2}\over {\partial\, \Xi_{2_{|_1}} }}\ , \ \cdot \cdot \cdot\,, \ {\bf V}_{2_{\,n}} \ = \ {{\partial\, {\bf V}_2}\over {\partial \,\Xi_{2_{|_n}} }}\ ;\\[0.1in]
& \cdot &\\[0.05in]
& \cdot &\\[0.05in]
{\bf V}_{\flat_{\,o}} & = & {{\partial\, {\bf V}_{\flat}}\over {\partial\, \Lambda_{\,\flat}}}\ , \ \ \ {\bf V}_{\flat_{\,1}} \ = \ {{\partial\, {\bf V}_{\,\flat}}\over {\partial\, \Xi_{\,\flat_{|_1}} }}\ , \ \cdot \cdot \cdot\,, \ {\bf V}_{\flat_{\,n}} \ = \ {{\partial\, {\bf V}_{\flat}}\over {\partial \,\Xi_{\,\flat_{|_n}} }}\ . \ \ \ \ \ \  \ \ \ \ \ \   \ \ \ \ \ \\
\end{eqnarray*}
Via the result in \cite{I}\,,\, or a direct calculation\,,\,  we have
$$
\langle\  {\bf V}_{l_{\,j}}\,, \ {\bf V}_{l_{\,j}}\,\rangle_{\,\btd} \ = \ C_j \ > \ 0\,, \ \ \ \  \ \ {\mbox{where}} \ \ \ 1 \ \le \ l \ \le \ \flat\,, \ \ \  j \ = \ 0\,, \ 1\,, \ 2\,, \cdot \cdot \cdot\,, \ n \ .  \leqno (A.2.1)
$$
For $\,j \ \not=\ {\tilde j}\,$,\,  if we control the ``\,interaction\,"\, term [\, see part (\,b\,) in {\bf \S\,A\,3\,.\,i} \,]\,:
$$
\langle\  {\bf V}_{l_{\,j}} \,, \ {\bf V}_{{\tilde l}_{\,{\tilde j}}}\,\rangle_{\,\btd} \ = \ -\, \int_{\R^n} [\ \Delta\,(\, {\bf V}_{l_{\,j}}\,)\ ] \cdot  {\bf V}_{{\tilde l}_{\,{\tilde j}}}
\ = \ n\,(\,n\,-\,2\,)\cdot {{n\,+\,2}\over {\,n\,-\ 2\,}} \cdot \int_{\R^n} \left( {\bf V}_{\,l}^{4\over {n\,-\,2}} \cdot  {\bf V}_{l_{\,j}}\  \right) \cdot  {\bf V}_{{\tilde l}_{\,{\tilde j}}} \ ,
$$
and make it small, that is,
$$
[\,{\,\flat\,(n\,+\,1)\,-\,1}\,] \cdot \langle\, {\bf V}_{l_{\,j}} \,, \ {\bf V}_{{\tilde l}_{\,{\tilde j}}}\,\rangle_{\,\btd} \ = \ o_{\,{\bar\lambda}_{\,\,\flat}}\,(\,1) \ \ \ \  {\mbox{for}} \ \ \ \ l \, \not=\ {\tilde l} \ \ \ {\mbox{or}} \ \ \ j \ \not=\ {\tilde j}\,, \leqno (A.2.2)
$$
where $\, 0  \ \le  \ j\ \not= \ \tilde j \ \le n\,$,\, $\, 1  \ \le  \ l \ \le \ \flat\,$,\, then $\,\{ \ {\bf V}_{\,l_{\,j}} \,\}\,$ is a collection of linearly independent ``\,vectors\,"\,.\, Refer to (\,A.2.3\,) below. In the above
$$
o_{\,{\bar\lambda}_{\,\,\flat}}\,(\,1) \ \to \ 0 \ \ \ \ {\mbox{as}} \ \ \ {\bar\lambda}_{\,\,\flat} \,\to\,0\,.  \leqno (\,A.2.3\,)
$$
Based on (\,A.3.29\,) and (\,A.3.30\,) in {\bf \S\,A\,3.i\,(\,b\,)}\,,\,  (\,A.2.2\,) happens when
$$
(\,n\,-\,2\,)\cdot \gamma \ >  \ \sigma \ \ \  \ {\mbox{and}} \ \ \  \ {\bar{\lambda}}_{\,\flat} \ \ \,{\mbox{is \ \ small \ \ enough}}\ . \leqno (\,A.2.4\,)
$$
See (\,1.22\,) and (\,1.32\,) in the original text\,,\, for $\,\gamma\,$ and $\,\sigma\,$,\, respectively\,.\,
Introduce
\begin{eqnarray*}
{\bf {\cal u}}_{\,1} & = & {\bf V}_{1_{\,o}} \,,\\[0.15in]
{\bf {\cal u}}_{\,2} &  = &  {\bf V}_{1_{\,1}}\ - \ \left(\ {{\langle\  {\bf {\cal u}}_{\,1}\,, \ {\bf V}_{1_{\,1}}\,\rangle_{\,\btd} }\over { \langle\  {\bf {\cal u}}_{\,1}\,, \ {\bf {\cal u}}_{\,1}\,\rangle_{\,\btd}  }}\ \right) \cdot {\bf {\cal u}}_{\,1} \,,\\[0.1in]
& \ & \hspace*{1.2in}[\ \uparrow \ o\,(\,1)\ ]\\[0.2in]
{\bf {\cal u}}_{\,3}  &  = &  {\bf V}_{1_{\,2}} \ - \ \left(\ {{\langle\, {\bf {\cal u}}_{\,1}\,, \  {\bf V}_{1_{\,2}}\,\rangle_{\,\btd} }\over { \langle\ {\bf {\cal u}}_{\,1}\,, \ {\bf {\cal u}}_{\,1}\,\rangle_{\,\btd}  }}\ \right) \cdot {\bf {\cal u}}_{\,1} \ - \ \left(\ {{\langle\, {\bf {\cal u}}_{\,2}\,, \  {\bf V}_{1_{\,2}}\,\rangle_{\,\btd} }\over { \langle\  {\bf {\cal u}}_{\,2}\,, \ {\bf {\cal u}}_{\,2}\,\rangle_{\,\btd}  }}\ \right) \cdot {\bf {\cal u}}_{\,2}  \,, \ \ \ \ \ \  \ \ \\[0.1in]
& \ & \hspace*{1in}[\ \uparrow \ o\,(\,1) \ ]\\[-0.1in]
& : & \\[-0.1in]
& : & \\[0.1in]
(\,A.2.5\,) \ \ \ \ \ \ \ \ \ \  {\bf {\cal u}}_{\,\flat\,(n\,+\,1) }  &  = &   {\bf V}_{\flat_{\,n}} \ - \ \sum_{m\,=\,1}^{\,\flat\,(n\,+\,1)\,-\,1}  \left(\ {{\langle\, {\bf {\cal u}}_{\,m}\,, \  {\bf V}_{\flat_{\,n}}\,\rangle_{\,\btd} }\over { \langle\  {\bf {\cal u}}_{\,m}\,, \ {\bf {\cal u}}_j\,\rangle_{\,\btd}  }}\ \right) \cdot {\bf {\cal u}}_{\,m}\,.
\end{eqnarray*}
Under condition (\,A.2.4\,)\,,\, we conclude that  $\,\{ \ {\bf {\cal u}}_{\,j} \,\}\,$ is a collection of orthogonal vectors\,,\, which can be normalized\,.\, \{\,Cf. Lemma 2.27 in \cite{III}\,.\}

\vspace*{0.5in}


{\bf \large \S\,A\,3.} \ \ {\bf \large Solving the linear inhomogeneous equation -- the case of}\\[0.1in]
\hspace*{0.65in}  {\bf \large multiple bubbles\,.}\\[0.2in]
%
%
%
{\bf \large {\bf  \S\,A\,3.a.}}   \ \
From the single bubble case in {\bf \S\,A\,1}\,,\, we move on to consider  the linear inhomogeneous equation  with the number of bubbles tending to infinity\,:
\begin{eqnarray*}
(A.3.1) \ \ \ \ \ \ \ \ \ \ \ \ \ \ \ \ \ \ \Delta \,\phi \ + \ {\cal B}_{\,\flat}\cdot \phi & = &  \hbar \ + \ {\mbox{P}}_{  \flat_{\,\parallelsum}}\ , \ \ \ \ \ \ \ \ \ \ \phi \ \in \ {\cal D}^{\,1,\,\,2}_{\flat_{\,\perp}}\  ,\\[0.2in]
 {\mbox{where}} \ \ \ \ \ \ \ \ \ \ \ \ \,  \ \ \ \ \ \ \ \ \ \ \ \  \ \ \ \ \ \ \ \ \  {\cal B}_{\,\flat}  & = & \left(\ {{n\,+\,2}\over {\,n\,-\ 2\,}}\,\right)  \cdot ( \,c_n \cdot K)\cdot W_{\,\,\flat}^{4\over{n\,-\,2}} \ . \ \ \ \ \ \ \ \ \ \ \ \   \ \ \ \ \ \ \ \ \ \ \ \   \ \ \ \ \ \ \ \ \ \ \ \ \ \ \
\end{eqnarray*}
Here $\,\hbar \ \in \  {\cal D}^{1\,,\ 2}\ $ is given (\,with further conditions being introduced later\,)\,,\, and
$$
{\mbox{P}}_{\flat_{\,\parallelsum}} \ = \ \sum_{l\,=\,1}^\flat  \,\left\{ \  a_{\,l} \cdot  V_{\,l}^{4\over {n\,-\,2}} \cdot  \left[\, (\,\lambda_{\,\,l} \ \partial_{\lambda_{\,\,l}}\,) \,V_{\,l}\,\right] \ + \ \sum_{j ,=\,1}^n  b_{\,l\,,\,j} \cdot   V_{\,l}^{4\over {n\,-\,2}} \cdot  \left[\, \left(\,\lambda_{\,\,l} \ \partial_{\,\xi_{\,l_{|_{\,j}}}}  \right) V_{\,l}\,\right]   \ \right\}\ , \leqno (A.3.2)
$$
with the (\,to\,-\,be\,-\,determined\,)  coefficients $\,a_{\,l}\,$ and $\,b_{\,l\,,\,j}\,$.

\newpage

\hspace*{0.5in}We show that when the bubbles are ``\,well\,-\,separated\,"\,,\, similar existence result as Proposition A.1.21 can be established.

\vspace*{0.2in}

{\bf \large {\bf  \S\,A\,3.b.}}   \ \ {\bf Rescaling\,: from $\,V_{\,l}\,$ to $\,{\bf V}_{\,l}\,.$ } \ \ Emphasizing on the quasi\,-\,hyperbolic  distance between two bubbles, we introduce the rescaling
\begin{eqnarray*}
\Phi\,(\,Y) \  \ = \  {\bar\lambda}_{\,\,\flat}^{\,{{n\,-\,2 }\over 2}}\  \cdot \phi\,(\,{\bar\lambda}_{\,\,\flat}\cdot Y\,)  \ \ \ \ \ \ \ \  \left(\ Y \ = \ {y\over {\bar\lambda}_{\,\,\flat}} \ \ \Longleftrightarrow  \ \ \ y \ = \ {\bar\lambda}_{\,\,\flat}\cdot Y\right)\,.
\end{eqnarray*}
With $\,\phi\,$ as in (\,A.3.1\,)\,,\, we check that
$$
\Delta_Y \,\Phi \ + \ {\cal {\bf {\cal B}}}_{\,\flat}  \cdot \Phi \ = \ {\bf  H}_{\,\flat} \ + \ {{\bf {\cal P}}}_{\flat_{\,\parallelsum}}\ .   \leqno (A.3.3)
$$
Here
\begin{eqnarray*}
(A.3.4) \ \ \ \ \ \ {\cal {\bf {\cal B}}}_{\,\flat} & = &   \left(\ {{n\,+\,2}\over {\,n\,-\ 2\,}} \,\right) \cdot (\, c_n \cdot {\bf K})   \cdot {\bf W}_\flat^{4\over{n\,-\,2}} \\[0.2in]
   {\bf W}_{\,\flat}  & = & \sum_{{\,l} \ =\,1}^\flat {\bf V}_{\,l}\ , \\[0.2in]
    {\bf V}_{\,l}\,(\,Y\,)  &  = &\!\! \left(\  {{  \Lambda_{\,l} }\over {  \Lambda_{\,l}^2 \ + \ \Vert\,Y \ - \ \Xi_{\,l}\,\Vert^{\,2} }}\ \right)^{\!\!{{\,n\,-\,2\,}\over 2}  }, \ \ \ {\mbox{where}} \ \ \Lambda_{\,l}\  = \  {{\lambda_{\,\,l}}\over {{\bar\lambda}_{\,\,\flat}}}\ , \ \ \ \ \Xi_{\,l} \ = \ {{\xi_{\,l}}\over {{\bar\lambda}_{\,\,\flat}}}\ ,\\[0.2in]
   (\,c_n\cdot {\bf K}\,)\,(\,Y)  & = &   (\,c_n \cdot K\,)\,(\,y)   \ , \ \ \ \  \ \ \ \   \ \ \ \  \ \ \ \    \ \ \ \  \ \ \ \   {\mbox{where}} \ \ \ y \ = \ {\bar\lambda}_{\,\,\flat}  \cdot Y \ , \\[0.2in]
{\bf  H}_{\,\flat}\ (\,Y\,) & = &  {\bar\lambda}_{\,\,\flat}^{\,{{ n\,+\,2 }\over { 2 }}} \cdot \hbar_{\,\flat}\,(\,{\bar\lambda}_{\,\,\flat} \cdot Y)\ ,  \\[0.2in]
{\mbox{and}} \ \  \ \  \ \ \ \  \, { {\bf {\cal P}}}_{\,\flat_{\,\parallelsum}}\!\!& = &\!\! \sum_{l\ =\,1}^\flat\, \left\{ \, \sum_{j\ =\,1}^n  b_{\,l\,,\,j} \cdot   {\bf V}_{\,l}^{4\over {n\,-\,2}} \cdot  \left[\, (\,\Lambda_{\,l} \ \partial_{\,\Xi_{\,l_{|_{\,j}}}})\, {\bf V}_{\,l}\,\right]  +     a_{\,l} \cdot  {\bf V}_{\,l}^{4\over {n\,-\,2}} \cdot  \left[\, (\,\Lambda_{\,l} \ \partial_{\Lambda_{\,l}}\,) \,{\bf V}_{\,l}\,\right] \, \right\}.
\end{eqnarray*}
The orthogonal condition [\,cf. (\,A.1.18\,)\ ] is transformed to\\[0.1in]
(\,A.3.5\,)
$$
\bigg\langle {\bf V}_{\,l}^{4\over {n\,-\,2}} \cdot  \left[\, (\,\Lambda_{\,l} \ \partial_{\Lambda_{\,l}}\,)\, {\bf V}_{\,l}\,\right]\,, \ \  \Phi  \bigg\rangle_{\!\int} \ = \ 0 \ \ {\mbox{and}} \ \  \bigg\langle {\bf V}_{\,l}^{4\over {n\,-\,2}} \cdot  \left[\, (\,\Lambda_{\,l} \ \partial_{\,\Xi_{\,l_{|_{\,j}}}}\,) \,{\bf V}_{\,l}\,\right]\,, \ \  \Phi  \bigg\rangle_{\!\int} \ \,= \ 0
$$
for $\,j \ = \ 1\,,\ 2\,, \ \cdot \cdot \cdot\,, \ n\, $\  and  $\,{\,l} \ = \ 1\,,\ 2\,, \ \cdot \cdot \cdot\,, \ \flat\,$.

\newpage

{\bf \large {\bf  \S\,A\,3\,.\,c.}}   \ \  {\bf  The weighted sup\,-\,norm\,.\,}\ \
As in \cite{Wei-Yan}\,,\, according to the rescaled centers $\,\{ \ \Xi_{\,l}\ \}\,$ of the rescaled bubbles\,,\,  let us consider
functions $\,F\,$ and $\,G\,,\,$  continuous  on $\,\R^n\,,\,$ and set\\[0.1in]
(A.3.6) \
\begin{eqnarray*}
 & \ & \hspace*{-1in} {\bf W}_*^{\,\flat} \ = \ \left\{\  F \ \in \ C^o\,(\,\R^n) \ \,\Bigg\vert \ \ \Vert \,F\,\Vert_{\,*_{\,Y}} \ := \   \sup_{Y\,\in\,\R^n}    \left( \  {{ |\,F\,(\,Y\,)\,|}\over { {\displaystyle{\sum_{{\,l}\,=\,1}^{\,\flat}}} \left(\,{1\over {   \,  1\ + \ \Vert\,Y \,-\ \Xi_{{\,l}}\,\Vert \,  }}\, \right)^{\!\! {{n\,-\,2}\over 2 } \ +\ \tau_{\,>1} }
 } }     \ <  \ \infty \right) \  \right\}\ {\bf ,}  \\[0.1in]
(A.3.7) \ \ \  & \ & \\[0.1in] & \ & \hspace*{-1in}  {\bf W}_{**}^{\,\flat}  \ = \ \left\{\  G \ \in \ C^o\,(\,\R^n)  \ \,\Bigg\vert \ \ \Vert \,G\,\Vert_{\,**_{\,Y}} \   := \  \sup_{Y\,\in\,\R^n}  \left( \   {{ |\,G\,(\,Y\,)\,|}\over { {\displaystyle{\sum_{{\,l}\,=\,1}^{\,\flat}}} \left(\, {1\over {   \,  1\ + \ \Vert\,Y \,-\ \Xi_{{\,l}}\,\Vert\,  }}\,\right)^{\!\! {{n\,-\,2}\over 2 }\ +\ 2 \ +\  \tau_{\,>1} }
 } }   \right)  \ <  \ \infty\  \right\}\ {\bf .}
\end{eqnarray*}
Likewise,
\begin{eqnarray*}
& \ & \Vert \,f\,\Vert_{\,*_y} \ := \   \sup_{y\,\in\,\R^n} \left(    {{ |\,f\,(\,y\,)\,|}\over { {\displaystyle{\ \sum_{{\,l}\,=\,1}^{\,\flat}}} \left(\,{1\over {   \,  1\ + \ \Vert\,y \,-\ \xi_{{\,l}}\,\Vert \,  }}\, \right)^{\!\! {{n\,-\,2}\over 2 } \ +\ \tau_{\,>1} } \
 } }   \right)\     {\bf ,}  \\[0.2in]
 {\mbox{and}} \ \ \  \ \ \ \  \  \ \ \  \ \ \ \  \ \  \ & \ & \Vert \,g\,\Vert_{\,**_y} \   := \  \sup_{y\,\in\,\R^n}  \left(  {{ |\,g\,(\,y\,)\,|}\over {\  {\displaystyle{\sum_{{\,l}\,=\,1}^{\,\flat}}} \left(\, {1\over {   \,  1\ + \ \Vert\,y \,-\ \xi_{{\,l}}\,\Vert\,  }}\,\right)^{\!\! {{n\,-\,2}\over 2 }\ +\ 2 \ +\  \tau_{\,>1} } \
 } }  \right) \ {\bf .}  \ \ \  \ \ \ \  \ \   \ \ \  \ \ \ \  \ \   \ \ \  \ \ \ \  \ \
\end{eqnarray*}
With
$$
F \,(\,Y\,) \ = \ {\bar\lambda}_{\,\,\flat}^{{\,n\,-\,2\,}\over 2} \cdot f\,(\,y\,)\ \ \ \ \ \ \ \ (\,y \ = \ {\bar\lambda}_{\,\,\flat} \cdot Y\ )\,,\,
$$
we have
$$
\Vert \,F\,\Vert_{\,*_{\,Y}} \ = \ {1\over {{ (\,\bar\lambda}_{\,\flat }\,)^{\,\tau_{\epsilon}} } }  \cdot \Vert \,f\,\Vert_{\,*}\ , \leqno (A.3.8)
$$
and
$$
G \,(\,Y\,) \ = \ {\bar\lambda}_{\,\,\flat}^{{n\,+\,2}\over 2} \cdot {\tt g}\,(\,y\,) \ \ \Longrightarrow \ \
\Vert \,G\,\Vert_{\,**_{\,Y}} \ = \ {1\over {{ (\,\bar\lambda}_{\,\flat }\,)^{\,\tau_{\epsilon}} } }  \cdot   \Vert \,{\tt g}\,\Vert_{\,**_{\,Y}}\ . \leqno (A.3.9)
$$

Note that
$$
F \,\in\,{\bf W}_*^{\,\flat} \ \ \Longrightarrow \ \ |\,F\,(\,Y\,)\,| \ \le \ {1\over { \ \Vert\,Y\,\Vert^{{{ n\,-\,2}  \over 2}\ + \ \tau_{\,>1} \ - \ o_{\,+}\,(\,1\,) } \  }} \ \ \  \ \mfor \  \ \Vert \,Y\,\Vert \ \gg \ 1\ . \leqno (A.3.10)
$$
Define
\begin{eqnarray*}
(A.3.11) \ \ \ \ \ \ {\bf W}_{*_\perp}^{\,\flat} \!\!\!& = &\!\! \Bigg\{ \ F \,\in\,{\bf W}_{*}^{\,\flat} \ \ \bigg\vert \ \
\bigg\langle V_{\Lambda_{\flat_{\,l}}\,,\,\Xi_{\flat_{\,l}}}^{4\over {n\,-\,2}} \cdot \left(\ \Lambda_{\flat_{\,l}} \cdot {{\partial  V_{\Lambda_{\flat_{\,l}}\,,\,\Xi_{\flat_{\,l}}} }\over {\partial\, \Lambda_{\flat_{\,l}}}}\right)\,, \ F  \bigg\rangle_{\!\int} \ = \ 0\,, \\[0.2in]
& \ & \hspace*{1.2in}{\mbox{and}} \ \ \ \ \bigg\langle V_{\Lambda_{\flat_{\,l}}\,,\ \Xi_{\flat_{\,l}}}^{4\over {n\,-\,2}} \cdot \left(\ \Lambda_{\flat_{\,l}} \cdot {{\partial\, V_{\Lambda_{\flat_{\,l}}\,,\,\Xi_{\flat_{\,l}}} }\over {\partial\, \Xi_{{\flat_{\,l}}_{|_m}} }}\right)\,, \ F  \bigg\rangle_{\!\int} \ \,= \ 0 \\[0.15in]
& \ & \ \ \ \ \ \ \ \ \ \ \ \ \ \ \ \ \   {\mbox{for}}  \ \ \ \,{\,l} \ = \ 1\,,\ 2\,, \ \cdot \cdot \cdot\,, \ \flat\,, \ \  {\mbox{ and}} \ \   \,m \ = \ 1\,,\ 2\,, \ \cdot \cdot \cdot\,, \ n\,. \, \ \ \
 \Bigg\} \ .
\end{eqnarray*}

\vspace*{0.2in}

{\bf \large {\bf  \S\,A\,3.d.    }}   \ \  {\bf Preliminary estimates.} \\[0.1in]
{\it Separation Lemma.}  \ \ Refer to in Lemma B.1 in \cite{Wei-Yan}\,.\, Given numbers $\,\sigma\,$,\, $\,\alpha\,$ and  $\,\beta\,$ such that $$\,0 \ < \ \sigma \ \le \mbox{Min}  \ \{\,\alpha\,, \ \beta\,\}\,,\, \leqno (A.3.12)
$$ we have
\begin{eqnarray*}
(A.3.13) & \ &      \left(\ {{ 1}\over { 1\ + \ \Vert\,Y\ - \ \Xi_{\,\,1}\,\Vert  }}\right)^{\!\!  \alpha  } \cdot \left(\ {{ 1}\over { 1\ + \ \Vert\,Y \ - \ \Xi_{\,2}\,\Vert }}\right)^{\!\!  \beta  }
\\[0.15in]
& \le & C_1 \cdot    {{1}\over {  \Vert\,\Xi_{\,\,1} \ - \ \Xi_{\,2}\,\Vert^{\,\sigma} }}  \cdot \left\{ \left(\ {{ 1}\over { 1\ + \ \Vert\,Y\ - \ \Xi_{\,\,1}\,\Vert  }}\right)^{\!\!  \alpha \,+\, \beta \,-\, \sigma  } \ \ + \ \left(\ {{ 1}\over { 1\ + \ \Vert\,Y \ - \ \Xi_{\,2}\,\Vert  }}\right)^{\!\!  \alpha \,+\, \beta \,-\, \sigma  } \ \right\}
\end{eqnarray*}
for (\,all\,) $\,Y \,,\, \ \Xi_{\,\,1} \,\not=\,\Xi_{\,2} \,\in\, \R^n\,.\ $
Here the positive constant $\,C_1\,$ depends only on $\,\sigma\,$,\, $\,\alpha\,$ and  $\,\beta\,$ only.

\vspace*{0.2in}


{\it Condensation Lemma.}  \ \  Refer to in Lemma B.2 in \cite{Wei-Yan}\,.\,  For a fixed number $\,\tilde\varsigma\ $ such that
$$
2 \ < \ \tilde\varsigma \ < \ {{n } }\,,\, \leqno (A.3.14)
$$
we have
$$
   \int_{\R^n} \left(\ {1\over { 1\,+\,\Vert\,Z\,\Vert  }}\  \right)^{\! {  \tilde\varsigma} }
   \cdot  {1\over {  \Vert \,Y_c \,-\,Z\,\Vert^{n\,-\,2}  }} \ d\,Z
\   \le \   C_2 \cdot  \left(\ {1\over { 1\,+\,\Vert\,Y_c\,\Vert   }}\,\right)^{\!   \tilde\varsigma\ -\,2 } \ \ .  \leqno (A.3.15)
$$
Here $\,Y_c \ \in \ \R^n\,$ is fixed\,,\,
and the positive constant $\,C_2\,$ depends only on $\,\tilde\varsigma\,$ and $\,n\,.\,$

\vspace*{0.2in}

{\it ``\,Down\,-\,grading\,"\, Lemma.}  \ \  Refer to in Lemma B.1 in \cite{Wei-Yan}\,.\, For an integer $\,\flat \,\ge\,2\,,\,$ together with numbers $\,\kappa \, \in \, (\,0\,,\ 2)\,$ and $\,\theta \ >  \ 0\,$ (\,small enough\,,\, fixed\,)\,,\, we have
\begin{eqnarray*}
(A.3.16) \ \ & \  &   \int_{\R^n} \left[\ \sum_{{\,l} \ =\,1}^\flat \left(\ {1\over { 1 \ + \ \Vert\,Z \,-\,\Xi_{\,l}\,\Vert }}\ \right)^{\!\! {{\,n\,-\,2\,}\over 2}\ +\  \kappa  } \ \right]  \cdot {1\over { \Vert \,Y \,-\,Z\, \Vert^{\,n\,-\,2}  }}  \cdot  \left[\ {\bf W}_\flat(Z)\ \right]^{4\over {n\,-\,2}}   \ d\,Z\ \\[0.2in]
& \le & C_3    \cdot \left[\ \sum_{{\,l} \ =\,1}^\flat \left(\ {1\over {1 \,+\,\Vert\,Y \,-\,\Xi_{\,l}\,\Vert  }}\ \right)^{\!\!{{\,n\,-\,2\,}\over 2} \,+\  {\kappa }\,+\, {\theta }} \ \right] \ \ \ \ \ \ \ \ \ \ \ \ \ \ \ {\mbox{for \ \ all}} \ \ \ \ Y \,\in\,\R^n\,.
\end{eqnarray*}
 Here the positive constant $\,C_3\,$ depends only on $\,\kappa\,$,\, $\theta\,$ and $\,n\,.\,$\bk
  For the proofs of these three statements\,,\, see \cite{Wei-Yan}\,.

\vspace*{0.2in}

We often make use of the following estimates too (\,via similar calculation as in \S\,3 in [\,{\bf I}\,]\ )\,.\\[0.1in]
(A.3.17)
\begin{eqnarray*}
\bigg\vert \ \left(\ \lambda \cdot {{\partial }\over {\partial\, \lambda}}\  \right) V_{\,\lambda \,,\  \xi}\ \bigg\vert & \le &  C \cdot V_{\,\lambda \,,\  \xi}\ , \ \ \ \ \ \bigg\vert \ \left(\ \lambda \cdot {{\partial }\over {\partial\, \xi_{\,|_{\,j}} }}\  \right) V_{\,\lambda \,,\  \xi}\ \bigg\vert \ \le \   C \cdot V_{\,\lambda \,,\  \xi} \\[0.2in]
\bigg\vert \ \left(\ \lambda^2 \cdot {{\partial^2 }\over {\partial \,\lambda^2}}\  \right) V_{\,\lambda \,,\  \xi}\ \bigg\vert & \le &  C \cdot V_{\,\lambda \,,\  \xi}\ , \ \ \ \ \ \bigg\vert \ \left(\ \lambda^2 \cdot {{\partial^2 }\over {\partial \,\xi_{\,|_{\,j}}^2 }}\  \right) \  V_{\,\lambda \,,\  \xi}\ \bigg\vert \ \le \   C \cdot V_{\,\lambda \,,\  \xi}  \ ,\\[0.2in]
\bigg\vert \ \left(\ \lambda^2 \cdot {{\partial^2 }\over {\partial\, \lambda \ \partial\, \xi_{\,|_{\,j}} }}\  \right)    V_{\,\lambda \,,\  \xi}\ \bigg\vert & \le &   C \cdot V_{\,\lambda \,,\  \xi} \ , \ \ \ \ \  \bigg\vert \ \left(\ \lambda^2 \cdot {{\partial^2 }\over {\partial \,\xi_{\,|_{\,j}} \ \partial \, \xi_{\,|_{\,k}} }}\  \right)    V_{\,\lambda \,,\  \xi}\ \bigg\vert \  \le \   C \cdot V_{\,\lambda \,,\  \xi} \\[0.1in]
& \ & \hspace*{2.3in}  \ \ \ \ \ \ \ \ \ \ \ {\mbox{for }}\ \ 1 \ \le \  j\,, \ k \ \le  \ n\ .
\end{eqnarray*}
In many places we deal with integrals of the following type.
\begin{eqnarray*}
 (A.3.18) \ \ \ \ \ \ \  \ \int
 {{  {{\bar\lambda}}_{\,\flat}^B \cdot r^{\,D}  }\over {(\,{{\bar\lambda}}_{\,\flat}^2 \ + \  r^{\,2})^{{A }\over 2} }} \cdot r^{n\,-\,1} \cdot dr
  & = &   \int              {{ {{\bar\lambda}}_{\,\flat}^{\,B\,+\,D\,+\,n} \cdot   \left(\ {r\over {{{\bar\lambda}}_{\,\flat} }}\right)^{\!D}   }\over {{{\bar\lambda}}_{\,\flat}^{\,A} \cdot \left[\ 1  \ + \  \left(\ {r\over {{{\bar\lambda}}_{\,\flat} }}\right)^{\!2}\ \right]^{{A }\over 2} }}  \cdot {{ r^{n\,-\,1} \cdot dr}\over {{{\bar\lambda}}_{\,\flat}^n}}\\[0.2in]
  & = &  {{\bar\lambda}}_{\,\flat}^{\,B\,+\,D\,-\,A\,+\,n} \cdot   \int           {{ R^{\ D\,+\,n\,-\,1}   }\over {\left[\ 1  \ + \  R^2\ \right]^{{A }\over 2} }}  \cdot dR\ . \ \ \ \ \ \ \ \ \ \ \ \ \ \ \
\end{eqnarray*}
Here $\,A\,,\, \ B\,$ and $\,D\,$ are (\,fixed\,)  positive numbers.
Thus
\begin{eqnarray*}
A \ >  \ D \ + \ n  \ \ \Longrightarrow \ \ \int_{\rho_\nu}^{\,\infty}
 {{  {{\bar\lambda}}_{\,\flat}^B \cdot r^{\,D}  }\over {(\,{{\bar\lambda}}_{\,\flat}^2 \ + \  r^{\,2})^{{A }\over 2} }} \cdot r^{n\,-\,1} \cdot dr & \le  & C_1 \cdot {{\bar\lambda}}_{\,\flat}^{\,B\,+\,D\,-\,A\,+\,n} \cdot  {1\over { \left(\ {{ \rho_{\,\nu} }\over {  {{\bar\lambda}}_{\,\flat}  }}\right)^{A\,-\,(\,D\,+\,n\,)} }}\\[0.2in]
 & \le & C_1\cdot {{\bar\lambda}}_{\,\flat}^{\,B\,+\,D\,-\,A\,+\,n} \cdot  (\,{{\bar\lambda}}_{\,\flat}^{1\,-\ \nu}\,)^{A\,-\,(\,D\,+\,n\,)}\ ,
\end{eqnarray*}
where
$$
    {\tilde \lambda}_{\,\flat}  \ = \  \sqrt[^\flat\,]{\,\,\lambda_{\,1} \cdot \cdot \cdot \lambda_{\,\flat}\,\, } \ \ \  [ \ = \ o_{\,+}\,(\,1\,)\ ] \ \  \ \ \ \ \ \ \ {\mbox{and}} \ \ \ \  \ \  \rho_{\,\nu} \ = \ {\bar\lambda}_{\,\,\flat}^{\,\nu}\ .
$$
Likewise,
\begin{eqnarray*}
& \ & \\
A \ < \ C \ + \ n  \ \ \Longrightarrow \ \ \int^{\,\rho_\nu}_0
 {{  {{\bar\lambda}}_{\,\flat}^B \cdot r^{\,C}  }\over {(\,{{\bar\lambda}}_{\,\flat}^2 \ + \  r^{\,2})^{{A }\over 2} }} \cdot r^{n\,-\,1} \cdot dr & \le  & C_1 \cdot {{\bar\lambda}}_{\,\flat}^{\,B\,+\,C\,-\,A\,+\,n} \cdot  \left(\ {{ \rho_{\,\nu} }\over {  {{\bar\lambda}}_{\,\flat}  }}\right)^{ (\,C\,+\,n\,)\,-\,A}  \\[0.2in]
 & \le & C_2\cdot {{\bar\lambda}}_{\,\flat}^{\,B\,+\,C\,-\,A\,+\,n} \cdot  {1\over { \ (\,{{\bar\lambda}}_{\,\flat}^{1\,-\ \nu}\,)^{\,[\,(\,C\,+\,n\,)\,-\,A\,]}\, }} \\[0.2in]
 & = & C_2\cdot {{\bar\lambda}}_{\,\flat}^{\,B\,+\,\nu \,\cdot\, (\,C\,-\,A\,+\,n\,)} \ .
\end{eqnarray*}
Also\,,\, \\[0.1in]
(A.3.19)
\begin{eqnarray*}
& \ & \int_{B_{\,\xi_{\,1}}(\,\rho_{\,\nu}\,)} \left(\ {{{\bar\lambda}}_{\,\flat}\over {{{\bar\lambda}}_{\,\flat}^2 \ + \ |\,y\,-\,\xi_{\,1}\,|^2 }}\  \right)^{{\,n\,-\,2\,}\over 2} \ d \,y\ \le \ C \cdot \int_0^{\ \rho_{\nu}}  \left(\ {{{\bar\lambda}}_{\,\flat}\over {\ {{\bar\lambda}}_{\,\flat}^2 \ + \ r^2 \ }}\  \right)^{{\,n\,-\,2\,}\over 2} \ r^{n\,-\,1} \cdot dr\\[0.2in]
& = & C \cdot \int_0^{\ \rho_{\nu}}  \left(\ {{{\bar\lambda}}_{\,\flat}^2\over {\ \lambda^2_1 \ + \ r^2 \ }}\  \right)^{{\,n\,-\,2\,}\over 2} \ \cdot \ {1\over { {{\bar\lambda}_{\,\,\flat}}^{{\,n\,-\,2\,}\over 2}   }} \ \cdot \  r^{n\,-\,1} \cdot dr\\[0.2in]
& = & C  \cdot \int_0^{\ \rho_{\nu}} {{  1 }\over {\  \left[\  1 \ + \ \left(\  {r\over {{ {\bar\lambda}_{\,\,\flat} }}}\  \right)^{\,2} \ \right]^{\,{{\,n\,-\,2\,}\over 2}  \  }} }\,\cdot \ { {\bar\lambda}_{\,\,\flat} }^{{n\,+\,2}\over 2} \, \cdot \, {{ \,r^{n\,-\,1} \cdot dr\,}\over { { {\bar\lambda}_{\,\,\flat} }^n  }} \\[0.2in]
& = & C \cdot { {\bar\lambda}_{\,\,\flat} }^{\! {{n\,+\,2}\over 2} } \cdot \int_0^{\ \rho_{\nu}}\  {{\ R^{\, n\,-\,1 } \cdot dR \ }\over {\  \  \left[\  1 \ + \ R^2 \ \right]^{\,{{\,n\,-\,2\,}\over 2}  } \ \ } }  \ \le \ C_1 \cdot { {\bar\lambda}_{\,\,\flat} }^{\! {{n\,+\,2}\over 2} } \cdot \left(\ {{\rho_{\,\nu}}\over { {\bar\lambda}_{\,\,\flat} }}\  \right)^{\!\!2}\ .\\
\end{eqnarray*}
In the same spirit,
\begin{eqnarray*}
  & \ & \int_{\rho_{\,\nu}}^\infty
\left[   \  \left(\ {{{ {\bar\lambda}_{\,\,\flat} }}\over {\lambda^2_1 \,+\, r^2}}\  \right)^{\!\!2}\!\!
\cdot {{  r^2  }\over {(\,{ {\bar\lambda}_{\,\,\flat} }^2 \ + \  r^{\,2})^{{n }\over 2} }} \ \right] \cdot r^{n\,-\,1} \cdot dr \\[0.2in]
& = & \int_{\rho_{\,\nu}}^\infty
\left[   \  \left(\ {{{ {\bar\lambda}_{\,\,\flat} }^2 }\over {\lambda^2_1 \,+\, r^2}}\  \right)^{\!\!2}\!\!
\cdot \left( {{  { {\bar\lambda}_{\,\,\flat} }^2 }\over {\,{ {\bar\lambda}_{\,\,\flat} }^2 \ + \  r^{\,2} }}\  \right)^{{n }\over {\,2}}  \ \right] \cdot {{ \ r^{n\,+\,1} \cdot dr \ }\over {{ {\bar\lambda}_{\,\,\flat} }^{n\,+\,2} }} \\[0.2in]
  & = &  \int_{{\rho_{\,\nu}}\over {{ {\bar\lambda}_{\,\,\flat} }}}^\infty\  {{R^{n\,+\,1}\ d\,R}\over { \ (\,1 \ + \ R^2 )^{{n\,+\,4}\over 2} \   }}\  \ \ \ \ \ \  \ \ \ \ \ \  \ \ \ \ \ \  \ \ \ \ \ \  \ \ \ \ \ \  \ \ \ \ \  \left( \ R  \ = \ {r\over {{ {\bar\lambda}_{\,\,\flat} }}} \ \right) \\[0.15in]
&\le  & C_2 \int_{{\rho_{\,\nu}}\over {{ {\bar\lambda}_{\,\,\flat} }}}^\infty  {{R^{n\,+\,1}\cdot d\,R}\over { R^{n\,+\,4} }} \  \le\  C_3 \cdot \left(\ {{{ {\bar\lambda}_{\,\,\flat} } }\over {\rho_{\,\nu} }}\  \right)^{\,2} \ = \ O\,\left(\ {\bar{\lambda}}_{\,\flat}^{\,2\,(\,1\,-\ \nu\,)}\  \right)\  \\[0.2in]
\end{eqnarray*}
and
\begin{eqnarray*}
& \ & \\
& \ & \int_{\rho_{\,\nu}}^\infty
\left[    \left(\ {{{ {\bar\lambda}_{\,\,\flat} }}\over {\lambda^2_1 \,+\, r^2}}\  \right)^{\!\!2}\!\!
\cdot {{  { {\bar\lambda}_{\,\,\flat} }^2  }\over {(\,{ {\bar\lambda}_{\,\,\flat} }^2 \ + \  r^{\,2})^{{n }\over 2} }}\right] \cdot r^{n\,-\,1} \cdot dr \\[0.2in]
& = &  \int_{\rho_{\,\nu}}^\infty
\left[   \, \left(\ {{{ {\bar\lambda}_{\,\,\flat} }^2 }\over {\lambda^2_1 \,+\, r^2}}\  \right)^{\!\!2}\!\!
\cdot \left( {{  { {\bar\lambda}_{\,\,\flat} }^2 }\over {\,{ {\bar\lambda}_{\,\,\flat} }^2 \ + \  r^{\,2} }}\  \right)^{{n }\over 2}   \right] \cdot {{ \ { {\bar\lambda}_{\,\,\flat} }^2 \cdot r^{n\,-\,1} \cdot dr \ }\over {{ {\bar\lambda}_{\,\,\flat} }^{n\,+\,2} }} \\[0.2in]
\\[0.15in]
& =   &   \int_{{\rho_{\,\nu}}\over {{ {\bar\lambda}_{\,\,\flat} }}}^\infty\  {{R^{n\,-\,1}\ d\,R}\over { \ (\,1 \ + \ R^2 )^{{n\,+\,4}\over 2} \   }}\ \le \ C_4 \int_{{\rho_{\,\nu}}\over {{ {\bar\lambda}_{\,\,\flat} }}}^\infty \ {{\ \ R^{n\,-\,1}\  d\,R\ \ }\over { R^{n\,+\,4} }} \  \le\  C_5 \cdot \left(\ {{{ {\bar\lambda}_{\,\,\flat} } }\over {\rho_{\,\nu} }}\  \right)^4 \ = \ O\,\left(\ {\bar{\lambda}}_{\,\flat}^{\,4\,(\,1\,-\ \nu\,)}\  \right)\ .\\[0.2in]
\end{eqnarray*}

\vspace*{0.5in}

{\bf \large {\bf  \S\,A\,3\,.\,e\,.    }}   \ \  Following  Lemma 2.1 in \cite{Wei-Yan}\,,\, we present the  following lemma in a version which fits into our general setting here. It basically says that if the first term in the right hand side of (\,A.3.3\,) is ``\,small\,", then the solution is ``\,small\," too\,.\,\bk
{\it For simplicity sake, we recognize that the numbers  $\,\{\, \lambda_{\,\,l}\, \}_{\,l\,=\,1}^{\,\flat} \,$ and the points  $\,\{\,\xi_{\,l}\, \}_{\,l\,=\,1}^{\,\flat} \,$ may vary with\,} $\flat\,.\,$ In case of a sequence of numbers $\,\{\,\flat_{\,i}\,\}\,$,\,  even though $\,\flat_{\,j} \ = \ \flat_{\,k\,}$,\, we accept  that    $\,\{\, \lambda_{\,\,l}\, \}_{\,l\,=\,1}^{\,\flat_{\,j}}  \,$ and  $\ \{\, \xi_{\,l}\, \}_{\,l\,=\,1}^{\,\flat_{\,j}} \,$  may be different from $\,\{\, \lambda_{\,\,l}\, \}_{\,l\,=\,1}^{\,\flat_{\,k\,}}  \,$ and  $\ \{\, \xi_{\,l}\, \}_{\,l\,=\,1}^{\,\flat_{\,k\,}} \,$ \,.\,

\newpage

{\bf (\,Smallness\,) Lemma A.3.20\,.\,} \ \ {\it For $\,n \, \ge \, 6\,,\,$\, consider a sequence of numbers $$\,\{\,\flat_{\,i}\,\}\,\subset\,[\ 2\,, \ \infty\,)\,,$$ and the corresponding sequences   $\,\{\, \lambda_{\,\,l}\, \}_{\,l\,=\,1}^{\,\flat_{\,i}}  \,$ and  $\ \{\, \xi_{\,l}\, \}_{\,l\,=\,1}^{\,\flat_{\,i}} \ $,}\,
{\it with}
$$ \ \ \ \ \ \ \ \ \  \ \ \ \  \ \
{\bar\lambda}_{\,\,\flat_{\,i}} \,\to\,0^{\,+} \ \ \ \  \ as \ \ \ \ i \ \to \ \infty \ \ \ \ \ \ \ \ \ \  \left(\ {\it{ recall \ \ that}} \ \  {\tilde \lambda}_{\,\flat_{\,i}}  \ = \  \sqrt[^{\flat_i}]{\,\,\lambda_{\,1} \cdot \cdot \cdot \lambda_{\,\flat_{\,i}} } \ \  \right).
$$
{\it Assume that the bubble parameters  $\,\{\, \lambda_{\,\,l}\, \}_{\,l\,=\,1}^{\,\flat_{\,i}}  \,$ and  $\ \{\, \xi_{\,l}\, \}_{\,l\,=\,1}^{\,\flat_{\,i}} \,$ satisfy the conditions\,} (\,1.4\,)\,,  (\,1.8\,)\,, (\,1.22\,)\,, (\,1.24\,)\,--\,(\,1.28\,)\,,\,  \,(\,1.30\,)\, {\it and}  \,(\,1.32\,)\, {\it of the main text}\,.\,
{\it Suppose that for each $\,\flat_{\,i}\,,\,$}
$$\,\Phi_{\,\flat_{\,i}} \, \in \, {\bf W}_{*_\perp}^{\,\flat_{\,i}}\ $$
{\it solves  equation\,}  \,(\,A.3.3\,)\  {\it \,for $\ {\bf  H}_{\,\flat_{\,i}}\, \in \, {\bf W}_{**}^{\,\flat_{\,i}} \ $,\, where}
$$\ \Vert\
 {\bf  H}_{\,\flat_{\,i}}\,\Vert_{\,**_{\,Y}} \ \to \ 0^{\,+}\ \ \ \ \, {\it{as}} \ \ \ \ {\bar\lambda}_{\,\,\flat_{\,i}} \,\to\,0^{\,+}\ .\leqno (\,A.3.21\,) $$
{\it Then we also have}
 $$\ \ \  \Vert\,\Phi_{\,\flat_{\,i}}\,\Vert_{\,*_{\,Y}} \ \to \ 0^{\,+}\ \ \ \ {\it{as}} \ \ \ \  {\bar\lambda}_{\,\,\flat_{\,i}} \,\to\,0^{\,+}\,.
 $$

\vspace*{0.3in}

{\it Remarks.} \ \ One may also replace condition (\,1.26\,) by
$$
\ \ \ \ \ \ \ \
 |\,( \,{\tilde c}_n\,K)\,(\,y\,) \ - \ n\,(\,n\,-\,2\,)\ | \ \le \  {\bar\lambda}_{\,\,\flat}^{\,\zeta}  \ \ \ \ \  \ {\mbox{for}} \ \ \ y\,\in\, B_{ \,\xi_{\,l} }\,(\,\rho_\nu \,)\   \ \ \ \ \ \ (\ l \ = \ 1\,,\, \ \cdot \cdot \cdot\,, \ \flat\,)\ , \leqno (\,A.3.22\,)
$$
together with the conditions (\,1.29\,) and (\,1.30\,) in the main text\,.\,   Concerning $\,\{\, \lambda_{\,\,l}\, \}_{\,l\,=\,1}^{\,\flat_{\,i}}  \,$ and  $\ \{\, \xi_{\,l}\, \}_{\,l\,=\,1}^{\,\flat_{\,i}} \ $,\, see the clarification after (\,1.6\,) of the main text\,.\,  \bk
 We follow the proof of Lemma 2.1 in \cite{Wei-Yan}\,,\, placing suitable adjustment based on the conditions presented in this article\,.\, It begins with a consideration on the opposite situation\,.\, Assume that the conclusion of Lemma A.3.20 does not hold\,,\, that is\,,\,
 $$\ \ \  \Vert\,\Phi_{\,\flat_{\,i}}\,\Vert_{\,*_{\,Y}} \ \,\not\!\!\!\longrightarrow \ 0^{\,+}\ \ \ \ {\mbox{as}}  \ \ \ \  {\bar\lambda}_{\,\,\flat_{\,i}} \,\to\,0^{\,+}\,.
 $$
 By a normalization and selecting a subsequence (\,if necessary\,)\,,\, we may accept that
 $$\ \ \  \Vert\,\Phi_{\,\flat_{\,i}}\,\Vert_{\,*_{\,Y}} \ = \ 1\ \ \ \ \ \ \ {\mbox{for}}  \ \ \ \  i \ = \ 1\,, \ 2\,, \ \cdot \cdot \cdot \ . \leqno (\,A.3.23\,)
 $$


\newpage

{\bf \large {\bf  \S\,A\,3.f.    }}   \ \  {\bf  Representation of the linear equation\,.\,}\ \ Using Green's function of the form
$$
- \, {{ 1 }\over {\ (\,n\,-\,2\,)\cdot \omega_n\ }}  \cdot {1\over { \  \Vert\,Y\ - \ Z\,\Vert^{\,n\,-\,2} \ }} \ \ \ \ \mfor \ \ Y \,\not=\,Z \ \in \ \R^n\,,
$$
it follows from equation (\,A.3.3\,) that
\begin{eqnarray*}
(\,A.3.24\,) \ \ \ \Delta_{\,Y} \,\Phi & = &  - \ \left( {{n\,+\,2}\over {\,n\,-\ 2\,}}\  \right) \cdot ( \, c_n \cdot {\bf K}\,)   \cdot {\bf W}_\flat^{4\over{n\,-\,2}}\cdot \Phi \ + \  {\bf H}_{\,\flat} \ + \ {{\bf {\cal P}}}_{\flat_{\,\parallelsum}}\  \ \  \ \ \ \   \left(\ \Phi\,\in\, {\cal D}^{\,1\,,\,2}\   \right) \\[0.3in]
\Longrightarrow \ \ \Phi\,(\,Y\,) & = &  {{1 }\over {(\,n\,-\,2\,)\cdot \omega_n}} \cdot \left( {{n\,+\,2}\over {\,n\,-\ 2\,}}\  \right) \cdot  \int_{\R^n} \,( \,c_n \cdot {\bf K}\,)\,( Z) \cdot [\,{\bf W}_\flat(Z)\,]^{4\over {n\,-\,2}}   \cdot  \!\Phi\,(Z) \,*\\[0.2in]
& \ & \hspace*{1.5in} *\, \left(\ \  {1\over { \  \Vert\,Y\ - \ Z\,\Vert^{\,n\,-\,2} \ }}\right) d\,Z  \ \ \cdot \cdot \cdot \cdot \cdot \  {\bf{(\,i)}}_{\,(\,A.3.24\,)}   \\[0.3in]
& \ &    \ \ \  - \, {{ 1 }\over {(\,n\,-\,2\,)\cdot \omega_n}} \cdot \left[\ \int_{\R^n} {\bf H}_{\,\flat}\,(Z)  \cdot\! \left(\  {1\over { \  \Vert\,Y\ - \ Z\,\Vert^{\,n\,-\,2} \ }}\right) d\,Z \    \right]\\[0.2in]
 & \ & \hspace*{3in}       \uparrow  \cdot \cdot \cdot \cdot \cdot \cdot\cdot \cdot \cdot \cdot\cdot\ {\bf{(ii)}}_{\,(\,A.3.24\,)}    \\[0.3in]
& \ &  \ \ \ \ \ \ \ \  - \, {{ 1 }\over {(\,n\,-\,2\,)\cdot \omega_n}} \cdot \left[\  \int_{\R^n}  {{\bf {\cal P}}}_{\flat_{\,\parallelsum}} (Z)\cdot\! \left(\ \ {1\over { \  \Vert\,Y\ - \ Z\,\Vert^{\,n\,-\,2} \ }}\right) d\,Z \ \right] \ .\, \\[0.1in]
 & \ & \hspace{3in}       \uparrow  \cdot \cdot \cdot \cdot \cdot \cdot\cdot \cdot \cdot \cdot\cdot\ {\bf{(iii)}}_{\,(\,A.3.24\,)}\\
\end{eqnarray*}
Here $\,\omega_n\,$ is the volume (measure) of the unit sphere in $\,\R^n\,.\,$ See\,,\, for example, pp. 84 \ \& \  109\, in \cite{McOwen}\,.


\newpage

{\bf \large {\bf  \S\,A\,3.\,h.    }}   \ \   {\bf The order.} \ \  The first term $\,{\bf{(i)}}_{\,(\,A.3.24\,)} \,$  in the right hand side of (\,A.3.24\,) can be estimated by
\begin{eqnarray*}
 (\,A.3.25\,) & \ & \bigg\vert\ \int_{\R^n} \,( c_n \cdot {\bf K})\,( Z) \cdot [\,{\bf W}_\flat(Z)\,]^{4\over {n\,-\,2}}   \cdot  \!\Phi\,(Z) \cdot\! \left(\ \ {1\over { \  \Vert\,Y\ - \ Z\,\Vert^{\,n\,-\,2} \ }}\right) d\,Z  \ \bigg\vert\\[0.15in]
 & \le & C \int_{\R^n} \  [\,{\bf W}_\flat(Z)\,]^{4\over {n\,-\,2}}   \cdot \left[\ {{ |\,\Phi\,(Z)\,| }\over {   \displaystyle{    \sum_{l\,=\,1}^\flat   \left(\  {{1}\over {     1\ + \ \Vert\,Z \,-\,\Xi_{\,l}\,\Vert }}\  \right)^{\!\! {{n\,-\,2}\over 2 } \,+\,\tau_{\epsilon}} \
  }} }\ \right] * \\[0.4in]
   & \ & \ \ \ \ \ \ \ \ \ \ \ \  *\   \left[\  \sum_{l\,=\,1}^\flat   \left(\  {{1}\over {     1\ + \ \Vert\,Z \,-\,\Xi_{\,l}\,\Vert }}\  \right)^{\!\! {{n\,-\,2}\over 2 } \,+\,\tau_{\epsilon}}  \, \right] \times  \left(\ {1\over { \  \Vert\,Y\ - \ Z\,\Vert^{\,n\,-\,2} \ }}\right)d\,Z\\[0.3in]
 & \le & C\! \cdot \!\Vert \,\Phi\,\Vert_{\,*_Y} \cdot \int_{\R^n}  [\,{\bf W}_\flat(Z)\,]^{4\over {n\,-\,2}}  \cdot \left[\   \sum_{l\,=\,1}^\flat    \left(\  {{1}\over {     1\ + \ \Vert\,Z \,-\,\Xi_{\,l}\,\Vert }}\  \right)^{\!\! {{n\,-\,2}\over 2 } \,+\,\tau_{\epsilon}} \  \right] \,* \\[0.3in]
 & \ & \hspace*{3in}\ \ \  \ \ \ \ *\, \ \left( {1\over { \  \Vert\,Y\ - \ Z\,\Vert^{\,n\,-\,2} \ }}\right) d\,Z\\[0.2in]
& \le &  C \cdot \left[\  \displaystyle{\sum_{l\ =\,1}^\flat }  \left(\  {1\over {    1\ + \ \Vert\,Y \,-\ \Xi_{\,l}\,\Vert  }}\ \right)^{\!  { {{n\,-\,2}\over 2 } \,+\,\tau_{\epsilon}\,+\,\theta}} \, \right]  \cdot  \Vert \,\Phi\,\Vert_{\,*_Y} \\[0.2in]
& \ & \hspace*{2.5in} \ \ \ \ \ \ \  \    [\  {\mbox{Down\,-\,grading \ \ Lemma\,,\ \ \ (\,A.3.16\,)}}\ ]\,.
\end{eqnarray*}

\newpage

{\bf \large {\bf  \S\,A\,3.h.    }}   \ \   {\bf The second and}\, (\,{\bf half of\,}\,) {\bf the third term.} \ \
Based on the same token (\,as in {\bf\S\,A\,3.\,g}\,)
\begin{eqnarray*}\emph{\emph{}}
 \ \ \ \ \ \ \ \ \ & \ &   \int_{\R^n}|\,{\bf H}_{\,\flat}\,(Z)\,  | \cdot {1\over { \Vert\,Y\ - \ Z\,\Vert^{\,n\,-\,2} }}\ \,d\,Z\\[0.15in]
& \le & C \!\cdot\! \Vert \,{\bf H}_{\,\flat} \Vert_{**_{\,Y}}  \int_{\R^n}   \left[\  \displaystyle{\sum_{l\,=\,1}^k }  \left(\ {1\over {    1\ + \ \Vert\,Z \,-\ \Xi_{\,l}\,\Vert  }}\ \right)^{\!\! { {{n\,+\,2}\over 2 } \ +\ \tau_{\epsilon}}} \, \right] \cdot  {1\over { \  \Vert\,Y\ - \ Z\,\Vert^{\,n\,-\,2} \ }}\ \,d\,Z\\[0.15in]
& \le & C' \cdot \Vert \,{\bf H}_{\,\flat}\, \Vert_{**_{\,Y}}  \cdot \left[\  \displaystyle{\sum_{l\,=\,1}^k }  \left(\ {1\over {    1\ + \ \Vert\,Y \,-\ \Xi_{\,l}\,\Vert  }}\ \right)^{\!\! { {{n\,+\,2}\over 2 } \,+\,\tau_{\epsilon}\,-\,2}} \, \right]\\[0.15in]
& \le & C' \cdot \Vert \,{\bf H}_{\,\flat}\, \Vert_{**_{\,Y}}  \cdot \left[\  \displaystyle{\sum_{l\,=\,1}^k }  \left(\ {1\over {    1\ + \ \Vert\,Y \,-\ \Xi_{\,l}\,\Vert  }}\ \right)^{\!\! { {{n\,-\,2}\over 2 } \,+\ \tau_{\,>1} }} \, \right]\ ,\\
\end{eqnarray*}
\begin{eqnarray*}
{\mbox{and}} \ \ \ \ & \ & \int_{\R^n} \left\{ \ \sum_{l\,=\,1}^\flat \left(\   {{\bf V}}_{\,l}^{4\over {n\,-\,2}} \cdot  \bigg\vert\, \Lambda_{\,l} \cdot  {{\partial {{\bf V}}_{\,l} }\over {\partial\, \Lambda_{\,l} }}    \,\bigg\vert\  \right)  \right\}\cdot {1\over { \Vert\,Y\ - \ Z\,\Vert^{\,n\,-\,2} }}\ \,d\,Z\\[0.15in]
& \le & C \cdot   \int_{\R^n}   \left[\  \displaystyle{\sum_{l\,=\,1}^\flat }  \left(\ {1\over {    1\ + \ \Vert\,Z \,-\ \Xi_{\,l}\,\Vert  }}\ \right)^{\!\!   n\,+\,2} \, \right] \cdot {1\over { \Vert\,Y\ - \ Z\,\Vert^{\,n\,-\,2} }}\ \,d\,Z\\[0.15in]
&  \le &  C' \cdot    \left[\  \displaystyle{\sum_{l\,=\,1}^\flat }  \left(\ {1\over {    1\ + \ \Vert\,Y \,-\ \Xi_{\,l}\,\Vert  }}\ \right)^{\!\! n } \, \right]\\[0.15in]
&  \le &  C' \cdot    \left[\  \displaystyle{\sum_{l\,=\,1}^\flat }  \left(\ {1\over {    1\ + \ \Vert\,Y \,-\ \Xi_{\,l}\,\Vert  }}\ \right)^{\!\! { {{n\,-\,2}\over 2 } \ +\ \tau_{\,>1} }} \, \right] \\[0.2in]
& \ & \hspace*{2in} \ \ \ \ \ \ \ \ \ \ \ \  \left(\ {\mbox{note \ \ that}} \ \ {{n\,+\,2}\over 2} \ > \ \tau_{\,>1} \ > 1 \ \right)\,.  \ \ \ \   \ \ \ \   \ \ \ \
\end{eqnarray*}
Recall that $\,\tau_{\,>1}\,$ is slightly {\it bigger}\, than one.
Here we use the Condensation Lemma (\,A.3.15\,)\,.\,
Note that, via (\,1.8) in the original text\,,\,
$$
{1\over {\,{\bar C}_{\,1}\,}} \ \le \ \Lambda_{\,l} \ \le \ {\bar C}_{\,1} \ \ \ \ {\mbox{for}} \ \ \ l \ = \ 1\,,\ 2\,, \ \cdot \cdot \cdot\,, \ \flat\,. \leqno (A.3.26)
$$


\newpage

{\bf \large {\bf  \S\,A\,3.i.    }}   \ \   {\bf Toward the smallness of the coefficients in} $\,{\tilde{\bf P}}_{\flat_{//}}\,.\,$ \ \
Multiplying  both sides of equation (A.3.3) by $\displaystyle{ \ \left(\Lambda_{\,1} \cdot {{ \partial\, {{\bf V}}_1\   }\over {\partial\, \Lambda_{\,1}}}\  \right)}$\ ,\, upon integration, we obtain
\begin{eqnarray*}
 (A.3.27) \ \ \ \ \ \ \ \ & \ & \bigg\langle \left(\ \Delta\,\Phi \ + \ \left\{\,{\tilde C}_{\,n} \cdot [\,c_n \cdot {\tilde K}(\,Y\,)\,] \cdot {\bf W}_\flat^{4\over {n\,-\,2}}\  \right\} \cdot   \Phi\right), \ \left(\Lambda_{\,1} \cdot {{ \partial\, {{\bf V}}_1,\   }\over {\partial\, \Lambda_{\,1}}}\  \right) \ \bigg\rangle_{\int}\ \ \ \ \ \ \ \ \ \ \ \ \ \ \ \ \ \ \ \ \   \\[0.2in]
  & = &   \bigg\langle \, {\bf H}_\flat\,,\  \left(\Lambda_{\,1} \cdot {{ \partial\, {{\bf V}}_1 \   }\over {\partial\, \Lambda_{\,1}}}\  \right) \ \bigg\rangle_{\int}  \ + \ \bigg\langle \ {\tilde{\bf P}}_{\flat_{//}}\,,\   \left(\Lambda_{\,1} \cdot {{ \partial\, {{\bf V}}_1,\   }\over {\partial\, \Lambda_{\,1}}}\  \right) \ \bigg\rangle_{\int}\ .\\[0.1in]
  & \ & \ \ [ \ (\,{\bf i}\,)_{\,(A.3.27\,)}\uparrow \ ] \ \ \  \ \  \ \ \ \  \ \  \ \ \ \ \ \ \ \ \ \ \ \ \ \ [ \ \uparrow \  (\,{\bf ii}\,)_{\,(A.3.27\,)}\ ]
\end{eqnarray*}

\vspace*{0.2in}

{\bf (\,a\,)} \ \ Making use of (A.3.17) and arguing as in (A.3.25)\,,\, we have\\[0.1in]
(A.3.28)
\begin{eqnarray*}
& \ & | \ (\,{\bf i}\,)_{\,(A.3.27\,)}\ |\\[0.15in]
& = & \bigg\vert  \ \bigg\langle \, {\bf H}_{\,\flat} \,,\  \left(\Lambda_{\,1} \cdot {{ \partial\, {{\bf V}}_1\   }\over {\partial\, \Lambda_{\,1}}}\  \right)\ \bigg\rangle_{\int} \ \bigg\vert\\[0.2in]
& \le & C_1\cdot  \int_{\R^n}|\,{\bf H}_{\,\flat}\,(Z)\,  | \cdot  \,{{{\bf V}}_1\,} \, \ \,d\,Z \ \ \ \ \ \ \ \  \  \ \ \ \ \ \ \ \  \ \ \ \ \ \ \  \  \ \ \ \ \ \ \ \  \ \ \ \ \ \ \  \  \ \ \ \ \ \ \ \  \   [\,{\mbox{via \ \ (\,A.3.17\,)}}\,]\\[0.2in]
& \le & C_2 \cdot \Vert \,{\bf H}_{\,\flat}\Vert_{**_{\,Y}} \cdot \int_{\R^n}   \left[\  \displaystyle{\sum_{l\,=\,1}^\flat }  \left(\ {1\over {    1\ + \ \Vert\,Z \,-\ \Xi_{\,l}\,\Vert  }}\ \right)^{\!\! { {{n\,-\,2}\over 2 }\,+\,2 \,+\ \tau_{\,>1}}} \  \right] * \\[0.2in] & \ & \hspace*{3in} * \,\left(\ {1\over {    1\ + \ \Vert\,Z \,-\ \Xi_{\,\,1}\,\Vert  }}\ \right)^{\! {2 \cdot {{n\,-\,2}\over 2 } }}  \ \,d\,Z\\[0.2in]
& = & C_3 \cdot \Vert \,{\bf H}_{\,\flat} \Vert_{**_{\,Y}} \cdot \int_{\R^n}   \left(\ {1\over {    1\ + \ \Vert\,Z \,-\ \Xi_{\,\,1}\,\Vert  }}\ \right)^{\!\! { {{n\,+\,2}\over 2 } \,+\ \tau_{\,>1}}} \cdot \left(\ {1\over {    1\ + \ \Vert\,Z \,-\ \Xi_{\,\,1}\,\Vert  }}\ \right)^{{\!\! n\,-\,2 }}  \ \,d\,Z \ + \ \\[0.2in]
& \ & \ \  + \ C_4 \cdot \Vert \,{\bf H}_{\,\flat}\Vert_{**_Y} \cdot \int_{\R^n}   \left[\  \displaystyle{\sum_{l\,\not=\,1} }  \left(\ {1\over {    1\ + \ \Vert\,Z \,-\ \Xi_{\,l}\,\Vert  }}\ \right)^{\!\! { {{n\,+\,2}\over 2 } \,+\ \tau_{\,>1}}} \, \right]* \\[0.2in] & \ & \hspace*{3in} * \,\left(\ {1\over {    1\ + \ \Vert\,Z \,-\ \Xi_{\,\,1}\,\Vert  }}\ \right)^{\,n\,-\,2}  \ \,d\,Z\\[0.7in]
& \le & C_5 \cdot \Vert \,{\bf H}_{\,\flat} \Vert_{**_{\,Y}}\ + \hspace*{1.5in} [ \ \downarrow \ \ {\mbox{Seperation \ \ Lemma}} \ (\,A.3.13\,)\ ]\\[0.15in]
 & \ &  \ \ \ \ \  \ \ \ + \ \ C_6\cdot \Vert \,{\bf H}_{\,\flat} \Vert_{**_Y} \cdot \int_{\R^n}   \left\{ \  \displaystyle{ \sum_{l\,\not=\,1} }  {1\over {\Vert \, \Xi_{\,l}\,-\, \Xi_{\,\,1}\,\Vert^{{n \,-\,2}\over 2 }   }} \cdot \left[\ \left(\ {1\over {    1\ + \ \Vert\,Z \,-\ \Xi_{\,l}\,\Vert  }}\ \right)^{\!\! { {{n  } \,+\ \tau_{\,>1}}}} \ + \right. \right. \\[0.2in]
 & \ & \hspace*{3in} \left. \left.+ \   \left(\ {1\over {    1\ + \ \Vert\,Z \,-\ \Xi_{\,\,1}\,\Vert  }}\ \right)^{\!\! {  {{n\,+\ \tau_{\,>1}}  } }}\  \right] \ \right\}   \,d\,Z\\[0.25in]
& \le & C_7 \cdot \Vert \,{\bf H}_{\,\flat} \Vert_{**_{\,Y}}\cdot \left[\ 1\ + \  \sum_{l\,\not=\,1}   \  {1\over {\Vert \, \Xi_{\,l}\,-\, \Xi_{\,\,1}\,\Vert^{{n \,-\,2}\over 2 }   }}  \  \right] \ \le \ \ C' \cdot \Vert \,{\bf H}_{\,\flat} \Vert_{**_Y}\cdot \left[\ 1\ + \ O\, \left(\ {\bar\lambda}_{\,\,\flat}^{ {{\,n\,-\,2\,}\over 2} \cdot \gamma  } \right)  \  \right]\\[0.2in]
&  \le &   C_8 \cdot \Vert \,{\bf H}_{\,\flat}\Vert_{**_{\,Y}}  \hspace*{3in}  \ \ \ \ \left[  \,{\mbox{refer \ \ to \ \ (\,1.22\,)}}\,\right]\ .\\
\end{eqnarray*}
Likewise, we estimate  terms involving  derivatives in $\,\xi\,$--\,directions.


\vspace*{0.2in}

{\bf (\,b)} \ \     \ \ {\it Weak interaction.} \ \ As for $\,(\,{\bf ii}\,)_{\,(\,A.3.27\,)}$\,,\,  recall that
$$
 { {\bf P}}_{\,\flat_{\,\parallelsum}} \ = \ \sum_{l\ =\,1}^\flat\, \left\{ \, \sum_{j\ =\,1}^n  b_{\,l\,,\,j} \cdot   {\bf V}_{\,l}^{4\over {n\,-\,2}} \cdot  \left[\, (\,\Lambda_{\,l} \ \partial_{\,\Xi_{\,l_{|_{\,j}}}})\, {\bf V}_{\,l}\,\right]  +     a_{\,l} \cdot  {\bf V}_{\,l}^{4\over {n\,-\,2}} \cdot  \left[\, (\,\Lambda_{\,l} \ \partial_{\Lambda_{\,l}}\,) \,{\bf V}_{\,l}\,\right] \, \right\}\ .
$$

 As in (\,A.3.28\,)
\begin{eqnarray*}
(A.3.29) \ \ \ \ \ \ \ \ & \ & \Bigg\vert  \ \Bigg\langle  \left(\ {{\bf V}}_{\,l}^{4\over {n\,-\,2}}  \cdot  {{\partial \, {{\bf V}}_{\,l}  }\over {\partial\, \Lambda_{\,l}}}\ \right) \,,\  \, \left(\ \Lambda_{\,1} \cdot  {{\partial {{{\bf V}}_1\,}  }\over {\partial\, \Lambda_{\,1}}}\  \right) \ \Bigg\rangle_{\int} \ \Bigg\vert \ \ \ \  \ \ \ \ \ \ \ \ \ \  \ \ \ \ (\,l\,\not=\,1\,)\ \ \ \ \ \ \ \ \ \ \ \ \\[0.15in]
& \le & C\cdot  \int_{\R^n} {{\bf V}}_{\,l}^{{\,n\,+\,2\,}\over {n\,-\,2}}  \cdot  \,{{{\bf V}}_1\,} \, \ \,d\,Z\\[0.15in]
& \le & C \cdot  \int_{\R^n}    \left(\ {1\over {    1\ + \ \Vert\,Z \,-\ \Xi_{\,{\,l}}\,\Vert^{\,2}  }}\ \right)^{\! {  {{n\,+\,2}\over 2 } }}   \cdot \left(\ {1\over {    1\ + \ \Vert\,Z \,-\ \Xi_{\,\,1}\,\Vert^{\,2}  }}\ \right)^{\! {  {{n\,-\,2}\over 2 } }}  \ \,d\,Z\\[0.15in]
& \le & C \cdot    {1\over {\ \Vert \, \Xi_{{\,l}}\,-\, \Xi_{\,\,1}\,\Vert^{n\,-\,2}  \ }}
 \ \ \ \ \ \ \  \  \ \ \ \ \ \ \  \  (\,l \ \not= \ 1\,)\,.
\end{eqnarray*}
Similarly we estimate
\begin{eqnarray*}
& \ & \Bigg\vert  \ \Bigg\langle  \left(\ {{\bf V}}_{\,l}^{4\over {n\,-\,2}}  \cdot  {{\partial \, {{\bf V}}_{\,l}  }\over {\partial\, \Xi_{\,l_{\,|_{\,j}}} }}\ \right) \,,\  \,  \left(\ \Lambda_{\,1} \cdot  {{\partial {{{\bf V}}_1\,}  }\over {\partial\, \Lambda_{\,1}}} \ \right)  \ \Bigg\rangle_{\int} \ \Bigg\vert\ .
\end{eqnarray*}
Also\,,\, fixing $\,l\,=\,1\,,\,$
$$
\Bigg\langle  \left(\ {{\bf V}}_{\,1}^{4\over {n\,-\,2}}  \cdot  {{\partial \, {{\bf V}}_{\,1}  }\over {\partial\, \Xi_{\,1_{\,|_j}}}}\ \right) \,,\  \,  \left(\ \Lambda_{\,1} \cdot  {{\partial {{{\bf V}}_1\,}  }\over {\partial\, \Lambda_{\,1}}} \ \right)  \ \Bigg\rangle_{\int} \ = \ 0 \ \ \ \  \ \mfor \ \ j \ =  \ 1\,,\, \ \cdot \cdot \cdot\,, \ n\,.
$$
[\,Cf. (\,A.1.12\,) and (\,A.10.11\,)\,.\,]   After summing up, the contributions amount to\\[0.1in]
(\,A.3.30\,)
\begin{eqnarray*}
& \ & \hspace*{-0.5in}\Bigg\langle \sum_{l\,=\,1}^\flat\, \left\{ \, \sum_{j\ =\,1}^n  b_{\,l\,,\,j} \cdot   {\bf V}_{\,l}^{4\over {n\,-\,2}} \cdot  \left[\, (\,\Lambda_{\,l} \ \partial_{\,\Xi_{\,l_{|_{\,j}}}})\, {\bf V}_{\,l}\,\right]  +     a_{\,l} \cdot  {\bf V}_{\,l}^{4\over {n\,-\,2}} \cdot  \left[\, (\,\Lambda_{\,l} \ \partial_{\Lambda_{\,l}}\,) \,{\bf V}_{\,l}\,\right] \, \right\}, \ \ \left(\ \Lambda_{\,1} \cdot  {{\partial {{{\bf V}}_1\,}  }\over {\partial\, \Lambda_{\,1}}} \ \right)  \ \Bigg\rangle_{\int} \\[0.2in]&  = & C\,(\,n\,) \cdot a_1 \\[0.2in]
& \ & \hspace*{-1in} \ \ \ \ \ \ \ \ \ \ \ \ + \  O \left(\,{\bar\lambda}_{\,\,\flat}^{\,(\,n\,-\,2\,)\cdot \, \gamma}\ \right) \cdot \,[ \  \max \ \{\  | \, b_{\,l\,,\,j}\,| \ \ | \ \flat \, \le  \, l \, \le 1\,, \ \ 1 \, \le  \, j \, \le \, n\, \}  \, + \,  \max \ \{\, | \, a_{\,l} \, | \ \ | \ 1 \, < \, j \, \le \, n \, \} \ ]\,.
\end{eqnarray*}
Here
$$
C\,(\,n\,) \ = \ \int_{\R^n} \left(\ \Lambda_{\,1} \cdot  {{\partial {{{\bf V}}_1\,}  }\over {\partial\, \Lambda_{\,1}}} \ \right)^{\!\!2}
$$
We finally conclude the smallness of the coefficient in {\bf \S\,A\,3\,.\,m}\,.\, See (\,A.3.49\,)\,.\,
Cf. also {\bf \S\,A\,6} and {\bf \S\,9\,.\,b}\,.

\vspace*{0.35in}

{\bf \large {\bf  \S\,A\,3.\,j\,.    }}   \ \ {\bf The Laplacian term.} \ \  Based on the perpendicular condition {\bf \S\,A\,3.e}\,,\, we have [ \ cf. (\,A.1.10\,)\,]
\begin{eqnarray*}
0 & =  & \bigg\langle \ \Delta\,\Phi\,, \ \  \left(\ \Lambda_{\,1} \cdot  {{\partial {{{\bf V}}_1\,}  }\over {\partial\, \Lambda_{\,1}}}\  \right)   \ \bigg\rangle_{\int}\ = \ \Lambda_{\,1} \cdot \bigg\langle \ \Phi\,, \ \  \Delta\,{{\partial {{{\bf V}}_1\,}  }\over {\partial\, \Lambda_{\,1}}} \ \bigg\rangle_{\int} \ = \ \Lambda_{\,1} \cdot \bigg\langle \ \Phi\,, \ \  {{\partial (\,\Delta_Y {{{\bf V}}_1\,} \,) }\over {\partial\, \Lambda_{\,1}}} \ \bigg\rangle_{\int}\\[0.15in]
& = & -\,{{n\,+\,2}\over {\,n\,-\ 2\,}} \cdot n\,(\,n\,-\,2\,) \cdot\Lambda_{\,1} \cdot \bigg\langle \ \Phi\,, \ \ \left(\  [\,{{{\bf V}}_1\,}  ]^{4\over {n\,-\,2}} \cdot {{\partial {{{\bf V}}_1\,}  }\over {\partial\, \Lambda_{\,1}}}  \right) \bigg\rangle_{\int} \ = \ 0 \ .\\
\end{eqnarray*}


Likewise,
\begin{eqnarray*}
(A.3.31) \ \ \ \ \ \ \ \ &  \ &  \bigg\langle \ (\,c_n\cdot {\bf K}\,)   \cdot {\bf W}_\flat^{4\over {n\,-\,2}}\   \cdot   \Phi\, , \ \  \left(\Lambda_{\,1} \cdot {{ \partial\, {{\bf V}}_1,\   }\over {\partial\, \Lambda_{\,1}}}\  \right) \ \bigg\rangle_{\int}\\[0.2in]
 & = & \bigg\langle \ [\  (\,c_n\cdot {\bf K}\,) \ - \ n\,(\,n\,-\,2\,) \ ]   \cdot {\bf W}_\flat^{4\over {n\,-\,2}}\   \cdot   \Phi\, , \ \  \ \ \left(\Lambda_{\,1} \cdot {{ \partial\, {{\bf V}}_1 \   }\over {\partial\, \Lambda_{\,1}}}\  \right) \ \bigg\rangle_{\int} \\[0.2in]
  & \ & \ \ \ \ \ \ \ \ \ \ \ \ \ \ \ \    + \ \,\bigg\langle \ n\,(\,n\,-\,2\,)  \cdot \left[\  {\bf W}_\flat^{4\over {n\,-\,2}}\ - \ [\,{{{\bf V}}_1\,}  ]^{4\over {n\,-\,2}}  \ \right]  \cdot   \Phi\, , \ \ \, \left(\Lambda_{\,1} \cdot {{ \partial\, {{\bf V}}_1,\   }\over {\partial\, \Lambda_{\,1}}}\  \right) \ \bigg\rangle_{\int}\ .
\end{eqnarray*}
Here we use the perpendicular condition one more time.
\newpage
\hspace*{0.5in}Let us consider [\ with similar argument as shown in (\,A.3.25\,)\ ]
\begin{eqnarray*}
 (A.3.32)   & \   & \bigg\vert \   \bigg\langle \  \left[\  {\bf W}_\flat^{4\over {n\,-\,2}}\ - \ [\,{{{\bf V}}_1\,}  ]^{4\over {n\,-\,2}}  \ \right]  \cdot   \Phi\, , \ \, \left(\Lambda_{\,1} \cdot {{ \partial\, {{\bf V}}_1 \   }\over {\partial\, \Lambda_{\,1}}}\  \right) \ \bigg\rangle_{\int} \ \bigg\vert\\[0.2in]
 & \le & C \cdot \int_{\R^n} |\,\Phi\,(Z)\,| \cdot \bigg\vert \ {\bf W}_\flat^{4\over {n\,-\,2}}\ - \ [\,{{{\bf V}}_1\,}  ]^{4\over {n\,-\,2}}  \ \bigg\vert \cdot \left( {1\over {[\ 1 \ + \Vert\,Z \, - \, \Xi_{\,\,1}\,\Vert\ ]^{\,n\,-\,2}}} \,\right) \ d\,Z\\[0.15in]
 & \le & C \cdot \Vert\,\Phi\Vert_{\,*_Y} \cdot  \int_{\R^n}   \bigg\vert \ {\bf W}_\flat^{4\over {n\,-\,2}}\ - \ [\,{{{\bf V}}_1\,}  ]^{4\over {n\,-\,2}}  \ \bigg\vert     * \\[0.2in]
 & \ & \ \ \ \  \ \ \ \  \   * \,\left( {1\over {(\,1 \ + \Vert\,Z \, - \, \Xi_{\,\,1}\,\Vert\,)^{\,n\,-\,2}}}\  \right)\cdot \left(\ \ \sum_{l\ =\,1}^\flat {1\over {   (  1\ + \ \Vert\,Z \,-\Xi_{\,l}\,\Vert\,)^{ {{n\,-\,2}\over 2 } \,+\ \tau_{\,>1} } }}\  \right) d\,Z\ .
\end{eqnarray*}
Note that
\begin{eqnarray*}
  \ \ & \ & \bigg\vert \  (\,a \,+\,s\,)^{4\over{n\,-\,2}} \ - \ a^{4\over {n \,-\,2}}  \ \bigg\vert
\ \le \ C_2 \cdot  s^{{4}\over {\,n\,-\,2\,}}\ \ \ \mfor \ \ a \ > \ 0\,, \
\ s \ >  \ 0\,, \ \ \    (\,n \ \ge \ 6\,)\\[0.15in]
{\mbox{and}} \ \ \ \ \ \ \ \,& \ & \ \
(\, a \ + \ b)^{\,\tau} \ \le \  a^{\,\tau} \ + \ b^{\,\tau} \ \ \mfor \ \ a\,, \ \,b \ \ge \ 0 \ \ \ \ {\mbox{and}} \  \ \ \tau \ \in \ (\,0\,, \ 1\,]\,.
\end{eqnarray*}
It follows that\\[0.1in]
(A.3.33)
\begin{eqnarray*}
   \left(\ {\bf W}_\flat^{4\over {n\,-\,2}} \ - \ {\bf V}_1^{4\over {n\,-\,2}}\ \right) & \le &  C  \cdot \left[\ {\bf V}_2^{4\over{n\,-\,2}} \ + \cdot \cdot \cdot \ + \ {\bf V}_\flat^{4\over{n\,-\,2}} \  \right] \\[0.15in]
&   \le &  C \cdot   \left[\ \left(\ {{\Lambda_{\,2} }\over {  \Lambda_{\,2}^2 \ + \ \Vert\,Z \ - \ \Xi_{\,2} \,\Vert^{\,2} }} \ \right)^{\!2 }  \ + \ \cdot \cdot \cdot \ + \  \left(\ {{\Lambda_{\,\flat} }\over {  \Lambda_{\,\flat}^2 \ + \ \Vert\,Z \ - \ \Xi_{\,\flat} \,\Vert^{\,2} }} \ \right)^{\!2 } \  \right]\\[0.15in]
&   \le & \ C_2 \cdot    \left[\ \left(\ {{1 }\over {  1 \ + \ \Vert\,Z \ - \ \Xi_{\,2} \,\Vert  }} \ \right)^{\!4 }  \ + \ \cdot \cdot \cdot \ + \  \left(\ {{1 }\over {  1 \ + \ \Vert\,Z \ - \ \Xi_{\,\flat} \,\Vert  }} \ \right)^{\!4 } \  \right]\ .
\end{eqnarray*}

\newpage

With (\,A.3.33\,)\,,\, we continue from (A.3.32)\\[0.1in]
%
%
(A.3.34)
\begin{eqnarray*}
 & \  &   \int_{\R^n}   \bigg\vert \ {\bf W}_\flat^{4\over {n\,-\,2}}\ - \ [\,{{{\bf V}}_1\,}  ]^{4\over {n\,-\,2}}  \ \bigg\vert     \cdot {1\over {(\,1 \ + \Vert\,Z \, - \, \Xi_{\,\,1}\,\Vert\,)^{\,n\,-\,2}}} \cdot \left(\ \ \sum_{l\ =\,1}^\flat {1\over {   (  1\ + \ \Vert\,Z \,-\Xi_{\,l}\,\Vert\,)^{ {{n\,-\,2}\over 2 } \,+\ \tau_{\,>1} } }}
\ \right) d\,Z\\[0.2in]
 & \le  &   C \cdot \int_{\R^n}  {1\over {(\,1 \ + \Vert\,Z \, - \, \Xi_{\,\,1}\,\Vert\,)^{\,n\,-\,2}}} \,*\,\\[0.15in]
  & \ & \   \ \ \ \ \ \ \ \ \ \ \ \ \ \ \ *\,\left(\ \ \sum_{l\ =\,1}^\flat {1\over {   (  1\ + \ \Vert\,Z \,-\Xi_{\,l}\,\Vert\,)^{ {{n\,-\,2}\over 2 } \,+\ \tau_{\,>1} } }}
\ \right) \cdot \left(\ \ \sum_{l\ =\,2}^\flat {1\over {   (  1\ + \ \Vert\,Z \,-\Xi_{\,l}\,\Vert\,)^{ 4} }}
\ \right) \,dZ \\[0.2in]
& = & C \cdot \int_{\R^n} {1\over {(\,1 \ + \Vert\,Z \, - \, \Xi_{\,\,1}\,\Vert\,)^{\,n\,-\,2\,+\, {{n\,-\,2}\over 2 } \,+\ \tau_{\,>1}  }}}  \,*\,\\[0.15in]
  & \ & \   \ \ \ \ \ \ \ \ \ \ \ \ \ \ \ *\, \left[\ \left(\ {{1 }\over {  1 \ + \ \Vert\,Z \ - \ \Xi_{\,2} \,\Vert  }} \ \right)^{\!4 }  \ + \ \cdot \cdot + \, \left(\, {{1 }\over {  1 \ + \ \Vert\,Z \ - \ \Xi_{\,\flat} \,\Vert^{\,2} }} \, \right)^{\!4 } \, \right] \,dZ  \\[0.2in]
& \ & \ \  + \ C \cdot \int_{\R^n} {1\over {(\,1 \ + \Vert\,Z \, - \, \Xi_{\,\,1}\,\Vert\,)^{\,n\,-\,2}}} \,*\,\\[0.15in]
  & \ & \   \ \ \ \ \ \ \ \ \ \ \ \ \ \ \ *\, \left(\ \ \sum_{l\ =\,2}^\flat {1\over {   (  1\ + \ \Vert\,Z \,-\Xi_{\,l}\,\Vert\,)^{ {{n\,-\,2}\over 2 } \,+\ \tau_{\,>1} } }}
\ \right) \cdot \left(\ \ \sum_{l\ =\,2}^\flat {1\over {   (  1\ + \ \Vert\,Z \,-\Xi_{\,l}\,\Vert\,)^{4} }}
\ \right) \,dZ  \\[0.2in]
& \le & C_1 \cdot \left[\ \sum_{l \ = \,2}^\flat \ {1\over { \,\Vert\,\Xi_{\,l} \ - \  \Xi_{\,\,1} \Vert^{\,4} }}  \ \right] \cdot \int_0^\infty {1\over {(\,1 \ + \ R\,)^{\,n\,-\,2\,+\, {{n\,-\,2}\over 2 } \,+\ \tau_{\,>1} } \ }} \  \,dR\ + \ \ \ \ \ \ \ \ [\ \leftarrow \ \cdot \cdot \cdot \ \ \mbox{(\,i)} \ ]\\[0.2in]
& \ &  \ \ \ \ \ \ \ \ \ \   + \ \, C \cdot \int_{\R^n}  {1\over {(\,1 \ + \Vert\,Z \, - \, \Xi_{\,\,1}\,\Vert\,)^{\,n\,-\,2}}} \,\,* \ \ \ \ \  \  \ \ \ \ \ \  \ \ \ \ \ \ \ [ \ {\mbox{via \ \ (\,A.3.13\,)}}   \ ] \\[0.2in]
& \ & \hspace*{-0.5in} *\, \left[\ \sum_{j\,=\,2}^\flat \  \sum_{l \ \not=\,1\,,\ j}^\flat \ {1\over { \,\Vert\,\Xi_{\,l} \ - \  \Xi_j \Vert^{\,4} }}  \left( {1\over {   (  1\ + \ \Vert\,Z \,-\Xi_{\,l}\,\Vert\,)^{ {{n\,-\,2}\over 2 } \,+\ \tau_{\,>1} } }}\,  +   {1\over {   (  1\ + \ \Vert\,Z \,-\Xi_{\,j}\,\Vert\,)^{ {{n\,-\,2}\over 2 } \,+\ \tau_{\,>1} } }}  \ \right) \right] \,dZ \\[0.1in]
& \ & \ \ \ \  \ \ \ \ \  \ \ \ \ \  \  [\ \uparrow  \ \ j\,=\,1 \ \ {\mbox{absorbed \ \ in \ \ the  \ \ first \ \ term \ \ (i)}} \ ]
\end{eqnarray*}

\newpage

\begin{eqnarray*}
& \le & C_2 \cdot\emph{} \left[\ \sum_{l \ = \,2}^\flat \ {1\over { \,\Vert\,\Xi_{\,l} \ - \  \Xi_{\,\,1} \Vert^{\,4\,} }}  \ \right]  \\[0.2in]
& \ &  \ \ \ \ \    + \ C_3 \cdot     \left[\ \sum_{j\,=\,2}^\flat \ \  \sum_{l \ \not= \,1\,; \, j}^\flat \ {1\over { \,\Vert\,\Xi_{\,l} \ - \  \Xi_j \Vert^{\,4\,} }}  \ \left(\  {1\over { \,\Vert\,\Xi_{\,\,1} \ - \  \Xi_j \Vert^{\,\tau_{\,>1}} }}  \ + \   {1\over { \,\Vert\,\Xi_{\,\,1} \ - \  \Xi_{\,l} \Vert^{\,\tau_{\,>1}} }}  \ \right) \right] \,*\,\\[0.15in]
  & \ & \   \ \ \ \ \ \ \ \ \ \ \ \ \ \ \ \ \ \ \ \ *\, \int_0^\infty {1\over {(\,1 \ + \ R\,)^{\,n\,-\,2\,+\, {{n\,-\,2}\over 2 } \,+\ \tau_{\,>1} } \ }} \  \,dR\ \\[0.2in]
& \le & C_4 \cdot  \left[\ O\,\left( \ {\bar\lambda}_{\,\,\flat}^{\ 4 \,\cdot \ \gamma } \   \ \right) \ + \ O\,\left(\  {\bar\lambda}_{\,\,\flat}^{(\,4\,+\ \tau_{\,>1}\,) \,\cdot \ \gamma } \   \ \right) \  \right]   \ .
\end{eqnarray*}
Note that $\,\displaystyle{n\, - \, 2 \ +\ {{\ n\,-\,2\ }\over 2} \ + \ \tau_{\,>1}\ >  \ n  \ \ \Longrightarrow \ \ n \ \ge \ 4}\,$.\,
In the above we make use of\\[0.1in]
(\,A.3.35)
\begin{eqnarray*}
& \ &  \sum_{j\,=\,2}^\flat \  \sum_{l \ \not= \,1\,; \  j}^\flat \ {1\over { \,\Vert\,\Xi_{\,l} \ - \  \Xi_j \Vert^{\,4\,} }}   \cdot   {1\over { \,\Vert\,\Xi_{\,\,1} \ - \  \Xi_{\,l} \Vert^{\,\tau_{\,>1}} }}\\[0.2in] & \le & \sum_{\,l\,=\,2}^\flat \left(\ \  \sum_{j \ \not= \,1\,; \ l}^\flat  {1\over { \,\Vert\,\Xi_{\,l} \ - \  \Xi_j \Vert^{\,4\,} }}  \ \right) \cdot   {1\over { \,\Vert\,\Xi_{\,\,1} \ - \  \Xi_{\,l} \Vert^{\,\tau_{\,>1}} }} \\[0.2in]
  & = &    {\bar\lambda}_{\,\,\flat}^{\ 4\, \cdot \ \gamma}   \cdot  \sum_{l \ \not= \,1\,; \, j}^\flat {1\over { \,\Vert\,\Xi_{\,\,1} \ - \  \Xi_{\,l} \Vert^{\,\tau_{\,>1}} }} \ \ \ \ \ \ \ \ \ (\,\tau_{\,>1} \ > \ 1\,)\\[0.2in]
      &= & O\,\left( \ {\bar\lambda}_{\,\,\flat}^{( \ 4\,\,+\, \tau_{\,>1}\ ) \,\cdot \ \gamma}  \ \right)\ .
\end{eqnarray*}

We are led to
$$
\bigg\vert \  \bigg\langle \ \Phi\ , \ \  \left[\,{\bf W}_\flat^{4\over {n\,-\,2}}\ - \ [\,{{{\bf V}}_1\,}  ]^{4\over {n\,-\,2}} \, \right] \cdot\left( \Lambda_{\,1} \cdot  {{\partial \,{{\bf V}}_1}\over {\partial\, \Lambda_{\,1}}} \ \right)  \ \bigg\rangle_{\int} \ \bigg\vert \ \le\ O\,\left( \, {\bar\lambda}_{\,\,\flat}^{\ 4 \,\cdot \, \gamma}  \ \right)\, \cdot\,  \Vert\,\Phi\Vert_{\,*_Y} \ . \leqno (A.3.36)
$$

\newpage

{\bf \large {\bf  \S\,A\,3.\,k\,.}}   \ \ [\ Refer also to {\bf {\bf \large {\bf  \S\,A\,3.\,j}}}\,.\ ] \ \
Next, let us consider\\[0.1in]
(A.3.37)
\begin{eqnarray*}
& \ & \bigg\vert \ \bigg\langle \left(\   \left\{ \ [\,c_n \cdot{\bf K}\,(\,Z\,)\,] \ - \ n\,(n\,-2)\,\right\} \cdot {\bf W}_\flat^{4\over {n\,-\,2}}  \cdot  \Phi\ \right)\,,  \ \ {{\partial {{{\bf V}}_1\,}  }\over {\partial\, \Lambda_{\,1}}}  \ \bigg\rangle_{\int}\ \bigg\vert  \\[0.15in]
& = & \bigg\vert \, \bigg\langle \left(\,   \left\{ \ [\,c_n \cdot{\bf K}\,(\,Z\,)\,] \, - \, n\,(n\,-2)\,\right\} \cdot \left[ \, \left(\ {\bf W}_\flat^{4\over {n\,-\,2}}  \,-\,  [\,{{{\bf V}}_1\,}  ]^{4\over {n\,-\,2}}  \ \right) \,+\,  [\,{{{\bf V}}_1\,}  \,]^{4\over {n\,-\,2}}  \ \right]    \cdot  \Phi\, \right)\,,  \ \ {{\partial\, {{{\bf V}}_1\,}  }\over {\partial\, \Lambda_{\,1}}}  \,\bigg\rangle_{\int} \bigg\vert  \\[0.15in]
 & \le & C \cdot \int_{\R^n} |\,\Phi\,(Z)\,| \cdot \bigg\vert \ {\bf W}_\flat^{4\over {n\,-\,2}}\ - \ [\,{{{\bf V}}_1\,}  ]^{4\over {n\,-\,2}}  \ \bigg\vert \cdot {1\over {(\,1 \ + \Vert\,Z \, - \, \Xi_{\,\,1}\,\Vert\,)^{\,n\,-\,2}}} \ \,d\,Z \ + \ \\[0.15in]
 & \ & \ \ \ \ \ \ \ \ \ \ \ \ \ \ + \ C \cdot \int_{\R^n}  \bigg\vert \, [\,c_n \cdot{\bf K}\,(\,Z\,)\,] \ - \ n\,(n\,-2) \,\bigg\vert \cdot [\,{{{\bf V}}_1\,}  ]^{4\over {n\,-\,2}}   \cdot  |\,\Phi\,| \cdot {1\over {(\,1 \,+\, |\, Z \,-\,\Xi_{\,\,1}\,|\,)^{\,n\,-\,2} }}\\[0.2in]
& \le &  C   \cdot \left[\ O\,\left(\, {\bar\lambda}_{\,\,\flat}^{\, 4\cdot \ \gamma}   \ \right)\  \right] \cdot \Vert\,\Phi\Vert_{\,*_Y} \ + \ \ \ \ \ \ \  \ \  \ \ \ \ \ \ \  \ \   [\, {\mbox{estimated \ \ as \ \ in \ \ (\,A.3.36\,)}}\,]\\[0.2in]
 & \ &   + \ C \cdot \Vert\,\Phi\,\Vert_{\,*_Y} \cdot \int_{\R^n}  \bigg\vert \, [\,c_n \cdot{\bf K}\,(\,Z\,)\,] \ - \ n\,(n\,-2)\,\bigg\vert   \cdot {1\over {(\,1 \,+\, \Vert\, Z \,-\,\Xi_{\,\,1}\,\Vert\,)^{\,n\,+\,2} }}\,*\\[0.1in]
 & \ & \hspace*{3in} * \ \left(\ \sum_{l\ =\,1}^\flat {1\over {   (  1\ + \ \Vert\,Z \,-\Xi_{\,l}\,\Vert\,)^{ {{n\,-\,2}\over 2 } \,+\ \tau_{\,>1} } }}
\ \right) d\,Z\\[0.15in]
& \le &  C   \cdot \left[\ O\,\left(\  {\bar\lambda}_{\,\,\flat}^{\, 4\cdot \ \gamma}  \ \right) \  \right]  \cdot \Vert\,\Phi\Vert_{\,*_Y}  \ + \ C \cdot \left[\ O\,\left(\ {\bar\lambda}_{\,\,\flat}^{( {{n\,+\,2}\over 2} \,+\ \tau_{\,>1}\,-\,\varepsilon \,) \cdot \ \gamma}  \ \right)   \  \right]  \cdot  \Vert\,\Phi\Vert_{\,*_Y}  \\[0.1in]
& \ & \ \ \ \ \ \ \ \ \ \ \ \ \ \   \ \ \ \ \ \   \ \ \ \ \ \  \ \  \ \ \ \ \ \ \ \ \ \ \ \  \ \ \ \ \ \   \ \ \ \ \ \    [\ \uparrow \ \ {\mbox{estimated \ \ as \ \ in \ \  {\bf{\S\,A\,3.g}}}}\ ]\\[0.2in]  & \ &   \ \ \ \ \ \ \ \ \ \  + \ C \cdot \Vert\,\Phi\,\Vert_{\,*_Y} \cdot \int_{\R^n}  \bigg\vert \, [\,c_n \cdot{\bf K}\,(\,Z\,)\,] \ - \ n\,(n\,-2)\,\bigg\vert   \cdot {1\over {(\,1 \,+\, \Vert\, Z \,-\,\Xi_{\,\,1}\,\Vert\,)^{\,n\,+\, {{n\,+\,2}\over 2 } \,+\ \tau_{\,>1} } }}
 \  d\,Z\ .
\end{eqnarray*}

\vspace*{0.5in}

{\bf \large {\bf  \S\,A\,3.\,l\,.}}  }  { \bf Estimate on the integral of the form} \,\{\,distant assumption -- refer to (\,1.26\,) of  the main text\,\}\,.\, This part is refined later in {\bf \S\,A\,9}\,.\\[0.1in]
(\,A.3.38\,)
\begin{eqnarray*}
& \ & \int_{\R^n}  \bigg\vert \, [\,c_n \cdot{\bf K}\,(\,Y\,)\,] \ - \ n\,(n\,-2)\,\bigg\vert   \cdot {1\over {(\,1 \,+\, \Vert\, Z \,-\,\Xi_{\,\,1}\,\Vert^{\,2}\,)^{\,q} }}\ d\,Z \ \ \ \ \ \left(\,q \ > \ {n\over 2}\,\ \right)\\[0.2in]
& = & C \cdot  {{\lambda_{\,1}^q}\over {\lambda_{\,1}^n}} \cdot \int_{\R^n}  \bigg\vert \, [\,c_n \cdot K\,(z)\,] \ - \ n\,(n\,-2)\,\bigg\vert   \cdot \left(\  {{\lambda_{\,1}}\over {\lambda_{\,1}^2  \,+\, \Vert\, z \,-\,\xi_{\,1}\,\Vert^{\,2}\,  }}\ \right)^{\!\!q} \ dz \ \ \  \ \ \ \  \left(\ z \ = \ \lambda \cdot Z\ \right)\ .
\end{eqnarray*}
Let $\,{\bf p_1} \, \in \, {\cal H}\,$ be such that
$$
 \Vert\ \xi_{\,1} \ - \ {\bf p}_1\,\Vert \ = \ {\mbox{dist}}\,(\,\xi_{\,1}\,,\,{\cal H}\,) \ = \ O\, (\,{\bar\lambda}_{\,\,\flat}^{1\,+\,\kappa}\,)\ . \leqno (\,A.3.39\,)
$$
As the integral is invariant under   translations and rotations, we may assume that
$$
{\bf p}_1 \ = \ 0 \ \ \ \ {\mbox{and \   the \   hyperplane   \ defined \   by }} \ \ \{ \,z_n \ = \ 0\, \} \ \ {\mbox{is \   the \   tangant \    space \   of \ \ {\cal H} \ \ at \ \ }} 0\,.
$$
For a fixed
$$
r \ ( \ > \  \Vert\,\xi_{\,1} \ - \ {\bf p}_1\,\Vert\ )  \ \ {\mbox{small \   enough\   so \   that \ \ }} {\mbox{dist}}\,(\,z\,,\,{\cal H}\,)  \ \le \ C \cdot |\,z_n\,| \ \ \ \ {\mbox{for}} \ \ \ \ z \ = \ (z_1\,,\, \cdot \cdot \cdot\,, \ z_n) \,\in\, B_o\,(r)\,,
$$
we have
\begin{eqnarray*}
(\,A.3.40\,) \ \ & \ & {{\lambda_{\,1}^q}\over {\lambda_{\,1}^n}} \cdot \int_{\R^n}  \bigg\vert \,[\,c_n \cdot{\bf K}\,(\,Y\,)\,] \ - \ n\,(\,n\,-\,2\,)\,\bigg\vert   \cdot \left(\  {{\lambda_{\,1}}\over {\lambda_{\,1}^2  \,+\, |\, z \,-\,\xi_{\,1}\,|^{\,2}\,  }}\ \right)^{\!\!q} \ dz\\[0.2in]
& \le &   {{\lambda_{\,1}^q}\over {\lambda_{\,1}^n}} \cdot \int_{B_{o}\,(r)}  \bigg\vert \,[\,c_n \cdot{\bf K}\,(\,Y\,)\,] \ - \ n\,(\,n\,-\,2\,)\,\bigg\vert   \cdot \left(\  {{\lambda_{\,1}}\over {\lambda_{\,1}^2  \,+\, |\, z \,-\,\xi_{\,1}\,|^{\,2}\,  }}\ \right)^{\!\!q} \ dz \\[0.2in]
& \ & \ \ \ \ +  \ C_1 \cdot  {{\lambda_{\,1}^q}\over {\lambda_{\,1}^n}} \cdot \int_{\R^n \setminus \,B_o(r)}    \left(\  {{\lambda_{\,1}}\over {\lambda_{\,1}^2  \,+\, |\, z \,-\,\xi_{\,1}\,|^{\,2}\,  }}\ \right)^{\!\!q} \ dz \ \ \ \ \left[\ = \ O\,\left(\ {\bar\lambda}_{\,\,\flat}^{2q \,-\,n}\ \ \right) \ \right]\\[0.2in]
& = &  C_2 \cdot  {{\lambda_{\,1}^q}\over {\lambda_{\,1}^n}} \cdot \left\{\int_{B_{{\bf p}_1}\,(r)} |\,z_n\,|^{\,\ell} \left(\  {{\lambda_{\,1}}\over {\lambda_{\,1}^2  \,+\, |\, z \,-\,\xi_{\,1}\,|^{\,2}\,  }}\ \right)^{\!\!q } \ \right\} \, \cdot [\,1\,+\,o_{\,+}\,(\,1\,)\,]  \ \ \   + \ O\,\left(\ {\bar\lambda}_{\,\,\flat}^{2q \ -\,n}\ \ \right)  \\[0.1in]
& \ & \hspace*{2.7in} [\ {\mbox{applying \ \ (\,1.26\,) \ \ in \ \ the \ \ main \ \ text}} \ ]\\[0.1in]
& \le &    C_3  \cdot {\bar\lambda}_{\,\,\flat}^\ell \cdot \int_{\R^n} |\,Z_n|^{\,\ell} \cdot \left(\ {1\over {1 \ + \ |\,Z \ - \ \Xi_{\,\,1}\,|^{\,2}}} \ \right)^{\!\!q} d\,Z   \ \ \ + \ O\,\left(\ {\bar\lambda}_{\,\,\flat}^{2q \ -\,n}\ \ \right)  \\[0.2in]
& = & O\,(\ {\bar\lambda}_{\,\,\flat}^{\ell}\ )   \ + \ O\,\left(\ {\bar\lambda}_{\,\,\flat}^{2q \ -\,n} \ \right)\\[0.1in]
& \ &  \ \ \   \  \ \ \ \ \ \ \ \ \ \ \ \ \ \ \  \left(\ 2\,q \ =  \ n \ + \ {{n\,+\,2}\over 2} \,+\ \tau_{\,>1} \ \ \Longrightarrow \ \ 2\,q \ - \ n \ = \ {{n\,+\,2}\over 2} \,+\ \tau_{\,>1}\  \ \right)\\[0.2in]
& = & O\,(\ {\bar\lambda}_{\,\,\flat}^{\,\ell}\ )   \ + \ O\,\left(\ {\bar\lambda}_{\,\,\flat}^{ \,{{n\,+\,2}\over 2} \ +\,\tau_{\,>1}   }\ \ \right) \ .\\
\end{eqnarray*}

\newpage

{\bf \large {\bf  \S\,A\,3.\,l\,'\,.    }}    \ \ { \bf Estimate on the integral of the form} [\,uniform assumption (\,A.3.21\,)\,]\,.\\[0.1in]
We begin with
\begin{eqnarray*}
(\,A.3.41\,)  & \ & \int_{\R^n}  \bigg\vert \  [\,c_n \cdot{\bf K}\,(\,Z\,)\,] \ - \ n\,(n\,-2)\,\bigg\vert   \cdot {1\over {(\,1 \,+\, \Vert\, Z \,-\,\Xi_{\,\,1}\,\Vert^{\,2}\,)^{\,q} }}\ d\,Z \\[0.1in]
& \ & \hspace*{4in} \ \ \ \ \ \ \ \ \ \ \ \ \ \  \ \left(\,q \ > \ {n\over 2}\,\ \right)\\[0.1in]
& = & C \cdot  {{\lambda_{\,1}^q}\over {\lambda_{\,1}^n}} \cdot \int_{\R^n}  \bigg\vert \, [\,c_n \cdot K\,(z)\,] \ - \ n\,(n\,-2)\,\bigg\vert   \cdot \left(\  {{\lambda_{\,1}}\over {\lambda_{\,1}^2  \,+\, \Vert\, z \,-\,\xi_{\,1}\,\Vert^{\,2}\,  }}\ \right)^{\!\!q} \ dz\\[0.1in]
& \ & \hspace*{4in} \ \ \ \ \ \ \ \ \ \ \ \ \ \  \  \left(\, z \ = \ \lambda \cdot Z\ \ \right)\ .
\end{eqnarray*}
Recall the assumption in (\,A.3.21\,)\,:
$$
|\,(\,c_n \cdot K\,)\,(\,y\,)\,  \ - \ n\,(\,n\,-\,2\,)\,|\  \le \ {\bar\lambda}_{\,\,\flat}^{\,\zeta} \ \ \  \ \mfor \ \ y \,\in\,B_{\xi_{\,l}}\,(\,\rho_\nu\,)\ \ \ \ \ \ \ \ \ \  (\ l \ = \ 1\,, \ \cdot \cdot \cdot\,, \ \flat\,)\,, \leqno (\,A.3.42\,)
$$
where
$$
\rho_{\,\nu}\ = \ O\,\left(\,{\bar\lambda}_{\,\,\flat}^{ \,\nu\,}\ \right)\ \ \ \ {\mbox{and}} \ \ \ \ \Vert \,\xi_{\,l}\,\Vert \ = \ O\,\left(\,{\bar\lambda}_{\,\,\flat}^{ \,1\,+\,\kappa\, } \ \ \right)\ ,
$$
$$
1 \ >  \ \gamma \ >  \ \nu\ >  \ {1\over 2} \cdot (\,1\ + \ \kappa\,) \ .\leqno{\mbox{and}}
$$
Refer to Remark 1.23  regarding  the leak of a constant in the right hand side of (\,A.3.42\,)\,.\,  We are led to
\begin{eqnarray*}
(\,A.3.43\,) \ \ \ & \ & {{\lambda_{\,1}^q}\over {\lambda_{\,1}^n}} \cdot \int_{\R^n}  \bigg\vert \,[\,c_n \cdot{\bf K}\,(z)\,] \ - \ n\,(\,n\,-\,2\,)\,\bigg\vert   \cdot \left(\  {{\lambda_{\,1}}\over {\lambda_{\,1}^2  \,+\, \Vert\, z \,-\,\xi_{\,1}\,\Vert^{\,2}\,  }}\ \right)^{\!q} \ dz\\[0.2in]
& \le &   {{\lambda_{\,1}^q}\over {\lambda_{\,1}^n}} \cdot \int_{B_{\xi_{\,l}}\,(\,\rho_\nu\,)}\   \bigg\vert \,[\,c_n \cdot{\bf K}\,(z)\,] \ - \ n\,(\,n\,-\,2\,)\,\bigg\vert   \cdot \left(\  {{\lambda_{\,1}}\over {\lambda_{\,1}^2  \,+\, \Vert\, z \,-\,\xi_{\,1}\,\Vert^{\,2}\,  }}\ \right)^{\!\!q} \ dz \\[0.2in]
& \ &\ \ \ \ \ \ \ \  \ \ \ \ \ \ \ \   \ \ \ \ \   \ \ \ \ \ \ \ \ \ +  \ C_1 \cdot  {{\lambda_{\,1}^q}\over {\lambda_{\,1}^n}} \cdot \int_{\R^n\, \setminus \,B_{\xi_{\,l}}\,(\,\rho_\nu\,)}    \left(\  {{\lambda_{\,1}}\over {\lambda_{\,1}^2  \,+\, \Vert\, z \,-\,\xi_{\,1}\,\Vert^{\,2}\,  }}\ \right)^{\!\!q} \ dz \\[0.2in]
& \le  &  C_2 \cdot  {{\lambda_{\,1}^{q\,+\,\zeta}}\over {\lambda_{\,1}^n}} \cdot \int_{B_{\xi_{\,l}}\,(\,\rho_\nu\,)} \left(\  {{\lambda_{\,1}}\over {\lambda_{\,1}^2  \,+\, \Vert\, z \,-\,\xi_{\,1}\,\Vert^{\,2}\,  }}\ \right)^{\!\!q } \ dz \ \   + \ O \left(\ \ {\bar\lambda}_{\,\,\flat}^{\,(\,2\,q \ -\,n\,)\,\cdot\, (\,1\,-\ \nu\,) } \ \ \right) \\[0.2in]
& = & O\,(\,{\bar\lambda}_{\,\,\flat}^{\,\zeta}\,)   \ + \ O \left(\ \ {\bar\lambda}_{\,\,\flat}^{\,(\,2\,q \ -\,n\,)\,\cdot\, (\,1\,-\ \nu\,) } \ \ \right)\ .
\end{eqnarray*}
Note that in the present situation [\ refer to (\,A.3.37\,)\ ]\,,
$$
 2\,q \ =  \ n \ + \ {{n\,+\,2}\over 2} \,+\ \tau_{\,>1} \ \ \Longrightarrow \ \ 2\,q \ - \ n \ = \ {{n\,+\,2}\over 2} \,+\ \tau_{\,>1}  \ .
$$


\newpage

It follows that
\begin{eqnarray*}
(A.3.44) \ \ & \ & \bigg\vert \ \bigg\langle \left(\  {\tilde C}_{\,n}\cdot \left\{  \ [\,c_n\cdot {\bf K}\,(\,Y\,)\,] \ - \ n\,(\,n\,-\,2\,)\,\right\}  \cdot {\bf W}_\flat^{4\over {n\,-\,2}}  \cdot  \Phi\ \right)\,,  \ \ {{\partial {{{\bf V}}_1\,}  }\over {\partial\, \Lambda_{\,1}}}  \ \bigg\rangle_{\int}\ \bigg\vert \ \ \ \ \ \  \ \ \ \  \ \ \ \ \ \ \ \ \
\end{eqnarray*}
$$
\le\  \cases{ \ \ { \displaystyle{       \Vert\,\Phi\,\Vert_{\,*_Y} \cdot  \left[\  O\,\left(\, {\bar\lambda}_{\,\,\flat}^{ {\ 4 \,\cdot\, \gamma }}\, \ \right)  \    + \ O\,\left(\,{\bar\lambda}_{\,\,\flat}^{\,\ell}\,\ \right)   \ + \ O\,\left(\ {\bar\lambda}_{\,\,\flat}^{ \left({{n\,+\,2}\over 2} \,+\ \tau_{\,>1} \ \right) \,\cdot\, (\,1 \ - \ \nu\,)  }\ \ \right) \ \right]     }  }  \ \ & \ \ \ \ \ \ \ \ \ \ \ \ \ \ \ \ \ \ \ \ \ \ \ \ \ \ \ \ \ \ \ \ \ \ \ \ \ \  \ \ \ \cr
& \cr
\ & \hspace*{-3.5in}\ \ \ \ \ \{ \ $\uparrow$ \ \ {\mbox Under \ \ condition (\,1.26\,) \ \ of \ \ the\ \  main\ \  text\ \  [\,8\,]\ .}\ \}\cr
\ & \cr
\ \   \displaystyle{       \Vert\,\Phi\,\Vert_{\,*_Y} \cdot  \left[\  O\,\left(\, {\bar\lambda}_{\,\,\flat}^{ {\ 4 \,\cdot\, \gamma }}\, \ \right)  \    + \ O\,\left(\,{\bar\lambda}_{\,\,\flat}^{\, \zeta}\,\ \right)   \ + \ O\,\left(\ {\bar\lambda}_{\,\,\flat}^{ \left({{n\,+\,2}\over 2} \,+\ \tau_{\,>1} \ \right) \,\cdot\, (\,1 \ - \ \nu\,)  }\ \ \right) \ \right]  \ .    }   &  \cr}
$$
 \ \ \ \ \ \ \ \ \ \ \ \ \ \ \ \ \ \ \ \ \ \ \ \ \ \ \ \ \ \ \ \ \ \ \ \ \ \ \ \ \ \ \ \ \ \ \ \ \ \  [\  $\uparrow$ \ \   Under condition (\,A.3.4\,)\,. \ ]

\vspace*{0.35in}

{\bf \large {\bf  \S\,A\,3.\,m\,.    }}    \ \ From (A.2.1)
\begin{eqnarray*}
\bigg\langle  {{{\bf V}}_1\,}^{4\over {n\,-\,2}} \cdot \left(\ \Lambda_{\,1} \cdot {{\partial\,{{{\bf V}}_1\,}  }\over {\partial\, \Lambda_{\,1}}}\ \right),\ \  \left(\ \Lambda_{\,1} \cdot {{\partial\,{{{\bf V}}_1\,}  }\over {\partial \, \Lambda_{\,1}}} \ \right) \ \bigg\rangle_{\int} &  = & C\,(n)\,, \\[0.2in]
& \ & \hspace*{-4.3in}\!\!\!{\mbox{and \ \ via \ \ symmetric \ \ cancellation \ \ [\,cf. (\,5.3\,)\,]}} \\[0.2in] \ \ \ \  \ \ \ \ \bigg\langle  {{{\bf V}}_1\,}^{4\over {n\,-\,2}} \cdot \left(\ \Lambda_{\,1} \cdot {{\partial\,{{{\bf V}}_1\,}  }\over  {\partial\,\Xi_{{\,1}_{\,|_{\,j}}} }} \ \right),\ \ \left(\ \Lambda_{\,1} \cdot {{\partial\,{{{\bf V}}_1\,}  }\over {\partial\, \Lambda_{\,1}}}\ \right) \ \bigg\rangle_{\int}  & = & 0\,. \ \ \ \ \ \ \ \ \ \ \ \ \ \ \ \  \ \ \ \ \ \ \ \ \ \ \ \ \ \ \ \
\end{eqnarray*}
Here $\,j \ = \ 1\,,\, \ 2\,,\,\ \cdot \cdot \cdot\,, \ n\,.\,$ Whereas from {\bf (b)}\,,\, {\bf \S\,A\,3.\,i}  (\,$l\,\not=\,1$\,)\,,
\begin{eqnarray*}
\bigg\langle  {{{\bf V}}_{\,l}\,}^{4\over {n\,-\,2}} \cdot \left(\ \Lambda_{\,l} \cdot {{\partial \,{{{\bf V}}_{\,l}\,}  }\over {\partial\, \Lambda_{\,l}}}\ \right),\ \  \left(\Lambda_{\,1} \cdot {{\partial\,{{{\bf V}}_1\,}  }\over {\partial \, \Lambda_{\,1}}}\ \right) \ \bigg\rangle_{\int}  \ = \ O\,\left( \ {1\over {\Vert\,\Xi_{\,\,1}\,-\,\Xi_{{\,l}} \,\Vert^{\,n\,-\,2} }} \ \right) \ = \ O\,\left({\bar\lambda}_{\,\,\flat}^{(\,n\,-\,2\,)\cdot \gamma}\ \right) \,.
\end{eqnarray*}
Recall that
$$
{\tilde{\bf P}}_{\flat_{\,\parallelsum}} \ = \ \sum_{l\,=\,1}^\flat   a_{\,l}  \cdot  {{\bf V}}_{\,l}^{4\over {n\,-\,2}} \cdot \left(\ \Lambda_{\,l} \cdot  {{\partial\,{{\bf V}}_{\,l} }\over {\partial\,\Lambda_{\,l} }}  \ \right) \ + \  \sum_{l\,=\,1}^\flat\,\sum_{j\,=\,1}^n b_{\,l\,, \ j} \cdot   {{\bf V}}_{\,l}^{4\over {n\,-\,2}} \cdot  \left(\ \Lambda_{\,l} \cdot  {{\partial\,{{\bf V}}_{\,l} }\over {\partial\,\Xi_{{l}_{\,|_{\,j}}} }} \ \right)\ ,
$$
Thus\\[0.1in]
(\,A.3.45\,)
$$
\bigg\langle \ {\tilde{\bf P}}_{\flat_{//}}\ ,\ \  \left(\ \Lambda_{\,1}\cdot {{\partial\, {{{\bf V}}_l\,}  }\over {\partial \, \Lambda_{\,1}}} \ \right) \ \bigg\rangle_{\int} \ =\  C\,(\,n\,) \cdot  a_{\,l}  \ + \   \sum_{l\,=\,2}^\flat\,\sum_{j\,=\,1}^n \ \left[\ b_{\,l\,, \ j}\cdot O\,\left(\lambda^{\,(\,n\,-\ 2\,)\,\cdot\, \gamma}\ \right) \ \right]  \ .
$$
Similar expression for
$$
\bigg\langle \ {\tilde{\bf P}}_{\flat_{//}}\ ,\ \  \left(\ \Lambda_{\,l} \cdot  {{\partial\,{{\bf V}}_{\,l} }\over {\partial\,\Xi_{\,{l}_{\,|_{\,{1_{\,j}}}}}}} \ \right) \ \bigg\rangle_{\int}  \ .
$$


\newpage

{\bf \large {\bf  \S\,A\,3.\,n\,.    }}    \ \ {\bf Matrix equation and estimates on the coefficients.} \ \ Matching with   equation (A.3.3)\,,\, we are led to \\[0.1in]
(A.3.46)
\begin{eqnarray*}
\bigg\langle \ {\tilde{\bf P}}_{\flat_{//}}\,,\   \left(\Lambda_{\,1} \cdot {{ \partial\, {{\bf V}}_1,\   }\over {\partial\, \Lambda_{\,1}}} \ \right) \ \bigg\rangle_{\int} &  = & \bigg\langle \left(\ \Delta\,\Phi \ + \ \left\{\,{\tilde C}_{\,n} \cdot [\,c_n \cdot {\tilde K}(\,Y\,)\,] \cdot {\bf W}_\flat^{4\over {n\,-\,2}}\  \right\} \cdot   \Phi\ \right), \ \left(\Lambda_{\,1} \cdot {{ \partial\, {{\bf V}}_1,\   }\over {\partial\, \Lambda_{\,1}}} \ \right) \ \bigg\rangle_{\int} \\[0.2in]
& \ & \ \ \ \ \ \ \ \ \  \ \ \ \ \ \ \ \ \ - \  \bigg\langle \, {\bf H}\,,\  \left(\Lambda_{\,1} \cdot {{ \partial\, {{\bf V}}_1,\   }\over {\partial\, \Lambda_{\,1}}} \ \right) \ \bigg\rangle_{\int} \ .
\end{eqnarray*}
Similar expression for
$$
\bigg\langle \ {\tilde{\bf P}}_{\flat_{//}}\ ,\ \  \left(\ \Lambda_{\,l} \cdot  {{\partial\,{{\bf V}}_{\,l} }\over {\partial\,\Xi_{\,{l}_{\,|_{\,{1_{\,j}}}}}}} \ \right) \ \bigg\rangle_{\int}  \ .
$$
Let
$$
 \vec{\,V}  \ = \ (\,a_1\,,\ b_{1\,,\ 1}\,,\,\cdot \cdot \cdot\,,\, b_{1\,, \ n}\ ; \    \ \cdot \cdot \cdot\,; \  a_{\,\flat}\,,\ b_{\,\flat\,,\ 1}\,,\,\cdot \cdot \cdot\,,\, b_{\,\flat\,, \ n}\,) \ \in  \ \R^{\,\flat\,\cdot\, (\,n\,+\,1\,)}\ .
$$
Consider the result in (\,A.3.36\,)\,,\, etc.\,,\, in the following  matrix equation\,:
$$
\tilde{\bf M}\, (\,\vec{\,c}\ ) \ = \  (\ \overrightarrow{{\mbox{R.H.S.}}}\,)\ . \leqno (\,A.3.47\,)
$$
Here $\,(\ \overrightarrow{{\mbox{R.H.S.}}}\,)\,$ is a $\,1 \,\times\,\{\ \flat\,\cdot\, (\,n\,+\,1\,)\,\}\,$ matrix with entries of the following order [\ from (\,A.3.28\,) \,and\, (\,A.3.44\,)\ ]\,:\\[0.1in]
(\,A.3.48\,)
\begin{eqnarray*}
& \ & O\,\left(\, \Vert\,{\bf H}\,\Vert_{**_Y} \,\ \right) \ +  \ \Vert\,\Phi\,\Vert_{\,*_Y} \cdot  \left[\  O\,\left(\, {\bar\lambda}_{\,\,\flat}^{ {\ 4\,\cdot\, \gamma }}\, \ \right)  \    + \ O\,\left({\bar\lambda}_{\,\,\flat}^{\,\ell}\,\ \right)   \ + \ O\,\left(\ {\bar\lambda}_{\,\,\flat}^{ \left({{n\,+\,2}\over 2} \,+\ \tau_{\,>1} \ \right) \,\cdot\, (\,1 \ - \ \nu\,)  }\ \ \right) \ \right]\ . \
\end{eqnarray*}
Note that $
 \Vert\,\Phi_{\,\flat_{\,i}}\,\Vert_{\,*_Y} \ = \ 1\,, \ \ {\bar\lambda}_{\,\,\flat_{\,i}} \ \to \ 0^{\,+} \ \ \ \ {\mbox{and}} \ \ \Vert\,{\bf H}_{\,\flat_{\,i}}\,\Vert_{**_Y}\ \to \ 0^{\,+}.\,$
Here  $\ o\,(\,1 \,)\ \to \ 0^{\,+}\  $ as $\ {\bar{\lambda}}_{\,\flat} \ \to \ 0^{\,+}\,$.\,
Whereas $\,\tilde{\bf M}\,$ is an ``\,almost\," diagonal matrix with   $\,O\,(\,1)\,$ diagonal entries\,,\,   $\,o_{\,+}\,(\,1\,)\,$ non\,-\,diagonal entries\,.\, In particular, the matrix is invertible\,.\, Cramer's rule shows that the inverse $\,{\tilde{\bf M}}^{\,-\,1}\,$ has similar form, that is, having $\,O\,(\,1\,)\,$ diagonal entries\,,\, and  $\,o_{\,+}\,(\,1\,)\,$ non\,-\,diagonal entries\,.\,
Applying the inverse  $\,{\tilde{\bf M}}^{\,-\,1}\,$  to both sides of equation (\,A.3.47\,)\,,\, we see that  each of the coefficient in $\,{\tilde{\bf P}}_{\flat_{//}}\,$ is  bounded from above by the following\,.

(\,A.3.49\,)
\begin{eqnarray*}
& \ & \max \ \left\{  \ |\,a_1\,|\,,\ |\,b_{1\,,\ 1}\,|\,,\,\cdot \cdot \cdot\,,\ |\,b_{1\,,\ n}\,|\ ; \    \ \cdot \cdot \cdot\,; \  |\,a_{\,\flat}\,|\,,\ |\,b_{\,\flat\,,\ 1}\,|\,,\,\cdot \cdot \cdot\,,\ |\,b_{\,\flat\,, \ n} \,|\ \right\}\\[0.2in]
& = &  O\,\left(\, \Vert\,{\bf H}\,\Vert_{**_Y} \,\ \right) \ +  \ 1 \cdot  \left[\  O\,\left(\, {\bar\lambda}_{\,\,\flat}^{ {\ 4 \,\cdot\, \gamma }}\, \ \right)  \    + \ O\,\left({\bar\lambda}_{\,\,\flat}^{\,\ell}\,\ \right)   \ + \ O\,\left(\ {\bar\lambda}_{\,\,\flat}^{\, \left({{n\,+\,2}\over 2} \,+\ \tau_{\,>1} \ \right) \,\cdot\, (\,1\ - \ \nu\,)  }\ \ \right) \ \right] \\[0.2in]
& = &  \ o_{\,{\bar\lambda}_{\,\,\flat}}(\,1\,) \ \ \ (\, \to \ 0^{\,+} \ \ \ \  {\mbox{as}} \  \ {\bar\lambda}_{\,\,\flat} \ \to \ 0^+\ )\ .\,
\end{eqnarray*}

It follows that [\,recall from (\,A.3.16\,) that $\,\theta\ >\,0$ is small\ ] \\[0.1in]
(A.3.50)
\begin{eqnarray*}
|\ \Phi_{\,\flat_{\,i}}\,(\,Y\,)\,| & \le &  C \cdot  \left(\ \ \sum_{l\ =\,1}^\flat \left(\ {1\over {      1\ + \ \Vert\,Y \,-\,\Xi_{\,l}\,\Vert   }} \ \right)^{\!\!{ {{n\,-\,2}\over 2 } \,+\ \tau_{\,>1} \,+\ \theta}} \ \ \ \right)   \\[0.15in]
& \ & \ \ \ \ \ \ \ \ \ \ \ \ \ + \ \, C' \cdot \Vert \,{\bf H}_{\,\flat_{\,i}}\, \Vert_{**_Y}  \cdot \left(\ \  \displaystyle{\sum_{l\,=\,1}^\flat }  \left(\ {1\over {    1\ + \ \Vert\,Y \,-\ \Xi_{\,l}\,\Vert  }}\ \right)^{\!\! { {{n\,-\,2}\over 2 } \,+\ \tau_{\,>1} }} \ \  \ \right)\\[0.15in]
& \ &  \ \ \ \ \ \ \ \ \ \ \ \ \ \ \ \ \ \ \ \ \ \ \ \ \ \ + \ o\,(\,1 \,)  \cdot  \left(\ \  \displaystyle{\sum_{l\,=\,1}^\flat }  \left(\ {1\over {    1\ + \ \Vert\,Y \,-\ \Xi_{\,l}\,\Vert  }}\ \right)^{\!\! { {{n\,-\,2}\over 2 } \,+\ \tau_{\,>1} }} \ \  \ \right)  \\[0.25in]
&\ &\hspace*{-1.2in} \Longrightarrow  \ \ \  {{ |\,\Phi_{\,\flat_{\,i}}\,(\,Y\,)\,|  }\over {    \left(\  \  \displaystyle{\sum_{l\,=\,1}^\flat }  \left(\ {1\over {    1\ + \ \Vert\,Y \,-\ \Xi_{\,l}\,\Vert  }}\ \right)^{\!\! { {{n\,-\,2}\over 2 } \,+\ \tau_{\,>1} }} \ \  \ \right) }}\\[0.2in]
& \le &   \overline{C} \ \cdot {{ \left(\ \  \displaystyle{\sum_{l\,=\,1}^\flat }  \left(\ {1\over {    1\ + \ \Vert\,Y \,-\ \Xi_{\,l}\,\Vert  }}\ \right)^{\!\! { {{n\,-\,2}\over 2 } \,+\ \tau_{\,>1}\,+\ \theta }} \ \  \ \right)  }\over { \left(\ \  \displaystyle{\sum_{l\,=\,1}^\flat }  \left(\ {1\over {    1\ + \ \Vert\,Y \,-\ \Xi_{\,l}\,\Vert  }}\ \right)^{\!\! { {{n\,-\,2}\over 2 } \,+\ \tau_{\,>1} }} \, \ \right)   }}\ \ + \ o_{\,+}\,(\,1\,)\  \ \ \ \ \ \  {\mbox{for}} \ \ Y\,\in\,\R^n\,.\\
\end{eqnarray*}
Here  $\ o\,(\,1 \,)\ \to \ 0\, $  as $\,{\bar{\lambda}}_{\,\flat_{\,i}} \ \to \ 0^{\,+}\,$ (\, that is\,,\, as $\,i\,\to\,\infty\,$)\,.\, Cf. (\,2.12\,) in \cite{Wei-Yan}\,.\, Due to the ``\,advantage\," of $\,\theta\,$ in the first term in the left hand side of the above, one can make the term to become small by controlling the distance $\,\Vert\,Y\,-\,\Xi_{\,l}\,\Vert\,.$\, Refer to $\,{\bf (\,i\,)}_{(\,A.3.55\,)}$\,.

\vspace*{0.35in}

{\bf \large {\bf  \S\,A\,3.\,o \,.    }}    \ \ Following \cite{Wei-Yan}\,,\,
  we seek to establish the following.\\[0.1in]
${\bf (\,*\,)_{\,(\,A.3.51\,)}}$ \ \  There exist (\,a fixed large number\,) $\,R_{\,o} \ >  \ 0\,$,\,  and (\,a fixed small number\,) $\,\varepsilon_o \ >  \ 0\,$,\,  both independent of $\, i \ \gg \ 1\,,\,$ such that for each $\,i\,\gg\,1\,,\,$  there is at least one integer $\,I_{\,i}\,$ with
$$
1 \ \le  \ I_{\,i} \ \le \ \flat_{\,i}\ , \leqno (\,A.3.51\,)
$$
(\,that integer $\,I_{\,i}\,$ may change as $\,i \ \to \ \infty\,$)\,,\, fulfilling
$$
\sup\ \bigg\{ \ |\,\Phi_{\,{\flat_{\,i}}_{\,i}}\,(\,Y\,)\ | \ \ \  \bigg\vert \ \ \ \  Y\, \in \,  B_{\,\Xi_{\,I_{\,i}}} (\,R_{\,o})   \ \bigg\} \ \,\ge \ \varepsilon_o \ > \  0\,. \leqno (\,A.3.52\,)
$$

\newpage
We demonstrate this by contradiction. Suppose that (A.3.51) does not hold. We first fix $\,R_{\,o}\,$ to be large enough so that
$$
  {{\overline{C}}\over { R_{\,o}^{\,\theta}}}  \ <  \ {1\over 2}   \ ,\leqno (\,A.3.53\,)
$$
where $\,\overline{C}\,$ and $\,\theta\,$ appear in (A.3.50). Then we choose $\,\varepsilon_o\,> \,0\,$ to be small enough so  that
$$
\varepsilon_o \cdot (\,1\,+\, R_{\,o}\,)^{ {{n\,-\,2}\over 2 } \ +\ \tau_{\,>1} }   \ < \   {1\over 2} \ . \leqno (\,A.3.54\,)
$$

Suppose that (\,for contradiction sake, and for simplicity sake, we rename the subsequence as if it is the original sequence\,) for  all $\,i\  \gg \ \ 1$\, we have
$$
\sup\ \bigg\{ \ |\,\Phi_{\flat_{\,i}}\,(\,Y\,)\ | \ \ \ \bigg\vert \ \  \ Y\, \in  \, B_{\,\Xi_{\,l }} (\,R_{\,o})  \  \bigg\} \ < \ \varepsilon_o \ \ \ \ \  \mfor (\,{\mbox{all}}\,) \ \ \  \ \ l \ = \ 1\,, \ \cdot \cdot \cdot \,, \ \flat_{\,i} \ . \leqno (\,A.3.55\,)
$$
For $\,i \ \gg \ 1\,,\,$ $\,{\bar\lambda}_{\,\,\flat_{\,i}}\,$ is small\,,\, and $\,\Vert\,\Xi_{\,k} \ - \ \Xi_{\,l}\,\Vert\,\gg\,1$ for $\,k\,\not=\,l\,,\,$ leading to
$$
 B_{\,\Xi_{\,k}}  (\,R_{\,o}\,) \ \cap \ B_{\,\Xi_{\,l}}  (\,R_{\,o}\,) \ = \ \emptyset\ \ \ \ \ \ \  \ \ \ \ \  (\, 1 \ \le \ l \ \not= \ k \ \le \ \flat_{\,i}\,)\ .
$$

 \vspace*{0.1in}

 $\,{\bf (\,i\,)}$ \ \ Outside the balls\,:
$$
Y \ \in \ \R^n  \bigg\backslash\, \left[\  {\displaystyle{\bigcup_{l \,=\, 1}^{\flat_{\,i}} B_{\,\Xi_{\,l}}  \,(R_{\,o})}}\,\right] \ \ \ \ \bigg(  \Longrightarrow \ \ \Vert\, Y\,-\,\Xi_{\,l}\,\Vert \ \ge \ R_{\,o} \ \ \ \ {\mbox{for \ \ all}} \ \ 1 \ \le \ l\  \le \ \flat_{\,i}\ \bigg) \,,\,
$$
we have\\[0.1in]
(\,A.3.56\,)
\begin{eqnarray*}
& \ & {1\over {   (  1\ + \ \Vert\,Y \,-\,\Xi_{\,l}\,\Vert\,)^{ {{n\,-\,2}\over 2 } \,+\ \tau_{\,>1} \,+\ \theta} }}\\[0.2in] & = & {1\over {   (  1\ + \ \Vert\,Y \,-\,\Xi_{\,l}\,\Vert\,)^{ {{n\,-\,2}\over 2 } \,+\ \tau_{\,>1}  } }} \cdot {1\over {   (  1\ + \ \Vert\,Y \,-\,\Xi_{\,l}\,\Vert\,)^{  \,\theta} }}\\[0.2in]
& \le & {1\over {   (  1\ + \ \Vert\,Y \,-\,\Xi_{\,l}\,\Vert\,)^{ {{n\,-\,2}\over 2 } \,+\ \tau_{\,>1}  } }} \cdot {1\over {   (  1\ + \ R_{\,o}\,)^{  \,\theta} }}\
 \le\ {1\over {   (  1\ + \ \Vert\,Y \,-\,\Xi_{\,l}\,\Vert\,)^{ {{n\,-\,2}\over 2 } \,+\ \tau_{\,>1}  } }} \cdot {1\over {  R_{\,o}^{  \,\theta} }}\\[0.2in]
\Longrightarrow & \ & {{\  \left(\  {\displaystyle{\sum_{l\,=\,1}^{\flat_{\,i}} \,  {1\over {   (  1\ + \ \Vert\,Y \,-\,\Xi_{\,l}\,\Vert\,)^{ {{n\,-\,2}\over 2 } \,+\ \tau_{\,>1} \,+\ \theta} }} }} \ \
\ \right)\  }\over { \left(\ \ {\displaystyle{\sum_{l\,=\,1}^{\flat_{\,i}} \, {1\over {   (  1\ + \ \Vert\,Y \,-\,\Xi_{\,l}\,\Vert )^{ {{n\,-\,2}\over 2 } \,+\ \tau_{\,>1}  } }} }}\ \
\ \right) }}  \ \le \   {1\over { R_{\,o}^{\,\theta}}} \\[0.1in]
& \ & \hspace*{2.5in}\ \ \ \ \ \ \  \mfor \ \ Y \ \in \ \R^n  \bigg\backslash\, \left[\  {\displaystyle{\bigcup_{l \,=\, 1}^{\flat_{\,i}} B_{\,\Xi_{\,l}}  \,(R_{\,o})}}\,\right]\ .
\end{eqnarray*}
It follows from (\,A.3.50\,) and (\,A.3.53\,) that
$$
{{ |\,\Phi_{\,\flat_{\,i}}\,(\,Y\,)\,|  }\over {    \left(\  \displaystyle{\sum_{l\,=\,1}^{\flat_{\,i}} }  \left(\ {1\over {    1\ + \ \Vert\,Y \,-\ \Xi_{\,l}\,\Vert  }}\ \ \right)^{\!\! { {{n\,-\,2}\over 2 } \,+\ \tau_{\,>1} }} \  \ \right) }}  \ \le \ \overline{C}  \cdot\, {1\over { R_{\,o}^{\,\theta}}}  \ \ + \ \ \overline{o_{\,+}\,(\,1\,)} \ <   \ {1\over 2}\  \leqno (A.3.57)
$$
for  $\,i\,\gg\,1\,$ and $\displaystyle{\ Y \ \in \ \R^n  \bigg\backslash\, \left[\  {\displaystyle{\bigcup_{l \,=\, 1}^{\flat_{\,i}} B_{\,\Xi_{\,l}}  \,(R_{\,o})}}\,\right]}$\ .


\vspace*{0.2in}
 $\,{\bf (\,ii\,)}$ \ \ For what is inside, say, \\[0.1in]
(\,A.3.58\,)
\begin{eqnarray*}
 \ \ \ Y \, \in \,    B_{\,\Xi_{\,1 }}  (R_{\,o})  \ \ \ \ \bigg[ \Longrightarrow \ \   \Vert\, Y\,-\,\Xi_{\,l}\,\Vert \ < \ R_{\,o} \ \ \ {\mbox{and \ \ hence}} \ \ \  |\,\Phi_{\,\flat_{\,i}}\,(\,Y\,)\,| \ \le \ \varepsilon_o & \ & \\[0.1in]
 & \ & \hspace*{-1in} {\mbox{based \ \ on}} \ \ (\,A.3.55\,)\ \bigg]\ {\bf ,} \\[-0.2in]
\end{eqnarray*}
we have\\[0.1in]
(\,A.3.59\,)
\begin{eqnarray*}
& \ & \sum_{l\,=\,1}^{\flat_{\,i}}    \left(\ {1\over {    1\ + \ \Vert\,Y \,-\ \Xi_{\,l}\,\Vert  }}\ \ \right)^{\!\! { {{n\,-\,2}\over 2 } \ +\ \tau_{\,>1} }}  \\[0.2in]
&  = &  \left(\ {1\over {    1\ + \ \Vert\,Y \,-\ \Xi_{\,\,1}\,\Vert  }}\ \ \right)^{\!\! { {{n\,-\,2}\over 2 } \ +\ \tau_{\,>1} }}  + \ \  \sum_{l\,=\,2}^{\flat_{\,i}}    \left(\ {1\over {    1\ + \ \Vert\,Y \,-\ \Xi_{\,l}\,\Vert  }}\ \ \right)^{\!\! { {{n\,-\,2}\over 2 } \ +\ \tau_{\,>1} }}\\[0.2in]
& \ & \ \ \ \ \  \ \ \    \ \ge \ \left(\ {1\over {    1\ + \ R_{\,o}  }}\ \ \right)^{\!\! { {{n\,-\,2}\over 2 } \ +\ \tau_{\,>1} }} \\[0.2in]
& \Longrightarrow & {{ |\,\Phi_{\,\flat_{\,i}}\,(\,Y\,)\,|  }\over {   \  \left[\  \displaystyle{\sum_{l\,=\,1}^{\flat_{\,i}} }  \left(\ {1\over {    1\ + \ \Vert\,Y \,-\ \Xi_{\,l}\,\Vert  }}\ \right)^{\!\! { {{n\,-\,2}\over 2 } \ +\ \tau_{\,>1} }} \, \right]\  }} \ \le \  \varepsilon_o\cdot (\,1\,+\, R_{\,o}\,)^{ {{n\,-\,2}\over 2 } \ +\ \tau_{\,>1} }  \ < \   {1\over 2} \,.\\
\end{eqnarray*}
Here we use (\,A.3.50\,) one more time\,,\, also (\,A.3.54\,)\,.\,
Together with  $\,{\bf (\,i\,)}_{\,(\,A.3.56\,)}$\,,\, we obtain  $$\,\Vert\,\Phi_{\,\flat_{\,i}}\,\Vert_{\,*_{\,Y}} \ \le  \ {1\over 2} \ \ \ \ \  \mfor  \ \ i \ \gg \ 1\ .\,$$ But this contradicts $\,\Vert\,\Phi_{\,\flat_{\,i}}\,\Vert_{\,*_{\,Y}} \ = \ 1\,$ for all $\,i\,\gg\,1\,.\,$  Hence the claim in  $\,{\bf (\,*\,)_{\,(\,A.5.51\,)}}\,$ stands.

\newpage

{\bf \large {\bf  \S\,A\,3.\,p \,.}}    \ \ {\bf Bound on} $\,|\,\Phi_{\,\flat_{\,i}}\,(\,Y\,)\,|\,.$ \ \ We proceed from  the normalization (\,A.3.23\,)\,,\, that is\,,\,
$$
\,\Vert\,\Phi_{\,\flat_{\,i}}\,\Vert_{\,*_{\,Y}} \ =   \ 1 \ \ \Longleftrightarrow \ \ \sup_{Y\,\in\,\R^n}\ \left\{ \ {{ |\,\Phi_{\,\flat_{\,i}}\,(\,Y\,)\,|  }\over {   \  \left[\  \displaystyle{\sum_{l\,=\,1}^{\flat_{\,i}} }\   \left(\ {1\over {    1\ + \ \Vert\,Y \,-\ \Xi_{\,l}\,\Vert  }}\ \right)^{\!\! { {{n\,-\,2}\over 2 } \ +\ \tau_{\,>1} }} \, \right]\  }} \ \right\} \ = \ 1\ .
$$
Consider the expression
 $$
 \sum_{l\,=\,1}^{\flat_{\,i}} \  \left(\ {1\over {    1\ + \ \Vert\,Y \,-\ \Xi_{\,l}\,\Vert  }}\ \right)^{\!\! { {{n\,-\,2}\over 2 } \ +\ \tau_{\,>1} } }\ .
 $$
 For $\,Y \ = \ \Xi_{\,\,1}\,,\,$
 \begin{eqnarray*}
 (\,A.3.60\,) \ \ \ \   & \ &  \left[ \ \,\sum_{l\,=\,1}^{\flat_{\,i}}  \left(\ {1\over {    1\ + \ \Vert\,Y \,-\ \Xi_{\,l}\,\Vert  }}\ \right)^{\!\! { {{n\,-\,2}\over 2 } \ +\ \tau_{\,>1} } } \ \right]_{\ Y \ = \ \Xi_{\,\,1} }\\[0.2in]
& = & 1  \ \ + \  \  \displaystyle{\sum_{l\,=\,2}^{\flat_{\,i}} }  \left(\ {1\over {    1\ + \ \Vert\,\Xi_{\,\,1} \,-\ \Xi_{\,l}\,\Vert  }}\ \ \right)^{\!\! { {{n\,-\,2}\over 2 } \,+\ \tau_{\,>1} }} \\[0.2in]
& \le  & 1 \ \ + \  \  C_1 \cdot \displaystyle{\sum_{l\,=\,2}^{\flat_{\,i}} }  \left(\ {1\over {     \Vert\,\Xi_{\,\,1} \,-\ \Xi_{\,l}\,\Vert  }}\ \ \right)^{\!\! { {{n\,-\,2}\over 2 } \,+\ \tau_{\,>1} }}\\[0.2in]
&   \le  &  1 \ + \  C_2  \cdot {\bar\lambda}_{\ {\flat_{\,i}}}^{ \left(\ { {{n\,-\,2}\over 2 } \,+\ \tau_{\,>1} } \ \right) \ \cdot \,\gamma} \ \ \ \ \ \ \ \ \ \   \\[0.2in]
& =  & 1 \ + \ o_{\,\lambda_{\flat_{\,i}}}\,(\,1\,) \ \ \ \ \ \ \ \  \ \ [ \ o_{\,\lambda_{\flat_{\,i}}}\,(\,1\,) \ \to\ 0^{\,+} \ \ \ \ {\mbox{as}} \ \ {\bar\lambda}_{\ {\flat_{\,i}}} \ \to \ 0^{\,+}\ , \ \ \ {\mbox{i.e.}}  \ \ i \ \to \ \infty  \ ]  \ \ \ \ \  \\[0.2in]
& \le  & 2  \ \ \ \ \ \ \ \ \ \ \ \ \ \ \  \ \ \ \ \ \ \mfor \ \ i \ \gg \ 1\,.\\
\end{eqnarray*}
{\it Minimization\,.\,} \ \ As for a general $\,Y\,\in\,\R^n\,,\,$ consider the (\,finite\,) collection of points (\,the centers of the bubbles\,)
$$
\Xi_{\,\,1}\,, \ \ \,\Xi_{2}\,, \    \,\cdot\,\cdot\,\cdot \ , \ \ \,\Xi_{\,l}\,, \ \,\cdot\,\cdot \cdot\,,\ \  \Xi_{\,\flat_{\,i}}\ .
$$
By finiteness\,,\, there is at least one point in the above collection\,,\, say $\,\Xi_{\,\,1}\,,\,$ such that
$$
\Vert\, Y \ - \ \Xi_{\,\,1}\,\Vert
$$
achieves the minimal value among the collection of non\,-\,negative numbers
$$
\{\ \Vert\, Y \ - \ \Xi_{\,\,1}\,\Vert\,, \ \, \Vert\, Y \ - \ \Xi_{\,2}\,\Vert\,,\ \cdot\,\cdot \cdot\,, \ \Vert\, Y \ - \ \Xi_{\,\flat}\,\Vert\ \}\ .
$$
That is\,,\,
$$
\Vert\, Y \ - \ \Xi_{\,\,1}\,\Vert \ \le \ \Vert\, Y \ - \ \Xi_{\,l}\,\Vert  \ \ \  \ \mfor \ \ l\ \not= \ 1\ .
$$
This leads to the decomposition
$$ \R^n \ = \ {\displaystyle{\bigcup_{\,j\,=\,1}^{\,k} }} \ \Upsilon_j\ , $$
where\\[0.1in]
(\,A.3.61\,)
$$
\Upsilon_1 \ :=\  \left\{ \, Y
 \ \in \ \R^n \ \ | \ \ \,\Vert\, Y \ - \ \Xi_{\,l}\,\Vert \ \ge \ \Vert\, Y \ - \ \Xi_{\,\,1}\,\Vert \ \ \ \ {\mbox{for}} \ \ \ \ l \ = \ 2\,,\, \ 3\,,\,\,\cdot\,\cdot \cdot\,, \ \flat\,\right\}\,, \ \ \ \ {\mbox{etc.}}
$$
(\,Note that $\,\Upsilon_{\,1}\,$ is not smooth\,.\, In case we desire integral estimates, as in {\bf \S\,A\,4\,.i\,},\, it is better off to replace $\,\Upsilon_{\,1}\,$ by a metric ball\,.\, Cf. as in {\bf \S\,A\,4\,.j}\,.  \,)\,.\,
Consider first the case that $\,Y \,\in \ \Upsilon_1 \,.\,$  We have
$$
Y \,\in \ \Upsilon_1 \ \ \Longrightarrow \ \ \Vert\, Y  \ - \ \Xi_{\,l} \,\Vert \ \ge \ {1\over 2} \cdot \Vert\,  \Xi_{\,\,1}  \ - \ \Xi_{\,l} \,\Vert \ \ \  \ \mfor \ \ l\ \not= \ 1\ . \leqno (\,A.3.62\,)
$$
For if this is not the case, we can obtain a contradiction via
\begin{eqnarray*}
\Vert\, \Xi_{\,\,1}   \ - \ \Xi_{\,l} \,\Vert & \le & \Vert\, \Xi_{\,\,1}   \ - Y \ + \ Y \ - \ \Xi_{\,l} \,\Vert \ \le \ \Vert\, Y \ - \  \Xi_{\,\,1} \,\Vert   \ + \ \Vert \,Y \ - \ \Xi_{\,l} \,\Vert \\[0.2in]
 & \le & 2\,\cdot \Vert \,Y \ - \ \Xi_{\,l} \,\Vert  \ \ \ \ \ \  \ \ \ \ \ \ \ ( \ {\mbox{as}} \ \ Y \,\in \ \Upsilon_1 \ )\\[0.2in]
 & < & 2\,\cdot {1\over 2} \cdot \Vert\,  \Xi_{\,\,1}  \ - \ \Xi_{\,l} \,\Vert \ \ \ \ \ \ \ \ (\ {\mbox{contradiction}} \ ) \ .
\end{eqnarray*}
As in  (\,A.3.61\,)\,,\, for $\,Y\,\in\,\Upsilon_{\,1}\,,\,$ we have
\begin{eqnarray*}
& \ & \displaystyle{\sum_{l\,=\,1}^{\flat_{\,i}} }  \left(\ {1\over {    1\ + \ \Vert\,Y \,-\ \Xi_{\,l}\,\Vert  }}\ \ \right)^{\!\! { {{n\,-\,2}\over 2 } \,+\ \tau_{\,>1} }} \hspace*{1.5in} ( \ \Xi_{\,l} \ = \ 0\ ) \\[0.2in]
& = & \left(\ {1\over {    1\ + \ \Vert\,Y \,-\ \Xi_{\,l}\,\Vert  }}\ \ \right)^{\!\! { {{n\,-\,2}\over 2 } \,+\ \tau_{\,>1} }}  \ \ + \  \  \displaystyle{\sum_{l\,=\,2}^{\flat_{\,i}} }  \left(\ {1\over {    1\ + \ \Vert\,Y \,-\ \Xi_{\,l}\,\Vert  }}\ \ \right)^{\!\! { {{n\,-\,2}\over 2 } \,+\ \tau_{\,>1} }} \\[0.2in]
& \le  & 1 \ \ + \  \  C_3 \cdot \displaystyle{\sum_{l\,=\,2}^{\flat_{\,i}} }  \left(\ {1\over {  \    \Vert\ \Xi_{\,\,1} \,-\ \Xi_{\,l}\ \Vert \  }}\ \ \right)^{\!\! { {{n\,-\,2}\over 2 } \,+\ \tau_{\,>1} }} \ \ \ \ \ \ \ \ \ \ \ \ \  [ \ {\mbox{via}} \ \ (\,A.3.62\,) \ ]\\[0.2in]
& \le  & 1 \ + \  C_4  \cdot {\bar\lambda}_{\ {\flat_{\,i}}}^{ \left(\ { {{n\,-\,2}\over 2 } \ +\ \tau_{\,>1} } \ \right) \ \cdot \,\gamma}  \\[0.2in]
& \le  & 1 \ + \ o_{\,\lambda_{\flat_{\,i}}}\,(\,1\,) \ \ \ \ \ \ \ \  \ \ [ \ o_{\,\lambda_{\flat_{\,i}}}\,(\,1\,) \ \to\ 0^{\,+} \ \ \ \ {\mbox{as}} \ \ {\bar\lambda}_{\ {\flat_{\,i}}} \ \to \ 0^{\,+} \ , \ \ {\mbox{i.e.}}  \ \ i \ \to \ \infty  \ ]   \\[0.2in]
& \le  & 2  \ \ \ \ \ \ \ \ \mfor \ \ i \ \gg \ 1 \ \  (\ {\mbox{and \ \ for}} \ \ \,Y\,\in\,\Upsilon_{\,1}\ )\ .
\end{eqnarray*}
Likewise for the other cases in the regions $\,\Upsilon_{\,2}\,,\, \ \cdot \cdot \cdot\,,\, \ \Upsilon_{\,\flat_{\,i}}\ .$\, Hence\\[0.1in]
(\,A.3.63\,)
 $$\,\Vert\,\Phi_{\,\flat_{\,i}}\,\Vert_{\,*_{\,Y}} \ =   \ 1 \ \ \Longrightarrow \  \ |\,\Phi_{\,\flat_{\,i}}\,(\,Y\,)\,| \ \le \ 2 \ \ \ \ \  \mfor  \ \ i \ \gg \ 1 \ \ \ {\mbox{and \ \ for \ \ all }} \ \ Y \,\in\,\R^n\ .\,$$

\vspace*{.25in}

{\bf \large {\bf  \S\,A\,3.\,q \,.}}    \ \ {\bf Proof of Lemma} \,A.3.20\,. \ \
Via the  translation
$$
Y \ \to \  {\bar Y}  \ := \ Y \ - \ \Xi_{\,l_{\,I_{\,i}}}  \ \ \Longleftrightarrow \ \ Y \ = \ {\bar Y}  \ + \ \Xi_{\,l_{\,I_{\,i}}} \,,\leqno (\,A.3.64\,)
$$
we set
$$
{\tilde\Phi_{\,I_{\,i}}} \,(\ {\bar Y}\,) \
 := \ \Phi_{\,I_{\,i}}\,(\ {\bar Y}\ + \ \Xi_{\,l_{\,I_{\,i}}} \,) \ \ \ \ \mfor \ \ {\bar Y} \,\in\,\R^n\ .
 $$
 Here $\,I_i\,$ is found in claim ${\bf(\,*\,)}_{\,(\,A.3.51\,)}\,.\, $
 The relevant $\,\perp\,$--\ \,conditions are expressed via\\[0.1in]
(\,A.3.65\,)
\begin{eqnarray*}
\bigg\langle \ V_{\Lambda_{I_{\,i}}\,,\ \Xi_{I_{\,i}}}^{4\over {n\,-\,2}} \cdot \left(\ \Lambda_{I_{\,i}} \cdot {{\partial\, V_{\Lambda_{I_{\,i}}\,,\ \Xi_{I_{\,i}}} }\over {\partial\, \Xi_{{I_{\,i}}_{|_m}} }}\right)\, \bigg\vert_{\,\Xi_{\,l_{\,I_{\,i}}}}\ , \ \ \  \Phi_{\,I_{\,i}} \  \bigg\rangle_{\!\int} &= & 0\\[0.2in]
\Longleftrightarrow \ \ \  \bigg\langle \ V_{\Lambda_{I_{\,i}}\,,\ 0}^{4\over {n\,-\,2}} \cdot \left(\ \Lambda_{I_{\,i}} \cdot {{\partial\, V_{\Lambda_{I_{\,i}}\,,\ \Xi} }\over {\partial\, \Xi_{\,|_{\,m}} }}\right)\, \bigg\vert_{\,\Xi \ = \ 0}\ , \ \ \ {\bar\Phi}_{\,I_{\,i}} \  \bigg\rangle_{\!\int} & = & 0 \ \ \ \ \mfor \ \ m \ = \ 1\,, \ \cdot \cdot \cdot\,, \ n\ ,\\
\end{eqnarray*}
and
\begin{eqnarray*}
\bigg\langle \ V_{\Lambda_{I_{\,i}}\,,\ \Xi_{I_{\,i}}}^{4\over {n\,-\,2}} \cdot \left(\ \Lambda_{I_{\,i}} \cdot {{\partial\, V_{\Lambda_{I_{\,i}}\,,\ \Xi_{I_{\,i}}} }\over {\partial\, \Lambda_{I_{\,i}} }}\right)\,\ , \ \ \  \Phi_{\,I_{\,i}} \  \bigg\rangle_{\!\int} &= & 0\\[0.2in]
\Longleftrightarrow \ \ \ \ \ \  \bigg\langle \ V_{\Lambda_{I_{\,i}}\,,\ 0}^{4\over {n\,-\,2}} \cdot \left(\ \Lambda_{I_{\,i}} \cdot {{\partial\, V_{\Lambda_{I_{\,i}}\,,\ 0} }\over {\partial\,  \Lambda_{I_{\,i}} }}\right)\ , \ \ \ {\bar\Phi}_{\,I_{\,i}} \  \bigg\rangle_{\!\int} & = & 0\ .
\end{eqnarray*}
Recall that
\begin{eqnarray*}
  {{\partial\, V_{\Lambda\,, \, \Xi}}\over {\partial\, \Lambda}} \ = \ {{\partial}\over {\partial\, \Lambda}} \left[\ \left({\Lambda\over {\Lambda^2 + \Vert\,\, Y \ - \ \Xi\,\Vert^{\,2}}}\  \right)^{{n - 2}\over 2} \ \right] & = &\!-\,{{n \,- \,2}\over 2} \cdot \Lambda^{{n\, -\, 4}\over 2} \cdot {{ (\Lambda^2 - \Vert\,\, Y \ - \ \Xi\,\Vert^{\,2})}\over {(\,\Lambda^2 \ + \ \Vert\,\, Y \ -\  \Xi\,\Vert^{\,2}\,)^{{n }\over 2} }}\ ,\\[0.2in]
  {{\partial\, V_{\Lambda\,, \, \Xi}}\over {\partial\, \Xi_j}} \ = \ {{\partial}\over {\partial\, \Xi_j}} \left[\ \left({\Lambda\over {\Lambda^2 + \Vert\,\, Y \ - \ \Xi\,\Vert^{\,2}}}\  \right)^{{n - 2}\over 2} \  \right] & = & -\,{{n\, -\, 2}\over 2} \cdot \Lambda^{{n - 2}\over 2} \cdot {{ 2\,(\ \Xi_j\, - \,Y_j\,  )}\over {(\,\Lambda^2 \ +\  \Vert\,\, Y - \Xi\,\Vert^{\,2}\,)^{{n }\over 2} }}\\[0.2in]
& \ & \ \ \ \ \ \ \ \ \ \  \ \ \ \ \ \ \mfor \ m \ = \ 1\,, \ \cdot \cdot \cdot\,,\, \ n\,. \ \ \ \
\end{eqnarray*}
Conditions  (\,1.22\,) \ $\Longrightarrow \ \Vert \ \Xi_j  \ - \ \Xi_k\ \Vert \to \ \infty\ $\, as $\,i \ \to \ \infty\,$\,,\, where $1\ \le \ j \ \not= \ k \ \le \ \flat_i$\ . Hence for any fixed $\,R\ > \ 0\,$ and $\,i\ \gg\ 1\,,$\, we have\\[0.1in]
(\,A.3.66\,)
\begin{eqnarray*}
 {\bf W}_{\,\flat_{\,i}}(\ Y\,)   & = &   \left(\  { {\Lambda_{\,I_i} } \over {\ \Lambda_{\,I_i}^2 + \Vert\,\, Y \ - \ \Xi_{\,I_i} \,\Vert^{\,2}\  }}\  \right)^{\!\!{{n - 2}\over 2}} \ + \ \sum_{l \ \not= \ I_i}  \left(\  { {\Lambda_{\,I_i} } \over {\ \Lambda_{\,I_i}^2 + \  \Vert\,\ Y \ - \ \Xi_{\,l} \Vert^{\,2}\  }}\  \right)^{\!\!{{n \,- \,2}\over 2}}\\[0.2in]& \le  &   \left(\  { {\Lambda_{\,I_i} } \over {\ \Lambda_{\,I_i}^2 + \Vert\,\, Y \ - \ \Xi_{\,I_i} \,\Vert^{\,2}\  }}\  \right)^{\!\!{{n - 2}\over 2}} \\[0.2in]
  & \ &  \hspace*{1in} + \ \sum_{l \ \not= \ I_i}  \left(\  { {\Lambda_{\,I_i} } \over {\ \Lambda_{\,I_i}^2 + \ [ \   \Vert\ \Xi_{\,I_i}  \ - \ \Xi_{\,l}  \,\Vert \ - \ \Vert\,\ Y \ - \ \Xi_{\,I_i} \Vert \ ]^{\,2}\  }}\  \right)^{\!\!{{n \,- \,2}\over 2}}\\[0.2in]
  & \le & \left(\  { {\Lambda_{\,I_i} } \over {\ \Lambda_{\,I_i}^2 + \Vert\,\, Y \ -\  \Xi_{\,I_i} \,\Vert^{\,2}\  }}\  \right)^{\!\!{{n - 2}\over 2}} \ \  + \ O\,\left(\  \sum_{l \ \not= \ I_i}  \left[\  {{ 1 } \over {\  \left[\ {1\over 2} \cdot \Vert\,\Xi_{\,I_i}  \ - \ \Xi_{\,l}\,\Vert\  \right]^{\,2}\  }}\  \right]^{\,{{n\, - \,2}\over 2}} \ \right)\\[0.2in]
  & \ & \hspace*{1in} \left[ \ \Vert\,Y \ - \ \Xi_{\,I_{\,i}} \,\Vert \ < \ R \ \ \Longrightarrow \ \ \Vert\,Y \ - \ \Xi_{\,I_{\,i}} \,\Vert  \ \le \ {1\over 2} \cdot  \Vert\ \Xi_{\,I_i}  \ - \ \Xi_{\,l}  \,\Vert\ \right] \\[0.1in]
  & =  &  \left(\  { {\Lambda_{\,I_i} } \over {\ \Lambda_{\,I_i}^2 + \Vert\,\, Y \ - \ \Xi_{\,I_i} \,\Vert^{\,2}\  }}\  \right)^{\!\!{{n - 2}\over 2}} \ \  + \ O\,\left(\   {\bar\lambda}_{\ {\flat_{\,i}}}^{ { {{n\,-\,2}\over 2 }  \, \cdot \,\gamma} }\ \right)   \\[0.25in]
  & \to &  \left(\  { {\Lambda_{\,I_i} } \over {\ \Lambda_{\,I_i}^2 + \Vert\,\, Y \ - \ \Xi_{\,I_i} \,\Vert^{\,2}\  }}\  \right)^{\!\!{{n - 2}\over 2}} \ \ \ \ \ \ {\mbox{uniformly \ \ in \ \ }} B_{ \ \Xi_{\,I_i}}\,(\,R\,) \ \ \ {\mbox{as}} \ \ i \ \to \ \infty \ .
\end{eqnarray*}

\newpage

Together with conditions (\,1.26\,)\,,\, (\,1.28\,) of the main text\,,\, and  the result in (\,A.3.49\,)\,,\, standard theory in elliptic partial differential equations can be used to show that there is a subsequence of $\,\{ \, {\tilde\Phi_{\,I_{\,i}}} \, \}\,$ which converges   uniformly in $\,B_o\,(R)\,$ (\,for any fixed $\,R\,>0\,$)  to $\,\Phi_{\infty} \, \in \, {\cal D}^{\,1\,,\,2}\,$,\, which is a solution of the equation
$$
\Delta\,\Phi_{\infty} \ + \ n\,(\,n\,+\,2\,) \cdot V_{\Lambda\,, \ 0}^{4\over {n\,-\,2}} \cdot\, \Phi_{\infty} \ = \ 0 \ \ \ \ \ \ {\mbox{in}} \ \ \ \R^n \ \ \ \ \ \ \  (\,\Lambda \ > \ 0 \ \ {\mbox{is \ \ fixed}} \ )\ .\leqno (\,A.3.67\,)
$$
Here
$$
V_{\,\Lambda\,, \ 0} \ = \ \left({\Lambda\over {\Lambda^2 + \Vert\,\, Y  \,\Vert^{\,2}}}\  \right)^{{n - 2}\over 2} \ \ .
$$
Via (\,A.3.52\,)\,,\, we have
$$
\sup\ \bigg\{ \ |\,\Phi_\infty\,(\ {\bar Y}\,)\ | \ \ \  \bigg\vert \ \ \ \  Y\, \in \,  B_o\, (\,R_{\,o})   \ \bigg\} \ \,\ge \ \varepsilon_o \ > \  0\,. \leqno (\,A.3.68\,)
$$
Hence $\,\Phi_{\infty}\ \not\equiv \ 0\,.\,$ From the classification of
non\,-\,negative  solutions of equation (\,A.3.67\,) \cite{Progress-Book}, we obtain
$$  \Phi_{\infty} \ = \  A_o \cdot \left(\ {{V_{ \Lambda\,,\ \Xi}}\over {\partial \,\Lambda}}\ \ \right) \, \bigg\vert_{\ (\,\Lambda\,,\ 0\,)} \, + \ \ \sum_{j\,=\,1}^n \ B_j  \cdot \left(\ {{\partial\, V_{ \Lambda\,,\ \Xi} }\over {\partial \ \Xi_{\,|_{\,j}} }} \ \ \right) \, \bigg\vert_{\ (\,\Lambda\,,\ 0\,)} \ {\bf ,} \leqno (\,A.3.69\,) $$
where $\,A_o\,,\, \ B_1\,, \ \cdot \cdot \cdot\,, \ B_n\,$ are coefficients. Due to (\,A.3.67\,)\,,\, at least one of the coefficient is non\,-\,zero\,.\bk
From (\,A.3.63\,) and (\,A.3.65\,) [\,with $\,\Xi \ = \ 0$\ ]\,,\, we have
$$
\bigg\vert \ V_{\Lambda_{I_{\,i}}\,,\ 0}^{4\over {n\,-\,2}} \cdot \left(\ \Lambda_{I_{\,i}} \cdot {{\partial\, V_{\Lambda_{I_{\,i}}\,,\ 0} }\over {\partial\,  \Lambda_{I_{\,i}} }}\right)\,\cdot\, {\bar\Phi}_{\,I_{\,i}} \  \bigg\vert_{\,Y} \le \, {C\over { \ (\ 1 \ + \ R\ )^{\ n\,+\,2}}} \ \ \ \ \  \mbox{for \ \ all} \ \ i \ \gg\ 1 \ \ \  (\,R \ = \ \Vert\,Y\,\Vert \ \ge \ 0 \,)\  .
$$
Likewise for the partial derivatives in $\,\Xi\,.\,$ Applying the Lebesque dominated convergence theorem to the integral in (\,A.3.65\,)\,,\, we see that
\begin{eqnarray*}
\bigg\langle \ V_{\Lambda\,,\ 0}^{4\over {n\,-\,2}} \cdot \left(\ \Lambda \cdot {{\partial\, V_{\Lambda\,,\ 0} }\over {\partial\, \Lambda}}\right)\ \,, \ \ \ {\bar\Phi}_{\,I_{\,i}} \  \bigg\rangle_{\!\int\ } & = & 0\\[0.2in]
{\mbox{and}} \ \ \ \ \ \ \ \ \bigg\langle \ V_{\Lambda\,,\ 0}^{4\over {n\,-\,2}} \cdot \left(\ \Lambda \cdot {{\partial\, V_{\Lambda\,,\ \Xi} }\over {\partial\, \Xi_{\,|_m} }}\right)\ \bigg\vert_{\,\Xi \ = \ 0}\,, \ \ \ {\bar\Phi}_{\,I_{\,i}} \  \bigg\rangle_{\!\int\ } & = &  0  \ \ \ \  \ \  \mfor \ m \ = \ 1\,, \ \cdot \cdot \cdot, \ n\,.\\
\end{eqnarray*}
These imply that the coefficients in (\,A.3.68\,) are all equal to zero, that is,
$$
A_o \ = \ B_ 1 \ = \ B_2 \ = \ \cdot \cdot \cdot \ = \ B_n \ = \ 0 \ \ \ \Longrightarrow \ \ \Phi_\infty \ \equiv \ 0 \ \ \ {\mbox{in}} \  \ \R^n\,.
$$
But this contradicts (\,A.3.68\,)\,.\, Hence the claim ${\bf(\,*\,)}_{\,(\,A.3.51\,)}\,$ must hold.


\newpage


{\large{\bf \S\,A\,4\,.  Solving the  linear inhomogeneous equation \,-\, }\\[0.1in] \hspace*{0.7in}{\large{\bf{\bf multiple bubbles case\,.\,}} }}\\[0.2in]
In this section we demonstrate that if the bubbles are well\,-\,separated [\,as expressed in the Smallness Lemma A.3.20\,]\,,\, then many aspects of the one bubble situation (\,{\bf \S\,A 1}\,) can be generalized to the multiple bubbles case. (\,See \,{\bf \S\,A 4\,b\,}\, for some of the differences.\,) \\[0.2in]
{\bf \S\,A\,4\,.\,a\,.} \ \ {\bf Set up in the Hilbert Space} $\ {\cal D}^{\,1\,,\,2} \,$. \ \ We are interested in finding a solution  of the linear inhomogeneous equation\\[0.1in]
(A.4.1)
$$
\Delta_Y \,\Phi \ + \  \left(\ {{n\,+\,2}\over {\,n\,-\ 2\,}} \cdot [\,c_n \cdot K\,] \cdot {\bf W}_\flat^{4\over{n\,-\,2}}\ \right)\cdot \Phi \ = \ \left(V_{1\,,\  o}^{4\over {n\,-\,2}} \cdot {\cal H} \right)  \ \ \ \ \ \  {\mbox{in}} \ \ \R^n\ \ \ \ \ (\,{\cal H}\,\in\, {\cal D}^{\,1\,,\,2}\,)\,.
$$
As in the above,
$$
V_{1\,,\  o}\,(\,Y\,) \ = \ \left(\ {1\over {\ 1\ + \ \Vert\,Y\,\Vert^2\   }}\ \right)^{\!\!{{\,n\,-\,2\,}\over 2} } \ \ \  \ \mfor \ \ Y\,\in\,\R^n\ .
$$
In weak form, (\,A.4.1\,) can be expressed as \\[0.1in]
(A.4.2)
$$
\int_{\R^n} \left\{\,  \Delta \,\Phi \, + \, \left(\ \ {{n\,+\,2}\over {\,n\,-\ 2\,}} \cdot  [\,c_n \cdot K\,] \cdot {\bf W}_\flat^{4\over{n\,-\,2}}\ \right) \cdot \Phi \ - \ \left({V_{1\,,\  o}^{4\over {n\,-\,2}}} \cdot {\cal H}   \right) \  \right\} \,\cdot\, \Psi \ = \ 0
$$ for all $\,\Psi \ \in \ {\cal D}^{\,1\,, \ 2}.\
$
Similar to {\bf \S\,A\,1}\,,\, the problem  can be cast in the form of investigating  the (functional) equation
\begin{eqnarray*}
(A.4.3) \ \ \ \ \ \ \ \ \ \ \ \ \ \ \ \ \ \ \  \ \ \ \ \ \ \ \ \ \ \ \ \ \ \ \ \ \ \  (\,{\bf I} \ + \ {\bf K}_{\,\flat}\,)\!\!\! & : & {\cal D}^{\,1\,, \ 2} \ \rightarrow \ {\cal D}^{\,1\,, \ 2}\\[0.2in]
(\,{\bf I} \ + \ {\bf K}_{\,\flat}\,) \,(\Phi) & = &  {\overline {\cal H}}   \ \ \ \ \ \ \ \ \ \ \left(\ {\overline {\cal H}} \,\in\  {\cal D}^{\,1\,, \ 2}\,\right) \,. \ \ \ \ \ \ \ \ \ \ \ \ \ \ \ \ \ \ \  \ \ \ \ \ \ \ \ \ \ \ \ \ \ \ \ \ \ \
\end{eqnarray*}
${\overline {\cal H}}$\  is related to $\,{\cal H}\,$ via
  the Riesz Representation Theorem\,:
$$
 \langle\  {\overline {\cal H}} \,, \ \Psi  \,\rangle_\btd \  = \  -\,\int_{\R^n}    \left(\ {V_{1\,,\  o}^{4\over {n\,-\,2}}} \cdot {\cal H} \right)   \cdot \Psi \ \ \ \ {\mbox{for \ \ all}} \ \ \Psi \ \in \ {\cal D}^{\,1\,, \ 2}. \leqno (A.4.4)
$$
The factor  $\,\displaystyle{\left( V_{1\,,\  o}^{4\over {n\,-\,2}}\  \right)}\,$ is implemented so as to make the right hand side of (\,A.4.4\,) into a {\it bounded} linear function [\,via the Sobolev inequality\,;\, cf. the argument in (\,A.1.12\,)\,]\,.\, Moreover, we observe that
\begin{eqnarray*}
(A.4.5) \ \ \ \ \ \ \ \   \ \ \ \ \ \ \ \   \langle\  {\overline {\cal H}} \,, \ \,{\overline {\cal H}}   \,\rangle_\btd & = & -\,\int_{\R^n}    \left(\ {V_{1\,,\  o}^{4\over {n\,-\,2}}} \cdot {\cal H} \right)   \cdot {\overline {\cal H}}  \ \le \  C\,(n) \cdot \Vert\, {\cal H} \,\Vert_{\,\btd} \cdot \Vert\, {\overline {\cal H}} \,\Vert_{\,\btd} \ \ \ \ \ \ \ \   \ \ \ \ \ \ \ \  \\[0.2in]
 \Longrightarrow \ \ \ \  \Vert\, {\overline {\cal H}} \,\Vert_{\,\btd} & \le & C\,(n) \cdot \Vert\, {\cal H} \,\Vert_{\,\btd} \ .
\end{eqnarray*}
Here we apply the argument as in (\,A.1.12\,)\,.\, See also (\.A.4.11\,) below\,.\,  \bk
%
%
%
%
Moving forward\,,\,  in (A.4.4)\,,\, $\,{\bf I}\,$ is the identity map, and
$$
 {\bf K}_{\,\flat} \ : \  {\cal D}^{\,1\,, \ 2} \ \rightarrow \ {\cal D}^{\,1\,, \ 2} \leqno (\,A.4.6\,)
$$
is a linear operator given by\\[0.1in](A.4.7)
$$
 \langle\  {\bf K}_{\,\flat} \,(\Phi)\,, \ \Psi  \,\rangle_\btd \ = \  -\,\int_{\R^n}   \left(\  {{n\,+\,2}\over {\,n\,-\ 2\,}} \cdot  [\,c_n \cdot K\,]  \cdot  {\bf W}_\flat^{4\over{n\,-\,2}}\ \right) \cdot \Phi   \cdot \Psi \ \ \ \ {\mbox{for \ \ all}} \ \ \ \Psi \ \in \ {\cal D}^{\,1\,, \ 2} \ {\bf .}
$$
As in {\bf \S\,A\,1}\,,\, for fixed centers of the bubbles in $\,{\bf W}_\flat\,$,\, we obtain  the decay of $\,{\bf W}_\flat\,(\,Y\,)\,$ near infinity in the order of $\,\Vert\, Y\,\Vert^{-\,(\,n\,-\,2\,) }$\, for $\,\Vert\,Y\,\Vert \ \gg \ 1\,.\,$ It  can be then checked that $\,
 {\bf K}_{\,\flat} \ : \  {\cal D}^{\,1\,, \ 2} \ \rightarrow \ {\cal D}^{\,1\,, \ 2}
$ is a compact linear operator   [\,via the Rellich\,-\,Kondrachov Compactness Theorem\,,\, cf. \cite{Olsen-Holden} \ ]\ )\ .

\vspace*{0.5in}

{\bf \S\,A\,4\,.\,b\,.} \   {\bf Differences between single and multiple bubble cases.} \ \  Firstly, unlike the case of one bubble (\,that is\,,\, $\flat \ = \ 1\,)\,,\,$  for
$$
[\,c_n \cdot K\,]  \ \equiv \ n\,(\,n\,-\,2) \ \ \ \ \ {\mbox{in}} \ \ \R^n
$$
and $\,\flat \,\ge \,2\,,\,$
$$
(\,{\bf I} \ + \ {\bf K}_{\,\flat}\,) \,( \, \partial_{\,\Xi_{\,l_{|_{\,j}}}}  \! {\bf V}_{\,l} \,)\ \not\equiv \ 0\,.
$$
Instead, one has
$$
\Vert\, (\,{\bf I} \ + \ {\bf K}_{\,\flat}\,) \ ( \, \partial_{\,\Xi_{\,l_{|_{\,j}}}}  \! {\bf V}_{\,l} \,)\ \Vert_\btd \ \to \ 0 \ \ \ \ {\mbox{as}} \ \ \ \ {\bar\lambda}_{\,\,\flat} \ \to \ 0^{\,+}\ .
$$
This can
 be verify by using the Riese Representation Theorem and similar calculations as in {\bf \S\,A\,4\,.\,g}\,.\, We  are not at a easy  position to precisely determine  the null space of $\,(\,{\bf I} \ + \ {\bf K}_{\,\flat}\,)\,$,\, even when $
[\,c_n \cdot K\,]  \ \equiv \ n\,(\,n\,-\,2)
$\,,\,   as opposite to the  one bubble situation [\,cf. (\,A.1.14\,)\,]\,.\,\bk
 Secondly\,,\, concerning the $\,\perp$\,-\,conditions [ \ cf. (\,A.1.18\,)\ ]\,,\,
\begin{eqnarray*}
\bigg\langle  \  {\overline {\cal H}}\,, \ {{\partial\, V_{\,\Lambda \,,\  \Xi}}\over {\partial \,\Lambda}}\ \bigg\rangle_\btd \ = \ 0 &  \Longleftrightarrow & \bigg\langle \left(\ V_{1\,,\  o}^{4\over {n\,-\,2}}\cdot {\cal H}\right)\,, \ \,{{\partial\, V_{\,\Lambda \,,\  \Xi}}\over {\partial\, \Lambda}}\, \bigg\rangle_{\!{\int}} \ = \ 0\\[0.2in]
&  \not\!\!\Longrightarrow & \bigg\langle  \,  {\cal H}\,, \ {{\partial\, V_{\,\Lambda \,,\  \Xi}}\over {\partial \,\Lambda}}\,\bigg\rangle_\btd  \ = \ 0 \ \ \ \ \ \ \ \left(\  {\mbox{unless}} \ \  V_{\,\Lambda \,,\  \Xi} \ = \ V_{1\,,\  o}\,\right)\,.
\end{eqnarray*}
 Instead\,,\, we consider the decomposition (\,valid under the conditions in Proposition 2.7 of the main text\,)
 $$
 {\cal D}^{\,1\,, \ 2}_{{\flat}} \ = \ {\cal D}^{\,1\,, \ 2}_{{\flat}}\,(\,\perp\,) \ \oplus \ {\cal D}^{\,1\,, \ 2}_{{\flat}}\,(\,/\!/\,)\ . \leqno (A.4.8)
 $$

 \newpage

 Here\\[0.1in]
 (A.4.9)
 $$
{\cal D}^{\,1\,, \ 2}_{{\flat}}\,(\,/\!/\,) \,:=\,{\mbox{Span}} \ \ \bigg\{ \ (\,\Lambda_{\,l} \,\cdot \, \partial_{\,\Lambda_{\,l}}\,)\,{\bf V}_l\,, \ \ \ (\,\Lambda_{\,l}\, \cdot \, \partial_{\,\Xi_{\,l_{|_{\,j}}}}\,)\,{\bf V}_l \,, \   \ \ \  j \ = \ 1\,,\ \cdot \cdot \cdot\,, \ n   \ \ \ \ {\mbox{and}} \ \ \ l \ = \ 1\,,\ \cdot \cdot \cdot\,, \ \flat \  \bigg\} \ ,
$$
 We first look at the $\,\perp$\,-\,space $\,{\cal D}^{\,1\,, \ 2}_{{\flat}}\,(\,\perp\,)\,.$

\vspace*{0.3in}

{\bf \S\,A\,4\,.\,c\,.} \ \  {\bf The $\,\perp$\,-\,space}\, $(\ {\cal D}^{\,1\,, \ 2}_{{\flat}}\,(\,\perp\,)\,,\ \ \langle \ \ \ \rangle_{\,\btd}\ )\,$ {\bf and Fredholm Alternative\,.}

\vspace*{0.2in}

As  $$\ (\,{\cal D}^{\,1\,, \ 2}_{{\flat}}\,(\,\perp\,)\,,\ \ \langle \ \ \ \rangle_{\,\btd}\ ) $$ is itself  a complete Hilbert space\,,\, we consider the restriction of the homogeneous equation [\,similarly defined as in (\,A.4.7\,)\,]
 \begin{eqnarray*}
(\,A.4.10\,) \ \ \ \  \ \ \ \   \ \ \ \   \ \ \ \   \ \ \ \   (\,{\bf I} \ + \ {\bf K}_{\,\flat}\,) & : & {\cal D}^{\,1\,, \ 2}_{{\flat}}\,(\,\perp\,) \ \rightarrow \ {\cal D}^{\,1\,, \ 2}_{{\flat}}\,(\,\perp\,)\\[0.2in]
(\,{\bf I} \ + \ {\bf K}_{\,\flat}\,) \,(\Phi) & = & \vec{\,0}  \ \ \ \ \ \ \ \ \ \ \ \ \ \ \  \ \ \ \ \ \ \ \ \ \ \  \ \ \left(\,\Phi \,\in\ {\cal D}^{\,1\,, \ 2}_{{\flat}}\,(\,\perp\,)\ \right)\,. \ \ \ \ \ \ \ \ \  \ \ \ \   \ \ \ \
\end{eqnarray*}
Assume that it has a nontrivial solution  $\, \Phi_{\,\flat_{\,i_{\,o}}}    \in \, {\cal D}^{\,1\,, \ 2}_{{\flat}}\,(\,\perp\,)\, \setminus \, \{\,0\,\}\,$ for   $\,{\bar\lambda}_{\,\flat_i}\ > \ 0$\,,\, where $\,{\bar\lambda}_{\,\flat_i}\ \to \ 0^{\,+}$\,.\,    That is,
\begin{eqnarray*}
& \ &  \langle\ \Phi_{\,\flat_{\,i_{\,o}}}   \,, \ \Psi  \,\rangle_\btd  \ + \ \langle\, {\bf K}_{\,\flat}  \,(\Phi_{\,\flat_{\,i_{\,o}}})\,, \ \Psi  \,\rangle_\btd \ = \ 0  \ \ \ \ {\mbox{for \ \ all}} \ \ \Psi \ \in \ {\cal D}^{\,1\,, \ 2}_{{\flat}}\,(\,\perp\,)\\[0.2in]
\Longrightarrow & \ & \Phi_{\,\flat_{\,i_{\,o}}} \ + \  {\bf K}_{\,\flat}  \,(\,\Phi_{\,\flat_{\,i_{\,o}}})  \ \in \ {\cal D}^{\,1\,, \ 2}_{{\flat}_{\,\parallelsum}}\ .
\end{eqnarray*}
It follows that \\[0.1in](\,A.4.11\,)
\begin{eqnarray*}
 & \ &  \Phi_{\,\flat_{\,i_{\,o}}} \ + \  {\bf K}_{\,\flat}  \,(\,\Phi_{\,\flat_{\,i_{\,o}}}) \ =  \ \sum_{l\,=\,1}^\flat \ \left\{ \   a_{\,l} \cdot  \left[\, (\,\Lambda_{\,l} \ \partial_{\Lambda_{\,l}}\,) \, {\bf V}_{\,l}\ \right] \ + \ \sum_{j\,=\,1}^n  \ b_{\,l\,,\,j} \cdot    \left[\, (\,\Lambda_{\,l} \ \partial_{\,\Xi_{\,l_{|_{\,j}}}}\,) {\bf V}_{\,l}\,\right]     \ \right\} \\[0.15in]
& \ & \hspace*{4.5in}\ \ \ [\ {\mbox{via \ \ (A.4.9)}}\,]\ .
\end{eqnarray*}


\newpage

{\bf \S\,A\,4\,.\,d\,.} \ \
{\bf  Tracing back the situation. } \ \  We highlight the correspondence in (\,A.4.4\,) for the following special case.\\[0.1in]
(\,A.4.12\,)
\begin{eqnarray*}
  \bigg\langle \,\left( {{\partial \,V_{\Lambda\,,\ \Xi} }\over {\partial \Lambda}}\  \right)   \,, \ \Psi \bigg\rangle_\btd &  = &   \int_{\R^n}  \bigg\langle \,\btd\,{{\partial \,V_{\Lambda\,,\ \Xi} }\over {\partial \Lambda}}   \,, \ \btd\,\Psi \bigg\rangle  \  = \ -\,  \int_{\R^n}  \Delta\,\left( {{\partial \,V_{\Lambda\,,\ \Xi} }\over {\partial \Lambda}}\  \right)  \cdot \,\Psi\\[0.2in]
   & = &
     -\,   \int_{\R^n}  \left(\ {{\partial\, (\,\Delta\, V_{\Lambda\,,\ \Xi}\,) }\over {\partial \Lambda}}\  \right)  \cdot \,\Psi \ = \   n\,(\,n\,-\,2\,) \cdot \int_{\R^n}   \left(\ {{\partial \,V_{\Lambda\,,\ \Xi}^{{\,n\,+\,2\,}\over {n\,-\,2}}  }\over {\partial \Lambda}}\  \right)  \cdot \,\Psi\\[0.2in]
          &  & \hspace*{-1.2in} = \ n\,(n\,+\,2) \cdot    \int_{\R^n}  V_{\Lambda\,,\ \Xi}^{4\over {n\,-\,2}} \cdot  \left(\ {{\partial   V_{\Lambda\,,\ \Xi} }\over {\partial \Lambda}}\  \right)  \cdot \,\Psi \ = \ n\,(n\,+\,2) \cdot \bigg\langle \left(V_{\Lambda\,,\ \Xi} ^{4\over {n\,-\,2}} \cdot {{\partial \,V_{\Lambda\,,\ \Xi} }\over {\partial \Lambda}}  \right), \ \Psi \bigg\rangle_{\!\int}\ . 
\end{eqnarray*}
Coming back to equation (A.4.10) -  one traces back the relations in (A.4.4)\,-\,(A.4.7) and shows  that $\,\Phi_{\,\flat_{\,i_{\,o}}}\,$ satisfies the linear inhomogeneous equation\,:
\\[0.1in]
(A.4.13)
\begin{eqnarray*}
& \ & \Delta \, \Phi_{\,\flat_{\,i_{\,o}}}\ + \  {{n\,+\,2}\over {\,n\,-\ 2\,}} \cdot (\,c_n \cdot K\,) \cdot W_{\flat}^{4\over{n\,-\,2}}\cdot  \Phi_{\,\flat_{\,i_{\,o}}} \\[0.2in]
& = &  0 \ + \  \sum_{l\,=\,1}^\flat \  \left\{ \ \sum_{j\,=\,1}^n  \ \beta_{\,l\,,\,j} \cdot   {\bf V}_{\,l}^{4\over {n\,-\,2}} \cdot  \left[\, (\,\Lambda_{\,l} \ \partial_{\,\Xi_{\,l_{|_{\,j}}}}\,)\, {\bf V}_{\,l}\,\right]  \ + \     \alpha_{\,l} \cdot  {\bf V}_{\,l}^{4\over {n\,-\,2}} \cdot  \left[\, (\,\Lambda_{\,l} \ \partial_{\Lambda_{\,l}}\,)\,{\bf V}_{\,l}\,\right]\  \right\}\ .\\[-0.1in]
\end{eqnarray*}
Furthermore,
$$
\Vert \, \Phi_{\,\flat_{\,i_{\,o}}}\,\Vert_{\,*_Y} \ \ \ \  \ {\mbox{makes \ \ sense}} \ \ \{\, {\mbox{see}} \ \ (A.3.50\emph{}) \ \ {\mbox{and}} \ \ [\,9\,]\,\}\,, \ \ {\mbox{and \ \ we \ \ can \ \ set \ \ }} \Vert \, \Phi_{\,\flat_{\,i_{\,o}}}\,\Vert_{\,*_Y} \ = \ 1\,.
$$

\vspace*{0.1in}

{\bf \S\,A\,4\,.\,e\,.} \ \
{\bf Isomorphism. } \ \
Applying Smallness Proposition A.3.20, we obtain
$$ \ \ \ \ \  \ \ \ \ \
\Vert \, \Phi_{\,\flat_{\,i_{\,o}}} \,\Vert_{\,*_Y} \ \to \ 0 \ \ \ \ \  \ {\mbox{as}} \ \ \ \ {\bar\lambda}_{\,\,\flat_{\,i}} \ \to \ 0^{\,+}\,.
$$
But this is a contradiction with (A.4.13)\,.\, Hence (A.4.12) has no non\,-\,trivial solution when $\,{\bar\lambda}_{\,\,\flat}  \ (\,> \ 0\,)\,$ is small enough. That is\,,\,  for all $\,{\bar\lambda}_{\,\,\flat}  \ (\,> \ 0\,)\,$ small enough\,,\,  the operator
$$
(\,{\bf I} \ + \ {\bf K}_{\,\flat}\,) \,: \ {\cal D}^{\,1\,, \ 2}_{{\flat}}\,(\,\perp\,) \ \to  \ {\cal D}^{\,1\,, \ 2}_{{\flat}}\,(\,\perp\,) \leqno (A.4.14)
$$
is an isomorphism\, (\,using the self adjoint property and the  Fredholm alternative method as in {\bf \S\,A\,1}\,)\,.\, Moreover, via the Open Mapping Theorem,
$$
(\,{\bf I} \ + \ {\bf K}_{\,\flat}\,) \,(\Phi) \ = \ {\cal H}_{\,\perp}\ \ \Longrightarrow \ \  \Vert\, \Phi \,\Vert_{\,\btd} \ \le \ C^*_1 \cdot   \Vert\,  {\cal H}_{\,\perp} \,\Vert_{\,\btd} \ . \leqno (A.4.15)
$$
Note that $\,C^*_1\,$  may depend  on $\,{\bar\lambda}_{\,\,\flat}\,.\,$ We provide a  ``\,stable\," estimate in Proposition A.4\.20 and in  {\bf{\,\S\,A\,4\,.\,g}}\,.\, 

\newpage

{\bf \S\,A\,4\,.f\,.} \ \
{\bf  Back to the full space} $\ {\cal D}^{\,1\,, \ 2}\,.\,$ \  \ \  Consider the (full) equation
\begin{eqnarray*}
(\,{\bf I} \ + \ {\bf K}_{\,\flat}\,) & : & {\cal D}^{\,1\,, \ 2} \ \rightarrow \ {\cal D}^{\,1\,, \ 2}\\[0.2in]
(\,{\bf I} \ + \ {\bf K}_{\,\flat}\,) \,(\Phi) & = &  {\overline {\cal H}}   \ \  \ \ \ \ \ \ \ \ \  \ \ \ \ \ \ \ \ \ (\ {\overline {\cal H}}  \ \in\  {\cal D}^{\,1\,, \ 2}\,)\ .
\end{eqnarray*}
Here $\,{\overline {\cal H}}\ \,$ is given by (A.4.5)\,.
Via the decomposition in (A.4.8)\,,\, we write
$$
{\overline {\cal H}}  \ \,=\, \ {\overline {\cal H}} _\perp \ +\ ( -\,{\overline H}_{/\!/}\, ) \ \ \Longleftrightarrow \ \ {\overline H}_\perp \ = \ {\overline H} \  + \ {\overline H}_{/\!/}\ , \leqno (A.4.16)
$$
where
$$
{\overline {\cal H}}_\perp  \ \in \ {\cal D}^{\,1\,, \ 2}_{{\flat}}\,(\,\perp\,)\ ,
$$
and
\begin{eqnarray*}
(A.4.17) \ \ \ \ \ \  & \ & {\overline {\cal H}}_{/\!/}  \ \in \ {\cal D}^{\,1\,, \ 2}_{{\flat}}\,(\,{/\!/}\,)\\[0.2in]
\Longrightarrow & \ &
 {\overline {\cal H}}_{/\!/} \ = \ \sum_{l\,=\,1}^\flat \ \left\{ \ \sum_{j\,=\,1}^n  B_{\,l\,,\,j} \cdot    \left[\, (\,\Lambda_{\,l} \ \partial_{\,\Xi_{\,l_{|_{\,j}}}}\,)\, {\bf V}_{\,l}\,\right]  \ + \     A_{\,l} \cdot  \left[\, (\,\Lambda_{\,l} \ \partial_{\Lambda_{\,l}}\,)\, {\bf V}_{\,l}\ \right] \ \right\} \ . \ \ \ \ \ \ \ \ \ \ \ \
\end{eqnarray*}
Via (A.4.14), the equation
$$
(\,{\bf I} \ + \ {\bf K}_{\,\flat}\,) \,(\Phi) \ = \  {\overline {\cal H}}_\perp  \ \ \ \ \ \ {\mbox{has \ \ an  \ \ unique \ \ solution \ \  }} \Phi \ \in \ {\cal D}^{\,1\,, \ 2}_{{\flat}}\,(\,\perp\,)\ . \leqno (A.4.18)
$$
Running through the transformations in  (\,A.4.4\,)\, and \,(\,A.4.7\,)\, one more time, we arrive at the conclusion that $\, \Phi \ [\ \in \ {\cal D}^{\,1\,, \ 2}_{{\flat}}\,(\,\perp\,)\ ]\,$ satisfies

\begin{eqnarray*}
\Delta_Y \,\Phi \ + \  {{n\,+\,2}\over {\,n\,-\ 2\,}} \cdot (\,c_n \cdot K) \cdot {\bf W}_\flat^{4\over{n\,-\,2}}\cdot \Phi & = &  \left(V_{1\,,\  o}^{4\over {n\,-\,2}} \cdot {\cal H} \right) \ + \ \hspace*{0.3in}\ [\,{\mbox{refer \ \ to  \ \ (A.4.3)}}\,]\\[0.2in]
& \ & \hspace*{-2.5in} + \  \sum_{l\,=\,1}^\flat\  \left\{ \ \sum_{j\,=\,1}^n  B_{\,l\,,\,j} \cdot {\bf V}_{\,l}^{4\over {n\,-\,2}} \cdot    \left[\, (\,\Lambda_{\,l} \ \partial_{\,\Xi_{\,l_{|_{\,j}}}}\,) \,{\bf V}_{\,l}\,\right]  \ + \     A_{\,l} \cdot {\bf V}_{\,l}^{4\over {n\,-\,2}} \cdot  \left[\, (\,\Lambda_{\,l} \ \partial_{\Lambda_{\,l}}\,) \,{\bf V}_{\,l}\,\right]\  \right\}  \ + \ \\[0.1in]
& \ & \hspace*{1.8in}[\,{\mbox{refer \ \ to  \ \ (\,A.4.17\,)}}\,]\\[0.1in]
& \ & \hspace*{-2in} + \ \sum_{l\,=\,1}^\flat\  \left\{ \ \sum_{j\,=\,1}^n  \beta_{\,l\,,\,j} \cdot   {\bf V}_{\,l}^{4\over {n\,-\,2}} \cdot  \left[\, (\,\Lambda_{\,l} \ \partial_{\,\Xi_{\,l_{|_{\,j}}}}\,) \,{\bf V}_{\,l}\,\right]  \ + \     \alpha_{\,l} \cdot  {\bf V}_{\,l}^{4\over {n\,-\,2}} \cdot  \left[\, (\,\Lambda_{\,l} \ \partial_{\Lambda_{\,l}}\,) \,{\bf V}_{\,l}\,\right]\  \right\}\\[0.1in]
& \ & \hspace*{1.7in}\ [\,{\mbox{refer \ \ to  \ \ (\,A.4.13\,)}}\,]
\end{eqnarray*}

\newpage

\begin{eqnarray*}
\Longrightarrow \ \  \Delta_Y \,\Phi \ + \  {{n\,+\,2}\over {\,n\,-\ 2\,}} \cdot (\,c_n \cdot K) \cdot {\bf W}_\flat^{4\over{n\,-\,2}}\cdot \Phi & = &  \left(V_{1\,,\  o}^{4\over {n\,-\,2}} \cdot {\cal H} \right) \ + \ \\[0.2in]
& \ & \hspace*{-2.9in} + \  \sum_{l\,=\,1}^\flat\  \left\{ \ \sum_{j\,=\,1}^n  b_{\,l\,,\,j} \cdot {\bf V}_{\,l}^{4\over {n\,-\,2}} \cdot    \left[\, (\,\Lambda_{\,l} \ \partial_{\,\Xi_{\,l_{|_{\,j}}}}\,)\, {\bf V}_{\,l}\,\right]  \ + \     a_{\,l} \cdot {\bf V}_{\,l}^{4\over {n\,-\,2}} \cdot  \left[\, (\,\Lambda_{\,l} \ \partial_{\Lambda_{\,l}}\,)\, {\bf V}_{\,l}\,\right]\  \right\} \ ,
\end{eqnarray*}
where
$$
a_{\,l} \ = \ A_{\,l} \ + \ \alpha_{\,l} \ \ \ \ {\mbox{and}} \ \ \ \ b_{\,l} \ = \ b_{\,l} \ + \ \beta_{\,l}\ .\leqno (A.4.19)
$$

\vspace*{0.15in}

{\bf \S\,A\,4\,.g\,.} \ \ {\bf Existence, uniqueness and estimates.}\\[0.2in]
{\bf Proposition A.4.20.} \ \ {\it For $\,\flat \ \ge \ 2\,,\,$ assume that the bubble parameters $\,\{\, \lambda_{\,\,l}\, \}_{\,l\,=\,1}^{\,\flat}  \,$ and  $\ \{\, \xi_{\,l}\, \}_{\,l\,=\,1}^{\,\flat} \,$ satisfy the conditions  in}  (\,1.4\,)\,,  (\,1.8\,)\,, (\,1.22\,)\,, (\,1.24\,)\,--\,(\,1.28\,)\,,\,  \,(\,1.30\,)\, {\it and}  \,(\,1.32\,) {\it \,of the main text\,}\,.\, {\it There is  a small positive number $\,{{\underline\lambda}}_{\,\, \epsilon}\,$ so that for all $\,{\bar\lambda}_{\,\,\flat} $\, satisfying}\,
$$0 \ < \ {\bar\lambda}_{\,\,\flat} \ \le \ {{\underline\lambda}}_{\,\, \epsilon}\,,$$
{\it  and for any given  $\,{\cal H}\,\in\,{\cal D}^{1\,\,,2}\,,\,$    the  equation}
\begin{eqnarray*}
(A.4.21)  \!\!\!\!\!\!    & \ & \Delta \,\Phi \ + \ {{n\,+\,2}\over {\,n\,-\ 2\,}} \cdot (\,c_n \cdot K) \cdot {\bf W}_{\flat}^{4\over{n\,-\,2}}\cdot \Phi \ = \  \left(\,V_{1\,,\  o}^{4\over {n\,-\,2}} \cdot {\cal H} \right) \ + \ {  {\bf P}}_{\flat_{\,\parallelsum}} \ , \\[0.2in]
{\it{where}}  & \ &\!\!\!\!\!\!\!\!  {  {\bf P}}_{\flat_{\,\parallelsum}}\ = \   \sum_{l\,=\,1}^\flat \ \left\{ \ \sum_{j\,=\,1}^n \  b_{\,l\,,\ j} \cdot {\bf V}_{\,l}^{4\over {n\,-\,2}} \cdot    \left[\, \left(\,\Lambda_{\,l} \ \partial_{\,\Xi_{\,l_{|_{\,j}}}}\,\right) \,{\bf V}_{\,l}\,\right]  \ + \     a_{\,l} \cdot {\bf V}_{\,l}^{4\over {n\,-\,2}} \cdot  \left[\, (\,\Lambda_{\,l} \ \partial_{\Lambda_{\,l}}\,) \,{\bf V}_{\,l}\,\right]\  \right\}\ ,
\end{eqnarray*}
{\it has a  unique solution $\,\Phi \, \in \, {\cal D}^{\,1\,, \ 2}_{{\flat}}\,(\,\perp\,)\ ,\, $} {\it and the coefficients $\,a_{\,l}\,$ and $\ b_{\,l\,,\ j}\,$ } {\it are uniquely determined by $ \,{\cal H}\,$} [\,{\it{\,refer to\,}} \,(A.4.19)\,]\,.
{\it Moreover}\,,\, {\it if}
$$
{\bf H} \, := \, \left(\,V_{1\,,\  o}^{4\over {n\,-\,2}} \cdot {\cal H} \right) \ \in\ C^o \,(\,\R^n)\,  \ \ \ \ {\it and}  \ \ \ \ \Vert \, {\bf H} \,\Vert_{\,**_Y} \ \ \ \mbox{ {\it is \ \ well\,-\,defined\,{\bf ,}\,}}\, \leqno (\,A.4.22\,)
$$
{\it   then}  $\,\Phi\,$ is in $\,C^2\,(\R^n)\,,\,$  {\it and}  $\ \Vert \, \Phi  \,\Vert_{\,*_Y}\,$ {\it is also well\,-\,defined\,,\, satisfying the following estimates}\, [\,(\,A.4.23\,) {\,\it and\,} (\,A.4.24\,)\,,\, {\it with a possibly smaller choice of} $\,{{\underline\lambda}'}_{\, \epsilon}$\ ]\,.
\begin{eqnarray*}
(A.4.23) \ \ \ \ \  \ \ \ \ \ \ \ \ \ \ \Vert \, \Phi  \,\Vert_{\,*_Y} & \le &  {\overline{C}}_1 \cdot \Vert \, {\bf H} \,\Vert_{\,**_Y}\ .\\[0.2in]
(A.4.24) \ \ \ \ \ \  \ \ \ \ \ \ \ \ \ \  \Vert\,\Phi\,\Vert_\btd^2 & \le & {\overline{C}}_2 \cdot  \flat    \cdot  \, \Vert\,{\bf H}\,\Vert_{\,**_Y}^2 \ \ \ \  \ \ \ \ \ \ \ \ \ \ \  \left[\   {\it{recall\ \ that \ \ \  }}  \gamma \ >  \ \sigma \ \right]\ . \ \ \ \ \ \ \ \ \ \ \ \ \ \ \
\end{eqnarray*}
{\it{Here}} $\ {\overline{C}}_1\ $ and $\ {\overline{C}}_2\ $ {\it are independent on $\,{\bar\lambda}_{\,\,\flat}\,$} [\,{\it once the conditions stated in this proposition are satisfied}\ ]\,.

\vspace*{0.2in}

{\bf Proof.} \ \ We establish the existence in  {\bf A\,4\,.\,f}\,.\, \\[0.1in]
${\bf *}_1$ \ \ The regularity of $\,\Phi\,$ follows from standard Elliptic Theory, see, for example, \cite{Gilbarg-Trudinger}\,.\,\\[0.1in]
${\bf *}_2$ \ \  Uniqueness can be checked via a substraction and the uniqueness of the isomorphism [\,cf. (\,A.4.14\,)\ ]
 \begin{eqnarray*}
(\,{\bf I} \ + \ {\bf K}_{\,\flat}\,) & : & {\cal D}^{\,1\,, \ 2}_{{\flat}}\,(\,\perp\,) \ \rightarrow \ {\cal D}^{\,1\,, \ 2}_{{\flat}}\,(\,\perp\,)\\[0.2in]
(\,{\bf I} \ + \ {\bf K}_{\,\flat}\,) \,(\Phi) & = & 0  \ \ \ \ \ \ \Longrightarrow \ \ \Phi \ \equiv  \ 0\ \ \ \ \  {\mbox{in}} \ \ \R^n .
\end{eqnarray*}
${\bf *}_3$ \ \ \,Under the condition in (\,A.4.22\,)\,,\, and using (\,A.3.50\,)\,, $\,
\Vert \, \Phi  \,\Vert_{\,*_Y}\,$ is well\,-\,defined. \\[0.1in]
${\bf *}_4$ \ \  As for the estimate in (A.4.22)\,,\,  suppose that the contrary holds, that is, a sequence $\,\{\,\Phi_{\,\flat_{\,i}}\,\}\,$ exists with $\,\bar\lambda_{\,\flat_{\i}} \ \to \ 0^{\,+}\,$,\, so that
$$
{{\Vert \, \Phi_{\,\flat_{\,i}}  \,\Vert_{\,*_Y}}\over { \,\Vert \,{\bf H}_{\,i} \,\Vert_{\,**_Y}  }} \ \to \ \infty \ \ \ \ {\mbox{for \ \ a \ \ sequence\  }} \ \ \bar\lambda_{\,\flat_{\i}} \ \to \ 0^{\,+}\,.
$$
Via a rescaling\,,\, we may assume that $\,\Vert \, {\bf H}_{\,i}  \,\Vert_{\,**_Y}  \to \,0^{\,+}\,$,\,  whereas $\,\Vert \, \Phi_{\,\flat_{\,i}} \,\Vert_{\,*_Y} \ = \ 1\,$ for all $\,i \ \gg \ 1\,.\,$  But this contradicts the Smallness Lemma A.3.20\,.\,\\[0.1in]
${\bf *}_5$ \ \ \,Finally\,,\, in the following subsection (\  {\bf \S\,A\,4\,.\,h}\ )\,,\, we  establish the estimate on $\,\Vert\,\Phi\,\Vert_\btd^2\,$  as detailed in  (\,A.4.23\,)\,.\qed
%

%
%
%
{\bf \S\,A\,4\,.\,h\,.} \ \
{\bf Estimate on} $\,\Vert\,\Phi\,\Vert_\btd^2\,$ - {\bf proof of (\,A.4.24\,)}\ .\\[0.15in]
From equation (A.4.21) and condition (\,A.4.22\,)\,,\, together with  $\,\Phi \, \in \, {\cal D}^{\,1\,, \ 2}_{{\flat}}\,(\,\perp\,)\ ,\, $\, we have\\[0.1in]
(A.4.25)
\begin{eqnarray*}
& \ & \Delta_Y \,\Phi \ + \ \left(\,{{n\,+\,2}\over {\,n\,-\ 2\,}}\,\right) \cdot ( \,c_n \cdot K)\cdot {\bf W}_\flat^{4\over{n\,-\,2}}\cdot \Phi \ \ =\  {\bf H} \ + \ {\bf P}_{\flat_{/\!/}} \\[0.2in]
\Longrightarrow & \ &  \int_{\R^n} (\Delta \,\Phi) \cdot \Phi \ \,+ \ \left(\,{{n\,+\,2}\over {\,n\,-\ 2\,}}\,\right)  \cdot \int_{\R^n} ( \,c_n \cdot K)\cdot {\bf W}_\flat^{4\over{n\,-\,2}}\cdot \Phi^2 \ =  \  \int_{\R^n} {\bf H} \cdot \Phi \ + \int_{\R^n} {\bf P}_{\flat_{/\!/}} \cdot \Phi\\[0.2in]
\Longrightarrow & \ &   \int_{\R^n} |\,\btd\,\Phi\,|^{\,2} \ \le \ \left(\,{{n\,+\,2}\over {\,n\,-\ 2\,}}\,\right) \cdot \int_{\R^n} ( \,c_n \cdot K)\cdot {\bf W}_\flat^{4\over{n\,-\,2}}\cdot \Phi^2 \ + \ \int_{\R^n} |\,{\bf H}| \cdot |\,\Phi\,|\\[0.1in]
& \ & \hspace*{2in} \ \ \ \ \ \ \  \left(\,\Phi \ \in \  {\cal D}^{\,1\,, \ 2}_{{\flat}}\,(\,\perp\,) \ \ \Longrightarrow \ \ \int_{\R^n} {\bf P}_{\flat_{/\!/}} \cdot \Phi \ = \ 0\ \right) \\[0.2in]
& \ & \hspace*{1.1in}\le \ C \cdot \left(\ \int_{\R^n} {\bf W}_\flat^{{2n}\over {\,n\,-\ 2\,}}\  \right)^{\!\!{2\over n}} \cdot \left(\ \int_{\R^n} |\,\Phi\,|^{{{2n}\over {\,n\,-\ 2\,}} } \right)^{\!\!\!{{n\,-\,2}\over {n}}} \\[0.2in]
  & \ & \hspace*{2in} \ \ \ \ \ \ \ + \  \left(\ \int_{\R^n} |\,{\bf H}|^{{2n}\over {n\,+\,2}}\  \right)^{\!\!{{n\,+\,2}\over {2n}}} \cdot \left(\ \int_{\R^n} |\,\Phi\,|^{{{2n}\over {\,n\,-\ 2\,}} } \right)^{\!\!\!{{n\,-\,2}\over {2n}}}\ .
\end{eqnarray*}
We establish the following in {\bf \S\,A\,4\,.\,j}\,:
\begin{eqnarray*}
(\,A.4.26\,) \ \ \ \ \ \ \ \  \    \int_{\R^n} {\bf W}_\flat^{{2n}\over {\,n\,-\ 2\,}} & = & \int_{\R^n} \left[ \  \sum_{l\,=\,1}^\flat   \left(\ {1\over {     1\ + \ \Vert\,Y \,-\ \Xi_{\,l}\,\Vert  }}\ \right)^{\!\! {{n\,-\,2}\over 2 }   }  \ \right]^{\,{{2n}\over {\,n\,-\ 2\,}} }  \ \ \ \ \ \ \ \  \ \ \ \ \ \ \ \  \\[0.2in]
& = & {\bar C}_{\,1} \cdot  \flat \ + \ o_{\,{\bar{\lambda}_{\,\flat}}}\,(\,1) \ , \ \ \ \ \ {\mbox{where}} \ \ \ {\bar C}_{\,1} \ := \  \int_{\R^n} V_{1\,,\  o}^{{2n}\over {\,n\,-\ 2\,}}\ {\bf .} \ \ \ \ \ \ \ \ \ \ \ \ \
\end{eqnarray*}
Here $\,o_{\,{\bar{\lambda}_{\,\flat}}}\,(\,1)   \,\to\, 0^{\,+}\,$ as ${\bar{\lambda}}_{\,\flat}  \ \to \ 0^{\,+}\,.$\, (\,Intuitively\,,\, it says that when the bubbles are well\,-\,separated, the integral can be approximated by the sum of the individual integration of a bubble\,.\,)
From definitions (A.3.6) and (A.3.7)\,,\, we see that
\begin{eqnarray*}
 |\,\Phi\,(\,Y\,)\,|  & \le & \Vert\,{\Phi}\Vert_{\,*_{\,Y}} \cdot \sum_{l\,=\,1}^\flat   \left(\ {1\over {     1\ + \ \Vert\,Y \,-\ \Xi_{\,l}\,\Vert  }}\ \right)^{\!\! {{n\,-\,2}\over 2 } \ +\  \tau_{\,>1} } \ {\bf ,} \ \ \ \ \ \ \ \ \ \ \ \ \ \\[0.2in]
 |\,{\bf H}\,(\,Y\,)\,| & \le &  \Vert\,{\bf H}\Vert_{\,**_{\,Y}} \cdot \sum_{l\,=\,1}^\flat   \left(\ {1\over {     1\ + \ \Vert\,Y \,-\ \Xi_{\,l}\,\Vert  }}\ \right)^{\!\! {{n\,+\,2}\over 2 } \ +\ \tau_{\,>1} }\ {\bf .}
\end{eqnarray*}
Note that\\[0.1in]
(A.4.27)
\begin{eqnarray*}
\left[\ {1\over {(\,1\ + \ R)^{{{\,n\,-\,2\,}\over 2} \,+\ \tau_{\,>1}}  }}\, \right]^{\,{{2n}\over {\ n\,-\,2\ }}  } & = & {1\over {(\,1\ + \ R)^{n \,+\ {{2n}\over {\,n\,-\ 2\,}} \,\cdot\, \tau_{\,>1}}  }} \ , \\[0.2in]
& \ &\hspace*{-0.7in}    \left\{\ {\bar C}_{\,2} \ := \  \int_{\R^n} {{d\,Y}\over {(\,1\ + \ R)^{\,n \,+\,{{2n}\over {\,n\,-\ 2\,}} \,\cdot\, \tau_{\,>1}}  }}  \ <  \ \infty  \ \right\} \ \ \ \ \ \ \ (\ R \ = \ \Vert \,Y\,\Vert\ )\ ,\\[0.2in]
\left[\ {1\over {(\,1\ + \ R)^{{{n\,+\,2}\over 2} \,+\ \tau_{\,>1}}  }}\, \right]^{\,{{2n}\over {\ n\,+\,2\ }}  } & = & {1\over {(\,1\ + \ R)^{n \,+\ {{2n}\over {n\,+\,2}} \,\cdot\, \tau_{\,>1}}  }} \ , \\[0.2in]
& \ &\hspace*{0.7in}    \left\{ \ {\bar C}_{\,3} \ := \  \int_{\R^n} {{d\,Y}\over {(\,1\ + \ R)^{\,n \,+\,{{2n}\over {n\,+\,2}} \,\cdot\, \tau_{\,>1}}  }}  \ <  \ \infty \ \right\}\,.\\
\end{eqnarray*}
Similar to (\,A.4.26\,)\,,\, we obtain (\ cf. {\bf \S\,A\,4\,.\,i}\,)
\begin{eqnarray*}
 (A.4.28) \ \ \ \ \ \ \ \  \ \ \  \int_{\R^n} \left[ \  \sum_{l\,=\,1}^\flat   \left(\ {1\over {     1\ + \ \Vert\,Y \,-\ \Xi_{\,l}\,\Vert  }}\ \right)^{\!\! {{n\,-\,2}\over 2 }   \ +\  \tau_{\,>1} }  \ \right]^{\,{{2n}\over {\,n\,-\ 2\,}} }
& = & {\bar C}_{\,2} \cdot  \flat  \ + \ o_{\,{\bar{\lambda}_{\,\flat}}}\,(\,1) \ {\bf ,} \ \ \ \ \ \ \ \  \ \  \ \ \ \ \ \ \ \  \ \   \\[0.2in]
(A.4.29)\ \ \ \ \ \ \ \ \ \ \
 \int_{\R^n} \left[ \  \sum_{l\,=\,1}^\flat   \left(\ {1\over {     1\ + \ \Vert\,Y \,-\ \Xi_{\,l}\,\Vert  }}\ \right)^{\!\! {{n\,+\,2}\over 2 }   \ +\  \tau_{\,>1} }  \ \right]^{\,{{2n}\over {n\,+\,2}} }
& = & {\bar C}_{\,3} \cdot \flat  \ + \ o_{\,{\bar{\lambda}_{\,\flat}}}\,(\,1) \ {\bf .}
\end{eqnarray*}
(\,A.4.26\,)\,,\, (\,A.4.27\,)\, and (\,A.4.28\,)\, hold under the conditions in Proposition {\bf A.4.20}\,.\, It follows that
\begin{eqnarray*}
\int_{\R^n} {\Phi}_\flat^{{2n}\over {\,n\,-\ 2\,}} & \le & C_4\cdot \flat \cdot \Vert\,{\Phi}\Vert_{\,*_{\,Y}}^{{n\,-\,2}\over {2\,n}}  \ ,\\[0.2in]
  \int_{\R^n} |\,{\bf H}\,|^{{2n}\over {\,n\,-\ 2\,}} & \le & C_5  \cdot \flat \cdot \Vert\,{\bf  H}\Vert_{\,**_{\,Y}}^{{2n}\over {\,n\,-\ 2\,}}\ .
\end{eqnarray*}
As a consequence\,,\, we obtain
\begin{eqnarray*}
  \left(\ \int_{\R^n} {\bf W}_\flat^{{2n}\over {\,n\,-\ 2\,}}\  \right)^{\!\!{2\over n}} \cdot \left(\ \int_{\R^n} |\,\Phi\,|^{{{2n}\over {\,n\,-\ 2\,}} } \right)^{\!\!\!{{n\,-\,2}\over {n}}}
& \le & C_6 \cdot  \left(\ \flat \ \right)^{{2\over n}} \cdot \left(\ \flat \cdot  \Vert\,\Phi\Vert_{\,*_Y}^{{2n}\over {\,n\,-\ 2\,}}\  \right)^{\!\!\!{{n\,-\,2}\over {n}}}\\[0.2in]
  & \le & C_6\cdot \flat^{2\over n} \cdot \flat^{{n\,-\,2}\over n} \cdot \Vert\,\Phi\Vert_{\,*_Y}^2 \ = \  C_6\cdot\flat \cdot \Vert\,\Phi\Vert_{\,*_Y}^2\,;\\[0.3in]
   \left(\ \int_{\R^n} |\,{\bf H}|^{{2n}\over {n\,+\,2}}\  \right)^{\!\!{{n\,+\,2}\over {2n}}} \cdot \left(\ \int_{\R^n} |\,\Phi\,|^{{{2n}\over {\,n\,-\ 2\,}} } \right)^{\!\!\!{{n\,-\,2}\over {2n}}}
& \le & C_7 \cdot  \left(\ \flat \cdot \Vert\,{\bf H}\Vert_{\,**_Y}^{{2\,n}\over {n\,+\,2}}  \ \right)^{{n\,+\,2}\over  {2\,n}} \cdot \left(\ \flat \cdot  \Vert\,\Phi\Vert_{\,*_Y}^{{2n}\over {\,n\,-\ 2\,}}\  \right)^{\!\!\!{{n\,-\,2}\over {2\,n}}}\\[0.2in]& \le & C_7\cdot \flat^{{n\,+\,2}\over {2n}} \cdot \flat^{{n\,-\,2}\over {2n}}\cdot \Vert\,{\bf H}\Vert_{\,**_Y}  \cdot \Vert\,\Phi\Vert_{\,*_{\,Y}}  \\[0.2in]
  & = & C_8\cdot \flat \cdot \Vert\,{\bf H}\Vert_{\,**_Y}  \cdot \Vert\,\Phi\Vert_{\,*_Y} \ .\\
\end{eqnarray*}

Combining with (A.4.24), we finally arrive at
\begin{eqnarray*}
\Vert\,\Phi\,\Vert_\btd^2 & \le & C_9 \cdot \flat \cdot \left(\,\Vert\,\Phi\Vert_{\,*_Y}^2 \ + \   \Vert\,{\bf H}\Vert_{\,**_Y}  \times \Vert\,\Phi\Vert_{\,*_Y} \,\right) \\[0.2in]
&  \le &  C_4 \cdot \flat \cdot \Vert\,{\bf H}\Vert_{\,**_Y}^2\ \ \ \ \ \ [\ {\mbox{via \ \ (A.4.23)}}\ ]\\[0.2in]
 & \  &  \hspace*{2in}\ \left[\ {\mbox{recall \ \ that \ \ }} {\bf H}\ = \ \left(\,V_{1\,,\  o}^{4\over {n\,-\,2}} \cdot {\cal H} \right)\ \right]\ .
\end{eqnarray*}
Thus we conclude the proof of (A.4.24)\,.

\newpage


{\bf \S\,A\,4\,.\,i.} \ \ {\bf Estimate of}\\[0.1in]
(\,A.4.30\,)
\begin{eqnarray*}
& \ &   [\ {\bf W}_\flat(\,Y\,) \ ]^{\,{{2n}\over {\,n\,-\ 2\,}}} \ - \ \ \sum_{l \,=\,1}^\flat\,  [\ {{\bf V}}_{l}\,(\,Y\,)\ ]^{\,{{2n}\over {\,n\,-\ 2\,}}} \\[0.2in]  &  = &  \left[\ \sum_{l \,=\,1}^\flat \left(\ {{\Lambda_{\,1}}\over {\Lambda^2_l \ + \ \Vert\,Y \ - \ \Xi_{\,l}\,\Vert^{\,2}}}\ \right)^{\!\!{{\,n\,-\,2\,}\over 2}} \ \right]^{ {{2n}\over {\,n\,-\ 2\,}}  } - \ \  \sum_{l \,=\,1}^\flat \left(\ {{\Lambda_{\,1}}\over {\ \Lambda^2_l \ + \ \Vert\,Y \ - \ \Xi_{\,l}\,\Vert^{\,2}\ }}\ \right)^{\!n}\ \ \ \ \ \ \ \ \ \ \\[0.1in]
& \ & \hspace*{4.5in}\ \ \mfor \ \ Y \,\in\,\R^n.
\end{eqnarray*}
Cf. {\bf \S\,4.\,b} of the original text, where the ``\,interaction term\,"\, is considered in detail\,.\, \{ \ The so called interaction term can be simply illustrated by
$$
(\,a \ + \ b\,)^{\,2} - \ (\,a^2 \ + \ b^2\,) \ = \ 2\,a\,b \ .
$$
\hspace*{4in}\   (  $\uparrow \ \ $ interaction.\,) \\[0.1in]
As the bubbles are far apart from each other (\,measured with the ``\,unit\," \, $\,{\bar\lambda}_{\,\,\flat}\,$)\,,\, the interaction term can be shown to be weak [\ cf. the Separation Lemma (\,A.3.13\,)\ ]\,.\,\,\}\bk
Recall that
 $$
\rho \ = \  {\bar\lambda}_{\,\,\flat}^\nu \ \ \ \ \ \ \ \ \  {\mbox{with}} \ \ \ \ \gamma \ + \ \nu \ >  \ 1\,.
$$
It satisfies
  (\,cf. Remark A.4.32\ )
$$
  \rho_\nu \ = \ o_{\,+}\,(\,1\,)\,\cdot\,\mbox{Min}  \ \{  \ \Vert \,\xi_{\,l} \ - \ \xi_{\,k} \Vert \ \ | \ \ 1 \ \le \ l\ \not= \ k\ \le \ \flat \ \}\ . \leqno (\,A.4.31\,)
$$
Here $\,\gamma\,$ appears in (\,1.22\,) of the original text.
%
%

\vspace*{0.1in}

{\it Remark.} A.4.32. \ \ Consider the case
\begin{eqnarray*}
 {{1}\over {  \ \Vert \,\Xi_{\,l} \ - \ \Xi_{\,k} \Vert \  }} & \approx &  {\bar\lambda}_{\,\,\flat}^{\  \gamma  } \ \ \ \ \mfor \ \ l \ \not= \ k\\[0.2in]
\Longleftrightarrow \ \  {1\over { \left(\ \ {{\Vert \,\xi_{\,l} \ - \ \xi_{\,k} \Vert}\over {\lambda }} \ \right)   }} & \approx &  {\bar\lambda}_{\,\,\flat}^{\  \gamma  }\\[0.2in]
\Longrightarrow \ \  \Vert \,\xi_{\,l} \ - \ \xi_{\,k} \Vert & \approx &  {\bar\lambda}_{\,\,\flat}^{\,1 \ - \ \,\gamma  }\\[0.2in]
(\,A.4.31\,) \ \  \Longrightarrow \ \  \nu \ > \ 1 \ - \ \gamma \ \ \Longleftrightarrow \ \ \gamma & > & 1 \ - \ \nu \ \ \Longleftrightarrow \ \ \nu \ + \ \gamma \ > \ 1\,.
\end{eqnarray*}
With this, we go on to the rescaled parameters introduced in {\bf \S\,A\,3\,.b}\,:
\begin{eqnarray*}
(A.4.33) \ \ \ \ \ \ \ \ \ & \ & B_{\,\Xi_{\,l\,}}(\, R_{\,\nu}\,)\ , \ \ \ \ \ \ {\mbox{where}} \ \ \ \ R_\nu \ := \ {\bar\lambda}_{\,\,\flat}^\nu \ \slash \  {\bar\lambda}_{\,\,\flat}   \ = \ {\bar\lambda}_{\,\,\flat}^{\ \nu \ - \ 1} \  \ \ \ \ \ \ \  \ \ \ \ \ \ \  \\[0.2in]
 \Longrightarrow & \ & R_{\,\nu} \  = \ o_{\,+}\,(\,1\,) \cdot\, \mbox{Min}  \ \left\{ \ \Vert \ \Xi_{\,l} \ - \ \Xi_{\,l\,'} \ \Vert\  \ \bigg\vert \ \ \  \ 1 \ \le \ l\, \not= \, \l\,' \ \le \ \flat\ \ \right\} \ .\ \ \ \ \ \ \ \ \ \ \ \ \ \ \ \ \ \
\end{eqnarray*}
In particular\,,\,
$$\ \ \  Y\,\in\,B_{\,\Xi_{\,1\,}}(\, R_{\,\nu}\,) \ \ \Longrightarrow \ \
\Vert\,Y\ - \ \Xi_{\,l}\,\Vert \ \ge \ {1\over 2} \cdot \Vert\,\Xi_{\,\,1}\ - \ \Xi_{\,l}\,\Vert \ \ \mfor \ \ l \ = \ 2\,, \ 3\,, \ \cdot \cdot \cdot\,, \ \flat\ . \leqno (\,A.4.34\,)
$$
In $\,B_{\,\Xi_{\,1\,}}(\, R_{\,\nu}\,)\,,\,$
\begin{eqnarray*}
(A.4.35) \ \ \ \ \ \ \ \ \   {\bf V}_1\,(\,Y\,) &  = & \left(\ {{\Lambda_{\,1}}\over { \Lambda_{\,1}^2 \ + \ \Vert\,Y\ - \ \Xi_{\,\,1}\,\Vert^{\,2}}}\  \right)^{\!\! {{\,n\,-\,2\,}\over 2}} \ \ge \  \left(\ {{\Lambda_{\,1}}\over {   \Lambda_{\,1}^2 \ + \  R_{\,\nu}^2}}\  \right)^{\!\! {{\,n\,-\,2\,}\over 2}} \\[0.2in]
& \ge & C_1 \cdot {\bar\lambda}_{\,\,\flat}^{\,(\,n\, -\ 2\,) \,\cdot\,(\,1\ - \ \nu\,) } \\[0.1in]
& \ &   \hspace*{-0.5in} [\ {\mbox{condition}} \ \  (\,1.8\,) \ \ \Longrightarrow \ \ \Lambda_{\,1}\ = \ O\,(\,1) \ ; \ \ \  {\mbox{\ \ for}} \ \  Y \,\in \,B_{\,\Xi_{\,\,1}}\,(\, R_{\,\nu}\,)\ ]\ .\ \ \ \ \ \ \ \ \ \ \ \ \ \ \ \ \ \
\end{eqnarray*}
As for the other bubbles\,,\, for $\,Y \,\in \,B_{\,\Xi_{\,\,1}}\,(\, R_{\,\nu}\,)\,,\,$
\begin{eqnarray*}
(A.4.36) \ \ \ \ \ \ \   & \ & {\bf V}_2\,(\,Y\,)  \ + \ \cdot \cdot \cdot \ + \ {\bf V}_\flat\,(\,Y\,)  \\[0.2in]
& = & \left(\ {{\Lambda_2}\over { \Lambda_2^2 \ + \ \Vert\,Y\ - \ \Xi_{\,2}\,\Vert^{\,2}}}\  \right)^{\!\! {{\,n\,-\,2\,}\over 2}}  \ + \ \cdot \cdot \cdot \ + \  \left(\ {{{\Lambda}_{\,\flat}}\over { {\Lambda}_{\,\flat}^2 \ + \ \Vert\,Y\ - \ \Xi_{\,\flat}\,\Vert^{\,2}}}\  \right)^{\!\! {{\,n\,-\,2\,}\over 2}} \\[0.2in]
& \le & C_2   \left[\  \left(\ {{1}\over { 1\ + \ \Vert\,\Xi_{\,\,1}\ - \ \Xi_{\,2}\,\Vert^{\,2}}}\  \right)^{\!\! {{\,n\,-\,2\,}\over 2}}  \ + \ \cdot \cdot \cdot \ + \  \left(\ {{1}\over { 1\ + \ \Vert\ \Xi_{\,\,1}\ - \ \Xi_{\,\flat}\,\Vert^{\,2}}}\  \right)^{\!\! {{\,n\,-\,2\,}\over 2}}  \ \right] \\[0.2in]
& \le & C_3  \cdot  {\bar\lambda}_{\,\,\flat}^{\,(n\,-\,2\,)\,\cdot\, \gamma} \ \ \ \  \ \   \ \ \ \ \ \ \  \ \   \ \ \  \ \ \ \  \ \   \ \ \  [ \ \uparrow \ \ {\mbox{using \ \ (\,A.4.34\,)}}\ ]\\[0.2in]
\Longrightarrow\ \ \ & \ & \!\!\!\!\!\!\!\!\!\!{{{\bf V}_2\,(\,Y\,)  \ + \ \cdot \cdot \cdot \ + \ {\bf V}_\flat\,(\,Y\,)   }\over {{\bf V}_1\,(\,Y\,)  }} \ \ = \ O\,\left(\, {\bar\lambda}_{\,\,\flat}^{\,(\,n\,-\,2\,)\,\cdot\, [\ (\ \gamma \ + \ \nu\ ) \ - \ 1\ ]}\ \ \right)\\[0.1in]
& \ &  \hspace*{1in}
 [  \ {\mbox{recall \  \ that}} \ \ \gamma \ + \ \nu \ > \ 1\ ; \ \ {\mbox{hold \ \ for}} \ \  Y \,\in \,B_{\,\Xi_{\,\,1}}\,(\, R_{\,\nu}\,)\ ]\ .
\end{eqnarray*}


\newpage

{\it Focusing on the leading term.} \ \ For $\,n \ \ge\ 6\,,\,$ applying Taylor expansion\,,\, we obtain
\begin{eqnarray*}
& \ &
(\,1 \ + \ t\,)^{{2n}\over {\,n\,-\ 2\,}} \ \le \ 1 \ + \ \left(\ {{2n}\over {\,n\,-\ 2\,}}\ \right) \cdot t \ +\ \left( \  {1\over 2} \cdot {{2n}\over {\,n\,-\ 2\,}} \cdot {{ n\,+\,2}\over {n\,-\,2}} \ \right) \cdot t^2 \ + \ C_1  \cdot t^{{2n}\over {\,n\,-\ 2\,}} \\[0.2in]
& \ & \hspace*{4in} \ \ \  \mfor \ \ 0 \ \le \ t\  \le \ 1\,.
\end{eqnarray*}
Here the positive constant $\,C_1\,$ depends on the dimnsion $\,n\,$ only\,.\,  When $\,t\ \ge \ 1\,,\,$ we have
$$
(\,1 \ + \ t\,)^{{2n}\over {\,n\,-\ 2\,}} \ \le  \ (\,t \ + \ t\,)^{{2n}\over {\,n\,-\ 2\,}}  \ \le \ 2^{{2n}\over {\,n\,-\ 2\,}}  \cdot t^{{2n}\over {\,n\,-\ 2\,}} \ \ \ \ \ \ \ \ {\mbox{for}} \ \ t \ \ge   \ 1\,. \leqno (A.4.37)
$$
Hence\\[0.1in]
(A.4.38)
 \begin{eqnarray*}
& \ &\!\!\!\!\!\!\!\!\ \ \
(\,1 \ + \ t\,)^{{2n}\over {\,n\,-\ 2\,}} \ \le \ 1 \ + \ \left(\ {{2n}\over {\,n\,-\ 2\,}}\ \right) \cdot t \ +\ \left( \  {1\over 2} \cdot {{2n}\over {\,n\,-\ 2\,}} \cdot {{ n\,+\,2}\over {n\,-\,2}} \ \right) \cdot t^2 \ + \ C_2  \cdot t^{{2n}\over {\,n\,-\ 2\,}} \\[0.2in]
& \ & \hspace*{4in} \ \ \  \mfor \ \ {\mbox{all}} \ \   t\  \ge \ 0\,.
\end{eqnarray*}

\vspace*{-0.35in}

Here
$$
C_2 \ = \ \max \ \left\{ \ 2^{{2n}\over {\,n\,-\ 2\,}} \ , \ \ C_1 \ \right\}\ .
$$
Accordingly\,,\, inside $\,B_{\,\Xi_{\,\,1}}\,(\, R_{\,\nu}\,)\,,\,$
\begin{eqnarray*}
(A.4.39)  & \   & \left[ \  \left(\ \, {{\bf V}}_1 \,+\, {{\bf V}}_2\,+\,\cdot \cdot \cdot\,+\, {{\bf V}}_\flat\right)\,\bigg\vert_{\,Y} \ \right]^{\,{{2n}\over{n\,-\,2}} }  \\[0.15in]
&=  &  [\ {{\bf V}}_1\,(\,Y\,)\ ]^{{2n}\over{n\,-\,2}} \ \cdot \  \left(\ 1\ + \  {{\  {{\bf V}}_2  \,+\, {{\bf V}}_3 \,+\,\cdot \cdot \cdot\,+\, {{\bf V}}_{\,\flat}  \ }\over {{\bf V}_1}}\ \Bigg\vert_{\,Y}\  \right)^{{2n}\over{n\,-\,2}}\\[0.2in]
&\le & [\ {{\bf V}}_1\,(\,Y\,)\ ]^{{2n}\over{n\,-\,2}} \ \cdot \  \left[\ 1\ + \  \left(\  {{2n}\over{n\,-\,2}}\ \right)\cdot \left(\  {{\  {{\bf V}}_2  \,+\, {{\bf V}}_3 \,+\,\cdot \cdot \cdot\,+\, {{\bf V}}_{\,\flat}  \ }\over {{\bf V}_1}}\ \Bigg\vert_{\,Y}\  \right) \ \right.\\[0.2in]
& \ & \ \ \ \ \  + \ \left( \  {1\over 2} \cdot {{2n}\over {\,n\,-\ 2\,}} \cdot {{ n\,+\,2}\over {n\,-\,2}} \ \right) \cdot  \left(\  {{\  {{\bf V}}_2  \,+\, {{\bf V}}_3 \,+\,\cdot \cdot \cdot\,+\, {{\bf V}}_{\,\flat}  \ }\over {{\bf V}_1}}\ \Bigg\vert_{\,Y}\  \right)^{2}\\[0.2in]
& \ &  \ \ \ \ \ \ \  \ \ \ \ \  \left. + \ \ C_2 \cdot  \left(\  {{\  {{\bf V}}_2  \,+\, {{\bf V}}_3 \,+\,\cdot \cdot \cdot\,+\, {{\bf V}}_{\,\flat}  \ }\over {{\bf V}_1}}\ \Bigg\vert_{\,Y}\  \right)^{{2n}\over {n\, - \, 2}} \ \right] \\[0.15in]
& \ & \hspace*{3in} \ \ \ \ \  \ \ \ \ \   \ \ \ \ \   {\mbox{for}} \ \ Y \,\in \,B_{\,\Xi_{\,\,1}}\,(\, R_{\,\nu}\,)\ .
\end{eqnarray*}

\newpage

It follows that\,,\,  for $\,Y \,\in \,B_{\,\Xi_{\,\,1}}\,(\, R_{\,\nu}\,)\,,\,$  we have  \\[0.1in]
(A.4.40)
\begin{eqnarray*}
0 & <  &  \left[ \  \left(\ \, {{\bf V}}_1 \,+\, {{\bf V}}_2\,+\,\cdot \cdot \cdot\,+\, {{\bf V}}_\flat\right)\,\bigg\vert_{\,Y} \ \right]^{\,{{2n}\over{n\,-\,2}} }  \ -\ [\ {{\bf V}}_1\,(\,Y\,)\ ]^{{2n}\over{n\,-\,2}} \\[0.25in]
& \le & [\ {{\bf V}}_1\,(\,Y\,)\ ]^{{2n}\over{n\,-\,2}} \ + \  \left(\ {{2n}\over{n\,-\,2}}\ \right)\cdot [\ {{\bf V}}_1\,(\,Y\,)\ ]^{{n\,+\,2}\over{n\,-\,2}} \cdot \left(\  {{\bf V}}_2  \,+\, {{\bf V}}_3 \,+\,\cdot \cdot \cdot\,+\, {{\bf V}}_{\,\flat} \right)\,\bigg\vert_{\,Y} \\[0.2in]
& \ & \ \ \ \ \ \ + \ \left( \  {1\over 2} \cdot {{2n}\over {\,n\,-\ 2\,}} \cdot {{ n\,+\,2}\over {n\,-\,2}} \ \right) \cdot [\ {{\bf V}}_1\,(\,Y\,)\ ]^{{4}\over{n\,-\,2}} \cdot \left[ \  \left(\  {{\bf V}}_2  \,+\, {{\bf V}}_3 \,+\,\cdot \cdot \cdot\,+\, {{\bf V}}_{\,\flat} \right)\,\bigg\vert_{\,Y} \ \right]^{2}\\[0.2in]
& \ & \ \ \ \ \ \ \ \ \ \  \ \  + \ C_2 \cdot \left[ \  \left(\  {{\bf V}}_2  \,+\, {{\bf V}}_3 \,+\,\cdot \cdot \cdot\,+\, {{\bf V}}_{\,\flat} \right)\,\bigg\vert_{\,Y} \ \right]^{{2n}\over{n\,-\,2}} \ \ \ \   \ \ \ \ \ \ \ \  \ \ \ \   \ \ \ \ \ \ \ \  [\,{\mbox{via \ \ (\,A.4.39\,)\,}}] \\[0.15in]
& \  & \ \  \ \ \ \ \ \ \ \ \ \ \ \ \ \ \ \ \ \ \ \ \ \  \ -\ [\ {{\bf V}}_1\,(\,Y\,)\ ]^{{2n}\over{n\,-\,2}}\\[0.25in]
& \le &  \left(\  {{2n}\over{n\,-\,2}} \ \right) \cdot  [\ {{\bf V}}_1\,(\,Y\,)\ ]^{{n\,+,2}\over{n\,-\,2}} \cdot \left(\  \ {{\bf V}}_2  \,+\, {{\bf V}}_3 \,+\,\cdot \cdot \cdot\,+\, {{\bf V}}_{\,\flat} \ \right)\,\bigg\vert_{\,Y}  \\[0.2in]
& \ & \ \ \ \ \ \ + \ \left(\ {1\over 2} \cdot {{2n}\over {\,n\,-\ 2\,}} \cdot {{ n\,+\,2}\over {n\,-\,2}}\ \right) \cdot  [\ {{\bf V}}_1\,(\,Y\,)\ ]^{{4}\over{n\,-\,2}} \cdot   \left[ \  \left(\  {{\bf V}}_2  \,+\, {{\bf V}}_3 \,+\,\cdot \cdot \cdot\,+\, {{\bf V}}_{\,\flat} \right)\,\bigg\vert_{\,Y} \ \right]^{2}\\[0.2in]
& \ & \ \ \ \ \  \ \ \ \ \ \ \ \ \  \ \ \ \ \  \  \ + \ \ C_2 \cdot   \left[ \  \left(\  {{\bf V}}_2  \,+\, {{\bf V}}_3 \,+\,\cdot \cdot \cdot\,+\, {{\bf V}}_{\,\flat} \right)\,\bigg\vert_{\,Y} \ \right]^{{2n}\over{n\,-\,2}} \\[0.25in]
& \le & C_1 \!\cdot \left(\  {1\over {1 \ + \ \Vert\, Y \ - \ \Xi_{\,\,1}\,\Vert}}\ \right)^{{n\,+\,2}\over 2} \cdot \sum_{l\,=\,2}^\flat \left(\ {1\over {1 \ + \ \Vert\, Y \ - \ \Xi_j\,\Vert}}\ \right)^{n\,-\,2}\  \ \ \cdot\,\cdot\,\cdot\,\cdot\,\cdot\,\cdot\,\cdot\,\cdot \,\cdot  \ \ (\ {\bf A}\ )_{\,(\,A.4.40\,)} \\[0.2in]
& \ &\!\!\!\!\!\!\!\! +  \  \left( \ {1\over 2} \cdot {{2n}\over {\,n\,-\ 2\,}} \cdot {{ n\,+\,2}\over {n\,-\,2}} \ \right) \cdot  [\ {{\bf V}}_1\,(\,Y\,)\ ]^{{4}\over{n\,-\,2}} \cdot   \left[ \  \left(\  {{\bf V}}_2  \,+\, {{\bf V}}_3 \,+\,\cdot \cdot \cdot\,+\, {{\bf V}}_{\,\flat} \right)\,\bigg\vert_{\,Y} \ \right]^{2}\,\cdot\,\cdot\,\cdot\ \ (\ {\bf B}\ )_{\,(\,A.4.40\,)} \\[0.2in]
& \ & \ \ \ \ \ \ \ \ \ \ \ \  \ \ \ \ \ \ \ \ \ \ \ \   \ \ \ \ \ \,   + \ C_2 \cdot\! \left[\  \sum_{l\,=\,2}^\flat \left(\ {1\over {1 \ + \ \Vert\, Y \ - \ \Xi_j\,\Vert}}\ \right)^{n\,-\,2} \ \right]^{\,{{2n}\over {\,n\,-\ 2\,}}} \ . \,\cdot\,\cdot\,\cdot\,\cdot\,\cdot\,\cdot \ \ (\ {\bf C}\ )_{\,(\,A.4.40\,)} \\[0.15in]
& \ & \hspace*{3.5in} [ \ {\mbox{holds \ \ for }} \ \  Y \,\in \,B_{\,\Xi_{\,\,1}}\,(\, R_{\,\nu}\,)\ ]\ .
\end{eqnarray*}

\newpage




\underline{$(\ {\bf A}\ )_{\,(\,A.4.40\,)} $}\,.\,  \ \ Via the Separation Lemma [\,see (\,A.3.13\,)\ ]\,,\, we proceed with\\[0.1in]
(A.4.41)
\begin{eqnarray*}
& \  &  \left(\ {1\over {1 \ + \ \Vert\, Y \ - \ \Xi_{\,\,1}\,\Vert}}\ \right)^{{n\,+\,2}\over 2} \cdot \sum_{l\,=\,2}^\flat \left(\ {1\over {1 \ + \ \Vert\, Y \ - \ \Xi_j\,\Vert}}\ \right)^{n\,-\,2}\\[0.15in]
& \le & C_1 \cdot  \sum_{l\,=\,2}^\flat \left(\ {1\over {(\,1 \ + \ \Vert\, Y \ - \ \Xi_{\,\,1}\,\Vert\,)^{\ n \,+\,\varepsilon}}} \ + \ {1\over {(\,1 \ + \ \Vert\, Y \ - \ \Xi_{\,l}\,\Vert\,)^{\ n \,+\,\varepsilon}}}\  \right) \cdot {1\over { \Vert\,\Xi_{\,l}\,-\,\Xi_{\,\,1}\,\Vert^{\,{{\,n\,-\,2\,}\over 2} \,-\,\varepsilon} }}\\[0.15in]
& \le & 2 \cdot  C_1 \cdot {1\over {(\,1 \ + \ \Vert\, Y \ - \ \Xi_{\,\,1}\,\Vert\,)^{\ n \,+\,\varepsilon} }}\, \cdot\  \sum_{l\ =\ 2}^\flat {1\over { \Vert\,\Xi_{\,l}\,-\,\Xi_{\,\,1}\,\Vert^{ \,{{ n\,-\,2 }\over { 2}}\,-\,\varepsilon } }} \\[0.1in]
& \ & \hspace*{3.5in} \ \ \ \ \ \ \  [ \ {\mbox{via \ \ (\,A.4.33\,) \ \& \ (\,A.4.34\,)}} \ ]  \\[0.15in]
& \le & C_2 \cdot {1\over {(\,1 \ + \ \Vert\, Y \ - \ \Xi_{\,\,1}\,\Vert\,)^{\,n\,+\,\varepsilon}}}\  \cdot \, {\bar\lambda}_{\,\,\flat}^{\,\left( {{\,n\,-\,2\,}\over 2} \,-\,\varepsilon \right)\,\cdot\,\gamma}\ .\\
\end{eqnarray*}
The index $\,\varepsilon \ > \ 0\,$ is inserted to make the integration over $\,\R^n\,$ finite\,,\, as
$$
\int_{\R^n}  {{d\,Y}\over {(\,1 \ + \ \Vert\, Y \,\Vert\,)^{\,n}}}\ = \ \infty\,,\, \ \ \ \ {\mbox{but}} \ \ \  \int_{\R^n}  {{d\,Y}\over {(\,1 \ + \ \Vert\, Y \,\Vert\,)^{\,n\,+\,\varepsilon}}} \ = \  {{ C_3\,(\,n\,) }\over {\varepsilon\,}} \ <  \ \infty\ . \leqno (A.4.42)
$$

\vspace*{0.3in}

\underline{$(\ {\bf C}\ )_{\,(\,A.4.40\,)} $}\,. \ \ We proceed to estimate $(\ {\bf C}\ )_{\,(\,A.4.40\,)} $\, first before $(\ {\bf B}\ )_{\,(\,A.4.40\,)} $\,.\,  For  $\,Y \,\in \,B_{\,\Xi_{\,\,1}}\,(\, R_{\,\nu}\,)\,$ and $\,l \, \not= \, 1\,,\,$ we have\\[0.1in]
(A.4.43)
\begin{eqnarray*}
& \ & \\[-0,3in]
& \ & {1\over {(\,1 \ + \ \Vert\, Y \ - \ \Xi_{\,l}\,\Vert\,)^{\,n\,-\,2}}} \\[0.15in]
  & \le & {1\over {(\,1 \ + \ \Vert\, Y \ - \ \Xi_{\,\,1}\,\Vert\,)^{ {{\,n\,-\,2\,}\over 2}  }}} \cdot {1\over {(\,1 \ + \ \Vert\, Y \ - \ \Xi_{\,l}\,\Vert\,)^{ {{\,n\,-\,2\,}\over 2}  }}}   \\[0.15in]
& \le & C \cdot {1\over {\Vert\,\Xi_{\,\,1}\ - \ \Xi_{\,l}\,\Vert^{ \, \left(\ {{\,n\,-\,2\,}\over 2} \,-\,\varepsilon \,\right)  \,\cdot \,\gamma } }} \cdot  {1\over {(\,1 \ + \ \Vert\, Y \ - \ \Xi_{\,\,1}\,\Vert\,)^{ {{\,n\,-\,2\,}\over 2} \,+\,\varepsilon }}} \\[0.15in]
  &\ & \hspace*{-0.7in}   \left[\,{\mbox{argue \ \ as \ \ in \ \ (A.4.41)}}\,, \ \ \ {\mbox{via \ \ (\,A.4..34\,)\ : }} \ \  \  {1\over {(\,1 \ + \ \Vert\, Y \ - \ \Xi_{\,\,1}\,\Vert\ )}} \ \ge \  {1\over {(\,1 \ + \ \Vert\, Y \ - \ \Xi_{\,l}\,\Vert\ )}}  \ \right]
  \end{eqnarray*}

  \newpage

  \begin{eqnarray*}
  \Longrightarrow & \ &  \sum_{l\,=\,2}^\flat \ {1\over {(\,1 \ + \ \Vert\, Y \ - \ \Xi_{\,l}\,\Vert\,)^{\,n\,-\,2}}} \\[0.15in]
   & \le &   C_4 \cdot \left(\ \ \sum_{l\,=\,2}^\flat\  {1\over {\ \Vert\,\Xi_{\,\,1}\ - \ \Xi_{\,l}\,\Vert^{ \left(\ \,{{\,n\,-\,2\,}\over 2} \,-\,\varepsilon \,\right)  \,\cdot \,\gamma }\  }}  \right) \cdot  {1\over {(\,1 \ + \ \Vert\, Y \ - \ \Xi_{\,\,1}\,\Vert\,)^{ {{\,n\,-\,2\,}\over 2} \,+\,\varepsilon }}}\\[0.15in]
  & \ \le  & C_5 \cdot {\bar\lambda}_{\,\,\flat}^{\left( {{\,n\,-\,2\,}\over 2} \,-\,\varepsilon \right)\,\cdot\,\gamma} \cdot\, {1\over {(\,1 \ + \ \Vert\, Y \ - \ \Xi_{\,\,1}\,\Vert\,)^{ {{\,n\,-\,2\,}\over 2} \,+\,\varepsilon }}}\\[0.2in]
  \Longrightarrow & \ & \left(\ \ \sum_{l\,=\,2}^j \ {1\over {(\ 1 \ + \ \Vert\, Y \ - \ \Xi_{\,l}\,\Vert\ )^{\,n\,-\,2}}}\  \right)^{{2n}\over {\,n\,-\ 2\,}} \, \le \,  C_6 \cdot {\bar\lambda}_{\,\,\flat}^{\ \left(\, n\,-\,\varepsilon' \right)\,\cdot\,\gamma} \,\cdot\,  {1\over {\ (\,1 \ + \ \Vert\, Y \ - \ \Xi_{\,\,1}\,\Vert\,)^{\, n\,+\,\varepsilon' }\ }}\ .
\end{eqnarray*}
Here
$$
\varepsilon' \ = \ {{2n}\over {\,n\,-\ 2\,}} \cdot\varepsilon \ .
$$

\vspace*{0.1in}

\underline{$(\ {\bf B}\ )_{\,(\,A.4.40\,)} $}\,. \ \ Likewise [\,the following expressions are to be evaluated at $\,Y \ \in \ B_{\,\Xi_{\,l}}\,(\,R_{\,\nu}\,)\,$]\,,\\[0.1in]
(A.4.44)
\begin{eqnarray*}
& \ &  \left(\ \, {{\bf V}}_2  \,+\, {{\bf V}}_3 \,+\,\cdot \cdot \cdot\,+\, {{\bf V}}_{\,\flat} \right)^{2}\\[0.2in]
& \le  &  {{\bf V}}_2   \cdot \left(\ \, {{\bf V}}_2  \,+\, {{\bf V}}_3 \,+\,\cdot \cdot \cdot\,+\, {{\bf V}}_{\,\flat} \right) \ + \  {{\bf V}}_3   \cdot \left(\ \, {{\bf V}}_2  \,+\, {{\bf V}}_3 \,+\,\cdot \cdot \cdot\,+\, {{\bf V}}_{\,\flat} \right) \ +  \\[0.2in]
& \ & \ \ + \ \cdot \cdot \cdot   \ + \  {{\bf V}}_{\,\flat}   \cdot \left(\ \, {{\bf V}}_2  \,+\, {{\bf V}}_3 \,+\,\cdot \cdot \cdot\,+\, {{\bf V}}_{\,\flat} \right)\\[0.2in]
& \le  &  {{\bf V}}_1   \cdot \left(\ \, {{\bf V}}_2  \,+\, {{\bf V}}_3 \,+\,\cdot \cdot \cdot\,+\, {{\bf V}}_{\,\flat} \right)  \ + \  {{\bf V}}_1   \cdot \left(\ \, {{\bf V}}_2  \,+\, {{\bf V}}_3 \,+\,\cdot \cdot \cdot\,+\, {{\bf V}}_{\,\flat} \right) \ + \ \\[0.2in]
& \ & \ \ \ + \ \cdot \cdot \cdot   \ + \  {{\bf V}}_1   \cdot \left(\ \, {{\bf V}}_2  \,+\, {{\bf V}}_3 \,+\,\cdot \cdot \cdot\,+\, {{\bf V}}_{\,\flat} \right)\ \ \ \ \ \ \ \ \ \ \ \ \ \ \ \ [\  {\mbox{argue \ \ as \ \ in \ \ (\,A.4.43\,)}} \ ] \\[0.2in]
\Longrightarrow  & \ & {{\bf V}}_1^{4\over {n\,-\,2}} \cdot \left(\ \, {{\bf V}}_2  \,+\, {{\bf V}}_3 \,+\,\cdot \cdot \cdot\,+\, {{\bf V}}_{\,\flat} \right)^{2} \ \bigg\vert_{\,Y} \\[0.15in]
& \le & \flat \cdot  {{\bf V}}_1^{4\over {n\,-\,2}} \cdot {{\bf V}}_1 \cdot  \left(\ \, {{\bf V}}_2  \,+\, {{\bf V}}_3 \,+\,\cdot \cdot \cdot\,+\, {{\bf V}}_{\,\flat} \right) \ \bigg\vert_{\,Y}\\[0.2in]
&    \le  & \ \ C_1 \cdot \flat \cdot {\bar\lambda}_{\,\,\flat}^{\,(\ n\,-\ \varepsilon\,)\,\cdot\, \gamma} \,\cdot\,  {1\over {\ (\,1 \ + \ \Vert\, Y \ - \ \Xi_{\,\,1}\,\Vert\,)^{\, n\,+\,\varepsilon' }\ }} \\[0.15in] & \le &   C_2 \cdot {\bar\lambda}_{\,\,\flat}^{\,(\ n\,-\ \varepsilon\,)\,\cdot\, \gamma \ - \ \sigma} \,\cdot\,  {1\over {\ (\,1 \ + \ \Vert\, Y \ - \ \Xi_{\,\,1}\,\Vert\,)^{\, n\,+\,\varepsilon }\ }} \\[0.1in]
& \ &     \bigg[\ {\mbox{argue \ \ as \ \ in \ \ (A.4.41)}}\ , \ \ {\mbox{using \ \ condition  \ \ (\,1.32\,)}}\ : \ {\mbox{namely}}\,, \ \ \ 2 \ \le \ \flat \ \le\  {1\over { { {\bar\lambda}_{\,\,\flat}   }^{\,\sigma} }}  \  \bigg]\ .\\
\end{eqnarray*}
{\it Combined\,.}\, Consider the following order (\,note that $\,\varepsilon\,$ and hence $\,\varepsilon'\,$ can be made as small as we like\,)\,:
 $$
 \left(\, {{\,n\,-\,2\,}\over 2} \,-\,\varepsilon\, \right)\,\cdot\,\gamma \ \le \ \,(\ n\,-\ \varepsilon\,)\,\cdot\, \gamma \ - \ \sigma \  \le \ (\ n\,-\ \varepsilon'\,)\,\cdot\, \gamma
$$
\hspace*{2.6in}$\uparrow \ \ \gamma \ >  \ \sigma \ \ \Longrightarrow \ \
 \left(\, {{n\,+\,2}\over 2}\,\right) \cdot \gamma \ > \ \sigma$\\[0.1in]
 for $\,\varepsilon\,$ to be small enough\,.\, Combining the estimates in $(\ {\bf A}\ )_{\,(\,A.4.40\,)} $\,, $(\ {\bf B}\ )_{\,(\,A.4.40\,)} $\, and $(\ {\bf C}\ )_{\,(\,A.4.40\,)} $\,,\, we obtain \\[0.1in]
 (\,A.4.45\,)
$$
 \left[\  \left(\  {\bf V}_1  \,+\, {{\bf V}}_2 \,+\,\cdot \cdot \cdot\,+\, {{\bf V}}_\flat \ \right)\,\bigg\vert_{\,Y} \ \right]^{\ {{2n}\over{n\,-\,2}} }\\[0.2in]
\ = \ [ \ {\bf V}_1\,(\,Y\,)\ \,]^{{2n}\over{n\,-\,2}}  \ \ + \ \ {\bf E}_{\,(A.4.45)\,} \ \ \ \ \mfor \ \ Y \ \in \  B_{\,\Xi_{\,l}}\,(\,R_{\,\nu}\,)\,(\,Y\,)\,,
$$
where
$$
|\, {\bf E}_{\,(A.4.45)\,}\,(\,Y\,)\,| \ \le   {{C\,(\,n\,)}\over {\ (\,1 \ + \ \Vert\, Y \ - \ \Xi_{\,\,1}\,\Vert\,)^{\ n\,+\,\varepsilon }\ }} \ \cdot \ {\bar\lambda}_{\,\,\flat}^{\,\left( {{\,n\,-\,2\,}\over 2} \,-\,o_{\,+}\,(\,1\,) \right)\,\cdot\,\gamma}\  \ \ \ \mfor \ \ Y \ \in \  B_{\,\Xi_{\,l}}\,(\,R_{\,\nu}\,)\,(\,Y\,)\ .
$$
Here $\,o_{\,+}\,(\,1\,)\,$ is an arbitrary but fixed small positive number, and the factor $\,1\,/\,\varepsilon\,$ in (\,A.4.42\,)  is absorbed in $\,o_{\,+}\,(\,1\,)\,$ with $\,{\bar\lambda}_{\,\,\flat} \ > \ 0$\, close to zero\,.

\vspace*{0.3in}


{\bf \S\,A\,4\,.\,j .} \ \ {\bf Proof of (\,A.4.26\,)\,. }\\[0.1in]
Upon obtaining similar estimates as in (\,A.4.45\,) insider $\,B_{\,\Xi_{\,l}}\,(\,R_{\,\nu}\,)\,(\,Y\,)\,$,\, $\,l \ = \ 2\,,\,\ 3\,, \ \cdot \cdot \cdot\,,\,\ \flat\,,\,$ and performing the integration
$$
 \int_{\R^n}
{\bf W}_\flat^{{2n}\over {\,n\,-\ 2\,}} \ = \  {\int}_{\displaystyle{ \bigcup_{l\,=\,1}^\flat B_{\,\Xi_{\,l}}\,(\,R_{\,\nu}\,) }}\hspace*{-0.6in}
{\bf W}_\flat^{{2n}\over {\,n\,-\ 2\,}}  \ \ \ \  + \ \ \ \  \int_{\displaystyle{\R^n \, \Big\backslash \,\bigcup_{l\,=\,1}^\flat B_{\,\Xi_{\,l}}\,(\,R_{\,\nu}\,) }}\hspace*{-0.6in}
{\bf W}_\flat^{{2n}\over {\,n\,-\ 2\,}} \ \ \ ,
$$

 \newpage

we first obtain\\[0.1in]
(A.4.46)
\begin{eqnarray*}
& \ & {\int}_{\displaystyle{  \bigcup_{l\,=\,1}^\flat B_{\,\Xi_{\,l}}\,(\,R_{\,\nu}\,) }}\hspace*{-0.6in}
{\bf W}_\flat^{{2n}\over {\,n\,-\ 2\,}}  \ \ \ = \  \sum_{l\,=\,1}^\flat \ \int_{B_{\,\Xi_{\,l}}\,(\,R_{\,\nu}\,)}
{\bf W}_\flat^{{2n}\over {\,n\,-\ 2\,}} \\[0.1in]
& \ & \ \ \ \ \ \ \ \ \ \ \ \ \ \ \ \ \ \ \ \ \ \ \ \ \ \ \ \ \ \ \ \  \ \ \ \ \ \ \ \ \ \ \ \ \ \ \ \ \left[ \ B_{\,\Xi_{\,l}}\,(\,R_{\,\nu}\,) \ \cap \ B_{\,\Xi_{\,l}}\,(\,R_{\,\nu}\,) \ = \ \emptyset \ \ \ \ \mfor  \  l \ \not= \ k  \ \right]\\[0.1in]
&  = &     \sum_{l\,=\,1}^\flat \ \int_{B_{\,\Xi_{\,l}}\,(\,R_{\,\nu}\,)}
  {{\bf V}}_l^{{2n}\over{n\,-\,2}}  \ + \ C\,(\,n\,) \cdot \flat \cdot  {\bar\lambda}_{\,\,\flat}^{\,\left( {{\,n\,-\,2\,}\over 2} \,-\,o_{\,+}\,(\,1\,)  \right)\,\cdot\,\gamma}  \\[0.2in]
& =  &     \sum_{l\,=\,1}^\flat \ \int_{\R^n}
  {{\bf V}}_l^{{2n}\over{n\,-\,2}}  \ + \ \left[ \ O \left(\  {\bar\lambda}_{\,\,\flat}^{\,\left( {{\,n\,-\,2\,}\over 2} \,-\,o_{\,+}\,(\,1\,) \right)\,\cdot\,\gamma }  \ \right) \ + \ O \left(\  {\bar\lambda}_{\,\,\flat}^{ \ n\, \cdot\, \left( \ 1 \ - \ \nu \ \right) } \ \right) \ \right] \cdot \flat \\[0.1in]
& \  &    \hspace*{2.5in} \left[ \ \int_{\R^n\ \setminus \, B_{\,\Xi_{\,l}}\,(\,R_{\,\nu}\,) }
  {{\bf V}}_l^{{2n}\over{n\,-\,2}} \ = \  O \left(\  {\bar\lambda}_{\,\,\flat}^{\ n \,\cdot \,\left( \ 1 \ - \ \nu \ \right) } \ \right) \ \right]\ \\[0.1in]
   & = & \flat \cdot \ \int_{\R^n}
{\bf V}_o^{{2n}\over {\,n\,-\ 2\,}} \  + \  \ + \ O \left(\  {\bar\lambda}_{\,\,\flat}^{\,\left( {{\,n\,-\,2\,}\over 2} \,-\,o_{\,+}\,(\,1\,) \right)\,\cdot\,\gamma \ - \ \sigma } \ \right)  \ + \ O \left(\  {\bar\lambda}_{\,\,\flat}^{ \ n\, \cdot\, \left( \ 1 \ - \ \nu \ \right)  \ - \ \sigma  } \ \right) \ .
\end{eqnarray*}
Recall
$$
{\bf V}_o \,(\,Y\,) \ = \ \left(\  {1\over {1\ + \ \Vert\,Y\Vert^{\,2}}}\  \right)^{\!\! {{\,n\,-\,2\,}\over 2}}   \ .
$$

\vspace*{0.2in}

{\it Estimate on}
$$
\int_{\displaystyle{\R^n \,\Big\backslash \,\bigcup_{l\,=\,1}^\flat B_{\,\Xi_{\,l}}\,(\,R_{\,\nu}\,) }}\hspace*{-0.6in}
{\bf W}_\flat^{{2n}\over {\,n\,-\ 2\,}} \ .\leqno (\,A.4.47\,)
$$
We have
$$
 {\bf W}_\flat^{{2n}\over {\,n\,-\ 2\,}} \ = \ {\bf W}_{\,\flat}  \cdot {\bf W}_\flat^{{\,n\,+\,2\,}\over {n\,-\,2}}\ \ \ \ {\mbox{and}} \ \ \ \
n \ \ge \ 6 \ \ \Longrightarrow \  \ {1\over 2} \cdot {{n\,+\,2}\over {\,n\,-\ 2\,}} \ \le \ 1\,. \leqno (\,A.4.4\,8)
$$
One can estimate (\,A.4.47\,) based on (\,A.4.48\,) and the argument as in {\bf \S\,A\,4\,.i}\,.\,
As $\, {{n\,+\,2}\over {\,n\,-\ 2\,}}\,$ may not be an integer\,,\,  the following consideration gives a cleaner picture [\,cf. the summation in (\,A.4.56\,)\,]\,.\, It begins with (\,the following expressions are to be evaluated at $\,Y\,\in\,\R^n$\,)\\[0.1in]

\newpage

(\,A.4.49\,)
\begin{eqnarray*}
\left[\  {\bf W}_{\,\flat}  \ \right]^{{\,n\,+\,2\,}\over {n\,-\,2}}  & = &   \left[\  {\bf V}_1\ + \  {\bf V}_2  \ +  \cdot \cdot \cdot  + \ {\bf V}_{\,\flat}  \ \right]^{\,{1\over 2} \cdot {{n\,+\,2}\over {\,n\,-\ 2\,}}} \ \cdot   \left[\  {\bf V}_1\ + \  {\bf V}_2  \ +  \cdot \cdot \cdot  + \ {\bf V}_{\,\flat}  \ \right]^{\,{1\over 2} \cdot {{n\,+\,2}\over {\,n\,-\ 2\,}}} \\[0.2in]
& \le &    \left[\  {\bf V}_1^{\,{1\over 2} \cdot {{n\,+\,2}\over {\,n\,-\ 2\,}}}\ + \  {\bf V}_2^{\,{1\over 2} \cdot {{n\,+\,2}\over {\,n\,-\ 2\,}}}  \ +  \cdot \cdot \cdot  + \ {\bf V}_\flat^{\,{1\over 2} \cdot {{n\,+\,2}\over {\,n\,-\ 2\,}}}  \ \right] \, * \\[0.2in]
& \ &  \ \ \ \ \  \ \ \ \ \ \ *\, \left[\  {\bf V}_1^{\,{1\over 2} \cdot {{n\,+\,2}\over {\,n\,-\ 2\,}}}\ + \  {\bf V}_2^{\,{1\over 2} \cdot {{n\,+\,2}\over {\,n\,-\ 2\,}}} \ +  \cdot \cdot \cdot  + \ {\bf V}_\flat^{\,{1\over 2} \cdot {{n\,+\,2}\over {\,n\,-\ 2\,}}}  \ \right] \\[0.1in]
& \ & \left[ \ {1\over 2} \cdot {{n\,+\,2}\over {\,n\,-\ 2\,}} \ \le \ 1 \ , \ \ \ \ (\,a \ + \ b\,)^\tau \ \le \ a^\tau \ + \ b^\tau \ \ \ \ \mfor \ \ a\,, \ b \ \ \mbox{and} \ \ \tau \ > \ 0\  \right]\\[0.1in]
& = & \left[\  {\bf V}_1^{{\,n\,+\,2\,}\over {n\,-\,2}} \ + \  {\bf V}_2^{{\,n\,+\,2\,}\over {n\,-\,2}} \ + \cdot \cdot \cdot  + \ {\bf V}_\flat^{{\,n\,+\,2\,}\over {n\,-\,2}} \ \right] \ + \ \\[0.2in]
& \ & \ \ \ \ \  + \ {\bf V}_1^{\,{1\over 2} \cdot {{n\,+\,2}\over {\,n\,-\ 2\,}}} \cdot \left[\    {\bf V}_2^{\,{1\over 2} \cdot {{n\,+\,2}\over {\,n\,-\ 2\,}}}  \ +  \cdot \cdot \cdot  + \ {\bf V}_\flat^{\,{1\over 2} \cdot {{n\,+\,2}\over {\,n\,-\ 2\,}}}  \ \right]\ + \ \\[0.2in]
& \  &  \hspace*{0.6in} \ \ \ \  + \  {\bf V}_2^{\,{1\over 2} \cdot {{n\,+\,2}\over {\,n\,-\ 2\,}}} \cdot \left[\    {\bf V}_1^{\,{1\over 2} \cdot {{n\,+\,2}\over {\,n\,-\ 2\,}}} \ +  \cdot \cdot \cdot  + \ {\bf V}_\flat^{\,{1\over 2} \cdot {{n\,+\,2}\over {\,n\,-\ 2\,}}}  \ \right]\ + \ \\[0.1in]
& \  &  \hspace*{1.6in} \ : \ \\[0.1in]
& \  &  \hspace*{0.6in} \ \ \ \ \ \ \ \ + \  {\bf V}_\flat^{\,{1\over 2} \cdot {{n\,+\,2}\over {\,n\,-\ 2\,}}} \cdot \left[\    {\bf V}_1^{\,{1\over 2} \cdot {{n\,+\,2}\over {\,n\,-\ 2\,}}}  \ +  \cdot \cdot \cdot  + \ {\bf V}_{\flat\,-\,1}^{\,{1\over 2} \cdot {{n\,+\,2}\over {\,n\,-\ 2\,}}} \ \right]  \ . \\
\end{eqnarray*}
Hence (\,the following expressions are to be evaluated at $\,Y\,\in\,\R^n$\,)
\begin{eqnarray*}
(\,A.4.50\,) \ \ \ \ \ \left[\  {\bf W}_{\,\flat}  \ \right]^{{2\,n}\over {n\,-\,2}}  & = &    \left[\  {\bf V}_1\ + \  {\bf V}_2  \ +  \cdot \cdot \cdot  + \ {\bf V}_{\,\flat}  \ \right] \cdot {\bf W}_\flat^{{\,n\,+\,2\,}\over {n\,-\,2}}\\[0.2in]
& \le &  \left[\  {\bf V}_1^{{2n}\over {\,n\,-\ 2\,}} \ + \ {\bf V}_1 \cdot  {\bf V}_2^{{\,n\,+\,2\,}\over {n\,-\,2}} \ + \cdot \cdot \cdot  + \ {\bf V}_1 \cdot  {\bf V}_\flat^{{\,n\,+\,2\,}\over {n\,-\,2}} \ \right]\ + \  \\[0.2in]
& \ &  \ \ \ \   \ + \ {\bf V}_1 \cdot {\bf V}_1^{\,{1\over 2} \cdot {{n\,+\,2}\over {\,n\,-\ 2\,}}} \cdot \left[\    {\bf V}_2^{\,{1\over 2} \cdot {{n\,+\,2}\over {\,n\,-\ 2\,}}}  \ +  \cdot \cdot \cdot  + \ {\bf V}_\flat^{\,{1\over 2} \cdot {{n\,+\,2}\over {\,n\,-\ 2\,}}}  \ \right]\ + \ \\[0.1in]
& \ &  \ \ \ \   \ + \ {\bf V}_1 \cdot  {\bf V}_2^{\,{1\over 2} \cdot {{n\,+\,2}\over {\,n\,-\ 2\,}}} \cdot \left[\    {\bf V}_1^{\,{1\over 2} \cdot {{n\,+\,2}\over {\,n\,-\ 2\,}}} \ +  \cdot \cdot \cdot  + \ {\bf V}_\flat^{\,{1\over 2} \cdot {{n\,+\,2}\over {\,n\,-\ 2\,}}}  \ \right] \ + \ \\[0.1in]
& \ &  \hspace*{1.5in}\ :\\[0.1in]
& \ &  \ \ \ \    \ + \  {\bf V}_1 \cdot  {\bf V}_\flat^{\,{1\over 2} \cdot {{n\,+\,2}\over {\,n\,-\ 2\,}}} \cdot \left[\    {\bf V}_1^{\,{1\over 2} \cdot {{n\,+\,2}\over {\,n\,-\ 2\,}}}  \ +  \cdot \cdot \cdot  + \ {\bf V}_{\flat\,-\,1}^{\,{1\over 2} \cdot {{n\,+\,2}\over {\,n\,-\ 2\,}}} \ \right] \ + \\[0.2in]
& \ & \hspace*{1.5in} \cdot \cdot \cdot \ + \ \cdot \cdot \cdot
\\[0.1in]
& \ & \hspace*{1.5in} \cdot \cdot \cdot \ + \ \cdot \cdot \cdot
\\[0.1in]
& \  &  + \  \left[\ {\bf V}_\flat\cdot   {\bf V}_1^{{\,n\,+\,2\,}\over {n\,-\,2}} \ + \  {\bf V}_\flat\cdot   {\bf V}_2^{{\,n\,+\,2\,}\over {n\,-\,2}} \ + \cdot \cdot \cdot  + \   {\bf V}_\flat^{{2n}\over {\,n\,-\ 2\,}} \ \right] \ + \\[0.2in]
& \ &  \ \ \ \   \ + \ {\bf V}_\flat\cdot {\bf V}_1^{\,{1\over 2} \cdot {{n\,+\,2}\over {\,n\,-\ 2\,}}} \cdot \left[\    {\bf V}_2^{\,{1\over 2} \cdot {{n\,+\,2}\over {\,n\,-\ 2\,}}}  \ +  \cdot \cdot \cdot  + \ {\bf V}_\flat^{\,{1\over 2} \cdot {{n\,+\,2}\over {\,n\,-\ 2\,}}}  \ \right]\ + \\[0.1in]
& \ &  \ \ \ \   \ + \ {\bf V}_\flat\cdot {\bf V}_2^{\,{1\over 2} \cdot {{n\,+\,2}\over {\,n\,-\ 2\,}}} \cdot \left[\    {\bf V}_1^{\,{1\over 2} \cdot {{n\,+\,2}\over {\,n\,-\ 2\,}}} \ +  \cdot \cdot \cdot  + \ {\bf V}_\flat^{\,{1\over 2} \cdot {{n\,+\,2}\over {\,n\,-\ 2\,}}}  \ \right]\ + \\[0.1in]
& \ &  \hspace*{1.5in} : \ \\[0.1in]
& \ &  \ \ \ \    \ + \  {\bf V}_\flat\cdot  {\bf V}_\flat^{\,{1\over 2} \cdot {{n\,+\,2}\over {\,n\,-\ 2\,}}} \cdot \left[\    {\bf V}_1^{\,{1\over 2} \cdot {{n\,+\,2}\over {\,n\,-\ 2\,}}}  \ +  \cdot \cdot \cdot  + \ {\bf V}_{\flat\,-\,1}^{\,{1\over 2} \cdot {{n\,+\,2}\over {\,n\,-\ 2\,}}} \ \right] \  . \ \ \ \ \ \ \ \ \ \ \ \ \  \ \ \ \ \ \ \ \ \ \ \ \ \
\end{eqnarray*}
{\it Three types of terms in}\, (\,A.4.50\,)\,.\,\\[0.1in]
(\,I\,) \ \ \ Unmixed term like $\,{\bf V}_1^{{2n}\over {\,n\,-\ 2\,}}\,.$\\[0.1in]
(\,II\,) \  \,Simple mixed terms like $\, {\bf V}_1  \cdot {\bf V}_1^{\,{1\over 2} \cdot {{n\,+\,2}\over {\,n\,-\ 2\,}}} \cdot  {\bf V}_j^{\,{1\over 2} \cdot {{n\,+\,2}\over {\,n\,-\ 2\,}}}\ \ (\,j\,\not=\,1\,)\,$ and $\ {\bf V}_1  \cdot {\bf V}_j^{ {{n\,+\,2}\over {\,n\,-\ 2\,}}} \,.$\\[0.1in]
(\,III\,) \ Triple mixed terms like $\, {\bf V}_1  \cdot {\bf V}_1^{\,{1\over 2} \cdot {{n\,+\,2}\over {\,n\,-\ 2\,}}} \cdot  {\bf V}_j^{\,{1\over 2} \cdot {{n\,+\,2}\over {\,n\,-\ 2\,}}}\,, \ \ j \ \not= \ k\,.$\\[0.1in]

{\it Estimate of the unmixed term. }
\begin{eqnarray*}
\int_{\displaystyle{\R^n \,\Big\backslash \,\bigcup_{l\,=\,1}^\flat B_{\,\Xi_{\,l}}\,(\,R_{\,\nu}\,) }}\hspace*{-0.6in}
{\bf V}_1^{{2n}\over {\,n\,-\ 2\,}}\ \  & \le &  \int_{\R^n \ \setminus \  B_{\,\Xi_{\,\,1}}\,(\,R_{\,\nu}\,) }
{\bf V}_1^{{2n}\over {\,n\,-\ 2\,}}  \\[0.2in]
& \le  & C\,\cdot\,{1\over {R_{\,\nu}^{\,n } }} \ = \  O\left(  \ {\bar\lambda}_{\,\,\flat}^{n\,(\,1\ - \ \nu\,) } \ \emph{} \right) \ .
\end{eqnarray*}
Thus
\begin{eqnarray*}
(A.4.51) \ \ \ \ \ \ \ \  \ \ \ \ \ \ \ \ & \ &
\int_{\displaystyle{\R^n \,\Big\backslash \,\bigcup_{l\,=\,1}^\flat B_{\,\Xi_{\,l}}\,(\,R_{\,\nu}\,) }}\ \left[\  {\bf V}_1^{{2n}\over {\,n\,-\ 2\,}}\ + \  {\bf V}_2^{{2n}\over {\,n\,-\ 2\,}}  \ +  \cdot \cdot \cdot  + \ {\bf V}_\flat^{{2n}\over {\,n\,-\ 2\,}}  \ \right] \ \ \ \ \ \ \ \ \ \ \ \ \ \ \ \ \\[0.2in]
&  = & \flat \cdot   O\left(  \ {\bar\lambda}_{\,\,\flat}^{\ n\,(\,1\ - \ \nu\,) } \  \right) \ = \  O\left( \  {\bar\lambda}_{\,\,\flat}^{n\,(\,1\ - \ \nu\,) \ - \ \sigma }\  \right) \ .\\
\end{eqnarray*}


\newpage

{\it Estimate of   double mixed terms\,.}\\[0.1in]
(A.4.52)
\begin{eqnarray*}
& \ & {\bf V}_1  \cdot {\bf V}_1^{\,{1\over 2} \cdot {{n\,+\,2}\over {\,n\,-\ 2\,}}} \cdot  {\bf V}_j^{\,{1\over 2} \cdot {{n\,+\,2}\over {\,n\,-\ 2\,}}} \ \ \ \ \ \ \ \ \ \ \ \ \ \  \ \ \ \ \ \ \   (\,j\,\not=\,1\,) \\[0.2in]
& = & \left(\ {{\Lambda_{\,1}}\over {\Lambda_{\,1}^2 \ + \ \Vert\,Y\ - \ \Xi_{\,\,1}\,\Vert^{\,2} }}\ \right)^{\!\!{{n\,+\,2}\over 4} \ + \ {{\,n\,-\,2\,}\over 2} }  \cdot \left(\ {{\Lambda_j}\over {\Lambda_j^2 \ + \ \Vert\,Y\ - \ \Xi_j\,\Vert^{\,2} }}\ \right)^{\!\!{{n\,+\,2}\over 4}} \ \ \ \ \ \ \ \ \ \ \ \ \ \   (\,j\ \not= \ 1\,)\\[0.2in]
& \le & C_1 \cdot   \left(\ {{1}\over {1 \ + \ \Vert\,Y\ - \ \Xi_{\,\,1}\,\Vert }}\ \right)^{\!\!{{n\,+\,2}\over 2} \ + \ (\,n\,-\,2\,)  }   \cdot \left(\ {{1}\over {1 \ + \ \Vert\,Y\ - \ \Xi_j\,\Vert  }}\ \right)^{\!\!{{n\,+\,2}\over 2} }\\[0.2in]
& \le & C_2 \cdot  {1\over {\Vert\, \Xi_{\,\,1} \ - \ \Xi_j\,\Vert^{{{n\,+\,2}\over 2}\,-\,\epsilon} }} \cdot \bigg[\ \left(\ {{1}\over {1 \ + \ \Vert\,Y\ - \ \Xi_{\,\,1}\,\Vert }}\ \right)^{\!\!{{n\,+\,2}\over 2} \ + \ (\,n\,-\,2\,) \,+\,\varepsilon }\\[0.2in]& \ &
\hspace*{3in} \ + \ \left(\ {{1}\over {1 \ + \ \Vert\,Y\ - \ \Xi_j\,\Vert }}\ \right)^{\!\!2\,(\,n\,-\,2\,)\,+\,\epsilon} \ \ \bigg]\\[0.1in]
& \ & \ \ \ \ \ \ \ \ \ \ \ \ \ \ \ \ \   \left[\ \,{\mbox{note \ \ that}} \ \  n \ \ge \ 6 \ \ \Longrightarrow \ \ 2\,(\,n\,-\,2\,) \ \ge \ {{n\,+\,2}\over 2} \ + \ (\,n\,-\,2\,) \,+\,\varepsilon \ \right]\ .
\end{eqnarray*}
Also
\begin{eqnarray*}
(A.4.53) \ \ \ \ \ \ \ \ \ \ \ \ \ & \ &
\int_{\,\R^n \setminus \,  \,B_{\xi_{\,l} } ({_{\,R\nu}})}    \left(\ {{1}\over {1 \ + \ \Vert\,Y\ - \ \Xi_{\,\,1}\,\Vert }}\ \right)^{\!\!{{n\,+\,2}\over 2} \ + \ (\,n\,-\,2\,) \,+\,\epsilon }\ d\,V\\[0.2in]
& \le & C  \int_{ {{\rho_{\,\nu}}\over {\Lambda_{\,1}}} }^\infty\  {{R^{n\,-\,1} }\over {\  R^{{{n\,+\,2}\over 2} \ + \ (\,n\,-\,2\,) \,+\,\epsilon  } \ }}\\[0.2in]& \le & C_1 \cdot {1\over {\  R^{{{\,n\,-\,2\,}\over 2} \ + \ \epsilon  } \ }} \ \Bigg\vert^{\,\infty}_{\,{{1\over { \ {\bar\lambda}_{\,\,\flat}^{\,1 \, - \, \nu} \  }}}} \ \ \le \  C_2 \cdot   {\bar\lambda}_{\,\,\flat}^{ \left(\ {{ n\,-\,2}\over 2}\,+\,\epsilon \right) \ \cdot \ (\,1\ - \ \nu\,) } \ . \ \ \ \ \ \ \ \ \ \ \ \ \ \ \ \ \ \ \
\end{eqnarray*}

\newpage

It follows that
\begin{eqnarray*}
 \ \ \ \ \ \ \ \ \ \ \ & \ &    \int_{\displaystyle{\R^n \,\Big\backslash \,\bigcup_{l\,=\,1}^\flat B_{\,\Xi_{\,l}}\,(\,R_{\,\nu}\,) }}\   {\bf V}_1 \cdot  {\bf V}_1^{\,{1\over 2} \cdot {{n\,+\,2}\over {\,n\,-\ 2\,}}} \cdot \left[\    {\bf V}_2^{\,{1\over 2} \cdot {{n\,+\,2}\over {\,n\,-\ 2\,}}}  \ +  \cdot \cdot \cdot  + \ {\bf V}_\flat^{\,{1\over 2} \cdot {{n\,+\,2}\over {\,n\,-\ 2\,}}}  \ \right]\ \ \ \ \ \ \ \ \ \ \\[0.2in]
& \le & C \left(\ \sum_{j\ =\,2}^\flat \ {1\over {\ \Vert\, \Xi_{\,\,1} \ - \ \Xi_j\,\Vert^{{{n\,+\,2}\over 2}\,-\,\epsilon}\  }}  \ \right) \cdot {\bar\lambda}_{\,\,\flat}^{ \left(\ {{ n\,-\,2}\over 2}\,+\,\epsilon \right) \ \cdot \ (\,1\ - \ \nu\,) } \\[0.2in]
&  \le  & C_1  \cdot {\bar\lambda}_{\,\,\flat}^{\  \left(\ {{ n\,+\,2}\over 2}\,-\,\epsilon \right) \ \cdot \,\gamma \ + \  \left(\ {{ n\,-\,2}\over 2}\,+\,\epsilon \right) \ \cdot \ (\,1\ - \ \nu\,) }\ .\\
\end{eqnarray*}
Summing up, we obtain
\begin{eqnarray*}
(A.4.54)  & \ &    \int_{\displaystyle{\R^n \,\Big\backslash \,\bigcup_{l\,=\,1}^\flat B_{\,\Xi_{\,l}}\,(\,R_{\,\nu}\,) }}\  \left\{ \  {\bf V}_1 \cdot  {\bf V}_1^{\,{1\over 2} \cdot {{n\,+\,2}\over {\,n\,-\ 2\,}}} \cdot \left[\    {\bf V}_2^{\,{1\over 2} \cdot {{n\,+\,2}\over {\,n\,-\ 2\,}}}  \ +  \cdot \cdot \cdot  + \ {\bf V}_\flat^{\,{1\over 2} \cdot {{n\,+\,2}\over {\,n\,-\ 2\,}}}  \ \right] \ + \ \cdot \cdot \cdot  \right.\ \ \ \ \ \ \ \ \ \ \\[0.2in]
& \ &  \ \ \ \ \ \ \ \ \ \ \ \ \ \ \ \ \ \ \ \ \ \ \ \ \ \ \ \ \ + \  \left.  {\bf V}_{\,\flat} \cdot  {\bf V}_{\,\flat}^{\,{1\over 2} \cdot {{n\,+\,2}\over {\,n\,-\ 2\,}}} \cdot \left[\    {\bf V}_1^{\,{1\over 2} \cdot {{n\,+\,2}\over {\,n\,-\ 2\,}}}  \ +  \cdot \cdot \cdot  + \ {\bf V}_{\,\flat\,-\,1}^{\,{1\over 2} \cdot {{n\,+\,2}\over {\,n\,-\ 2\,}}}  \ \right] \ \right\}\\[0.2in]
&  \le  & C_1  \cdot {\bar\lambda}_{\,\,\flat}^{\  \left(\ {{ n\,+\,2}\over 2}\,-\,\epsilon \right) \ \cdot \,\gamma \ + \  \left(\ {{ n\,-\,2}\over 2}\,+\,\epsilon \right) \ \cdot \ (\,1\ - \ \nu\,) } \ \times \,\flat\\[0.2in]
&  \le  & C_2  \cdot {\bar\lambda}_{\,\,\flat}^{\  \left(\ {{ n\,+\,2}\over 2}\,-\,\epsilon \right) \ \cdot \,\gamma \ + \  \left(\ {{ n\,-\,2}\over 2}\,+\,\epsilon \right) \ \cdot \ (\,1\ - \ \nu\,) \ - \ \sigma} \ .\\
\end{eqnarray*}
Likewise,
\begin{eqnarray*}
(A.4.55) & \ & {\bf V}_1  \cdot {\bf V}_j^{ {{n\,+\,2}\over {\,n\,-\ 2\,}}} \\[0.2in]
& = & \left(\ {{\Lambda_{\,1}}\over {\Lambda_{\,1}^2 \ + \ \Vert\,Y\ - \ \Xi_{\,\,1}\,\Vert^{\,2} }}\ \right)^{\!\! {{\,n\,-\,2\,}\over 2} }  \cdot \left(\ {{\Lambda_j}\over {\Lambda_j^2 \ + \ \Vert\,Y\ - \ \Xi_j\,\Vert^{\,2} }}\ \right)^{\!\!{{n\,+\,2}\over 2}} \ \ \ \ \ \ \ \ \ \ \   (\,j\ \not= \ 1\,)\\[0.2in]
& \le & C_1 \cdot   \left(\ {{1}\over {1 \ + \ \Vert\,Y\ - \ \Xi_{\,\,1}\,\Vert }}\ \right)^{\!\! n\,-\,2\ }   \cdot \left(\ {{1}\over {1 \ + \ \Vert\,Y\ - \ \Xi_j\,\Vert  }}\ \right)^{\!\!n\,+\,2}\\[0.2in]
& \ &\hspace*{-0.75in}\le \ C_2 \cdot  {1\over {\Vert\, \Xi_{\,\,1} \ - \ \Xi_j\,\Vert^{{{\ n\,-\,2}}} }} \cdot \left[\ \left(\ {{1}\over {1 \ + \ \Vert\,Y\ - \ \Xi_{\,\,1}\,\Vert }}\ \right)^{\!\!n\,+\,2}  \ + \ \left(\ {{1}\over {1 \ + \ \Vert\,Y\ - \ \Xi_j\,\Vert }}\ \right)^{\!\!n\,+\,2} \ \ \right]
\end{eqnarray*}
\begin{eqnarray*}
\Longrightarrow  & \ &    \int_{\displaystyle{\R^n \,\Big\backslash \,\bigcup_{l\,=\,1}^\flat B_{\,\Xi_{\,l}}\,(\,R_{\,\nu}\,) }}\   {\bf V}_1 \cdot   \left[\    {\bf V}_2^{ {{n\,+\,2}\over {\,n\,-\ 2\,}}}  \ +  \cdot \cdot \cdot  + \ {\bf V}_\flat^{  {{n\,+\,2}\over {\,n\,-\ 2\,}}}  \ \right]\ \ \ \ \ \ \ \ \ \ \\[0.2in]
& \le & C_3 \left(\ \sum_{j\ =\,2}^\flat \ {1\over {\ \Vert\, \Xi_{\,\,1} \ - \ \Xi_j\,\Vert^{\ n\,-\,2}\  }}  \ \right) \cdot {\bar\lambda}_{\,\,\flat}^{\ 2\, \cdot \, (\,1\ - \ \nu\,) } \\[0.2in]
&  \le  & C_4  \cdot {\bar\lambda}_{\,\,\flat}^{\  \left(\ n\,-\,2\ \right) \ \cdot \,\gamma \ + \  2 \, \cdot \ (\,1\ - \ \nu\,) }\\[0.2in]
\Longrightarrow & \ &    \int_{\displaystyle{\R^n \,\Big\backslash \,\bigcup_{l\,=\,1}^\flat B_{\,\Xi_{\,l}}\,(\,R_{\,\nu}\,) }}\  \left\{ \  {\bf V}_1 \cdot  {\bf V}_1^{\,{1\over 2} \cdot {{n\,+\,2}\over {\,n\,-\ 2\,}}} \cdot \left[\    {\bf V}_2^{\,{1\over 2} \cdot {{n\,+\,2}\over {\,n\,-\ 2\,}}}  \ +  \cdot \cdot \cdot  + \ {\bf V}_\flat^{\,{1\over 2} \cdot {{n\,+\,2}\over {\,n\,-\ 2\,}}}  \ \right] \ + \ \cdot \cdot \cdot  \right.\ \ \ \ \ \ \ \ \ \ \\[0.2in]
& \ &  \ \ \ \ \ \ \ \ \ \ \ \ \ \ \ \ \ \ \ \ \ \ \ \ \ \ \ \ \ \ \ \ \ \ \ \ \ \ \ \ \ \ \ \ \ \ \ + \  \left.  {\bf V}_{\,\flat} \cdot   \left[\    {\bf V}_1^{  {{n\,+\,2}\over {\,n\,-\ 2\,}}}  \ +  \cdot \cdot \cdot  + \ {\bf V}_{\,\flat\,-\,1}^{  {{n\,+\,2}\over {\,n\,-\ 2\,}}}  \ \right] \ \right\}\\[0.2in]
&  \le  & C_5  \cdot {\bar\lambda}_{\,\,\flat}^{\  \left(\ n\,-\,2 \right) \ \cdot \,\gamma \ + \  2\, \cdot \ (\,1\ - \ \nu\,) } \ \times \,\flat
\  \le  \ C_6    \cdot {\bar\lambda}_{\,\,\flat}^{\  \left(\ n\,-\,2 \right) \ \cdot \,\gamma \ + \  2\, \cdot \ (\,1\ - \ \nu\,) \ - \ \sigma }\ \ \ .\\[0.2in]
\end{eqnarray*}

{\it Estimate of  triple mixed  terms.} \ \ Let us take a look at, for example\,,\, \\[0.2in]
(A.4.56)
\begin{eqnarray*}
& \ & {\bf V}_1  \cdot {\bf V}_2^{\,{1\over 2} \cdot {{n\,+\,2}\over {\,n\,-\ 2\,}}} \cdot  {\bf V}_3^{\,{1\over 2} \cdot {{n\,+\,2}\over {\,n\,-\ 2\,}}} \\[0.2in]
& \le &  C   \left(\ {{1}\over {1 \ + \ \Vert\,Y\ - \ \Xi_{\,\,1}\,\Vert }}\ \right)^{\!\! n\,-\,2 }  \cdot \left(\ {{1}\over {1 \ + \ \Vert\,Y\ - \ \Xi_{\,2}\,\Vert  }}\ \right)^{\!\!{{n\,+\,2}\over 2} }\cdot \left(\ {{1}\over {1 \ + \ \Vert\,Y\ - \ \Xi_3\,\Vert  }}\ \right)^{\!\!{{n\,+\,2}\over 2} }\\[0.2in]
& \le  &    \ C_2  \cdot {1\over {\Vert\, \Xi_{\,\,1} \ - \ \Xi_{\,2}\,\Vert^{{{n\,+\,2}\over 2}\,-\,\epsilon} }} \cdot \left[\ \left(\ {{1}\over {1 \ + \ \Vert\,Y\ - \ \Xi_{\,\,1}\,\Vert }}\ \right)^{\!\!(\,n\,-\,2\,)\,+\,\epsilon}  \ + \right.\\[0.2in]
& \ & \left. \ \ \ \ \ \ \ \ \ \ \ + \  \left(\ {{1}\over {1 \ + \ \Vert\,Y\ - \ \Xi_{\,2}\,\Vert }}\ \right)^{\!\!(\,n\,-\,2\,)\,+\,\epsilon} \ \right] \ \,\times
 \,\left(\ {{1}\over {1 \ + \ \Vert\,Y\ - \ \Xi_3\,\Vert  }}\ \right)^{\!\!{{n\,+\,2}\over 2}} \\[0.2in]
& \le  &   C_3  \cdot {1\over {\ \Vert\, \Xi_{\,\,1} \ - \ \Xi_{\,2}\,\Vert^{{{n\,+\,2}\over 2}\,-\,\epsilon}\ }}\cdot {1\over {\ \Vert\, \Xi_{\,\,1} \ - \ \Xi_{\,3}\,\Vert^{{{\,n\,-\,2\,}\over 2}}\  }} \cdot \left[\ \left(\ {{1}\over {1 \ + \ \Vert\,Y\ - \ \Xi_{\,\,1}\,\Vert }}\ \right)^{\!\!n\,+\,\epsilon} \ \ +  \right. \\[.2in]
& \ &\hspace*{3.3in} \left. \ + \ \left(\ {{1}\over {1 \ + \ \Vert\,Y\ - \ \Xi_{\,3}\,\Vert }}\ \right)^{\!\!n\,+\,\epsilon} \ \right] \ + \  \\[0.2in]& \  &  \\[0.2in]
& \  &  \ \ \ \ \ \  + \ C_4   \cdot {1\over {\Vert\, \Xi_{\,\,1} \ - \ \Xi_{\,2}\,\Vert^{{{n\,+\,2}\over 2}\,-\,\epsilon}}}\cdot {1\over {\Vert\, \Xi_{\,2} \ - \ \Xi_{\,3}\,\Vert^{{{\,n\,-\,2\,}\over 2}\,} }} \cdot \left[\ \left(\ {{1}\over {1 \ + \ \Vert\,Y\ - \ \Xi_{\,2}\,\Vert }}\ \right)^{\!\!n\,+\,\epsilon}   \ \ + \right. \\[.2in]
& \ &\hspace*{3.3in} \left. \ + \  \left(\ {{1}\over {1 \ + \ \Vert\,Y\ - \ \Xi_{\,3}\,\Vert }}\ \right)^{\!\!n\,+\,\epsilon} \ \right] \ .\\
\end{eqnarray*}
Upon integration,
summing up and making use of condition (\,1.22\,) in the main text\, , we have\\[0.1in]
(\,A.4.57\,)
\begin{eqnarray*}
 & \ &  \int_{\displaystyle{\R^n \,\Big\backslash \,\bigcup_{l\,=\,1}^\flat B_{\,\Xi_{\,l}}\,(\,R_{\,\nu}\,) }}\ \left\{ \
  {\bf V}_1 \cdot   {\bf V}_2^{\,{1\over 2} \cdot {{n\,+\,2}\over {\,n\,-\ 2\,}}} \cdot  \left[\   {\bf V}_{3 }^{\,{1\over 2} \cdot {{n\,+\,2}\over {\,n\,-\ 2\,}}} \ + \ \cdot \cdot \cdot \ + \  {\bf V}_{\flat}^{\,{1\over 2} \cdot {{n\,+\,2}\over {\,n\,-\ 2\,}}} \ \right] \ \ +  \right. \\[0.2in]
& \ &\hspace*{1in} \ \ \ \ \ \ \ \ \ \ \ \  + \ \cdot \cdot \cdot      \ + \  {\bf V}_1 \cdot  {\bf V}_\flat^{\,{1\over 2} \cdot {{n\,+\,2}\over {\,n\,-\ 2\,}}} \cdot \left[\    {\bf V}_2^{\,{1\over 2} \cdot {{n\,+\,2}\over {\,n\,-\ 2\,}}}  \ +  \cdot \cdot \cdot  + \ {\bf V}_{\flat\,-\,1}^{\,{1\over 2} \cdot {{n\,+\,2}\over {\,n\,-\ 2\,}}} \ \right] \ + \\[0.2in]
& \ & \hspace*{2.3in} \cdot \cdot \cdot \ + \ \cdot \cdot \cdot
\\[0.1in]
& \ & \hspace*{2.3in} \cdot \cdot \cdot \ + \ \cdot \cdot \cdot
\\[0.1in]
& \ & \hspace*{1.3in}   \ \ \ \   \ + \ {\bf V}_\flat\cdot {\bf V}_1^{\,{1\over 2} \cdot {{n\,+\,2}\over {\,n\,-\ 2\,}}} \cdot \left[\    {\bf V}_2^{\,{1\over 2} \cdot {{n\,+\,2}\over {\,n\,-\ 2\,}}}  \ +  \cdot \cdot \cdot  + \ {\bf V}_{\flat\,-\,1}^{\,{1\over 2} \cdot {{n\,+\,2}\over {\,n\,-\ 2\,}}}  \ \right]\ + \\[0.1in]
& \ & \hspace*{1.3in}  \ \ \ \   \ + \ {\bf V}_\flat\cdot {\bf V}_2^{\,{1\over 2} \cdot {{n\,+\,2}\over {\,n\,-\ 2\,}}} \cdot \left[\    {\bf V}_1^{\,{1\over 2} \cdot {{n\,+\,2}\over {\,n\,-\ 2\,}}} \ +  \cdot \cdot \cdot  + \ {\bf V}_{\flat\,-\,1}^{\,{1\over 2} \cdot {{n\,+\,2}\over {\,n\,-\ 2\,}}}  \ \right]\ + \\[0.1in]
& \ &  \hspace*{1.5in} : \ \\[0.1in]
& \ & \hspace*{1.3in} \ \ \ \   \left. \ + \  {\bf V}_\flat\cdot  {\bf V}_{\flat\,-\,1}^{\,{1\over 2} \cdot {{n\,+\,2}\over {\,n\,-\ 2\,}}} \cdot \left[\    {\bf V}_1^{\,{1\over 2} \cdot {{n\,+\,2}\over {\,n\,-\ 2\,}}}  \ +  \cdot \cdot \cdot  + \ {\bf V}_{\flat\,-\,1}^{\,{1\over 2} \cdot {{n\,+\,2}\over {\,n\,-\ 2\,}}} \ \right] \ \right\}\\[0.2in]
   &\  \le &  C  \cdot {\bar\lambda}_{\,\,\flat}^{\ n\,\cdot\, \gamma \ + \   \epsilon \cdot\, [\,1\,-\ (\,\gamma \ + \ \nu\,)\,] }\ .\\
 \end{eqnarray*}
Summing up the estimates in  (\,A.4.51\,)\,,\,  (\,A.4.54\,)\,,\,  (\,A.4.55\,)\, and \, (\,A.4.57\,)\,,\, we arrive at \\[0.1in]
(A.4.58)
\begin{eqnarray*}
\int_{\displaystyle{\R^n \,\Big\backslash \,\bigcup_{l\,=\,1}^\flat B_{\,\Xi_{\,l}}\,(\,R_{\,\nu}\,) }}\ {\bf W}_\flat^{{2n}\over {\,n\,-\ 2\,}}
       &  \le &  C_1\cdot  {\bar\lambda}_{\,\,\flat}^{\,n\,(\,1\ - \ \nu\,) \ - \ \sigma } \ \ + \\[-0.2in]
        & \ & \ \ \ \ + \ C_2  \cdot \,{\bar\lambda}_{\,\,\flat}^{\  \left(\ {{ n\,+\,2}\over 2}\,-\,\epsilon \right) \ \cdot \,\gamma \ + \  \left(\ {{ n\,-\,2}\over 2}\,+\,\epsilon \right) \ \cdot \ (\,1\ - \ \nu\,)\ - \ \sigma }\ \\[0.2in]
        & \ &  \hspace*{-1in}+ \ C_3 \cdot {\bar\lambda}_{\,\,\flat}^{\  \left(\ n\,-\,2 \right) \ \cdot \,\gamma \ + \  2\, \cdot \ (\,1\ - \ \nu\,) \ - \ \sigma } \ \
        + \     C_4 \cdot   {\bar\lambda}_{\,\,\flat}^{\ n\,\cdot\, \gamma \ + \   \epsilon \cdot\, [\,1\,-\ (\,\gamma \ + \ \nu\,)\,] } \ .
\end{eqnarray*}

Combining with the estimate ``\,inside\,"\, [\ that is\,,\, (A.4.46)\ ] we finally get\\[0.1in]
(A.4.59)
\begin{eqnarray*}
\int_{ \R^n }
{\bf W}_\flat^{{2n}\over {\,n\,-\ 2\,}}  & = &  \flat \cdot \ \int_{\R^n}
{\bf V}_o^{{2n}\over {\ n\,-\ 2\ }} \  + \ O\,\left( \ {\bar\lambda}_{\,\,\flat}^{\ \left( {{\,n\,-\,2\,}\over 2} \,-\,o_{\,+}\,(\,1\,) \right)\,\cdot\,\gamma \ - \ \sigma}\,\right)\ + \ O\,\left( \   {\bar\lambda}_{\,\,\flat}^{\ n\,(\,1\ - \ \nu\,) \ - \ \sigma } \,\right) \ . \\[0.1in]
& \ & \hspace*{-1.3in}\left[\,{\mbox{with \ \ the \ \ last \ \ three \ \ terms \ \ in\ \  (\,A.4.58\,) \ \ absorbed}}  \ \, \uparrow\  \right]
\\[0.2in]
 &  &  \ \ \ \ \ \ \ \
\end{eqnarray*}

{\bf \S\,A\,4\,.\,k\,.} \ \
{\it Things we know about}\, $\Phi\,.\,$ \\[0.1in]
{\bf (\,a)} \ \ By equation (\,A.4.1\,)\,,\, suppose that ${\,\cal H}\,$ is continuous\,,\, via the Sobolev embedding theorem, Schauder theory, and a boot\,-\,strap argument similar to the proof of Lemma 2.4 (\,pp. 203\,) in \cite{Book-Schoen-Yau}\,,\, we obtain $\,\Phi\,\in \,C^{\,2}\ .$\\[0.1in]
{\bf (\,b)} \ \ Note that the equation as well as $\,{\tilde{\bf P}}_{\flat_{\,\parallelsum}}\,$  can be pulled back to $\,S^n\,\setminus \,{\bf N}\,$ via the stereographic projection (the pull back of $\,{\tilde{\bf P}}_{\flat_{\,\parallelsum}}\,$ to $\,S^n\,$ is smooth)\,.\, Here $\,{\bf N}\,$ denotes the north pole of $\,S^n.\,$  As
$$
\Phi \,\in\, {\cal D}^{\,1\,,\,2} \ \ \Longrightarrow \ \ \int_{\R^n} |\,\Phi\,|^{{2n}\over {\,n\,-\ 2\,}} \ <  \ \infty\,,
$$
the pull\,-\,back of $\,\Phi\,$ has a removable singularity at the north pole (at least in $\,C^o\,$ sense\,,\, similar to a result obtained by Brezis\,-\,Kato \cite{Brezis-Kato}\,)\,.\, Transferring back to $\,\R^n\,,\,$ we have
$$
\Phi\, (\,Y\,) \ = \ O\left(\ {1\over { R^{\,n\,-\,2} }}\ \right) \ \ \ \ \ \mfor \ \ R \, = \, \Vert\,Y\,\Vert \ \gg \ 1\,. \leqno (A.4.60)
$$




\newpage




\newpage

{\large{\bf \S\,A\,5\,. \ \ Solving the  ``\,perpendicular\,"  equation\,--\, a proof of }}\\[0.1in]
\hspace*{0.75in}  {\large{\bf Proposition 2.7\,.\,}} \\[0.2in]
Consider solving the equation (\,for $\,\Phi\,$)
$$
\Delta\,(\,{\bf W}_{\,\flat} \ + \ \Phi\,)  \ + \ (\,c_n \,\cdot\, K)\,\cdot\,({\bf W}_{\,\flat} \ + \ \Phi)_+^{{\,n\,+\,2\,}\over {n\,-\,2}} \ = \ {\tilde{\bf P}}_{\flat_{\,\parallelsum}}\,, \ \ \ \ \ {\mbox{where}} \ \ \ \Phi \ \in \  {\cal D}^{\,1\,, \ 2}_{\,\flat_{\,\perp}}\,. \leqno (A.5.1)
$$
Recall that\\[0.1in]
(\,A.5.2\,)
$$
{\tilde{\bf P}}_{\flat_{\,\parallelsum}} \ = \  \sum_{l\ =\,1}^\flat   c_{\ l}\,\cdot\, {{\bf V}}_{\,l}^{4\over {n\,-\,2}}\,\cdot\,\left(\ \Lambda_{\,l}\,\cdot\, {{\partial \, {{\bf V}}_{\,l} }\over {\partial \Lambda_{\,l} }}  \right) \ \,+\, \ \sum_{l\,=\,1}^\flat  \left(\,\sum_{j\ =\,1}^n c_{\ l\,,\,j}\,\cdot\,  {{\bf V}}_{\,l}^{4\over {n\,-\,2}}\,\cdot\, \left[\  \Lambda_{\,l}\,\cdot\, {{\partial \,{{\bf V}}_{\,l} }\over {\partial \,\Xi_{{\,l}_{|_j}} }}\  \right]\ \right)\ .
$$
Making use of the relation
$$
\Delta\,{{\bf V}}_{\,l} \ + \ n\,(\,n\,-\,2\,) \,{{\bf V}}_{\,l}^{{\,n\,+\,2\,}\over {n\,-\,2}} \ = \ 0 \ \ \ \ \Longrightarrow \ \ \ \ \Delta\,{\bf W}_{\,\flat} \ = \ -\,n\,(\,n\,-\,2\,)\,\left(\ \sum_{l\,=\,1}^\flat \  {{\bf V}}_{\,l}^{{\,n\,+\,2\,}\over {n\,-\,2}} \ \right) \ ,
$$
 equation  (\,A.5.1\,) can be rewritten as
 $$
\Delta \,\Phi \ + \ \left\{ \ {{n\,+\,2}\over {\,n\,-\ 2\,}}\,\cdot\,( \,c_n\,\cdot\,K)\cdot {\bf W}_\flat^{4\over{n\,-\,2}} \right\}\,\cdot \Phi \ \ =\  \left[\ {\bf N}\,(\Phi) \ + \ {\cal l}_{\ \flat}  \ \right]  \ + \ {\bf P}_{\flat_{/\!/}}\ , \leqno (\,A.5.3\,)
$$
where\\[0.1in]
(\,A.5.4\,)
\begin{eqnarray*}
{\bf N}\,(\Phi) & = & \ - \ \left\{ \ (\,c_n\,\cdot\,K)\,\cdot\,\left[\ ({\bf W}_{\,\flat} \ + \ \Phi)_+^{{\,n\,+\,2\,}\over {n\,-\,2}} \ -\,{\bf W}_\flat^{{n\,+\,2}\over{n\,-\,2}} \ - \ {{n\,+\,2}\over {\,n\,-\ 2\,}} \,\cdot\,{\bf W}_\flat^{4\over{n\,-\,2}}\cdot \Phi   \ \right] \ \right\}\,, \\[0.2in]
{\cal l}_{\ \bf\flat}  & = & \ - \ \left\{\  (\,c_n\,\cdot\,K)\,\cdot\, {\bf W}_\flat^{{\,n\,+\,2\,}\over {n\,-\,2}} \ -\, n\,(\,n\,-\,2)\,\cdot\,\sum_{l\,=\,1}^\flat  {{\bf V}}_{\,l}^{{\,n\,+\,2\,}\over {n\,-\,2}}  \ \right\}.\\
\end{eqnarray*}

{\it Remark\,} A.5.5\,.\, \ \ As in {\bf S\,A\,4\,.k}\,,\,
Via elliptic regularity theory \cite{Gilbarg-Trudinger}, any weak solution $\,\Phi\,\in \,{\cal D}^{\,1\,, \ 2}_{\,\flat_{\,\perp}}\,$ of\, (\,A.5.1\,) automatically becomes $\,\Phi\,\in\,C^2\,(\,\R^n)\,.\,$
Note that equation (\,A.5.1\,) as well as $\,{\tilde{\bf P}}_{\flat_{\,\parallelsum}}\,$  can be pulled back to $\,S^n\,$ via the stereographic projection (\,the pull back of $\,{\tilde{\bf P}}_{\flat_{\,\parallelsum}}\,$ is smooth)\,.\, As
$$
\Phi \,\in\, {\cal D}^{\,1\,,\,2} \ \ \Longrightarrow \ \ \int_{\R^n} |\,\Phi\,|^{{2n}\over {\,n\,-\ 2\,}} \ <  \ \infty\,,
$$
the pull\,-\,back of $\,\Phi\,$ has a removable singularity at the north pole (at least in $\,C^o\,$ sense\,,\, similar to a result obtained by Brezis\,-\,Kato \cite{Brezis-Kato}\,)\,.\, Transferring back to $\,\R^n\,,\,$ we have
$$
\Phi\, (\,Y\,) \ = \ O\left( {1\over { R^{\,n\,-\,2} }}\ \right) \ \ \ \ \ \mfor \ \ R \, = \, \Vert\ Y\,\Vert \ \gg \ 1\,.
$$

\newpage

\hspace*{0.5in}Based on the scheme presented in \cite{Wei-Yan}\,,\, we discuss a proof of Proposition 2.5 (of the original text)\,,\, divided into the following Lemma {\bf A.5.6}\,,\, {\bf \S\,A\,5.a}\, and {\bf \S\,A\,5.b}\,.\, We follow closely the proofs of  Lemmas 2.4 and 2.5 in \cite{Wei-Yan}\,.

\vspace*{0.2in}

{\bf Lemma A.5.6.} \ \ {\it For $\,n \,\ge\,6\,,\,$ under the conditions in Proposition} {\bf A.4.20}\,,\,    {\it assume also that $\,\Phi\,\in\,{\bf W}_*^\flat\,$  \,}   [\,{\it{so \ that \ }} {\it $\Vert \ \Phi\,\Vert_{*_Y}\,$ makes sense\,}\,,\,  {\it{refer to}} (\,A.3.6\,)\ ]\,.\,
{\it Then we have the following.}\\[0.1in]
${\bf (i)}_{\,(\,5.6)}$ \  \,$\,{\bf N}\,(\Phi)\,$ and $\ {\cal l}_{\ \flat}\,$ {\it can be written in the form }
$$
 {\cal F}\,\cdot\,h \ ,
$$
{\it where}
 $\,h\,\in\,{\cal D}^{1\,\,,\ 2}\,,\ $  {\it and } ${\cal F} \,\in\,L^{n\over 2}\,(\,\R^n\,)$\ .
 \\[0.2in]
${\bf (ii)}_{\,(\,5.6)}$  \  \,$\displaystyle{\Vert\,{\bf N}\,(\Phi)\, \Vert_{**_Y} \  \le \  C\,\cdot\,\flat^{\,{{4}\over {n\,-\,2}} } \,\cdot\,  \Vert \,\Phi\,\Vert_{*_Y}^{{\,n\,+\,2\,}\over {n\,-\,2}} }$ \ .  \\[0.2in]
${\bf (iii)}_{\,(\,5.6)}$  \  $\displaystyle{\Vert\, {\cal l}_{\ \flat}\, \Vert_{**_Y} \  \le \  C\,\cdot\, {\bar\lambda}_{\,\,\flat}^{\,{\vartheta}_{\,(\,A.5.7\,)} }\ .}$ \\[0.1in]
{\it Here}
$$
{\vartheta}_{\,(\,A.5.7\,)}
 \ = \ \mbox{Min}  \, \left\{ \     {{n}\over 2} \,\cdot\,(\,1\, - \ \nu\,)\,, \  \ \ \  \ell    \cdot  \nu \  \right\}\ \,-\ o_{\,+}\,(\,1\,)\,.\leqno (\,A.5.7\,)
$$

\vspace*{0.2in}

{\bf Proof.} \ \ ${\bf (i)}_{\,(\,5.6)}$\,. \ \ Check  that $\,{\bf N}\,(\Phi)\,$ and $\,{\cal l}_{\ \flat}\,$ can be written in the form
\begin{eqnarray*}
 \ \ \ \ \ \ \ \ \ \ \ \  & \ &  \ \ \ \ \ \ \ \ \  {\cal F}\cdot h \,,\, \ \ \ \
 \ \ {\mbox{where}} \ \ \ \  h\,\in\,{\cal D}^{\,1\,,\,2} \ \ \\[0.2in]
  {\mbox{and}} \ \ \ \ \ \ \ \ \   & \ &   |\,{\cal F}\,(\,Y\,)\,| \ \le \ {C\over {R^{\,2\,+\,\varepsilon} }} \ \  \ \ \  {\mbox{for}} \ \ \ \  R \ = \ \Vert\,Y\,\Vert \ \gg \ 1 \ , \ \ \ \ {\mbox{and}} \ \  {\cal F} \,\in\,L^{n\over 2}\,(\,\R^n\,) \ .  \ \ \ \ \ \ \ \ \ \\
\end{eqnarray*}
Here $\,\epsilon\ > \ 0\,$ is a small number\,.\,  [\,Cf. the note after (\,A.1.12\,)\,.\,]\,
We begin with the observation that
\begin{eqnarray*}
[\ {\bf V}_l\,(\,Y\,)\ ] ^{{\,n\,+\,2\,}\over {n\,-\,2}} \,  & = & [\ V_{1\,,\  o}\,(\,Y\,)\ ] ^{{4}\over {\,n\,-\,2\,}} \,\cdot\,\left[\ {{ [\ {\bf V}_l\,(\,Y\,)\ ] ^{{\,n\,+\,2\,}\over {n\,-\,2}}  }\over { [\ V_{1\,,\  o}\,(\,Y\,)\ ] ^{{4}\over {\,n\,-\,2\,}}   }}\ \right]\\[0.2in]
 & = &  \Lambda_{\,1}^{{n\,+\,2}\over 2} \,\cdot\,\left[\  \left(\ {1\over {1\ + \ R^2}}\  \right)^{\!\!2}    \ \right]\,\cdot\,\left[\  {{ (\,1\,+\,\Vert\,Y\,\Vert^{\,2}\,)^{\,2}}\over {(\,\Lambda_{\,1}^2\,+\,\Vert\,Y\,-\,\Xi_{\,l}\,\Vert^{\,2}\,)^{{n\,+\,2}\over 2}   }}  \ \right]  \ .
 \end{eqnarray*}
Here $\,l \ = \ 1\,,\, \ \cdot \cdot \cdot\,, \ \flat\,.\,$  We may take
$$ {\cal F}\,(\,Y\,) \ = \   \left(\ {1\over {1\ + \ R^2}}\  \right)^{\!\!2}\,,\, \ \ \ \ {\mbox{and}} \ \ \ \ h\,(\,Y\,) \ = \  \left[\  {{ (\,1\,+\,\Vert\,Y\,\Vert^{\,2})^{\,2}}\over {(\,\Lambda_{\,1}^2\,+\,\Vert\,Y\,-\,\Xi_{\,l}\,\Vert^{\,2}\,)^{{n\,+\,2}\over 2}   }} \ \right] \  \in   \ {\cal D}^{\,1\,,\,2}\ ,
$$
where $\,R \ = \ \Vert\,Y\,\Vert\,.\,$\
Similarly\,,\, we deal with the term\,:
$$
{\bf W}^{{\,n\,+\,2\,}\over {n\,-\,2}}\ .
$$
As $\,{\bf W}_{\,\flat}\,$ is the sum of a finite number of  bubbles\,,\, with centers $\,\{\,\Xi_{\,l}\,\}_{\,l\,=\,1}^\flat\ $  spread out in $\,\R^n\,,\,$  we see that there is a number $\,R_{\,\flat} \ \gg \ 1\,$  [\,depends on the collection of points $\,\{\,\Xi_{\,l}\,\}_{\,l\,=\,1}^\flat\ $]\ ,\, such that
$$
{\bf W}_\flat\,(\,Y\,) \  \le \    {{1}\over {R^{\, n\,-\,2 } }} \ \ \ \ \mfor \ \  R \ = \ \Vert\,Y\,\Vert \ \ge \ R_{\,\flat}\,,
$$
Next,
$$
\Vert \,\Phi\,\Vert_{*_Y} \ \le \ C \ \ \Longrightarrow \ \ |\,\Phi\,(\,Y\,)\,| \  \le \  {{C_2}\over {\ R^{ {{\,n\,-\,2\,}\over 2}\  +\,\tau_{\,>1} }\  }} \ \ \ \ \mfor \ \ \Vert\,Y\,\Vert \ = \ R \ \ge \ R_{\,\flat}\,.$$
(\,We may let $\,R_{\,\flat}\,$ to be larger if necessary.\,)
It follows that\\[0.1in]
(\,A.5.8\,)
\begin{eqnarray*}
& \ & [\,{\bf W}_\flat\,(\,Y\,)  \,+\,\Phi\,(\,Y\,) \,]_+^{{4}\over {\,n\,-\,2\,}}\\[0.2in] & \le & [\,{\bf W}_\flat\,(\,Y\,)  \,+\,|\,\Phi\,(\,Y\,) \,|\,]^{\,{{4}\over {n\,-\,2}}}\\[0.2in]
  &\le & [\,{\bf W}_\flat\,(\,Y\,) \,]^{{4}\over {\,n\,-\,2\,}} \,+\,|\,\Phi\,(\,Y\,) \,|^{{4}\over {\,n\,-\,2\,}}  \hspace*{1in}  \left(\,n \ \ge \ 6 \ \ \Longrightarrow \ \  {{4}\over {n\,-\,2}}  \ \le \ 1 \right)\\[0.2in]
& \le &  {C \over {(\,1 \,+\,R)^{\,2 \ +\ {4\over {\,n\,-\,2\,}}\,\cdot\,\tau_{\,>1} } }} \ \ \ \ \ \ \ \ \ \ \ \ \ \mfor \ \  R \ \ge \ R_{\,\flat}\,.
\end{eqnarray*}
\begin{eqnarray*}
\Longrightarrow \ \ [\ {\bf W}_\flat\,(\,Y\,) \,+\,\Phi\,(\,Y\,) \ ]_+^{\,{{n\,+\,2}\over {\,n\,-\ 2\,}}} &  = & [\ {\bf W}_\flat\,(\,Y\,)  \,+\,\Phi\,(\,Y\,) \ ]_+^{{4}\over {\,n\,-\,2\,}}\,\cdot\,(\,{\bf W}_{\,\flat} \,+\,\Phi\,)_+\\[0.15in]
& = &  O\,\left( \ {1\over {(\,1 \,+\,R\,)^{\,2 \ +\ {4\over {n\,+\,2}} \,\cdot \ \tau_{\,>1} } }}\  \right) \,\cdot\, (\,{\bf W}_{\,\flat} \,+\,\Phi\,)_+
\end{eqnarray*}
for $\,R \ \ge \ R_{\,\flat} \ \gg \ 1\,,$\, Also\,,\, \emph{}
$$
(\,{\bf W}_{\,\flat} \,+\,\Phi\,) \ \in \ {\cal D}^{\,1\,,\,2} \ \ \Longrightarrow \ \ (\,{\bf W}_{\,\flat} \,+\,\Phi\,)_+ \ \in \ {\cal D}^{\,1\,,\,2}\ .
$$
This in this case we may take
$$ {\cal F}\,(\,Y\,) \ = \   \ {1\over {(\,1 \,+\,R\,)^{\,2 \ +\ {4\over {n\,+\,2}} \,\cdot \ \tau_{\,>1} } }}\  \ ,\, \ \ \ \ {\mbox{and}} \ \ \ \ h  \ = \ (\,{\bf W}_{\,\flat} \,+\,\Phi\,)_+   \  \in   \ {\cal D}^{\,1\,,\,2}\ ,
$$
where $\,R \ = \ \Vert\,Y\,\Vert\,.\,$

\newpage

Similarly we work on the term
$$
{\bf W}_\flat^{\,{{4}\over {n\,-\,2}}}\,\cdot\,\Phi\ .
$$
Finally, together with the Minkowski inequality\,:
$$
\Vert \, F\ + \ G\,\Vert_{L^2} \ \le \ \Vert \, F \,\Vert_{\,L^2} \ + \ \Vert \,  G\,\Vert_{\,L^2} \ \ \  \ \ \mfor \ \ F\ \,\&\, \ G \ \in \ L^2\,(\,\R^n)\,,
$$
we conclude the proof of \,${\bf (\,i)}_{\,(\,5.6\,)}$ \,.

\vspace*{0.2in}

${\bf (\,ii\,)}_{\,(\,5.6\,)}$\,.\ \   Similar to (\,A.4.37\,)\,,\, we have
$$
(\,1 \,+\,T\,)^{\,{{n \,+\,2}\over {n\,-\,2}}}  \  = \  1 \ + \ \left( \ {{n\,+\,2}\over {\,n\,-\ 2\,}}\ \right)
\cdot T \ + \ O \left(\ \,|\,T\,|^{{n\,+\,2 }\over {  n\,-\,2}}\ \right) \ , \leqno (\,A.5.9\,)
$$
which holds when $|\,T\,|$ is small  (\,via a form of Taylor expansion with remainder\,)\,,\, and also  when $\, T\,> \,0\,$ is large\,,\,  by incorporating a larger constant in front of $\,O \left(\, |\,T\,|^{{n\,+\,2 }\over {  n\,-\,2}}\ \right) \,$.\,  It follows that (\,for the positive constant $\,C_1\,$ to be large enough\,)\\[0.1in]
(\,A.5.10\,)
\begin{eqnarray*}
& \ & \bigg\vert \ (\,1 \,+\,T)^{\,{{n \,+\,2}\over {n\,-\,2}}} \  - \left(\ 1 \ + \ {{n\,+\,2}\over {\,n\,-\ 2\,}}
\cdot T \right) \ \bigg\vert \ \le \ C_1\,\cdot\, |\,T\,|^{\,{{n\,+\,2}\over {\,n\,-\ 2\,}} }\ \ \ \ \ \ \ \ \ {\mbox{for}} \ \ -\,1 \ < \ T \ < \ \infty\\[0.2in]
\Longrightarrow  & \ & (A \ + \ B)^{{\,n\,+\,2\,}\over {n\,-\,2}} \ = \  A^{{\,n\,+\,2\,}\over {n\,-\,2}} \ + \ {{n\,+\,2}\over {\,n\,-\ 2\,}}\,\cdot\,A^{4\over {n\,-\,2}}\cdot B \ + \ O\left(\, |\,B\,|^{{\,n\,+\,2\,}\over {n\,-\,2}} \ \right)  \ \  \ \ \  \ \  \  \\[0.15in]
& \ &\hspace*{1.7in} [\ n \ \ge \ 6\,,\, \ \ A \ >  \ 0 \,, \ \ \ \ -\,A \ < \ B \ \ \Longleftrightarrow \ \ A \ + \ B \ >  \ 0\ ]\   .
\end{eqnarray*}
To prepare for (\,A.5.13\,) below\,,\, we note that, via H\"older's inequality,
\begin{eqnarray*} (\,A.5.11\,)
 & \ &  \sum_{i \,=\,1}^\flat a_{\,i} \cdot\,b_{\,i}  \ \le \ \left(\ \ \sum_{i \,=\,1}^\flat\ a_{\,i}^{\,p} \ \right)^{ \!\!{1\over p} }\,\cdot\,\left(\ \ \sum_{i \,=\,1}^\flat\ b_{\,i}^{\,q}  \ \right)^{ \!\!{1\over q} } \\[0.1in]
 & \ & \ \ \ \ \ \ \ \ \ \ \ \  \ \ \ \ \ \  \ \ \ \ \ \  \ \ \ \ \ \   \left( \, a_{\,i} \,,\ \,  b_{\,i} \ >  \ 0 \,; \ \ p\,,\ \,q \ > \ 1\,,\, \ \ {1\over p} \ + \ {1\over q} \ = \ 1\ \right) \\[0.1in]
 \Longrightarrow  & \ &\!\!\!\!\!\!\!  \left(\ \ \sum_{i \,=\,1}^\flat\ a_{\,i} \,\cdot\,b_{\,i}  \ \right)^{\!\!p}  \ \le \  \left(\ \ \sum_{i \,=\,1}^\flat\ b_{\,i}^{\,q}  \ \right)^{ \!\!{p\over {\,q}} }\,\cdot\, \left(\ \ \sum_{i \,=\,1}^\flat\ a_{\,i}^{\,p} \ \right)
 \\[0.15in]
 & \ & \ \ \ \ \ \ \ \ \ \ \ \  \ \ \ \ \ \  \ \ \ \ \ \ \ \ (\,\downarrow \ \ {\mbox{by \ \ taking \ \ }} b_{\,i}  \ = \ 1 \ \ {\mbox{for}} \ \ i \ = \ 1\,,\ 2\,, \ \cdot \cdot \cdot\,, \ \flat\,)
\\[0.1in]
 \Longrightarrow  & \ &\!\!\!\!\!\!\!   \left(\  \sum_{i \,=\,1}^\flat\ a_{\,i} \cdot\,1 \ \right)^{\!\!p}  \ \le \   \flat^{\,p\,-\,1}\,\cdot\, \left(\ \ \sum_{i \,=\,1}^\flat\ a_{\,i}^{\,p} \ \right) \ \ \  \left( \ {1\over p} \ + \ {1\over q} \ = \ 1 \ \,   \Longrightarrow   \ \, {p\over q} \ = \ p \ - \ 1 \ \right)\,.
\end{eqnarray*}
Consider the case when
$$
\Phi\,(\,Y\,) \ > \ -\, {\bf W}_\flat\,(\,Y\,) \ \ \  \ \Longleftrightarrow \ \ (\,{\bf W}_{\,\flat} \,+\,\Phi\,)\,(\,Y\,) \ > \ 0 \,.
$$
As such, set
$$
U_{\,>} \ := \ \left\{\  Y \,\in\,\R^n \ \, \bigg\vert \ \ (\,{\bf W}_{\,\flat} \,+\,\Phi\,)\,(\,Y\,) \ > \ 0 \ \right\}\ .\leqno (\,A.5.12\,)
$$
Applying (\,A.5.10\,)\,,\, we have \\[0.1in]
(\,A.5.13\,)
\begin{eqnarray*}
& \ & \bigg\vert \ (\,{\bf W}_{\,\flat} \,+\,\Phi\,)_+^{{\,n\,+\,2\,}\over {n\,-\,2}} \ - \  {\bf W}_\flat^{{\,n\,+\,2\,}\over {n\,-\,2}} \ - \ {{n\,+\,2}\over {\,n\,-\ 2\,}}\,\cdot\, {\bf W}_\flat^{{4}\over {\,n\,-\,2\,}} \,\cdot\,\Phi\ \bigg\vert_{\,Y} \ \ \ \ \ \ \ \  \ (\  \mbox{for} \ \ Y\,\in\,U_{\,>} \ )\\[0.2in]
& \le & C\cdot\!|\,\Phi\,(\,Y\,)\,|^{{\,n\,+\,2\,}\over {n\,-\,2}}
\\[0.2in]
& \le & C\!\cdot\!\left[\ |\,\Phi\,(\,Y\,)\,|\,\cdot\,\ \left(\, \,\sum_{l\,=\,1}^\flat {1\over {   (  1\ + \ \Vert\,Y \,-\,\Xi_{\,l}\,\Vert\,)^{ {{n\,-\,2}\over 2 } \,+\ \tau_{\,>1} } }}\  \right)^{\!\!-\,1} \ \right]^{\,{{n\,+\,2}\over {\,n\,-\ 2\,}}}\ * \\[0.2in]
& \ & \left.\hspace*{2.75in} * \ \left(\  \sum_{l\,=\,1}^\flat {1\over {   (  1\ + \ \Vert\,Y \,-\,\Xi_{\,l}\,\Vert\,)^{ {{n\,-\,2}\over 2 } \,+\ \tau_{\,>1} } }}\  \right)^{{\,n\,+\,2\,}\over {n\,-\,2}}\ \right]\\[0.15in]
& \le & C    \!\cdot\Vert\,\Phi\,\Vert_{*_Y}^{{\,n\,+\,2\,}\over {n\,-\,2}}  \left[\  \flat^{\,{{4}\over {n\,-\,2}} }\,\cdot\,  \sum_{l\,=\,1}^\flat \left(\ {1\over  {   (  1\ + \ \Vert\,Y \,-\,\Xi_{\,l}\,\Vert\,)^{ {{\,n\,-\,2\,}\over 2}\,\cdot\,{{n\,+\,2}\over {\,n\,-\ 2\,}}   \,\ +\ {{n\,+\,2}\over {\,n\,-\ 2\,}} \,\cdot \, {{n\,+\,2\,}\over {\,n\,-\,2\,}}  } }}\  \right)\ \right]\\[0.1in]
& \ & \hspace*{4.5in} \ \ \ [\ {\mbox{via\ \ (\,A.5.11\,)}}\ ]\\[0.1in]
& \le & C\,\cdot\,\flat^{\,{{4}\over {n\,-\,2}} }\,\cdot\,\Vert\,\Phi\,\Vert_{*_Y}^{{\,n\,+\,2\,}\over {n\,-\,2}}\,*\\[0.1in]
& \ & \ \ \ \ \ \ \ \ \ \ \ \ \ * \ \left[\  \sum_{l\,=\,1}^\flat \left(\, \, {1\over {   (  1\ + \ \Vert\,Y \,-\,\Xi_{\,l}\,\Vert\,)^{ {{n\,+\,2}\over 2}     \,\ +\ \,\tau_{\,>1}  }} }\,\cdot\,{1\over {   (  1\ + \ \Vert\,Y \,-\,\Xi_{\,l}\,\Vert\,)^{  \,{{4}\over {n\,-\,2}} \,\cdot\, \tau_{\,>1}  }}  }\right)\ \right] \\[0.2in]
& \ & \ \ \ \ \ \ \ \ \ \ \ \ \ \ \ \ \ \ \ \ \ \ \ \ \  \ \ \ \ \ \ \ \ \ \ \ \ \ \ \ \  \ \ \ \ \ \ \left[\  \downarrow \ \ {\mbox{as}} \ \ \  \left(\, {1\over {   (  1\ + \ \Vert\,Y \,-\,\Xi_{\,l}\,\Vert\,)^{  \,{{4}\over {n\,-\,2}} \,\cdot\, \tau_{\,>1}  }}  }\right) \ \le \ 1 \ \right] \\[0.1in]
& \le & C\,\cdot\,\flat^{\,{{4}\over {n\,-\,2}} }  \,\cdot\,\Vert\,\Phi\,\Vert_{*_Y}^{{\,n\,+\,2\,}\over {n\,-\,2}}\,\cdot\, \left[\ \sum_{l\,=\,1}^\flat \ {1\over {   (  1\ + \ \Vert\,Y \,-\,\Xi_{\,l}\,\Vert\,)^{ {{\,n\,-\,2\,}\over 2}   \,+\,2  \,\ +\ \,\tau_{\,>1}  }} }\ \right]  \ \ \ \ \  \ \mfor \ \ Y\,\in\,U_{\,>}\  .
\end{eqnarray*}
Hence we obtain \,${\bf (ii)}_{(\,5.6\,)}$  when (\,A.5.12\,) holds\,.\bk
In case
$$
\Phi\,(\,Y\,) \ \le  \ -\, {\bf W}_\flat\,(\,Y\,) \ \ \  \ \Longrightarrow \ \  |\,\Phi\,(\,Y\,)\,| \ \ge \ {\bf W}_\flat\,(\,Y\,) \,.
 $$
As in (\,A.5.12\,)\,,\, set
 $$
L_{\,\le} \ := \ \left\{\  Y \,\in\,\R^n \ \, \bigg\vert \ \ (\,{\bf W}_{\,\flat} \,+\,\Phi\,)\,(\,Y\,) \ \le \ 0 \ \right\}\ .\leqno (\,A.5.14\,)
$$
It follows that
\begin{eqnarray*}
  & \ & \bigg\vert \ (\,{\bf W}_{\,\flat} \,+\,\Phi\,)_+^{{\,n\,+\,2\,}\over {n\,-\,2}} \ - \  {\bf W}_\flat^{{\,n\,+\,2\,}\over {n\,-\,2}} \ - \ {{n\,+\,2}\over {\,n\,-\ 2\,}}\,\cdot\, {\bf W}_\flat^{{4}\over {\,n\,-\,2\,}} \,\cdot\,\Phi\ \bigg\vert_{\,Y} \ \ \ \ \ \ \  (\ \mbox{for} \ \ Y\,\in\,L_{\,\le}\ ) \\[0.15in]
& \le &  0 \ + \ [\,{\bf W}_\flat\,(\,Y\,)\,]^{{\,n\,+\,2\,}\over {n\,-\,2}} \ + \ {{n\,+\,2}\over {\,n\,-\ 2\,}}\,\cdot\, [\,{\bf W}_\flat\,(\,Y\,)\,]^{{4}\over {\,n\,-\,2\,}} \,\cdot\, |\,\Phi\,(\,Y\,)\,|  \ \ \ \ ( \ [ \ {\mbox{-ve}} \ ]_+ \ = \ 0 \ )  \\[0.15in]
& \le & C\,\cdot\,  |\,\Phi\,(\,Y\,)\,|^{{\,n\,+\,2\,}\over {n\,-\,2}} \ \ \ \ \ \ \ \ \ \ \ \ \ \ \ \ \ \ \ \ \ \mbox{for} \ \ Y\,\in\,L_{\,\le}\  \ \ \ \ \ \ \ \ \ \ \ \ \ \ \ \ \ \ \ \  [ \ {\mbox{via}} \ \ (\,A.5.14\,)\ ]\ .
\end{eqnarray*}
Arguing as in (\,A.5.13\,), we come to the conclusions in  \,${\bf (\,ii)}_{\,(\,5.6\,)}$\,. \\[0.3in]
%
%
%
%
%
%
%
\noindent${\bf (\,iii\,)_{\,(\,5.6\,)}}\,.$\ \ Write\\[0.1in]
(\,A.5.15\,)
\begin{eqnarray*}
& \ &  \left[ \ (\,c_n\,\cdot\,K)\,\cdot\, {\bf W}_\flat^{{\,n\,+\,2\,}\over {n\,-\,2}} \ -\, n
\,(\,n\,-\,2)\,\cdot\, \sum_{l\,=\,1}^\flat  {{\bf V}}_{\,l}^{{\,n\,+\,2\,}\over {n\,-\,2}}\ \right]_{\,Y}\\[0.2in]
 & \  & \!\!\!\!\!\!\!\!\! =  \left[(\,c_n\,\cdot\,K)\,\cdot\, {\bf W}_\flat^{{\,n\,+\,2\,}\over {n\,-\,2}} - \,(\,c_n\,\cdot\,K)  \sum_{l\,=\,1}^\flat  {{\bf V}}_{\,l}^{{\,n\,+\,2\,}\over {n\,-\,2}} \ \right]_{\,Y}  +    \left\{\,   [\,(\,c_n\,\cdot\,K) \, - \, n\,(\,n\,-\,2)\ ] \sum_{l\,=\,1}^\flat  {{\bf V}}_{\,l}^{{\,n\,+\,2\,}\over {n\,-\,2}} \ \right\}_{\,Y} \ .\\[0.1in]
 & \ & \ \leftarrow \ \ \ \hspace{0.8in}{\bf (\,a\,)}_{\,(\,5.15\,)}  \hspace{0.8in} \rightarrow \ \ \ \  \  \ \ \ \leftarrow \ \ \hspace{0.7in} {\bf (\,b\,)}_{\,(\,5.15\,)} \hspace{0.7in} \rightarrow
\end{eqnarray*}
In can be discerned in $\,{\bf (\,b\,)}_{\,(\,5.15\,)}\,$ that the flatness of $\,K\,$ is involved.

\vspace*{0.2in}

${\bf (\,a\,)}_{\,(\,5.15\,)}$\ \  Estimate of
\begin{eqnarray*}
& \ &
 \left[ \  {\bf W}_\flat^{{\,n\,+\,2\,}\over {n\,-\,2}} \ - \  \ \sum_{l \,=\,1}^\flat  {{\bf V}}_{\,l}^{{\,n\,+\,2\,}\over {n\,-\,2}} \ \right]_{\,Y} \\[0.1in]   & = &  \left[\ \sum_{l \,=\,1}^\flat \left(\ {{\Lambda_{\,1}}\over {\Lambda^2_l \ + \ \Vert\,y \ - \ \Xi_{\,l}\,\Vert^{\,2}}}\ \right)^{\!\!{{\,n\,-\,2\,}\over 2}} \ \right]^{ {{n\,+\,2}\over {\,n\,-\ 2\,}}  } \ - \ \  \sum_{l \,=\,1}^\flat \left(\ {{\Lambda_{\,1}}\over {\Lambda^2_l \ + \ \Vert\,y \ - \ \Xi_{\,l}\,\Vert^{\,2}}}\ \right)^{\!\!{{n\,+\,2}\over 2}} .
\end{eqnarray*}
Here we deal with pointwise estimate (\,cf. {\bf \S\,A\,4\,.i} and {\bf \S\,A\,4\,.i} for integral estimates\,)\,,\, as in  {\bf \S\,A\,3\,.p}\,,\, we
consider the (finite) collection of points (\,the centers of the bubbles\,)
$$
\Xi_{\,\,1}\,, \ \ \,\Xi_{2}\,, \    \,\cdot\,\cdot\,\cdot \ , \ \ \,\Xi_{\,j}\,, \ \,\cdot\,\cdot \cdot\,,\ \  \Xi_{\ \flat}\ .
$$
As in (\,A.3.61\,)\,,\, let \\[0.1in]
(\,A.5.16\,)
$$
\Upsilon_j \,:=\ \left\{ \,y \ \in \ \R^n \ \ | \ \ \,\Vert\, y \ - \ \Xi_{\,l}\,\Vert \ \ge \ \Vert\, y \ - \ \Xi_j\,\Vert \ \ \ \ {\mbox{for}} \ \ 1 \ \le \  l \ \not= \ j \ \le  \ \flat\,\right\}\,.
$$
Recall that
$$ \R^n \ = \ {\displaystyle{\bigcup_{\,j\ =\,1}^{\,k} }} \,\Upsilon_j\ . $$

\vspace*{0.15in}

{\it Focusing on the leading term.} \ \ Cf. {\bf \S\,A\,4\,.i}\,.\  The following expressions are to be evaluated at  $\,Y \ \in \ \Upsilon_1\,$.\, \\[0.1in]
(\,A.5.17\,)
\begin{eqnarray*}
0 & <  & \left(\ \, {{\bf V}}_1 \,+\, {{\bf V}}_2\,+\,\cdot\,\cdot\,\cdot\,+\, {{\bf V}}_\flat\right)^{\,{{n\,+\,2}\over{n\,-\,2}} }\ -\ \left(\  {{\bf V}}_1^{{n\,+\,2}\over{n\,-\,2}} \,+\, {{\bf V}}_2^{{n\,+\,2}\over{n\,-\,2}}\,+\,\cdot\,\cdot\,\cdot\,+\, {{\bf V}}_\flat^{{n\,+\,2}\over{n\,-\,2}}\ \right)\\[0.15in]
& \le &  {{\bf V}}_1^{{n\,+\,2}\over{n\,-\,2}} \ + \  {{n\,+\,2}\over{n\,-\,2}}\cdot {{\bf V}}_1^{4\over{n\,-\,2}}\,\cdot\,\left(\  {{\bf V}}_2  \,+\, {{\bf V}}_3 \,+\,\cdot\,\cdot\,\cdot\,+\, {{\bf V}}_{\,\flat} \right) \ +  \\[0.2in]
& \ & \  \ \  \ \ \ \ \ \ \ \ \ \ \ \    + \ C_1\,\cdot\,\left(\ \, {{\bf V}}_2  \,+\, {{\bf V}}_3 \,+\,\cdot\,\cdot\,\cdot\,+\, {{\bf V}}_{\,\flat} \right)^{{n\,+\,2}\over{n\,-\,2}} \\[0.15in]
& \  & \ \  \ \ \ \ \ \ \ \ \ \ \ \ \ \ \ \ \ \ \ \ \ \  - \ \left(\ {{\bf V}}_1^{{n\,+\,2}\over{n\,-\,2}} \,+\, {{\bf V}}_2^{{n\,+\,2}\over{n\,-\,2}}\,+\,\cdot\,\cdot\,\cdot\,+\, {{\bf V}}_\flat^{{n\,+\,2}\over{n\,-\,2}}\ \right) \ \ \ \ \ \ \    [\,{\mbox{via \ \ (\,A.5.10\,)\,}}]\\[0.15in]
& \le &  {{n\,+\,2}\over{n\,-\,2}}\cdot  {{\bf V}}_1^{4\over{n\,-\,2}}\,\cdot\,\left(\  {{\bf V}}_2  \,+\, {{\bf V}}_3 \,+\,\cdot\,\cdot\,\cdot\,+\, {{\bf V}}_{\,\flat} \right) \ + \ C_2\,\cdot\,\left(\  {{\bf V}}_2  \,+\, {{\bf V}}_3 \,+\,\cdot\,\cdot\,\cdot\,+\, {{\bf V}}_{\,\flat} \right)^{{n\,+\,2}\over{n\,-\,2}}\\[0.2in]
& \le & C_3 \!\cdot \left(\ {1\over {1 \ + \ \Vert\, Y \ - \ \Xi_{\,\,1}\,\Vert}}\ \right)^4\,\cdot\, \sum_{j\,=\,2}^\flat \left(\ {1\over {1 \ + \ \Vert\, Y \ - \ \Xi_j\,\Vert}}\ \right)^{n\,-\,2}\ + \  \ \ \    \cdot \cdot  \cdot \cdot \cdot  {\bf (\,a\,)}_{\,(\,A.5.17\,)} \\[0.1in]
 & \ & \ \ \ \ \ \ \ \ \ \ \ \ \ \ \ \ \ \ \ \ \ \ \ \ \   \   +\ C_4 \cdot \left[ \  \sum_{j\,=\,2}^\flat \left(\ {1\over {1 \ + \ \Vert\, Y \ - \ \Xi_j\,\Vert}}\ \right)^{n\,-\,2} \,\right]^{{\,n\,+\,2\,}\over {n\,-\,2}}. \ \ \cdot \cdot \cdot \cdot\  {\bf (\,b\,)}_{\,(\,A.5.17\,)}
\end{eqnarray*}


\newpage

{\it Separation.} \ \ For the term \,${\bf (\,a\,)}_{\,(\,A.5.17\,)}$\,,\, applying the Separation Lemma [\,see (\,A.3.12\,) and (\,A.3.13\,)\ ]\,,\, we proceed with\\[0.1in]
(\,A.5.18\,)
\begin{eqnarray*}
& \  &  \left(\ {1\over {1 \ + \ \Vert\, Y \ - \ \Xi_{\,\,1}\,\Vert}}\ \right)^4\,\cdot\,\sum_{j\,=\,2}^\flat \left(\ {1\over {1 \ + \ \Vert\, Y \ - \ \Xi_j\,\Vert}}\ \right)^{n\,-\,2}\\[0.15in]
& \le & C_1\,\cdot\, \sum_{j\,=\,2}^\flat \left[ \ \left(\ {1\over {(\,1 \ + \ \Vert\, Y \ - \ \Xi_{\,\,1}\,\Vert\,)^{{{n\,+\,2}\over 2} \, +\ \tau_{\,>1}}}} \ + \ {1\over {(\,1 \ + \ \Vert\, Y \ - \ \Xi_j\,\Vert\,)^{{{n\,+\,2}\over 2} \,+\ \tau_{\,>1}}}}\  \right)\,* \right.\\[0.2in]
& \ & \hspace*{3.7in} \ \ \ \ \left.*\ {1\over { \ \Vert\,\Xi_j\,-\,\Xi_{\,\,1}\,\Vert^{\,{{n\,+\,2}\over 2} \,-\ \tau_{\,>1}}\  }} \ \right] \\[0.15in]
& \le & C_2\,\cdot\,{1\over {(\,1 \ + \ \Vert\, Y \ - \ \Xi_{\,\,1}\,\Vert\,)^{\,{{n\,+\,2}\over 2} \,+\ \tau_{\,>1}}}}\,\cdot\, \left(\ \sum_{j\,=\,2}^\flat\  {1\over {\  \Vert\,\Xi_j\,-\,\Xi_{\,\,1}\,\Vert^{ \,{{ n\,+\,2 }\over { 2}}\,-\ \tau_{\,>1} } \ }}\ \right) \\[0.2in]
& \ & \hspace*{3.1in} \ \ \ \ \ [\ {\mbox{via}} \ \ (\,A.3.61\,) \ \  {\mbox{for}} \ \ Y \ \in \ \Upsilon_1\ ] \\[0.15in]
& \ & \hspace*{-0.37in}\left[ \ \le \  C_3 \cdot\,{1\over {(\,1 \ + \ \Vert\, Y \ - \ \Xi_{\,\,1}\,\Vert\,)^{\,{{n\,+\,2}\over 2} \,+\ \tau_{\,>1}}}}\,\cdot\,\, {\bar\lambda}_{\,\,\flat}^{\left(\, {{n}\over 2} \,-\ o_{\,+}\,(\,1\,)\, \right)\,\cdot\,\gamma} \ \ \ \ \  \ \ \ \ \ {\mbox{for}} \ \ Y \ \in \ \Upsilon_1\ \right] \ .\\[0.2in]
\end{eqnarray*}
{\it Note}\, A.5.19\,. \ \ We have
$$
 {{n\,+\,2}\over 2} \,-\ \tau_{\,>1}\ =\  {{n\,+\,2}\over 2} \,-\ [\  1 \ + \ o_{\,+}\,(\,1\,)\ ] \ =\  {n\over 2} \,-\,o_{\,+}\,(\,1\,)\ ,
$$
as $\,\tau_{\,>1}\,$ is a number slightly bigger than one.\\[0.2in]
\hspace*{0.5in}Looking at the second term \,${\bf (\,b\,)_{\,(\,A.5.17\,)}}$\,,\, for $\,Y \ \in \ \Upsilon_{\,1}\,$,\, we have
\begin{eqnarray*}
& \ & {1\over {(\,1 \ + \ \Vert\, Y \ - \ \Xi_{\,l}\,\Vert\,)^{\,n\,-\,2}}} \\[0.15in]
  & \le & {1\over {(\,1 \ + \ \Vert\, Y \ - \ \Xi_{\,\,1}\,\Vert\,)^{ {{\,n\,-\,2\,}\over 2}  }}}\,\cdot\,{1\over {(\,1 \ + \ \Vert\, Y \ - \ \Xi_{\,l}\,\Vert\,)^{ {{\,n\,-\,2\,}\over 2}  }}}  \ \ \ \ \ [\ {\mbox{via}} \ \ (\,A.3.61\,) \ \  {\mbox{for}} \ \ Y \ \in \ \Upsilon_1\ ] \\[0.15in]
& \le & C_1\,\cdot\,{1\over {\Vert\,\Xi_{\,\,1}\ - \ \Xi_{\,l}\,\Vert^{  {{\,n\,-\,2\,}\over 2} \ -\ {{n\,-\,2}\over {n\,+\,2}} \,\cdot\  \tau_{\,>1}  } }}\,\cdot\, {1\over {(\,1 \ + \ \Vert\, Y \ - \ \Xi_{\,\,1}\,\Vert\,)^{ {{\,n\,-\,2\,}\over 2} \,+\, {{n\,-\,2}\over {n\,+\,2}} \,\cdot\  \tau_{\,>1} }}}\\[0.2in]
  &\ & \hspace*{0.3in}    \left[\,{\mbox{argue \ \ as \ \ in \ \ (\,A.5.18\,)}}\,, \ \  {1\over {(\,1 \ + \ \Vert\, Y \ - \ \Xi_{\,\,1}\,\Vert\,)}} \ \ge \  {1\over {(\,1 \ + \ \Vert\, Y \ - \ \Xi_{\,l}\,\Vert\,)}}  \ \right]\\[0.15in]
  \Longrightarrow & \ &  \sum_{l\,=\,2}^\flat \ {1\over {(\,1 \ + \ \Vert\, Y \ - \ \Xi_{\,l}\,\Vert\,)^{\,n\,-\,2}}}
   \  \le \   C_2\,\cdot\,\left(\  \sum_{l\,=\,2}^\flat \ {1\over {\Vert\,\Xi_{\,\,1}\ - \ \Xi_{\,l}\,\Vert^{  {{\,n\,-\,2\,}\over 2} \,-\, {{n\,-\,2}\over {n\,+\,2}} \,\cdot\  \tau_{\,>1}  } }}  \right)\,*\\[0.1in]
   & \ & \hspace*{3in}    *\,{1\over {\ (\ 1 \ + \ \Vert\, Y \ - \ \Xi_{\,\,1}\,\Vert\,)^{ {{\,n\,-\,2\,}\over 2} \,+\, {{n\,-\,2}\over {n\,+\,2}} \,\cdot\  \tau_{\,>1} }\ }}\\[0.15in]
   & \ & \hspace*{1in}\ \ \  \le \ C_3\,\cdot\,{1\over {\ (\ 1 \ + \ \Vert\, Y \ - \ \Xi_{\,\,1}\,\Vert\,)^{ {{\,n\,-\,2\,}\over 2} \,+\, {{n\,-\,2}\over {n\,+\,2}} \,\cdot\  \tau_{\,>1} }\ }}\,\cdot\,  {\bar\lambda}_{\,\,\flat}^{\, \left(\, {{\,n\,-\,2\,}\over 2} \,-\, {{n\,-\,2}\over {n\,+\,2}} \,\cdot\  \tau_{\,>1} \,\right)\,\cdot \,\gamma }\,  \\[0.2in]
  \Longrightarrow & \ &\hspace*{-0.4in} \left(\ \sum_{l\,=\,2}^\flat \ {1\over {(\,1 \ + \ \Vert\, Y \ - \ \Xi_{\,l}\,\Vert\,)^{\ n\,-\,2}}}\  \right)^{\!\!{{n\,+\,2}\over {\,n\,-\ 2\,}} }\!\!\!\! \le \     C_4 \cdot {\bar\lambda}_{\,\,\flat}^{\, \left(\, {{n}\over 2} \,-\,o_{\,+}\,(\,1\,) \,\right) \, \cdot  \, \gamma } \cdot\,{1\over {\ (\ 1 \ + \ \Vert\, Y \ - \ \Xi_{\,\,1}\,\Vert\,)^{ {{n\,+\,2}\over 2} \,+\, \tau_{\,>1}  } }} \ . \\[0.15in]
\end{eqnarray*}
See Note A.5.19\,.\,
 Summing up
$$ \R^n \ = \ \bigcup_{j\,=\,1}^\flat \Upsilon_j\ ,$$
we obtain\\[0.1in]
(\,A.5.20\,)
\begin{eqnarray*}
& \ & \left[ \ {\bf W}_\flat^{{\,n\,+\,2\,}\over {n\,-\,2}} \ - \  \sum_{l \,=\,1}^\flat\   {{\bf V}}_l^{{\,n\,+\,2\,}\over {n\,-\,2}} \ \right]_{\,Y} \ \ \ \ \ \ \ \left( \ \downarrow \ \ \  {{n\,+\,2}\over 2} \,+\, \tau_{\,>1}  \ = \ {{\,n\,-\,2\,}\over 2} \,+\, 2\,+\, \tau_{\,>1}  \ \right) \\[0.2in]
& \le  &  C\,\cdot\,{\bar\lambda}_{\,\,\flat}^{\ \left( {{n}\over 2} \ -\ o_{\,+}\,(\,1\,)  \right)\,\cdot\,\gamma} \,\cdot\, \left(\ \sum_{l\,=\,1}^\flat \  {1\over {(\ 1 \ + \ \Vert\, Y \ - \ \Xi_{\,l} \,\Vert\ )^{ {{\,n\,-\,2\,}\over 2}\ +\, 2 \,+\, \tau_{\,>1} }}}\  \right)\ \ \ \ {\mbox{for \ \ all}} \ \ \,Y\,\in\,\R^n\,.\,\
\end{eqnarray*}

\newpage

{\it Note.}\, A.5.21\,. \ \  The above consideration also furnish the reason for the bound\,:
$$
{\bf W}_{\,\flat} \,(\,Y\,) \ = \ \sum_{l\,= \,1}^\flat \left( {{\Lambda_{\,l}}\over {\ \Lambda^2_{\,l} \ + \ \Vert\ Y \,- \ \Xi_{\,l}\,\Vert^{\,2}\ }}\  \right)^{\!\!{{\,n\,-\,2\,}\over 2}} \ \ \le \ C_o \ \ \ \ \mfor \ \ Y \ \in \ \R^n.
$$
Here $\,C_0\,$ is independent on $\,\flat\,$ and the centers $\,\{ \ \Xi_{\,l}\ \}\,$ [\,see (\,A.5.21\,\,)\ ]  as long as conditions (\,1.8) and (\,1.22) in the original text are fulfilled\,,\, together with
$$
 {\bar{\lambda}}_{\,\flat} \ \le \ {\underline{\lambda}}_{\,\epsilon}\ .
$$
Intuitively, as
$$
0 \ <  \ {1\over { {\bar C}_1^{\ } }}  \ \le \ \Lambda_{\,l} \ \le \ {\bar C}_1 \ \ \ \ \ \ \ {\mbox{for}} \ \ \ l \ = \ 1\,,\, \ 2\,, \ \cdot \cdot \cdot\,, \ \flat\,,\,
$$
and the centers
$$
\Xi_{\,\,1}\,, \ \ \,\Xi_{2}\,, \    \,\cdot\,\cdot\,\cdot \ , \ \ \,\Xi_{\,j}\,, \ \,\cdot\,\cdot \cdot\,,\ \  \Xi_{\ \flat}\
$$
are ``\,widely\," separated from each other\,,\, the sum in (\,A.5.20\,) cannot be too large\,.\,\bk
Formally\,,\,  for $\,Y \,\in\, \Upsilon_1\,,\,$ we have
\begin{eqnarray*}
\Vert\ \Xi_{\,\,1} \ - \ \Xi_{\,l}\,\Vert & \le & \ \Vert\ Y \ - \ \Xi_{\,\,1}\,\Vert \ + \ \Vert\ Y \ - \ \Xi_{\,l}\,\Vert \ \le \  2 \,\Vert\ Y \ - \ \Xi_{\,l}\,\Vert\\[0.2in]
\Longrightarrow \ \ \Vert\ Y \ - \ \Xi_{\,l}\,\Vert & \ge & {1\over 2} \cdot  \Vert\ \Xi_{\,\,1} \ - \ \Xi_{\,l}\,\Vert \ \ \ \ \ \ \ \ \ \ \ ( \ {\mbox{for }} \ \ Y \,\in\, \Upsilon_1 \ ; \ \ l \ = \ 2\,, \ \cdot \cdot \cdot \ \flat\ )\,.
\end{eqnarray*}
Thus for $\,Y \,\in\, \Upsilon_1\,,\,$ we obtain
\begin{eqnarray*}
{\bf W}_{\,\flat}\,(\,Y\,) &  = & \sum_{l\ = \ 1}^\flat \left( {{\Lambda_{\,l}}\over {\ \Lambda^2_{\,l} \ + \ \Vert\ Y \ - \ \Xi_{\,l}\,\Vert^{\,2}\ }}\  \right)^{\!\!{{\,n\,-\,2\,}\over 2}} \\[0.2in]
 &  = & \left( {{\Lambda_{\,1}}\over {\ \Lambda^2_{\,1} \ + \ \Vert\ Y \ - \ \Xi_{\,\,1}\,\Vert^{\,2}\ }}\  \right)^{\!\!{{\,n\,-\,2\,}\over 2}} \ + \ \sum_{l\ = \ 2}^\flat \left( {{\Lambda_{\,l}}\over {\ \Lambda^2_{\,l} \ + \ \Vert\ Y \ - \ \Xi_{\,l}\,\Vert^{\,2}\ }}\  \right)^{\!\!{{\,n\,-\,2\,}\over 2}} \\[0.2in]
& \le & \left( {{1}\over {\ \Lambda_{\,1}  \ }}\  \right)^{\!\!{{\,n\,-\,2\,}\over 2}} \ + \  C_1 \cdot \sum_{l\ = \ 2}^\flat \left( {{1}\over {\ \Vert\ \Xi_{\,\,1} \ - \ \Xi_{\,l}\,\Vert^{\,2}\ }}\  \right)^{\!\!{{\,n\,-\,2\,}\over 2}}\\[0.2in]
& \le & C_2 \ + \ C_3 \cdot {\bar{\lambda}}_{\,\flat}^{\, \, (\,n\ - \ 2\,)\,\cdot\, \gamma} \ \le \ C_o\ \ \ \ \ \ \ \ \ \ \ ( \ {\mbox{for \ \ all }} \ \ Y \,\in\, \Upsilon_1 \ )\,.
\end{eqnarray*}
Similarly we work on $\,\Upsilon_{\,2}\,, \ \cdot \cdot \cdot \,, \ \Upsilon_{\,\flat}\ .$\, Here we see that $\,C_o\,$ depends on $\,{\underline{\lambda}}_{\,\epsilon}\,$,\, $\,{\bar C}_1\,$ and $\,\gamma\,.\,$

\newpage

${\bf (\,b\,)}_{\,(\,A.5.17\,)}$ \ \ {\it Flatness of $\,K\,.\,$ } \ \ Let us consider the term  (\ cf. \S\,A\,4\,g\,)
$$
 \bigg\vert \  (\,c_n\,K)\left(     \lambda \,\cdot\,Y   \right) \ - \ n\,(\,n\,-\,2)\ \bigg\vert \cdot\left( {{ \Lambda_{\,1}}\over { \Lambda_{\,1}^2 \ + \ \Vert\,Y \ - \ \Xi_{\,\,1}\,\Vert^{\,2}  }}\ \right)^{{n\,+\,2}\over 2}\ .
$$
It turns out that the following condition is more nature in this setting.

\vspace*{0.2in}

 {\bf  Uniform\ \  assumption}\,. \\[0.1in]
 (\,A..5.22\,)
$$ \ \
|\,(\,c_n\,\cdot\,K\,)\,(\,y\,)\,  \ - \ n\,(\,n\,-\,2)\,|\  \le\  C \ {\bar\lambda}_{\,\,\flat}^\zeta \ \ \  \ \ \ {\mbox{for  \ \ all}} \ \ y \,\in\,B_{\,\xi_{\,l}}\,(\,\rho_\nu\,)\ \ \ \ (\ l \,= \, 1\,, \ \cdot\,\cdot \cdot\,, \ \flat\ )\,.
$$

\vspace*{0.1in}

{\underline{\it Inside}}\,. \ \   Consider
$$
\Vert\,y \ - \ \xi_{\,1}\,\Vert \ < \ \rho_{\,\nu}\ \ \ \ \Longleftrightarrow \ \  \Vert\,Y \ - \ \Xi_{\,\,1}\,\Vert \ <  \    {\rho_{\,\nu}\over {  {\lambda_{\,\flat}} }}\,.\leqno (\,A.5.23\,)
$$
Under the uniform condition, we have
  \begin{eqnarray*}(\,A.5.24\,) \ \ \ \
  & \ & |\,(\,c_n\,\cdot\,K\,)\,(\,y\,)\,  \ - \ n\,(\,n\,-\,2)\,|\  \le C \ {\bar\lambda}_{\,\,\flat}^\zeta \ \ \  \ \mfor \ \ y \,\in\,B_{\xi_{\,1}}\,(\,\rho_\nu\,)\\[0.2in]
\Longrightarrow \  \  & \ & \bigg\vert \, (\,c_n\,K)\left(     \lambda \,\cdot\,Y   \right) \ - \ n\,(\,n\,-\,2)\,\bigg\vert \cdot\left( {{ \Lambda_{\,1}}\over { \Lambda_{\,1}^2 \ + \ \Vert\,Y \ - \ \Xi_{\,\,1}\,\Vert^{\,2}  }}\ \right)^{{n\,+\,2}\over 2}\\[0.2in]
& \le & C_1\,\cdot\,\bigg\vert \  (\,c_n\,K)\left(     \lambda \,\cdot\,Y   \right) \ - \ n\,(\,n\,-\,2)\,\bigg\vert \cdot\left( {{ \Lambda_{\,1}}\over { \Lambda_{\,1}^2 \ + \ \Vert\,Y \ - \ \Xi_{\,\,1}\,\Vert   }}\ \right)^{{n\,+\,2} } \ \ \ \ \ \ \ \ \  \ \  \\[0.2in]
& \le & C_2 \,\cdot\,   {\bar\lambda}_{\,\,\flat}^\zeta   \,\cdot\, {{1}\over { (\,1 \ + \ \Vert\,Y \ - \ \Xi_{\,\,1}\,\Vert \,)^  {  {{n\,+\,2}\over 2} \,-\ \tau_{\,>1}}  }}\,\cdot\, {{1}\over { (\,1 \ + \ \Vert\,Y \ - \ \Xi_{\,\,1}\,\Vert \,)^  {  {{n\,+\,2}\over 2} \,+\ \tau_{\,>1}}  }}\\[0.15in]
& \le &  C_3\,\cdot\,  {\bar\lambda}_{\,\,\flat}^\zeta  \,\cdot\, {{1}\over { (\,1 \ + \ \Vert\,Y \ - \ \Xi_{\,\,1}\,\Vert \,)^  {  {{n\,+\,2}\over 2} \ +\,\tau_{\,>1}}  }} \\[0.15in]
& \le &  C_3\,\cdot\,  {\bar\lambda}_{\,\,\flat}^{\,\zeta}  \,\cdot\, {{1}\over { (\,1 \ + \ \Vert\,Y \ - \ \Xi_{\,\,1}\,\Vert \,)^  {  {{n\,+\,2}\over 2} \ +\ \tau_{\,>1}}  }} \ \ \  \ \ \ \ \ \ \ \mfor \ \ y \,\in\,B_{\,\xi_{\,1}}\,(\,\rho_\nu\,)\ .\\
\end{eqnarray*}

\newpage

{\underline{\it Outside}}\,. \ \  Next\,,\, consider the situation
$$
\Vert\,y \ - \ \xi_{\,1}\,\Vert \ >  \ \cdot\,\rho_\nu\ \ \ \ \Longleftrightarrow \ \  \Vert\,Y \ - \ \Xi_{\,\,1}\,\Vert \ > \   {{\rho_{\,\nu}}\over {  { {\bar\lambda}_{\,\,\flat}} }} \ \ (\ = \  {\bar\lambda}_{\,\,\flat}^{\,\nu \ - \ 1}\ \gg \ 1\ )\,. \leqno (\,A.5.25\,)
$$
We proceed with
\begin{eqnarray*}
(\,A.5.26\,) \ \ \ \  & \ &  \bigg\vert \, (\,c_n\,K)\left(     \lambda \,\cdot\,Y   \right) \ - \ n\,(\,n\,-\,2)\,\bigg\vert \cdot\left( {{ \Lambda_{\,1}}\over { \Lambda_{\,1}^2 \ + \ \Vert\,Y \ - \ \Xi_{\,\,1}\,\Vert^{\,2}  }}\ \right)^{{n\,+\,2}\over 2}\\[0.15in]
& \le & C_1\,\cdot\,{{1}\over { (\,1 \ + \ \Vert\,Y \ - \ \Xi_{\,\,1}\,\Vert \,)^{\,n\,+\,2} }}\\[0.1in]& \ &  \ \ \  \  \ \ \ \ \ \ \ \ \  [\ \,|\,(\,c_n\,\cdot\,K\,)\,(\,y\,)\,  \ - \ n\,(\,n\,-\,2)\,|\  \le\  C \ \ {\mbox{for \  \ all}} \ Y \,\in\,\R^n \ ]\ \ \ \  \ \ \ \  \ \ \ \  \\[0.1in]
& \le & C_1 \,\cdot\, {{1}\over { (\,1 \ + \ \Vert\,Y \ - \ \Xi_{\,\,1}\,\Vert \,)^  {  {{n\,+\,2}\over 2} \,-\ \tau_{\,>1}}  }}\,\cdot\, {{1}\over { (\,1 \ + \ \Vert\,Y \ - \ \Xi_{\,\,1}\,\Vert \,)^  {  {{n\,+\,2}\over 2} \,+\ \tau_{\,>1}}  }}\\[0.15in]
& \le & C_1\,\cdot\,\left(\, {{\bar\lambda}_{\,\,\flat}\over \rho_{\,\nu}} \,\right)^{\!\!{ {{n\,+\,2}\over 2} \,-\ \tau_{\,>1} }}  \times   {{1}\over { (\,1 \ + \ \Vert\,Y \ - \ \Xi_{\,\,1}\,\Vert \,)^  {  {{\,n\,-\,2\,}\over 2} \ +\,2\,+\ \tau_{\,>1}}  }} \\[0.15in]
& \le & C_1 \cdot{\bar\lambda}_{\,\,\flat}^{\,{\left(  {{n\,+\,2}\over 2} \,-\ \tau_{\,>1}\right)\,\cdot\,(\,1 \ - \ \nu) }}  \times   {{1}\over { (\,1 \ + \ \Vert\,Y \ - \ \Xi_{\,\,1}\,\Vert \,)^  {  {{\,n\,-\,2\,}\over 2} \ +\,2\,+\ \tau_{\,>1}}  }}\ \\[0.2in]
& \le & C_1 \cdot{\bar\lambda}_{\,\,\flat}^{\,{\left(  {{n}\over 2} \,-\ o_{\,+}\,(\,1\,)\right)\,\cdot\,(\,1 \ - \ \nu) }}  \times   {{1}\over { (\,1 \ + \ \Vert\,Y \ - \ \Xi_{\,\,1}\,\Vert \,)^  {  {{\,n\,-\,2\,}\over 2} \ +\,2\,+\ \tau_{\,>1}}  }}\ \\[0.2in]
& \ & \ \ \ \  \ \ \ \  \ \ \ \  \ \ \ \ \ \ \ \ \ \ \ \ \ \ \ \ \ \mfor Y\,\in\,\R^n \ \ {\mbox{satisfying}} \ \ \Vert\,Y \ - \ \Xi_{\,\,1}\,\Vert \ > \   {{\rho_{\,\nu}}\over {  { {\bar\lambda}_{\,\,\flat}} }}\ .
\end{eqnarray*}
Refer to Note A.5.19\,.\,
 Likewise\,,\, we consider similar situation where the center $\ \Xi_{\,\,1}\ $ is replaced by $\ \Xi_{\,l}\,,\,$  for $\ l \ = \ 2\,, \ 3\,, \ \cdot\,\cdot\,\cdot\,\ , \ \flat\ .$\,\bk
Summing up the results in (\,A.5.24\,) and (\,A.5.26\,)\,,\, we obtain\\[0.1in]
(\,A.5.27\,)
\begin{eqnarray*}
|\ (\,c_n\,\cdot\,K) \ - \ n\,(\,n\,-\,2)\,|\,\cdot\, \left(\ \sum_{l\ =\,1}^\flat  V^{{\,n\,+\,2\,}\over {n\,-\,2}}_{\,l } \ \right)  & \le & C \,\cdot\,{\bar\lambda}_{\,\,\flat}^{\,{\vartheta}_{\,(\,A.5.7\,)} }\,\cdot\,  \left(\ \ \sum_{j\,=\,1}^\flat {1\over {   (  1\ + \ \Vert\,Y \,-\,\Xi_{\,l}\,\Vert\,)^{ {{n\,+\,2}\over 2}     \,\ +\ \,\tau_{\,>1}  }} }\right)  \ .
\end{eqnarray*}
Here
\begin{eqnarray*}
{\vartheta}_{\,(\,A.5.7\,)}  
 & = &  \mbox{Min}  \ \left\{ \    \left[\  {{n}\over 2} \,-\ o_{\,+}\,(\,1\,)\,\right]\,\cdot\, (\,1\ - \ \nu\,) \,, \ \ \ \zeta \  \right\}\ .
\end{eqnarray*}

\newpage

{\bf  Distance assumption.} \ \ Recall (\,1.26\,) in the main text\,\,:\\[0.1in]
(\,A.5.28\,)
$$
(\,{\tilde c}_n\,K\,)\,(\,y\,) \,=\,n\,(\,n\,-\,2\,) \ - \ C\,(\,{\bf p}_{\,y}\,)\cdot \Vert\,y \,-\,{\bf p}_{\,y}\,\Vert^{\,\ell} \ + \  {\bf R}_{\,\ell \  + \ 1}(\,y\,) \ \  \mfor \ \ y \ \in \ {\cal O}\,.
$$
Here
$$
 {\bf R}_{\,\ell \  + \ 1}(\,y\,)  \ \le \ C \,\Vert\,y \,-\,{\bf p}_{\,y}\,\Vert^{\,\ell\ +\,1}\ \ \  \mfor \ \ y \ \in \ {\cal O}\,,
$$
and the integer $\,\ell\,$ satisfies
$$\,\ell \ \in \ [\ 2\,, \ \, n \ - \ 2\ )\,. $$
Moreover, the open set $\,{\cal O}\,$ is chosen to be ``\,staying close\," to $\,{\cal H}\,$ so that the ``\,projection\," of $\,y\,$ unto $\,{\cal H}\,$,\, denoted by $\,{\bf p}_{\,y}\,\in\,{\cal H}\,,\,$ is uniquely defined\,:
$$
\,{\bf p}_{\,y} \ \in  \,{\cal H} \ \ \ {\mbox{so \ \  that \ \ Dist}}\  (\,y\,,\ {\cal H}\,) \ = \ \Vert  \,y \ - \ {\bf p}_{\,y}\ \Vert \,.\,
$$
The center of each bubble (\,that is\,,\, \,$\xi_{\,l}$\,) is located around $\,{\cal H}\,$,\, properly,
$$
   {\mbox{dist}}\,(\,\xi_{\,l}\,,\ {\cal H}\,) \ \le  \   {\bar\lambda}_{\,\,\flat}^{\,1\ + \ \kappa}\ \ \  \ \ {\mbox{for}} \ \ \  l \ = \ 1\,, \ 2\,,\cdot \cdot \cdot\,, \ \flat\,, \ \ \ \ \ \   {\mbox{where}} \ \ \kappa \, \in  \, (\,0\,,\   1\,) \ \ {\mbox{is \ \ fixed }}\,.\leqno (\,A.5.29\,)
$$
By choosing $\,{\bar{\lambda}}_{\,\flat}\,$ smaller (\,hence smaller $\,\rho_{\,\nu}\,$)\,,\, we may take it that
$$
\Vert\ y \ - \ \xi_{\,1}\ \Vert \ < \ \rho_{\,\nu} \ \ \Longrightarrow \ \ \Vert\ {\bf p}_{\,y} \ - \ \xi_{\,1}\ \Vert \ < \ 2\,\rho_{\,\nu}\ . \leqno (\,A.5.30\,)
$$
Similar expressions for $\xi_{\,l}$\,,\, $l \ = \   2\,,\cdot \cdot \cdot\,, \ \flat\,.\,$
Refer to \,{\bf \S\,4\,c}\, in the main text\, . \bk
On account of the Triangle Inequality\,,\, we obtain
$$
\Vert\ y \ - \ {\bf p}_{\,y}\ \Vert \ \le \ \Vert\ y \ - \ \xi_{\,1}\ \Vert \ + \ \Vert\ {\bf p}_{\,y} \ - \ \xi_{\,1}\ \Vert  \ \le \  3\,\cdot\,\rho_{\,\nu}\ .
$$
It follows that
\begin{eqnarray*}
(\,A.5.31\,) \ \ \   & \ & \bigg\vert \  (\,c_n\,K\,)\left(     \lambda \,\cdot\,Y   \right) \ - \ n\,(\,n\,-\,2)\,\bigg\vert\ = \  \bigg\vert \  (\,c_n\,K)\left(     y  \right) \ - \ n\,(\,n\,-\,2)\,\bigg\vert \ \ \ \ \ \ \ \ \ \ \ \\[0.2in]
  &  \le &  C_1\,\cdot\,\Vert\   y \ - \ {\bf p}_{\,y}\,\Vert^\ell\\[0.2in]
  &  \le &   C_2\,\cdot\,\rho_{\,\nu}^\ell \ \le \  C_2 \cdot {\bar{\lambda}}_{\,\flat}^{\, \ell\  \cdot \ \nu}\\[0.1in]
   & \ & \ \ \ \ [\ {\mbox{for}} \ \ {\bar{\lambda}}_{\,\flat} \ \ {\mbox{to \ \ small\,;\ \ depending \ \ on \ \ the \ \ error \   \ term \ \ in (\,A.5.28\,)}}\ ] \\[0.1in]
 & \ & \ \ \ \ \  \ \ \ \ \ \ \ \ \ \ \ \ \ \ \ \ \ \ \ \  \ \ \ \ \ \ \ \ \ \  \mfor \ \ {\bar\lambda}_{\,\,\flat} \cdot Y \ = \ y \,\in\,B_{\,\xi_{\,1}}\,(\,\rho_{\,\nu}\,)\ .
  \end{eqnarray*}
  Similar expressions for $\xi_{\,l}$\,,\, $l \ = \   2\,,\cdot \cdot \cdot\,, \ \flat\,.\,$ Thus we can take $$\,\zeta \ = \ \ell \,\cdot\, \nu\ \ \ \ \  {\mbox{in \ \ \ (\,A.5.27\,)}}\ . $$
This leads to
$$
{{\vartheta}_{\,(\,A.5.7\,)} }
 \ = \ \mbox{Min}  \, \left\{ \    {{n}\over 2}  \,\cdot\,(\,1\, - \ \nu\,)\,, \  \ \ \  \ell \cdot \nu \  \right\} \ \ - \ o_{\,+}\,(\,1\,)\,. \leqno (\,A.5.32\,)
$$
Finally we obtain \,${\bf (ii)}_{(\,5.6)}\,$ by combining \,(\,A.5.27\,)\, and  \,(\,A.5.31\,)\, with (\,A.5.20\,)\,,\, after  observing
\begin{eqnarray*}
  \gamma \ + \ \nu \ >  \ 1 \ \ \Longrightarrow \ \ \gamma \ >  \ 1 \ - \ \nu \ \ \
\Longrightarrow  \ \   {{n}\over 2}  \cdot \gamma \ > \   {{n}\over 2} \cdot  (\,1\ - \ \nu\,)\ .
\end{eqnarray*}
This comes to the end of the proof of Lemma A.\,5.6\,. \qedwh
{\it Note}\, A.5.33\,.  \ \ It can be perceived that
$$
{{\vartheta}_{\,(\,A.5.7\,)} } \ \le  \ \ell \cdot \nu \ \ \Longrightarrow \ \ {{\vartheta}_{\,(\,A.5.7\,)} } \ <  \ \ell \ \ \ \ \ (\,{\mbox{as}} \ \ \nu \ <  \ 1\,)\ .
$$

 \vspace*{0.5in}

{\bf \S\,A\,5\,.\,b. \ \ Contraction mapping.}   \ \ See the proof of Proposition 2.3 in \cite{Wei-Yan}\,.\, Based on Proposition A.4.20\,,\, for  a given $\,{\cal H}\,\in\,{\cal D}^{1\,,\ 2}\,,\,$ the equation
\begin{eqnarray*}
   & \ & \left[ \ \Delta  \ + \ {{n\,+\,2}\over {\,n\,-\ 2\,}} \cdot (\,c_n \cdot K) \cdot {\bf W}_{\flat}^{4\over{n\,-\,2}}\ \right]\, \Phi \ = \  \left(\,V_{1\,,\  o}^{4\over {n\,-\,2}} \cdot {\cal H} \right) \ + \ {  {\bf P}}_{\flat_{\,\parallelsum}} \
\end{eqnarray*}
has a  unique solution $\,\Phi \, \in \, {\cal D}^{\,1\,, \ 2}_{{\flat}}\,(\,\perp\,)\ ,\, $
\begin{eqnarray*}
{\mbox{where}}  & \ &  {  {\bf P}}_{\flat_{\,\parallelsum}}\ = \   \sum_{l\,=\,1}^\flat \ \left\{ \ \sum_{j\,=\,1}^n  b_{\,l\,,\ j} \cdot {\bf V}_{\,l}^{4\over {n\,-\,2}} \cdot    \left[\, (\,\Lambda_{\,l} \ \partial_{\,\Xi_{\,l_{|_{\,j}}}}\,) \,{\bf V}_{\,l}\,\right]  \ + \     a_{\,l} \cdot {\bf V}_{\,l}^{4\over {n\,-\,2}} \cdot  \left[\, (\,\Lambda_{\,l} \ \partial_{\Lambda_{\,l}}\,) \,{\bf V}_{\,l}\,\right]\  \right\}\ ,
\end{eqnarray*}
and the coefficients $\,a_{\,l}\,$ and $\ b_{\,l\,,\ j}\,$ are uniquely determined by $ \,{\cal H}\,$ [\,refer to  \,(\,A.4.19\,)\,]\,.\,  With these notations, we define the relation $\,{\cal L}_{\,\flat}\,$ given by
$$
\Phi \ = \ {\cal L}_{\,\flat}\,(\,{\cal H}\,)\,. \leqno (\,A.5.34\,)
$$
 ${\cal L}_{\,\flat}\,$ can be considered as a kind of ``\,(\,left\,) inverse\," of the linear operator
$$\left[ \ \Delta \ + \   \left({{n\,+\,2}\over {\,n\,-\ 2\,}}\ \right) \,\cdot\,( \,c_n\,\cdot\,K)\cdot {\bf W}_{\flat}^{4\over{n\,-\,2}}  \ \right]  \,. $$
Moreover, from Proposition {\bf A.4.20} again\,,\,
\begin{eqnarray*}
(\,A.5.35\,) \ \ \ \ \ \ \ \ \ \ \ \ \ \ \ \ \ \ \ \   & \ &   \Vert\,{\cal L}_{\,\flat}\,(\,{\cal H}\,)\,\Vert_{\,*} \ = \ \Vert\,\Phi\,\Vert_{\,*}\ \le \ C\,\cdot\,\Vert\,H\,\Vert_{\,**_{\,Y}}\\[0.2in]
 \ \ \ {\mbox{and}} \ \ \ \ & \ &  \Vert \, {\cal L}_{\,\flat}\,(\,{\cal H}\,) \,\Vert_\btd  \ = \ \Vert\,\Phi\,\Vert_\btd\ \le \ C\,\cdot\,\sqrt{\flat}\,\cdot\, \Vert\,{\cal H}\,\Vert_{\,**_{\,Y}} \ . \ \ \ \  \ \ \ \ \ \ \ \ \ \ \ \ \ \ \   \\
\end{eqnarray*}
We go ahead and described a subset $\,{\cal E} \,$ of a Banach space\,,\,  in which
\begin{eqnarray*}
(\,A.5.36\,) \ \ \ \ \ \ \ \ \ \ \ \ \  \ \   {\cal A}\  : \   {\cal E}  & \to &  {\cal E} \ \ \ \ \ \ \ \  \\[0.2in]
\Phi &\mapsto & {\cal L}_{\ \flat} \,(\,N\,(\,\Phi\,)\,) \ + \ {\cal L}_{\ \flat}\, (\, {\cal l}_{\ \flat}\,)   \ \ \ \ \ [ \ {\mbox{refer \ \ to \ \ (\,A.5.4\,)}} \ ] \ \ \ \ \ \ \ \ \ \ \ \ \ \ \ \ \ \ \ \ \ \ \ \ \ \  \ \ \ \ \ \ \ \ \ \ \ \ \
\end{eqnarray*}
is  well\,-\,defined\,,\, and is a {\it contraction mapping}\,.\bk
With this in mind, consider the Banach space
 $$
 {\cal B}  \   :=
  \  \left(\, {\bf W}_{* }^\flat \ \cap \ {\cal D}^{\,1\,,\,2}_{\flat_\perp}\, \right)
 $$
 with the norm given by
$$
\Vert\,\bullet\,\Vert_\oplus   :=    \Vert \,\bullet\,\Vert_{\,*} \ + \ \Vert \,\bullet\,\Vert_{\,\btd} \ . \leqno (\,A.5.37\,)
$$
Define the (\,fixed\,) index $\,{{\varpi}_{\,(\,A.5.38\,)} }\,$  by
$$
  {{\varpi}_{\,(\,A.5.38\,)} } \ = \ \left(\,{\vartheta}_{\,(\,A.5.7\,)}  \ - \ {1\over 2}\,\cdot\,\sigma\,\right) \ - \ o_{\,+}\,(\,1\,) \ . \leqno (\,A.5.38\,)
$$
We restrict ourselves to the metric (\,closed\,) ball
$$
{\cal E} \ := \ \left\{ \ \Psi \ \in \ {\cal B}  \  \ \bigg\vert \ \ \Vert \,\Psi\,\Vert_{*_Y} \ + \ \Vert\,\Psi\Vert_\btd \ \le \ {\bar\lambda}_{\ \flat}^{{{{\varpi}_{\,(\,A.5.38\,)} }} } \ \right\}\ \leqno (\,A.5.39\,)
$$
with respect to the norm $\,\Vert\,\bullet\,\Vert_\oplus\,$.\, Under the condition that
$$
\flat^{\,1\ +  \  {{n\,-\,2}\over 8}   }  \cdot \, \lambda^{ \,{{{\varpi}_{\,(\,A.5.38\,)} }} }_{\ \flat}  \ = \ o_{\,+}\,(\,1\,) \ , \leqno (\,A.5.40\,)
$$
As
$$
2 \ \le \ \flat \ \le\  {1\over { { {\bar\lambda}_{\,\,\flat}   }^{\,\sigma} }} \ ,
$$
the following condition is sufficient for (\,A.5.40\,)\,:\\[0.1in]
 (\,A.5.41\,)
$$
\left(\,{\vartheta}_{\,(\,A.5.7\,)}  \ - \ {1\over 2}\,\cdot\,\sigma\,\right) \ > \ \left( 1\ +  \  {{n\,-\,2}\over 8}  \right)\cdot \sigma \ \ \ \  \Longleftrightarrow \ \ \ \ {\vartheta}_{\,(\,A.5.7\,)}   \ > \ \left(\  {{ n \ + \ 10 }\over 8}\ \right)\, \cdot \sigma \ .
$$
We claim that [\,under the conditions in Lemma A.5.6 and (\,A.5.41\,)\ ] \\[0.1in]
(\,A.5.42\,)\\[0.1in]
$\, {\cal A} \,: \ {\cal E}  \ \to \ {\cal E}\,$ {\it is  well\,-\,defined and a contraction mapping with respect to the norm}\, $\,\Vert\,\bullet\,\Vert_\oplus\,$\emph{}\,.\\[0.1in]
 Once this is established\,,\, it follows that  there is a unique fixed point $\,\Phi\,\in\,{\cal E}\,$:
$$
 {\cal A} \,(\,\Phi\,) \ = \ \Phi\,.  \leqno (\,A.5.43\,)
$$
Reading backward [\,that is\,,\, via (\,A.5.3\,)\,,\, (\,A.5.4\,)\, and \,(\,A.5.34\,)\ ]\,,\,  $\,\Phi\,$ solves the equation
$$
\Delta\,(\,{\bf W}_{\ \flat} \ + \ \Phi\,)  \ + \ (\,c_n \,\cdot\, K)\,\cdot\,({\bf W}_{\ \flat} \ + \ \Phi)_+^{{\,n\,+\,2\,}\over {n\,-\,2}} \ = \ {\tilde{\bf P}}_{\flat_{\,\parallelsum}}\,, \ \ \ \ \ {\mbox{where}} \ \ \ \Phi \ \in \  {\cal D}^{\,1\,, \ 2}_{\,\flat_{\,\perp}}\,.
$$

\newpage

{\it Proof of the claim\,.\,  Step} 1\,.\ \ We begin with\\[0.1in]
(\,A.5.44\,)
\begin{eqnarray*}
\Vert \,{\cal A} \,(\Phi)\,\Vert_{\,*} & \le & \Vert \,{\cal L}_{\ \flat} \,({\bf N}\,(\Phi)) \ + \  {\cal L}_{\ \flat}\, ( {\cal l}_{\ \flat}) \,\Vert_{\,*} \\[0.2in]
& \le & \Vert \,{\cal L}_{\ \flat} \,({\bf N}\,(\Phi)) \,\Vert_{\,*}  \ + \ \Vert \,{\cal L}_{\ \flat}\, ( {\cal l}_{\ \flat}) \,\Vert_{\,*} \\[0.2in]
& \le & C_1\,\cdot\,[\ \Vert \,{\bf N}\,(\Phi) \,\Vert_{\,**} \ + \ \Vert \, {\cal l}_{\ \flat} \,\Vert_{\,**} \ ]  \ \ \ \ \ \ \ \ \ \ \ \  \ \ \ \  \ \ \ \ \ \ \ \ \ \ \ \ \ \ \ \ \ \ \ \ \ [ \ {\mbox{via}} \ \ (\,A.5.35\,) \ ] \\[0.2in]
& \le & C_2\,\cdot\,\left[\ \flat^{4\over {\ n\,-\ 2\ }}\,\cdot\,\Vert \, \Phi  \,\Vert_{\,*}^{{\ n\,+\,2\ }\over {n\,-\,2}} \ + \ {\bar\lambda}_{\ \flat}^{ \,{\vartheta}_{\,(\,A.5.7\,)} } \ \right]  \ \ \ \ \ \ \  \left[\ {\mbox{using}} \ \ {\bf (\,ii\,)}_{\mbox(\,5.6\,)} \ \ {\mbox{and}} \ \  {\bf (\,iii\,)}_{\mbox(\,5.6\,)} \ \right]\\[0.2in]
& \le & C_2\,\cdot\,\left( \ \flat\,\cdot\ \Vert \, \Phi  \,\Vert_{\,*} \ \right)^{4\over {n\ - \ 2}} \cdot  \Vert \, \Phi  \,\Vert_{\,*} \ + \  C_2 \cdot {\bar\lambda}_{\,\,\flat}^{ \,{\vartheta}_{\,(\,A.5.7\,)} }  \\[0.2in]
& \le & C_3\,\cdot\,\left( \ \flat\,\cdot\ \lambda^{ \,{{{\varpi}_{\,(\,A.5.38\,)} }} }_{\ \flat} \ \right)^{4\over {n\ - \ 2}} \cdot  \lambda^{ \,{{{\varpi}_{\,(\,A.5.38\,)} }} }_{\ \flat}  \ + \  C_2 \cdot {\bar\lambda}_{\ \flat}^{ {\vartheta}_{\,(\,A.5.7\,)} }   \ \ \ \ \ \ \ [\ {\mbox{via}} \ \ (\,A.5.40\,) \ ] \\[0.2in]
 & \le &   \lambda^{ \,{{{\varpi}_{\,(\,A.5.38\,)} }} }_{\ \flat} \ \ \ \mfor \ \ \Phi \ \in \  {\cal E} \ \ \ {\mbox{and}} \ \ {\bar\lambda}_{\ \flat} \ \ \ {\mbox{small \ \ enough}}\ .
\end{eqnarray*}
Here we apply the first part of (\,A.5.40\,)\,.\, Likewise\,,\,\\[0.1in]
(\,A.5.45\,)
\begin{eqnarray*}
\Vert \,{\cal A} \,(\Phi)\,\Vert_\btd  & \le & \Vert \  {\cal L}_{\ \flat} \,({\bf N}\,(\Phi)) \ + \ {\cal L}_{\ \flat}\, ( {\cal l}_{\ \flat}) \ \Vert_\btd
\\[0.2in] & \le & C_1\,\cdot\,\left[\ \Vert \,{\cal L}_{\ \flat} \,({\bf N}\,(\Phi)) \,\Vert_\btd   \ + \ \Vert \,{\cal L}_{\ \flat}\, ( {\cal l}_{\ \flat}) \,\Vert_\btd   \ \right]\\[0.2in]
& \le & C_2\,\cdot\,\sqrt{\flat\,}\,\cdot\,[\ \Vert \,{\bf N}\,(\Phi) \,\Vert_{\,**}  \ + \ \Vert \, {\cal l}_{\ \flat} \,\Vert_{\,**}  \ ]\\[0.2in]
& \le & C_3\,\cdot\,\sqrt{\flat\,}\,\cdot\,\left[\  \flat^{4\over {\ n\,-\ 2\ }}\,\cdot\Vert \, \Phi  \,\Vert_{\,*}^{  {{n\,+\,2}\over {\,n\,-\ 2\,}} } \ + \ {\bar\lambda}_{\ \flat}^{ \,{\vartheta}_{\,(\,A.5.7\,)} } \ \right] \\[0.2in]
& \le & C_4\,\cdot\,\left[ \ \sqrt{\,\flat\,} \,\,\cdot\  \left( \ \flat\,\cdot\ \Vert \, \Phi  \,\Vert_{\,*} \ \right)^{4\over {n\ - \ 2}} \ \right]\cdot  \Vert \, \Phi  \,\Vert_{\,*} \ + \  C_2 \cdot\sqrt{\,\flat\,} \,\,\cdot {\bar\lambda}_{\ \flat}^{ \,{\vartheta}_{\,(\,A.5.7\,)} }  \\[0.2in]
& \le & C_5\,\cdot\,\left( \ \flat^{ {{n\,-\,2}\over 8} \ + \ 1 }  \cdot \, \lambda^{ \,{{{\varpi}_{\,(\,A.5.38\,)} }} }_{\ \flat}  \ \right)^{4\over {n\ - \ 2}}  \cdot  \lambda^{ \,{{{\varpi}_{\,(\,A.5.38\,)} }} }_{\ \flat}    \ + \  C_3 \cdot{\bar\lambda}_{\ \flat}^{\  o_{\,+}\,(\,1\,)} \cdot{\bar\lambda}_{\ \flat}^{ {{\varpi}_{\,(\,A.5.38\,)} }} \ \ \ \ \ \ \ [\ {\mbox{via}} \ \ (\,A.5.40\,) \ ] \\[0.2in]
 & \le &   \lambda^{ {{{\varpi}_{\,(\,A.5.38\,)} }} }_{\ \flat} \ \  \mfor \ \Phi \ \in \ {\cal E} \ \  {\mbox{and}} \ \  {\bar\lambda}_{\ \flat} \ \ {\mbox{small \ \ enough}}\\[0.1in]
   & \ & \hspace*{1.8in} \left[\  {\mbox{recall \ \ that}} \ \ {{\varpi}_{\,(\,A.5.38\,)} } \, = \, \left(\,{\vartheta}_{\,(\,A.5.7\,)}  \ - \ {1\over 2}\,\cdot\,\sigma\,\right) \ - \ o_{\,+}\,(\,1\,) \ \right] \, {\bf .}
\end{eqnarray*}

\newpage

Recall from (\,A.5.3\,)   that
\begin{eqnarray*}
{\bf N}\,(\Phi) & = & \ - \ \left\{ (\,c_n\,\cdot\,K)\,\cdot\,\left[\ ({{\bf W}}_{\ \flat} \ + \ \Phi)^{{\,n\,+\,2\,}\over {n\,-\,2}} \ -\,{{\bf W}}_\flat^{{n\,+\,2}\over{n\,-\,2}} \ - \ \left(\ {{n\,+\,2}\over {\,n\,-\ 2\,}}\ \right)  \,\cdot\,{{\bf W}}_\flat^{4\over{n\,-\,2}}\cdot \Phi   \ \right] \ \right\}\,, \\[0.2in]
{\cal l}_{\ \flat}  & = & \ - \ \left\{ (\,c_n\,\cdot\,K)\,\cdot\, {{\bf W}}_\flat^{{\,n\,+\,2\,}\over {n\,-\,2}} \ -\, n\,(\,n\,-\,2)\,\cdot\,\sum_{l\,=\,1}^\flat \  {{\bf V}}_l^{{\,n\,+\,2\,}\over {n\,-\,2}}  \ \right\}.\\
\end{eqnarray*}
Hence we demonstrate that $\,{\cal A}\,$ as in (\,A.5.36\,) is well\,-\,defined.

\vspace*{0.2in}

{\it Proof of the claim\,.\,  Step} 2\,.\ \  For $\, n \ \ge \ 6\,,\,$ and $\,a \ > \ 0\,$ fixed\,,\, consider
$$
F\,(\,t\,) \ := \ (\,a \ + \ t\,)^{{\,n\,+\,2\,}\over {n\,-\,2}} \ - \ a^{{\,n\,+\,2\,}\over {n\,-\,2}} \ - \ \left(\ {{n\,+\,2}\over {\,n\,-\ 2\,}}\ \right)   a^{{4}\over {\,n\,-\,2\,}}\,\cdot\,t \ \ \ \ \mfor \ \ |\,t\,| \ = \ o_{\,+}\,(\,1\,)\ .
$$
We have
$$
|\,F'(t)\,| \ \le \ C\,\cdot\,|\,t\,|^{4\over {n\,-\,2}} \ \ \ \ \mfor \ \ |\,t\,| \ = \ o_{\,+}\,(\,1\,)\,. \leqno  (\,A.5.46\,)
$$
Cf. also (\,A.5.10\,)\,.\, As observed in Note A.5.21\,,\,  that
$$
(\, 0 \ < \ ) \ \ {{\bf W}}_{\ \flat}\,(\,Y\,) \ \le \ C_o \ \ \ \ \mfor\ \ {\mbox{all}} \ \ \ \ Y \ \in \ \R^n\,.  \leqno (\,A.5.47\,)
$$
Here we see that $\,C_o\,$ depends on $\,{\underline{\lambda}}_{\,\epsilon}\,$,\, $\,{\bar C}_1\,$ and $\,\gamma\,.\,$\,
Likewise (\,argue as in Note A.5.21\,)
$$
\Vert \, \Phi  \,\Vert_{\,*} \ \ \ \ {\mbox{is \ \ small \ \ }} \Longrightarrow \ \ |\, \Phi\,(\,Y\,)\,| \ \ \ \ {\mbox{is \ \ uniformly \ \  small  }} \mfor\ \ {\mbox{all}} \ \ Y \ \in \ \R^n\ .  \leqno (\,A.5.48\,)
$$

Recall that
$$
\Vert \, {\bf N}\,(\Phi_1) \ - \  {\bf N}\,(\Phi_2)\,\Vert_{\,**} \ = \ \sup_{Y\,\in\,\R^n}   \left\{ \, {{ |\,{\bf N}\,(\Phi_1) \ - \  {\bf N}\,(\Phi_2)\,|_{\,Y}}\over {\
  \ {\displaystyle{\sum_{l\,=\,1}^\flat}} \left(\  {1\over {  \  (  \,1\ + \ \Vert\ Y \,-\ \Xi_{\,l}\ \Vert\,)^{ {{\,n\,-\ 2\,}\over 2 } \, +\ 2\ +\ \tau_{\,>1}}\  }}\ \right)\,
 \  }} \,\right\}\ .
$$

\newpage

We continue with\\[0.1in]
(\,A.5.49\,)
\begin{eqnarray*}
& \ &  {{ |\ {\bf N}\,(\Phi_1) \ - \  {\bf N}\,(\Phi_2)\,|_{\,Y} }\over{\
  \ {\displaystyle{\sum_{l\,=\,1}^\flat}} \left(\  {1\over {  \  (  \ 1\ + \ \Vert\,Y \,-\,\Xi_{\,l}\,\Vert\,)^{ {{\,n\,-\ 2\,}\over 2 } \, +\ 2\,+\ \tau_{\,>1}}\  }}\ \right)\,
 \  }} \\[0.25in]
  & \le &  {  { \ C_1\,\cdot\,\left[\ |\, \Phi_1\,(\,Y\,) \,|^{4\over {n\,-\,2}} \ + \  |\, \Phi_2 \,(\,Y\,) \,|^{4\over {n\,-\,2}} \ \right\}\,\cdot\, |\, \Phi_1\,(\,Y\,) \ - \ \Phi_2\,(\,Y\,) \,| } \over{\
  \ {\displaystyle{\sum_{l\,=\,1}^\flat}} \left(\  {1\over {  \  (  \ 1\ + \ \Vert\,Y \,-\,\Xi_{\,l}\,\Vert\,)^{ {{n\,-\,2}\over 2 } \, +\ 2\,+\ \tau_{\,>1}}\  }}\ \right)\,
 \  }}\\[0.1in]
  & \ & \hspace*{1.75in} \ \ \ \ \ \ [ \ {\mbox{via \ \ (\,A.5.46\,)\,, \ \ (\,A.5.47\,)  \ \ and \ \ (\,A.5.48\,)}} \ ]\\[0.1in]
& \le & C\,\cdot\,\left[\ \Vert \, \Phi_1 \,\Vert_{\,*}^{4\over {n\,-\,2}} \ + \  \Vert\, \Phi_2 \,\Vert_{\,*}^{4\over {n\,-\,2}} \ \right]\,\cdot\, \Vert\, \Phi_1 \ - \ \Phi_2 \,\Vert_{\,*} * \\[0.25in]
& \ & \hspace*{2.5in} * \, {{\  \left\{\ {\displaystyle{\sum_{l\,=\,1}^\flat}} \left(\ {1\over {   (  \ 1\ + \ \Vert\,Y \,-\,\Xi_{\,l}\,\Vert\,)^{ {{n\,-\,2}\over 2 } \,+\ \tau_{\,>1}} }}\ \right)\,
\right\}^{{\,n\,+\,2\,}\over {n\,-\,2}} \ \  }\over{\
  \ {\displaystyle{\sum_{l\,=\,1}^\flat}} \left(\  {1\over {  \  (  \ 1\ + \ \Vert\,Y \,-\,\Xi_{\,l}\,\Vert\,)^{ {{\,n\,-\ 2\,}\over 2 } \, +\ 2\,+\ \tau_{\,>1}}\  }}\ \right)\,
 \  }}\\[0.25in]
& \le & C\,\cdot\,\left[\ \Vert \, \Phi_1 \,\Vert_{\,*}^{4\over {n\,-\,2}} \ + \  \Vert\, \Phi_2 \,\Vert_{\,*}^{4\over {n\,-\,2}} \ \right]\,\cdot\, \Vert\, \Phi_1 \ - \ \Phi_2 \,\Vert_{\,*} * \\[0.1in]
 & \ & \hspace*{3in} \left[\  \downarrow \ \ {\mbox{see \ \ (\,A.5.11\,)\,; \ \ }} {\mbox{cf. \ \ (\,A.5.13\,)}}  \ \right] \\[0.1in]
& \ & \hspace*{2in} * \, {{ \  \left\{\ {\displaystyle{\flat^{4\over {n\,-\,2}}\,\cdot\,\sum_{l\,=\,1}^\flat}} \left(\ {1\over {   (  1\ + \ \Vert\,Y \,-\,\Xi_{\,l}\,\Vert\,)^{ {{n\,+\,2}\over 2 } \ +\  {{n\,+\,2}\over {\,n\,-\ 2\,}}\,  \cdot \, \tau_{\,>1}} }}\ \right)\,
\right\}  \  }\over{\
  \ {\displaystyle{\sum_{l\,=\,1}^\flat}} \left(\  {1\over {  \  (  \ 1\ + \ \Vert\ Y \,-\,\Xi_{\,l}\,\Vert\,)^{ {{\,n\,-\ 2\,}\over 2 } \, +\ 2\,+\ \tau_{\,>1}}\  }}\ \right)\,
 \  }}\\[0.2in]
 & \le & C_1\cdot \flat^{\ {4\over {n\,-\,2}}}\,\cdot\,\left[\ \Vert \, \Phi_1 \,\Vert_{\,*}^{4\over {n\,-\,2}} \ + \  \Vert\, \Phi_2 \,\Vert_{\,*}^{4\over {\,n\,-\ 2\,}} \ \right]\,\cdot\, \Vert\, \Phi_1 \ - \ \Phi_2 \,\Vert_{\,*} \\[0.2in]
  & \le & C_2\,\cdot\,\left( \ \flat\,\cdot\ \lambda^{ \,{{{\varpi}_{\,(\,A.5.38\,)} }} }_{\ \flat} \ \right)^{4\over {n\ - \ 2}}\,\cdot\, \Vert\, \Phi_1 \ - \ \Phi_2 \,\Vert_{\,*} \ .
\end{eqnarray*}
As a result, when $\,\lambda_{\ \flat} \,$ is small enough\,,\, together with (\,A.5.40\,)\,,\, we have
$$
  \Vert \,{\cal L}_{\ \flat} \,({\bf N}\,(\Phi_1)) \ - \ {\cal L}_{\ \flat} \,({\bf N}\,(\Phi_2)) \,\Vert_{\,*} \ < \ {\bar c}\,\cdot\, \Vert\, \Phi_1 \ - \ \Phi_2 \,\Vert_{\,*}\ \ \ \ \ ( \, {\bar c} \ <  \ 1\,)\ . \leqno (\,A.5.50\,)
$$
Likewise [\ for the gradient estimate\,;\, cf. (\,A.5.45\,) \ ]\\[0.1in]
(\,A.5.51\,)
\begin{eqnarray*}
\Vert \,{\cal L}_{\ \flat} \,({\bf N}\,(\Phi_1)) \ - \ {\cal L}_{\ \flat} \,({\bf N}\,(\Phi_2)) \,\Vert_{\,\btd} & = &  \Vert \,{\cal L}_{\ \flat} \,[\,{\bf N}\,(\Phi_1) \ - \  {\bf N}\,(\Phi_2) \,]  \,\Vert_{\,\btd} \\[0.2in]
& \le & C\,\cdot\,\sqrt{\flat\,}\,\cdot\, \Vert \, {\bf N}\,(\Phi_1) \ - \  {\bf N}\,(\Phi_2)\,\Vert_{\,**} \\[0.15in]
 & \le & C'\,\cdot\,\sqrt{\flat\,}\,\cdot\,\left[\ \Vert \, \Phi_1 \,\Vert_{\,*}^{4\over {n\,-\,2}} \ + \  \Vert\, \Phi_2 \,\Vert_{\,*}^{4\over {n\,-\,2}} \ \right]\,\cdot\, \Vert\, \Phi_1 \ - \ \Phi_2 \,\Vert_{\,*}\\[0.2in]
 & < & {\bar c}\,\cdot\, \Vert\, \Phi_1 \ - \ \Phi_2 \,\Vert_{\,*}\ \ \ \ \   [\ {\bar c} \ <  \ 1\,,\, \ \ \ \ {\mbox{via \ \ in}}  \ \ (\,A.5.4\,)\ ]\ .
\end{eqnarray*}
Here we use (\,A.5.40\,)\,.
Combining (\,A.5.49\,) and (\,A.5.51\,)\,,\, we obtain
$$
\Vert\, {\cal A} \,(\Phi_1) \ - \ {\cal A} \,(\Phi_2)\,\Vert_{\,*} \ + \ \Vert\, {\cal A} \,(\Phi_1) \ - \ {\cal A} \,(\Phi_2)\,\Vert_\btd \ < \  {\bar c}\,\cdot\,\left[\ \Vert\,  \Phi_1  \ - \  \Phi_2 \,\Vert_{\,*} \ + \ \Vert\,  \Phi_1  \ - \  \Phi_2 \,\Vert_\btd \, \right]\ .
$$
 Thus we arrive at a proof to claim A.5.42\,.

\vspace*{0.2in}

{\bf \S\,A\,5\,.\,d\,.} \ \    {\bf Estimate on the coefficient $\,c_{\,l}\,$ and  $\,c_{\,l\,, \ j}\,$  } \ \   From (\,A.5.36\,) and (\,A.5.43\,)\,,\, Lemma A.5.6 and condition (\,A.5.40\,)\,,\, we have
\begin{eqnarray*}
\Phi & = & {\cal L}_{\ \flat} \,(\,{\bf N}\,(\Phi)\,) \ + \ {\cal L}_{\,\flat}\, (\, {\cal l}_{\ \flat}\,)\ , \\[0.2in]
{\mbox{and}} \ \ \ \ \ \Vert\,\Phi\,\Vert_{*_Y} & \le &  {\bar\lambda}_{\ \flat}^{{\varpi}_{\,(\,A.5.38\,)} } \
\ \ \ \Longrightarrow  \ \ \
\Vert\,{\bf N}\,(\Phi))\,\Vert_{**_Y} \  =  \ O\,(\,{\bar\lambda}_{\ \flat}^{{\varpi}_{\,(\,A.5.38\,)} })
 \ \ \ \ {\mbox{and}} \ \ \ \
\Vert\, {\cal l}_{\ \flat} \,\Vert_{**_Y} \ = \  O\,(\,{\bar\lambda}_{\ \flat}^{{\varpi}_{\,(\,A.5.38\,)} } \,)\ .\\
\end{eqnarray*}
$\,{\cal L}_{\ \flat} \,({\bf N}\,(\Phi))\,$ and $\, {\cal L}_{\,\flat}\, ( {\cal l}_{\ \flat})\,$ solve the following equations\,,\, respectively\,,\,\\[0.1in]
(\,A.5.52\,)
\begin{eqnarray*}
& \ & \Delta_Y \,(\ {\cal L}_{\ \flat} \,({\bf N}\,(\,\Phi\,))\ ) \ + \ \left\{ \  {{n\,+\,2}\over {\,n\,-\ 2\,}}\,\cdot\,(\,c_n\,\cdot\,K)\,\cdot\,{\bf W}_\flat^{4\over{n\,-\,2}} \ \right\}\,\cdot\,(\,{\cal L}_{\ \flat} \,({\bf N}\,(\,\Phi\,)\ )\\[0.2in]
 & = &   {\bf N}\,(\Phi) \ + \  \sum_{l\,=\,1}^\flat\  \left\{ \ \alpha_{\,l}\,\cdot\, {\bf V}_{\,l}^{4\over {n\,-\,2}}\,\cdot\, \left[\, (\,\Lambda_{\,l} \ \partial_{\Lambda_{\,l}}\,) {\bf V}_{\,l}\,\right] \ + \ \sum_{j\,=\,1}^n  \beta_{\,l\,,\,j}\,\cdot\,  {\bf V}_{\,l}^{4\over {n\,-\,2}}\,\cdot\, \left[\, (\,\Lambda_{\,l} \ \partial_{\,\Xi_{\,l_{|_{\,j}}}}\,)\, {\bf V}_{\,l}\,\right] \,     \right\}\ ,\\
\end{eqnarray*}

\newpage

(\,A.5.53\,)
\begin{eqnarray*}
& \ & \Delta_Y \,(\ {\cal L}_{\,\flat}\, (\, {\cal l}_{\ \flat}\,)\ ) \ + \  \left\{ \ {{n\,+\,2}\over {\,n\,-\ 2\,}}\,\cdot\, (\,c_n\,\cdot\,K)\,\cdot\,{\bf W}_\flat^{4\over{n\,-\,2}} \ \right\}\,\cdot\,(\ {\cal L}_{\,\flat}\, ( \  {\cal l}_{\ \flat}\,)\ )\\[0.2in]
 & = &   {\cal l}_{\ \flat} \ + \  \sum_{l\,=\,1}^\flat \ \left\{ \ {\tilde\alpha}_{\,l}\,\cdot\, {\bf V}_{\,l}^{4\over {n\,-\,2}}\,\cdot\, \left[\, (\,\Lambda_{\,l} \ \partial_{\Lambda_{\,l}}\,) {\bf V}_{\,l}\,\right] \ + \ \sum_{j\,=\,1}^n  {\tilde\beta}_{\,l\,,\,j}\,\cdot\,  {\bf V}_{\,l}^{4\over {n\,-\,2}}\,\cdot\, \left[\, (\,\Lambda_{\,l} \ \partial_{\,\Xi_{\,l_{|_{\,j}}}}\,)\, {\bf V}_{\,l}\,\right] \,     \right\}\ .
 \end{eqnarray*}
 Adding up (\,A.5.52\,) and (\,A.5.53\,)\,, together with the first half of (\,A.3.48\,)\,,\, we obtain
$$
|\,c_{\,l} \,| \ + \ |\,c_{\,l\,, \ j}\,| \ \le \ |\,\alpha_{\,l}\,|  \ + \ |\,{\tilde\alpha}_{\,l}\,| \ + \ |\,\beta_{\,l\,, \ j}\,|  \ + \ |\,\tilde\beta_{\,l\,, \ j}\,|  \ \le \ C\,\cdot\, {\bar\lambda}_{\ \flat}^{{\varpi}_{\,(\,A.5.38\,)} }\ . \leqno (\,A.5.54\,)
$$


\hspace*{0.5in}Combining the discussion in {\bf \S\,A\,5\,.\,a}\,--\,{\bf d}\,,\,
we are now ready to ascertain our next result (\ cf. Proposition 2.3 in \cite{Wei-Yan}\ ) \ .

\vspace*{0.2in}

{\bf Proposition A.5.55\,.} \ \ {\it Under the conditions in Proposition\,} {\bf A.4.20}\,,\,  {\it together with}
$$
 {\vartheta}_{\,(\,A.5.7\,)}  \ > \  \left(\  {{ n \ + \ 10 }\over 8}\ \right)\, \cdot \sigma  \ \ \left( \  \Longrightarrow  \ \ {\vartheta}_{\,(\,A.5.7\,)}  \ >  \ 2\, \sigma\,, \ \ as \  \ n \ \ge \ 6 \ \right) \ . \leqno (\,A.5.56\,)
$$
{\it Here}
$$
{{\vartheta}_{\,(\,A.5.7\,)} }
 \ = \ \mbox{Min}  \ \left\{ \   {{n}\over 2}  \,\cdot\,(\,1\, - \ \nu\,)\,, \  \ \ \  \ell    \cdot  \nu \  \right\} \ -\ \  o_{\,+}\,(\,1\,)\,, \ \ \ \ \ {\it{where}} \ \  \   \gamma \ > \ \nu \ \   {\it{and}} \ \ \gamma \ + \ \nu \ >  \ 1\,.
$$

 {\it There exist a positive number  $\,{\underline{\lambda}}_{\  \epsilon}\,$} (\ {\it made smaller\,, if necessary\,,\, than the one in Proposition}\, {\bf A.4.20}\,) {\it  and a positive constant $\,C\,$} [\,{\it independent on}\, $\,{\underline{\lambda}}_{\ \epsilon}\,,$\,  $\,\ell$ \,{\it and} \,\,$\nu$\ ]\ ,\, {\it such that for all $\,{\bar\lambda}_{\ \flat}\,$ satisfying}
 $$
   0\ <  \ {\bar\lambda}_{\ \flat} \ \le \  {\underline{\lambda}}_{\ \epsilon} \ ,
    $$
    {\it equation\,} (\,A.5.1\,) {\it has a unique} ({\it small\,}) {\it solution\,} $\,\Phi  \ \in \ {\bf W}^{\, \flat}_*\,\cap\,{\cal D}^{\,1\,,\,2}_\perp\,.\, $ {\it Moreover, }
$$
\Vert \, \Phi \,\Vert_\btd \ + \ \Vert \, \Phi \,\Vert_{\,*} \ \le \  C\,\cdot\,{\bar\lambda}_{\ \flat}^{{{{\varpi}_{\,(\,A.5.38\,)} }} }
$$
{\it and}
$$
|\ c_{\,l} \,| \ + \ |\ c_{\,l\,, \ j}\,|  \ \le \ C\,\cdot\, {\bar\lambda}_{\ \flat}^{{\varpi}_{\,(\,A.5.38\,)} }\ \ \ \ \ \ \ \  {\it for} \ \ \  l \ = \ 1\,, \ \cdot \cdot \cdot\,, \ \flat\,, \ \ \ \ j \ = \ 1\,, \ \cdot \cdot \cdot\,, \ n\,.
$$
{\it In addition\,,\, $\Phi$ depends on the bubble parameters} $\,(\,\Lambda_{\,l}\ ,\ \Xi_{\,l}\,)_{\ 1\ \le\ l\ \le \ \flat}$\, [\ {\it fulflling the stated conditions}\ ] {\it \,in a $\,C^{\,1}$\,-\,manner\,.}

\vspace*{0.2in}

\hspace*{0.5in}For the $C^1$\,-\,conclusion\,,\, we refer the readers to  the argument in the proof of Proposition 4.2 in \cite{Pino}\,.\,

\newpage

{\bf \S\,A\,5\,.e\,.} \ \    {\it Some properties we know about}\, $\Phi\,.\,$ See also {\bf \S\,A\,4\,.k} and Remark A.5.5\,. \\[0.1in]
(\,a) \ \ $\displaystyle{\Phi \,\in\, {\bf W}^{\, \flat}_* \ \ \Longrightarrow \ \ \Phi\,(\,Y\,) \ =  \ O\left( \ {1\over { R^{ {{\,n\,-\,2\,}\over 2} \ +\, 1 \ }  }}\ \right)}\ $ \ \ for \ \ $\,R\,=\,\Vert\,Y\,\Vert \ \gg \ 1\,.\,$\\[0.1in] Recall that $$\tau_{\,>\,1} \ = \ 1\  + \ o_{\,+}\,(\,1\,)$$
(\,b) \ \ By applying the Sobolev embedding theorem, Schauder theory, and a boot\,-\,argument\,,\, we obtain $\,\Phi\,\in \,C^{\,2\,,\,\alpha}\ .$ Cf. See {\bf \S\,A\,4\,.k} and Remark A.5.5\,.\\[0.1in]
(\,c) \ \ Suppose that the arrangement of the bubbles is invariant under a rigid motion $\,{\bf{\cal T}}\,$,\, then (\,via uniqueness\,) $\,\Phi\,$ is also invariant under the same rigid motion\,.\, That is\,,\,
$$
\Phi\,(\,Y\,) \ = \ \Phi\,(\,{\bf{\cal T}}\,(\,Y\,)\,) \ \ \ \mfor\ \ Y \,\in\,\R^n\,.
$$

\vspace*{0.5in}

{\bf \S\,A\,5\,.\,f\,.} \ \  {\bf Transferring back.} \ \ Via the transformation on {\bf{\S\,A\,3\,.b}}\,,\, we obtain from Proposition A.5.55 the following (\ cf. Proposition 2.3 in \cite{Wei-Yan}\ )\,.\\[0.2in]
 {\bf Proposition A.5.57.} \ \ {\it Under the conditions in Proposition}\, {\bf A.5.55}\,  \,,\,
 {\it   equation\,} (2.1) {\it in the main text\,} [\,8\,] {\it  has a unique\,} (\,{\it small\,}) {\it \emph{}solution\,} $\,\phi_{\ \flat}  \ \in \ W^\flat_*\,\cap\,{\cal D}^{\,1\,,\,2}_\perp\,$ {\it satisfying}
 $$
 \Vert \, \phi_{\ \flat}  \,\Vert_\btd  \ \le \  C\,\cdot\,{\bar\lambda}_{\ \flat}^{{{\varpi}_{\,(\,A.5.38\,)} }}\,, \leqno{\bf
 {(i)_{\,(\,5.57\,)}}}
 $$
  $$
\Vert \, \phi_{\ \flat}  \,\Vert_{\,*} \ \le \  C\,\cdot\,{\bar\lambda}_{\ \flat}^{{{\varpi}_{\,(\,A.5.38\,)} }\ -\ \tau_{\,>1}}  \,, \leqno{\bf
 {(ii)_{\,(\,5.57\,)}}}
 $$
   $$ |\ c_{\,l} \,| \ + \
|\,c_{\ l_{\,,\,j}}\,|  \ \le \  C\,\cdot\,{\bar\lambda}_{\ \flat}^{{{\varpi}_{\,(\,A.5.38\,)} }}  \leqno{\bf
 {(iii)_{\,(\,5.57\,)}}}
 $$ {\it for}\, $\,j \ = \ 0\,, \ 1\,,\ \cdot\,\cdot \cdot\,\ n\,,$\,   {\it{and}} \ $l \ = \ 1\,,\ \cdot\, \cdot \cdot\,, \ \flat\,.\,$ {\it In addition\,,\, $\phi_{\ \flat}\,$  depends on the parameters} $\,(\,\lambda_{\,\,l}\ ,\ \xi_{\,l}\,)_{\,1\ \le\ l\ \le \ \flat}$\, [\ {\it fulfilling the condition in this Proposition \,}] {\it \,in a $\,C^1$\,-\,manner\,.}\\[0.2in]
 Recall from (\,A.5.38\,) that
 $$
 {{\varpi}_{\,(\,A.5.38\,)} } \ = \
 \mbox{Min}  \ \left\{ \     {{n}\over 2}  \,\cdot\,(\,1\, - \ \nu\,)\,, \  \ \ \  \ell    \cdot  \nu \  \right\} \ \ - \ {1\over 2} \cdot \sigma \ -\ \ o_{\,+}\,(\,1\,)\ .
 $$

\newpage

 {\large{\bf \S\,A\,6\,.\,} \   {\bf Reduced functional and a proof of Lemma 3.1.}}\\[0.2in]
Recall from (\,1.10\,) in the main text  the  functional
$$
\ {\bf I}\,(\,f\,)\ = \ {1\over 2}\ \int_{\R^n}\ \langle\,\btd\,f\,,\,\btd\,f\,\rangle\ -\ \left(\ {{n\,-\,2}\over {2n}}\,\right)\,\cdot\,\int_{\R^n}\,(\,{\tilde c}_n\!\cdot K\,)\,f_+^{{2n}\over {\,n\,-\ 2\,}}\ \ \ \mfor \ \  f\,\in\,{\cal D}^{\,1,\,\,2}\,, \leqno (\,A.6.1\,)
$$
and its the first Fr\'echlet derivative
$$
\ \ \ \ {\bf I}\,'\,(\,f\,)\,[\,h\,]\ = \ \int_{\R^n} \!\left[\ \langle\,\btd\,f\,,\,\btd\,h\,\rangle\ -\ (\,{\tilde c}_n\!\cdot K\,)\,f^{{\,n\,+\,2\,}\over{n\,-\,2}}_+\cdot\,h\ \right]\ \ \ \ \ \ \ \ \mfor\ \ h\ \in\ {\cal D}^{\,1,\,\,2}\,. \leqno (\,A.6.2\,)
$$
Let us recapture the following\,.\,
$$
\,W_{\,\,\flat}\,(\,\lambda_{\,1}\,, \ \cdot \cdot \cdot\,, \ \lambda_{\ \flat}\,; \ \xi_{\,1}\,, \ \cdot \cdot \cdot\,, \ \xi_{\ \flat}\,) \ \left( \ = \ \sum_{l \ = \ 1}^\flat V_{\,l} \ \right) \ = \ \sum_{l \ = \ 1}^\flat \left( \ {{\lambda_{\,\,l} }\over { \lambda_{\,\,l}^2 \ + \ \Vert\, y \ - \ \xi_{\,l} \,\Vert^2 }}\ \right)^{\!\!{{\,n\,-\,2\,}\over 2}} {\bf .}$$
Here the bubbles parameters are given by
 $$
 {\bf P} _{\ \flat} \ = \ (\,\lambda_{\,1}\,, \ \cdot \cdot \cdot\,, \ \lambda_{\ \flat}\,; \ \xi_{\,1}\,, \ \cdot \cdot \cdot\,, \ \xi_{\ \flat}\,)\,, \leqno (\,A.6.3\,)
 $$

 \vspace*{-0.12in}

 $$
(\,\lambda_{\,1}\,, \ \cdot \cdot \cdot\,, \ \lambda_{\ \flat}\, )\ \in \ \R^+ \times \cdot \cdot \cdot \times \R^+\  \ \ \ \ {\mbox{and}} \ \ \ \
(\,\xi_{\,1}\,, \ \cdot \cdot \cdot\,, \ \xi_{\ \flat}\,) \ \in \ \R^n\, \times \cdot \cdot \cdot \times \R^n\,.
$$
{\it They are supposed to fulfill the conditions in Proposition}\, A.5.57\,.\, Let $$\,\phi_{\,\,\flat}  \, \in \, W^\flat_{*} \,\cap \, {\cal D}^{1\,,\  2}_\perp\,$$  be the (\,unique\,) \,``\,small\,"\, solution  of equation (\,2.1\,) in the main text\, \,,\, granted by Proposition A.5.57\,.\, That is\,,\,
$$
\Delta\, (\,W_{\,\,\flat} \ + \ \phi_{\ \flat}\,) \ +  \ (\,{\tilde c}_n\!\cdot K\,)\,\cdot\,(\,W_{\,\,\flat} \ + \ \phi_{\ \flat}\,)^{{\,n\,+\,2\,}\over {n\,-\,2}} \ = \ {\bf P}_{/\!/} \ ,
$$
where $\,{\bf P}_{/\!/} \,$ is given by (\,A.5.2\,)\,,\, that is\,,\,\\[0.1in]
(\,A.6.4\,)
$$
{\tilde{\bf P}}_{\flat_{\,\parallelsum}} \ = \  \sum_{l\ =\,1}^\flat   c_{\ l}\,\cdot\, {{\bf V}}_{\,l}^{4\over {n\,-\,2}}\,\cdot\,\left(\ \Lambda_{\,l}\,\cdot\, {{\partial \, {{\bf V}}_{\,l} }\over {\partial \Lambda_{\,l} }}  \right) \ \,+\, \ \sum_{l\,=\,1}^\flat  \left(\,\sum_{j\ =\,1}^n c_{\ l\,,\,j}\,\cdot\,  {{\bf V}}_{\,l}^{4\over {n\,-\,2}}\,\cdot\, \left[\  \Lambda_{\,l}\,\cdot\, {{\partial \,{{\bf V}}_{\,l} }\over {\partial \,\Xi_{{\,l}_{|_j}} }}\  \right]\ \right)\ .
$$
Note that if all the coefficients $c_{\ l}\,$ and $\,c_{\ l\,,\,j}\,$ vanish (\,for $\ l \, = \, 1\,, \ \cdot \cdot \cdot\,, \ \flat\,,\,$ and $\ j \ = \ 1\,, \ \cdot \cdot \cdot\,, \ n$\,)\,,\,   then for that particular $\,W_{\,\,\flat} \ + \ \phi_{\ \flat}\,$,\, it is a solution of the equation
$$
\Delta\, (\,W_{\,\,\flat} \ + \ \phi_{\ \flat}\,) \ +  \ (\,{\tilde c}_n\!\cdot K\,)\,\cdot\,(\,W_{\,\,\flat} \ + \ \phi_{\ \flat}\,)^{{\,n\,+\,2\,}\over {n\,-\,2}} \ = \ 0\,.\, \leqno (\,A.6.5\,)
$$

\newpage

Introduce the reduced functional\\[0.1in]
(\,A.6.6\,)
\begin{eqnarray*}
 & \ &  {\bf I}_{\,\cal R} \,(\,\lambda_{\,1}\,, \ \cdot \cdot \cdot\,, \ \lambda_{\ \flat}\,; \ \xi_{\,1}\,, \ \cdot \cdot \cdot\,, \ \xi_{\ \flat}\,) \\[0.2in]
  & = & {\bf I}\,(\ W_{\,\,\flat} \ + \ \phi_{\ \flat}\ ) \\[0.2in]
  & = & {1\over 2}\ \int_{\R^n}\ \langle\,\btd\,(\,W_{\,\,\flat} \ + \ \phi_{\ \flat}\,)\,,\ \ \ \btd\,(\,W_{\,\,\flat} \ + \ \phi_{\ \flat}\,)\ \rangle\\[0.2in]
   & \ & \hspace*{2in}-\ \left(\ {{n\,-\,2}\over {2n}}\,\right)\,\cdot\,\int_{\R^n}\,(\,{\tilde c}_n\!\cdot K\,)\,(\,W_{\,\,\flat} \ + \ \phi_{\ \flat}\,)_+^{{2n}\over {\,n\,-\ 2\,}} \ {\bf .}
\end{eqnarray*}
Note that $\,{\bf I}_{\,\cal R}\,$
 depends only on the finite number of  bubbles parameters $\, {\bf P} _{\ \flat} \,$ as shown in (\,A.6.3\,)\,,\, with the understanding that the parameters fulfilling the conditions in Proposition A.5.57\,.\, We compute\\[0.1in]
(\,A.6.7\,)
\begin{eqnarray*}
D_{\,{l_{\,|j}}}\,  {\bf I}_{\,\cal R}
     & =  &    \int_{\R^n} \!\left[\ \langle\ \btd\,(W_{\,\,\flat} \ + \ \phi_{\ \flat})\,,\,\btd\,D_{\,{l_{\,|j}}}\, (\ V_l\,\,+\,\,\phi_{\ \flat}\ )\ \rangle\ \right.\\[0.2in]
     & \ & \hspace*{2in}   \left. -\ (\,{\tilde c}_n\!\cdot K\,)\,(W_{\,\,\flat} \ + \ \phi_{\ \flat})^{{\,n\,+\,2\,}\over{n\,-\,2}}_+\,\cdot\,D_{\,{l_{\,|j}}}\, (\ V_l\,\,+\,\,\phi_{\ \flat}\ )\ \right] \ \ \ \ \ \ \ \ \ \ \ \ \ \ \ \ \ \ \ \ \ \ \\[0.2in]
  & = &  -\, \int_{\R^n} \!\left[\ \Delta\,(W_{\,\,\flat} \ + \ \phi_{\ \flat}) \ + \ (\,{\tilde c}_n\!\cdot K\,)\,(W_{\,\,\flat} \ + \ \phi_{\ \flat})^{{\,n\,+\,2\,}\over{n\,-\,2}}_+ \  \right]\,\cdot\,{\tilde h} \\[0.in]
  & \ & \hspace*{3.5in} \,[ \ {\tilde h} \ = \  D_{\,{l_{\,|j}}}\, (\ V_l\,\,+\,\,\phi_{\ \flat}\ )\ ]\\[0.1in]
  & = &  -\, \int_{\R^n} {\bf P}_{/\!/}   \,\cdot\,{\tilde h} \hspace*{2.5in} [\ {\mbox{via \ \ equation}} \ \ (\,A.6.4\,) \ ]\\[0.2in]
  & = &  {1\over { n\,(\,n\, + \, 1\,)}} \cdot  \Bigg\langle  \left\{ \ \sum_{l\,'\ =\,1}^\flat \ \left[\  c_{\,l}\,\cdot\,\lambda_{\,l\,'}\,\cdot\,{{\partial\, V_{\,l\,'}  \over {\partial \lambda_{\,l\,'}}} } \ + \    \sum_{j\,=\,1}^n c_{j_{\,l\,}}  \,\cdot\,\lambda_{\,l\,'}\,\cdot\,{{\partial \,V_{\,l\,'}   \over {\partial\, \xi_{{\,l\,'|_{\,j} }}}} } \ \right] \  \right\}  \ , \ \,    {\tilde h}\ \Bigg\rangle_{\,\btd} \ .
\end{eqnarray*}
for $\ l \ = \ 1\,, \ \cdot \cdot \cdot\,, \ \flat\,,\,$ and $\ j \ = \ 1\,, \ \cdot \cdot \cdot\,, \ n\,.\,$
Here we apply  the Riesz Representation Theorem\,,\, Proposition A.5.57\,,\,  (\,A.6.2\,)\, and similar observation as in (2.2) in the main text\, \,.

\newpage

{\bf Proposition A.6.8\,.} \ \   {\it For a smaller choice of} $\ {\underline{\lambda}}_{\ \epsilon}\,$ {\it    in Proposition}  A.5.57\,.\, {\it suppose that the bubble parameters}
$$
{\bf P} _{\,(\,\flat\,)} \ = \ (\  \lambda_{\,1}\,, \ \cdot \cdot \cdot\,, \ \lambda_{\,\,\flat}\,; \ \  \xi_{\,1}\,, \ \cdot \cdot \cdot\,, \ \xi_{\,\flat}\,)
$$
{\it  satisfy  the conditions in Proposition}  A.5.57\,.\, {\it  Assume also that}
$$
   \gamma \ > \ {{\sigma}\over {\,n\ - \ 2\,  }} \ .
   \leqno (\,A.6.9\,)
$$
{\it If} $\ {\bf P} _{\,(\,\flat\,)}\,$  {\it is a critical point of the reduced functional $\, {\bf I}_{\,\cal R}$,\,} [\ {\it as defined in}\, (\,A.6.6\,)\ ]\,,\,  {\it that is}
$$
 {{\partial \,{\bf I}_{\,\cal R}}\over {\partial\, \lambda_{\,\,l} }}\,\Bigg\vert_{\, {\bf P} _{\,(\,\flat\,)}}= \ {{\partial \,{\bf I}_{\,\cal R} }\over {\partial \,\xi_{{\,l \,|_{\,j}}} }} \,\Bigg\vert_{\,{\bf P} _{\,(\,\flat\,)}} \ \,= \ 0 \ \ \ \ \ {\it{for}}  \ \ \ l\  = \ 1\,,\, \cdot \cdot \cdot\,, \   \flat\,, \ \ \ \ j \ = \ 1\,, \ 2\,, \ \cdot \cdot \cdot \,, \ n\,, \leqno  (\,A.6.10\,)
$$
{\it then\,}  $\,W_{\,\,\flat} \ + \ \phi_{\,\,\flat}\,$ [\,{\it corresponding to} $\,{\bf P} _{\,(\,\flat\,)}\,$]\, {\it is also a critical point of}\, (\,{\it the full functional}\ )\, $\,{\bf I} \ $.\, {\it That is\,,\,}
$$
{\bf I}\,' \, (\ W_{\,\,\flat} \ + \ \phi_{\ \flat}\ ) \ = \ 0\,.
$$

\vspace*{0.15in}


{\it Proof\,.} \ \
The proof of Proposition A.6.8 is similar to the proof of Theorem 2.8 in   \cite{III} (\,see also  Theorem 2.12 in \cite{Progress-Book}\ )\,.\, More care has to be taken here as the number of bubble (\,denoted by $\,\flat\,$) may tend to $\,\infty\,.\,$ \bk
A  direct calculation as in  \cite{I} shows that
$$
\bigg\Vert\  {{\partial \,V_{\,l} }\over {\partial \lambda_{\,\,l}}}\ \bigg\Vert_\btd \ = \ {{ {\bar C}_{o}}\over { \lambda_{\,\,l}}} \ \ \ \ \ {\mbox{and}} \ \ \ \  \ \bigg\Vert \ {{\partial \,V_{\,l} }\over {\partial \,\xi_{\,l_{\,|\,j}}}}\  \bigg\Vert_\btd \ = \ {{{\bar C}_{\,1}}\over { \lambda_{\,\,l}}}\  \leqno (\,A.6.11\,)
$$
for  $\,  l \,= 1\,, \ \cdot\,\cdot \cdot\,, \ \flat\,; \ \  j\, = 0\,, \ 1\,,\,\cdot\,\cdot \cdot\,, \ n\ .\,$
Here the dimensional positive constants $\,{\bar C}_o\,$ and $\,{\bar C}_1\,$ are independent on $\,\lambda_{\,\,l}\,$ and $\,\xi_{\,l_{\,|\,j}}\,$ [\,and $\,{\bar C}_1\,$ is independent on $\ j \ = \ 1\,, \ 2\,,\,\cdot \cdot \cdot\,, \ n$\,]\,.\,
Let\\[0.1in]
(\,A.6.12\,)
\begin{eqnarray*}
{\tilde Q}_{\,{l_{\,o}}} & = & {{\partial \,V_{\,l} }\over {\partial \lambda_{\,\,l}}}\ \bigg/\,\bigg\Vert \,{{\partial \,V_{\,l} }\over {\partial \lambda_{\,\,l}}}\,\bigg\Vert_\btd \ = \ {1\over { {\bar C}_{\,o} }}\,\cdot\, \lambda_{\,\,l}\,\cdot\, {{\partial \,V_{\,l} }\over {\partial \lambda_{\,\,l}}} \ \ \Longleftrightarrow \ \ \lambda_{\,\,l}\,\cdot\, {{\partial \,V_{\,l} }\over {\partial \lambda_{\,\,l}}}  \ = \   {\bar C}_{\,o} \cdot {\tilde Q}_{\,{l_{\,o}}} \,,\\[0.2in]
{\mbox{and}} \ \ \ \ \   {\tilde Q}_{\,{l_{\,|j}}} & = & {{\partial \,V_{\,l} }\over {\partial \,\xi_{\,{l_{\,|j}}} }}\ \bigg/\,\bigg\Vert \, {{\partial \,V_{\,l} }\over {\partial \,\xi_{\,{l_{\,|j}}} }}\, \bigg\Vert_\btd \ = \  {1\over {{\bar C}_{\,1} }}\,\cdot\,\lambda_{\,\,l} \cdot {{\partial \,V_{\,l} }\over {\partial \,\xi_{\,{l_{\,|j}}} }}  \ \ \Longleftrightarrow \ \  \lambda_{\,\,l} \cdot {{\partial \,V_{\,l} }\over {\partial \,\xi_{\,{l_{\,|j}}} }} \ = \ {\bar C}_{\,1} \cdot   {\tilde Q}_{\,{l_{\,|j}}}  \  \ \ \
\end{eqnarray*}
for  $\,  l \,= 1\,, \ \cdot\,\cdot \cdot\,, \ \flat\,; \ \  j\, = 0\,, \ 1\,,\,\cdot\,\cdot \cdot\,, \ n\ .\,$
Fixed any $\,l\,$ ($\,1 \ \le \ l \ \le \ \flat\,$)\,,\, a direct calculation as in \cite{I}  shows that  \\[0.1in]
(\,A.6.13\,)
\begin{eqnarray*}
 \langle\ {\tilde Q}_{\,{l_{|\,o}}}\,, \ {\tilde Q}_{\,{l_{|\,o}}}\  \rangle_{\,\btd} & = &  1 \ , \ \ \ \ \ \ \ \ \langle\ {\tilde Q}_{\,{l_{\,|j}}}\,, \ {\tilde Q}_{\,{l_{\,|j}}}\  \rangle_{\,\btd} \ = \ 1\ , \\[0.2in]
  \langle\ {\tilde Q}_{\,{l_{\,|j}}}\,, \ {\tilde Q}_{\,{l_{\,|j'} }}\ \rangle_{\,\btd} & = &  0 \   \ \ \ \ \ \  \mfor \ \ j^{\,'}   \ \not=\   j\ ,  \ \ \ \ \ {\mbox{and}} \ \ \  \   \langle\ {\tilde Q}_{\,{l_{\,|\,o}}}\,, \ {\tilde Q}_{\,{l_{\,|j } }}\ \rangle_{\,\btd} \ = \  0 \ .
\end{eqnarray*}
Using the weak interaction and integration by parts as in {\bf \S\,A\,4\,.\,j}\,,\,  we have [\ see also Note A.6.25 below\ ]
$$
 \langle\ {\tilde Q}_{\,{l_{\,|j}}}\,, \ {\tilde Q}_{\,{{l\,'}_{|\,j\,'}}}\,\rangle_{\,\btd} \ = \ O\,\left( \, {\bar\lambda}_{\ \flat}^{\,(\,n\,-\,2\,)\,\cdot \, \gamma}\, \right) \ \ \ \ \mfor \ \   l\ \not=\  l\,'\  \ \ (\,1 \ \le \ j\,, \ \ j\,' \ \le \ n\,)\,. \leqno (\,A.6.14\,)
$$
In particular [\,as discussed in {\bf \S\,A\,2}\,]\,,\, the collection  $$\,\{\ {\tilde Q}_{\,l_{\,j}} \ \,| \ \   l \,= 1\,, \ \cdot\,\cdot \cdot\,, \ \flat\,; \ \  j\, = \,0\,, \ 1\,,\,\cdot\,\cdot \cdot\,, \ n\ \}$$  forms a basis for the orthogonal complement of $\,{\cal D}^{\,1\,,\,2}_{\flat_\perp}\ $ .\, In addition\,,\, for a fixed set of $\,(\,j\,,\ l\,)$\, [ \ $1\ \le  \ j \ \le n\,, \ \ 1 \ \le \ l  \ \le \ \flat $\ ]\,,\,
$$
 \Vert \ D_{\,{l_{\,|_{\,j\,'}}}}\,{\tilde Q}_{\,{l_{\,|j}}}\,\,  \Vert_{\,\btd}  \ \le \  {C \over { {\bar\lambda}_{\,\,\flat}  }} \ \ \ \mfor \ \ j\,' \ = \ 0\,, \ 1\,,\,\cdot\,\cdot \cdot\,, \ n\ .\leqno (\,A.6.15\,)
 $$
 This can be checked via the equation for $\,V_{\,l}\,$ and  the integration by parts formula.
 See {\bf \S\,3\,c} in \cite{I}\,.\bk
Fixed  $\,(\,j\,,\ l\,)\,,\,$ via (\,A.6.7\,)\,,\, we obtain\\[0.1in]
(\,A.6.16\,)
\begin{eqnarray*}
``\,(\,A.6.10\,)\," \!\!\!\!\!\!\!\!\!& \ & \Longrightarrow\ \  D_{\,{l_{\,|j}}}\, {\bf I}\,(\, W_{\,\,\flat}\,+\,\phi_{\ \flat}\,) \ = \ 0 \ \ \mfor  \ \ \ l \,= \,1\,, \ \cdot\,\cdot \cdot\,, \ \flat\,; \ \ \ j\, = \,0\,, \ 1\,,\,\cdot\,\cdot \cdot\,, \ n\,\\[0.2in]
 & \ & \hspace*{-0.8in} \Longrightarrow \ \  \Bigg\langle\, \left\{ \ \sum_{l\,'\ =\,1}^\flat \ \left[\  c_{o_{\,l\,'}}\,\cdot\,\lambda_{\,l\,'}\,\cdot\,{{\partial\, V_{\,l\,'}  \over {\partial \lambda_{\,l\,'}}} } \ + \    \sum_{j\,'\,=\,1}^n c_{j\,'_{\,l\,'\,}}  \,\cdot\,\lambda_{\,l\,'}\,\cdot\,{{\partial \,V_{\,l\,'}   \over {\partial\, \xi_{{\,l\,'|_{\,j\,'} }}}} } \ \right] \  \right\}\ , \\[0.1in]
  & \ & \hspace*{3in}\ \ \ \  [\ \lambda_{\,\,l}\,\cdot D_{\,{l_{\,|\,j}}}\, (\,V_l\,\,+\,\,\phi_{\ \flat}\ )\ ] \ \Bigg\rangle_\btd \! = \ 0\\[0.15in]
\\[0.2in]
 & \ & \hspace*{-0.8in} \Longrightarrow \ \  \Bigg\langle\,   \left(\   \sum_{l\,'\,=\,1}^\flat \, \left[\  [\ C_o\,\cdot\,c_{o_{\,l\,'}}\ ]\cdot {\tilde Q}_{\,{l\,'_{\,o}}}   +     \sum_{j\,'\,=\,1}^n \ [\  C_1\,\cdot\,c_{j\,'_{\ l\,'}}\ ]   \cdot {\tilde Q}_{\,{l\,'_{\,j\,'}}}   \, \right]   \, \right), \   \\[0.1in]
  & \ & \hspace*{3in}\ \ \ \  [\ \lambda_{\,\,l}\,\cdot D_{\,{l_{\,|\,j}}}\, (\,V_l\,\,+\,\,\phi_{\ \flat}\ )\ ] \ \Bigg\rangle_\btd \! = \ 0 \\[0.1in]
 & \ & \hspace*{4.4in} [ \ {\mbox{via}} \ \ (\,A.6.12\,) \ ]\\[0.1in]
 & \ &\hspace*{2.2in}\mfor\ \  l \,= \,1\,, \ \cdot\,\cdot \cdot\,, \ \flat\,; \ \ \ j\, = \,0\,, \ 1\,,\,\cdot\,\cdot \cdot\,, \ n\,.
\end{eqnarray*}
Using (\,A.6.13\,)\,,\, the last expression in (\,A.6.16\,)  can be written as\\[0.1in]
(\,A.6.17\,)
$$ A_{\,{l_{\,|\,j}}}\,  \ +  \, \sum_{(j'\,, \, l') \ \not=\  (j\,,\,l) }\!\!\!A_{\,{l'_{\,|j'} } } \,\langle \ {\tilde Q}_{l\,'_{\,|j\,'}}\,, \ {\tilde Q}_{l_{\,|j}}\ \rangle_{\,\btd}\emph{} \ +\   \sum_{l\,'\,,\ j\,' }    \,A_{\,{l\,'_{\,|j\,'} } } \,\langle \ {\tilde Q}_{l'_{\,|j'}}\,, \ \ (\,\lambda_{\,\,l}\,\cdot\,D_{\,l_{\,|\,j}}\,)\, \phi_{\ \flat}\ \rangle_{\,\btd} \ = \ 0\,,
$$
where
\begin{eqnarray*}
(\,A.6.18\,) \ \ \ \ \ \ \ \ \ \ \ \ \ \ \ \ \ \   \  \ \ \ \ \ \ \ \  \ \ \ \ \ \ \ \ \  \ \ \ \ \ \ \ \ \    A_{\,{l_{\,|\,o}}} &  = &  {\bar C}_o^2\,\cdot\,c_{\,l}\ ,\  \ \ \ \ \ \ \ \  \ \ \ \ \ \ \ \ \  \ \ \ \ \ \ \ \ \    \  \ \ \ \ \ \ \ \  \ \ \ \ \ \ \ \ \  \ \ \ \ \ \ \ \ \    \  \ \ \ \ \ \ \ \  \ \ \ \ \ \ \ \ \  \ \ \ \ \ \ \ \ \    \\[0.2in]
{\mbox{and}} \ \ \ \ \ \ \ \ \ \ \ \ \ \ \ \ \ \ \ \ \ \ \ \   \ \ \ \   \ \ \ \ \ \ \ \ \ \ \ \ \ \ \ \ \ \ \   \ \ \ \ A_{\,{l_{\,|\,j}}} &  = &  {\bar C}_1^2\,\cdot\,c_{\,l_{\,j}}\
  \end{eqnarray*}
 for $\   l \,= 1\,, \ \cdot\,\cdot \cdot\,, \ \flat\,; \ \  j\, = 0\,, \ 1\,,\,\cdot\,\cdot \cdot\,, \ n\ .\,$ Our goal is to show that
\begin{eqnarray*}
& \ & A_{\,{l_{\,|\,j}}} \ = \ 0 \ \ \ \  \Longrightarrow \ \  c_{\,{l_{\,|\,j}}} \ = \ 0 \ \ \ \ \mfor\ \  \,  l \,=\,1\,, \ \cdot\,\cdot \cdot\,, \ \flat\,; \ \ \ j\, =\,\emph{} 0\,, \ 1\,,\,\cdot\,\cdot \cdot\,, \ n\,,\\[0.2in]
& \Longrightarrow &  \ \ \
{\bf P}_{/\!/}  \ = \ 0 \ \ \Longrightarrow \ \ {\bf I}\,'\,(\,W_{\,\,\flat} \ + \ \phi_{\ \flat}\,)  \ = \ 0\,.
\end{eqnarray*}
To this end, we observe that $\,\phi_{\ \flat} \, \in \  {\cal D}^{1\,\ 2}_{\flat_\perp}\,$,\,   where the perpendicular condition leads to,\\[0.1in]
(\,A.6.19\,)
\begin{eqnarray*}
  & \ & \ \ \langle \,  \phi_{\ \flat}\,, \ \ {\tilde Q}_{\,{l_{\,|j}}}\,\rangle_{\,\btd} = 0 \ \ \ \ \ \ \ \ \ \ \ \  \mfor \ \ \ \   l \,= 1\,, \ \cdot\,\cdot \cdot\,, \ \flat\,; \ \ \ j\, = 0\,, \ 1\,,\,\cdot\,\cdot \cdot\,, \ n\,,\ \ \ \ \ \ \ \ \\[0.15in]
& \Longrightarrow & \ \ \langle \,D_{\,{l_{\,|_{\,j}}}} \,\phi_{\ \flat}\,, \ \ {\tilde Q}_{\,{l_{\,|j}}}\,\rangle_{\,\btd}\  +\  \langle \,  \phi_{\ \flat}\,, \ \ D_{\,{l_{\,|_{\,j}}}}\,{\tilde Q}_{\,{l_{\,|j}}}\,\rangle_{\,\btd} \ \ =\ 0\\[0.075in]
& \Longrightarrow & \ \   |\,\langle \,D_{\,{l_{\,|_{\,j}}}}\,  \phi_{\ \flat}  \,, \ \ {\tilde Q}_{\,{l_{\,|j}}}\,\,\rangle_{\,\btd}\ | \  \le \   \Vert \   \phi_{\ \flat} \ \Vert_{\,\btd} \,\,*\,\Vert \ D_{\,{l_{\,|_{\,j}}}}\,{\tilde Q}_{\,{l_{\,|j}}}\,\,  \Vert_{\,\btd} \ \le \  {{C'}\over {  \lambda_{\,\,l}}}\,\cdot\  \Vert\,   \phi_{\ \flat}\ \Vert_{\,\btd}\\[0.1in]
& \ & \hspace*{3.5in} \ \ \ \ \ \ [\ {\mbox{via}} \ (\,A.6.15\,)\ ]\\[0.11in]
& \Longrightarrow & \ \   |\,  \langle \  (\,\lambda_{\,\,l}\,\cdot\, D_{\,{l_{\,|_{\,j}}}}\,)\,  \phi_{\ \flat} \  \, \ \ {\tilde Q}_{\,{l_{\,|j}}}\,, \ \ \rangle_{\,\btd}\,|\ \le \ C \cdot \,{\bar\lambda}_{\ \flat}^{{{\varpi}_{\,(\,A.5.38\,)} }}\\[0.1in]
& \ & \hspace*{2.7in} \ \ \ \ \ \ [\ {\mbox{via}} \ {\bf (i)}_{{\bf{\,(\,5.6\,)}}} \ {\mbox{in \ \ Proposition A.5.57}}\ ]\,.
\end{eqnarray*}
On account of this,  equation (\,A.6.13\,) can be written as
$$
A_{\,{l_{\,|_{\,j}}}} \ + \ \sum_{l\,,\ j }      \left(\ A_{\,{l'_{\,|\,{j'}}}} \,\cdot\,\clubsuit^{ {{l}_{\,|j}} }_{ \,{{l'}_{\,|j'}} }\right) \ = \ 0 \ \ \ \ \mfor\ \  \,  l \,=\,1\,, \ \cdot\,\cdot \cdot\,, \ \flat\,; \ \ \ j\, =\,\emph{} 0\,, \ 1\,,\,\cdot\,\cdot \cdot\,, \ n\,, \leqno (\,A.6.20\,)
$$
where the numbers
 $$
 \bigg\vert \,\clubsuit^{ {{l}_{\,|\,j}} }_{ \,{{l'}_{ |\,j'}} }\,\bigg\vert \ = \ O\ \left(\,{\bar\lambda}_{\ \flat}^{\,{{\varpi}_{\,(\,A.5.38\,)} }} \,\right) \ + \   O\ \left(\,{\bar\lambda}_{\ \flat}^{\,(\ n\,-\,2\ )\,\cdot\, \gamma} \,\right)\ . \leqno (\,A.6.21\,)
  $$
  Here we apply (\,A.6.13\,)\,,\, (\,A.6.14\,)\,  and (\,A.6.19\,)\,.\, \bk
Equations (\,A.6.20\,) can be expressed as a  system linear equations with unknowns $\,A_{\,{l_{\,|\,j}}} \,$,\, whose  matrix $\,{\bf M}\,$ (\,of dimension $\,[\ \flat\cdot (\,n\,+\,1)\ ]\,\times [\ \flat\cdot (\,n\,+\,1)\ ] \ ;= \ m \,\times\, m$\ )\, has entries $$\ \delta^{ {{l}_{\,|\,j}} }_{ {{l\,'}_{\,|\,j\,'}} }  \ + \  \clubsuit^{ {{l}_{\,|\,j}} }_{ \,{{l\,'}_{\,|\,j\,'}} } \,,$$ where $\, \delta^{ {{l}_{\,|\,j}} }_{ {{l\,'}_{\,|\,j\,'}} }\,$ forms the identity matrix $\,{\bf I}\,.\,$ We estimate the determinant via the Leibniz formula\,:
$$\displaystyle{ {\mbox{Det}} \,\,{\bf M} \ = \ \sum_{\tau \,\in\,S_m}\! \left(\ {\mbox{sgn}}\,(\,\tau\,) \prod_{I\,=\,1}^m A_{\,I\,,\ \tau\,(\,I\,)} \right)} \ .$$
In the above, $\,S_m\,$ is the permutation group of a set with $\,m\,$ elements\,,\, and $\,A_{\,I\,, \, J}$\, is the $\,(\,I\,,\ J\,)$\,-\,th entry of the matrix \,{\bf M}\,.
This yields the estimate\\[0.1in]
(\,A.6.22\,)
$$
|\,{\mbox{Det}} \,{\bf M}\,| \ \ge \ 1 \ - \ [\ \flat\cdot (\,n\,+\,1)\ ]\,!\,\cdot\,O\ \left(\ \left[\ O\ \left(\,{\bar\lambda}_{\ \flat}^{\,{{\varpi}_{\,(\,A.5.38\,)} }} \,\right) \ + \   O\ \left(\,{\bar\lambda}_{\ \flat}^{\,(\ n\,-\,2\ )\,\cdot\, \gamma} \,\right)\ \right]^{\ \flat\cdot (\,n\,+\,1)\,-\,1 }\  \right)\ {\bf .}
$$
We know that
$$
( \ 2\ \le \ ) \ \flat \ \le \ {1\over { {\bar{\lambda}}_{\,\flat}^{\,\sigma} }}\ .
$$
For any integer  $\,k \ \ge \ 2\,$,\, via Mathematical Induction, we obtain
$$
k\,! \ \le  \ k^{\ k\,-\,1}
$$
It follow that
\begin{eqnarray*}
(\,A.6.23\,) \ \ \ \ \ \ \ \ \ \ \ \ \ \ \ \ \ \ \ \  [\ \flat\cdot (\,n\,+\,1)\ ]\,! & \le &  [\ \flat\cdot (\,n\,+\,1)\ ]^{\,(\,n\,+\,1\,)\, \cdot\, \flat \ - \ 1 }
 \\[0.2in]
  & \le  & C_4 \cdot \left[\ {1\over { {\bar\lambda}_{\ \flat}^{ \,\sigma \ +  \  o_{\,+}\,(\,1\,)} }}\ \right]^{\ \flat\,\cdot\, (\,n\,+\,1\,)\ - 1 } \ .\ \ \ \ \ \ \ \ \ \ \ \ \ \ \ \ \ \ \ \ \ \ \ \ \ \ \ \ \ \
\end{eqnarray*}

Here
$$
{1\over { \lambda_{\ \flat}^{\,o_{\,+}\,(\,1\,)} }}\  \ge \ n\,+\,1 \ \ \ \ \ \ \mfor \ \  \lambda_{\ \flat}\ \ \,{\mbox{small \ \ enough}}\,.
$$
From condition (\,A.5.56\,) in Proposition A.5.55 (\,and hence automatically a  condition in Proposition A.5.57\,)\,,\,
$$
 {\vartheta}_{\,(\,A.5.7\,)}  \ > \  \left(\  {{ n \ + \ 10 }\over 8}\ \right)\, \cdot \sigma  \ \ \ \ \ \ {\mbox{and}} \ \ \ \  \  \ n \ \ge \ 6 \ \ \Longrightarrow \ \  {\vartheta}_{\,(\,A.5.7\,)}  \ > \   2\,\cdot\, \sigma   \ \ \Longrightarrow \ \  {{\varpi}_{\,(\,A.5.38\,)} }\ > \     \sigma  \ .
$$
Likewise, condition (\,A.6.9\,) shows that
$$
(\,n \ - \ 2\,) \cdot \gamma \ > \ \sigma\,.
$$
Thus if $\ {\bar{\lambda}}_{\,\flat}$\, is small enough\,,\,
then we get
$$
|\ {\mbox{Det}} \,\,{\bf M}\,| \ \ge \ 1 \ - \ o_{\,+}\,(\,1\,)\ ,
$$
where $\,o_{\,+}\,(\,1\,) \ \to \ 0\,$ as $\ {\bar{\lambda}}_{\ \flat} \ \to \ 0^{\,+}\,.$\,
$\,{\bf M}\,$ is then shown to be invertible (\,when $\, {\bar\lambda}_{\ \flat}\,$ is made small enough\,)\,.\, \bk Hence  zero solutions (\,or the trivial solutions\,) are the only option (\,that is, $\,A_{\,{l_{\,|{j}}}}  \ = \ 0\,$ for $\, l \,= 1\,, \ \cdot\,\cdot \cdot\,, \ \flat\,$ and $\,j\  =\  0\,,\, 1\,, \ \cdot\,\cdot \cdot\,, \   n$\,)\,.\, This, in turn [\,via (\,A.6.18\,)\,]\,,\, implies that $\,c_{j_{\ l}} \ =\ 0\,$ for $\, l \,= 1\,, \ \cdot\,\cdot \cdot\,, \ \flat\,$ and $\,j\, = 0\,, \ 1\,,\,\cdot\,\cdot \cdot\,, \ n\,.$\, It follows that
$$\ \ \ \ \ \ \ \ \ \ \ \ \ \  \ \ \ \ \ \ \ \ \ \  \ \ \ \ \ \ \ \ \  \ \ \ \ \  \  \ \ \ \ \  \
{\bf I}\,'\,(\ W_{\,\,\flat} \ + \ \phi_{\ \flat}\ )  \ \equiv \ 0\,.\ \ \ \ \ \ \ \ \ \ \ \ \ \  \ \ \ \ \ \ \ \ \ \  \  \ \ \ \ \  \  \ \ \ \ \  \  \ \ \ \ \  \square
$$

\vspace*{0.3in}

{\it Note\,} A.6.24\,. \ \  [\,Refer to (\,A.6.24\,)\,.\, See also {\bf \S\,A\,4\,.\,j}\ .\,]\, For $\,1 \ \le \ l \ \not=\ l\,' \ \le \ \flat\,,\,$ we compute\,:
\begin{eqnarray*}
   \langle\ {\tilde Q}_{\,{l_{\,|j}}}\,, \ {\tilde Q}_{\,{{l\,'}_{|\,j\,'}}}\,\rangle_{\,\btd}
& = & \bigg\vert \,\int_{\R^n} \bigg\langle \btd\!\left(  {\tilde Q}_{\,{l_{\,|j}}}\,\, \right)\,, \ \ \btd\!\left(\,{\tilde Q}_{\,{l\,'_{\,|j\,'}}}\,\,\right) \bigg\rangle\,\bigg\vert \\[0.2in]
& \le &   C_1 \cdot  \lambda_{\ \!l}\,\cdot\,\lambda_{\,l\,'}\,   \bigg\vert \, \int_{\R^n} \bigg\langle \btd\! \left(\ D_{\,{l_{\,|_{\,j}}}}\, V_{\,l} \right)\  , \ \ \btd\! \left(\ D_{\,{l\,'_{\,|_{\,j\,'}}}}\,V_{\,l\,'} \right)\,\bigg\rangle\,\,\bigg\vert \\[0.2in]
& \le &   C_1 \cdot  \lambda_{\ \!l}\,\cdot\,\lambda_{\,l\,'}\,   \bigg\vert \, \int_{\R^n} \left[\  \Delta \left(\ D_{\,{l_{\,|_{\,j}}}}\, V_{\,l} \right)\  \right]\,{\bf \cdot}\,\left(\ D_{\,{l\,'_{\,|_{\,j\,'}}}}\,V_{\,l\,'} \right)\,\bigg\vert\\[0.1in]
& \ & \hspace*{3in} [ \ {\mbox{via}} \ \ (\,A.6.12\,) \ ]\\[0.1in]
& = &  C_1 \cdot  \lambda_{\ \!l}\,\cdot\,\lambda_{\,l\,'}\,\cdot\,   \bigg\vert \, \int_{\R^n} \left[  \ D_{\,{l_{\,|_{\,j}}}}\,  \,[\,\Delta\, V_{\,l}\,]  \ \right]\,\cdot\,\,\left(\ D_{\,{l\,'_{\,|_{\,j\,'}}}}\,V_{\,l\,'} \right)\,\,\bigg\vert\\[0.2in]
& \le &  C_3\cdot\,   \lambda_{\ \!l}\,\cdot\,\lambda_{\,l\,'}\,\cdot\,  {1\over {\,\lambda_{\ \!l}\,}}\,\cdot\, {1\over { \,\lambda_{\,l\,'}\,}} \,\cdot\,\int_{\R^n} V_{\,l}^{{\,n\,+\,2\,}\over {n\,-\,2}}\,\cdot\,V_{\,l\,'} \\[0.1in]
& \ & \hspace*{3in} [ \ \uparrow \ \ {\mbox{via}} \ \ (\,A.1.6\,) \ ]\\[0.1in]
& \le & C\,\cdot\,{1\over {\  {\bf d}_{\ l\,, \ l\,'}^{\ n\,-\ 2}\ }} \hspace*{2in}\ \ \ \ \   [ \ {\mbox{see}} \ \ (\,A.3.29\,) \ ] \\[0.2in]
& = & O\,\left(\ {\bar\lambda}_{\ \flat}^{\ (\,n\,-\,2)\, \cdot\, \gamma} \ \right) \ \ \ \ \ \ \ \left(\ \mbox{recall \ \ that} \ \ {\bf d}_{\,j\,,\,j\,'} \ := \ {{\ \Vert\, \,\xi_{\,j}\ - \ \xi_{\,j\,'}\,\Vert \  }\over
{\sqrt{\,\lambda_{\,j} \cdot \,\lambda_{\,j\,'}\,}}}  \ \right)\ \,{\bf .}
\end{eqnarray*}

\newpage

 {\large{\bf \S\,A\,7\,.\,} \   {\bf Isolating the key terms in the reduced functional\,. }}\\[0.2in]
For direct reference, let us recall that
\begin{eqnarray*}
(\,A.7.1\,)\ \ \ \ \ \ \ \ \ \ \ \ \ \ \  {\bf I}\,(\,f\,) & = &  {1\over 2}\,\int_{\R^n}\,\langle\,\btd\,f\,,\,\btd\,f\,\rangle\ -\ \left(\,{{n\,-\,2}\over {2n}}\,\right)\,\cdot\,\int_{\R^n}\,(\,{\tilde c}_n\!\cdot K\,)\,f_+^{{2n}\over {\,n\,-\ 2\,}} \,,\\[0.2in]
(\,A.7.2\,)\ \ \ \  \ \ \ \ \ {\bf I}\,'\,(\,f\,)\,[\,\varphi\,] & = &
\int_{\R^n} \!\left[\ \langle\,\btd\,f\,,\,\btd\,\varphi\,\rangle\ -\ (\,{\tilde c}_n\!\cdot K\,)\,f^{{\,n\,+\,2\,}\over{\,n\,-\ 2\,}}_+\,\cdot\,\varphi\ \right] \,, \\[0.2in](\,A.7.3\,)\ \ \
(\,{\bf I}\,'' (\,f\,)\, [\,\psi] \  \varphi\,)&  = & \int_{\R^n} \!\left[\,  \langle \,\btd \,\psi\,,\, \btd \, \varphi \, \rangle - \left( \ {{n\,+\,2}\over {\,n\,-\ 2\,}} \ \right)\cdot  (\,{\tilde c}_n\!\cdot K\,)\,f^{{4}\over {\,n\,-\ 2\,}}_+\cdot \psi\cdot \varphi\, \,\right] \,,\\[0.2in]
 & \ &\hspace*{2in}\ \ \    {\mbox{where}} \ \  f\,, \ \varphi \ \ \&\  \ \psi \  \in \  {\cal D}^{\,1,\,\,2}\,.\ \ \ \
\end{eqnarray*}
As in the above, via Proposition A.5.57\,,\,
$\,\phi_{\ \flat}   \in \, W^\flat_{*} \,\cap \, {\cal D}^{1\,,\, 2}_\perp\,$  is the (\,unique\,) \,``\,small\,"\, solution  of equation (\,2.1\,)\,,\, associated with $\,W_{\ \flat}$\,,\, in which the bubble parameters satisfies the conditions in Proposition A.5.57\,.\,  That is\,,\,
\begin{eqnarray*}
 (\,A.7.4\,) \ \ \ \  \ \ \ \ \ \  \ \  \ \ \ \  \ \   & \ & \Delta\, (\ W_{\,\,\flat} \ + \ \phi_{\ \flat}\,) \ +  \ (\,{\tilde c}_n\!\cdot K\,)\,\cdot\,(\,W_{\,\,\flat} \ + \ \phi_{\ \flat}\,)^{{\,n\,+\,2\,}\over {n\,-\,2}} \ = \ {\bf P}_{/\!/} \ \ \ \  \ \ \ \ \ \  \ \ \ \ \ \  \ \  \\[0.2in]
   \Longrightarrow & \ &    {\bf I}\,'\,(\,W_{\,\,\flat} \ + \ \phi_{\ \flat})  \ [\ \phi_{\ \flat} \ ]   \ = \ 0 \  \ \ \ \ \ \ {\mbox{as}} \ \ \ \phi_{\ \flat}  \,\in\, {\cal D}^{1\,\, 2}_\perp\,.
\end{eqnarray*}
[ \ ${\bf P}_{/\!/} \,$ is given by (\,A.5.2\,)\,.\, Here we apply the perpendicular condition and (\,A.1.18\,)\,.\,]

\vspace*{0.3in}

{\bf \S\,A\,7.\,a\,.} \ \  {\bf  Heuristic} ``\,{\bf picture}\,". \\[0.1in]
\emph{}Via a form of Taylor expansion\,, we see that\\[0.1in]
(\,A.7.5\,)
\begin{eqnarray*} & \ &
{\bf I}_{\,\cal R}\, (\,\lambda_{\,1}\,, \ \cdot\,\cdot \cdot\,, \ \lambda_{\ \flat}\,; \ \xi_{\,1}\,, \ \cdot\,\cdot \cdot\,, \ \xi_{\ \flat}\,)\\[0.2in] &  = & {\bf I}\,(\,W_{\,\,\flat} \ + \ \phi_{\ \flat}\,) \ \ \  \ \ \ \ \ \ \ \ \ \ \ \ \ \ \  \ \ \ \ \ \ \ \ \ \left[\ \,W_{\ \flat}  \ = \   \sum_{l \ = \ 1}^\flat \left( \ {{\lambda_{\,\,l} }\over { \ \lambda_{\,\,l}^2 \ + \ \Vert\, y \ - \ \xi_{\,l} \,\Vert^2\  }}\ \right)^{{\,n\,-\,2\,}\over 2} \ \right] \\[0.15in] & = & {\bf I}\,(\,W_{\,\,\flat}) \ + \  {\bf I}\,'\,(\,W_{\,\,\flat}) \, [\ \phi_{\ \flat} \ ]   \ + \ {1\over 2}\,\cdot\, \left(\  {\bf I}\,''\,(\,W_{\,\,\flat} \ + \ s\cdot \phi_{\ \flat}) \,   [\  \phi_{\ \flat}\, ] \ \phi_{\ \flat}\ \right) \ \\[0.1in]
  & \ & \hspace*{4in} \ \ \  \ \ \ \ \  \ \ \ \ \  \ \ \ \ \ [\ s \ \in \ (\,0\,,\ 1)\,]
  \end{eqnarray*}

  \newpage

  \begin{eqnarray*}
  & = & {\bf I}\,(\,W_{\,\,\flat}) \ + \
  {\bf I}\,'\,(\ W_{\,\,\flat} \ + \ \phi_{\ \flat}) \, [\ \phi_{\ \flat} \ ]   \ + \  {1\over 2}\,\cdot\, \left(\  {\bf I}\,''\,(\,W_{\,\,\flat} \ + \ t\cdot \phi_{\ \flat}) \,   [\  \phi_{\ \flat}\, ] \ \phi_{\ \flat}\ \right)\\[0.1in]
  & \ &  \ \ \ \ \ \  \ \ \ \ \ \ \ \ \{\ \ \uparrow \ \ {\bf I}\,'\,(\,W_{\,\,\flat} \ + \ \phi_{\ \flat})  \ [\ \phi_{\ \flat} \ ]   \ = \ 0 \ \}   \ \ \ \  \ \ \ \ \  \ \  [\ t \ \in \ (\,0\,,\ 1)\,] \\[0.1in]
& = & {\bf I}\,(\,W_{\,\,\flat}) \ + \ {1\over 2} \cdot\,  \left(\, {\bf I}\,''\,(\,W_{\,\,\flat} \ + \ t\cdot \phi_{\ \flat}) \,   [\  \phi_{\ \flat}\, ] \ \phi_{\ \flat}\ \right)  \ .
\end{eqnarray*}

\vspace*{0.3in}

{\it Note.} A.7.6. \ \ Refer to (\,A.7.3\,)\,.\, Using Proposition A.5.57\,,\, we obtain
\begin{eqnarray*}
\int_{\R^n} \langle\,\btd\,\phi_{\ \flat}\,,\,\btd\,\phi_{\ \flat}\,\rangle & = &  \Vert \, \phi_{\ \flat} \,\Vert_{\,\btd}^2 \ \le \ C_1 \cdot {\bar\lambda}_{\ \flat}^{2\,{{\varpi}_{\,(\,A.5.38\,)} }}\,,\\[0.2in]
\int_{\R^n}W_{\,\,\flat}^{{4}\over {\,n\,-\,2\,}} \phi_{\ \flat}^{\,2}  & \le  & \left( \ \int_{\R^n}W_{\,\,\flat}^{{2\,n}\over {n\,-\,2}}   \ \right)^{\!\!{2\over n}   } \cdot \left( \ \int_{\R^n} |\,\phi_{\ \flat}\,|^{{2\,n}\over {n\,-\,2}} \ \right)^{\!\!{{n\,-\,2}\over n}   }\\[0.2in]
& \le & C_2 \cdot \flat^{2\over n} \cdot  \Vert \, \phi_{\ \flat} \,\Vert_{\,\btd}^2 \\[0.1in]
& \ & \ \ \ \ \ \ \ \ \  \ \ \ \ [\ {\mbox{via \ \ Sobolev \ \ inequality\,,\, \ \ and \ \ (\,A.4.59\,)}}\ ]\\[0.2in]
& \le & C_3\cdot {\bar\lambda}_{\ \flat}^{2\,{{\varpi}_{\,(\,A.5.38\,)} } \ - \ {{2}\over n} \cdot \,\sigma }\ .
\end{eqnarray*}

Our goal is to show that
\begin{eqnarray*}
 {\bf I}\,(\,W_{\,\,\flat} \ + \ \phi_{\ \flat}\,) &  = &  {\bf I}\,(\,W_{\,\,\flat}) \ + \  ``\,(\,{\mbox{smaller \ \ order \ \ term}})\,"\,,\\[0.1in]
& \ &  \hspace*{-1.3in} \left[ \ {\bf I}\,(\,W_{\,\,\flat})   \ = \   {1\over 2}\,\int_{\R^n}\,\langle\,\btd\,W_{\,\,\flat}\,,\,\btd\,W_{\,\,\flat}\,\rangle\ -\ \left(\,{{n\,-\,2}\over {2n}}\,\right)\,\cdot\,\int_{\R^n}\,(\,{\tilde c}_n\!\cdot K\,)\,W_{\,\,\flat}^{{2n}\over {\,n\,-\ 2\,}} \ \right]\\[0.2in]
{{\partial  }\over {\partial\, \lambda_{\,1}}}  \  {\bf I}\,(\,W_{\,\,\flat} \ + \ \phi_{\ \flat}\,) &  = &  {{\partial \,{\bf I}\, (\,W_{\,\,\flat}) }\over {\partial\, \lambda_{\,1}}}  \ + \  ``\,(\,{\mbox{smaller \ \ order \ \ term}})\," , \\[0.2in]
{{\partial  }\over {\partial\, \xi_{1_{\,|_1}} }}  \  {\bf I}\,(\,W_{\,\,\flat} \ + \ \phi_{\ \flat}\,) &  = &  {{\partial \,{\bf I} \, (\,W_{\,\,\flat}) }\over {\partial\, \xi_{1_{\,|_1}} }}  \ + \  ``\,(\,{\mbox{smaller \ \ order \ \ term}})\,"\,,... \ \ \ \  \ {\mbox{etc.}}  \\
\end{eqnarray*}

In what follows we would give a precise description on the smaller order terms. The arguments runs in  similar ways as those in  {\bf \S\,A\,4}\, and \,{\bf \S\,A\,5}\,.\,\bk {\bf\  In the remaining discussion (\ \S\,A\,7\,--\ \S\,A\,12\ )\,,\, we always assume that  the conditions in Proposition A.6.8 hold\,.}


\newpage

{\bf \S\,A\,7.\,b\,.} \ \  {\bf  Estimating the } ``\,{\bf small order terms}\,"\,. \\[0.1in]
We begin with
\begin{eqnarray*} \ (\,A.7.7\,) \ \ \ \ \ \ \ \ \
{\bf I}\,(\,W_{\,\,\flat} \ + \ \phi_{\ \flat}) & = &  {1\over 2}\,\int_{\R^n}\,\langle\,\btd\,(\,W_{\,\,\flat} \ + \ \phi_{\ \flat})\,,\,\btd\,(\,W_{\,\,\flat} \ + \ \phi_{\ \flat}\ )\,\rangle\\[0.2in]
& \ &  \ \ \ \ \ \  \ \ \ \ \ \ \ \ \ -\ \left(\,{{n\,-\,2}\over {2n}}\,\right)\,\cdot\,\int_{\R^n}\,(\,{\tilde c}_n\!\cdot K\,)\,(\,W_{\,\,\flat} \ + \ \phi_{\ \flat})_+^{{2n}\over {\,n\,-\ 2\,}}\ \,. \ \ \ \ \ \ \ \ \
\end{eqnarray*}
Let us observe that\,,\, via integration by parts\,,
$$
{{\ \partial \,{\bf I}\, (\,W_{\,\,\flat})\  }\over {\partial\, \lambda_{\,1}}}    \ = \  -\, {1\over 2}\,\int_{\R^n}\,\Delta\,(\,W_{\,\,\flat}\,)\,\cdot \,{{\partial\,V_1}\over { \partial\, \lambda_{\,1} }} \ -\ \int_{\R^n}\,(\,{\tilde c}_n\!\cdot K\,)\,W_{\,\,\flat}^{{\,n\,+,2\,}\over {n\,-\,2}}\,\cdot \,{{\partial\,V_1}\over { \partial\, \lambda_{\,1} }} \ , \leqno (\,A.7.8\,)
$$
which is the main object we seek to extract from (\,A.7.7\,)\,.\,
It follows that\\[0.1in]
(\,A.7.9\,)
\begin{eqnarray*}
 & \ &  \lambda_{\,1}\cdot{{\partial \, {\bf I}_{\,\cal R} }\over {\partial\, \lambda_{\,1}}} \ = \ \lambda_{\,1}\cdot {{\partial \,{\bf I}}\over {\partial\, \lambda_{\,1}}} \, (\,W_{\,\,\flat} \ + \ \phi_{\ \flat})\\[0.2in]
    &= & \int_{\R^n}\,\bigg\langle\,\btd\,( W_{\,\,\flat} \ + \ \phi_{\ \flat})\,,\, \ \  \ \btd\,\left( \lambda_{\,1}\cdot{ { \partial \, V_1 } \over {\partial\, \lambda_{\,1}}}   \ + \ \lambda_{\,1}\cdot{ { \partial\, \phi_{\ \flat} } \over {\partial\, \lambda_{\,1}}}\  \right)\,\bigg\rangle\\[0.2in]
& \ & \ \ \ \ \ \ \ \ \ \  \ \ \ \ \ \ \ \ \ \ \ \ \ \  \ \ \ \ \ \ \ \ \ \ \ \  -\  \int_{\R^n}\,(\,{\tilde c}_n\!\cdot K\,)\,(\,W_{\,\,\flat} \ + \ \phi_{\ \flat})_+^{{\,n\,+\,2\,}\over {n\,-\,2}}\,\cdot\,\left(\ \lambda_{\,1}\cdot { { \partial\, V_1 } \over {\partial\, \lambda_{\,1}}}   \ + \ \lambda_{\,1}\cdot { { \partial \, \phi_{\ \flat} } \over {\partial\, \lambda_{\,1}}}\  \right) \\[0.2in]
& = &  {\bf I}\,' \,(\,W_{\,\,\flat} \ + \ \phi_{\ \flat}\,)  \left[\  \hslash \ \right]\\[0.1in]
& \ &  \ \ \ \  \ \ \ \ \  \  \ \ \ \  \  \ \ \ \  \   \left(\  \ {\mbox{in \ \ the \ \ style \ \ of  \ \ (\,A.7.2\,)\,,\ \ where}} \ \   \hslash \ = \ \lambda_{\,1}\cdot { { \partial \,V_1 } \over {\partial\, \lambda_{\,1}}}   \ + \ \lambda_{\,1}\cdot { { \partial \,\phi_{\ \flat} } \over {\partial\, \lambda_{\,1}}}  \right)\\[0.1in]
& = & \bigg\langle \ {\bf I}\,' \,(\,W_{\,\,\flat} \ + \ \phi_{\ \flat})\,, \ \ \ \hslash\ \bigg\rangle_\btd \\[0.2in]
& \ &  \ \ \ \ \  \ \ \ \  (\ \uparrow \ \ {\mbox{via \ \ Riesz's \ \ Representation \ \ Theorem}}\, )\\[0.2in]
& = & \bigg\langle \ {\bf I}\,' \,(\,W_{\,\,\flat} \ + \ \phi_{\ \flat})\,, \ \ \left(\ \lambda_{\,1}\cdot { { \partial \,V_1 } \over {\partial\, \lambda_{\,1}}}  \right) \, \bigg\rangle_\btd   \ + \ \ \bigg\langle \ {\bf I}\,' \,(\,W_{\,\,\flat} \ + \ \phi_{\ \flat})\,, \ \left(\ \lambda_{\,1}\cdot { { \partial \phi_{\ \flat} } \over {\partial\, \lambda_{\,1}}}  \right) \, \bigg\rangle_\btd\\[0.2in]
& = &  {\bf I}\,' \,(\,W_{\,\,\flat} \ + \ \phi_{\ \flat}\,)  \left[\   \lambda_{\,1}\cdot { { \partial \,V_1 } \over {\partial\, \lambda_{\,1}}}   \ \right]  \ + \ \bigg\langle\, {\bf I}\,'\, (\,W_{\,\,\flat} \ + \ \phi_{\ \flat})\,, \  \ \left(\ \lambda_{\,1}\cdot { { \partial \phi_{\ \flat} } \over {\partial\, \lambda_{\,1}}}  \right)\, \bigg\rangle_\btd
\end{eqnarray*}

\newpage

\begin{eqnarray*}
& = & \int_{\R^n}\,\bigg\langle\,\btd\,(\, W_{\,\,\flat} \ + \ \phi_{\ \flat})\,,\, \ \ \ \btd\,\left( \lambda_{\,1}\cdot { { \partial \,V_1  } \over {\partial\, \lambda_{\,1}}}   \right)\,\bigg\rangle \ \ \ \ \cdot\,\cdot\,\cdot\,\cdot\,\cdot\,\cdot\,\cdot\, \cdot\,\cdot\,\cdot\,\cdot\,\cdot\,\cdot\,\cdot\,\cdot\,\cdot\, \cdot\,\cdot\,\cdot\,\cdot\,\ \ ({\bf I})_{\,(\,A.7.8\,)}\\[0.1in]
 & \ & \hspace*{4in} [\ \updownarrow \ \ {\mbox{via  \ \ (\,A.7.2\,)}} \ ]\\[0.1in]
 & \ & \ \ \ \ \ \  \ \ - \  \int_{\R^n}\,(\,{\tilde c}_n\!\cdot K\,)\,(\,W_{\,\,\flat} \ + \ \phi_{\ \flat})_+^{{\,n\,+\,2\,}\over {n\,-\,2}}\,\cdot\,\left(\ \lambda_{\,1}\cdot { { \partial \,V_1  } \over {\partial\, \lambda_{\,1}}}\ \right)   \ \ \ \ \cdot\,\cdot\,\cdot\,\cdot\,\cdot\,\cdot\,\cdot\,\cdot\,\cdot\,\cdot\,\cdot \cdot\,\cdot\,\ \ ({\bf II})_{\,(\,A.7.8\,)}\\[0.2in]
 & \ & \ \ \ \ \ \  \ \   \    \ \ \ \ \ \ \ \ \ \ \  \ + \ \bigg\langle\, {\bf I}\,'\, (\,W_{\,\,\flat} \ + \ \phi_{\ \flat})\,, \  \ \left(\ \lambda_{\,1}\cdot { { \partial \phi_{\ \flat} } \over {\partial\, \lambda_{\,1}}}  \right)\, \bigg\rangle_\btd  \ \ \ \ \cdot\,\cdot\,\cdot\,\cdot\,\cdot\,\cdot\,\cdot\,\cdot\,\cdot\,\cdot\,\cdot\,\cdot \cdot\, \ \ ({\bf III})_{\,(\,A.7.8\,)} .
\end{eqnarray*}


\vspace*{0.35in}

{\bf \S\,A\,7.\,c\,.} \ \
{\bf Estimate of ($\,{\bf III})_{\,(\,A.7.8\,)}\,$}\,.  \\[0.2in]
Using (\,A.5.2\,) and (\,A.7.2\,)\,,\, we have
\begin{eqnarray*}
& \ & \bigg\langle \ {\bf I}\,'\  (\ W_{\,\,\flat} \ + \ \phi_{\ \flat}\ )\,, \ \ \left(\ \lambda\,\cdot\, { { \partial \,\phi_{\ \flat} } \over {\partial\, \lambda_{\,1}}}  \right) \ \bigg\rangle_{\!\btd} \\[0.2in]
&  = &  {1\over {n\,(\,n\,+\,2\,) }} \cdot\,\Bigg\langle \   \sum_{l\,=\,1}^\flat \, \left[\  c_{{\,l}}\,\cdot\,\lambda_{\ \!l}\,\cdot\,{{\partial \,V_{\,l} }\over {\partial \lambda_{\,\,l} }}  \ + \    \sum_{j\ =\,1}^n c_{\,l_{\,j}}  \,\cdot\,\lambda_{\,\,l}\,\cdot\, {{\partial \, V_{\,l} }\over {\partial\, \xi_{\,{l\,|_{\,j} }}}} \ \right]\,, \ \  \left(\ \lambda_{\,1}\,\cdot\, { { \partial \,\phi_{\ \flat} } \over {\partial\, \lambda_{\,1}}}  \right)  \ \Bigg\rangle_{\!\btd}.
\end{eqnarray*}
From the $\,\perp$\,-\,condition $\, \phi_{\ \flat}  \,\in\, {\cal D}^{1\,\, 2}_\perp\,$,\, [\,cf. (\,A.1.16\,) and (\,1.14\,)\,]\,,\,   we have \\[0.1in]
(\,A.7.10\,)
\begin{eqnarray*}
& \ & \bigg\langle \  {{\partial \, V_{\,l} }\over {\partial\, \xi_{\,{l\,|_{\,j} }}}}   \ , \ \ \phi_{\ \flat}\, \bigg\rangle_\btd \ = \ 0 \\[0.2in]
& \Longrightarrow & \bigg\langle  \ {{\partial^2\,V_{\,l}}\over {\partial \lambda_{\,1} \,\partial\, \xi_{\,{l\,|_{\,j} }}}} \ , \ \ \phi_{\ \flat} \bigg\rangle_\btd \ + \ \bigg\langle\   {{\partial \, V_{\,l} }\over {\partial\, \xi_{\,{l\,|_{\,j} }}}}  \ , \ \ {{\partial\, \phi_{\ \flat}}\over {\partial \,\lambda_{\,1}}} \,\bigg\rangle_\btd \ = \ 0 \\[0.2in]
& \Longrightarrow &  \bigg\vert \  \bigg\langle\, \left(\ \lambda_{\,1}\,\cdot\, { { \partial \,V_{\,l} } \over {\partial\, \xi_{1_{\,|\,1} }}}  \right)    \,, \ \  \left(\ \lambda_{\,1}\,\cdot\, { { \partial \,\phi_{\ \flat} } \over {\partial\, \lambda_{\,1}}}  \right)  \  \bigg\rangle_\btd \ \bigg\vert \ = \ \bigg\vert \   \bigg\langle  \ \lambda_{\,1}^2 \cdot {{\partial^2\,V_{\,l}}\over {\partial \lambda_{\,1} \,\partial\, \xi_{\,{l\,|_{\,j} }}}} \ , \ \ \phi_{\ \flat} \bigg\rangle_\btd  \ \bigg\vert \  \\[0.2in]
& \  & \ \ \ \  \le \ C \cdot  \lambda_{\,1}^2 \,\cdot\,\bigg\vert \ \int_{\R^n} \, \left( \ \Delta\, \left[\ {{\partial^{\,2} \,V_{\,l} }\over {\partial\, \lambda_{\,1} \,\partial\, \xi_{1_{\,|\,1} }}} \ \right] \ \right)\,\cdot\,\phi_{\ \flat}  \ \bigg\vert \\[0.1in]
& \ & \ \ \ \  \ \ \ \ \  \ \ \ \ \  \  \ \ \ \  \  \ \ \ \  \   \ \ \ \  \  \left( \ \uparrow \ \   \mbox{via \ \ the \ \ integration \ \ by \ \ parts \ \ formula} \ \right) \\[0.1in]
& \  & \ \ \ \  \le \ C \cdot  \lambda_{\,1}^2 \,\cdot\,\bigg\vert \ \int_{\R^n}  \left[\ {{\partial^{\,2} \,(\,\Delta\,V_{\,l}\,) }\over {\partial\, \lambda_{\,1} \,\partial\, \xi_{1_{\,|\,1} }}} \ \right]\,\cdot\,\phi_{\ \flat}  \ \bigg\vert \\[0.2in]
& \  & \ \ \ \  \le \ C_1 \cdot \lambda_{\,1}^2\,\cdot\, \bigg\vert \ \int_{\R^n} \left[\ {{\partial^{\,2}\! }\over {\partial\, \lambda_{\,1} \,\partial \,\xi_{1_{\,|\,1} }}}\,\left(\ \, V_{\,l}^{{\,n\,+\,2\,}\over {n\,-\,2}} \,\right) \ \right]\,\cdot\,\phi_{\ \flat} \ \bigg\vert
 \ \ \ \ \ \  \left[\ {\mbox{applying \ \ equation \ \ (\,A.1.6\,)}}\ \right]\\[0.1in]
& \ &  \ \ \ \ \le \ C_2\cdot  \lambda_{\,1}^2\,\cdot\,{1\over {\lambda^2_1}}\,\cdot\,\int_{\R^n} V_1^{{\,n\,+\,2\,}\over {n\,-\,2}}\,\cdot\,|\ \phi_{\ \flat}\,| \\[0.1in]
& \ &  \hspace*{-0.5in}   \left(\ \ {\mbox{as}} \ \  \bigg\vert \ {{\partial \,V_{\,l} }\over {\partial\, \lambda_{\,1} }} \ \bigg\vert  \ \le \ C_3\,\cdot\,{1\over {\lambda_{\,1}}}\,\cdot\,V_{\,l} \,, \ \  \bigg\vert \ {{\partial \,V_{\,l} }\over { \partial \,\xi_{1_{\,|\,1} }}}\ \bigg\vert  \ \le \ C_4\,\cdot\,{1\over {\lambda_{\,1}}}\,\cdot\,V_{\,l}  \ \ {\mbox{and}} \ \   \bigg\vert \ {{\partial^2 \,V_{\,l} }\over {\partial\, \lambda_{\,1} \,\partial \,\xi_{1_{\,|\,1} }}}\ \bigg\vert  \ \le \ C_5\,\cdot\,{1\over {\lambda_{\,1}^2}}\,\cdot\,V_{\,l} \ \right) \\[0.2in]
& \ &   \le \  C_6 \,\cdot\,\left(\  \int_{\R^n}   V^{{2n}\over {\,n\,-\ 2\,}}\  \right)^{\! {{n\,+\,2}\over {2n}} }\,\cdot\,\,\left( \int_{\R^n}    |\,\phi_{\ \flat}\,|^{{2n}\over {\,n\,-\ 2\,}}  \right)^{\! {{n\,-\,2}\over {2n}} } \ \,\le \ C_7\,\cdot\,\Vert\,\phi_{\ \flat}\,\Vert_\btd \  .\\
\end{eqnarray*}

Here we invoke the Sobolev inequality\,\, and (\,A.3.17\,)\,.\, It follows that\\[0.1in]
(\,A.7.11\,)
\begin{eqnarray*}
  \bigg\vert \ \bigg\langle \,{\bf I}\,'\,(W_{\,\,\flat} \ + \ \phi_{\ \flat})\,, \  \left(\ \lambda_{\,1}\,\cdot\,{ { \partial\, \phi_{\ \flat} } \over {\partial\, \lambda_{\,1}}}\  \right)  \ \bigg\rangle_\btd \ \bigg\vert    & \le  & C \cdot \flat\cdot\,[\ \max \ \{ \  |\,c_{\,l}\,|\,, \ \ |\, c_{\,l_{\,j}}\,|\ \}  \ ]\,\cdot\, \Vert\,\phi_{\ \flat}\,\Vert_\btd  \\[0.2in]& \le &\  C_1\cdot {\bar\lambda}_{\ \flat}^{\,2\,{{\varpi}_{\,(\,A.5.38\,)} } \ - \ \sigma}  \ .\
\end{eqnarray*}
Recall from Proposition A.5.23 that
$$
 \Vert \, \phi_{\ \flat} \,\Vert_\btd  \ \le \  C\,\cdot\,{\bar\lambda}_{\ \flat}^{{\,{{\varpi}_{\,(\,A.5.38\,)} }}}  \ , \ \ \ \ \max \ \{ \ |\,c_{\,l}\,|\,, \ \ |\, c_{\,l_{\,j}}\,|\ \}  \ \le \  C\,\cdot\, {\bar\lambda}_{\ \flat}^{{\,{{\varpi}_{\,(\,A.5.38\,)} }} } \ \ \ \ {\mbox{and}} \ \ \ \
\flat \ \le \  {1\over {{\bar\lambda}_{\ \flat}^{\,\sigma}  }}  \  .
$$

 \vspace*{0.1in}

 {\it Note.}  A.7.12. \ \ For this part of the argument, one can replace $\, \left( \ \lambda_{\,1}\,\cdot\,{ { \partial\, V_1 } \over {\partial\, \lambda_{\,1}}}\  \right)\,\,$ by $\,\, \left( \ \lambda_{\,1}\,\cdot\,{ { \partial\, V_1} \over {\partial\, \xi_{1_{\,j}}}}\  \right)\,.$ via the inequality $\,\, \vert \ \lambda_{\,1}\,\cdot\,{ { \partial\, V_1} \over {\partial\, \xi_{1_{\,j}}}} \vert \ \le \ C \,V_1\,,\ ....$
etc.

\newpage

{\bf \S\,A\,7.\,d\,.} \ \
{\bf Estimate of ($\,{\bf I})_{\,(\,A.7.8\,)}\,$}\,.  \\[0.2in]   We have
\begin{eqnarray*}
(\,A.7.13\,) \ \ \ \  & \ & \int_{\R^n}\,\bigg\langle\btd\,(\, W_{\,\,\flat} \ + \ \phi_{\ \flat}\,)\,,\,\ \btd\,\left( \lambda_{\,1}\cdot{ { \partial \,W_{\,\,\flat} } \over {\partial\, \lambda_{\,1}}}   \right)\,\bigg\rangle\\[0.2in]
   & = & \int_{\R^n}\,\bigg\langle\,\btd  W_{\,\,\flat}\,,\ \,\btd\,\left(\,\lambda_{\,1}\cdot { { \partial \,V_1 } \over {\partial\, \lambda_{\,1}}}   \right)\,\bigg\rangle \ + \ \int_{\R^n}\,\bigg\langle\,\btd\,\phi_{\ \flat}\,,\ \btd\,\left(\,\lambda_{\,1}\cdot{ { \partial \,V_1} \over {\partial\, \lambda_{\,1}}}   \right)\,\bigg\rangle \\[0.1in]
   & \ & \hspace*{2.7in} (\,\uparrow \ \ = \ 0\ \ \ {\mbox{via \ \ the  \ \ }}\perp\,-\,{\mbox{condition}}\,)\ \ \ \ \  \ \ \ \ \ \ \ \  \ \ \  \\[0.1in]
& = &  -\ \int_{\R^n} (\Delta\,W_{\,\,\flat})\,\cdot\,\left(\ \lambda_{\,1}\cdot { { \partial \,V_1} \over {\partial\, \lambda_{\,1}}}   \right) \ \ \ (\,\leftarrow \,\ {\mbox{part \ \ of \ \ the \ \ key \ \ term}}\,) \  {\bf .}\\
\end{eqnarray*}

\vspace*{0.15in}

{\bf \S\,A\,7.\,e\,.} \ \ {\bf Estimate of ($\,{\bf II})_{\,(\,A.7.8\,)}\,$}\,.  \\[0.2in]
We start with
\begin{eqnarray*}
(\,A.7.14\,)   & \ & \int_{\R^n}\,(\,{\tilde c}_n\!\cdot K\,)\,(W_{\,\,\flat} \ + \ \phi_{\ \flat})_+^{{\,n\,+\,2\,}\over {n\,-\,2}}\,\cdot\,\left(\,\lambda_{\,1}\cdot { { \partial \,V_1 } \over {\partial\, \lambda_{\,1}}}\ \right)\\[0.15in]
& = &  \int_{\R^n}\,(\,{\tilde c}_n\!\cdot K\,)\,(W_{\,\,\flat} )^{{\,n\,+\,2\,}\over {n\,-\,2}}\,\cdot\,\left(\ \lambda_{\,1}\cdot { { \partial\,V_1 } \over {\partial\, \lambda_{\,1}}}\ \right) \ +  \\[0.1in]
& \ &  \ \ \ \ \ \ \ \ \ (\,\uparrow \,\ {\mbox{part \ \ of \ \ the \ \ key \ \ term}}\,) \\[0.2in]
& \ &  \ \  + \  \int_{\R^n}\,(\,{\tilde c}_n\!\cdot K\,)\,\left[\ (W_{\,\,\flat} \ + \ \phi_{\ \flat})_+^{{\,n\,+\,2\,}\over {n\,-\,2}} - \ (W_{\,\,\flat} )^{{\,n\,+\,2\,}\over {n\,-\,2}} \ \right] \cdot\! \left(\ \lambda_{\,1}\cdot { { \partial \,V_1} \over {\partial\, \lambda_{\,1}}}\ \right) \cdot \cdot \cdot \,{\bf (\,E\,)}_{ (\,A.7.13\,)}\ .
\end{eqnarray*}
Summing up (\,A.7.9\,)\,,\, (\,A.7.13\,)\,  and \,(\,A.7.14\,)\,,\, we obtain\\[0.1in]
(\,A.7.15\,)
\begin{eqnarray*}
 \lambda_{\,1}\cdot{{\partial \, {\bf I}_{\,\cal R} }\over {\partial\, \lambda_{\,1}}}  & = &  -\int_{\R^n} (\Delta\,W_{\,\,\flat})\,\cdot\,\left(\ \lambda_{\,1}\cdot { { \partial \,V_1} \over {\partial\, \lambda_{\,1}}}   \right) \ - \  \int_{\R^n}\,(\,{\tilde c}_n\!\cdot K\,)\,(W_{\,\,\flat} )^{{\,n\,+\,2\,}\over {n\,-\,2}}\,\cdot\,\left(\ \lambda_{\,1}\cdot { { \partial\,V_1 } \over {\partial\, \lambda_{\,1}}}\ \right)\\[0.2in]
 & \ & \ \ \ \  \ + \ O\,\left( \,{\bar\lambda}_{\ \flat}^{\,2\,{{\varpi}_{\,(\,A.5.38\,)} } \ - \ \sigma}  \, \right) \ +  \ {\bf (\,E\,)}_{ (\,A.7.13\,)}\\[0.2in]
 & = & \lambda_{\,1}\cdot{{\partial \, {\bf I}\,(\,W_{\ \flat} \,) }\over {\partial\, \lambda_{\,1}}}   \ + \ O\,\left( \,{\bar\lambda}_{\ \flat}^{\,2\,{{\varpi}_{\,(\,A.5.38\,)} } \ - \ \sigma}  \, \right) \ +  \ {\bf (\,E\,)}_{ (\,A.7.13\,)}\ .
\end{eqnarray*}

 \vspace*{0.1in}

 {\it Note.}  A.7.16. \ \ This part of the argument works in a similar fashion for  $\,\, \left( \ \lambda_{\,1}\,\cdot\,{ { \partial\, V_1} \over {\partial\, \xi_{1_{\,j}}}}\  \right)\,,....$,\, etc.

\vspace*{0.3in}

{\bf \S\,A\,7.\,f\,.} \ \ {\bf Estimate of $\,{\bf (\,E\,)}_{ (\,A.7.13\,)}\,$}\,.  \\[0.2in]
Consider the expression
$$
\left[\ (W_{\,\,\flat} \ + \ \phi_{\ \flat})_+^{{\,n\,+\,2\,}\over {n\,-\,2}} \ - \ (W_{\,\,\flat} )^{{\,n\,+\,2\,}\over {n\,-\,2}} \ \right] \,\cdot\,\left(\ \lambda_{\,1}\cdot { { \partial \,V_1} \over {\partial\, \lambda_{\,1}}}\ \right)\ .
$$
As
$$
W_{\,\,\flat}\,(\,y\,) \ + \ \phi_{\ \flat}\,(\,y\,)  \ < \ 0 \ \ \Longleftrightarrow \ \   (\ 0 \ < \ ) \ \ W_{\,\,\flat}\,(\,y\,) \ < \  -\,\phi_{\ \flat}\,(\,y\,)      \ \ \Longrightarrow \ \ W_{\,\,\flat}\,(\,y\,)  \ < \  |\,\phi_{\ \flat}\,(\,y\,) \,|\,.
$$
It follows that

\vspace*{-0.35in}

\begin{eqnarray*}
\Omega_{\ -} \!\!\!\! & \subset &\!\!\!\!  \Omega_{\,\,||}\ , \\[0.2in]
{\mbox{where}} \ \ \ \  \  \Omega_{\,\,-} & = & \{ \ y \,\in\,\R^n \ | \ \,W_{\,\,\flat}\,(\,y\,) \ + \ \phi_{\ \flat}\,(\,y\,)  \ < \ 0 \ \} \\[0.15in]
 {\mbox{and}} \ \   \ \ \Omega_{\,\,||} & = & \{ \  y \,\in\,\R^n \ | \ \,W_{\,\,\flat}\,(\,y\,) \ <  \  |\,\phi_{\ \flat}(\,y\,)\ | \ \} \ .
\end{eqnarray*}
Thus\\[0.1in]
(\,A.7.17\,)
\begin{eqnarray*}
& \ & \bigg\vert \ \int_{\Omega_{\,\,-}}\,(\,{\tilde c}_n\!\cdot K\,)\,\left[\ (\,W_{\,\,\flat} \ + \ \phi_{\ \flat}\,)_+^{{\,n\,+\,2\,}\over {n\,-\,2}} \ - \ (W_{\,\,\flat} )^{{\,n\,+\,2\,}\over {n\,-\,2}} \ \right] \,\cdot\, \left(\ \lambda_{\,1}\cdot { { \partial \,V_1} \over {\partial\, \lambda_{\,1}}}\ \right) \ \bigg\vert \\[0.15in]
& = & \bigg\vert \  \int_{\Omega_{\,\,-}}\,(\,{\tilde c}_n\!\cdot K\,)\,(W_{\,\,\flat} )^{{\,n\,+\,2\,}\over {n\,-\,2}}\,\cdot\, \left(\,\lambda_{\,1}\cdot { { \partial \,V_1 } \over {\partial\, \lambda_{\,1}}}\ \right) \ \bigg\vert \ \ \ \ \ \ \ \ \ \ \ \ \  [\ (\,W_{\,\,\flat} \ + \ \phi_{\ \flat}\,)_+ \ = \ 0 \ \ \ \ {\mbox{in}} \ \ \Omega_{\,\,-} \ ] \\[0.15in]& \le & \int_{\Omega_{\,||}}\,(\,{\tilde c}_n\!\cdot K\,)\,(W_{\,\,\flat} )^{{\,n\,+\,2\,}\over {n\,-\,2}}\,\cdot\,\bigg\vert\ \lambda_{\,1}\cdot { { \partial\,V_1 } \over {\partial\, \lambda_{\,1}}} \ \bigg\vert \\[0.15in]
& \le & C_1\, \int_{\Omega_{\,||}} |\,\phi_{\ \flat}\,|^{{\,n\,+\,2\,}\over {n\,-\,2}} \,\cdot  V_1  \ \ \ \  \ \ \ \ \  \ \ \ \  \left(\ \  W_{\,\,\flat} \ <  \  |\,\phi_{\ \flat}\,|  \ \ {\mbox{in}} \ \ \Omega_{\,\,||}\,, \ \ \ \  {\mbox{and}} \ \ \ \bigg\vert\ \lambda_{\,1}\cdot { { \partial\,V_1 } \over {\partial\, \lambda_{\,1}}} \ \bigg\vert  \ \le \ C \ V_1\ \right) \\[0.2in]
& \le &  C_2\,\cdot\,\left(\, \int_{\R^n} |\,\phi_{\ \flat}\,|^{{2\,n}\over {\,n \,- \,2\,}}  \right)^{{n + 2\,}\over {2\,n}} \,\cdot\,\left(\, \int_{\R^n} [\, V_1\,]^{{2\,n}\over {\,n\, - \,2\,}}  \right)^{{\,n \, - \,2}\over {2\,n}} \\[0.2in]
&  = & C_3\,\cdot\,\Vert\,\phi_{\ \flat}\,\Vert_{\,\btd}^{{\,n\,+\,2\,}\over {n\,-\,2}} = \ O\left( \ {\bar{\lambda}}_{\flat}^{  \, {{\varpi}_{\,(\,A.5.38\,)} } }\,\cdot\,{\bar{\lambda}}_{\flat}^{  \, {4\over {n\,-\,2}} \ \cdot\, {{\varpi}_{\,(\,A.5.38\,)} } }  \  \right)\ .
\end{eqnarray*}

Here we apply ${\bf (\,i\,)_{\,5.57} }$\,.\, Similar argument works under the condition $\, W_{\ \flat}\,(y\,) \ \le \ 2 \,\phi_{\ \flat}(\,y\,)\,\,$.\,
Next, consider the portion
$$
\Omega^{\,+} \ = \ \{ \ W_{\,\,\flat} \ + \ \phi_{\ \flat} \ \ge   \ 0 \ \} \ .
$$
In this case
\begin{eqnarray*}
& \ & \bigg\vert \ \int_{\Omega^{\,+}}\,(\,{\tilde c}_n\!\cdot K\,)\,\left[\ (W_{\,\,\flat} \ + \ \phi_{\ \flat})_+^{{\,n\,+\,2\,}\over {n\,-\,2}} \ - \ (W_{\,\,\flat} )^{{\,n\,+\,2\,}\over {n\,-\,2}} \ \right] \,\cdot\, \left(\ \lambda_{\,1}\cdot { { \partial \,V_1} \over {\partial\, \lambda_{\,1}}}\ \right) \ \bigg\vert \\[0.2in]& = & \bigg\vert \ \int_{\Omega^{\,+}}\,(\,{\tilde c}_n\!\cdot K\,)\,\left[\ (W_{\,\,\flat} \ + \ \phi_{\ \flat})^{{\,n\,+\,2\,}\over {n\,-\,2}} \ - \ (W_{\,\,\flat} )^{{\,n\,+\,2\,}\over {n\,-\,2}} \ \right] \,\cdot\,\left(\ \lambda_{\,1}\cdot { { \partial \,V_1} \over {\partial\, \lambda_{\,1}}}\ \right) \ \bigg\vert  .\\
\end{eqnarray*}

Using a version of Taylor expansion [\ cf. (\,A.5.9\,) and (\,A.5.10\,) \,]\,,\,  we obtain \\[0.1in]
(\,A.7.18\,)
\begin{eqnarray*}
& \ & \!\!\!\!\!\!\!\!\!(\,1 \,+\, a\,)^{{n\, +\, 2}\over {n\, - \,2}} \ = \ 1 \ + \ \left(\, {{n \,+\, 2}\over {n \, - \,2}} \,\right)\,\cdot\,a \ + \ O (\,1)\,\cdot\,\,a^{{n\, + \,2}\over {n \,-\, 2}} \ \ \ \  \ \ \ \ \ \ \ \ \ \ \ \ \ \ \ \ \ \  \ \, \mfor \ \ |\,a\,|\ \le \ {1\over 2}\\[0.2in]
& \Longrightarrow & \ \ \bigg\vert \ (\,1 \,+\, a\,)^{{n\, +\, 2}\over {n\, - \,2}} \ - \ 1 \ - \  \left(\, {{n \,+\, 2}\over {n \, - \,2}} \,\right)\,\cdot\,a \  \bigg\vert \ \le \ C\,\cdot\,a^{{n\, +\, 2}\over {n\, - \,2}} \ \ \ \ \ \ \ \ \ \ \ \ \ \ \ \ \ \ \mfor \ \ |\,a\,|\ \le\  {1\over 2}\,,\\[0.25in]
& \  & \!\!\!\!\!\!\!\!\!\!\!\!\!\!\!\!(\,b \,+\, c\,)^{{n \,+\, 2}\over {n \,-\, 2}} \ = \ \ b^{{n\, + \, 2}\over {n\, -\, 2}}\,\cdot\,\left(\, 1 \,+\, {c\over b}\right)^{{n \,+\, 2}\over {n \,-\, 2}} \ = \  b^{{n \, + \,2}\over {n\, -\, 2}}\,\cdot\,\left[\,1 \, + \, \left(\, {{n \,+\, 2}\over {n \, - \,2}} \,\right)\,\cdot\, {c\over b} \, + \, O (\,1)\,\cdot\,\left(\, {c\over b}\right)^{{n + 2}\over {n - 2}}\  \right]\\[0.1in]
& \ & \hspace*{4.8in}\ \  \mfor \ \ \  \bigg\vert \, {{\,c\,}\over b} \bigg\vert\  \le \ {1\over 2} \\[0.1in]
& \Longrightarrow & \ \ \bigg\vert\  (\,b \,+\, c\,)^{{n \,+\, 2}\over {n \,-\, 2}}  \ -\  b^{{n \,+\, 2}\over {n\, - \,2}} \ \bigg\vert \ =  \ \left(\, {{n + 2}\over {n - 2}}\,\right)\,\cdot\,b^{4\over {n - 2}}\,\cdot\,c\ + \ O \left(\ \,c^{{n\, + \,2}\over {n\, - \,2}}\,\right) \ \ \ \  \mfor \ \ \bigg\vert \, {{\,c\,}\over b}\, \bigg\vert \ \le\  {1\over 2}\,.\\
\end{eqnarray*}
Here $\,a\,$, $\,b \, > \, 0\,$ and $\,c\ (\,\pm\,)\,$ are numbers satisfying the indicated conditions\,.\, Let
$$
\Omega^{\,+}_{\le\,{1\over 2} } \ = \ \left\{ \ y\ \in \ \R^n \ \  \ \bigg\vert \ \ \ \  \bigg\vert \ {{\ \phi_{\ \flat} \,(\,y\,) \ }\over {W_{\,\,\flat}\,(\,y\,)}}\  \bigg\vert \ \le\  {1\over 2} \ \right\} \ . \leqno (\,A.7.19\,)
$$

\newpage

Applying (\,A.7.18\,) and (\,A.7.19\,)\,,\, we obtain\\[0.1in]
(\,A.7.20\,)
\begin{eqnarray*}
& \ &  \int_{\Omega^{\,+}_{\le\,{1\over 2} }}\,(\,{\tilde c}_n\!\cdot K\,)\,\left[\ (W_{\,\,\flat} \ + \ \phi_{\ \flat})_+^{{\,n\,+\,2\,}\over {n\,-\,2}} \ - \ (W_{\,\,\flat} )^{{\,n\,+\,2\,}\over {n\,-\,2}} \ \right] \,\cdot\, \left(\ \lambda_{\,1}\cdot { { \partial \,V_1} \over {\partial\, \lambda_{\,1}}}\ \right)   \\[0.2in]
 & = &  \ \left(\ {{n\,+\,2}\over {\,n\,-\ 2\,}}\  \right)\,\cdot\, \int_{\Omega^{\,+}_{\le\,{1\over 2} }}\, (\,{\tilde c}_n\!\cdot K\,)\,(W_{\,\,\flat} )^{{4}\over {\,n\,-\,2\,}}\,\cdot\,\left(\,\lambda_{\,1}\cdot { { \partial \,V_1 } \over {\partial\, \lambda_{\,1}}}\ \right)\,\cdot\,\phi_{\ \flat}  \   +\  O \left(\ \ \int_{\R^n} |\,\phi_{\ \flat}\,|^{{2n}\over {\,n\,-\ 2\,}}\  \right)\\[0.2in]
& \ & \hspace*{2in} \ \ \ \ \left[\ \ O \left(\, \Vert\,\phi_{\ \flat}\,\Vert_{\,\btd}^{{2n}\over {\,n\,-\ 2\,}}\  \right) \ = \ O\,\left(\  {\bar\lambda}_{\ \flat}^{\,{{2n}\over {\,n\,-\ 2\,}}  \cdot \,{{\varpi}_{\,(\,A.5.38\,)} }}\  \right) \  \uparrow  \ \right]\ .
\end{eqnarray*}
Let
$$
\Omega^{\,+}_{\ge\,{1\over 2} } \ = \ \left\{ \ y\ \in \ \R^n \ \  \bigg\vert \ \  \bigg\vert \ {{\phi_{\ \flat} \,(\,y\,) }\over {W_{\,\,\flat}\,(\,y\,)}}\  \bigg\vert \ \ge\  {1\over 2} \ \Longleftrightarrow \ \ |\,\phi_{\ \flat} \,(\,y\,)\,| \ \ge \ 2\,\cdot\,|\,W_{\,\,\flat}\,(\,y\,)\,| \  \right\} \ .
$$

Calculate as in (\,A.7.17\,)\,,\, we obtain
\begin{eqnarray*}
(\,A.7.21\,) \ \ \ \ & \ & \Bigg\vert \  \int_{\Omega^{\,+}_{\ge\,{1\over 2} }}\,(\,{\tilde c}_n\!\cdot K\,)\,\left[\ (W_{\,\,\flat} \ + \ \phi_{\ \flat})_+^{{\,n\,+\,2\,}\over {n\,-\,2}} \ - \ (W_{\,\,\flat} )^{{\,n\,+\,2\,}\over {n\,-\,2}} \ \right] \,\cdot\, \left(\ \lambda_{\,1}\cdot { { \partial \,V_1} \over {\partial\, \lambda_{\,1}}}\ \right) \ \Bigg\vert   \ \ \ \ \ \ \ \ \ \ \ \\[0.2in]
& \le &  C_2\,\cdot\,\left(\, \int_{\R^n} |\,\phi_{\ \flat}\,|^{{2\,n}\over {n \,- \,2}}  \right)^{{n \,+\, 2}\over {2\,n}} \,\cdot\,\left(\, \int_{\R^n} |\, V_1\,|^{{2\,n}\over {\,n\,-\,2\,}}  \right)^{{n - 2}\over {2n}}\\[0.2in]& = & C_3\,\cdot\,\Vert\,\phi_{\ \flat}\,\Vert_{\,\btd}^{{\,n\,+\,2\,}\over {n\,-\,2}}  \ = \ O\,\left( \ {\bar{\lambda}}_{\flat}^{  \, {{\varpi}_{\,(\,A.5.38\,)} } }\,\cdot\,{\bar{\lambda}}_{\flat}^{  \, {4\over {n\,-\,2}} \ \cdot\, {{\varpi}_{\,(\,A.5.38\,)} } }  \  \right)\ .\\
\end{eqnarray*}
Thus we obtain
\begin{eqnarray*}
& \ &  \int_{\R^n}\,(\,{\tilde c}_n\!\cdot K\,)\,\left[\ (W_{\,\,\flat} \ + \ \phi_{\ \flat})_+^{{\,n\,+\,2\,}\over {n\,-\,2}} \ - \ (W_{\,\,\flat} )^{{\,n\,+\,2\,}\over {n\,-\,2}} \ \right] \,\cdot\,\left(\ \lambda_{\,1}\cdot { { \partial \,V_1} \over {\partial\, \lambda_{\,1}}}\ \right)   \\[0.2in]
 & = &  \ \left(\ {{n\,+\,2}\over {\,n\,-\ 2\,}}\  \right)\,\cdot\, \int_{\Omega^{\,+}_{\le\,{1\over 2} }}\, (\,{\tilde c}_n\!\cdot K\,)\,(W_{\,\,\flat} )^{{4}\over {\,n\,-\,2\,}}\,\cdot\,\left(\,\lambda_{\,1}\cdot { { \partial \,V_1 } \over {\partial\, \lambda_{\,1}}}\ \right)\,\cdot\,\phi_{\ \flat}  \   +\   O\,\left( \ {\bar{\lambda}}_{\flat}^{ {{n\,+\,2}\over {\,n\,-\ 2\,}}\,\,\cdot\,\,{{\varpi}_{\,(\,A.5.38\,)} } } \  \right)\\[0.2in] & = &  \ \left(\ {{n\,+\,2}\over {\,n\,-\ 2\,}}\  \right)\,\cdot\, \int_{\R^n}\,(W_{\,\,\flat} )^{{4}\over {\,n\,-\,2\,}}\,\cdot\,\left(\,\lambda_{\,1}\cdot { { \partial \,V_1 } \over {\partial\, \lambda_{\,1}}}\ \right)\,\cdot\,\phi_{\ \flat}  \   +\   O\,\left( \ {\bar{\lambda}}_{\flat}^{ {{n\,+\,2}\over {\,n\,-\ 2\,}}\,\,\cdot\,\,{{\varpi}_{\,(\,A.5.38\,)} } } \  \right)\ .
\end{eqnarray*}
$\bigg\{$\ For the portions $\,\Omega_{\,-}\,$ and  $\,\,\Omega^{\,+}_{\ge\,{1\over {\,2}} }$\, [\, that is\,,\, $\, W_{\ \flat}\,(y\,) \ \le \ 2 \,\phi_{\ \flat}(\,y\,)\ ]\,$,\,  cf. (\,A.7.17\,)\ \bigg\}$\ . $\bk
Via the $\,\perp$\,--\,condition [\ cf. (\,A.4.11\,) \ ]\,,\,
$$
\int_{\R^n}\, V_1^{{4}\over {\,n\,-\,2\,}}\,\cdot\,\left(\,\lambda_{\,1}\cdot { { \partial \,V_1 } \over {\partial\, \lambda_{\,1}}}\ \right)\,\cdot\,\phi_{\ \flat} \ = \ {1\over {n\,(\,n \ + \ 2\,)  }} \cdot\,\bigg\langle \ \left(\,\lambda_{\,1}\cdot { { \partial \,V_1 } \over {\partial\, \lambda_{\,1}}}\ \right)\ , \ \ \phi_{\ \flat} \ \bigg\rangle_\btd \ = \ 0\,, \leqno (\,A.7.22\,)
$$
we may write\\[0.1in]
(\,A.7.23\,)
\begin{eqnarray*}
& \ & \int_{\R^n}\,(\,{\tilde c}_n\!\cdot K\,)\,(W_{\,\,\flat} )^{{4}\over {\,n\,-\,2\,}}\,\cdot\,\left(\ \lambda_{\,1}\cdot{ { \partial \,V_1 } \over {\partial\, \lambda_{\,1}}}\ \right)\,\cdot\,\phi_{\ \flat}\\[0.15in]
& = & \int_{\R^n}\,(\,{\tilde c}_n\!\cdot K\,)\,\left[\ (W_{\,\,\flat} )^{{4}\over {\,n\,-\,2\,}}\,\cdot\,\left(\ \lambda_{\,1}\cdot{ { \partial \,V_1 } \over {\partial\, \lambda_{\,1}}}\ \right) \ - \  V_1^{{4}\over {\,n\,-\,2\,}}\,\cdot\, \left(\,\lambda_{\,1}\cdot { { \partial \,V_1 } \over {\partial\, \lambda_{\,1}}}\ \right)\, \right]\,\cdot\,\phi_{\ \flat} \ \cdot\,\cdot\,\cdot\,\cdot\,\cdot\,\cdot \ \ ({\bf I}_{\,(\,A.7.23\,)}\ ) \\[0.15in]
& \ & \ \ \ \ \ \ \ + \ \int_{\R^n}\,[\ (\,{\tilde c}_n\!\cdot K\,)\ - \ n\,(\,n\,-\,2)\ ]\,\cdot\,\left[\    V_1^{{4}\over {\,n\,-\,2\,}}\,\cdot\,\left(\,\lambda_{\,1}\cdot  { { \partial \,V_1 } \over {\partial\, \lambda_{\,1}}}\ \right)\, \right]\,\cdot\,\phi_{\ \flat}\,.  \ \cdot\,\cdot\cdot\,\cdot\,\cdot\,\cdot\, \ \ ({\bf II}_{\,(\,A.7.23\,)}\ )\\[0.2in]
& \ & \ \ \ \ \  \ \  \ \  \ \  \ \  \ + \  n\,(\,n\,-\,2)\,\cdot\,\int_{\R^n}\, V_1^{{4}\over {\,n\,-\,2\,}}\,\cdot\, \left(\,\lambda_{\,1}\cdot { { \partial \,V_1 } \over {\partial\, \lambda_{\,1}}}\ \right)\,\cdot\, \phi_{\ \flat} \ \  ( \ \,= \ 0\,)\ .
\end{eqnarray*}
Before we proceed further, let us conclude this section with
 \begin{eqnarray*}
(\,A.7.24\,)\ \ \ \ \ {\bf (\,E\,)}_{ (\,A.7.13\,)}\,\ \bigg[ & = &   \int_{\R^n}\,(\,{\tilde c}_n\!\cdot K\,)\,\left[\ (W_{\,\,\flat} \ + \ \phi_{\ \flat})_+^{{\,n\,+\,2\,}\over {n\,-\,2}} - \ (W_{\,\,\flat} )^{{\,n\,+\,2\,}\over {n\,-\,2}} \ \right] \cdot\! \left(\ \lambda_{\,1}\cdot { { \partial \,V_1} \over {\partial\, \lambda_{\,1}}}\ \right)\ \bigg]\ \ \ \ \ \ \ \ \   \\[0.2in]
& = & (\ {\bf I}_{\,(\,A.7.23\,)}\ ) \ + \  (\ {\bf II}_{\,(\,A.7.23\,)}\ ) \   +\   O\,\left( \ {\bar{\lambda}}_{\flat}^{ {{n\,+\,2}\over {\,n\,-\ 2\,}}\,\,\cdot\,\,{{\varpi}_{\,(\,A.5.38\,)} } } \  \right)\ .
\end{eqnarray*}

 \vspace*{0.1in}

 {\it Note.}  A.7.25.\ \ This part of the argument works in a similar way for  $\,\, \left( \, \lambda_{\,1}\,\cdot\,{ { \partial\, V_1} \over {\partial\, \xi_{1_{\,j}}}} \ \right)\,$.

 \vspace*{0.3in}

{\bf \S\,A\,7.\,g\,.} \ \
{\bf Estimate of \,$(\,{\bf I}_{\,(\,A.7.23\,)}\,)$\ .\,} \ \  We begin with\\[0.1in]
 (\,A.7.26\,)
\begin{eqnarray*}
n \ \ge \ 6 \ \ \Longrightarrow & \ & {4\over {n \ - \ 2}} \ \le \ 1\\[0.2in]
\Longrightarrow & \ & (\ A_1\ + \ A_2 \ + \ \cdot \cdot \cdot \ + \ A_k\, )^{4\over {\,n \ -  \ 2\,}} \ \le \
 A_1^{4\over {\,n \ -  \ 2\,}}\ + \ A_2^{4\over {\,n \ -  \ 2\,}} \ + \ \cdot \cdot \cdot \ + \ A_k^{4\over {\,n \ -  \ 2\,}} \\[0.1in]
  & \ & \hspace*{1.2in}\ \ \ \mfor \ \ A_i \ \ge \ 0\,, \ \ \ k \,\in\,\N \ , \ \ \ i \ = \ 1\,,\, \ 2\,, \ \cdot \cdot \cdot\,, \ k \ \\[0.1in]
\Longrightarrow & \ &   [\ W_{\,\,\flat}\,(\,y\, )\ ]^{\,{{4}\over {n\,-\,2}}} \ \le \    \left[\ V_1^{{4}\over {\,n\,-\,2\,}}  \ + \ V_2^{{4}\over {\,n\,-\,2\,}} \,+\,V_3^{{4}\over {\,n\,-\,2\,}}  \ + \ \cdot\,\cdot\,\cdot\,\ + \ V_{\,\flat}^{{4}\over {\,n\,-\,2\,}}\ \right]_{\,y} \\[0.2in] \Longrightarrow & \ &[\  W_{\,\,\flat}\,(\,y\, )\ ]^{{4}\over {\,n\,-\,2\,}} \ - \ [\  V_1\,(\,y\, )\ ]^{{4}\over {\,n\,-\,2\,}} \ \le \  \left[\  \ V_2^{{4}\over {\,n\,-\,2\,}} \,+\,V_3^{{4}\over {\,n\,-\,2\,}}  \ + \ \cdot\,\cdot\,\cdot\,\ + \ V_{\,\flat}^{{4}\over {\,n\,-\,2\,}}\ \right]_{\,y}\,.
 \end{eqnarray*}
Next\,,\,\\[0.1in]
 (\,A.7.27\,)
\begin{eqnarray*}
``\,| \ (\ {\bf I}_{\,(\,A.7.23\,)}\ ) \,|\," & \le & \int_{\R^n} \bigg\vert\,(W_{\,\,\flat} )^{{4}\over {\,n\,-\,2\,}}\,\cdot\,\left(\,\lambda_{\,1}\cdot  { { \partial \,V_1 } \over {\partial\, \lambda_{\,1}}}\ \right) \ - \  V_1^{{4}\over {\,n\,-\,2\,}}\,\cdot\,\left(\,\lambda_{\,1}\cdot  { { \partial \,V_1 } \over {\partial\, \lambda_{\,1}}}\ \right) \,\bigg\vert \,\cdot\, |\,\phi_{\ \flat}\,| \\[0.2in]
& \le & \int_{\R^n} \left[\  (W_{\,\,\flat} )^{{4}\over {\,n\,-\,2\,}}  \ - \  V_1^{{4}\over {\,n\,-\,2\,}} \,\right]  \,\cdot\, \bigg\vert \ \lambda_{\,1}\cdot  { { \partial \,V_1 } \over {\partial\, \lambda_{\,1}}} \ \bigg\vert \,\cdot\, |\,\phi_{\ \flat}\,|  \\[0.2in]
& \le & \int_{\R^n} \left(\ \  \left[\ (W_{\,\,\flat} )^{{4}\over {\,n\,-\,2\,}}  \ - \  V_1^{{4}\over {\,n\,-\,2\,}} \,\right]\,\cdot\,V_1 \,\right)\,\cdot\, |\,\phi_{\ \flat}\,| \,dy \ \ \ \  \left( \ \bigg\vert \,\lambda_{\,1}\cdot  { { \partial \,V_1 } \over {\partial\, \lambda_{\,1}}}\,\bigg\vert  \ \le \ V_1 \ \right) \\[0.2in]
& \le & \left[\ \int_{\R^n}  \left(\   (W_{\,\,\flat} )^{{4}\over {\,n\,-\,2\,}} \ - \  V_1^{{4}\over {\,n\,-\,2\,}} \,  \right)^{{2n}\over {n\,+\,2}} \,\cdot\,V_1^{{2n}\over {n\,+\,2}} \ \right]^{{n\,+\,2}\over {2\,n}} \ \cdot\, \left(\ \  \int_{\R^n} |\,\phi_{\ \flat}\,|^{{2n}\over {\,n\,-\ 2\,}}\  \right)^{\!\!{{n\,-\,2}\over {2n}}}\\[0.2in]
& \le & C\,  \left[\ \int_{\R^n}  \left(\   (W_{\,\,\flat} )^{{4}\over {\,n\,-\,2\,}} \ - \  V_1^{{4}\over {\,n\,-\,2\,}} \,  \right)^{{2n}\over {n\,+\,2}} \,\cdot\,V_1^{{2n}\over {n\,+\,2}} \ \right]^{{n\,+\,2}\over {2\,n}} \ \cdot\,\Vert\,\phi_{\ \flat}\,\Vert_\btd\ . \\[0.1in]
& \ & \hspace*{3in} \ \ \ \ \ \ \ \ \ \ \ \ \ \ \ \  \left[\ \uparrow \ \  O\left(\ {\bar\lambda}_{\ \flat}^{\bar\omega} \  \right) \ \right] \ .
\end{eqnarray*}
We estimate the last integral in (\,A.7.27\,) in the following regions. C\,f. {\bf \S\,A\,4\,.\,j}\,.
$$
B_{\,\xi_{\,1}} (\,\rho_{\,\nu}) \,, \ \   B_{\xi_{\,2}} (\,\rho_{\,\nu}) \,, \ \ \ \cdot\,\cdot\,\cdot\,\,, \ \ \ B_{\xi_{\,\flat}} (\,\rho_{\,\nu}) \ \ \ \ \ \ \ {\mbox{and}} \ \ \ \ \  \R^n \,\setminus \, \left(\  \cup\ B_{\,\xi_{\,l}} \,(\,{\rho_{\,\nu}}) \, \right)\ .
$$
Here $\,\rho_{\,\nu} \ = \ {\bar\lambda}_{\,\,\flat}^{\,\nu}\,.\,$
We proceed with \\[0.1in]
(\,A.7.28\,)
\begin{eqnarray*}
& \ &  \left[\ \int_{\R^n}  \left(\   (W_{\,\,\flat} )^{{4}\over {\,n\,-\,2\,}} \ - \  V_1^{{4}\over {\,n\,-\,2\,}} \,  \right)^{{2n}\over {n\,+\,2}} \,\cdot\,V_1^{{2n}\over {n\,+\,2}} \ \right]^{{n\,+\,2}\over {2\,n}}\\[0.2in]
& = &  \left\{  \left[ \,   \int_{B_{\,\xi_{\,1}} (\,\rho_{\,\nu})}   \!\!\!\!+ \  \cdot \cdot \cdot  \ + \  \int_{B_{\xi_{\ \flat}} (\,\rho_{\,\nu})}  \!\!+ \  \int_{\R^n \,\setminus \, \left(\, \cup\ B_{\xi_l} \,(\,{\rho_{\,\nu}}) \right) }  \,\right]  \left[  \, (W_{\,\,\flat} )^{{4}\over {\,n\,-\,2\,}} \ - \  V_1^{{4}\over {\,n\,-\,2\,}} \,  \right]^{{2n}\over {n\,+\,2}} \!\!\cdot\,V_1^{{2n}\over {n\,+\,2}} \, \right\}^{\!{{n\,+\,2}\over {2\,n}}}\\[0.2in]
& \le &   \left[\    \int_{B_{\,\xi_{\,1}} (\,\rho_{\,\nu})}    \left(\   (W_{\,\,\flat} )^{{4}\over {\,n\,-\,2\,}} \ - \  V_1^{{4}\over {\,n\,-\,2\,}} \,  \right)^{{2n}\over {n\,+\,2}} \!\cdot\,V_1^{{2n}\over {n\,+\,2}} \ \right]^{{n\,+\,2}\over {2\,n}} \ + \ \\[0.in]
& \ & \hspace*{2in} \cdot \\[0.01in]
& \ & \hspace*{2in} \cdot \\[0.01in]
& \ & \hspace*{2in} \cdot \\[0.01in]
& \ & \ \ \ \  \ \ +  \left[\    \int_{B_{\xi_{\ \flat}} (\,\rho_{\,\nu})}    \left(\   (W_{\,\,\flat} )^{{4}\over {\,n\,-\,2\,}} \ - \  V_1^{{4}\over {\,n\,-\,2\,}} \,  \right)^{{2n}\over {n\,+\,2}} \,\cdot\,V_1^{{2n}\over {n\,+\,2}} \ \right]^{{n\,+\,2}\over {2\,n}}\\[0.2in]
& \ & \ \ \ \  \ \ \ \ \ \ \ \  \ \ +   \left[\    \int_{\R^n \,\setminus \, \left(\, \cup\ B_{\xi_l} \,(\,{\rho_{\,\nu}}) \right) }     \left(\   (W_{\,\,\flat} )^{{4}\over {\,n\,-\,2\,}} \ - \  V_1^{{4}\over {\,n\,-\,2\,}} \,  \right)^{{2n}\over {n\,+\,2}} \,\cdot\,V_1^{{2n}\over {n\,+\,2}} \ \right]^{{n\,+\,2}\over {2\,n}} \\[0.1in]
& \ & \ \ \ \ \ \ \ \ \ \ \     [\ {\mbox{on \ \ the \ \ base \ \ of}} \ \  (\,A.7.26\,)\ , \ \ {\mbox{note \ \ that \ \ }} n \ \ge \ 6\ \ \Longrightarrow \ \ {{2n}\over {n\,+\,2}} \ < \ 1 \  ] \ .
\end{eqnarray*}

\vspace*{0.2in}

{\bf (\,i\,)} \ \  In $
B_{\,\xi_{\,1}} (\,\rho_{\,\nu})\,.$ \ \ (\,Refer also to {\bf \S\,A\,4\,.j}\,.)\,  From (\,A.7.26\,)
and
\begin{eqnarray*}
(\,A.7.29\,) \ \ \  \ & \ &  V_2^{{4}\over {\,n\,-\,2\,}} \,+\,V_3^{{4}\over {\,n\,-\,2\,}}  \ + \ \cdot\,\cdot\,\cdot\,\ + \ V_{\,\flat}^{{4}\over {\,n\,-\,2\,}} \\[0.2in]
& \le &   C\,\cdot\,{1\over {{\bar\lambda}_{\ \flat}^2}}\,\cdot\,\left[\ \left(\ {1\over { 1 \ + \ \Vert\ \Xi_{\,\,1} \ -  \ \Xi_{\,2}\,\Vert  }}  \right)^{\!\!4} \ + \ \cdot\,\cdot\,\cdot\,\ + \ \left(\ {1\over { 1 \ + \ \Vert\ \Xi_{\,\,1} \ -  \ \Xi_{\ \flat} \,\Vert  }}  \right)^{\!\!4}  \ \right] \ \ \ \  \ \ \  \ \ \ \ \ \   \\[0.2in]
& \le & C\,\cdot\,{1\over {{\bar\lambda}_{\ \flat}^2}}\,\cdot\,{\bar\lambda}_{\ \flat}^{4\,\gamma} \ .
 \end{eqnarray*}
Thus
\begin{eqnarray*}(\,A.7.30\,)
& \ &  \int_{B_{\,\xi_{\,1}} (\,\rho_{\,\nu})}  \left(\   (W_{\,\,\flat} )^{{4}\over {\,n\,-\,2\,}} \ - \  V_1^{{4}\over {\,n\,-\,2\,}} \,  \right)^{{2n}\over {n\,+\,2}} \,\cdot\,V_1^{{2n}\over {n\,+\,2}} \   \\[0.2in]
& \le  & C_1 \int_{B_{\,\xi_{\,1}} (\,\rho_{\,\nu})}  \left(\    V_2^{{4}\over {\,n\,-\,2\,}} \,+\,V_3^{{4}\over {\,n\,-\,2\,}}  \ + \ \cdot\,\,\cdot\,\cdot \ + \ V_{\,\flat}^{{4}\over {\,n\,-\,2\,}}  \right)^{{2n}\over {n\,+\,2}} \,\cdot\, V_1^{{2n}\over {n\,+\,2}} \ \ \ \ \ [ \ {\mbox{via}} \ \ (\,A.7.26\,) \ ] \\[0.2in]
& \le & C_2 \cdot {1\over {{\bar\lambda}_{\ \flat}^{{ 4\,n}\over { n\,+\,2 }}  }} \cdot {\bar\lambda}_{\ \flat}^{{{ 2\,n}\over { n\,+\,2 }}\  \cdot\  4\,\gamma} \cdot \int_{B_{\,\xi_{\,1}} (\,\rho_{\,\nu})}    V_1^{{2n}\over {n\,+\,2}}   \\[0.2in]
& \le &  C_3  \cdot {1\over {{\bar\lambda}_{\ \flat}^{{ 4\,n}\over { n\,+\,2 }}  }} \cdot {\bar\lambda}_{\ \flat}^{{{ 2\,n}\over { n\,+\,2 }}\  \cdot\  4\,\gamma} \cdot  \int_0^{ \rho_{\,\nu} }\left( {\lambda_{\,1}\over {\lambda_{\,1}^2 \ + \ r^2}}\  \right)^{ {{n\,-\,2}\over {n\,+\,2}}\,\cdot\,n } \cdot\, r^{n\,-\,1}\ dr\\[0.2in]
& = & C_3  \cdot {1\over {{\bar\lambda}_{\ \flat}^{{ 4\,n}\over { n\,+\,2 }}  }} \cdot {\bar\lambda}_{\ \flat}^{{{ 2\,n}\over { n\,+\,2 }}\  \cdot\  4\,\gamma} \cdot  {{\lambda_{\,1}^n}\over {\ \lambda_{\,1}^{ {{n\,-\,2}\over {n\,+\,2}}\,\cdot\,n } \   }}\,\cdot\, \int_0^{ \rho_{\,\nu}\over \lambda_{\,1} }\left( {1\over {1\ + \ R^2}}\  \right)^{ {{n\,-\,2}\over {n\,+\,2}} \,\cdot\, n } \cdot R^{n\,-\,1}\ d\,R \ \ \ \  \ \ \  \ \ \ \ \ \\[0.2in]
& \le & C_4  \cdot {1\over { {\bar\lambda}_{\ \flat}^{{ 4\,n}\over { n\,+\,2 }}  }} \cdot {\bar\lambda}_{\ \flat}^{{{ 2\,n}\over { n\,+\,2 }}\  \cdot\  4\,\gamma} \cdot {\bar\lambda}_{\ \flat}^{{ 4\,n}\over { n\,+\,2 }} \cdot \left[\ \ln \left( { \rho_{\,\nu}\over \lambda_{\,1} } \right)  \ \right]
\ = \   C_4  \cdot {\bar\lambda}_{\ \flat}^{{{ 2\,n}\over { n\,+\,2 }}\,  \cdot\, 4\,\gamma \ - \ o_{\,{\bar\lambda}_{\,\,\flat}} \,(\,1)}  \ .  \end{eqnarray*}
Here
$$   o_{\,{\bar\lambda}_{\,\,\flat}} \,(\,1) \ \to \ 0^{\,+} \ \ \ \ \  {\mbox{as}} \ \ \ \ \ \ {\bar\lambda}_{\ \flat} \ \to \ 0^{\,+}\ .
$$

\vspace*{0.1in}


{\it Note.} \ \ A.7.31\,. \ \ Concerning the indices in (\,A.7.29\,)\,,\,
$$
n \ \ge \ 7 \ \ \Longleftrightarrow \ \ {{2\,(\,n \ - \ 2\,)}\over {n\,+\,2}} \cdot n \ >  \ n \ \ \Longrightarrow \ \ \int_0^{ \infty }\left( {1\over {1\ + \ R^2}}\  \right)^{ {{n\,-\,2}\over {n\,+\,2}} \ \cdot\  n } R^{n\,-\,1}\ d\,R \ <  \ \infty\,.
$$
$$
n \ =\ 6 \ \ \Longleftrightarrow \ \ {{2\,(\,n \ - \ 2\,)}\over {n\,+\,2}} \cdot n \ =  \ n \ \ \Longrightarrow \ \ \int_0^{ \rho_{\,\nu}\over \lambda_{\,1} } \left( {1\over {1\ + \ R^2}}\  \right)^{ {{n\,-\,2}\over {n\,+\,2}} \ \cdot\  n } R^{n\,-\,1}\ d\,R \ = \ O\,\left(\ \ln \left( { \rho_{\,\nu}\over \lambda_{\,1} } \right)  \ \right)\,.
$$
Moreover,
$$
 {{\lambda_{\,1}^n}\over {\ \lambda_{\,1}^{ {{n\,-\,2}\over {n\,+\,2}}\,\cdot\,n } \   }}\,\cdot\,{1\over {(\,\lambda_{\,1}^2\,)^{ {{2n}\over {n\,+\,2}}   } \   }}  \ = \ 1
$$
as
$$
n \ - \  {{n\,-\,2}\over {n\,+\,2}}\,\cdot\,n  \ - \ 2\cdot {{2n}\over {n\,+\,2}}    \ = \ {{\ n\,(\,n\,+\,2) \ - \ n\,(\,n\,-\,2) \,-\, 4\,n\  }\over {n\,+\,2}} \ = \ 0\,.
$$

\vspace*{0.2in}

After taking the power of $\ \displaystyle{\left(\,{{n\,+\,2}\over {2\,n}}\,\right)}\ $ into account, the contribution is given by \\[0.1in]
(\,A.7.32\,)
$$
 \left[\ \int_{B_{\,\xi_{\,1}} \,(\,\rho_{\,\nu})\,}  \left(\   (W_{\,\,\flat} )^{{4}\over {\,n\,-\,2\,}} \ - \  V_1^{{4}\over {\,n\,-\,2\,}} \,  \right)^{{2n}\over {n\,+\,2}} \,\cdot\,V_1^{{2n}\over {n\,+\,2}} \ \right]^{{n\,+\,2}\over {2\,n}} \ = \  O\left(\ {\bar\lambda}_{\ \flat}^{ 4\,\gamma\ - \ o_\lambda (\,1)}  \  \right)  \ .
 $$
\hspace*{3in}
(\ {\mbox{benchmark \ \ for \ \ this \ \ part}} \ \ $\uparrow \ $)\ .

\vspace*{0.2in}

{\bf (ii)} \ \ In $
B_{\xi_{\,2}} \,(\,\rho_{\,\nu})\,$.\, Following (\,A.7.29,)\,,\, we have\,,\,\\[0.1in]
(\,A.7.33\,)
\begin{eqnarray*}
 & \ &[\  W_{\,\,\flat}\,(\,y\, )\ ]^{{4}\over {\,n\,-\,2\,}} \ - \ [\  V_1\,(\,y\, )\ ]^{{4}\over {\,n\,-\,2\,}} \ \le \  \left[\  \ V_2^{{4}\over {\,n\,-\,2\,}} \,+\,V_3^{{4}\over {\,n\,-\,2\,}}  \ + \ \cdot\,\cdot\,\cdot\,\ + \ V_{\,\flat}^{{4}\over {\,n\,-\,2\,}}\ \right]_{\,y}\\[0.2in]
& \ & \hspace*{2in}\ \ \ \ \ \ \   \le \ [\  V_2\,(\,y\, )\ ]^{{4}\over {\,n\,-\,2\,}} \ + \  C\,\cdot\,{1\over {{\bar\lambda}_{\ \flat}^2}}\,\cdot\,{\bar\lambda}_{\ \flat}^{4\,\gamma} \ \  \mfor \ \ y \,\in\, B_{\,\xi_{\,2}} (\,\rho_{\,\nu})\,.\, \ \ \ \ \
 \end{eqnarray*}

 \newpage

Thus
\begin{eqnarray*}
(\,A.7.34\,) & \ &  \int_{B_{\,\xi_{\,2}} (\,\rho_{\,\nu})} \left[\ \left(\   (W_{\,\,\flat} )^{{4}\over {\,n\,-\,2\,}} \ - \  V_1^{{4}\over {\,n\,-\,2\,}} \,  \right)^{{2n}\over {n\,+\,2}} \,\cdot\,V_1^{{2\,n}\over {\,n \ + \ 2\,}}\ \right]_{\,y}\\[0.2in]
& \le  & C_1 \int_{B_{\,\xi_{\,2}} (\,\rho_{\,\nu})} \left(\  [\   V_2(\,y\,)\ ]^{{4}\over {\,n\,-\,2\,}} \,+\,  C\,\cdot\,{1\over {{\bar\lambda}_{\ \flat}^2}}\,\cdot\, {\bar\lambda}_{\ \flat}^{4\,\gamma}  \right)^{{2n}\over {n\,+\,2}} \,\cdot\, [\   V_1(\,y\,)\ ]^{{2\,n}\over {\,n \ + \ 2\,}} \\[0.2in]
& \le & C_2\,\cdot\, \int_{B_{\,\xi_{\,2}} (\,\rho_{\,\nu})}   \left[\ \left(\    [\   V_2(\,y\,)\ ]^{{4}\over {\,n\,-\,2\,}}    \right)^{{2n}\over {n\,+\,2}}  \ \
 + \ \left(\     {1\over {{\bar\lambda}_{\ \flat}^2}}\,\cdot\,{\bar\lambda}_{\ \flat}^{4\,\gamma}  \right)^{{2n}\over {n\,+\,2}} \ \right]\,\cdot\,[\   V_1(\,y\,)\ ]^{{2\,n}\over {\,n \ + \ 2\,}} \ .\ \ \ \ \ \  \ \ \ \ \ \ \\[0.1in]
& \ & \ \ \ \ \ \ \  \ \ \ \ \ \ \ \ \ \ \ \ \  \ \ \ \ \  \uparrow \ \ {\mbox{first \ \ term}} \ \ \ \ \ \ \  \uparrow \ \ {\mbox{ second \ \ term \  \ }}
\end{eqnarray*}

Consider the second term in (\,A.7.34\,) first.\\[0.1in]
(\,A.7.35\,)
\begin{eqnarray*}
 & \ & \int_{B_{\,\xi_{\,2}} (\,\rho_{\,\nu})} \left( \ {1\over {{\bar\lambda}_{\ \flat}^2}}\,\cdot\, {\bar\lambda}_{\ \flat}^{4\,\gamma}  \right)^{{2n}\over {n\,+\,2}} \,\cdot\,[\   V_1(\,y\,)\ ]^{{2\,n}\over {\,n \ + \ 2\,}} \\[0.2in]
& \le & C_1 \cdot \int_{B_{\,\xi_{\,2}} (\,\rho_{\,\nu})} \left( \ {1\over {{\bar\lambda}_{\ \flat}^2}}\,\cdot\, {\bar\lambda}_{\ \flat}^{\,4\,\gamma}  \right)^{{2n}\over {n\,+\,2}} \,\cdot\,\left(\ {\lambda_{\,1}\over {\lambda_{\,1}^2 \ + \ \Vert\,\xi_{\,1}\ - \ \xi_{\,2}\,\Vert^2}}\  \right)^{ {{n\,-\,2}\over {n\,+\,2}} \cdot\, n } \\[0.2in]
& \le & C_2\,\cdot\int_{B_{\,\xi_{\,2}} (\,\rho_{\,\nu})}  \left[\ \left(\     {1\over {{\bar\lambda}_{\ \flat}^2}}\,\cdot\,{\bar\lambda}_{\ \flat}^{4\,\gamma}  \right)^{{2n}\over {n\,+\,2}} \ \right] \,\cdot\,\left[\ {1\over {\lambda_{\,1}}} \cdot {1\over {1 \ + \ \Vert\ \Xi_{\,\,1}\ - \ \Xi_{\,2}\ \Vert^2}} \ \right]^{ \, {{n\,-\,2}\over {n\,+\,2}} \cdot\, n } \\[0.1in]
& \ & \hspace*{2.5in} \left[\ {\mbox{recall \ \ that \ \ }} \Xi_{\,\,1} \ = \ {{\xi_{\,1}}\over { {\bar\lambda}_{\ \flat} }} \ \ {\mbox{and}} \ \ \Xi_{\,2} \ = \ {{\xi_{\,2}}\over { {\bar\lambda}_{\ \flat} }} \ \right]\\[0.1in]
& \ \le & C_3 \cdot  \left[\ \left(\     {1\over {{\bar\lambda}_{\ \flat}^2}}\,\cdot\,{\bar\lambda}_{\ \flat}^{4\,\gamma}  \right)^{{2n}\over {n\,+\,2}} \ \right]\,\cdot\,\left[\ {1\over {\sqrt{\,\lambda_{\,1}\,}}} \cdot {1\over {\Vert\ \Xi_{\,\,1}\ - \ \Xi_{\,2}\ \Vert}} \ \right]^{ \, {{n\,-\,2}\over {n\,+\,2}} \cdot\, 2\,n }  \cdot  {\mbox{Vol}}\,(\ B_{\xi_{\,2}} \,(\,\rho_{\,\nu})\,)\ . \\
\end{eqnarray*}

\newpage

As for  the first term\,,\, we start with\\[0.1in]
(\,A.7.36\,)
\begin{eqnarray*}
& \ & \left[ \ V_2^{ {{4}\over {\,n \ - \ 2\,}}\  \cdot\  {{2\,n}\over {\,n \ + \ 2\,}} } \ \cdot\   V_1^{{2\,n}\over {\,n \ + \ 2\,}}\ \right]_{\,y}\\[0.2in]
& = & \left(\,{{\lambda_{\,2}}\over {\ \lambda_{\,2}^2 \ + \ \Vert\, y \ - \ \xi_{\,2}\,\Vert^2\ }} \,   \right)^{\!\! {{4\,n}\over {\,n \ + \ 2\,}} } \ \cdot \left(\,{{\lambda_{\,1}}\over {\ \lambda_{\,1}^2 \ + \ \Vert\, y \ - \ \xi_{\,1}\,\Vert^2\ }} \,   \right)^{\!\! {{\,n\,(\,n \ - \ 2\,)\,}\over {\,n \ + \ 2\,}} }\\[0.2in]
& \le & C_o \cdot \left(\,{{1}\over {\ \lambda_{\,2} \ }} \,   \right)^{\!\! {{4\,n}\over {\,n \ + \ 2\,}} } \ \cdot \left(\,{{1}\over {\ \lambda_{\,1} \ }} \,   \right)^{\!\! {{\,n\,(\,n \ - \ 2\,)\,}\over {\,n \ + \ 2\,}} } \cdot  \left(\,{{1}\over {\ 1 \ + \ \Vert\, Y \ - \ \Xi_{\,2}\,\Vert^2\ }} \,   \right)^{\!\! {{4\,n}\over {\,n \ + \ 2\,}} } \ *\\[0.2in]
& \ & \ \ \ \ \ \ \ \ \ \ \ \ \ \ \ \ \ \ \ \ \ \ \ \  \ \ \ \ \ \ \ \ \ \ \ \ \ *\,\left(\,{{1}\over {\ 1 \ + \ \Vert\, Y \ - \ \Xi_{\,\,1}\,\Vert^2\ }} \,   \right)^{\!\! {{\,n\,(\,n \ - \ 2\,)\,}\over {\,n \ + \ 2\,}} } \ \ \ \ \ \ \ \ \ \ \left[ \ Y \ = \ {y\over { {\bar\lambda}_{\,\,\flat} }}\  \right] \\[0.2in]
& \le & C_1 \cdot \left(\,{{1}\over {\ \lambda_{\,2} \ }} \,   \right)^{\!\! {{n^2 \ + \ 2\,n}\over {\,n \ + \ 2\,}} } \ \cdot  \left(\,{{1}\over {\ 1 \ + \ \Vert\, Y \ - \ \Xi_{\,2}\,\Vert\ }} \,   \right)^{\!\! {{8\,n}\over {\,n \ + \ 2\,}} } \ \cdot \left(\,{{1}\over {\ 1 \ + \ \Vert\, Y \ - \ \Xi_{\,\,1}\,\Vert\ }} \,   \right)^{\!\! {{\,2\,n\,(\,n \ - \ 2\,)\,}\over {\,n \ + \ 2\,}} }\\[0.2in]
& \le & C_2 \cdot \left(\,{{1}\over {\ \lambda_{\,2} \ }} \,   \right)^{\!\!n } \ \cdot  \left(\,{{1}\over {\  \ \Vert\, \Xi_{\,\,1} \ - \ \Xi_{\,2}\,\Vert\ }} \,   \right)^{\!\!{{8\,n}\over {\,n \ + \ 2\,}}}\,*\\[0.2in]
& \ & \ \ \ \ \  \ \ \ \ \ \ \ \ \ \,* \left[\  \left(\,{{1}\over {\ 1 \ + \ \Vert\, Y \ - \ \Xi_{\,2}\,\Vert\ }} \,   \right)^{\!\! {{\,2\,n\,(\,n \ - \ 2\,)\,}\over {\,n \ + \ 2\,}} }  \ + \ \  \left(\,{{1}\over {\ 1 \ + \ \Vert\, Y \ - \ \Xi_{\,\,1}\,\Vert\ }} \,   \right)^{\!\! {{\,2\,n\,(\,n \ - \ 2\,)\,}\over {\,n \ + \ 2\,}} } \ \right]\ .
\end{eqnarray*}

We continue with\\[0.1in]
(\,A.7.37\,)
 \begin{eqnarray*}
 & \ & \\[-0.35in]
& \ & \int_{B_{\xi_{\,2}} (\,\rho_{\,\nu}) } \left(\,{{1}\over {\ 1 \ + \ \Vert\, Y \ - \ \Xi_{\,2}\,\Vert\ }} \,   \right)^{\!\! {{\,2\,n\,(\,n \ - \ 2\,)\,}\over {\,n \ + \ 2\,}} }  \, dy\\[0.2in]
  &  \le &  C_3\,\cdot\,\int_0^{ \rho_{\,\nu} }\left(\,{{1}\over {\ 1 \ + \ R}} \,   \right)^{\!\! {{\,2\,n\,(\,n \ - \ 2\,)\,}\over {\,n \ + \ 2\,}} }  \, \cdot\, r^{n\,-\,1}\ dr \ \ \ \ \ \ \ \ \left[\ R \ = \ \Vert\, Y \ - \ \Xi_{\,2}\,\Vert\ \ = \ {r\over {{\bar\lambda}_{\,\,\flat}}} \ \right]\\[0.2in]
& = & C_4 \cdot \lambda_{\,2}^n \,\cdot\, \int_0^{ \rho_{\,\nu}\over \lambda_{\,1} }\left( {1\over {1\ + \ R}}\  \right)^{\!\! {{\,2\,n\,(\,n \ - \ 2\,)\,}\over {\,n \ + \ 2\,}} } \cdot\, R^{\,n\,-\,1}\ dR \\[0.15in]
& \le &  C_5,\cdot\, {\bar{\lambda}}_{\,\flat}^{\,n\ -\ o_\lambda\,(\,1\,)} \ \ \ \ \ \ \  \  \left( \ {{\,2\,n\,(\,n \ - \ 2\,)\,}\over {\,n \ + \ 2\,}} \ \ge \ n \ \ \ \Longleftrightarrow \ \ n \ \ge \ 6 \ ; \ \ \ {\mbox{cf. \ \ Note A.7.31}}\ \right)  .
 \end{eqnarray*}
Moreover\,,\, we have
 \begin{eqnarray*}
(\,A.7.38\,) \ \ \  \ \ \ \  \ \ \  \ \ \ \ & \ & \int_{B_{\xi_{\,2}} (\,\rho_{\,\nu}) } \left(\,{{1}\over {\ 1 \ + \ \Vert\, Y \ - \ \Xi_{\,\,1}\,\Vert\ }} \,   \right)^{\!\! {{\,2\,n\,(\,n \ - \ 2\,)\,}\over {\,n \ + \ 2\,}} }  \, dy\\[0.2in]
  &  \le &  C_6\,\left(\,{{1}\over {\ 1 \ + \ \Vert\, \Xi_{\,2} \ - \ \Xi_{\,\,1}\,\Vert\ }} \,   \right)^{\!\! {{\,2\,n\,(\,n \ - \ 2\,)\,}\over {\,n \ + \ 2\,}} } \cdot\,\int_0^{ \rho_{\,\nu} }  \, \cdot\, r^{n\,-\,1}\ dr\ \ \  \ \ \ \ \ \ \  \ \ \ \ \ \ \  \ \ \ \ \\[0.2in]
  & \le & C_7 \cdot {\bar{\lambda}}_{\ \flat}^{\,  \gamma\, \cdot\, n \cdot {{\,2\,(\,n\,-\,2\,)\,}\over {\,n\,+\,2\,}}  } \cdot  {\bar{\lambda}}_{\ \flat}^{\, n\,\nu}\ = \ C_7\cdot  {\bar{\lambda}}_{\ \flat}^{\, n\,(\,\gamma \ + \ \nu\,)} \cdot {\bar{\lambda}}_{\ \flat}^{\,  \gamma\, \cdot\, n \cdot {{\,n\,(\,n\,-\,6\,)\,}\over {\,n\,+\,2\,}}  }  \\[0.2in]
  & \le & C_8 \cdot  {\bar{\lambda}}_{\ \flat}^{\, n\,}  \ \ \ \ \ \  \ \ \ \ \ [\ (\,\gamma \ + \ \nu\,) \ > \ 1\ ]
 \end{eqnarray*}
Combining (\,A.7.35\,)\,--\,(\,A.7.38\,)\,,\, we obtain
 \begin{eqnarray*}
 (\,A.7.39\,) \ \ \  \ \ \ \  \ \ \ & \ & \left[\ \int_{B_{\xi_{\,2}} \,(\,\rho_{\,\nu})\,}  \left(\   (W_{\,\,\flat} )^{{4}\over {\,n\,-\,2\,}} \ - \  V_1^{{4}\over {\,n\,-\,2\,}} \,  \right)^{{2n}\over {n\,+\,2}} \,\cdot\,V_1^{{2n}\over {n\,+\,2}} \ \right]^{{n\,+\,2}\over {2\,n}} \\[0.2in]
 & \le &  C_9 \cdot  \left(\     {1\over {{\bar\lambda}_{\ \flat}^2}}\,\cdot\,{\bar\lambda}_{\ \flat}^{4\,\gamma}  \right) \,\cdot\, \left[\ \rho_{\,\nu}^n \ \right]^{{n\,+\,2}\over {2\,n}} \cdot  \ \left[\ {1\over {\lambda_{\,1}}} \cdot {1\over {\Vert\ \Xi_{\,\,1}\ - \ \Xi_{\,2}\ \Vert}} \ \right]^{ \, n\,-\,2 }\ \ \  \ \ \ \  \ \ \  \ \ \ \ \ \ \  \ \ \ \  \ \ \  \ \ \ \ \ \ \  \ \ \ \  \ \ \  \ \ \ \  \\[0.2in]
  & \ & \ \ \ \ \ + \ C_{\,10} \cdot \lambda_{\,2}^{-\ o_\lambda\,(\,1\,)} \cdot  \left[\   {1\over {\Vert\ \Xi_{\,\,1}\ - \ \Xi_{\,2}\ \Vert}} \ \right]^{ \,4}
 \end{eqnarray*}
Estimating in a similar fashion, then summing up the results, we obtain \\[0.1in]
(\,A.7.40\,)
\begin{eqnarray*}
& \ & \\[-0.35in]
& \ &   \left\{ \ \left[\   \int_{B_{\xi_{\,2}} (\,\rho_{\,\nu})}   \ + \  \cdot \cdot \cdot  \ + \  \int_{B_{\xi_{\, \flat}} (\,\rho_{\,\nu})}   \ \right] \  \left(\   (W_{\,\,\flat} )^{{4}\over {\,n\,-\,2\,}} \ - \  V_1^{{4}\over {\,n\,-\,2\,}} \,  \right)^{{2n}\over {n\,+\,2}} \,\cdot\,V_1^{{2n}\over {n\,+\,2}} \ \right\}^{{n\,+\,2}\over {2\,n}}\\[0.2in]
& \le &  C_{11} \cdot  \left(\     {1\over {{\bar\lambda}_{\ \flat}^2}}\,\cdot\,{\bar\lambda}_{\ \flat}^{4\,\gamma}  \right) \,\cdot\, {\bar\lambda}_{\ \flat}^{{{n\,+\,2}\over {2}}\, \cdot \ \nu} \cdot {\bar\lambda}_{\ \flat}^{- {{{n\,-\,2}\over {2}}} } \ \left[\ \,\sum_{j\,=\,2}^\flat\ {1\over {\Vert\ \Xi_{\,\,1}\ - \ \Xi_j\ \Vert^{ \, n\,-\,2 } }} \ \right] \ \ \ \cdot \cdot \cdot \ (\,{\bf *}\,)_{\,(\,A.7.40\,) } \\[0.1in]
  & \ & \ \ \ \ \ \ \ \ \ \  \  \ \ \ \ \ \ \ \   \left[\  {\mbox{see \ \ Note \ A.7.41 \ \ below}} \ ; \ \ \ \  \emph{}\uparrow \ [\ \cdot \cdot \cdot \ ] \ = \ O\left( {\bar\lambda}_{\ \flat}^{\ \gamma \cdot  {{\,n\,-\,2\,}\over 2}} \  \right) \  \right] \\[0.1in]
   & \ & \ \ \ \ \ + \ \ C_{\,12} \cdot {\bar\lambda}_{\ \flat}^{-\ o_\lambda\,(\,1\,)} \cdot  \left[\ \sum_{j\,=\,2}^\flat\   {1\over {\Vert\ \Xi_{\,\,1}\ - \ \Xi_j\ \Vert^{\,4}}} \ \right]\\[0.2in]
    & \le & C_{13} \cdot {\bar\lambda}_{\ \flat}^{\,4\,\gamma \ -\ o_\lambda\,(\,1\,)} \ \ \ \ \ \ \ \ \   \ (\,n \ \ge \ 6\,)\ .
 \end{eqnarray*}

\newpage

{\it Note.} \ \,A.7.41. \ \
Consider the indices in $(\,{\bf *}\,)_{\,(\,A.7.40\,) }\,:$
\begin{eqnarray*}
& \  & -\,2\ \ + \ 4\,\gamma \ + \ { {{\,n\,+\,2\,}\over 2} \cdot \,\nu} \ - \  {{\,n\,-\,2\,}\over 2}   \ + \ \left(   \,n\,-\,2\,  \right) \ \gamma \\[0.2in]
& = & 4\,\gamma \ + \ \left(  {{n\,-\,2}\over 2} \right) \cdot 2\,\gamma \ + \ { {{\,n\,+\,2\,}\over 2} \cdot \,\nu} \ - \  {{n\,-\,2}\over 2}  \ - \ 2\\[0.2in]
& = & 4\,\gamma \ + \  \left(  {{n\,-\,6}\over 2} \right) \cdot \gamma  \ + \  \left(  {{n\,+\,2}\over 2} \right) \cdot \gamma \ + \ { {{\,n\,+\,2\,}\over 2} \cdot \,\nu} \ - \  {{n\,+\,2}\over 2}\\[0.2in]
& \ge &  4\,\gamma \ + \   \left(  {{n\,+\,2}\over 2} \right) \cdot  [\ (\,\gamma \ + \ \nu\,) \ - \ 1 \ ]  \ \  \ \ \ \ \  \ \ \ \ \ \  \ \ \ \  \ ( \ n \ \ge \ 6\ ) \\[0.2in]
& > &4\,\gamma  \ \ \ \ \  \ \ \ \ \ \  \  \ \ \ \ \  \ \ \ \ \ \  \ \ \ \ \  \ \ \ \ \ \  \ \ \ \ \ \ \ \  \ \ \ \ \ \  \   \  \ \ \ \ \  \ \ \ \ \ \  \ [\ (\,\gamma \ + \ \nu \,) \ >  \ 1\ ]\ .
\end{eqnarray*}

 \vspace*{0.2in}

 {\it Note.}  A.7.42.\ \ This part of the argument works in a similar way for  $\,\, \left( \, \lambda_{\,1}\,\cdot\,{ { \partial\, V_1} \over {\partial\, \xi_{1_{\,j}}}}\  \right)\,$.

  \vspace*{0.3in}

{\bf \S\,A\,7.\,h\,.} \ \
{\bf Inside}  $\ \R^n \,\setminus \, \left\{\  \cup\ B_{\,\xi_{\,l}} (\,{\rho_{\,\nu}}) \, \right\}\,.$ \,\ \ Refer to (\,A.7.28\,)\,.
It follows from (\,A.7.26\,) that \\[0.1in]
(\,A.7.42\,)
\begin{eqnarray*}
& \ &   \left[\    \int_{\R^n \,\setminus \, \left(\, \cup\ B_{\xi_{\,l}} \,(\,{\rho_{\,\nu}}) \right) }     \left(\   (W_{\,\,\flat} )^{{4}\over {\,n\,-\,2\,}} \ - \  V_1^{{4}\over {\,n\,-\,2\,}} \,  \right)^{{2n}\over {n\,+\,2}} \,\cdot\,V_1^{{2n}\over {n\,+\,2}} \ \right]^{{n\,+\,2}\over {2\,n}} \\[0.1in]
& \le & \left[\ \int_{\R^n \,\setminus \, \left(\, \cup\ B_{\xi_{\,l}} (\,{\rho_{\,\nu}}) \, \right)}  \left(\    V_2^{{4}\over {\,n\,-\,2\,}} \,+\,V_3^{{4}\over {\,n\,-\,2\,}}  \ + \ \cdot\,\cdot\,\cdot\,\ + \ V_{\,\flat}^{{4}\over {\,n\,-\,2\,}}  \right)^{{2n}\over {n\,+\,2}}*\,V_1^{{2n}\over {n\,+\,2}} \ \right]^{{n\,+\,2}\over {2\,n}} \ .
\end{eqnarray*}
For what is inside\,,\, we continue with $\displaystyle{\left( {\mbox{note \ \ that \ \ }} {n\over {n\,+\,2}} \ <  \ 1\ , \ \ {\mbox{cf.}} \ \ {\bf {\S\,A\,4\,.\,g}}\ \right)}$\\[0.1in]
(\,A.7.43\,)
\begin{eqnarray*}
& \ & \left(\    V_2^{{4}\over {\,n\,-\,2\,}} \,+\,V_3^{{4}\over {\,n\,-\,2\,}}  \ + \ \cdot\,\cdot\,\cdot\,\ + \ V_{\,\flat}^{{4}\over {\,n\,-\,2\,}}  \right)^{{2n}\over {n\,+\,2}} \\[0.2in]
 & = &   \left(\    V_2^{{4}\over {\,n\,-\,2\,}} \,+\,V_3^{{4}\over {\,n\,-\,2\,}}  \ + \ \cdot\,\cdot\,\cdot\,\ + \ V_{\,\flat}^{{4}\over {\,n\,-\,2\,}}  \right)^{{ n}\over {n\,+\,2}}\,\cdot\,\left(\    V_2^{{4}\over {\,n\,-\,2\,}} \,+\,V_3^{{4}\over {\,n\,-\,2\,}}  \ + \ \cdot\,\cdot\,\cdot\,\ + \ V_{\,\flat}^{{4}\over {\,n\,-\,2\,}}  \right)^{{ n}\over {n\,+\,2}}  \\[0.2in]
& \le &    \left[\  V_2^{{4\,n}\over {(\,n\,+\,2\,) (\,n\,-\,2\,)   }}   \ + \,\cdot\,\cdot\,\cdot\, + \ V_{\,\flat}^{{4\,n}\over {(\,n\,+\,2\,) (\,n\,-\,2\,)   }}  \ \right]\,\cdot\,\left[\  V_2^{{4\,n}\over {(\,n\,+\,2\,) (\,n\,-\,2\,)   }}   \ + \,\cdot\,\cdot\,\cdot\, + \ V_{\,\flat}^{{4\,n}\over {(\,n\,+\,2\,) (\,n\,-\,2\,)   }}  \ \right]\\[0.2in]
& = & \left[\  V_2^{{8\,n}\over {(\,n\,+\,2\,) (\,n\,-\,2\,)   }}   \ + \,\cdot\,\cdot\,\cdot\, + \ V_{\,\flat}^{{8\,n}\over {(\,n\,+\,2\,) (\,n\,-\,2\,)   }}  \ \right] \\[0.2in]
& \  &   \ + \  V_2^{{4\,n}\over {(\,n\,+\,2\,) (\,n\,-\,2\,)   }}\,\cdot\, \left[\  V_3^{{4\,n}\over {(\,n\,+\,2\,) (\,n\,-\,2\,)   }}   \ + \,\cdot\,\cdot\,\cdot\, + \ V_{\,\flat}^{{4\,n}\over {(\,n\,+\,2\,) (\,n\,-\,2\,)   }}  \ \right]\\[0.1in]
& \  &    \ \cdot\,\ \\[0.1in]
& \  &   \ \cdot\,\ \\[0.1in]
& \  &   \ \ \ \ + \  V_{\,\flat}^{{4\,n}\over {(\,n\,+\,2\,) (\,n\,-\,2\,)   }}\,\cdot\, \left[\  V_2^{{4\,n}\over {(\,n\,+\,2\,) (\,n\,-\,2\,)   }}   \ + \,\cdot\,\cdot\,\cdot\, + \ V_{\flat\,-\,1}^{{4\,n}\over {(\,n\,+\,2\,) (\,n\,-\,2\,)   }}  \ \right]  \ . \\
\end{eqnarray*}
For the coupled terms, we proceed with \\[0.1in]
(\,A.7.44\,)
\begin{eqnarray*}
& \ & V_1^{{2\,n}\over { \,n\,+\,2\,      }} \,\cdot\,  V_2^{{8\,n}\over {(\,n\,+\,2\,) (\,n\,-\,2\,)   }}\\[0.2in]
& = & \left(\,{{\lambda_{\,2}}\over {\ \lambda_{\,2}^2 \ + \ \Vert\, y \ - \ \xi_{\,2}\,\Vert^2\ }} \,   \right)^{\!\! {{n\,(\,n \,-\,2\,)}\over {\,n \ + \ 2\,}} } \ \cdot \left(\,{{\lambda_{\,1}}\over {\ \lambda_{\,1}^2 \ + \ \Vert\, y \ - \ \xi_{\,1}\,\Vert^2\ }} \,   \right)^{\!\! {{4\,n\,}\over {\,n \ + \ 2\,}} }\\[0.2in]
& = &  \left(\,{{1}\over {\ \lambda_{\,1} \ }} \,   \right)^{\!\! {{\,n\,(\,n \ - \ 2\,)\,}\over {\,n \ + \ 2\,}} } \cdot \left(\,{{1}\over {\ \lambda_{\,2} \ }} \,   \right)^{\!\! {{4\,n\,}\over {\,n \ + \ 2\,}} }  *\\[0.1in]
& \ & \ \ \ \ \ \  \ \ \ \,* \,   \left(\,{{1}\over {\ 1 \ + \ \Vert\, Y \ - \ \Xi_{\,\,1}\,\Vert^2\ }} \,   \right)^{\!\! {{\,n\,(\,n \ - \ 2\,)\,}\over {\,n \ + \ 2\,}} } \cdot \left(\,{{1}\over {\ 1 \ + \ \Vert\, Y \ - \ \Xi_{\,2}\,\Vert^2\ }} \,   \right)^{\!\! {{4\,n}\over {\,n \ + \ 2\,}} } \\[0.2in]
& \le & C_1\,\cdot\,{1\over {{\bar\lambda}_{\ \flat}^n}}\,\cdot\,\left(\ {1\over {1 \ + \ \Vert\,Y \ - \ \Xi_{\,\,1}\,\Vert }}\  \right)^{\!\! {{2\,n\,(\,n\,-\,2)}\over {n\,+\,2}} } \ \cdot\, \  \left(\ {1\over {1 \ + \ \Vert\,Y \ - \ \Xi_{\,2}\,\Vert }}\  \right)^{\!\! {{8\,n}\over {n\,+\,2}} }\\[0.2in]
& \le & C_2\,\cdot\,  \left(\ {1\over {  \Vert\,\Xi_{\,\,1} \ - \ \Xi_{\,2}\,\Vert }} \ \right)^{\!\! {{8\,n}\over {n\,+\,2}}\,- \ \epsilon } *\\[0.1in]
& \ &  \ \ \   *\,{1\over {{\bar\lambda}_{\ \flat}^n}}\,\cdot\,\left[\ \left(\ {1\over {1 \ + \ \Vert\,Y \ - \ \Xi_{\,\,1}\,\Vert }} \ \right)^{\!\! {{2\,n\,(\,n\,-\,2)}\over {n\,+\,2}} \,+\,\epsilon } \ + \ \left(\ {1\over {1 \ + \ \Vert\,Y \ - \ \Xi_{\,2}\,\Vert }}\ \right)^{\!\! {{2\,n\,(\,n\,-\,2)}\over {n\,+\,2}} \,+\,\epsilon } \ \right]\\[0.2in]
& \ & \ \ \ \ \ \ \ \  \ \ \ \ \  \ \ \ \  \ \ \ \ \  \ \ \ \  \ \ \ \   \left\{ \ {\mbox{observe \ \ that}} \ \ n \ \ge \ 6 \ \ \Longrightarrow \  \ {{2\,n\,(\,n\,-\,2)}\over {n\,+\,2}} \ \ge \ {{8\,n }\over {n\,+\,2}}  \ \right\}\ .
\end{eqnarray*}
Here $\,\epsilon  \ >  \ 0\,$ is a fixed small number.
We come to
\begin{eqnarray*}
(\,A.7.45\,)     & \ &{1\over {{\bar\lambda}_{\ \flat}^n}}\,\cdot\,\int_{\,\R^n \setminus \,  \,B_{\xi_{\,2} } (\,{\rho_{\,\nu}})}    \left(\ {{1}\over {1 \ + \ \Vert\,Y\ - \ \Xi_{\,2}\,\Vert }}\ \right)^{\!\! {{2\,n\,(\,n\,-\,2)}\over {n\,+\,2}} \ + \ \epsilon }\ d\,V_y \ \ \ \ \   \left[\ Y \ = \ {y\over{{\bar\lambda}_{\ \flat}}} \ \right] \ \ \ \  \ \ \\[0.2in]
& \le & C_1\,\cdot\,  \int_{ {{\rho_{\,\nu}}\over { {\bar\lambda}_{\, \flat}}} }^\infty\  {{R^{n\,-\,1}\ }\over {\  R^{ {{2\,n\,(\,n\,-\,2)}\over {n\,+\,2}} \ + \ \epsilon  } \ }} \ dR\\[0.2in]& \le & C_2\,\cdot\, {1\over {\  R^{\,  (\,n\,-\,6\,) \,\cdot \,{{n }\over { n\,+\,2 }} \ + \ \epsilon  } \ }} \ \ \bigg\vert_{{{\rho_{\,\nu}}\over { {\bar\lambda}_{\, \flat}}}}^{\,\infty}\\[0.2in]& \le & C_3  \,\cdot\,  {\bar\lambda}_{\ \flat}^{ \,\left( (\,n\,-\,6\,)\,\,\cdot\, \, {{n }\over { n\,+\,2 }} \ + \ \epsilon  \, \right) \ \cdot\,\ (\,1\ - \ \nu\,) } \ \ \ \ \ [ \ {\mbox{see   \  \ also \ \ }} (\,A.7.31\,) \ ]
\end{eqnarray*}
Thus
\begin{eqnarray*}
(\,A.7.46\,) \ \   & \ &  \sum_{j\,=\,2}^\flat \left[\    \int_{\R^n \,\setminus \  \left\{\  \cup\ B_{\xi_{\,l}} \,(\,{\rho_{\,\nu}}) \ \right\} }     V_j^{{8\,n}\over {(\,n\,+\,2\,) (\,n\,-\,2\,) }}  \,\cdot\,V_1^{{2n}\over {n\,+\,2}} \ \right]^{{n\,+\,2}\over {2\,n}}  \\[0.2in]
& = & O\left( {\bar\lambda}_{\ \flat}^{ \,4\,\gamma \,-\ \epsilon'} \right) \,\cdot\, O\left( \  {\bar\lambda}_{\ \flat}^{ \left(\ (\ n\,-\,6\,)\,\cdot\,{{1 }\over { 2 }} \ + \ \epsilon'   \right) \ \cdot\,\ (\,1\ - \, \nu\,) }\ \right)  \ \ \ \ \ \  \ \ \ \ \  \left[\ \epsilon' \ = \ {{n\,+\,2}\over {2\,n}}\,\cdot\,\epsilon\ \right]\ .\\
\end{eqnarray*}
\hspace*{0.5in}Next, for the triple terms, we move on in a similar way as in with (\,A.7.44\,)\,.\, Cf. {\bf \S\,A\,4\,.\,j}\,,\, in particular\,,\, (\,A.4.56\,)\,. \\[0.1in]
(\,A.7.47\,)
\begin{eqnarray*}
& \ & \left[ \  V_1^{{2\,n}\over { \,n\,+\,2\,      }} \,\cdot\, V_2^{{4\,n}\over {(\,n\,+\,2\,) (\,n\,-\,2\,)   }}\,\cdot\, V_3^{{4\,n}\over {(\,n\,+\,2\,)\, (\,n\,-\,2\,)   }}\ \right]_{\,y}\\[0.2in]
& \le &  C \,\cdot\,{1\over {{\bar\lambda}_{\ \flat}^n}}\,\cdot\, \left(\ {{1}\over {1 \ + \ \Vert\,Y\ - \ \Xi_{\,\,1}\,\Vert }}\ \right)^{\!\! {{2\,n\,(\,n\,-\,2)}\over {n\,+\,2  }} } \,\cdot\,\left(\ {{1}\over {1 \ + \ \Vert\,Y\ - \ \Xi_{\,2}\,\Vert  }}\ \right)^{\!\!{{4\,n}\over {n\,+\,2}} }\ * \\[0.2in]
& \ &  \ \ \ \ \   * \  \left(\ {{1}\over {1 \ + \ \Vert\,Y\ - \ \Xi_{\,3}\,\Vert  }}\ \right)^{\!\!{{4\,n}\over {n\,+\,2}} } \ \ \  ( \ {\mbox{see \ \ note \ \ at \ \ the  \ \ end \ \ of \ \ this \ \ group }}\ )
\end{eqnarray*}

\newpage

\begin{eqnarray*}
& \le  &  C_2\,\cdot\,{1\over {{\bar\lambda}_{\ \flat}^n}}\,\cdot\,{1\over {\,\Vert\, \Xi_{\,\,1} \ - \ \Xi_{\,2}\,\Vert^{{{4\,n}\over {n\,+\,2}}\,-\,\epsilon} }}\,\cdot\,\left[\ \left(\ {{1}\over {1 \ + \ \Vert\,Y\ - \ \Xi_{\,\,1}\,\Vert }}\ \right)^{\!\!{{2\,n\,(\,n\,-\,2)}\over {n\,+\,2  }}\,+\,\epsilon}  \ + \ \right.\\[0.2in]
& \ & \left.\ \ \ \ \ \ \ \ \ \ \ \ \ \ \ \ \   + \  \left(\ {{1}\over {1 \ + \ \Vert\,Y\ - \ \Xi_{\,2}\,\Vert }}\ \right)^{\!\!{{2\,n\,(\,n\,-\,2)}\over {n\,+\,2  }})\,+\,\epsilon} \ \right] \ \,\times \  \left(\ {{1}\over {1 \ + \ \Vert\,Y\ - \ \Xi_{\,3}\,\Vert  }}\ \right)^{\!\!{{4\,n}\over {n\,+\,2}}}\\[0.2in]
& \le &  C_{\,3}\,\cdot\,{1\over {{\bar\lambda}_{\ \flat}^n}}\cdot  {1\over {\,\Vert\, \Xi_{\,\,1} \ - \ \Xi_{\,2}\,\Vert^{{{4\,n}\over {n\,+\,2}}\,-\,\epsilon} }}  \,\cdot\,{1\over {\,\Vert\, \Xi_{\,\,1} \ - \ \Xi_{\,3}\,\Vert^{{{4\,n}\over {n\,+\,2}}} }}\,\cdot\,\left[\ \left(\ {{1}\over {1 \ + \ \Vert\,Y\ - \ \Xi_{\,\,1}\,\Vert }}\ \right)^{\!\!{{2\,n\,(\,n\,-\,2)}\over {n\,+\,2  }}\ +\  \epsilon}  + \right.\\[0.2in]
& \ & \hspace*{2in} \ \ \ \ \ \ \ \ \  +  \  \left. \left(\ {{1}\over {1 \ + \ \Vert\,Y\ - \ \Xi_{\,2}\,\Vert }}\ \right)^{\!\!{{2\,n\,(\,n\,-\,2)}\over {n\,+\,2  }})\ +\ \epsilon} \ \right]   \\[0.2in]
& \ & \hspace*{3in}+  \  \left. \left(\ {{1}\over {1 \ + \ \Vert\,Y\ - \ \Xi_{\,3}\,\Vert }}\ \right)^{\!\!{{2\,n\,(\,n\,-\,2)}\over {n\,+\,2  }})\,+\,\epsilon} \ \right] \\[0.2in]
& \ & \!\!\!\!\!\!\!\!\left\{ \ {\mbox{note \ \ that}} \ \ \ \ {{\,n\,-\,2\,}\over 2} \cdot \left[\  {{2\,n}\over { \,n\,+\,2\,      }}  \ + \  {{4\,n}\over {(\,n\,+\,2\,) (\,n\,-\,2\,)   }} \ + \  {{4\,n}\over {(\,n\,+\,2\,) (\,n\,-\,2\,)   }} \ \right] \ = \ n \ \right\}\ .\\
\end{eqnarray*}
Following (\,A.7.44\,)\,,\, using the integral estimate as in (\,A.7.45\,)\,,\, after summing up and taking the power of $\ \displaystyle{\left( \ {{n\,+\,2}\over {2\,n}} \ \right)}\ $ into account [\,with the help of (\,A.7.26\,)\ ]\,,\, the contribution of the  terms added up to [ \ cf. (\,A.4.57\,)\ ]
\begin{eqnarray*}
 \ \   & \ &  \sum_{j\,,\ k \ \ge \ 2}^{j \,\not=\,k} \left[\    \int_{\R^n \,\setminus \  \left\{\  \cup\ B_{\xi_{\,l}} \,(\,{\rho_{\,\nu}}) \ \right\} }   V_j^{{4\,n}\over {(\,n\,+\,2\,) (\,n\,-\,2\,) }}  \,\cdot\,    V_k^{{4\,n}\over {(\,n\,+\,2\,) (\,n\,-\,2\,) }}  \,\cdot\,V_1^{{2n}\over {n\,+\,2}} \ \right]^{{n\,+\,2}\over {2\,n}}  \\[0.2in]
& = & O\left(\ {\bar\lambda}_{\ \flat}^{ \,4\,\gamma \,-\ {{n\,+\,2}\over 2} \cdot\,\epsilon'}\ \right) \,\cdot\, O\left( \  {\bar\lambda}_{\ \flat}^{ \left(\, (\ n\,-\,6\,)\,\cdot\,{{1 }\over { 2 }} \ + \ \epsilon'  \,  \right) \ \cdot\,\ (\,1\ - \, \nu\,) }\ \right)  \ \ \ \ \ \    \left[\ \epsilon' \ = \ {{n\,+\,2}\over {2\,n}}\,\cdot\,\epsilon\ \right] \ .
\end{eqnarray*}

\newpage

As a conclusion,\\[0.1in]
(\,A.7.48\,)
$$
\int_{\R^n} \bigg\vert \ (\,W_{\,\,\flat}\, )^{{4}\over {\,n\,-\,2\,}}\,\cdot\,\left(\,\lambda_{\,1}\cdot  { { \partial \,V_1 } \over {\partial\, \lambda_{\,1}}}\ \right) \ - \ \  V_1^{{4}\over {\,n\,-\,2\,}}\,\cdot\,\left(\,\lambda_{\,1}\cdot  { { \partial \,V_1 } \over {\partial\, \lambda_{\,1}}}\ \right) \,\bigg\vert \,\cdot\, |\,\phi\,| \,dy \ = \  O\left( \, {\bar\lambda}_{\ \flat}^{ 4\,\gamma\ - \ \tilde\varepsilon }\ \right)\,\cdot\, O\left(\,{\bar\lambda}_{\,\,\flat}^{{{\varpi}_{\,(\,A.5.38\,)} }} \, \right)\ .
$$
Here $\,\tilde\varepsilon\,$ is a fixed small positive number.

 \vspace*{0.1in}

 {\it Note.}  A.7.49.\ \ This part of the argument works in a similar way for  $\,\, \left( \ \lambda_{\,1}\,\cdot\,{ { \partial\, V_1} \over {\partial\, \xi_{1_{\,j}}}}\  \right)\,$.

\vspace*{0.3in}

{\bf \S\,A\,7.\,i\,.} \ \ {\bf Estimate of}  \ $\,{\bf II}_{\ (\,A.7.23\,)}\,$\,.\, \\[0.1in]
Consider the remaining term\,:\\[0.1in]
(\,A.7.50\,)
 \begin{eqnarray*}
& \ & \Bigg\vert \ \int_{\R^n}\,[\ (\,{\tilde c}_n\!\cdot K\,)\,(\,y\,)\ - \ n\,(\,n\,-\,2)\ ]\,\cdot\,\left[\    V_1^{{4}\over {\,n\,-\,2\,}}\,\cdot\,\left(\,\lambda_{\,1}\cdot  { { \partial \,V_1 } \over {\partial\, \lambda_{\,1}}}\ \right)\, \right]_{\,y}\,\cdot\,\phi\,(\,y\,) \ dy  \\[0.2in]
& \le &C_1\cdot \int_{\R^n}\,|\,(\,{\tilde c}_n\!\cdot K\,)\,(\,y\,)\ - \ n\,(\,n\,-\,2)\,|\,\cdot\,     [\ V_1\,(\,y\,)\ ]^{{\,n\,+\,2\,}\over {n\,-\,2}} \,\cdot\,|\, \phi\,(\,y\,)\,|\,dy\\[0.15in]
& \le & C_2\cdot \int_{\R^n}\,|\,(\,{\tilde c}_n\!\cdot {\bf K}\,(\,Y\,)\,)\ - \ n\,(\,n\,-\,2)\,|\,\cdot\,     {1\over { {\bar\lambda}_{\,\,\flat}^{{n\,+\,2}\over 2} }} \,*\,\\[0.2in]
& \ & \hspace*{1.7in} *\  \left(\ {1\over {1\ + \ \Vert\,Y\ - \ \Xi_{\,\,1}\Vert}}\  \right)^{\!\! n\,+\,2 }     \,\cdot\,{1\over { {\bar\lambda}_{\,\,\flat}^{{\,n\,-\,2\,}\over 2} }} \,\cdot\,|\, \Phi\,(\,Y\,)\,|\cdot {\bar\lambda}_{\,\,\flat}^n\, \,dY\\[0.1in]
& \ & \ \ \ \ \ \ \ \ \ \ \ \ \  \left[ \ Y \ = \ {y\over {{\bar\lambda}_{\,\,\flat}  }} \ , \ \ \ \ \ {\bf K}\,(\,Y\,) \ = \ K\,(\,{\bar\lambda}_{\,\,\flat} \cdot\, Y\,)\ , \ \ \ \Phi\,(\,Y) \  \ = \  {\bar\lambda}_{\,\,\flat}^{\,{{n\,-\,2 }\over 2}}\  \cdot \phi\,(\,{\bar\lambda}_{\,\,\flat}\cdot Y\,) \  \right] \\[0.1in]
& \le & C_2\cdot  \int_{\R^n}\,|\,(\,{\tilde c}_n\!\cdot K\,(\,Y\,)\,)\ - \ n\,(\,n\,-\,2)\,| \,\cdot\, \left(\ {1\over {1\ + \ \Vert\,Y\ - \ \Xi_{\,\,1}\Vert}}\  \right)^{\!\! n\,+\,2 }       \,\cdot\,|\, \Phi\,(\,Y\,)\,|\,\cdot\, dY\\[0.15in]
& \le & C_3\cdot  \Vert\,\Phi \,\Vert_*\,\cdot\,\int_{\R^n}\,|\,(\,{\tilde c}_n\!\cdot K\,(\,Y\,)\,)\ - \ n\,(\,n\,-\,2)\,|\cdot  \left(\ {1\over {1\ + \ \Vert\,Y\ - \ \Xi_{\,\,1}\Vert}} \ \right)^{\!\! n\,+\,2 }        \cdot\\[0.15in]
& \ & \ \ \  \,\cdot\,\left[\ \left(\ {1\over {1\ + \ \Vert\,Y\ - \ \Xi_{\,\,1}\Vert}} \  \right)^{\!\!{{\,n\,-\,2\,}\over 2}\,+\,\tau_{\,>1}}\ +\,\left(\ {1\over {1\ + \ \Vert\,Y\ - \ \Xi_{\,2}\Vert}} \  \right)^{\!\!{{\,n\,-\,2\,}\over 2}\,+\,\tau_{\,>1}}\ +\,\cdot\,\cdot\,\cdot\,\ \right.\\[0.1in]
& \
   & \hspace*{3in}\ \ \ \left. \ + \  \left(\ {1\over {1\ + \ \Vert\,Y\ - \ \Xi_{\,\flat}\Vert}} \ \right)^{\!\!{{\,n\,-\,2\,}\over 2} \,+\,\tau_{\,>1}} \ \  \right]\,\,dY\ .
\end{eqnarray*}

\newpage

We continue with the estimate as in  {\bf \S\,A\,3.\,l}\ , (\,A.3.40\,)\,:\\[0.1in]
(\,A.7.51\,)
\begin{eqnarray*}
& \  &  \int_{\R^n}\,|\,(\,{\tilde c}_n\!\cdot K\,(\,Y\,)\,)\ - \ n\,(\,n\,-\,2)\,|\cdot  \left(\ {1\over {1\ + \ \Vert\,Y\ - \ \Xi_{\,\,1}\Vert}} \  \right)^{\!\! (\,n\,+\,2)\ +\,{{\,n\,-\,2\,}\over 2}\,+\,\tau_{\,>1} } \ \ dY\\[0.2in]
& \le & O\,(\,{\bar\lambda}_{\ \flat}^{\,\ell}\,)   \ + \ O \left(\ \ {\bar\lambda}_{\ \flat}^{\,{{n \ +\ 2}\over 2} \ +\,\tau_{\,>1}\,  } \ \right) \ .\\
\end{eqnarray*}

As for the mixed terms\,,\, we estimate via
\\[0.1in]
(\,A.7.52\,)
\begin{eqnarray*}
& \ &   \left(\ {1\over {1\ + \ \Vert\,Y\ - \ \Xi_{\,\,1}\,\Vert}}\  \right)^{\!\! n\,+\,2 }       \,\cdot\,  \left(\ {1\over {1\ + \ \Vert\,Y\ - \ \Xi_{\,2}\,\Vert}}\  \right)^{\!\!{{\,n\,-\,2\,}\over 2}\,+\,\tau_{\,>1}}   \\[0.2in]
& \le &   {1\over { \Vert\,\Xi_{\,\,1} \ - \ \Xi_{\,2}\,\Vert^{{{n\,+\,2}\over 2} \,+\,\tau_{\,>1}\,-\,\epsilon}  }} \,\cdot\, \left[\  \left(\ {1\over {1\ + \ \Vert\,Y\ - \ \Xi_{\,\,1}\,\Vert}}\  \right)^{\!\! n\,+\,\epsilon } \ + \  \left(\ {1\over {1\ + \ \Vert\,Y\ - \ \Xi_{\,2}\,\Vert}}\  \right)^{\!\! n\,+\,\epsilon } \ \right]\ .
\end{eqnarray*}
Thus
 \begin{eqnarray*}
& \ &   \int_{\R^n}\,|\,(\,{\tilde c}_n\!\cdot {\bf K}\,(\,Y\,)\,)\ - \ n\,(\,n\,-\,2)\,|\cdot  \left(\ {1\over {1\ + \ \Vert\,Y\ - \ \Xi_{\,\,1}\Vert}}\  \right)^{\!\! n\,+\,2 }        \cdot\\[0.15in]
& \ & \ \ \ \  \ \ \ \ \,\cdot\,\left[\ \left(\ {1\over {1\ + \ \Vert\,Y\ - \ \Xi_{\,2}\Vert}}\  \right)^{\!\!{{\,n\,-\,2\,}\over 2}\,+\,\tau_{\,>1}} \ +\,\cdot\,\cdot\,\cdot\,\   \ + \  \left(\ {1\over {1\ + \ \Vert\,Y\ - \ \Xi_{\,\flat}\Vert}}\  \right)^{\!\!{{\,n\,-\,2\,}\over 2} \,+\,\tau_{\,>1}}\  \right]\,\,dY\\[0.2in]
& \le & C_1 \cdot\sum_{l\,=\,2}^\flat\   {1\over { \Vert\,\Xi_{\,\,1} \ - \ \Xi_{\,l}\,\Vert^{{{n\,+\,2}\over 2} \,+\,\tau_{\,>1}\,-\,\epsilon}  }} \\[0.2in]
& \le &  C_2\,\cdot\,  {\bar\lambda}_{\ \flat}^{ \ \left[\  {{n\,+\,2}\over 2} \,+\,\tau_{\,>1}\,-\,\epsilon\ \right]\, \cdot\,  \gamma \ } \ .
\end{eqnarray*}

\newpage

Hence
  \begin{eqnarray*}
(\,A.7.53\,) \ \ \ \  \ \ \ & \ & \Bigg\vert \ \int_{\R^n}\,[\ (\,{\tilde c}_n\!\cdot K\,)\ - \ n\,(\,n\,-\,2)\ ]\,\cdot\,\left[\    V_1^{{4}\over {\,n\,-\,2\,}}\,\cdot\,\left(\,\lambda_{\,1}\cdot  { { \partial \,V_1 } \over {\partial\, \lambda_{\,1}}}\ \right)\, \right]\,\cdot\,\phi \ \Bigg\vert \ \ \ \  \ \ \   \ \ \ \  \ \ \   \\[0.15in]
& \le & C\cdot  \Vert\,\Phi \,\Vert_*\,\cdot\,\left[\   {\bar\lambda}_{\ \flat}^{\,\ell} \ + \ {\bar\lambda}_{\ \flat}^{ \ \left[\  {{n\,+\,2}\over 2} \,+\,\tau_{\,>1}\,-\,\epsilon\ \right]\, \cdot\,  \gamma } \ \right] \\[0.2in]
& \le & C\cdot  {\bar\lambda}_{\ \flat}^{{{{\varpi}_{\,(\,A.5.38\,)} }}  }\,\cdot\,\left[\   {\bar\lambda}_{\ \flat}^{\,\ell} \ + \ {\bar\lambda}_{\ \flat}^{ \  \left[\  {{n\,+\,2}\over 2} \,+\,\tau_{\,>1}\,-\,\epsilon\ \right]\, \cdot\,  \gamma   } \ \right] \ .
\end{eqnarray*}

 \vspace*{0.1in}

 {\it Note.}  A.7.54.\ \ This part of the argument works in a similar way for  $\,\, \left( \ \lambda_{\,1}\,\cdot\,{ { \partial\, V_1} \over {\partial\, \xi_{1_{\,j}}}}\  \right)\,$.

 \vspace*{0.3in}

{\bf \S\,A\,7.\,j\,.} \ \   \ \ As a conclusion, we obtain\\[0.1in]
(\,A.7.55\,)
\begin{eqnarray*}
 \left(\ \lambda_{\,1}\cdot {{\partial }\over  {\partial\, \lambda_{\,1}}} \, \right)\,  {\bf I}_{\,R}
& = &  \left(\ \lambda_{\,1}\cdot {{\partial }\over  {\partial\, \lambda_{\,1}}} \, \right)\,  {\bf I}\,(\,W_{\,\,\flat}\,) \ + \   {\cal E}_{\,(\,A.{\bf 7}.55\,)} \\[0.2in]
   & = &  -\,\int_{\R^n} \left[\  (\Delta\,W_{\,\,\flat}) \ + \ (\,{\tilde c}_n\!\cdot K\,)\,(W_{\,\,\flat} )^{{\,n\,+\,2\,}\over {n\,-\,2}}  \ \right]\,\cdot\,\left(\,\lambda_{\,1}\cdot { { \partial \,V_1 } \over {\partial\, \lambda_{\,1}}}   \right)   \ + \   {\cal E}_{\,(\,A.{\bf 7}.55\,)}   \ .
 \end{eqnarray*}
 where\\[0.1in]
 (\,A.7.56\,)\\[0.1in]
 \hspace*{1.5in} \ \ \ \   \ \ \ \ \ \ \ \ \ (\,A.7.24\,) \ $\downarrow $ \ \ \ \ \ \ \ \ \ \  (\,A.7.11\,) \ $\downarrow $  \ \ \ \ \ \ \ \ \   \ \ \ \ (\,A.7.48\,) \ $\downarrow $
 $$
  {\cal E}_{\,(\,A.{\bf 7}.55\,)} \ = \ O\,\left(\, {\bar\lambda}_{\ \flat}^{{{{\varpi}_{\,(\,A.5.38\,)} }}  }\, \right)\,\cdot\,\Bigg[\  O\left( {\bar\lambda}_{\ \flat}^{{{4\over {\ n\,-\,2\ }} \,\cdot\ {{\varpi}_{\,(\,A.5.38\,)} }}  } \right) \  + \ \   O\left(\  {\bar\lambda}_{\ \flat}^{ {{\varpi}_{\,(\,A.5.38\,)} } \ - \ \sigma} \ \right) \   + \      O\left(\  {\bar\lambda}_{\ \flat}^{{ \, \,4\,\gamma \ -\ \tilde{\varepsilon}\,}  } \ \right)
 $$
$$
  \ \ \ \ \ \ \ \ \ \   \ \ \ \ \ \ \ \ \ \    \ \ \ \ \ \ \ \ \ \    \ \ \ \ \ \ \ \ \ \    \ \ \ \ \ \ \ \ \ \   \ \ \ \ \ \ \ \ \ \    \ \ \ \ \ \ \ \ \ \   \ \ \ \ \ \ \ \ \ \  +  \ \   O\left(\  {\bar\lambda}_{\ \flat}^{\,\ell}  \ \right)\ \  \Bigg]\ .
 $$
 \hspace*{4in}  \ \ \ \ \ \ \ \ \ \ \ (\,A.7.53\,) \,$\uparrow $\\[0.1in]
 Note that in (\,A.7.53\,)\,,\,
$$
n \ \ge  \ 6 \ \ \Longrightarrow \ \  O\left(\ {\bar\lambda}_{\ \flat}^{ \  \left[\  {{n\,+\,2}\over 2} \,+\,\tau_{\,>1}\,-\,\epsilon\ \right]\, \cdot\,  \gamma   } \ \right)  \ \ \ \ {\mbox{can \ \ be \ \ absorbed \ \ into}} \ \   \ \ O\left( \ {\bar\lambda}_{\ \flat}^{{ \, \,4\,\gamma \ -\ \tilde{\varepsilon}\,}  }\ \right)\ .
$$
On account of (\,A.5.38\,)\,:
$$
 {{\varpi}_{\,(\,A.5.38\,)} } \ = \
 \mbox{Min}  \ \left\{ \   {{n}\over 2}  \,\cdot\,(\,1\, - \ \nu\,)\,, \  \ \ \  \ell    \cdot  \nu \  \right\} \ \ - \ {1\over 2} \cdot \sigma \ - \ o_{\,+}\,(\,1\,)\ ,
 $$
as $\,\nu \ < \ 1\,,$\, it is proper to say that
$$
{{\varpi}_{\,(\,A.5.38\,)} } \ < \ l \ \ \Longrightarrow \ \ {{\varpi}_{\,(\,A.5.38\,)} } \ - \ \sigma \ < \ l \,.
$$
Hence
$$
   O\left(\ {\bar\lambda}_{\ \flat}^{ \  \ell  } \ \right)  \ \ \ \ {\mbox{can \ \ be \ \ absorbed \ \ into}} \ \   \ \ O\left( \ {\bar\lambda}_{\ \flat}^{{ \, \, {{\varpi}_{\,(\,A.5.38\,)} } \ - \ \sigma\,}  }\ \right)\ .
$$
We finally arrive at
 \begin{eqnarray*}
  (\,A.7.57\,)&\ &  {\cal E}_{\,(\,A.{\bf 7}.55\,)} \ = \ O\,\left(\, {\bar\lambda}_{\ \flat}^{{{{\varpi}_{\,(\,A.5.38\,)} }}  }\, \right)\,\cdot\,\Bigg[\  O\left(\   {\bar\lambda}_{\ \flat}^{{{4\over {n\,-\,2}} \,\cdot\ {{\varpi}_{\,(\,A.5.38\,)} }}  } \ \right)  \ \ + \ \ O\left(\    {\bar\lambda}_{\ \flat}^{{ \, \,4\,\gamma \ -\ o_{\,+}\,(\,1\,)\,}  }  \ \right) \\[0.2in]
   & \ &  \ \ \ \ \ \ \ \ \   \ \ \ \ \ \ \ \ \  \ \ \ \ \ \ \ \ \   \ \ \ \ \ \ \ \ \   \ \ \ \ \ \ \ \ \ \ \ \ \ \ \ \ \ \ \ \ \ \ \ \ \ \ \ + \ \   O\left(\   {\bar\lambda}_{\ \flat}^{ {{\varpi}_{\,(\,A.5.38\,)} } \ - \ \sigma}    \ \right)\ \  \Bigg]\\[0.1in]
  & = &  O\,\left(\, {\bar\lambda}_{\ \flat}^{{{{\varpi}_{\,(\,A.5.38\,)} \ + \ \mu_{\,(A.7.57\,)}}}  }\, \right) \ ,
 \end{eqnarray*}
where
$$
\mu_{\,(A.7.57\,)} \, = \, \mbox{Min} \   \left\{ \  {{4\over {n\,-\,2}} \,\cdot\ {\varpi}_{\,(\,A.5.38\,)}}  \ , \  \   4\,\gamma \ -\ o_{\,+}\,(\,1\,) \ , \ \   {\varpi}_{\,(\,A.5.38\,)} \ - \ \sigma\    \right\}\ . \leqno (\,A.7.58\,)
$$
\vspace*{0.2in}
Here we work under the assumptions in Proposition 6.8\,.\, Note that condition (\,1.34\,) of the main text\, , namely\,,
$$
\mbox{Min}  \, \left\{ \     {{n}\over 2} \,\cdot\,(\,1\, - \ \nu\,)\,, \  \ \ \  \ell    \cdot  \nu \  \right\}  \ \,-\ o_{\,+}\,(\,1\,) \ > \ \left( {{\,n\ + \ 10}\over 8}\  \right) \cdot \sigma \ \ge \ 2\cdot \sigma \ \ \ \ (\,n \ \ge \ 6\,)\ ,
$$
implies
$$
{\varpi}_{\,(\,A.5.38\,)} \ - \ \sigma \ > \ 0\ .
$$

\vspace*{0.2in}

{\it Note.} \ \ A\,.\,7\,.\,59\ . \ \
 Likewise,
 \begin{eqnarray*}
 \left(\ \lambda_{\,1}\cdot {{\partial }\over {\partial\, \xi_{1_{\,|_j}}}}  \, \right)\,  {\bf I}_R & = &  \left(\ \lambda_{\,1}\cdot {{\partial }\over  {\partial\, \xi_{1_{\,|_j}}}}  \, \right)\,  {\bf I}\,(\,W_{\,\,\flat}\,) \ + \   {\cal E}_{\,(\,A.{\bf 7}.55\,)} \\[0.2in] & = &  -\,\int_{\R^n} \left[\  (\Delta\,W_{\,\,\flat}) \ + \ (\,{\tilde c}_n\!\cdot K\,)\,(W_{\,\,\flat} )^{{\,n\,+\,2\,}\over {n\,-\,2}}  \ \right]\,\cdot\,\left(\,\lambda_{\,1}\cdot { { \partial \,V_1 } \over {\partial\, \xi_{1_{\,|_j}}}}   \right)   \ \ + \  \   {\cal E}_{\,(\,A.{\bf 7}.55\,)}  \ .
 \end{eqnarray*}
Here $\,j \ = \ 1\,, \ \cdot \cdot \cdot\,, \ n\,.$\, Along the same line\,,\,  we deal with other bubbles parameters $\,\xi_{\,2}\,, \ \cdot \cdot \cdot\,, \ \xi_{\,\flat}\ .$

\newpage

{\bf \large \S\,A\,8.} \ \ {\bf \large Extracting the key information in} ${\bf(\,A\,)}_{\,(\,A.8.1\,)}$ -- \ {\bf \large remaining }\\[0.1in]\hspace*{0.65in} {\bf \large estimates\ \  [\,refer to \S\,4\,b in the main text\ ]}\ . \\[0.2in]
We continue from (\,A.7.55\,)\,.\\[0.1in]
(\,A.8.1\,)
\begin{eqnarray*}
& \ & \left(\ \lambda_{\,1}\cdot {\partial\over {\partial\, \lambda_{\,1}}}\ \right)  \, {\bf I}_{\,\cal R}\\[0.2in]  & = &  - \int_{\R^n} \left[\  \Delta\,W_{\,\,\flat}  \ + \  (\,{\tilde c}_n\!\cdot K\,)\,(W_{\,\,\flat} )^{{\,n\,+\,2\,}\over {n\,-\,2}}  \ \right]\,\cdot\,\left[ \, \left(\ \lambda_{\,1}\cdot {\partial\over {\partial\, \lambda_{\,1}}}\ \right) V_1   \, \right]   \ \ \ + \ \    {\cal E}_{\,(\,A.{\bf 7}.55\,)}  \\[0.2in]
  & = & \ -\int_{\R^n} \left[\  \Delta\,W_{\,\,\flat} \ + \  n\,(\,n\,-\,2)\cdot\,(W_{\,\,\flat} )^{{\,n\,+\,2\,}\over {n\,-\,2}} \ \right]\,\cdot\,\left[ \, \left(\ \lambda_{\,1}\cdot {\partial\over {\partial\, \lambda_{\,1}}}\ \right) V_1   \, \right]    \ \  \,\cdot\,\,\cdot\,\cdot\,\cdot\,\,\cdot\,\cdot\,\cdot\, \cdot\,\cdot\,\cdot\ \   {\bf (A)}_{(\,A.8.1\,)}\\[0.2in]
 & \ &  \ \ \ \ \ \ \ \ \ \ \  \ \  - \    \int_{\R^n}\, [\,(\,{\tilde c}_n\!\cdot K\,) \ - \ n\,(\,n\,-\,2)\,]  \cdot\,(W_{\,\,\flat} )^{{\,n\,+\,2\,}\over {n\,-\,2}} \cdot\left[ \, \left(\ \lambda_{\,1}\cdot {\partial\over {\partial\, \lambda_{\,1}}}\ \right) V_1   \, \right]    \ \ \,\cdot\,\cdot\,\cdot\,\cdot   \ \ {\bf (B)}_{(\,A.8.1\,)}\\[0.2in]
 & \ & \ \ \ \ \ \ \    \ \ \ \ \ \ \   \ \ \ \ \ \ \  \ \ \    \ + \ \ {\cal E}_{\,(\,A.{\bf 7}.55\,)}  \ .
 \end{eqnarray*}
 Through the equation
 $$
 \ \ \  \ \ \ \ 
 -\,\Delta \,V_{\,l}\ = \ n\,(\,n\,-\,2)\,\cdot\, V_{\,l}^{{\,n\,+\,2\,}\over {n\,-\,2}}\ \ \  \ \ \ \  \ \ \  \ \ \ \  (\ l \ = \ 1\,,\ \cdot \cdot \cdot\ , \ \flat\ )\,,\,
 $$
 we obtain
 \\[0.1in]
(\,A.8.2\,)
  \begin{eqnarray*}
 {\bf (A)}_{(\,A.8.1\,)}  \ & = &  n\,(\,n\,-\,2)\,\cdot\,\int_{\R^n} \Bigg\{  \,\left[\  V_1^{{\,n\,+\,2\,}\over {n\,-\,2}} \ + \  V_2^{{\,n\,+\,2\,}\over {n\,-\,2}} \ +\,\cdot\,\cdot\,\cdot\, + \ V_{\,\flat}^{{\,n\,+\,2\,}\over {n\,-\,2}} \ \right] \\[0.2in]
    & \ & \ \ \ \ \ \  \ \ \ \ \ \  \ \ \ \ \ \  \ \ \ \ \
    - \ \left[\  V_1\ + \  V_2  \ + \,\cdot\,\cdot\,\cdot\, + \ V_{\,\flat}  \ \right]^{{\,n\,+\,2\,}\over {n\,-\,2}} \  \, \Bigg\} \cdot\left[ \, \left(\ \lambda_{\,1}\cdot {\partial\over {\partial\, \lambda_{\,1}}}\ \right) V_1   \, \right] \ .\\
\end{eqnarray*}
In order to discern the interaction between bubbles, we decomposite the integral in $\, {\bf (A)}_{\,(\,A.8.1\,)}\,$ as   \\[0.1in]
(\,A.8.3\,)
$$
\int_{\R^n} \ = \ \int_{B_{\,\xi_{\,1}} (\,{\rho_{\,\nu}}\,)}  \! + \ \ \sum_{l\,=\,2}^\flat \ \int_{B_{\xi_l} (\,{\rho_{\,\nu}})} \ + \ \  \int_{\R^n \,\setminus \, \left(\ \cup\ B_{\xi_l} (\,{\rho_{\,\nu}}) \, \right)}  \ =\,:\ {\bf A}^1_{\,(\,A.8.1\,)} \ + \ {\bf A}^2_{\,(\,A.8.1\,)}  \ + \ {\bf A}^3_{\,(\,A.8.1\,)} \ .
$$

\newpage

{\bf \S\,A\,8\,.\,a\,.}\ \   {\bf Extracting  the key information in} $\,{\bf A}^1_{\,(\,A.8.1\,)}\,$ [\,{\it  inside}\, $\,B_{\,\xi_{\,1}}(\,\rho_{\,\nu}\,)\,$\,]\,.\ \
In the ball $\,B_{\,\xi_{\,1}}(\,\rho_{\,\nu}\,)\,$\,,\,
$$
 {{\ V_2\,(\,y\,) \ + \ \cdot\,\cdot\,\cdot\,\ + \ V_{\,\flat}\,(\,y\,)\   }\over {V_1\,(\,y\,)  }} \ = \ O\,\left( \ {\bar\lambda}_{\ \flat}^{\ (\,n\,-\,2\,)\,\cdot\, [\ (\,\gamma\,+\,\nu\,) \ - \ 1\,]}\ \right)
\ \ \  \ \ \mfor \ \ \  y\,\in\,B_{\,\xi_{\,1}}(\,\rho_{\,\nu})\ .
$$
Similar to (\,A.4.33\,)\,--\,(\,A.4.36\,)\,,\,
we have\\[0.1in]
(\,A.8.4\,)
\begin{eqnarray*}
& \ &   \left[\  V_1\ + \  V_2  \ + \ \cdot\,\cdot\,\cdot\,\ + \ V_{\,\flat}  \ \right]^{{\,n\,+\,2\,}\over {n\,-\,2}}  \ = \ V_1^{{\,n\,+\,2\,}\over {n\,-\,2}}\,\cdot\,\left[\  1 \ + \ \bigg( \ {{  V_2  \ + \ \cdot\,\cdot\,\cdot\,\ + \ V_{\,\flat}}\over {V_1}} \ \bigg)  \, \right]^{{\,n\,+\,2\,}\over {n\,-\,2}}\\[0.15in]
& = & V_1^{{\,n\,+\,2\,}\over {n\,-\,2}}\,\cdot\,\left[\  1 \ + \ {{n\,+\,2}\over {\,n\,-\ 2\,}}\,\cdot\,\left(\ \ {{  V_2  \ + \ \cdot\,\cdot\,\cdot\,\ + \ V_{\,\flat}}\over {V_1}}  \ \right)\ + \ O \left(\ \  {{  V_2  \ + \ \cdot\,\,\cdot\,\cdot \ + \ V_{\,\flat}}\over {V_1}}\  \right)^{{\,n\,+\,2\,}\over {n\,-\,2}}\  \right]\\[0.15in]
& = & V_1^{{\,n\,+\,2\,}\over {n\,-\,2}} \ + \  {{n\,+\,2}\over {\,n\,-\ 2\,}} \,\cdot\,V_1^{4\over {n\,-\,2}}\,\cdot\, \left[\   V_2  \ + \ \cdot\,\cdot\,\cdot\,\ + \ V_{\,\flat}  \ \right] \ + \ O \left(\,  \left[ \, V_2  \ + \ \cdot\,\cdot\,\cdot\,\ + \ V_{\,\flat} \, \right]^{  {{n\,+\,2}\over {\,n\,-\ 2\,}}}\  \right) \\[0.1in]
& \ & \hspace*{5in} {\mbox{in}} \ \ \ B_{\,\xi_{\,1}} (\,{\rho_{\,\nu}}\,)\ .
\end{eqnarray*}

It follows from (\,A.8.3\,) and (\,A.8.4\,)  that\\[0.1in]
(\,A.8.5\,)
\begin{eqnarray*}
& \ &  n\,(\,n\,-\,2)\,\cdot\,\int_{B_{\,\xi_{\,1}} (\,{\rho_{\,\nu}}\,)}     \Bigg\{  \,\left[\  V_1^{{\,n\,+\,2\,}\over {n\,-\,2}} \ + \  V_2^{{\,n\,+\,2\,}\over {n\,-\,2}} \ +\,\cdot\,\cdot\,\cdot\, + \ V_{\,\flat}^{{\,n\,+\,2\,}\over {n\,-\,2}} \ \right]  \\[0.2in]
    & \ &   \ \ \ \ \ \  \ \ \ \ \ \  \ \ \ \ \ \  \ \ \ \ \ \  \ \ \ \ \ \  \ \ \
    - \ \left[\  V_1\ + \  V_2  \ + \,\cdot\,\cdot\,\cdot\, + \ V_{\,\flat}  \ \right]^{{\,n\,+\,2\,}\over {n\,-\,2}} \  \Bigg\}\,\cdot\,\left[ \, \left(\ \lambda_{\,1}\cdot {\partial\over {\partial\, \lambda_{\,1}}}\ \right) V_1   \, \right] \\[0.2in]
& = & -\,n\,(\,n\,-\,2)\,\cdot\int_{B_{\,\xi_{\,1}} (\,{\rho_{\,\nu}}\,)} \!\!\left\{ \  {{n\,+\,2}\over {\,n\,-\ 2\,}}  \,\cdot\, \left[\   V_2  \ + \ \cdot\,\cdot\,\cdot\,\ + \ V_{\,\flat}  \ \right] \   \right\}\cdot\,V_1^{4\over {n\,-\,2}} \cdot\left[ \, \left(\ \lambda_{\,1}\cdot {\partial\over {\partial\, \lambda_{\,1}}}\ \right) V_1   \, \right]\\[0.2in]
& \ & \hspace*{4in} + \ \ {\bf{error}}_{\,(\,A.8.5\,)}\ ,
\end{eqnarray*}
where
$$
|\,{\bf{error}}_{\,(\,A.8.5\,)}\ \,| \ \le \  C\,\cdot\, \int_{B_{\,\xi_{\,1}} (\,{\rho_{\,\nu}}\,)}\left(  \, V_2  \ + \ \cdot\,\cdot\,\cdot\,\ + \ V_{\,\flat} \, \right)^{  {{n\,+\,2}\over {\,n\,-\ 2\,}}}\,\cdot\, \bigg\vert \  \lambda_{\,1}\,\cdot\,{ { \partial \,V_1} \over {\partial\, \lambda_{\,1}}}\ \bigg\vert \ . \leqno
(\,A.8.6\,)
$$

\newpage

 {\bf \S\,A\,8\,.\,b\,.} \ \
{\bf Estimate of $\ |\,{\bf{error}}_{\,(\,A.8.5\,)}\,| \,$.}\,  \ \
Similar to (\,A.7.29\,)\,,\, in $\,B_{\,\xi_{\,1}}(\,\rho_{\,\nu}\,)\,$,\,\\[0.1in]
(\,A.8.7\,)
\begin{eqnarray*}
& \ & \left[\ \,V_2\,(\,y\,)   \ + \ \cdot\,\cdot\,\cdot\,\ + \ V_{\,\flat}\,(\,y\,)\,\right]^{{\,n\,+\,2\,}\over {n\,-\,2}}\\[0.2in]
& \le & C\,\cdot\,\left[\ \left(\ {\lambda_{\,2}\over {\lambda^2_2 \ + \ \Vert\ \xi_{\,1}\,-\,\xi_{\,2}\,\Vert^2 }}\  \right)^{{\,n\,-\,2\,}\over 2} \ \ + \ \cdot\,\cdot\,\cdot\,\ + \  \ \left(\ {{\bar\lambda}_{\ \flat}\over {\lambda^2_\flat \ + \ \Vert\ \xi_{\,1}\,-\,\xi_\flat\,\Vert^2 }}\  \right)^{{\,n\,-\,2\,}\over 2} \ \right]^{\,{{n\,+\,2}\over {\,n\,-\ 2\,}}} \\[0.2in]
& \le & C_1\,\cdot\,\left(\ {1\over {{\bar\lambda}_{\ \flat}   }}\  \right)^{{n\,+\,2}\over 2}\,\cdot\,\left(\  \sum_{l\,=\,2}^\flat \,{1\over { \Vert\ \Xi_{\,\,1} \ - \ \Xi_{\,l}\,\Vert^{\,n\,-\,2}  }}\  \right)^{\!\!{{n\,+\,2}\over {\,n\,-\ 2\,}}} \ \le \ C_2\,\cdot\,{{ {\bar\lambda}_{\ \flat}^{(\,n\,+\,2)\,\cdot\, \gamma}}\over { {\bar\lambda}_{\ \flat}^{{n\,+\,2}\over 2}  }} \ \ \ \ \ \ \ \ \ [\ y\,\in\,B_{\,\xi_{\,1}}(\,\rho_{\,\nu}\,)\ ]\ .\\
\end{eqnarray*} Derivative in the $\,\xi_{1_{\,|_n}}$-\,direction\,.
Moreover, via (\,A.3.17\,)\,,
\begin{eqnarray*}
 \int_{B_{\,\xi_{\,1}} (\,{\rho_{\,\nu}}\,)}  \bigg\vert \  \lambda_{\,1}\,\cdot\,{ { \partial \,V_1} \over {\partial\, \lambda_{\,1}}}\ \bigg\vert & \le & C_1  \int_{B_{\,\xi_{\,1}} (\,{\rho_{\,\nu}}\,)}  V_1\\[0.2in]
& \le & C_1 \int_{B_{\,\xi_{\,1}}(\,\rho_{\,\nu}\,)} \left(\ {\lambda_{\,1}\over {\lambda^2_1 \ + \ \Vert\ y\,-\,\xi_{\,1}\,\Vert^2 }}\  \right)^{{\,n\,-\,2\,}\over 2} \ d \,y\\[0.2in]
& \le &
C_2\,\cdot\,\lambda_{\,1}^{\! {{\,n\,+\,2\,}\over 2} }\,\cdot\,\left(\ {{\rho_{\,\nu}}\over \lambda_{\,1}} \right)^{\!\!2}\ .
\end{eqnarray*}
It follows that
\begin{eqnarray*}
(\,A.8.8\,) \ \ \ \ \ |\,{\bf{error}}_{\,(\,A.8.5\,)}\ \,|  & \le &  C_3\,\cdot\, {\bar\lambda}_{\ \flat}^{(\,n\,+\,2\,)\,\cdot\,\ \gamma}\,\cdot\, \left(\ {\rho_{\,\nu}\over {\bar\lambda}_{\ \flat}}\  \right)^{\!\!2} \\[0.2in]& = &  O\,\left( \  {\bar\lambda}_{\ \flat}^{(\,n\,+\,2\,)\,\cdot\,\ \gamma}\ \right)\,\cdot\,O\,\left(\ {1\over { {\bar{\lambda}}_{\ \flat}^{\,2\ (\,1\,-\ \nu\,)} }}\  \right)\\[0.2in]
& \le & O\,\left( \  {\bar\lambda}_{\ \flat}^{\,n\,\cdot\, \gamma}\ \right)\,\cdot\, O\,\left( \  {\bar\lambda}_{\ \flat}^{\ 2\ [\ (\,\gamma \ + \ \nu \,) \ - \ 1\ ]}\ \right)\ \ \ \ \ \ \ \ \   \ \ \ \ \ (\ \gamma \ + \ \nu \ >  \ 1\ )\,.
\end{eqnarray*}


\newpage

 {\bf \S\,A\,8\,.\,c\,.} \ \
{\bf Estimate of \,$\,{\bf T}\,.$}  \ \
Recall from (\,4.11\,) of the main text\,  that
$$
{\bf T} \ =\  {1\over { {\bf d}_{1\,,\,2}^{\,2}  }}\,\cdot\,\left(\ {{\lambda_{\,2}}\over {\lambda_{\,1}}}  \right) \ + \ {1\over  { {\bf d}_{1\,,\,2}^{\,2}  }}\,\cdot\,{{\Vert\,\,y\ - \ \xi_{\,1}\,\Vert^2}\over
{\lambda_{\,1}\,\cdot\,\,\lambda_{\,2}}}   \ + \ {1\over { {\bf d}_{1\,,\,2}^{\,2}  }}\cdot {{\,2\,(\,y\ - \ \xi_{\,1})\,\cdot\,(\,\xi_{\,1}\ - \ \xi_{\,2})\, }\over
{\lambda_{\,1}\,\cdot\,\,\lambda_{\,2}}}\ , \leqno (\,A.8.9\,)
$$
where
$$
{\bf d}_{1\,,\,2}^2 \ = \  {{\Vert\,\xi_{\,1}\ -\ \xi_{\,2}\,\Vert^{\,2}}\over{\, {\,\lambda_{\,1} \,\cdot\,\lambda_{\,2}\,}\,}}\ .
$$
The first term in the R.H.S. of (\,A.8.9\,) is estimated by
$$
  {1\over { {\bf d}_{1\,,\,2}^{\,2}  }}\,\cdot\,\left(\ {{\lambda_{\,2}}\over {\lambda_{\,1}}}  \right) \ \le \  C_1\,\cdot\, {1\over {  {{\Vert\,\xi_{\,1}\ -\ \xi_{\,2}\,\Vert^{\,2}}\over{\, {\,\lambda_{\,1} \,\cdot\,\lambda_{\,2}\,}\,}}   }}  \ \le  \ C_2\,\cdot\, {1\over { \Vert  \ \Xi_{\,\,1} \ - \ \Xi_{\,2} \ \Vert^{\,2} }}\ =  \  O\,\left(\, {\bar\lambda}_{\ \flat}^{\,2\,\gamma} \,\right)\ .  \leqno (\,A.8.10\,)
$$
Observe that
$$
 {1\over{ {\bf d}_{1\,,\,2}^{\,2}  }}\cdot {{\Vert\,\,y\ - \ \xi_{\,1}\ \Vert^2}\over
{\lambda_{\,1}\,\cdot\,\,\lambda_{\,2}}}   \  \le \    {1\over { {\bf d}_{1\,,\,2}^{\,2}  }}\,\cdot\,{{\rho^2_\nu}\over
{\lambda_{\,1}\,\cdot\,\,\lambda_{\,2}}}\ = \  O\,\left(\  {\bar\lambda}_{\ \flat}^{\ 2\,[\ (\,\gamma\,+\,\nu\,) \, - \ 1\,]} \ \right)  \leqno (\,A.8.11\,)
$$
for $\,y \, \in \, B_{\,\xi_{\,1}} (\,\rho_\nu\,)\,.$\,
The last term  in the R.H.S. of (\,A.8.9\,) is not so small \{\,when compared to the first two terms in the R.H.S. of (\,A.8.9\,)\,,\, but, being asymmetric, it can be integrated away [ cf. (\,A.8.19\,)\ ] \}\,:
 \begin{eqnarray*}
(\,A.8.12\,)   \ \ \ \ \ \ \  \ \ \ \ \ \  \  \ & \ &  \Bigg\vert\  {1\over { {\bf d}_{1\,,\,2}^{\,2}  }}\,\cdot\,{{\, (\,y\ - \ \xi_{\,1})\cdot(\,\xi_{\,1}\ - \ \xi_{\,2})\, }\over
{\lambda_{\,1}\,\cdot\,\,\lambda_{\,2}}}\ \Bigg\vert\\[0.2in]
& \le & {1\over{ {\bf d}_{1\,,\,2}^{\,2}  }}\cdot {{\,\Vert\,y\ - \ \xi_{\,1}\,\Vert\, }\over
{\sqrt{\lambda_{\,1}\,\cdot\,\,\lambda_{\,2}\,}}} \,\cdot\,{{\,\Vert\,\xi_{\,1}\ - \ \xi_{\,2}\,\Vert\, }\over
{\sqrt{\lambda_{\,1}\,\cdot\,\,\lambda_{\,2}\,}}} \\[0.2in]
& \le & {1\over{ {\bf d}_{1\,,\,2}^{\,2}  }}\,\cdot\,{{\,\rho_{\,\nu}\, }\over
{\sqrt{\lambda_{\,1}\,\cdot\,\,\lambda_{\,2}\,}}} \,\cdot\,{\bf d}_{1\,,\,2}  \ =  \  O\,\left(\  {\bar\lambda}_{\ \flat}^{\ (\,\gamma\,+\,\nu\,) \, - \ 1 } \ \right)\\[0.2in]
\Longrightarrow \ \ \    {\bf T} &  = & O\,\left(\ {\bar\lambda}_{\ \flat}^{2\,\gamma}\ \right) \ + \  O\,\left(\  {\bar\lambda}_{\ \flat}^{\ (\,\gamma\,+\,\nu\,) \, - \ 1 } \ \right)\\[0.2in]
& = & O\,\left(\  {\bar\lambda}_{\ \flat}^{\ (\,\gamma\,+\,\nu\,) \  - \ 1 } \ \right)  \ \ \  \mfor \ \ y\,\in\,{B_{\,\xi_{\,1}} (\,{\rho_{\,\nu}}\,)} \\[0.2in]
& \ &   [\ {\mbox{as}} \ \ 2\,\gamma \ >  \ (\,\gamma\,+\,\nu\,) \, - \ 1 \ \Longleftrightarrow \ \ 1 \ + \ \gamma \ >  \ \nu\ ; \ \ {\mbox{also}} \ \ \nu \ < \ 1\  ]\,. \ \ \ \ \ \ \  \ \ \ \ \ \  \ \ \ \ \ \ \  \ \ \ \ \ \  \
\end{eqnarray*}
\vspace*{0.1in}

\newpage

From (\,A.8.5\,)\,,\,
we continue with
 \begin{eqnarray*}
(\,A.8.13\,)
 & \ & -\,n\,(\,n\,-\,2)\,\cdot\,\int_{B_{\,\xi_{\,1}} (\,{\rho_{\,\nu}}\,)} \left\{ \, \left(\ \,{{n\,+\,2}\over {\,n\,-\ 2\,}} \,\right)\,\cdot\,V_1^{4\over {n\,-\,2}}\,\cdot\,\left[\   V_2   \ \right]  \  \right\}\,\cdot\, \left(\,\lambda_{\,1}\cdot { { \partial \,V_1} \over {\partial\, \lambda_{\,1}}}\ \right)\\[0.2in]
& =  &   +\ n\,(\,n\,-\,2) \cdot\left( \,{{n\,-\,2}\over { 2}} \,\right)\cdot  \left(\ \,{{n\,+\,2}\over {\,n\,-\ 2\,}} \,\right)\,\cdot\,\int_{B_{\,\xi_{\,1}}\,(\,{\rho_{\,\nu}})} V_1^{4\over {n\,-\,2}}
\cdot \left(\ {1\over {  \lambda_{\,1}^{{n\,-\,2}\over 2 }   }}\,\cdot\,{1\over {{\bf d}_{1\,,\,2}^{\,n\,-\,2}    }}\   \right) {\bf *} \\[0.2in]
& \ &   \hspace*{-0.4in}  {\bf *}
 \ \left\{  \  \lambda_{\,1}^{{n \,-\, 2}\over 2}
\cdot {{ \lambda_{\,1}^2\, -\, \Vert\, y \,-\, \xi_{\,1}\,\Vert^{\,2}}\over {(\,\lambda_{\,1}^2 \,+\, \Vert\, y \,-\, \xi_{\,1}\,\Vert^{\,2})^{{n }\over 2} }}\ \right\} \   \times \left\{ \  1 \ - \  \left(\,{{\,n\,-\,2\,}\over 2} \right)\,\cdot\,{\bf T}  \ + \  O\,(\,|\  {\bf T} \,|^{\,2}\,) \ \right\}\\[0.2in]
& =  &   +\ n\,(\,n\,-\,2) \cdot\left( \,{{n\,-\,2}\over { 2}} \,\right)\cdot  \left(\ \,{{n\,+\,2}\over {\,n\,-\ 2\,}} \,\right)\,\cdot\,\int_{B_{\,\xi_{\,1}}\,(\,{\rho_{\,\nu}})} {\bf[} \ \cdot\,\cdot \cdot\  |_{\,y} \  {\bf ]}_{\, (\,A.8.14\,)}   \ + \ \\[0.1in]
& \
 & \ \ \  \left[ \ \ \uparrow \ \  {\mbox{of \ \ order}} \ \   {{1}\over {{\bf d}_{\,1\,,\ 2}^{\,n\,-\ 2 }  }}   \ ; \ \ {\mbox{discussed \ \ in \ \ the \ \ main \ \ text\,, \ \ see \ \  {\bf \S\,4.\,b}}} \ \right]\\[0.1in]
& \ & \hspace*{-0.5in} +\ \ C_1\,\cdot \int_{B_{\,\xi_{\,1}}\,(\,{\rho_{\,\nu}})}  {\bf[} \ \cdot\,\cdot \cdot\  |_{\,y} \  {\bf ]}_{\, (\,A.8.14\,)}   \cdot {\bf T} \ + \ O \left( \  \int_{B_{\,\xi_{\,1}}\,(\,{\rho_{\,\nu}})}
{\bf[} \ \cdot\,\cdot \cdot\  |_{\,y} \  {\bf ]}_{\, (\,A.8.14\,)}  \cdot | \  {\bf T}\ |^2 \ \right)\ .
    \end{eqnarray*}
Here
 $$
 {\bf[} \ \cdot\,\cdot \cdot\  |_{\,y} \  {\bf ]}_{\, (\,A.8.14\,)}   \ = \  [\,V_1\,(\,y\,)\,]^{4\over {n\,-\,2}}\
\cdot\, \left(\ {1\over {  \lambda_{\,1}^{{n\,-\,2}\over 2 }   }}\,\cdot\,{1\over {{\bf d}_{1\,,\,2}^{\,n\,-\,2}    }}\   \right)\,\cdot\,\left\{  \  \lambda_{\,1}^{{n \,-\, 2}\over 2}
\cdot {{ \,\lambda_{\,1}^2\, -\, \Vert\, y \,-\, \xi_{\,1}\,\Vert^{\,2}\,}\over {\ (\,\lambda_{\,1}^2 \,+\, \Vert\, y \,-\, \xi_{\,1}\,\Vert^{\,2})^{{n }\over 2}\  }}\ \right\}\ .  \leqno (\,A.8.14\,)
 $$
From (\,A.8.9\,) and (\,A.8.10\,)\,,\, \\[0.1in]
(\,A.8.15\,)
 \begin{eqnarray*}
  \int_{B_{\,\xi_{\,1}} (\,{\rho_{\,\nu}}\,)}  {\bf[} \ \cdot\,\cdot \cdot\  |_{\,y} \  {\bf ]}_{\, (\,A.8.14\,)}    \cdot {\bf T} & = &   \left[\  \int_{B_{\,\xi_{\,1}}\,(\,{\rho_{\,\nu}})} {\bf[} \ \cdot\,\cdot \cdot\  |_{\,y} \  {\bf ]}_{\, (\,A.8.14\,)}  \ \right] \,\,\times\,O\,(\,{\bar{\lambda}}_{\ \flat}^{\,2\,\gamma}\,)\ \ + \ \\[0.2in]
  & \ & \  + \  \int_{B_{\,\xi_{\,1}}\,(\,{\rho_{\,\nu}})}{\bf[} \ \cdot\,\cdot \cdot\  |_{\,y} \  {\bf ]}_{\, (\,A.8.14\,)}   \,\,\times\,\left[\   {1\over{ {\bf d}_{1\,,\,2}^{\,2}  }}\cdot {{\Vert\,\,y\ - \ \xi_{\,1}\ \Vert^2}\over
{\lambda_{\,1}\,\cdot\,\,\lambda_{\,2}}} \ \right]\\[0.2in]
& \ &\hspace*{-1in}  + \  \int_{B_{\,\xi_{\,1}}\,(\,{\rho_{\,\nu}})}  {\bf[} \ \cdot\,\cdot \cdot\  |_{\,y} \  {\bf ]}_{\, (\,A.8.14\,)}   \,\,\times\,\left[\   {1\over { {\bf d}_{1\,,\,2}^{\,2}  }}\cdot {{\,2\,(\,y\ - \ \xi_{\,1})\,\cdot\,(\,\xi_{\,1}\ - \ \xi_{\,2})\, }\over
{\lambda_{\,1}\,\cdot\,\,\lambda_{\,2}}}\ \right]  \ .
    \end{eqnarray*}

\newpage

 {\bf \S\,A\,.\,8\,.\,d\,.} \ \
{\bf Estimate of first order teams in}\, (\,A.8.13\,)\,. \\[0.1in]
Let us consider the second team in the R.H.S. of (\,A.8.15\,)\,:\\[0.1in]
(\,A.8.16\,)
 \begin{eqnarray*}
 & \ &  \int_{B_{\,\xi_{\,1}}\,(\,{\rho_{\,\nu}})}{\bf[} \ \cdot\,\cdot \cdot\  |_{\,y} \  {\bf ]}_{\, (\,A.8.14\,)} \,\,\times\,\left[\   {1\over{ {\bf d}_{1\,,\,2}^{\,2}  }}\cdot {{\Vert\,\,y\ - \ \xi_{\,1}\ \Vert^2}\over
{\lambda_{\,1}\,\cdot\,\,\lambda_{\,2}}} \ \right] \ \ \ \ \ \ \ (\ {\mbox{take}} \ \ r \ = \ \Vert\,\,y\ - \ \xi_{\,1}\ \Vert \ )\\[0.2in]
& = & C \cdot   \left\{ \   {{1}\over {{\bf d}_{\,1\,,\ 2}^{\,n\,-\ 2 }  }}   \  \right\}\,\cdot\,\left\{ \,  {1\over { {\bf d}_{1\,,\,2}^{\,2}  }}\right\}\,\cdot\, \left[\ \int_0^{\rho_{\,\nu}}\!
\left[   \    {{\lambda_{\,1}^2}\over {(\,\lambda^2_1 \,+\, r^2\,)^{\,2}}}
\cdot {{ \lambda_{\,1}^2}\over {(\,\lambda_{\,1}^2 \ + \  r^{\,2})^{{n }\over 2} }}\right]\,\cdot\,r^{n\,-\,1}\,\cdot\, {{r^2}\over {\lambda_{\,1}\cdot \lambda_{\,2}}}\ dr   \ \right] \\[0.1in]
& \ & \ \ \ \ \ \ \ \ \  [ \ \uparrow \ \ {\mbox{unit \ \ for \ \ this \ \ part \ \ of \ \ the  \ \ calculation}} \ ]\\[0.1in]
& \ & \hspace*{-0.35in}\ \  - \ C \cdot  \left\{ \   {{1}\over {{\bf d}_{\,1\,,\ 2}^{\,n\,-\ 2 }  }}   \  \right\}\,\cdot\,\left\{ \,  {1\over { {\bf d}_{1\,,\,2}^{\,2}  }}\right\}\,\cdot\, \left[\ \int_0^{\rho_{\,\nu}}\!
\left[   \    {{  \lambda_{\,1}^2    }\over {(\,\lambda^2_1 \,+\, r^2\,)^{\,2}}}
\cdot {{ r^2    }\over {(\,\lambda_{\,1}^2 \ + \  r^{\,2})^{{n }\over 2} }}\right]\,\cdot\,r^{n\,-\,1}\,\cdot\, {{r^2}\over {\lambda_{\,1}\cdot \lambda_{\,2}}}\ dr   \ \right] \ . \\
  \end{eqnarray*}
We have\\[0.1in]
(\,A.8.17\,)
 \begin{eqnarray*}
&   & \left\{ \   {{1}\over {{\bf d}_{\,1\,,\ 2}^{\,n\,-\ 2 }  }}   \  \right\}\,\cdot\,\left\{ \,  {1\over { {\bf d}_{1\,,\,2}^{\,2}  }}\right\}\,\cdot\, \left[\ \int_0^{\rho_{\,\nu}}\!
\left[   \    {{\lambda_{\,1}^2}\over {(\,\lambda^2_1 \,+\, r^2\,)^{\,2}}}
\cdot {{ \lambda_{\,1}^2}\over {(\,\lambda_{\,1}^2 \ + \  r^{\,2})^{{n }\over 2} }}\right]\,\cdot\,r^{n\,-\,1}\,\cdot\, {{r^2}\over {\lambda_{\,1}\cdot \lambda_{\,2}}}\ dr   \ \right] \\[0.2in]
& \le &C_1 \cdot  \left\{ \  {{1}\over {{\bf d}_{\,1\,,\ 2}^{\,n\,-\ 2 }  }}   \  \right\}\,\cdot\,\left\{ \,  {1\over { {\bf d}_{1\,,\,2}^{\,2}  }}  \,\right\}\,\cdot\, \left[\ \int_0^{\rho_{\,\nu}}\!
    {{  \lambda_{\,1}^2\,\cdot\,r^{n\,+\,1}\,\cdot\,dr  }\over {\ (\,\lambda_{\,1}^2 \ + \  r^{\,2})^{{n\,+\,4 }\over 2} \  }}  \ \right]\\[0.2in]
    & = & C_2 \cdot  \left\{ \   {{1}\over {{\bf d}_{\,1\,,\ 2}^{\,n\,-\ 2 }  }}   \  \right\} \,*\,O\,\left(\, {\bar\lambda}_{\ \flat}^{\,2\,\gamma} \ \right)\ .
       \end{eqnarray*}
       \hspace*{1.5in}\ \ \ \ \ \ \ [ \ $\uparrow$ \ \  similar to (\,A.8.10\,)\ ]\\[0.1in]
     Likewise, for the second term in (\,A.8.16\,)\\[0.1in]
      (\,A.8.18\,)
 \begin{eqnarray*}
&   & \left\{ \   {{1}\over {{\bf d}_{\,1\,,\ 2}^{\,n\,-\ 2 }  }}   \  \right\}\,\cdot\,\left\{ \,  {1\over { {\bf d}_{1\,,\,2}^{\,2}  }}\right\}\,\cdot\, \left[\ \int_0^{\rho_{\,\nu}}\!
\left[   \    {{\lambda_{\,1}^2}\over {(\,\lambda^2_1 \,+\, r^2\,)^{\,2}}}
\cdot {{ r^2}\over {(\,\lambda_{\,1}^2 \ + \  r^{\,2})^{{n }\over 2} }}\right]\,\cdot\,r^{n\,-\,1}\,\cdot\, {{r^2}\over {\lambda_{\,1}\cdot \lambda_{\,2}}}\ dr   \ \right]\\[0.2in]
 & \le &C_3 \cdot  \left\{ \  {{1}\over {{\bf d}_{\,1\,,\ 2}^{\,n\,-\ 2 }  }}   \  \right\}\,\cdot\,\left\{ \,  {1\over { {\bf d}_{1\,,\,2}^{\,2}  }}  \,\right\}\,\cdot\, \left[\ \int_0^{\rho_{\,\nu}}\!
    {{   \,r^{n\,+\,3}\,\cdot\,dr  }\over {\ (\,\lambda_{\,1}^2 \ + \  r^{\,2})^{{n\,+\,4 }\over 2} \  }}  \ \right] \ \ \ \ [ \ {\mbox{similar  \ \ to \ \ (\,A.3.19\,)}}  \ ]\\[0.2in]
& = &   \left\{  \  {{1}\over {{\bf d}_{\,1\,,\ 2}^{\,n\,-\ 2 }  }}   \  \right\} \,\,\times\, O\,\left(\, {\bar\lambda}_{\ \flat}^{\,2\,\gamma\ - \ o_{\,\lambda}\,(\,1)} \,\right) \ \ \ \ \ \ \ \ \ \ \ \  \ \left( \ \int {{d\,R}\over {R}} \ = \ \ln\,R\ + \ C  \ \right) \ .\\
  \end{eqnarray*}
  Here
  $$
  o_{\,\lambda}\,(\,1) \ \to \ 0^{\,+} \ \ \ \ {\mbox{as}} \ \ \ \ {\bar\lambda}_{\ \flat} \ \to \ 0^{\,+}\ .
  $$
Via symmetric,\\[0.1in]
 (\,A.8.19\,)
 \begin{eqnarray*}
& \  &   \int_{B_{\,\xi_{\,1}}\,(\,{\rho_{\,\nu}})}  {\bf[} \ \cdot\,\cdot \cdot\  |_{\,y} \  {\bf ]}_{\, (\,A.8.14\,)}   \times \, {1\over { {\bf d}_{1\,,\,2}^{\,2}  }}\cdot {{\,2\,(\,y\ - \ \xi_{\,1})\,\cdot\,(\,\xi_{\,1}\ - \ \xi_{\,2})\, }\over
{\lambda_{\,1}\,\cdot\,\,\lambda_{\,2}}} \ dy   \\[0.2in]
&  =  &   \ C_2    \left(\ {1\over {{\bf d}_{1\,,\,2}^{\,n} }}   \right) \cdot\left(\,{{\lambda_{\,1}^{{\,n\,-\,2\,}\over 2}} \over {\lambda_{\,1}\,\cdot\,\lambda_{\,2}}}  \right)\,*\,\\[0.2in]
& \ & \ \
\,\times \int_{B_{\,o}\,(\,{\rho_{\,\nu}})}\!\!\left( {{\lambda_{\,1}}\over {\,\lambda_{\,1}^2 \,+\, \Vert\, {\bar y} \, \,\Vert^{\,2}  }}\ \right)^{\!\!2}
\left\{  \,
  {{ (\,\lambda_{\,1}^2\, -\, \Vert\, {\bar y} \, \,\Vert^{\,2})}\over {\ (\,\lambda_{\,1}^2 \,+\, \Vert\, {\bar y} \, \,\Vert^{\,2})^{{n }\over 2}\  }}\right\}\,\cdot\,  \left[\   \, {\bar y} \, \,\cdot\,\left(\ \xi_{\,1} \, - \, \xi_{\,2} \right)\ \right]   \ d \,\bar y \ \ = \ \ 0 \
\\[0.1in]
 & \ & \hspace*{1in} \ \ \ \ \ \ \  \ \uparrow \ \  \  {\mbox{ symmetric  }} \ \ \ \  \ \ \ \uparrow \ \  \ \ \ \ \ \ \ \ \ \ \ \ \ \ \ \ \ \   \uparrow \ \  {\mbox{anti\,-\,symmetric}}\ \\[0.1in]
 & \  &  \ \ \ \ \ \ \ \  \ \ \  \ \ \ \ \  \ \  \left\{\ \bar y \ = \ y \ - \ \xi_{\,1}\ , \ \ \ \ \  [\ -\, {\bar y} \ ]\, \,\cdot\,\left(\ \, \xi_{\,1} \, - \, \xi_{\,2} \right) \ = \ -\,\left\{\  {\bar y} \, \,\cdot\,\left(\ \, \xi_{\,1} \, - \, \xi_{\,2} \right)\ \right\} \    \right\} \ .
    \end{eqnarray*}
    Thus we obtain
     \begin{eqnarray*}
(\,A.8.20\,)
 & \ & -\,n\,(\,n\,-\,2)\,\cdot\,\int_{B_{\,\xi_{\,1}} (\,{\rho_{\,\nu}}\,)} \left\{ \, \left(\ \,{{n\,+\,2}\over {\,n\,-\ 2\,}} \,\right)\,\cdot\,V_1^{4\over {n\,-\,2}}\,\cdot\,\left[\   V_2   \ \right]  \  \right\}\,\cdot\, \left(\,\lambda_{\,1}\cdot { { \partial \,V_1} \over {\partial\, \lambda_{\,1}}}\ \right)\\[0.2in]
&  & \hspace*{-1in} \ = (\,+\,)\ n\,(\,n\,-\,2) \!\cdot\!\left(  {{n\,-\,2}\over { 2}}  \right)  \left(\ {{n\,+\,2}\over {\,n\,-\ 2\,}}  \right)\,\cdot\,\int_{B_{\,\xi_{\,1}}\,(\,{\rho_{\,\nu}})} \!\!\!  {\bf[} \ \cdot\,\cdot \cdot\  |_{\,y} \  {\bf ]}_{\, (\,A.8.14\,)}  \,{\bf *}\ \ \  [ \ \leftarrow \ \ {\mbox{see \ \ (\,A.8.13\,)}}\ ]\\[0.2in]
& \ & \hspace*{3in}{\bf *}\,\left\{ \  1 \  + O\,\left(\, {\bar\lambda}_{\ \flat}^{\,2\,\gamma\ - \ o_{\,\lambda}\,(\,1)} \,\right)   \ \right\}\ + \ \\[0.2in]
& \ & \ \ \ \  \ \ +\ \  \left( \  \int_{B_{\,\xi_{\,1}}\,(\,{\rho_{\,\nu}})}   {\bf[} \ \cdot\,\cdot \cdot\  |_{\,y} \  {\bf ]}_{\, (\,A.8.14\,)}    \cdot | \  {\bf T}\ |^{\,2} \ \right)\ .
    \end{eqnarray*}


\newpage

 {\bf \S\,A\,8\,.\,\,e\,.} \ \
{\bf Estimate of second order teams $\,{\bf T }^2\,$ in}\, (\,A.8.13\,)\,. \\[0.1in]
For the second order terms [\,in (\,A.8.13\,)\,]\,:\\[0.1in]
(\,A.8.21\,)
\begin{eqnarray*}
{\bf T }^2 & =  &   \left[\ {1\over {{\bf d}^{2} }}\,\cdot\,\left(\ {{\lambda_{\,2}}\over {\lambda_{\,1}}}  \right) \ \right]^{\,2} \ + \ \left[\ {1\over {{\bf d}^{2} }}\,\cdot\,{{\Vert\,\,y\ - \ \xi_{\,1}\,\Vert^2}\over
{\lambda_{\,1}\,\cdot\,\,\lambda_{\,2}}}    \ \right]^{\,2} \ + \  \left[\  {1\over {{\bf d}^{2} }}\,\cdot\,{{\,2\,(\,y\ - \ \xi_{\,1})\cdot(\,\xi_{\,1}\ - \ \xi_{\,2})\, }\over
{\lambda_{\,1}\,\cdot\,\,\lambda_{\,2}}} \ \right]^{\,2} \ + \ \\[0.2in]
& \ & \ \ \ \ \ + \ 2\,\cdot\, \left[\ {1\over {{\bf d}^{2} }}\,\cdot\,\left(\ {{\lambda_{\,2}}\over {\lambda_{\,1}}}  \right) \ \right]\,\cdot\,\left[\ {1\over {{\bf d}^{2} }}\,\cdot\,{{\Vert\,\,y\ - \ \xi_{\,1}\,\Vert^2}\over
{\lambda_{\,1}\,\cdot\,\,\lambda_{\,2}}}    \ \right] \\[0.2in]
& \ & \ \ \ \ \ \ \ \ \  + \ 2\,\cdot\, \left[\ {1\over {{\bf d}^{2} }}\,\cdot\,\left(\ {{\lambda_{\,2}}\over {\lambda_{\,1}}}  \right) \ \right]\,\cdot\,\left[\  {1\over {{\bf d}^{2} }}\,\cdot\,{{\,2\,(\,y\ - \ \xi_{\,1})\cdot(\,\xi_{\,1}\ - \ \xi_{\,2})\, }\over
{\lambda_{\,1}\,\cdot\,\,\lambda_{\,2}}} \ \right] \\[0.2in]
& \ & \ \ \ \ \ \ \ \ \ \  \ \ \ \ + \ 2\,\cdot\,\left[\ {1\over {{\bf d}^{2} }}\,\cdot\,{{\Vert\,\,y\ - \ \xi_{\,1}\,\Vert^2}\over
{\lambda_{\,1}\,\cdot\,\,\lambda_{\,2}}}    \ \right]\,\cdot\,\left[\  {1\over {{\bf d}^{2} }}\,\cdot\,{{\,2\,(\,y\ - \ \xi_{\,1})\cdot(\,\xi_{\,1}\ - \ \xi_{\,2})\, }\over
{\lambda_{\,1}\,\cdot\,\,\lambda_{\,2}}} \ \right] \ .
 \end{eqnarray*}
 Note that
\begin{eqnarray*}
  &  \  &    \bigg\vert \ 2\,\cdot\, \left[\ {1\over {{\bf d}^{2} }}\,\cdot\,\left(\ {{\lambda_{\,2}}\over {\lambda_{\,1}}}  \right) \ \right]\,\cdot\,\left[\ {1\over {{\bf d}^{2} }}\,\cdot\,{{\Vert\,\,y\ - \ \xi_{\,1}\,\Vert^2}\over
{\lambda_{\,1}\,\cdot\,\,\lambda_{\,2}}}    \ \right] \ \bigg\vert \\[0.2in]
& \le &  \left[\ {1\over {{\bf d}^{2} }}\,\cdot\,\left(\ {{\lambda_{\,2}}\over {\lambda_{\,1}}}  \right) \ \right]^{\,2} \ +  \ \left[\ {1\over {{\bf d}^{2} }}\,\cdot\,{{\Vert\,\,y\ - \ \xi_{\,1}\,\Vert^2}\over
{\lambda_{\,1}\,\cdot\,\,\lambda_{\,2}}}    \ \right]^{\,2} \ \ , \\[0.5in]
  &  \  &    \bigg\vert \ 2\,\cdot\, \left[\ {1\over {{\bf d}^{2} }}\,\cdot\,\left(\ {{\lambda_{\,2}}\over {\lambda_{\,1}}}  \right) \ \right]\,\cdot\,\left[\  {1\over {{\bf d}^{2} }}\,\cdot\,{{\,2\,(\,y\ - \ \xi_{\,1})\cdot(\,\xi_{\,1}\ - \ \xi_{\,2})\, }\over
{\lambda_{\,1}\,\cdot\,\,\lambda_{\,2}}} \ \right] \ \bigg\vert \\[0.2in]
& \le &  \left[\ {1\over {{\bf d}^{2} }}\,\cdot\,\left(\ {{\lambda_{\,2}}\over {\lambda_{\,1}}}  \right) \ \right]^{\,2} \ +  \ \left[\  {1\over {{\bf d}^{2} }}\,\cdot\,{{\,2\,(\,y\ - \ \xi_{\,1})\cdot(\,\xi_{\,1}\ - \ \xi_{\,2})\, }\over
{\lambda_{\,1}\,\cdot\,\,\lambda_{\,2}}} \ \right]^{\,2} \ \ , \\[0.5in]
 {\mbox{and}} \ \ \ \ \ \ \  \ \ \ \ \ \ \    &  \  &    \bigg\vert \ 2\,\cdot\,\left[\ {1\over {{\bf d}^{2} }}\,\cdot\,{{\Vert\,\,y\ - \ \xi_{\,1}\,\Vert^2}\over
{\lambda_{\,1}\,\cdot\,\,\lambda_{\,2}}}    \ \right]\,\cdot\,\left[\  {1\over {{\bf d}^{2} }}\,\cdot\,{{\,2\,(\,y\ - \ \xi_{\,1})\cdot(\,\xi_{\,1}\ - \ \xi_{\,2})\, }\over
{\lambda_{\,1}\,\cdot\,\,\lambda_{\,2}}} \ \right]   \ \bigg\vert \\[0.2in]
& \le &  \left[\ {1\over {{\bf d}^{2} }}\,\cdot\,{{\Vert\,\,y\ - \ \xi_{\,1}\,\Vert^2}\over
{\lambda_{\,1}\,\cdot\,\,\lambda_{\,2}}}    \ \right]^{\,2} \ + \  \left[\  {1\over {{\bf d}^{2} }}\,\cdot\,{{\,2\,(\,y\ - \ \xi_{\,1})\cdot(\,\xi_{\,1}\ - \ \xi_{\,2})\, }\over
{\lambda_{\,1}\,\cdot\,\,\lambda_{\,2}}} \ \right]^{\,2}  \ .  \ \ \ \ \ \ \   \ \ \ \ \ \ \   \ \ \ \ \ \ \  \\
 \end{eqnarray*}

\newpage

Thus it is enough to estimate the square terms\,;
$$
\left[\ {1\over {{\bf d}^{2} }}\,\cdot\,\left(\ {{\lambda_{\,2}}\over {\lambda_{\,1}}}  \right) \ \right]^{\,2} \ , \ \ \left[\ {1\over {{\bf d}^{2} }}\,\cdot\,{{\Vert\,\,y\ - \ \xi_{\,1}\,\Vert^2}\over
{\lambda_{\,1}\,\cdot\,\,\lambda_{\,2}}}    \ \right]^{\,2}\ \ \ \ {\mbox{and}} \ \ \ \left[\  {1\over {{\bf d}^{2} }}\,\cdot\,{{\,2\,(\,y\ - \ \xi_{\,1})\cdot(\,\xi_{\,1}\ - \ \xi_{\,2})\, }\over
{\lambda_{\,1}\,\cdot\,\,\lambda_{\,2}}} \ \right]^{\,2}
$$
Let us start with
$$
  \left[\ {1\over {{\bf d}^{2} }}\,\cdot\,\left(\ {{\lambda_{\,2}}\over {\lambda_{\,1}}}  \right) \ \right]^{\,2}  \ = \ O\,\left(\,{\bar\lambda}_{\ \flat}^{\,4\,\gamma}\,\right)\,. \leqno (\,A.8.22\,)
$$
Similar to (\,A.8.16\,) and (\,A.8.18\,)\,,\\[0.1in]
(\,A.8.23\,)
 \begin{eqnarray*}
 & \ &  \int_{B_{\,\xi_{\,1}}\,(\,{\rho_{\,\nu}})}{\bf[} \ \cdot\,\cdot \cdot\ {\bf ]}  \,\,\times\,\left[\   {1\over{ {\bf d}_{1\,,\,2}^{\,2}  }}\cdot {{\Vert\,\,y\ - \ \xi_{\,1}\ \Vert^2}\over
{\lambda_{\,1}\,\cdot\,\,\lambda_{\,2}}} \ \right]^{\,2}\\[0.2in]
& = & C \cdot   \left\{ \   {{1}\over {{\bf d}_{\,1\,,\ 2}^{\,n\,-\ 2 }  }}   \  \right\}\,\cdot\,\left\{ \,  {1\over { {\bf d}_{1\,,\,2}^{\,4}  }}\right\}\,\cdot\, \left[\ \int_0^{\rho_{\,\nu}}\!
\left[   \    {{\lambda_{\,1}^2}\over {(\,\lambda^2_1 \,+\, r^2\,)^{\,2}}}
\cdot {{ \lambda_{\,1}^2}\over {(\,\lambda_{\,1}^2 \ + \  r^{\,2})^{{n }\over 2} }}\right]\,\cdot\,r^{n\,-\,1}\,\cdot\, {{r^4}\over {\lambda_{\,1}^2\cdot \lambda_{\,2}^2}}\ dr   \ \right] \\[0.2in]
& \ &  \!\!\!\!|- \ C \cdot  \left\{ \   {{1}\over {{\bf d}_{\,1\,,\ 2}^{\,n\,-\ 2 }  }}   \  \right\}\,\cdot\,\left\{ \,  {1\over { {\bf d}_{1\,,\,2}^{\,4}  }}\right\}\,\cdot\, \left[\ \int_0^{\rho_{\,\nu}}\!
\left[   \    {{  \lambda_{\,1}^2    }\over {(\,\lambda^2_1 \,+\, r^2\,)^{\,2}}}
\cdot {{ r^2   }\over {(\,\lambda_{\,1}^2 \ + \  r^{\,2})^{{n }\over 2} }}\right]\,\cdot\,r^{n\,-\,1}\,\cdot\, {{r^4}\over {\lambda_{\,1}^2\cdot \lambda_{\,2}^2}}\ dr   \ \right] \ . \\
  \end{eqnarray*}
  We continue with \\[0.1in]
  (\,A.8.24\,)
\begin{eqnarray*}
& \ &  \left\{ \   {{1}\over {{\bf d}_{\,1\,,\ 2}^{\,n\,-\ 2 }  }}   \  \right\}\,\cdot\,\left\{ \,  {1\over { {\bf d}_{1\,,\,2}^{\,4}  }}\right\}\,\cdot\, \left[\ \int_0^{\rho_{\,\nu}}\!
\left[   \    {{\lambda_{\,1}^2}\over {(\,\lambda^2_1 \,+\, r^2\,)^{\,2}}}
\cdot {{ \lambda_{\,1}^2}\over {(\,\lambda_{\,1}^2 \ + \  r^{\,2})^{{n }\over 2} }}\right]\,\cdot\,r^{n\,-\,1}\,\cdot\, {{r^4}\over {\lambda_{\,1}^2\cdot \lambda_{\,2}^2}}\ dr   \ \right]\\[0.2in]
& = & O\,(\,1)\cdot  {{1}\over {{\bf d}_{\,1\,,\ 2}^{\,n\,-\ 2 }  }} \,\cdot\, {{1}\over {{\bf d}_{\,1\,,\ 2}^{\,4 }  }} \ \cdot\,\left[\ \int_0^{\rho_{\,\nu}}\!
    {{    r^{n\,+\,3}\,\cdot\,dr  }\over {(\,\lambda_{\,1}^2 \ + \  r^{\,2})^{{n\,+\,4 }\over 2} }}  \ \right] \ \ \ \ \ \ \ \ \ \ \  \ \ \ \ \ \ \ \ \ \ \  [ \ {\mbox{cf.}} \ \ (\,A.3.18\,) \ ]\\[0.2in]
& = & O\,(\,1)\cdot  {{1}\over {{\bf d}_{\,1\,,\ 2}^{\,n\,-\ 2 }  }} \,\,\cdot\,O\,\left(\, {\bar\lambda}_{\ \flat}^{\,4\,\gamma\ - \ o_{\,+}\,(\,1\,)} \,\right)\ .
\end{eqnarray*}
As for the second term in (\,A.8.23\,)\,,\,\\[0.1in]
(\,A.8.25\,)
\begin{eqnarray*}
& \ &  \left\{ \   {{1}\over {{\bf d}_{\,1\,,\ 2}^{\,n\,-\ 2 }  }}   \  \right\}\,\cdot\,\left\{ \,  {1\over { {\bf d}_{1\,,\,2}^{\,4}  }}\right\}\,\cdot\, \left[\ \int_0^{\rho_{\,\nu}}\!
\left[   \    {{\lambda_{\,1}^2}\over {(\,\lambda^2_1 \,+\, r^2\,)^{\,2}}}
\cdot {{ r^2}\over {(\,\lambda_{\,1}^2 \ + \  r^{\,2})^{{n }\over 2} }}\right]\,\cdot\,r^{n\,-\,1}\,\cdot\, {{r^4}\over {\lambda_{\,1}^2\cdot \lambda_{\,2}^2}}\ dr   \ \right]\\[0.2in]
& = & O\,(\,1)\cdot  {{1}\over {{\bf d}_{\,1\,,\ 2}^{\,n\,-\ 2 }  }} \,\cdot\, {{1}\over {{\bf d}_{\,1\,,\ 2}^{\,4 }  }} \ \cdot\,\left[\ \int_0^{\rho_{\,\nu}}\!
    {{    r^{n\,+\,5}\,\cdot\,dr  }\over {(\,\lambda_{\,1}^2 \ + \  r^{\,2})^{{n\,+\,4 }\over 2} }}  \ \right]\\[0.2in]
    & = & O\,(\,1)\cdot  {{1}\over {\ {\bf d}_{\,1\,,\ 2}^{\,n\,-\ 2 } \  }} \,\,\cdot\,O\,\left(\, {\bar\lambda}_{\ \flat}^{\,2\,\gamma} \,\right)\,\,\cdot\,O\,\left(\, {\bar\lambda}_{\ \flat}^{\ 2\,[\ (\ \gamma \ + \ \nu\ ) \ - \ 1\ ]} \,\right)\ .\\
\end{eqnarray*}
As for the last square term\,:
\begin{eqnarray*}
(\,A.8.26\,) \  \ \ \ \  \ \  \left[\  {1\over {{\bf d}^{2} }}\,\cdot\,{{\,2\,(\,y\ - \ \xi_{\,1})\cdot(\,\xi_{\,1}\ - \ \xi_{\,2})\, }\over
{\lambda_{\,1}\,\cdot\,\,\lambda_{\,2}}} \ \right]^{\,2} &  \le &   {1\over {{\bf d}^{4} }}\,\cdot\,2\,\cdot\,\Vert\,\,y\ - \ \xi_{\,1}\,\Vert^2\,\cdot\, {{\Vert\,\,\xi_{\,1}\ - \ \xi_{\,2}\,\Vert^2}\over
{\lambda_{\,1}^2\,\cdot\,\,\lambda_{\,2}^2}}\\[0.2in]
& \le & C \,\left[\ {1\over {{\bf d}^{2} }}\,\cdot\,  \Vert\,\,y\ - \ \xi_{\,1}\,\Vert^2 \ \right] \cdot   {{1}\over
{\lambda_{\,1}\,\cdot\,\,\lambda_{\,2}}}  \ .
 \end{eqnarray*}
Thus a similar calculation as in (\,A.8.16\,)\,--\, (\,A.8.18\,) leads to
 \begin{eqnarray*}
(\,A.8.27\,) \  \ \ \  \  \ \ \  & \ & \bigg\vert \  \int_{B_{\,\xi_{\,1}}\,(\,{\rho_{\,\nu}})}{\bf[} \ \cdot\,\cdot \cdot\ {\bf ]}  \,\,\times\,\left[\  {1\over {{\bf d}^{2} }}\,\cdot\,{{\,2\,(\,y\ - \ \xi_{\,1})\cdot(\,\xi_{\,1}\ - \ \xi_{\,2})\, }\over
{\lambda_{\,1}\,\cdot\,\,\lambda_{\,2}}} \ \right]^{\,2} \ \bigg\vert  \  \ \ \  \  \ \ \  \  \ \ \  \  \ \ \  \\[0.2in]
& = &  \left\{  \  {{1}\over {{\bf d}_{\,1\,,\ 2}^{\,n\,-\ 2 }  }}   \  \right\} \,\,\times\, O\,\left(\, {\bar\lambda}_{\ \flat}^{\,2\,\gamma\ - \ o_{\,\lambda}\,(\,1)} \,\right) \ .
  \end{eqnarray*}
Combining  (\,A.8.24\,)\,,\, (\,A.8.25\,) and (\,A.8.27\,)\,,\, we obtain
   \begin{eqnarray*}
 (\,A.8.28\,)   & \ & \bigg\vert \   \int_{B_{\,\xi_{\,1}}\,(\,{\rho_{\,\nu}})}   {\bf[} \ \cdot\,\cdot \cdot\ {\bf ]}_{\,y}    \cdot | \  {\bf T}\ |^{\,2} \ \bigg\vert\ = \  \left\{  \  {{1}\over {{\bf d}_{\,1\,,\ 2}^{\,n\,-\ 2 }  }}   \  \right\} \,\,\times\,O\,\left(\, {\bar\lambda}_{\ \flat}^{\,2\,\gamma\ - \ o_{\,\lambda}\,(\,1)} \,\right)\\[0.2in]
\Longrightarrow  & \ & -\,n\,(\,n\,-\,2)\,\cdot\,\int_{B_{\,\xi_{\,1}} (\,{\rho_{\,\nu}}\,)} \left\{ \, \left(\ \,{{n\,+\,2}\over {\,n\,-\ 2\,}} \,\right)\,\cdot\,V_1^{4\over {n\,-\,2}}\,\cdot\,\left[\   V_2   \ \right]  \  \right\}\,\cdot\, \left(\,\lambda_{\,1}\cdot { { \partial \,V_1} \over {\partial\, \lambda_{\,1}}}\ \right)\\[0.2in]
& \ & \hspace*{-0.75in}= \  +\,n\,(\,n\,-\,2) \cdot\left( \,{{n\,-\,2}\over { 2}} \,\right)\cdot  \left(\ \,{{n\,+\,2}\over {\,n\,-\ 2\,}} \,\right)\,\cdot\,\int_{B_{\,\xi_{\,1}}\,(\,{\rho_{\,\nu}})}   {\bf[} \ \cdot\,\cdot \cdot\ {\bf ]}_{\,y}  \, \times\,\left\{ \  1 \  + \,  O\,\left(\, {\bar\lambda}_{\ \flat}^{\,2\,\gamma\ - \ o_{\,\lambda}\,(\,1)} \,\right)   \ \right\}\ .
    \end{eqnarray*}

\newpage

 {\bf \S\,A\,.8\,.\,f\,.} \ \
{\bf Estimate of}  $\displaystyle{\ \ \ {\bf A}^2_{\,(\,A.8.1\,)} \ = \ \left( \ \ \sum_{l\,=\,2}^\flat \ \int_{B_{\xi_l} (\,{\rho_{\,\nu}})}\ \right)\ . }$  \\[0.1in]
Inside the ball $\,B_{\,\xi_{\,2}} (\,{\rho_{\,\nu}}) \,$,
\begin{eqnarray*}
 (\,A.8.29\,)    \ \ \ \ \ \ & \ &   \bigg\vert\  \int_{B_{\xi_{\,2}} (\,{\rho_{\,\nu}})}     \left\{  \  \sum_{l\,=\,1}^\flat \, V_{\,l}^{{\,n\,+\,2\,}\over {n\,-\,2}}   \
    - \ \left[\  W_{\,\,\flat}  \ \right]^{{\,n\,+\,2\,}\over {n\,-\,2}} \  \right\}\,\cdot\,\left(\ \lambda_{\,1}\cdot { { \partial \,V_1} \over {\partial\, \lambda_{\,1}}}\ \right) \ \bigg\vert \\[0.2in]
    & \le & C\,\cdot\, \int_{B_{\xi_{\,2}} (\,{\rho_{\,\nu}})}  V_2^{4\over {n\,-\,2}}\,\cdot\,\left[\   V_1  \ + \ V_3 \ + \ \,\cdot\,\cdot\,\cdot\,\ + \ V_{\,\flat}  \ \right]   \,\,\cdot\, \ V_1  \ \,+ \ {\bf{error}}_{\, (\,A.8.29\,)}\ , \ \ \ \ \ \  \ \ \ \ \ \   \ \ \ \ \ \
\end{eqnarray*}
where
$$
|\,{\bf{error}}_{\, (\,A.8.29\,)} \,| \ \le \  C\,\cdot\, \int_{B_{\xi_{\,2}} (\,{\rho_{\,\nu}})}\left(   V_1^{{\,n\,+\,2\,}\over {n\,-\,2}}  \ + \  V_3^{{\,n\,+\,2\,}\over {n\,-\,2}} \ + \ \cdot\,\cdot\,\cdot\,\ + \ V_{\,\flat}^{{\,n\,+\,2\,}\over {n\,-\,2}}\  \right)\,\cdot\,  \ V_1  \ ,
$$
which can be estimated as in \S\,A\,\,8\,.\,b\,,\, leading to
$$
|\,{\bf{error}}_{\, (\,A.8.29\,)} \,|  \ = \ O\,\left( \  {\bar\lambda}_{\ \flat}^{n\,\cdot\,\ \gamma}\ \right)  \,\times  O\,\left(\, {\bar\lambda}_{\ \flat}^{\,2\,[\ (\,\gamma\ + \ \nu\,) \ - \ 1 \ ] } \ \right)\ .  \leqno (\,A.8.30\,)
$$





\vspace*{0.3in}

 {\bf \S\,A\,.8\,.\,g\,.} \ \ {\bf Estimate of the first integral in the right hand side of} \  (\,A.8.29\,) \,.\\[0.1in]
In $\,B_{\xi_{\,2}}(\,\rho_\nu\,)\,$,\,\\[0.1in]
(\,A.8.31\,)
\begin{eqnarray*}
& \ &  [\ V_1   \ + \ V_3   \ + \ \cdot\,\cdot\,\cdot\,\ + \ V_{\,\flat}\ ]_{\,y} \\[0.2in]
& \le & C\,\cdot\,\left[\ \left(\ {\lambda_{\,1}\over {\lambda^2_1 \ + \ \Vert\ \xi_{\,2}\,-\,\xi_{\,1}\,\Vert^2 }}\  \right)^{{\,n\,-\,2\,}\over 2} \ + \ \left(\ {\lambda_3\over {\lambda^2_3 \ + \ \Vert\ \xi_{\,2}\,-\,\xi_3\,\Vert^2 }}\  \right)^{{\,n\,-\,2\,}\over 2} \ \ + \ \cdot\,\cdot\,\cdot\,\ + \right.\\[0.15in]
& \ & \hspace*{3.5in} \left. \ + \  \left(\ {{\bar\lambda}_{\ \flat}\over {\lambda^2_\flat \ + \ \Vert\ \xi_{\,2}\ - \ \xi_\flat\,\Vert^2  }}\  \right)^{{\,n\,-\,2\,}\over 2} \ \right]\\[0.2in]
& \le & C_1\,\cdot\,\left(\ {1\over {{\bar\lambda}_{\ \flat}   }}\  \right)^{{\,n\,-\,2\,}\over 2}\,\cdot\,\left(\  \sum_{l\,\not=\,2}\  {1\over { \Vert\ \Xi_{\,2} \ - \ \Xi_{\,l}\,\Vert^{\,n\,-\,2}  }}\  \right) \ \le \ C_2\,\cdot\,{{ {\bar\lambda}_{\ \flat}^{(\,n\,-\,2\,)\,\cdot\,\gamma}}\over { {\bar\lambda}_{\ \flat}^{{\,n\,-\,2\,}\over 2}  }} \ \ \ \ \ \ \ \ \ [\ y\,\in\,B_{\,\xi_{2}}(\,\rho_\nu\,)\ ]\ .\\
\end{eqnarray*}
Likewise,
$$
V_1\,(\,y\,) \ \le \ C_2\,\cdot\,{{ {\bar\lambda}_{\ \flat}^{(\,n\,-\,2)\,\cdot\,\gamma}}\over { {\bar\lambda}_{\ \flat}^{{\,n\,-\,2\,}\over 2}  }} \ \ \ \ \ \ \ \ \ [\ y\,\in\,B_{\,\xi_{2}}(\,\rho_\nu\,)\ ]\ .
$$
It follow that
$$
\left[\   V_1  \ + \ V_3 \ + \ \,\cdot\,\cdot\,\cdot\,\ + \ V_{\,\flat}  \ \right]   \,\,\cdot\,\, V_1  \ \le  \  C \,\cdot\,{{ {\bar\lambda}_{\ \flat}^{\,2\cdot (\,n\,-\,2\,)\,\cdot\,\gamma}}\over { {\bar\lambda}_{\ \flat}^{ n\,-\,2 }  }}  \ \ \ \ \ \ \ \ \mfor \ \  y\,\in\,B_{\,\xi_{2}}(\,\rho_\nu\,) \ . \leqno (\,A.8.32\,)
$$
Moreover,
\begin{eqnarray*}
& \ & \int_{B_{\xi_{\,2}}(\,\rho_{\,\nu}\,)} \left(\ {\lambda_{\,2}\over {\lambda^2_2 \ + \ \Vert\ y\,-\,\xi_{\,2}\,\Vert^2 }}\  \right)^{\,2} \ d\,y\\[0.3in]
& \le & C\,\cdot\,\int_0^{\rho_{\nu}}  \left(\ {\lambda_{\,2}\over {\lambda^2_2 \ + \ r^2 }}\  \right)^{\,2} \ r^{n\,-\,1}\,\cdot\,dr \ \ \ \ \ [\ \leftarrow \ \ {\mbox{estimated \ \ as \ \ in \ \ (\,A.3.19\,)}}\ ]\\[0.3in]
& = & C\,\cdot\,\int_0^{\rho_{\nu}} {{\lambda_{\,2}^{n\,-\,2 }\,\cdot\,\left(\  {r\over {\lambda_{\,2}}}\  \right)^{\! n\,-\,1 }\,\cdot\,{{dr}\over {\lambda_{\,2}}}  }\over { \left[\  1 \ + \ R^2 \ \right]^2 }} \\[0.3in]
& \le & C_1\,\cdot\,\lambda_{\,2}^{n\,-\,2}\,\cdot\,\int_0^{{\rho_{\,\nu}}\over \lambda_{\,2}}\ \  {{R^{n\,-\,1} \ dR }\over { R^4 }} \ \le \ C_2\,\cdot\,\lambda_{\,2}^{n\,-\,2 }\,\cdot\,\left(\ {{\rho_{\,\nu}}\over \lambda_{\,2}} \right)^{\!\!n\,-\,4}
\end{eqnarray*}
Thus\\[0.15in]
(\,A.8.33\,)
\begin{eqnarray*}
& \ &   \int_{B_{\xi_{\,2}} (\,{\rho_{\,\nu}})}  V_2^{4\over {n\,-\,2}}\,\cdot\,\left[\   V_1  \ + \ V_3 \ + \ \,\cdot\,\,\cdot\,\cdot \ + \ V_{\,\flat}  \ \right]   \, \cdot\,   V_1 \ \le  \  C \,\cdot\, {\bar\lambda}_{\ \flat}^{2\,\cdot\, (\,n\,-\,2\,)\,\cdot\,\gamma}   \,\cdot\,\left(\ {{1}\over {\lambda_{\flat}^{1\,-\ \nu} }}\  \right)^{\!\!n\,-\,4}\\[0.2in]
& \ & \hspace{3in} \ \ \ \ \ \ \ \ \ \le \ C_1\,\cdot\,{\bar\lambda}_{\ \flat}^{\  n\,\cdot\, \gamma} \\[0.2in]
 \Longrightarrow & \ & |\, {\bf A}^2_{\,(\,A.8.1\,)}\,| \ \le \  C_3  \,\cdot\,{\bar\lambda}_{\ \flat}^{\  n\,\cdot\, \gamma}\,\cdot\,(\,\flat\ - \ 1\,) \ = \ 0\,\left(\ \  {\bar\lambda}_{\ \flat}^{\  n\,\cdot\,\gamma\ - \ \sigma} \ \right) \ .
\end{eqnarray*}
Here we use
$$
\gamma \ + \ \nu \ >  \ 1\,.
$$



\newpage

 {\bf \S\,A\,8\,.\,h\,.} \ \
{\bf Estimate of}  $\displaystyle{\ \ \ {\bf A}^3_{(\,A.8.1\,) } \ \left( \ = \   \int_{\R^n \,\setminus \ \left(\, \cup\ B_{\xi_l} (\,{\rho_{\,\nu}}) \, \right)}\ \right)\ . }$  \ As in {\bf \S\,A 4\,.\,j}\,,\, we have \\[0.1in]
\begin{eqnarray*}
(\,A.8.34\,)\ \ \ \ \ \ \ \ \ & \ &   \bigg\vert \  \int_{\R^n \,\setminus \ \left(\, \cup\ B_{\xi_l} (\,{\rho_{\,\nu}}) \, \right)}  \left\{  \  \sum_{l\,=\,1}^\flat \, V_{\,l}^{{\,n\,+\,2\,}\over {n\,-\,2}}   \
    - \ \left[\  W_{\,\,\flat}  \ \right]^{{\,n\,+\,2\,}\over {n\,-\,2}} \  \right\}\,\cdot\,\left(\ \lambda_{\,1}\cdot { { \partial \,V_1} \over {\partial\, \lambda_{\,1}}}\ \right) \ \bigg\vert \ \ \ \ \ \ \ \ \ \ \ \ \ \   \\[0.2in]
    & \le &  C\,\cdot\,  \int_{\R^n \,\setminus \ \left(\, \cup\ B_{\xi_l} (\,{\rho_{\,\nu}}) \, \right)}  \left\{  \  \left[\  W_{\,\,\flat}  \ \right]^{{\,n\,+\,2\,}\over {n\,-\,2}} \ - \  \sum_{l\,=\,1}^\flat \, V_{\,l}^{{\,n\,+\,2\,}\over {n\,-\,2}}      \  \right\}\,\cdot\,V_1   \ .  \\[0.2in]
\end{eqnarray*}
Note that
$$
n \ \ge \ 6 \ \ \Longrightarrow \  \ {1\over 2}\,\cdot\,{{n\,+\,2}\over {\,n\,-\ 2\,}} \ \le \ 1\,.
$$
It follows that (\,similar to {\bf \S\,A\,7.\,h}\,)\,.\\[0.1in]
(\,A.8.35\,)
\begin{eqnarray*}
\left[\  W_{\,\,\flat}  \ \right]^{\,{{n\,+\,2}\over {\,n\,-\ 2\,}} }  & = &   \left[\  V_1\ + \  V_2  \ + \,\cdot\,\cdot\,\cdot\, + \ V_{\,\flat}  \ \right]^{\,{{1\over 2} \,\cdot \,{{n\,+\,2}\over {\,n\,-\ 2\,}}}} \ \cdot\, \, \left[\  V_1\ + \  V_2  \ + \,\cdot\,\cdot\,\cdot\, + \ V_{\,\flat}  \ \right]^{\,{{1\over 2} \,\cdot \,{{n\,+\,2}\over {\,n\,-\ 2\,}}}} \\[0.2in]
& \le &    \left[\  V_1^{\,{1\over 2}\,\cdot\,{{n\,+\,2}\over {\,n\,-\ 2\,}}}\ + \  V_2^{\,{1\over 2}\,\cdot\, {{n\,+\,2}\over {\,n\,-\ 2\,}}}  \ + \,\cdot\,\cdot\,\cdot\, + \ V_{\,\flat}^{\,{1\over 2}\,\cdot\,{{n\,+\,2}\over {\,n\,-\ 2\,}}}  \ \right]\,*\\[0.2in]
& \ & \hspace*{1in} *\,\left[\  V_1^{\,{1\over 2}\,\cdot\,{{n\,+\,2}\over {\,n\,-\ 2\,}}}\ + \  V_2^{\,{1\over 2}\,\cdot\,{{n\,+\,2}\over {\,n\,-\ 2\,}}} \ + \,\cdot\,\cdot\,\cdot\, + \ V_{\,\flat}^{\,{1\over 2}\,\cdot\, {{n\,+\,2}\over {\,n\,-\ 2\,}}}  \ \right]\\[0.2in]
& = & \left[\  V_1^{{\,n\,+\,2\,}\over {n\,-\,2}} \ + \  V_2^{{\,n\,+\,2\,}\over {n\,-\,2}} \ +\,\cdot\,\cdot\,\cdot\, + \ V_{\,\flat}^{{\,n\,+\,2\,}\over {n\,-\,2}} \ \right] \\[0.2in]
& \ & \ \ \ \ \  + \ V_1^{\,{1\over 2}\,\cdot\,{{n\,+\,2}\over {\,n\,-\ 2\,}}}\,\cdot\,\left[\    V_2^{\,{1\over 2}\,\cdot\,{{n\,+\,2}\over {\,n\,-\ 2\,}}}  \ + \,\cdot\,\cdot\,\cdot\, + \ V_{\,\flat}^{\,{1\over 2}\,\cdot\,{{n\,+\,2}\over {\,n\,-\ 2\,}}}  \ \right]\\[0.1in]
& \  &  \ \ \ \  \ \ + \  V_2^{\,{1\over 2}\,\cdot\,{{n\,+\,2}\over {\,n\,-\ 2\,}}}\,\cdot\,\left[\    V_1^{\,{1\over 2}\,\cdot\,{{n\,+\,2}\over {\,n\,-\ 2\,}}} \ + \,\cdot\,\cdot\,\cdot\, + \ V_{\,\flat}^{\,{1\over 2}\,\cdot\, {{n\,+\,2}\over {\,n\,-\ 2\,}}}  \ \right]\\[0.1in]
& \  &  \ \ \ \ \ \ \ \ \ \ \ \ \ \ \,\cdot\,\ \\[0.01in]
& \  &    \ \ \ \ \ \ \ \ \ \ \ \ \ \ \, \cdot\,\ \\[0.01in]
& \  &   \ \ \ \  \ \ + \
  V_{\,\flat}^{\,{1\over 2}\,\cdot\,{{n\,+\,2}\over {\,n\,-\ 2\,}}}\,\cdot\,\left[\    V_1^{\,{1\over 2}\,\cdot\,{{n\,+\,2}\over {\,n\,-\ 2\,}}}  \ + \,\cdot\,\cdot\,\cdot\, + \ V_{\flat\,-\,1}^{\,{1\over 2}\,\cdot\,{{n\,+\,2}\over {\,n\,-\ 2\,}}} \ \right]  \ . \\
\end{eqnarray*}

\newpage

Consider terms like (\,needed in {\bf \S\,A\,9\,}\,)\,\\[0.1in]
(\,A.8.36\,)
\begin{eqnarray*}
& \ & \left[\ V_1\,(\,y\,) \,\cdot\,  V_j\,(\,y\,)\ \right]^{ \, {{n\,+\,2}\over {\,n\,-\ 2\,}}} \ \ \ \ \ \ \ \ \ \ \ \  \ \ \ \ \ \  \ \ \ \ \ \  \ \ \ \ \ \  \ \ \ \ \ \   \ \ \ (\ j\ \not= \ 1\ ) \\[0.2in]
& = & \left(\ {{\lambda_{\,1}}\over {\lambda_{\,1}^2 \ + \ \Vert\ y\ - \ \xi_{\,1}\,\Vert^2 }}\ \right)^{\!\!{{\,n\,-\,2\,}\over 2} } \,\cdot\,\left(\ {{\lambda_{\,j}}\over {\lambda_{\,j}^2 \ + \ \Vert\ y\ - \ \xi_j\,\Vert^2 }}\ \right)^{\!\!{{n\,+\,2}\over 2}} \\[0.2in]
& \le & C_1\,\cdot\,{1\over {{\bar\lambda}_{\ \flat}^n}}\,\cdot\, \left(\ {{1}\over {1 \ + \ \Vert\ Y\ - \ \Xi_{\,\,1}\,\Vert^2 }}\ \right)^{\!\!{{\,n\,-\,2\,}\over 2} } \,\cdot\,\left(\ {{1}\over {1 \ + \ \Vert\ Y\ - \ \Xi_j\,\Vert^2 }}\ \right)^{\!\!{{n\,+\,2}\over 2}  } \\[0.2in]
& \le & C_2\,\cdot\,{1\over {{\bar\lambda}_{\ \flat}^n}}\,\cdot\, \left(\ {{1}\over {1 \ + \ \Vert\ Y\ - \ \Xi_{\,\,1}\,\Vert }}\ \right)^{\!\!n\,-\,2 }  \,\cdot\,\left(\ {{1}\over {1 \ + \ \Vert\ Y\ - \ \Xi_j\,\Vert  }}\ \right)^{\!\!n\,+\,2 } \\[0.2in]
& \le & C_3\,\cdot\,{1\over {{\bar\lambda}_{\ \flat}^n}}\,\cdot\,{1\over {\Vert\  \Xi_{\,\,1} \ - \ \Xi_j\,\Vert^{\,n\,-\,2\, } }}\,\cdot\,\left[\ \left(\ {{1}\over {1 \ + \ \Vert\ Y\ - \ \Xi_{\,\,1}\,\Vert }}\ \right)^{\!\! n\,+\,2  } \right. \\[0.2in]
& \ & \ \ \ \ \ \ \ \ \ \ \ \ \ \ \ \ \ \ \ \ \  \ \ \  \ \ \  \ \ \  \ \  \ \ \ \ \ \ \ \ \ \ \ \ \ \ \ \   \left.+ \ \left(\ {{1}\over {1 \ + \ \Vert\ Y\ - \ \Xi_j\,\Vert }}\ \right)^{\! n\,+\,2  } \ \  \right]\ .
\end{eqnarray*}
\hspace*{0.5in}\ \ \ \ \ \ $\uparrow$ \ \ combined in the integration factor.\\[0.1in]
\noindent Upon integration, we obtain
\begin{eqnarray*}
(\,A.8.37\,)   & \ &
{1\over {{\bar\lambda}_{\ \flat}^n}}\,\cdot\,  \int_{\,\R^n \setminus \,  \,B_{\xi_{\,1} } (\,{\rho_{\,\nu}})}    \left(\ {{1}\over {1 \ + \ \Vert\ Y\ - \ \Xi_{\,\,1}\,\Vert }}\ \right)^{\!\! (\,n\,+\,2\,)   } \ d\,V_y \ \ \ \ \ \ \      \left[\ Y \ = \ {y\over \lambda_{\,1}} \ \right] \ \ \ \ \ \ \ \ \ \ \\[0.2in]
& \le & C\,\cdot\,  \int_{ {{\rho_{\,\nu}}\over {\lambda_{\,1}}} }^\infty\  {{R^{n\,-\,1} \ dR}\over {\ \, R^{ \,(\,n\,+\,2\,) \, }  \ \, }}\\[0.2in]& \le & C \cdot \left( -\,{1\over 2}\right) \,\cdot\, {1\over {\  R^{\,2   } \ }}\ \ \bigg\vert_{ \,{{\rho_{\,\nu}}\over {\, \lambda_{\,1}\, }} }^{\,\infty}\  = \ {C\over 2}  \,\cdot\,  {\bar\lambda}_{\ \flat}^{ \ 2 \cdot\, (\,1\ - \ \nu\,) }\ .
\end{eqnarray*}

\newpage

Hence\\[0.1in]
(\,A.8.38\,)
\begin{eqnarray*}
& \ &       \int_{\R^n \,\setminus \ \left(\, \cup\ B_{\xi_l} (\,{\rho_{\,\nu}}) \, \right)}    V_1  \,\cdot\, \left[\    V_2^{ {{n\,+\,2}\over {\,n\,-\ 2\,}}}  \ + \,\cdot\,\cdot\,\cdot\, + \ V_{\,\flat}^{ {{n\,+\,2}\over {\,n\,-\ 2\,}}}  \ \right]\\[0.2in]
& \le & C_1 \left(\ \sum_{j\ =\,2}^\flat \ {1\over {\Vert\  \Xi_{\,\,1} \ - \ \Xi_j\,\Vert^{(\,n\,-\,2\,) } }}  \ \right)\,\cdot\,{\bar\lambda}_{\ \flat}^{ \ 2\,\cdot\, (\,1\ - \ \nu\,) }\\[0.2in]
& \le & C_2 \,\cdot\,{\bar\lambda}_{\ \flat}^{ \   (\,n\,-\,2\,)\cdot \gamma \ + \ 2\,\cdot\,(\,1\ - \ \nu\,)  } \\[0.2in]
& \le &  C_3 \,\cdot\,{\bar\lambda}_{\ \flat}^{ \   {{n\,+\,2}\over 2}\,\cdot\,\gamma \ + \ {{\,n\,-\,2\,}\over 2}\,\cdot\,(\,1\ - \ \nu\,)  } \ \ \ \ \ \ \ \ \ \ \ \ \ \  [\,{\mbox{recall \ \ that \ \ }} \gamma \ >  \ (\ 1\ - \ \nu\ ) \ ]\ .
\end{eqnarray*}
Next\,,\, consider terms like [\,cf. (\,A.4.52\,)\,]\\[0.1in]
(\,A.8.39\,)
\begin{eqnarray*}
& \ & V_1 \,\cdot\,V_1^{\,{1\over 2}\,\cdot\,{{n\,+\,2}\over {\,n\,-\ 2\,}}}\,\cdot\, V_j^{\,{1\over 2}\,\cdot\, {{n\,+\,2}\over {\,n\,-\ 2\,}}} \\[0.2in]
& = & \left(\ {{\lambda_{\,1}}\over {\lambda_{\,1}^2 \ + \ \Vert\ y\ - \ \xi_{\,1}\,\Vert^2 }}\ \right)^{\!\!{{\,n\,-\,2\,}\over 2} \ + \ {{n\,+\,4}\over 4} } \,\cdot\,\left(\ {{\lambda_{\,j}}\over {\lambda_{\,j}^2 \ + \ \Vert\ y\ - \ \xi_j\,\Vert^2 }}\ \right)^{\!\!{{n\,+\,2}\over 4}} \ \ \ \ \ \ \ \ (\,j\ \not= \ 1\,)\\[0.2in]
& \le & C\,\cdot\,{1\over {{\bar\lambda}_{\ \flat}^n}}\,\cdot\, \left(\ {{1}\over {1 \ + \ \Vert\ Y\ - \ \Xi_{\,\,1}\,\Vert^2 }}\ \right)^{\!\!{{n\,+\,2}\over 4} \ + \ {{\,n\,-\,2\,}\over 2} }  \,\cdot\,\left(\ {{1}\over {1 \ + \ \Vert\ Y\ - \ \Xi_j\,\Vert^2 }}\ \right)^{\!\!{{n\,+\,2}\over 4}  } \\[0.1in]
& \ & \ \ \ \ \ \ \ \ \ \ \ \ \ \ \ \ \ \ \ \ \ \ \ \ \ \ \ \  \ \ \ \ \ \ \ \ \ \ \ \ \ \left(\ \,{\mbox{note \ \ that}} \ \ {{n\,+\,2}\over 4} \ + \ {{\,n\,-\,2\,}\over 2} \ + \ {{n\,+\,2}\over 4}  \ = \ n\  \right)\\[0.15in]
& \le & C_1\,\cdot\,{1\over {{\bar\lambda}_{\ \flat}^n}}\,\cdot\, \left(\ {{1}\over {1 \ + \ \Vert\ Y\ - \ \Xi_{\,\,1}\,\Vert }}\ \right)^{\!\!{{n\,+\,2}\over 2} \ + \ (\,n\,-\,2\,)  }  \,\cdot\,\left(\ {{1}\over {1 \ + \ \Vert\ Y\ - \ \Xi_j\,\Vert  }}\ \right)^{\!\!{{n\,+\,2}\over 2} }\\[0.2in]
& \le & C_2\,\cdot\,{1\over {{\bar\lambda}_{\ \flat}^n}}\,\cdot\,{1\over {\Vert\  \Xi_{\,\,1} \ - \ \Xi_j\,\Vert^{{{n\,+\,2}\over 2}} }}\,\cdot\,\left[\ \left(\ {{1}\over {1 \ + \ \Vert\ Y\ - \ \Xi_{\,\,1}\,\Vert }}\ \right)^{\!\!{{n\,+\,2}\over 2} \ + \ (\,n\,-\,2\,)  } \right.\\[0.2in]
& \ & \ \ \ \ \ \ \  \ \ \  \ \ \ \ \ \ \ \ \ \ \ \ \ \ \ \ \ \ \ \ \ \ \ \ \ \ \ \  \ \ \ \ \  \ \ \ \ \ \ \ \ \ \ \ \left. \ + \ \left(\ {{1}\over {1 \ + \ \Vert\ Y\ - \ \Xi_j\,\Vert }}\ \right)^{\!\!{{n\,+\,2}\over 2} \ + \ (\,n\,-\,2\,)  }\ \right]\ .
\end{eqnarray*}

\newpage

Upon integration, we obtain
(\,A.8.40\,)
\begin{eqnarray*}
& \ &
{1\over {{\bar\lambda}_{\ \flat}^n}}\,\cdot\,  \int_{\,\R^n \,\setminus \,  \,B_{\xi_{\,l} } (\,{\rho_{\,\nu}})}    \left(\ {{1}\over {1 \ + \ \Vert\ Y\ - \ \Xi_{\,\,1}\,\Vert }}\ \right)^{\!\!{{n\,+\,2}\over 2} \ + \ (\,n\,-\,2\,) }\ d\,V_y \ \ \ \ \ \ \ \ \ \     \left[\ Y \ = \ {y\over \lambda_{\,1}} \ \right]\\[0.2in]
& \le & C\,\cdot\,  \int_{ {{\rho_{\,\nu}}\over {\lambda_{\,1}}} }^\infty\  {{R^{n\,-\,1} \ dR}\over {\  \, R^{{{n\,+\,2}\over 2} \ + \ (\,n\,-\,2\,)  } \ \, }}\ \ \ \ \ \ \ \ \ \  \ \ \ \ \  \ \ \ \ \ \ \ \ \ \  \ \ \ \ \  \ \ \ \ \ \ \ \ \ \  \ \ \ \ \ \ \ \ \ \ \ \   \left[\ R \ = \ \Vert\,Y\,\Vert\ \right]\\[0.2in]
& \le & C_1 \,\cdot\, {1\over {\  R^{{{\,n\,-\,2\,}\over 2} } \ }}\ \ \bigg\vert_{ \,{{\rho_{\,\nu}}\over {\lambda_{\,1}}} }^{\,\infty}\\[0.2in]& \le & C_2  \,\cdot\,  {\bar\lambda}_{\ \flat}^{ \ \left(\ {{ n\,-\,2}\over 2}\ \right) \ \cdot\,\ (\,1\, - \ \nu\,) } \ \ .\\
\end{eqnarray*}
It follows that
\begin{eqnarray*}
(\,A.8.41\,) \ \ \ \ \ & \ &       \int_{\R^n \,\setminus \ \left(\, \cup\ B_{\xi_l} (\,{\rho_{\,\nu}}) \, \right)}    V_1^{\,{1\over 2}\,\cdot\,{{n\,+\,2}\over {\,n\,-\ 2\,}}}\,\cdot\,\left[\    V_2^{\,{1\over 2}\,\cdot\,{{n\,+\,2}\over {\,n\,-\ 2\,}}}  \ + \,\cdot\,\cdot\,\cdot\, + \ V_{\,\flat}^{\,{1\over 2}\,\cdot\,{{n\,+\,2}\over {\,n\,-\ 2\,}}}  \ \right]\,\cdot\, V_1\ \ \ \ \ \ \ \ \ \ \ \ \ \  \\[0.2in]
& \le & C \left(\ \sum_{j\ =\,2}^\flat \ {1\over {\Vert\  \Xi_{\,\,1} \ - \ \Xi_j\,\Vert^{{{n\,+\,2}\over 2} } }}  \ \right)\,\cdot\,{\bar\lambda}_{\ \flat}^{ \ \left(\ {{ n\,-\,2}\over 2}\,+\,\epsilon \right) \,\,\cdot\,\, (\,1\ - \ \nu\,) }\\[0.2in]
& \le & C_1 \,\cdot\,{\bar\lambda}_{\ \flat}^{ \   {{n\,+\,2}\over 2} \,\cdot\, \gamma \ + \ {{\,n\,-\,2\,}\over 2} \,\cdot\, (\,1\ - \ \nu\,)  }\ .
\end{eqnarray*}

Finally\,,\, we consider terms like [\,cf. (\,A.4.56\,)\,]\\[0.1in]
(\,A.8.42\,)
\begin{eqnarray*}
& \ & V_1 \,\cdot\,V_2^{\,{1\over 2}\,\cdot\,{{n\,+\,2}\over {\,n\,-\ 2\,}}}\,\cdot\, V_3^{\,{1\over 2}\,\cdot\, {{n\,+\,2}\over {\,n\,-\ 2\,}}} \\[0.2in]
& \le &  C_1 \,\cdot\,{1\over {{\bar\lambda}_{\ \flat}^n}}\,\cdot\, \left(\ {{1}\over {1 \ + \ \Vert\ Y\ - \ \Xi_{\,\,1}\,\Vert }}\ \right)^{\!\! n\,-\,2 } \,\cdot\,\left(\ {{1}\over {1 \ + \ \Vert\ Y\ - \ \Xi_{\,2}\,\Vert  }}\ \right)^{\!\!{{n\,+\,2}\over 2} } *\\[0.2in]
& \ & \hspace*{3.8in} *\,\left(\ {{1}\over {1 \ + \ \Vert\ Y\ - \ \Xi_3\,\Vert  }}\ \right)^{\!\!{{n\,+\,2}\over 2} }\\[0.2in]
& \le  &   \ C_2\,\cdot\,{1\over {{\bar\lambda}_{\ \flat}^n}}\,\cdot\,{1\over {\Vert\  \Xi_{\,\,1} \ - \ \Xi_{\,2}\,\Vert^{{{n\,+\,2}\over 2}} }}\,\cdot\,\left[\ \left(\ {{1}\over {1 \ + \ \Vert\ Y\ - \ \Xi_{\,\,1}\,\Vert }}\ \right)^{\!\!(\,n\,-\,2)}\right.   \\[0.2in]
& \ & \hspace*{1in}\left. + \ \left(\ {{1}\over {1 \ + \ \Vert\ Y\ - \ \Xi_{\,2}\,\Vert }}\ \right)^{\!\!(\,n\,-\,2)} \ \right] \,\,\times\,  \left(\ {{1}\over {1 \ + \ \Vert\ Y\ - \ \Xi_3\,\Vert  }}\ \right)^{\!\!{{n\,+\,2}\over 2}} \\[0.2in]
& \le &   C_3\,\cdot\,{1\over {{\bar\lambda}_{\ \flat}^n}}\,\cdot\,{1\over {\Vert\  \Xi_{\,\,1} \ - \ \Xi_{\,2}\,\Vert^{{{n\,+\,2}\over 2}}}}\cdot {1\over {\Vert\  \Xi_{\,\,1} \ - \ \Xi_{\,3}\,\Vert^{{{\,n\,-\,2\,}\over 2}} }} \,\,*\\[0.2in]
& \ & \hspace*{1.2in}*\, \left[\ \left(\ {{1}\over {\ 1 \ + \ \Vert\ Y\ - \ \Xi_{\,\,1}\,\Vert \ }}\ \right)^{\!\!n}   + \ \left(\ {{1}\over {\ 1 \ + \ \Vert\ Y\ - \ \Xi_{\,3}\,\Vert \ }}\ \right)^{\!\!n} \ \right]    \\[0.2in]& \  &  \\[0.2in]
& \  &  \ \ + \ C_3\,\cdot\,{1\over {{\bar\lambda}_{\ \flat}^n}}\,\cdot\,{1\over {\Vert\  \Xi_{\,\,1} \ - \ \Xi_{\,2}\,\Vert^{{{n\,+\,2}\over 2}}}}\cdot {1\over {\Vert\  \Xi_{\,2} \ - \ \Xi_{\,3}\,\Vert^{{{\,n\,-\,2\,}\over 2}\,-\,\epsilon} }}  \,\,*\\[0.2in]
& \ & \hspace*{1.2in}*\, \left[\ \left(\ {{1}\over {\ 1 \ + \ \Vert\ Y\ - \ \Xi_{\,2}\,\Vert \ }}\ \right)^{\!\!n}  \ + \ \left(\ {{1}\over {\ 1 \ + \ \Vert\ Y\ - \ \Xi_{\,3}\,\Vert\  }}\ \right)^{\!\!n} \ \right]
\end{eqnarray*}
\begin{eqnarray*}
\Longrightarrow \ \ & \ &
 \int_{\R^n \,\setminus \, \left(\, \cup\ B_{\xi_l} (\,{\rho_{\,\nu}}) \, \right)}    V_2^{\,{1\over 2}\,\cdot\, {{n\,+\,2}\over {\,n\,-\ 2\,}}}\,\cdot\, \left[\   V_{3 }^{\,{1\over 2}\,\cdot\,{{n\,+\,2}\over {\,n\,-\ 2\,}}} \ + \ \cdot\,\cdot\,\cdot\,\ + \  V_{\flat}^{\,{1\over 2}\,\cdot\,{{n\,+\,2}\over {\,n\,-\ 2\,}}} \ \right]\,\cdot\,V_1   \\[0.2in]
  & \le &  C \,\cdot\,{\bar\lambda}_{\ \flat}^{\ n\,\cdot\, \gamma \ - \  o_{\,{\bar\lambda}_{\,\flat\,} } (\,1\,)}\ \ \ \ \ \ \\[0.1in]
   & \ &  \ \ \ \ \ \ \ \ \ [ \ {\mbox{similar \ \ to \ \ (\,A.8.41\,)\,, \ \ here}}  \ \ o_{\,{\bar\lambda}_{\,\flat\,} } (\,1\,) \ \to\ 0^+ \ \   {\mbox{as}} \ \ {\bar\lambda}_{\,\flat\,} \ \to \ 0^+ \ ] \ .\\[0.3in]
\end{eqnarray*}
 {\bf \S\,\S\, A\,8\,.\,i.\,.}\ \
Summing up, we have\\[0.1in]
(\,A.8.43\,)
\begin{eqnarray*}
|\,{\bf A}^3_{\,(\,A.8.1\,)}\,| & \le &     \int_{\R^n \,\setminus \ \left(\, \cup\ B_{\xi_l} (\,{\rho_{\,\nu}}) \, \right)}  \left\{  \  \left[\  W_{\,\,\flat}  \ \right]^{{\,n\,+\,2\,}\over {n\,-\,2}} \ - \  \sum_{l\,=\,1}^\flat \, V_{\,l}^{{\,n\,+\,2\,}\over {n\,-\,2}}      \  \right\}\,\cdot\,V_1  \\[0.2in]
    & \le &   C_1 \,\cdot\, {\bar\lambda}_{\ \flat}^{ \   {{n\,+\,2}\over 2}\,\cdot\,\gamma \ + \ {{\,n\,-\,2\,}\over 2}\,\cdot\,(\,1\ - \ \nu\,)  }\,\cdot\, {\bar\lambda}_{\ \flat}   \ + \  C_2\,\cdot\,{\bar\lambda}_{\ \flat}^{\ n\cdot \gamma \ - \  o_{\,{\bar\lambda}_{\,\flat\,} } (\,1\,) } \\[0.2in]
    & \le &   C_3 \,\cdot\, {\bar\lambda}_{\ \flat}^{ \   {{n\,+\,2}\over 2}\,\cdot\,\gamma \ + \ {{\,n\,-\,2\,}\over 2}\,\cdot\,(\,1\ - \ \nu\,)  }\,\cdot\, {\bar\lambda}_{\ \flat}^{  \   - \ o_{\,{\bar\lambda}_{\,\flat\,} } (\,1\,) } \ .  \\[0.2in]
\end{eqnarray*}
Note that (\,refer to {\bf \S\,A\,8.\,e} \ \& \ {\bf \S\,A\,8.\,f}\ )
$$
\gamma \ + \ \nu \ > \ 1 \ \ \Longrightarrow \  \ \gamma \ > \ (\,1\ - \ \nu\,) \ \ \Longrightarrow \  \  n \,\cdot\,\gamma \ > \  {{n\,+\,2}\over 2}\,\cdot\,\gamma \ + \ {{\,n\,-\,2\,}\over 2}\,\cdot\,(\,1\ - \ \nu\,)\,.
$$

 \vspace*{0.3in}

 {\bf \S\,\S\, A\,8\,.\,j\,.}\ \ {\bf Outside}. Refer to (4.15) in the main text\,,\,  where we introduce the error term $\,{\bf E}_{\,(\,4.15\,)}\,$,\, related to $\,{\bf{E}}_{\,(\,A.8.44\,)}\,$ below\,:
    \begin{eqnarray*}
(\,A.8.44\,)\ \ \ \ \ \ \ \ \    & \ & \int_0^{\rho_{\,\nu}}\!
\left[    \left(\ {{\lambda_{\,1}}\over {\lambda^2_1 \,+\, {\bar r}^2}}\  \right)^{\!\!2}\!\!
\cdot {{ (\,\lambda_{\,1}^2 \ -\  {\bar r}^{\,2})}\over {(\,\lambda_{\,1}^2 \ + \  {\bar r}^{\,2})^{{n }\over 2} }}\right]\,\cdot\, {\bar r}^{n\,-\,1}\,\cdot\,d\,{\bar r}\\[0.15in]
& = & \int_0^\infty
\left[    \left(\ {{\lambda_{\,1}}\over {\lambda^2_1 \,+\, {\bar r}^2}}\  \right)^{\!\!2}\!\!
\cdot {{ (\,\lambda_{\,1}^2 \ -\  {\bar r}^{\,2})}\over {(\,\lambda_{\,1}^2 \ + \  {\bar r}^{\,2})^{{n }\over 2} }}\right]\,\cdot\, {\bar r}^{n\,-\,1}\,\cdot\,d\,{\bar r}  \ + \ {\bf{E}}_{\,(\,A.8.44\,)} \ .\ \ \ \ \ \ \ \ \ \ \ \ \ \ \ \ \ \ \ \ \ \ \ \ \ \ \
\end{eqnarray*}

Here
\begin{eqnarray*}
 |\,{\bf{E\,'}}_{\,(\,A.8.44\,)} \,| & \le & \int_{\rho_{\,\nu}}^\infty
\left[    \left(\ {{\lambda_{\,1}}\over {\lambda^2_1 \,+\, {\bar r}^2}}\  \right)^{\!\!2}\!\!
\cdot {{ (\,\lambda_{\,1}^2 \ + \  {\bar r}^{\,2})}\over {(\,\lambda_{\,1}^2 \ + \  {\bar r}^{\,2})^{{n }\over 2} }}\right]\,\cdot\, {\bar r}^{n\,-\,1}\,\cdot\,d\,{\bar r}\\[0.2in]
  & = & \int_{\rho_{\,\nu}}^\infty \left(\ {\lambda_{\,1}\over { \,\lambda_{\,1}^2 \, +\  {\bar r}^{\,2}\, }}  \right)^{\!\!2}\,\cdot\,  {{1}\over {(\,\lambda_{\,1}^2 \, +\  {\bar r}^{\,2}\, )^{{{n \,- \ 2}\over 2 } }}}    \,\cdot\,{\bar r}^{\,n\,-\,1}\ d\,{\bar r}  \\[0.2in]
  & = &  \int_{  {{ \rho_{\,\nu} }\over {\lambda_{\,1}}} }^\infty \left(\ {1\over { \,1 \, +\   \left(\ {{ \bar r }\over {\lambda_{\,1}}}  \right)^{\!2}\, }}  \right)^{\!\!2}\,\cdot\,\left(\ {1\over { \,1 \, +\   \left(\ {{ \bar r }\over {\lambda_{\,1}}}  \right)^{\!2}\, }}  \right)^{\!\!{{\,n\,-\ 2\,}\over 2}}    \,\cdot\, \left(\ {{ \bar r }\over {\lambda_{\,1}}}  \right)^{\!n\,-\emph{}\,1} \ d\,\left(\ {{ \bar r }\over {\lambda_{\,1}}}  \right)\ \ \ \ \ \ \ \ \ \ \ \ \ \ \ \ \ \  \\[0.2in]
  & \le & C_1\,\cdot\, \int_{  {{ \rho_{\,\nu} }\over {\lambda_{\,1}}} }^\infty {{R^{n\,-\,1} \ dR }\over {R^{n\,+\,2}  }} \ \le \ C_2\,\cdot\,{{-1}\over {R^2}} \ \bigg\vert_{\,{  {{ \rho_{\,\nu} }\over {\lambda_{\,1}}} }}^{\,\infty} \ \ \ \ \ \ \  \ \ \ \ \ \ \  \ \ \ \ \ \ \  \ \ \   \left[\ {\bar R} \ = \ {{ \bar r }\over {\lambda_{\,1}}} \ \right] \\[0.2in]
    & = & O \left(\ {\bar \lambda}_{\,\flat}^{\ 2\,   (\,1 \,-\ \nu\,)}\ \right)\ .
  \end{eqnarray*}

\newpage

 {\bf \S\,A\,8\,.\,k\,.} \ \ {\bf Conclusion.} \ \
Combining the above discussion in this section\,,\, we finally arrive at\\[0.1in]
(\,A.8.45\,)
\begin{eqnarray*}
{\bf (A)}_{\,(\,A.8.1\,)} & = &  -\, {\hat C}_{\,1\,,\,2}   \left\{ \  \sum_{l\,=\ 2}^\flat \left(\ {{\lambda_{\,1}^{{\,n\,-\ 2\,}\over 2}\,\cdot\,\lambda_{\,\,l}^{{\,n\,-\ 2\,}\over 2}}\over {\Vert\ \xi_{\,1}\ - \ \xi_{\,l}\,\Vert^{\,n\,-\,2} }}\  \right)
 \cdot\left[\,1 \ + \ O\,\left(\ {\bar\lambda}_{\ \flat}^{\,2\,(\,1\ - \ \nu\,)}\,\right)
 \ \right]\,\,\right\} \ +\  \\[0.25in]
 & \ &     \ \ \ \ \ \ \ \ \ \ \ \ \ \ \ \  + \  \   O\left(\   {\bar\lambda}_{\ \flat}^{ \   {{n\,+\,2}\over 2}\,\,\cdot\,\,\gamma \ + \ {{\,n\,-\,2\,}\over 2}\,\cdot\,(\,1\ - \ \nu\,)  \ - \ o_{\,+}(\,1\,)} \ \right) \ + \ O\left( \  {\bar\lambda}_{\ \flat}^{ \,n \,\gamma \ - \ \sigma}\ \right)  \\[0.25in] & = &  -\, {\hat C}_{\,1\,,\,2}   \left\{ \  \sum_{l\,=\ 2}^\flat \left(\ {{\lambda_{\,1}^{{\,n\,-\ 2\,}\over 2}\,\cdot\,\lambda_{\,\,l}^{{\,n\,-\ 2\,}\over 2}}\over {\Vert\ \xi_{\,1}\ - \ \xi_{\,l}\,\Vert^{\,n\,-\,2} }}\  \right)
 \cdot\left[\,1 \ + \ O\,\left(\ {\bar\lambda}_{\ \flat}^{2\,(\,1\ - \ \nu\,)}\,\right)
 \ \right]\,\,\right\} \ +\  \\[0.25in]
 & \ &  \ \ \ \ \ \ \ \  + \  \   O\left(\ {\bar\lambda}_{\ \flat}^{ \ n\,\cdot\,(\,1\ - \ \nu\,)  \ - \ o_{\,+}(\,1\,)} \ \right) \ + \ O\left( \  {\bar\lambda}_{\ \flat}^{ \,n \,\gamma \ - \ \sigma}\ \right) \ ,
\end{eqnarray*}
as
$$
\gamma \ + \ \nu \ > \ 1 \ \ \Longrightarrow  \ \ \gamma \ > \ 1\ - \ \nu \ \ \Longrightarrow \   {{n\,+\,2}\over 2}\,\,\cdot\,\,\gamma \   >  \ {{\,n\,-\,2\,}\over 2}\,\cdot\,(\,1\ - \ \nu\,)  \ .
$$
In the above\,,\,
$$
\ \ {\hat C}_{\,1\,,\,2}\ = \ {{ (\,n\,-\,2)^{\,2} }\over 2}\,\cdot\,\omega_n\ . \leqno{\mbox{where}}
$$

Here we take into the account that
$$
 o_{\,{\bar\lambda}_{\,\flat\,} } (\,1\,) \ =  \  \varepsilon \,[\ (\,\gamma \ + \ \nu\,) \ - \ 1\ ]\ .
$$
Cf. (\,A.8.8\,)\,,\, (\,A.8.30\,)\,,\, (\,A.8.33\,)\,,\, (\,A.8.38\,)\, and \,(\,A.8.43\,)\,.\, \bk
Here we work under the assumptions in Proposition 6.8\,.\, Note that condition (\,1.35\,) of the main text\,,\, namely\,,
$$
\gamma \ > \ {\sigma\over {\,n \ - \ 2\,}} \ ,
$$
implies
$$
\gamma \ > \ {\sigma\over {\,n\,}} \ \ \Longrightarrow \ \ n\,\gamma \ -\ \sigma \ > \ 0\ .
$$

\newpage

{\large{\bf \S\,A\,\,9\,.  \ \ Estimate of}} $\,{\bf(B)}_{{\,(\,A.8.1\,)}}\,.$\\[0.1in]
Let us recall the expression of $\,{\bf(B)}_{{\,(\,A.8.1\,)}}\,:$
$$
 - \   \left[\ \int_{\R^n}\, [\,(\,{\tilde c}_n\!\cdot K\,) \ - \ n\,(\,n\,-\,2)\,]  \cdot\,(W_{\,\,\flat} )^{{\,n\,+\,2\,}\over {n\,-\,2}}\,\cdot\,\left(\,\lambda_{\,1}\,\cdot\,{ { \partial \,V_1 } \over {\partial\, \lambda_{\,1}}}\ \right)  \ \right] \ .
$$
From (\,A.8.4\,)\,,\, for $\,y\,\in\,B_{\,\xi_{\,1}}(\,\rho_{\,\nu}\,)\ ,\,$ \\[0.1in]
(\,A.9.1\,)
\begin{eqnarray*}
& \ &\left[\  (W_{\,\,\flat} )^{{\,n\,+\,2\,}\over {n\,-\,2}}\,\cdot\,\left(\ \lambda_{\,1}\,\cdot\,{ { \partial  V_1 } \over {\partial\, \lambda_{\,1}}}\ \right)\ \right]_{\,y}\\[0.2in]
& = &   \left[\ V_1^{{\,n\,+\,2\,}\over {n\,-\,2}} \ + \  {{n\,+\,2}\over {\,n\,-\ 2\,}} \,\cdot\,V_1^{4\over {n\,-\,2}}\,\cdot\,\left[\   V_2  \ + \ \cdot\,\cdot\,\cdot\,\ + \ V_{\,\flat}  \ \right] \ + \ O \left(\   V_2^{{\,n\,+\,2\,}\over {n\,-\,2}}  \ + \ \cdot\,\cdot\,\cdot\,\ + \ V_{\,\flat}^{{\,n\,+\,2\,}\over {n\,-\,2}}\  \right)\ \right]_{\,y}\!\!* \\[0.2in]
& \ & \hspace*{4.3in} \,*\,\left(\ \lambda_{\,1}\,\cdot\,{ { \partial \,V_1 } \over {\partial\, \lambda_{\,1}}}\ \right)_{y}\\[0.2in]
& = &\left( V_1^{{\,n\,+\,2\,}\over {n\,-\,2}}\,\cdot\,\lambda_{\,1}\,\cdot\, { { \partial \,V_1 } \over {\partial\, \lambda_{\,1}}}\ \right)\,\cdot\, \left[\ 1 \, + \,  {{n\,+\,2}\over {\,n\,-\ 2\,}} \,\cdot\,\left[\   {{V_2}\over {V_1}}  \ + \ \cdot\,\cdot\,\cdot\,\ + \  {{V_{\,\flat}}\over {V_1}}  \ \right] \right. \ + \ \\[0.2in]
& \ & \hspace*{3in} \left. \, + \  O \left(\   \left[\, {{V_2}\over {V_1}}\,\right]^{{\,n\,+\,2\,}\over {n\,-\,2}}  \, + \,\,\cdot\,\cdot\,\cdot\,\ + \ \left[\, {{V_{\,\flat}}\over {V_1}}\,\right]^{{\,n\,+\,2\,}\over {n\,-\,2}}\  \right)\ \right]\ .
\end{eqnarray*}
We continue with [\,see {\bf {\S\,A\,4\,.h}}\,,\, in particular (\,A.4.36\,)\ ]
\begin{eqnarray*}
  \left[\   {{V_2}\over {V_1}}  \ + \ \cdot\,\cdot\,\cdot\,\ + \  {{V_{\,\flat}}\over {V_1}}  \ \right]_{\,y} \ = \  O\,\left( \ {\bar\lambda}_{\ \flat}^{\ (\,n\ -\,2\,)\,\cdot\, [\ (\,\gamma \ + \ \nu \,)\ - \  1\,] }\ \right)\ \ \ \ \mfor y\,\in\,B_{\,\xi_{\,1}}(\,\rho_{\,\nu}\,)  \ .\\
\end{eqnarray*}
It follows that\\[0.1in]
(\,A.9.2\,)
\begin{eqnarray*}
&  \ &   -\,\int_{B_{\,\xi_{\,1}} (\,{\rho_{\,\nu}}\,)}\,[\,(\,{ c}_n\!\cdot K\,) \ - \ n\,(\,n\,-\,2)\,] \cdot\,(W_{\,\,\flat} )^{{\,n\,+\,2\,}\over {n\,-\,2}}\,\cdot\,\left(\ { { \partial \,V_1 } \over {\partial\, \lambda_{\,1}}}\ \right)  \\[0.2in]   &  = &  \left\{ \ -\, \int_{B_{\,\xi_{\,1}} (\,{\rho_{\,\nu}}\,)}\, [\,(\,{ c}_n\!\cdot K\,) \ - \ n\,(\,n\,-\,2)\,] \cdot\,V_1^{{\,n\,+\,2\,}\over {n\,-\,2}}\,\cdot\,\left(\ { { \partial \,V_1 } \over {\partial\, \lambda_{\,1}}}\ \right)  \  \right\}\,*\\[0.2in]
& \ & \hspace*{2.5in} *\ \left[\,1 \ + \   O\,\left( \ {\bar\lambda}_{\ \flat}^{\ (\,n\ -\,2\,)\,\cdot\, [\ (\,\gamma \ + \ \nu \,)\ - \  1\,] }\ \,\right) \ \right] \ .
\end{eqnarray*}

\newpage

{\underline{Outside}}\\[0.1in]
(\,A.9.3\,)
\begin{eqnarray*}
&  \ &
\bigg\vert \ \int_{\R^n \setminus\, B_{\,\xi_{\,1}} (\,{\rho_{\,\nu}}\,)  }\, [\,(\,{\tilde c}_n\!\cdot K\,) \ - \ n\,(\,n\,-\,2)\,]  \cdot\,(\,W_{\,\,\flat}\, )^{{\,n\,+\,2\,}\over {n\,-\,2}}\,\cdot\,\left(\,\lambda_{\,1}\,\cdot\,{ { \partial \,V_1} \over {\partial\, \lambda_{\,1}}}\ \right)\   \bigg\vert\\[0.2in]
& \le & C\,\cdot\,\int_{\R^n \setminus\, B_{\,\xi_{\,1}} (\,{\rho_{\,\nu}}\,)  }\, (W_{\,\,\flat} )^{{\,n\,+\,2\,}\over {n\,-\,2}}\,\cdot\,V_1\\[0.2in]
& \le & C\,\cdot\,   \left[\ {\bar\lambda}_{\ \flat}^{\,n\,(\,1\,-\ \nu\,)} \ + \     {\bar\lambda}_{\ \flat}^{ \   {{n\,+\,2}\over 2}\,\cdot\,\gamma \ + \ {{\,n\,-\,2\,}\over 2} \,\cdot\, (\,1\ - \ \nu\,)  } \,\cdot\, {\bar\lambda}_{\ \flat}^{  \   - \ 2\,\varepsilon \,[\ (\,\gamma \ + \ \nu\,) \ - \ 1\ ] }   \  \right]\\[0.1in]
& \ & \hspace*{-0.45in}  \left[\ {\mbox{estimate \ \ as \ \ \ in \ \ {\bf \S\,A\,8.\,h}}}\,,\, \ \ \ {\mbox{using}} \ \ \ \  (\,1\ +\ a\,)^{{\,n\,+\,2\,}\over {n\,-\,2}} \ \le \ C\cdot \left(\ 1\ + \ a^{{\,n\,+\,2\,}\over {n\,-\,2}} \ \right) \ \ \ \ {\mbox{for}} \ \ a \ \ge \ 0\ \right] \\[0.1in]
& \le & C_1\,\cdot\,  {\bar\lambda}_{\ \flat}^{\,n\,(\,1\,-\ \nu\,)} \,\cdot\, {\bar\lambda}_{\ \flat}^{  \   - \ 2\,\varepsilon \,[\ (\,\gamma \ + \ \nu\,) \ - \ 1\ ] }  \ \ \ \ \  \ \ \ \ \ \ \ \ \ \ \ \ \ \ \ [\ \mbox{as \ \ before}\,,\, \ \  \gamma \ >  \ (\ 1\ - \ \nu\ ) \ ]\ .
\end{eqnarray*}
We come to the conclusion that
\begin{eqnarray*}
(\,A.9.4\,) \ \ \ \ \ \ \ \ &  \ &
 - \   \left[\ \int_{\R^n}\, [\,(\,{\tilde c}_n\!\cdot K\,) \ - \ n\,(\,n\,-\,2)\,]  \cdot\,(W_{\,\,\flat} )^{{\,n\,+\,2\,}\over {n\,-\,2}}\,\cdot\,\left(\,\lambda_{\,1}\,\cdot\,{ { \partial \,W_{\,\,\flat} } \over {\partial\, \lambda_{\,1}}}\ \right)  \ \right]\ \ \ \ \ \ \ \ \ \ \ \ \ \ \ \ \ \ \ \ \ \ \ \ \ \ \ \ \ \ \ \
  \\[0.2in]
  &  = &  \left\{ -\, \int_{B_{\,\xi_{\,1}} (\,{\rho_{\,\nu}}\,)}\, [\ (\,{ c}_n\!\cdot K\,) \ - \ n\,(\,n\,-\,2)\ ] \cdot\,V_1^{{\,n\,+\,2\,}\over {n\,-\,2}}\,\cdot\,\left(\ { { \partial \,V_1 } \over {\partial\, \lambda_{\,1}}}\ \right)  \right\}\,*\\[0.2in]
  & \ & \ \ \ \ \ \ \ \ \ \ \ \ \ \  \ \ \ \ \ \ \  \ \ \ \ \ \ \  \ \ \ \ \ \ \  * \ \left[\ 1 \ + \ O\,\left( \ {\bar\lambda}_{\ \flat}^{\ (\,n\ -\,2\,)\,\cdot\, [\ (\,\gamma \ + \ \nu \,)\ - \  1\,] }\ \right) \  \right]\ + \\[0.2in]
    & \ & \ \ \ \ \ \ \ \ \ \ \ \ \ \ \ \ \ \ \ \  + \ O \left(\  {\bar\lambda}_{\ \flat}^{\ n\,(\,1\,-\ \nu\,) \ - \ o_{\,+}\,(\,1\,) }  \,\ \right)\ . \ \ \ \ \ \ \ \
\end{eqnarray*}

Gathering together (\,A.8.44\,) and (\,A.9.4\,)\,,\, we obtain
\begin{eqnarray*}
(\,A.9.5\,) & \ & \lambda_{\,1}\cdot {\partial\over {\partial\, \lambda_{\,1}}} \, {\bf I}\,(\,W_{\,\,\flat} \ + \ \phi\,)\\[0.2in]  & = &     -\, {\tilde C}_2\,(\,n)  \left\{ \  \sum_{l\,=\ 2}^\flat \left(\ {{\lambda_{\,1}^{{\,n\,-\ 2\,}\over 2}\,\cdot\, \lambda_{\ \!l}^{{\,n\,-\ 2\,}\over 2}}\over {\Vert\ \xi_{\,1}\ - \ \xi_{\,l}\,\Vert^{\,n\,-\,2} }}\  \right)
 \cdot\left[\,1 \ + \ O\,\left(\,{\bar\lambda}_{\ \flat}^{\,2\,\gamma}\,\right) \ \right]\,\,\right\} \ +\  \\[0.2in]
 & \ & \ \ \ \  +  \  \left\{\  -\, \int_{B_{\,\xi_{\,1}} (\,{\rho_{\,\nu}}\,)}\, [\,(\,{ c}_n\!\cdot K\,) \ - \ n\,(\,n\,-\,2)\,] \cdot\,V_1^{{\,n\,+\,2\,}\over {n\,-\,2}}\,\cdot\,\left(\,\lambda_{\,1}\,\cdot\, { { \partial \,V_1 } \over {\partial\, \lambda_{\,1}}}\ \right) \  \right\}\,*\ \ \ \ \ \ \ \\[0.2in]
 & \ & \hspace*{1in} \,\times \ \left[\,1 \ + \   O\,\left( \ {\bar\lambda}_{\ \flat}^{\ (\,n\ -\,2\,)\,\cdot\, [\ (\,\gamma \ + \ \nu \,)\ - \  1\,] }\ \right) \ \right]\\[0.1in]& \ &
  \hspace*{3in}\ \ \ \  \ \ \ \ \ \ \ \ \ \ \ \  + \ \ { \cal E}_{\,(\,A.9.5\,)\,} \ ,
 \end{eqnarray*}
 where
 \begin{eqnarray*}
{ \cal E}_{\,(\,A.9.5\,)\,} & = &   O\left(\   {\bar\lambda}_{\ \flat}^{ \   {{n\,+\,2}\over 2}\,\,\cdot\,\,\gamma \ + \ {{\,n\,-\,2\,}\over 2}\,\cdot\,(\,1\ - \ \nu\,) \ - \ o_{\,{\bar\lambda}_{\,\flat\,} } (\,1\,) } \right)\  \   +  \   O \left(\  {\bar\lambda}_{\ \flat}^{\ n\,(\,1\,-\ \nu)} \,\cdot\, {\bar\lambda}_{\ \flat}^{  \   - \ 2\,\varepsilon \,[\ (\,\gamma \ + \ \nu\,) - \ 1\ ] } \, \right) \\[0.2in]
& \ & \hspace*{4in} \ \ \ \ \ + \ \ {\cal E}_{\,(\,A.7.55\,)} \\[0.1in]
& = &   O \left(\  {\bar\lambda}_{\ \flat}^{\ n\,(\,1\,-\ \nu)} \,\cdot\, {\bar\lambda}_{\ \flat}^{  \   - \ 2\,\varepsilon \,[\ (\,\gamma \ + \ \nu\,) - \ 1\ ] } \, \right)  \ + \ \ {\cal E}_{\,(\,A.7.55\,)}  \ ,
 \end{eqnarray*}
 as
 $$
 \gamma \ >  \ 1 \ - \ \nu \ \ \Longrightarrow \ \   {{n\,+\,2}\over 2}\,\,\cdot\,\,\gamma \ + \ {{\,n\,-\,2\,}\over 2}\,\cdot\,(\,1\ - \ \nu\,) \ >\  n\,(\,1\,-\ \nu\,)\,.
 $$

\vspace*{0.3in}

 {\bf \S\,A\,9\,.a\,.}\ \
{\bf Estimate on } \\[0.1in]
(\,A.9.6\,)
\begin{eqnarray*}
& \ &  -\, \int_{B_{\,\xi_{\,1}} (\,{\rho_{\,\nu}}\,)}\, [\,(\,{ c}_n\!\cdot K\,) \ - \ n\,(\,n\,-\,2)\,] \cdot\,V_1^{{\,n\,+\,2\,}\over {n\,-\,2}}\,\cdot\,\left(\,\lambda_{\,1}\,\cdot\, { { \partial \,V_1 } \over {\partial\, \lambda_{\,1}}}\ \right)\\[0.2in]
& =  &    \int_{B_{\,\xi_{\,1}} (\,{\rho_{\,\nu}}\,)}\  [\,n\,(\,n\,-\,2) \ - \ (\,{ c}_n\!\cdot K\,) \,] \cdot\,V_1^{{\,n\,+\,2\,}\over {n\,-\,2}}\,\cdot\,\left(\,\lambda_{\,1}\,\cdot\, { { \partial \,V_1 } \over {\partial\, \lambda_{\,1}}}\ \right) \ \ \cdot\,\cdot\,\cdot\,\cdot\,\cdot\,\cdot\,\cdot\,\cdot\cdot\, \cdot\ \ (\,{{\bf B}}_{(\,A.9.6\,)}\,)\\[0.2in]
& =  & -\,{{n\ -\ 2}\over 2}\,\cdot\,\int_{B_{\,\xi_{\,1}} (\,{\rho_{\,\nu}}\,)}\, \,[\,n\,(\,n\,-\,2) \,-\, (\,{ c}_n\!\cdot K\,)  \ ] \,\cdot\, \left(\ {\lambda_{\,1}\over {\lambda^2_1 \ + \ \Vert\ y\,-\,\xi_{\,1}\,\Vert^2}}\  \right)^{\!\!{{n\,+\,2}\over 2} }* \\[0.2in]
& \ & \ \ \ \ \ \ \ \  \ \  \ \ \ \ \ \ \ \  \ \  \ \ \ \ \ \ \ \ \ \ \ \ \   * \  \left(\ {\lambda_{\,1}\over {\lambda^2_1 \ + \ \Vert\ y\,-\,\xi_{\,1}\,\Vert^2}}\  \right)^{\!\!{{\,n\,-\,2\,}\over 2} }  \cdot  \left(\ {{\lambda^2_1 \ - \ \Vert\ y\,-\,\xi_{\,1}\,\Vert^2}\over {\lambda^2_1 \ + \ \Vert\ y\,-\,\xi_{\,1}\,\Vert^2}}\  \right)\ dy\\[0.2in]
& = &  -\,{{n\ -\ 2}\over 2}\,\cdot\, \int_{B_{\,\xi_{\,1}} (\,{\rho_{\,\nu}}\,)} [\ \cdot\,\cdot\,\cdot\,\ ]_{\,(\,A.9.6\,)}  \ dy \ ,
\end{eqnarray*}
where\\[0.1in]
(\,A.9.7\,)
$$ [\ \cdot\,\cdot\,\cdot\,\ ]_{\,(\,A.9.6\,)}  \ = \   \,[\,n\,(\,n\,-\,2) \,-\, (\,{ c}_n\!\cdot K\,)\,(\,y\,)  \ ] \,\cdot\, \left(\ {{\lambda^2_1 \ - \ \Vert\ y\,-\,\xi_{\,1}\,\Vert^2}\over {\lambda^2_1 \ + \ \Vert\ y\,-\,\xi_{\,1}\,\Vert^2}}\  \right)\,\cdot\,\left(\ {\lambda_{\,1}\over {\lambda^2_1 \ + \ \Vert\ y\,-\,\xi_{\,1}\,\Vert^2}}\  \right)^{\!\!n }  \ .
$$
Recall that
there exists a fixed compact set $\,{\cal C}_{\,[]} \, \subset \, \R^n\,$   such that
$$
 \{\  \xi_{\,1}\,, \ \cdot \cdot \cdot\,, \ \xi_{\ \flat}  \ \} \ \subset \ {\cal C}_{\,[]} \ \ \ \ \ [ \ {\mbox{independent \   \ on}} \ \  \flat \ ] \ . \leqno (\,A.9.8\,)
$$
Now we take over from  (\,4.46\,)  in {\bf \S\,4\,c} of the main text\,:
$$ \int_{B_{\,\xi_{\,1}}\,(\,\rho_{\,\nu}\,)  }\,[\,n\,(\,n\,-\,2\,) \,-\, (\,{ c}_n\!\cdot K\,)  \ ]  \cdot  \widetilde{ [\ \cdot\,\cdot\,\cdot\,\ ]\,}_{(\,A.9.9\,)}  \ dy \ , \leqno (\,A.9.9\,)
$$
where
$$
\widetilde{ [\ \cdot\,\cdot\,\cdot\,\ ]\,}_{(\,A.9.9\,)}  \ = \ \left(\ {\lambda_{\,1}\over {\lambda_{\,1}^2 \ + \ \Vert\ y\,-\,\xi_{\,1}\,\Vert^2}}\  \right)^{\!\!n }  \,\cdot\, \left(\ {{\lambda_{\,1}^2 \ - \ \Vert\ y\,-\,\xi_{\,1}\,\Vert^2}\over {\lambda_{\,1}^2 \ + \ \Vert\ y\,-\,\xi_{\,1}\,\Vert^2}}\  \right)\ . \leqno (\,A.9.10\,)
$$
Recall that
$$
\xi_{1_{\,|_{\,n}}} \ \le \ {\bar\lambda}_{\ \flat}^{\ 1\,+\,\kappa} \ \ \  \ {\mbox{with}} \ \ \ \  2\,\nu \ >  \ 1 \ + \ \kappa\ .
$$

\vspace*{0.3in}


 {\bf \S\,A\,9\,.b\,.}\ \  {\bf Error terms in}\,  (\,4.46\,) {\bf \,of the main text [\,8\,]\,.} \ \ Applying the expansion in (\,4.44\,) and (\,4.45\,) in the main text\,,\, we continue with [\ cf. (\,4.46\,) of the main text \ ]\\[0.1in]
 (\,A.9.11\,)
\begin{eqnarray*}
&  \ &     \int_{B_{\,\xi_{\,1}} (\,{\rho_{\,\nu}}\,)}
\,[\,n\,(\,n\,-\,2\,) \,-\, (\,{ c}_n\!\cdot K\,)\,(\,y\,)  \ ] \,\cdot\, \widetilde{ [\ \cdot\,\cdot\,\cdot\,\ ]\,}_{(\,A.9.9\,)} \  dy \\[0.2in]
& = &     \int_{B_{\,\xi_{\,1}} (\,{\rho_{\,\nu}}\,)}
\,[\ C\,( \, {\bf p}_{\,y}\ )\cdot  \,
 \Vert\ y \ - \ {\bf p}_{\,y}\,\Vert^{\,\ell} \ - \ {\bf R}\,(\,y\,) \ ]  \,\cdot\, \widetilde{ [\ \cdot\,\cdot\,\cdot\,\ ]\,}_{(\,A.9.9\,)} \  dy \\[0.2in]
& = &  \int_{B_{\,\xi_{\,1}}\,(\,\rho_{\,\nu}\,)   }  C\,(\,{\bf p}_{\,y}\,) \,\cdot \Vert \,y  \ - \ {\bf p}_{\,y}\,\Vert^{\,\ell} \,\cdot\, \widetilde{ [\ \cdot\,\cdot\,\cdot\,\ ]\,}_{(\,A.9.9\,)} \  dy  \ + \ {\bf E}^{\,\mbox{Remainder}}_{\, (\,A.9.11\,) }
 \\[0.1in]
 & \ &  \hspace*{3in} \ [ \ {\mbox{see \ \ (\,A.9.12\,)}} \ \ \uparrow \ ]
 \\[0.1in]
 & = & C\,(\,{\bf p}_{\,1}\,) \,\cdot    \int_{B_{\,\xi_{\,1}}\,(\,\rho_{\,\nu}\,)   }    \Vert \,y  \ - \ {\bf p}_{\,y}\,\Vert^{\,\ell} \,\cdot\, \widetilde{ [\ \cdot\,\cdot\,\cdot\,\ ]\,}_{(\,A.9.9\,)} \  dy \ + \ \, {\bf E}^{\,\mbox{Remainder}}_{\, (\,A.9.11\,) } \ + \ {\bf E}^{\,\mbox{Shift}}_{\, (\,A.9.11\,) }
 \\[0.1in]
 & \ &  \hspace*{4in} \ [ \ {\mbox{see \ \ (\,A.9.13\,)}} \ \ \uparrow \ ]
 \\[0.1in]
  & = & C\,(\,{\bf p}_{\,1}\,) \,\cdot    \int_{B_{\,\xi_{\,1}}\,(\,\rho_{\,\nu}\,)   }   |\ y_n \ - \ \xi_{1_{\,|_{\,n}}}\,|^{\,\,\ell}\,\cdot\, \widetilde{ [\ \cdot\,\cdot\,\cdot\,\ ]\,}_{(\,A.9.9\,)} \  dy \\[0.2in]
   & \ & \ \ \ \ \ \ + \ C\,(\,{\bf p}_{\,1}\,) \,\cdot   {{ \ell\,\cdot\,(\,\ell\,-\,1\,) }\over {2}} \cdot \xi_{1_{\,|_{\,n}}}^2 \cdot  \int_{B_{\,\xi_{\,1}}\,(\,\rho_{\,\nu}\,)   }   |\ y_n \ - \ \xi_{1_{\,|_{\,n}}}\,|^{\,\,\ell \ - \ 2}\,\cdot\, \widetilde{ [\ \cdot\,\cdot\,\cdot\,\ ]\,}_{(\,A.9.9\,)} \  dy\ + \ \\[0.2in]
  & \ & \ \ \ \ \ \ \ \ \ \ \ \ \ \ \ \ \ \ \ + \ \, {\bf E}^{\,\mbox{Remainder}}_{\, (\,A.9.11\,) } \ + \ \, {\bf E}^{\,\mbox{Shift}}_{\, (\,A.9.11\,) } \ + \ {\bf E}^{\,\mbox{1st}}_{\, (\,A.9.11\,) } \ + \ {\bf E}^{\,\mbox{Lower}}_{\, (\,A.9.11\,) }  \ + \ {\bf E}^{\,\mbox{Thin}}_{\, (\,A.9.11\,) }
  \\[0.1in]
 & \ &  \hspace*{3.3in} \ [ \  \leftarrow \ \  {\mbox{see \ \ Note \ \ A.9.14\,}} \ \ \uparrow \ \ \ \ \rightarrow \ ]
 \\[0.1in]
& = &   C\,(\,{\bf p}_{\,1}\,) \,\cdot   \int_{B_{o}\,(\,\rho_{\,\nu}\,)}   |\,{\tilde y}_n  \,|^{\,\,\ell}    \times   \left(\, {\lambda_{\,1}\over {\lambda^2_1 \ + \ \Vert\ {\tilde y}\,\Vert^2}}\, \right)^{\!\!n }  \,\cdot\, \left(\ \, {{\lambda^2_1 \ - \ \Vert\ {\tilde y}\,\Vert^2}\over {\lambda^2_1 \ + \ \Vert\ {\tilde y}\,\Vert^2}} \, \right) \ d{\tilde y}\\[0.1in]
& \ &  \ \ \ \ \ \ \ \ \ \ \ \ \ \ \ \ \ \ \ \ \hspace*{2in} [\ \tilde y \ = \  y\,-\,\xi_{\,1} \ \ \Longrightarrow \ \ {\tilde y}_{\,n} \ = \ y_{\,n} \ - \ \xi_{1\,|_{\,n}} \ ]\\[0.1in]
   & \ & \ \ \ \ \ + \ C\,(\,{\bf p}_{\,1}\,) \,\cdot   {{ \ell\,\cdot\,(\,\ell\,-\,1\,) }\over {2}} \cdot \xi_{1_{\,|_{\,n}}}^2 *\\[0.1in]
     & \ & \ \ \ \  \ \ \  \ \ \ *\ \int_{B_{o}\,(\,\rho_{\,\nu}\,)}   |\,{\tilde y}_n  \,|^{\,\,\ell \ - \ 2}    \times   \left(\, {\lambda_{\,1}\over {\lambda^2_1 \ + \ \Vert\ {\tilde y}\,\Vert^2}}\, \right)^{\!\!n }  \,\cdot\, \left(\ \, {{\lambda^2_1 \ - \ \Vert\ {\tilde y}\,\Vert^2}\over {\lambda^2_1 \ + \ \Vert\ {\tilde y}\,\Vert^2}} \, \right) \ d{\tilde y} \ +  \\[0.2in]
  & \ & \ \ \ \ \ \ \ \ \ \ \ \ \ \ \ \ \ \ \ + \ \, {\bf E}^{\small\mbox{Remainder}}_{\, (\,A.9.11\,) } \ + \ \, {\bf E}^{\,\mbox{Shift}}_{\, (\,A.9.11\,) } \ + \ {\bf E}^{\,\mbox{1st}}_{\, (\,A.9.11\,) } \ + \ {\bf E}^{\,\mbox{Lower}}_{\, (\,A.9.11\,) }  \ + \ {\bf E}^{\,\mbox{Thin}}_{\, (\,A.9.11\,) }  \\[0.2in]
& = & \lambda^{\,\ell}_1\cdot C\,(\,{\bf p}_{\,1}\,) \,\cdot   \int_{B_o  \,\left(\,\lambda^{-1}\,\cdot\,{\rho_{\,\nu}}\ \right)} |\,Y_n\,|^{\,\,\ell}\,\cdot\,\left(\ {1\over {1 \ + \ \Vert\ Y\,\Vert^2}}\  \right)^{\!\!n }  \,\cdot\, \left(\ {{1 \ - \ \Vert\ Y\,\Vert^2}\over {1 \ + \ \Vert\ Y\,\Vert^2}}\  \right)\ d\,Y   \ + \\[0.2in]
& \ &\ \ \ \ \ + \ \lambda^{\,\ell \ - \ 2}_1\cdot C\,(\,{\bf p}_{\,1}\,)  \,\cdot   {{ \ell\,\cdot\,(\,\ell\,-\,1\,) }\over {2}} \,*\\[0.1in]
     & \ & \ \ \ \  \ \ \  \ \ \ \ \ \ \ \ *\,  \int_{B_o  \,\left(\,\lambda^{-1}\,\cdot\,{\rho_{\,\nu}}\ \right)} |\,Y_n\,|^{\,\,\ell \ - \ 2}\,\cdot\,\left(\ {1\over {1 \ + \ \Vert\ Y\,\Vert^2}}\  \right)^{\!\!n }  \,\cdot\, \left(\ {{1 \ - \ \Vert\ Y\,\Vert^2}\over {1 \ + \ \Vert\ Y\,\Vert^2}}\  \right)\ d\,Y   \ + \\[0.2in]
  & \ & \ \ \ \ \ \ \ \ \ \ \ \ \ \ \ \ \ \ \ + \ \, {\bf E}^{\,\mbox{Remainder}}_{\, (\,A.9.11\,) } \ + \ \, {\bf E}^{\,\mbox{Shift}}_{\, (\,A.9.11\,) } \ + \ {\bf E}^{\,\mbox{1st}}_{\, (\,A.9.11\,) } \ + \ {\bf E}^{\,\mbox{Lower}}_{\, (\,A.9.11\,) }  \ + \ {\bf E}^{\,\mbox{Thin}}_{\, (\,A.9.11\,) }    \\[0.2in]
& = & \lambda^{\,\ell}_1\cdot C\,(\,{\bf p}_{\,1}\,) \,\cdot    \left[\ -\,\int_{\R^n} |\,Y_n\,|^{\,\,\ell}\,\cdot\,\left(\ {1\over {1 \ + \ \Vert\ Y\,\Vert^2}}\  \right)^{\!\!n }  \,\cdot\, \left(\ {{\Vert\ Y\,\Vert^2 \ - \ 1   }\over {\Vert\ Y\,\Vert^2 \ + \ 1  }}\  \right)\ d\,Y  \ \right] \ + \ \\[0.2in]
& \ &\ \ \ \  \ \ \ + \  \lambda^{\,\ell \ - \ 2 }_1\cdot C\,(\,{\bf p}_{\,1}\,)  \,\cdot   {{ \ell\,\cdot\,(\,\ell\,-\,1\,) }\over {2}} *\\[0.1in]
     & \ & \ \ \ \  \ \ \  \ \ \ \ \ \ \ \ \ \ \ \ \  *  \left[\ -\,\int_{\R^n} |\,Y_n\,|^{\,\,\ell\ - \ 2}\,\cdot\,\left(\ {1\over {1 \ + \ \Vert\ Y\,\Vert^2}}\  \right)^{\!\!n }  \,\cdot\, \left(\ {{\Vert\ Y\,\Vert^2 \ - \ 1   }\over {\Vert\ Y\,\Vert^2 \ + \ 1  }}\  \right)\ d\,Y  \ \right] \ +  \\[0.2in]
  & \ & \ \ \ \ \ \ + \ \, {\bf E}^{\,\mbox{Remainder}}_{\, (\,A.9.11\,) } \ + \ \, {\bf E}^{\,\mbox{Shift}}_{\, (\,A.9.11\,) } \ + \ {\bf E}^{\,\mbox{1st}}_{\, (\,A.9.11\,) } \ + \ {\bf E}^{\,\mbox{Lower}}_{\, (\,A.9.11\,) } \ + \ {\bf E}^{\,\mbox{Thin}}_{\, (\,A.9.11\,) } \ + \ {\bf E}^{\,\mbox{Out}}_{\, (\,A.9.11\,) }   \  .   \\[0.1in]
 & \ &  \hspace*{4in} \ \ \ \ \ \  [ \ {\mbox{see \ \ (\,A.9.15\,) }} \ \ \uparrow   \ ]
\end{eqnarray*}

\newpage

In the above, \\[0.1in]
 (\,A.9.12\,)
$$
\bigg\vert \ {\bf E}^{\,\mbox{Remainder}}_{\, (\,A.9.11\,) }\,\bigg\vert \ \le \    \int_{B_{\,\xi_{\,1}} (\,{\rho_{\,\nu}}\,)}
\,|\  {\bf R}\,(\,y\,) \ | \cdot \Bigg\vert \  \, {{\lambda^2_1 \ - \ \Vert\ y\,-\,\xi_{\,1}\,\Vert^2}\over {\lambda^2_1 \ + \ \Vert\ y\,-\,\xi_{\,1}\,\Vert^2}} \ \, \Bigg\vert \,\cdot\,\left(\ {\lambda_{\,1}\over {\lambda^2_1 \ + \ \Vert\ y\,-\,\xi_{\,1}\,\Vert^2}}\  \right)^{\!\!n }\  dy\ ,
$$
which is estimated in in the next sub\,-\,section (\,that is\,,\, {\bf \S\,A\,9\,.c}\,)\,.\, Next\,,\, \\[0.1in]
 (\,A.9.13\,)
\begin{eqnarray*}
\bigg\vert\ {\bf E}^{\,\mbox{Shift}}_{\, (\,A.9.11\,) }\,\bigg\vert & \le &   \int_{B_{\,\xi_{\,1}}\,(\,{\rho_{\,\nu}})}  |\   C \,({\bf p}) \ - \ C \,({\bf p_1})  \ | \cdot \Vert \,y  \ - \ {\bf p}_{\,y}\,\Vert^\ell    \ \,*\\[0.2in]
& \ & \hspace*{1in} \,* \ \Bigg\vert \  \, {{\lambda^2_1 \ - \ \Vert\ y\,-\,\xi_{\,1}\,\Vert^2}\over {\lambda^2_1 \ + \ \Vert\ y\,-\,\xi_{\,1}\,\Vert^2}} \ \, \Bigg\vert \,\cdot\,\left(\ {\lambda_{\,1}\over {\lambda^2_1 \ + \ \Vert\ y\,-\,\xi_{\,1}\,\Vert^2}}\  \right)^{\!\!n }\  dy \ ,
\end{eqnarray*}

which is estimated in {\bf \S\,A\,9\,.d}\,.\,

\vspace*{0.2in}

{\it Note}  \,A.9.14\,\,. \ \  $\,{\bf E}^{\,\mbox{1st}}_{\, (\,A.9.11\,) }\,$  is related to the first order team in (\,4.44\,) of the main text\,,\, with self\,-\,cancellation\,.\, $\,{\bf E}^{\,\mbox{Lower}}_{\, (\,A.9.11\,) }\,$ comes from the lower order teams in  (\,4.44\,) of the main text\,,\, and $\,{\bf E}^{\,\mbox{Thin}}_{\, (\,A.9.11\,) }\,$ the thin layer, see {\bf \S\,A\,9.\,e}\,.\,

\vspace*{0.12in}

Lastly\,,\,\\[0.1in]
 (\,A.9.15\,)
\begin{eqnarray*}
\bigg\vert \ {\bf E}^{\,\mbox{Out}}_{\, (\,A.9.11\,) }  \,\bigg\vert & \le & \!\! C_1 \cdot  \lambda^{\,\ell}_1\cdot  \int_{\R^n \,\setminus \,B_o  \,\left(\,\lambda^{-1}\,\cdot\,{\rho_{\,\nu}}\ \right)} \!\! |\,Y_n\,|^{\,\,\ell}\,\cdot\,\left(\ {1\over {1 \ + \ \Vert\ Y\,\Vert^2}}\  \right)^{\!\!n }  \,\cdot\, \Bigg\vert \  {{1 \ - \ \Vert\ Y\,\Vert^2}\over {1 \ + \ \Vert\ Y\,\Vert^2}}  \ \Bigg\vert \  \ d\,Y+ \ \\[0.2in]
& \ & \hspace*{-1in} \ \ \ \ +  \
   C_2 \cdot \lambda^{\,\ell }_1\cdot  \int_{\R^n \,\setminus \,B_o  \,\left(\,\lambda^{-1}\,\cdot\,{\rho_{\,\nu}}\ \right)} \ \  |\,Y_n\,|^{\,\,\ell \ - \ 2}\,\cdot\,\left(\ {1\over {1 \ + \ \Vert\ Y\,\Vert^2}}\  \right)^{\!\!n }  \,\cdot\, \Bigg\vert \  {{1 \ - \ \Vert\ Y\,\Vert^2}\over {1 \ + \ \Vert\ Y\,\Vert^2}}  \ \Bigg\vert \  \ d\,Y\\[0.2in]
& \le &  C_1 \cdot {\bar\lambda}^{\,\ell}_{\ \flat} \cdot
   \int_{\lambda^{-1} \,\cdot\, {\rho_{\,\nu}}}^\infty {{R^{\,\ell}\,\cdot\, R^{\,n\,-\,1}}\over { R^{\,2n} }} \ \ dR \ + \
   C_2 \cdot  {\bar\lambda}^{\,\ell}_{\ \flat}\,\cdot\, \int_{{\rho_{\,\nu}}}^\infty {{R^{\,\ell \ - \ 2}\,\cdot\, R^{\,n\,-\,1}}\over { R^{\,2n} }}\\[0.1in]
     & \ & \hspace*{3.5in} \ \ \ \ \ \ \ \ \ \  \ [\ R \ = \ \Vert \,Y\,\Vert \ ] \\[0.1in]
   & \le & C_3\cdot {\bar\lambda}^{\,\ell}_{\ \flat}  \,\cdot \, {\bar\lambda}_{\ \flat}^{\,(\ n\,-\ \ell\,)\, \cdot\, (\,1 \ - \ \nu) } \ + \
   C_4 \cdot  {\bar\lambda}^{\,\ell}_{\ \flat}\,\cdot \, {\bar\lambda}_{\ \flat}^{\,(\ n \,+\,2\,-\ \ell\,)\, \cdot\, (\,1 \ - \ \nu) } \\[0.2in]
   & \le & C_5\cdot {\bar\lambda}^{\,\ell}_{\ \flat}  \,\cdot \, {\bar\lambda}_{\ \flat}^{\,(\ n\,-\ \ell\,)\, \cdot\, (\,1 \ - \ \nu) }  \\[0.1in]
    & \ & \hspace*{3in}\ \ \ \ \ \  [\ {\mbox{recall \ \ that}}\ \ \ell \ \le \ n\ - \ 2\ ]\ .
\end{eqnarray*}


The key term
$$
\int_{\R^n} |\,Y_n\,|^{\,\,\ell}\,\cdot\,\left(\ {1\over {1 \ + \ \Vert\ Y\,\Vert^2}}\  \right)^{\!\!n }  \,\cdot\, \left(\ {{\Vert\ Y\,\Vert^2 \ - \ 1   }\over {\Vert\ Y\,\Vert^2 \ + \ 1  }}\  \right)\ d\,Y  \leqno (\,A.9.16\,)
$$
is computed in the main text (\,see {\bf \S\,4\,c}\,)\,,\, and is shown to be a positive number. Likewise, the number
$$
\int_{\R^n} |\,Y_n\,|^{\,\,\ell\ - \ 2}\,\cdot\,\left(\ {1\over {1 \ + \ \Vert\ Y\,\Vert^2}}\  \right)^{\!\!n }  \,\cdot\, \left(\ {{\Vert\ Y\,\Vert^2 \ - \ 1   }\over {\Vert\ Y\,\Vert^2 \ + \ 1  }}\  \right)\ d\,Y
$$
is also positive (\,cf. {\bf \S\,4\,c}\,). In what follows we estimate the error terms.

\vspace*{0.3in}


 {\bf \S\,A\,9\,.c\,.}\ \  {\bf Estimate of $\, {\bf E}^{\,\mbox{Remainder}}_{\, (\,A.9.11\,) }\,$} {\bf \ \  coming from the  expansion}\,: \ \
 $$
(\,{\tilde c}_n\,K\,)\,(\,y\,) \,=\,n\,(\,n\,-\,2\,) \ - \ C\,(\,{\bf p}_{\,y}\,)\cdot \Vert\,y \,-\,{\bf p}_{\,y}\,\Vert^{\,\ell} \ + \  {\bf R}_{\,\ell \  + \ 1}(\,y\,) \ \  \mfor \ \ y \ \in \ {\cal O}\,.
$$
Referring  to (\,A.9.12\,)\,,\, we have \\[0.1in]
 (\,A.9.17\,)
 \begin{eqnarray*}
& \ & \bigg\vert \ {\bf E}^{\,\mbox{Remainder}}_{\, (\,A.9.11\,) }\  \bigg\vert \\[0.2in]
& \le &  {\bar C}_{\,R} \cdot    \int_{B_{\,\xi_{\,1}}\,(\,{\rho_{\,\nu}})}  \Vert \,y  \ - \ {\bf p}_{\,y}\,\Vert^{\,\ell\ +\,1}    \ | \cdot \Bigg\vert \  \, {{\lambda^2_1 \ - \ \Vert\ y\,-\,\xi_{\,1}\,\Vert^2}\over {\lambda^2_1 \ + \ \Vert\ y\,-\,\xi_{\,1}\,\Vert^2}} \ \, \Bigg\vert \,\cdot\,\left(\ {\lambda_{\,1}\over {\lambda^2_1 \ + \ \Vert\ y\,-\,\xi_{\,1}\,\Vert^2}}\  \right)^{\!\!n }\  dy \\[0.1in]
& \ &\hspace*{3in}\ \ \ \ \   (\ {\mbox{via \ \ (\,1.26\,) \ \ in \ \ the \ \ main \ \ text}}\ ) \\[0.1in]
 & \le &  C_2 \cdot    \int_{B_{\,\xi_{\,1}}\,(\,{\rho_{\,\nu}})}  \Vert \,y  \ - \ {\bf p}_{\,y}\,\Vert^{\,\ell\ +\,1}    \ |   \,\cdot\,\left(\ {\lambda_{\,1}\over {\lambda^2_1 \ + \ \Vert\ y\,-\,\xi_{\,1}\,\Vert^2}}\  \right)^{\!\!n }\  dy \\[0.2in]
  & \le &  C_2 \cdot    \int_{B_{\,\xi_{\,1}}\,(\,{\rho_{\,\nu}})}     \Vert\ y\,-\,\xi_{\,1}\,\Vert^{\,\ell\ +\,1}   \,\cdot\,\left(\ {\lambda_{\,1}\over {\lambda^2_1 \ + \ \Vert\ y\,-\,\xi_{\,1}\,\Vert^2}}\  \right)^{\!\!n }\  dy  \ + \ \,  {\bf E}_{\, (\,A.9.17\,) }  \\[0.2in]
  & = & C_2 \cdot \int_{B_{o}\,(\,\rho_{\,\nu}\,)} \ \   |\,{\tilde y}_n  \,|^{\,\,\ell \ + \ 1}    \times   \left(\, {\lambda_{\,1}\over {\lambda^2_1 \ + \ \Vert\ {\tilde y}\,\Vert^2}}\, \right)^{\!\!n }   \ d{\tilde y}    \ + \ \,  {\bf E}_{\, (\,A.9.17\,) }\\[0.1in]
     & \ & \hspace*{3in}\ \ \ \ \ \ \ \    [\ \tilde y \ = \  (\ y_1\,-\,\xi_{\,1}\,, \ y_2\,, \ \cdot \cdot \cdot\,, \ y_n\ ) \ ] \\[0.1in]
  & \ \le & C_2 \cdot \int_{\R^n}  \ |\,{\tilde y}_n  \,|^{\,\,\ell \ + \ 1}    \times   \left(\, {\lambda_{\,1}\over {\lambda^2_1 \ + \ \Vert\ {\tilde y}\,\Vert^2}}\, \right)^{\!\!n }   \ d{\tilde y} \ + \ \,  {\bf E}_{\, (\,A.9.17\,) }  \\[0.2in]
  & \le & C_3  \cdot {\bar\lambda}_{\ \flat}^{\,\ell \ + \ 1} \cdot  \int_0^\infty {{R^{\,\,\ell  \ + \ 1} \cdot R^{\,n\,-\,1} }\over {(\, 1 \ + \ R^2\,)^n } }\ \ dR  \ + \ \,  {\bf E}_{\, (\,A.9.17\,) }  \\[0.2in]
   & \ & \hspace*{3in} \ \left\{ \ R \ = \ \Vert \,Y\,\Vert \ = \ \bigg\Vert \ {y\over \lambda_{\,1}}  \bigg\Vert  \ \ge \ \bigg\vert \ {y_n\over \lambda_{\,1}}  \ \bigg\vert  \ \right\}\\[0.2in]
  & \le & C_4  \cdot {\bar\lambda}_{\ \flat}^{\,\ell \ + \ 1}   \ + \ \,  {\bf E}_{\, (\,A.9.17\,) }   \ \ \ \ \ \ \ \  \  \ \ \ \ \ \ \ \  \  \ \ \ \ \ \ \ \  \ ( \ \ell \ \le \ n\ - \ 2\ )\,.
\end{eqnarray*}
Here  $\,{\bf E}_{\, (\,A.9.17\,) }\,$  is related to the lower order teams mentioned in Note A.9.14\,,\, and the thin layer, which is estimated by [\,see (\,A.9.29\,) in \,{\bf \S\,A\,9\,.e}\ ]\,.
$$
| \ {\bf E}_{\, (\,A.9.17\,)} \ | \ \ \le \ C\,  {\bar\lambda}_{\ \flat}^{ \, (\,\ell \,+\,1\,) \cdot (\,1 \,+\,\kappa\,)} \cdot {\bar\lambda}_{\ \flat}^{\,\kappa} \  . \leqno (\,A.9.18\,)
$$

\vspace*{0.3in}


 {\bf \
 S\,A\,9\,.d\,.}\ \  {\bf Estimate of}\, $\,\displaystyle{ {\bf E}^{\,\mbox{Shift}}_{\, (\,A.9.11\,) }\, } \,$.\ \
Introducing the polar coordinates in $\,\R^{\,n\,-\,1}\,,\,$ we have [\,\,recall that $\,{\bf p}_1 \ = \ 0\,$,\, refer to (\,4.24\,)\,,\,  (\,4.33\,) and (\,4.44\,) of the main text\ ]
$$
|\  C \,(\,{\bf p}_{\,y}\,) \ - \ C \,(\,{\bf p_1}\,)  \ | \ \le \ {  C}_1 \,\cdot\,{{\bar r}} \ \ \ \ \ \ \ \left(\ \ {{\bar r}} \ = \ \sqrt{\ {\tilde y}_1\  + \ \cdot\,\cdot\,\cdot\,\ + \ {\tilde y}_{n\,-\,1}\ }\ \, \right)\leqno (\,A.9.19\,)
$$
for $\,y\,\in\,{B_{\,\xi_{\,1}}\,(\,{\rho_{\,\nu}})}\,.\,$ Here
$$
{\tilde y} \ =\ y \ - \ \xi_{\,1} \ .
$$
Due to the uniform boundness on $\,\Vert\,\xi_{\,l}\,\Vert\,$\, [\,\,see (\,1.4\,) of the main text \ ]\,:

$$
 \{\  \xi_{\,1}\,, \ \cdot \cdot \cdot\,, \ \xi_{\,\flat}  \ \} \ \subset \ {\cal C}_{\,[]} \ \ \ \ \ [ \ \leftarrow \ \  {\mbox{compact set\,, \ \ independent \   \ on}} \ \  \flat \ ] \ ,
$$

 the constant $\,{ C}_1\,$ can be chosen to be independent on $\,\xi_{\,1}\,.\,$
We continue with\\[0.1in]
(\,A.9.20\,)
\begin{eqnarray*}
& \ & \bigg\vert \  \int_{B_{\,\xi_{\,1}}\,(\,{\rho_{\,\nu}})}  [\   C \,(\,{\bf p}_{\,y} \ - \ C \,({\bf p_1})  \ ]\cdot \Vert \,y  \ - \ {\bf p}_{\,y}\,\Vert^\ell     \cdot\left( {\lambda_{\,1}\over {\lambda_{\,1}^2 \ + \ \Vert\ y\,-\,\xi_{\,1}\,\Vert^2}}\  \right)^{\!\!n }  \,*\\[0.1in]
& \ & \hspace*{3.4in} *\, \left(\ {{\lambda_{\,1}^2 \ - \ \Vert\ y\,-\,\xi_{\,1}\,\Vert^2}\over {\lambda_{\,1}^2 \ + \ \Vert\ y\,-\,\xi_{\,1}\,\Vert^2}}\  \right)\ dy \ \bigg\vert \\[0.2in]
& \le &    \int_{B_{\,\xi_{\,1}}\,(\,{\rho_{\,\nu}})}  |\  C \,(\,{\bf p}_{\,y} \ - \ C \,({\bf p_1})  \  |\cdot \Vert \,y  \ - \ {\bf p}_{\,y}\,\Vert^\ell   \cdot\left( {\lambda_{\,1}\over {\lambda_{\,1}^2 \ + \ \Vert\ y\,-\,\xi_{\,1}\,\Vert^2}}\  \right)^{\!\!n }    \ dy\\[0.2in]
& \le & C_1 \,\cdot\, \int_{B_{o}\,(\,{\rho_{\,\nu}})}  {{\bar r}} \,\cdot\,|\,{\tilde y}_n\,|^\ell     \cdot\left( {\lambda_{\,1}\over {\lambda_{\,1}^2 \ + \ \Vert\ {\tilde y}\,\Vert^2}}\  \right)^{\!\!n }    \ d\,{\tilde y}  \  \\[0.1in]
& \ & \hspace*{4in} \ \ \ \ \ \  [ \ {\tilde y} \ =\ y \ - \ \xi_{\,1} \ ]\\[0.1in]
& \le & C_2 \int_{ \ \{ \,{{\bar r}}  \ \le \ \rho_{\,\nu}\, \}   \ \cap \ \{ \,|\,{\tilde y}_n\,| \ \le \ \rho_{\,\nu}\, \}  }  {{\bar r}} \,\cdot\,|\,y_n\,|^\ell    \,\cdot\, \left(\ {\lambda_{\,1}\over {\lambda^2_1 \ + \ \Vert\,{\tilde y} \,\Vert^2}}\  \right)^{\!\!n }  \ d{\tilde y}  \ + \ \,  {\bf E}_{\, (\,A.9.20\,) } \\[0.2in]
& \le & C_3\,\cdot\,\int_{-\,\rho_{\,\nu}}^{\,\rho_{\,\nu}} \left\{ \ \int_0^{\, \rho_{\,\nu}  }  {{\bar r}}\,\cdot\, \left(\ \ {\lambda_{\,1}\over {[\ \lambda^2_1 \ + \ {\tilde y}_n^2\,] \ + \  {{\bar r}}^2}} \ \right)^{\!\!n }  {{\bar r}}^{\,(\,n\,-\,1)\,-\,1 }\ d {{\bar r}}\ \right\}\,\cdot\, |\,y_n\,|^\ell  \ d{\tilde y}_{\,n} \ + \ \,  {\bf E}_{\, (\,A.9.20\,) }  \\[0.2in]
& = & C_3\,\cdot\,\int_{-\,\rho_{\,\nu}}^{\,\rho_{\,\nu}} \left\{ \ \int_0^{\, \rho_{\,\nu}  }  \left(\ {\lambda_{\,1}\over {\,\Omega^2 \ + \  {{\bar r}}^2}\,} \right)^{\!\!n }  {{\bar r}}^{\  n\,-\,1   }\ d {{\bar r}}\ \right\}\,\cdot\, |\,y_n\,|^\ell  \ d{\tilde y}_{\,n}  \ + \ \,  {\bf E}_{\, (\,A.9.20\,) } \\[0.2in]
& \ & \hspace*{3.5in} \ \ \ \ \ \ \ \  \left[\ \Omega \ = \ \sqrt{\,\lambda^2_1 \ + \ {\tilde y}_{\,n}^{\,2} \ }\    \, \right] \\[0.1in]
& = & C_3\,\cdot\, \int_{-\,\rho_{\,\nu}}^{\,\rho_{\,\nu}}\ \left(\ {{\lambda_{\,1}}\over {\Omega}}\  \right)^{\!\!n}  \left\{ \ \int_0^{\ \rho_{\,\nu}  }  \left(\ {\Omega \over {\,\Omega^2 \ + \  {{\bar r}}^2}\,} \right)^{\!\!n }  {{\bar r}}^{\ n\,-\,1  }\ d {{\bar r}}\ \right\}\cdot  |\,y_n\,|^\ell  \ d{\tilde y}_{\,n}  \\[0.1in]
& \ &  \hspace*{4in} \ \ \ \ \ \ \ \ \ \ + \ \,  {\bf E}_{\, (\,A.9.20\,) }  \\[0.1in]
& = & C_3\,\cdot\, \int_{-\,\rho_{\,\nu}}^{\,\rho_{\,\nu}} \ \left(\ {{\lambda_{\,1}}\over {\Omega}}\  \right)^{\!\!n}\,\cdot\, \left\{ \ \int_0^{\, \rho_{\,\nu}  }  \left(\ {1\over {\,1 \ + \ \left(\ {{  {{\bar r}}  }\over {  \Omega}}\  \right)^{\!2} }\,} \right)^{\!\!n }  \left(\ {{  {{\bar r}}  }\over {  \Omega}}\  \right)^{\!  n\,-\,1 }\ d \left(\ {{  {{\bar r}}  }\over {  \Omega}}\  \right) \ \right\}\,\cdot\, |\,y_n\,|^\ell  \ d{\tilde y}_{\,n}\\[0.1in]
& \ & \hspace*{4.5in} \ + \ \,  {\bf E}_{\, (\,A.9.20\,) }  \\[0.1in]
& =& 2 \cdot C_3 \cdot \lambda_{\,1}^n  \!\cdot  \left(\, \int_0^{\,\rho_{\,\nu}} \left(\ \ {1\over {\lambda^2_1 \ + \ |\,{\tilde y}_n \,|^2}}\  \right)^{\!\!{{n}\over 2} } \,\cdot\, |\,y_n\,|^\ell  \ d{\tilde y}_{\,n}  \,  \right)\cdot\left(\, \int_0^{\, {{  \,\rho_{\,\nu} }\over {\Omega}}}   \left[\ {1\over {1 \ + \ R^2}} \ \right]^{\,n} \cdot R^{\,n\,-\,1}\,dR \, \right) \\[0.1in]
& \ &  \hspace*{3in} \ + \ \,  {\bf E}_{\, (\,A.9.20\,) }  \ \ \ \ \ \ \ \ \ \ [\ \uparrow \ = \ O\,(\,1)\ ] \, \\[0.1in]
& \le &  C_4 \,\cdot\,\lambda_{\,1}^n\,\cdot\,\int_0^{\,\rho_{\,\nu}} {{\lambda_{\,1}}\over { \lambda_{\,1}^{n}   }}\,\cdot\, \left(\ {1\over {1 \ + \ \left[\ {{ {\tilde y}_n}\over {\lambda_{\,1}}}  \ \right]^{\,2} }}\  \right)^{\!\!{{n}\over 2} }\,\cdot\,\left(\,{{ {\tilde y}_n} \over { \lambda_{\,1}}}\  \right)^{\!\ell}\, \,d\left(\,{{ {\tilde y}_n} \over { \lambda_{\,1}}}\  \right)   \ + \ \,  {\bf E}_{\, (\,A.9.20\,) }  \, \\[0.2in]
& = &     C_4 \,\cdot\,\lambda_{\,1}^{\,\ell\ + \ 1}\,\cdot\, \int_0^{\,{{ \rho_{\,\nu} }\over {\lambda_{\,1}}}} {{ {\tilde Y}_n^{\,\ell} \ \ d \,{\tilde Y}_n }\over {\  ( \, 1 \ + \ {\tilde Y}_n^2)^{n\over 2} \  }} \ + \ \,  {\bf E}_{\, (\,A.9.20\,) }  \ \ \ \ \ \ \ \ \ \ \ \ \ \  \ \left(\  {\tilde Y}_n \ = \ {{ {\tilde y}_n }\over {  \lambda_{\,1} }} \ \right) \\[0.2in]
& = &     C_4 \,\cdot\,\lambda_{\,1}^{\,\ell\ + \ 1  }  \ + \ \,  {\bf E}_{\, (\,A.9.20\,) }  \ \ \ \ \ \ \ \ \ \ \ \ \ \ \ \ \ \  \ \ \ \ \ \  \ \ \ \ \ \ \ \ \  \ (\ {\mbox{for}} \ \ \ell \ \le  \ n\ - \ 2\,)\,.
\end{eqnarray*}
Here  $\,{\bf E}_{\, (\,A.9.20\,) }\,$ is related to the lower order teams mentioned in Note A.9.14\,,\, and the thin layer, which can be estimated by (\,see \,{\bf \S\,A\,9\,.e}\,)
$$
| \ {\bf E}_{\, (\,A.9.20\,)} \ | \ \ \le \ C\, \left(\, {\bar\lambda}_{\ \flat}\, \right)^{ \ell \, \cdot \, (\\,1 \,+\ \kappa\,)} \,\cdot\, {\bar\lambda}_{\ \flat}^{\,\kappa} \ = \ C\, {\bar\lambda}_{\ \flat}^{\, \ell } \cdot {\bar\lambda}_{\ \flat}^{\,(\,\ell \ + \ 1\,) \,\cdot\,\kappa}\ . \leqno (\,A.9.21\,)
$$

 \newpage

 {\bf \S\,A\,9\,.e\,.}\ \  {\bf Estimate the error in the thin region\,. }
 As
 $$
\Bigg\vert\  {{\lambda^2_1 \ - \ \Vert\,y\,-\,\xi_{\,1}\,\Vert^2}\over {\lambda^2_1 \ + \ \Vert\,y\,-\,\xi_{\,1}\,\Vert^2}}\ \Bigg\vert  \ \le  \ 1 \ \ \ \ \ \ \mfor \  \ y \ \in  \ \R^n\,, \leqno (\,A.9.22\,)
 $$
 we have\,:\\[0.1in]
(\,A.9.23\,)
\begin{eqnarray*}
& \ & \!\!\!\!\!\!\!\!\!\!\!\Bigg\vert \ \int_{B_{\,\xi_{\,1}}\,(\,\rho_{\,\nu}\,) \ \cap \ \{ \,|\,y_n\,| \ \le \ C \,\Delta\, \}  }\,[\,n\,(\,n\,-\,2) \,-\, (\,{ c}_n\cdot\, K\,)  \ ] \,\cdot\,\left(\ {\lambda_{\,1}\over {\lambda^2_1 \ + \ \Vert\,y\,-\,\xi_{\,1}\,\Vert|^2}}\  \right)^{\!\!n }  \,*\ \ \ \ \ \ \ \ \ \ \ \ \ \ \ \ \\[0.2in]
& \ & \hspace*{3.5in} \ \ \ \ \ \   * \  \left(\ {{\lambda^2_1 \ - \ \Vert\,y\,-\,\xi_{\,1}\,\Vert^2}\over {\lambda^2_1 \ + \ \Vert\,y\,-\,\xi_{\,1}\,\Vert^2}}\, \right)\ dy \ \Bigg\vert\\[0.2in]
& \le & \int_{B_{\,\xi_{\,1}}\,(\,\rho_{\,\nu}\,) \ \cap \ \{ \,|\,y_n\,| \ \le \ C \,\Delta\, \}  }\,|\,n\,(\,n\,-\,2) \,-\, (\,{ c}_n\cdot\, K\,)  \ | \,\cdot\,\left(\ {\lambda_{\,1}\over {\lambda^2_1 \ + \ \Vert\,y\,-\,\xi_{\,1}\,\Vert|^2}}\  \right)^{\!\!n }\ dy\\
\end{eqnarray*}
where [\ recall (\,4.38\,) of the main text\ ]
$$
\Delta \ = \ O\,\left(\, \rho_{\,\nu}^2\, \right) \ + \ O\,\left( \ \xi_{1_{|_1}}\  \right) \ = \  O\,\left( \ {\bar\lambda}_{\ \flat}^{\,2\,\nu} \ \right) \ + \ O\,\left( \ {\bar\lambda}_{\ \flat}^{\,1\ + \ \kappa} \ \right)\ .
$$
Recall that we  accept the normalization
$$
\xi_{\,1} \ = \ (\,0\,, \ \cdot\,\cdot \cdot\,, \ 0 \,, \  \ \xi_{\,1_{\,|_{\,n}}}\,) \ \ \ \ \ \ \ {\mbox{with}} \ \ \ \ \ \xi_{\,1_{\,|_{\,n}}} \ > \ 0\,.
$$
Next\,,\, we observe that \\[0.1in]
(\,A.9.24\,)
\begin{eqnarray*}
\Vert\,y \ - \ {\bf p}_{\,y} \ \Vert & \le & |\,y_n\,|  \ + \ |\, {\bf p}_{\,y|_n} \,|  \ + \  O\,(\,\Vert\  y\,\Vert^2 \,)\ \ \ \ \ \ \,[ \ {\mbox{via \ \ (\,4.40\,) \ \ of \ \  the\ \  main\ \  text }} \ ]  \\[0.2in]
& \le &  C_1\, \Delta\ + \   O\,(\,\Vert\  y\,\Vert^2 \,) \ \ \ \ \ \ \ \ \ \ \ \  \ \ \ \ \ \  \ \,[ \ {\mbox{via \ \ (\,4.37\,) \ \ of \ \  the\ \  main\ \  text }} \ ]  \\[0.2in]
&  \le & C_2 \cdot \xi_{1_{|_1}} \ + \ C_1\cdot\,(\,\rho_{\,\nu}^2\, ) \\[0.2in]
&  \le & C_2 \cdot {\bar\lambda}_{\ \flat}^{\,1\ + \ \kappa} \ + \ C_1\cdot\,{\bar\lambda}_{\ \flat}^{\,2\,\nu}     \\[0.2in]
\Longrightarrow & \ &  | \ n\,(\,n\,-\,2) \,-\, (\,{ c}_n\cdot\, K\,)\,(\,y\,)  \ |  \  = \ O\,\left( \ {\bar\lambda}_{\ \flat}^{\,\ell\,\cdot\, (\ 1\ + \ \kappa \ )} \ \right)  \ +  \ O\,\left( \ {\bar\lambda}_{\ \flat}^{\,\ell\,\cdot\, (\ 2\,\nu \ )} \ \right)\\[0.1in]
& \ & \hspace*{1.5in}  {\mbox{for}} \ \ \ \  y  \in \ B_{\,\xi_{\,1}}\,(\,\rho_{\,\nu}\,) \ \cap \ \{ \ |\,y_n\,| \ \le \ C \,\Delta\  \} \ \ \ \ \ \ \ \ \ (\,\ell \ \ge \ 2\,)\ .
\end{eqnarray*}
Next\,,\, we estimate  the integral\,:
$$
\int_{B_{\,\xi_{\,1}}\,(\,\rho_{\,\nu}\,) \ \cap \ \{ \,|\,y_n\,| \ \le \ C \,\Delta\, \}  } \left(\ {\lambda_{\,1}\over {\lambda^2_1 \ + \ \Vert\,y\,-\,\xi_{\,1}\,\Vert^2}}\  \right)^{\!\!n }    \ dy\ .
$$
After moving the center to $\,\xi_{\,1}\,$,\, we continue with
\begin{eqnarray*}
(\,A.9.25\,) \ \ \ \ \ \ \ \ & \ & \int_{B_{\,\xi_{\,1}}\,(\,\rho_{\,\nu}\,) \ \cap \ \{ \,|\,y_n\,| \ \le \ C \,\Delta\, \}  } \left(\ {\lambda_{\,1}\over {\lambda^2_1 \ + \ \Vert\,y\,-\,\xi_{\,1}\,\Vert^2}}\  \right)^{\!\!n } \ \, dy \ \ \ \ \ \ \ \   \ \ \ \ \ \ \ \   \ \ \ \ \ \ \ \  \\[0.2in]
&  = &
\int_{ B_o \,(\,\rho_{\,\nu}\,) \ \cap \ \{ \,|\,{\tilde y}_n\,| \ \le \ {\tilde C}\cdot \Delta\, \}  } \left(\ {\lambda_{\,1}\over {\lambda^2_1 \ + \ \Vert\,{\tilde y} \,\Vert^2}}\  \right)^{\!\!n }  \ \ d{\tilde y}\ .
\end{eqnarray*}


 Here $\,{\tilde y} \ = \ y \ - \ \xi_{\,1}\,,\,$ and $$
{\tilde r} \ = \  \sqrt{\ {\tilde y}_1\  + \ \cdot\,\cdot\,\cdot\,\ + \ {\tilde y}_{n\,-\,1}\ } \ \ .
$$
Note that
\begin{eqnarray*}
(\,A.9.26\,) \ \ \ \ \ \ \ \ \ \ \ \ \ \ \ \ \ \ \ \  {\tilde y} & = & y \ - \ \xi_{\,1} \ = \ (\,y_1\,, \ \cdot\,\cdot \cdot\,, \ y_{\,n\,-\,1} \,, \ y_n \ - \ \xi_{1_{\,|_{\,n}}}\,) \ \ \ \ \ \ \ \ \ \ \ \ \ \ \ \ \ \ \ \ \ \ \ \ \  \ \ \ \ \   \\[0.2in]
\ \ \Longrightarrow   \ \ \ \ |\, {\tilde y}_{\,n} \,| & \le  & |\, y_{\,n} \,| \ + \ |\, \xi_{1_{\,|_{\,n}}}\,|  \ \le \ {\tilde C}\cdot \Delta \ \ \ \ {\mbox{for}}  \ \ \ \ \ |\,y_n\,| \ \le \ C \cdot\,\Delta \ .
\end{eqnarray*}
From (\,A.9.25\,)\,,\, we move on to\\[0.1in]
(\,A.9.27\,)
\begin{eqnarray*}
& \ & \int_{ B_o \,(\,\rho_{\,\nu}\,) \ \cap \ \{ \
|\,{\tilde y}_n\,| \ \le \ {\tilde C}\,\cdot\, \Delta            \  \} } \left(\ {\lambda_{\,1}\over {\lambda^2_1 \ + \ \Vert\,{\tilde y} \,\Vert^2}}\  \right)^{\!\!n }  \ \ d{\tilde y} \\[0.2in]
& \le & C_1\,\cdot\,\int_{-\,{\tilde C}\,\cdot\, \Delta }^{\ {\tilde C}\,\cdot\, \Delta} \left\{ \ \int_0^{\,\rho_{\,\nu}  }  \left(\ {\lambda_{\,1}\over {[\ \lambda^2_1 \ + \ {\tilde y}_n^2\,] \ + \  {\tilde r}^2}}\  \right)^{\!\!n  }\   {\tilde r}^{\,(\,n\,-\,1)\,-\,1 }\ d {\tilde r}\ \right\} \,d\,{\tilde y}_n  \\[0.2in]
& = & C_1 \,\cdot\,  \int_{-\,{\tilde C}\,\cdot\, \Delta }^{\ {\tilde C}\,\cdot\, \Delta}\left\{ \ \int_0^{\,\rho_{\,\nu}  }  \left(\ {\lambda_{\,1}\over {\,\Omega^2 \ + \  {\tilde r}^2}\,} \right)^{\!\!n  } \  {\tilde r}^{\ n\,-\,2 }\ d {\tilde r}\ \right\} \,d\,{\tilde y}_n \\[0.2in]
& \ & \!\!\!\!\!\!\!\! \Bigg[\ \Omega \ = \ \sqrt{\,\lambda^2_1 \ + \ {\tilde y}_n^2 \,}\ \approx \ \lambda_{\,1} \, \ \ {\mbox{inasmuch \ \ as}} \ \ \ {\tilde y}_n\  \, \ \le \ O\left( \rho_{\,\nu}^{2} \right) \ = \  \ O\left( \, {\bar\lambda}_{\,\,\flat}^{2\,\nu}\, \right) \ = \  o\,(\,\lambda_{\,1}\,) \ ,\\[0.1in]
& \ & \hspace*{1in}\ \  {\mbox{recall \ \ that \ \ }} 2\,\nu \ >  \ 1\  ; \ \ \ \ {\mbox{here}} \ \ {\tilde r} \ = \ \sqrt{\,{\tilde y}_1^2 \ + \ \cdot \cdot \cdot \ + \ {\tilde y}_{n\,-\,1}\,}\  \ \Bigg]\\[0.2in]
& = & C_1\,\cdot \int_{-\,{\tilde C}\,\cdot\, \Delta }^{\ {\tilde C}\,\cdot\, \Delta} \ \left(\ {{\lambda_{\,1}}\over {\Omega}}\  \right)^{\!\!n}  \left\{ \ \int_0^{\,\rho_{\,\nu}  }  \left(\ {\Omega \over {\,\Omega^2 \ + \  {\tilde r}^2}\,} \right)^{\!\!n }  {\tilde r}^{\,n\,-\,2 }\ d {\tilde r}\ \right\} \,d\,{\tilde y}_n\\[0.3in]
& = & C_1\,\cdot \int_{-\,{\tilde C}\,\cdot\, \Delta }^{\ {\tilde C}\,\cdot\, \Delta} \ \left(\ {{\lambda_{\,1}}\over {\Omega}}\  \right)^{\!\!n}\,\cdot\,\left(\ {{1}\over {\Omega}}\  \right) \,*\\[0.2in]
& \ &  \hspace*{2in} *\, \left\{ \ \int_0^{\,\rho_{\,\nu}  }  \left(\ {1\over {\,1 \ + \ \left(\ {{  {\tilde r}  }\over {  \Omega}}\  \right)^{\!2} }\,} \right)^{\!\!n }  \left(\ {{  {\tilde r}  }\over {  \Omega}}\  \right)^{\! n\,-\,2 }\ d \left(\ {{  {\tilde r}  }\over {  \Omega}}\  \right) \ \right\} \,d\,{\tilde y}_n\\[0.3in]
& =& C_1 \,\cdot  \lambda_{\,1}^{n}  \,\cdot\, \left(\ \int_0^{\,{\tilde C}\,\cdot\, \Delta} \left(\ {1\over {\lambda^2_1 \ + \ |\,{\tilde y}_n \,|^2}}\  \right)^{\!\!{{n\,+\,1}\over 2} }  \,d\,{\tilde y}_n \right)\,* \ \ \ \ \ \ \ \ \ \ \ \ \ \ \ \  \left( \ \tilde R \ = \    {{\tilde r  }\over {  \Omega}} \ \ \downarrow \ \right) \\[0.2in]
& \ &  \hspace*{3in} *\, \left(\ \ \int_0^{\, {{ \rho_{\,\nu} }\over {\Omega}}}   \left[\ {1\over {1 \ + \ {\tilde R}^2}} \ \right]^{\,n}\,\cdot\,{\tilde R}^{\,n\,-\,2}\,d\,\tilde R\ \right)\\[0.1in]
& \ & \hspace*{1.8in}   \left[\ \int_{\,0}^{\,\infty}   \left[\ {1\over {1 \ + \ {\tilde R}^2}} \ \right]^{\,n}\,\cdot\,{\tilde R}^{\,n\,-\,2}\,d\,\tilde R\ = \ O\,(\,1)\ \ \ \uparrow \ \right] \, \\[0.1in]
& \le &  C_2 \,\cdot\,\lambda_{\,1}^{\,n}\,\cdot\,\int_0^{\,{\tilde C}\,\cdot\, \Delta} {{\lambda_{\,1}}\over { \lambda_{\,1}^{\,n\,+\,1}   }}\,\cdot\,\left[\ {1\over {1 \ + \ \left(\, {{ {\tilde y}_n}\over {\lambda_{\,1}}}  \ \right)^{2} }}\  \right]^{\,\,{{n\,+\,1}\over 2} }  \,d\left(\,{{ {\tilde y}_n} \over { \lambda_{\,1}}}\  \right) \ \ \ \ \ \ \ \ \ \ \ \ \ \  \ \left(\  {\tilde Y}_n \ = \ {{ {\tilde y}_n }\over {  \lambda_{\,1} }} \ \right) \\[0.2in]
& \le &  C_3  \,\cdot\,  \int_0^{\,\lambda_{\,1}^{-1}\,\cdot\, (\ {\tilde C}\,\cdot\, \Delta\,)}\ d\,{\tilde Y}_n  \ \,\le \ \,C_4 \,\cdot\, (\,\lambda_{\,1}^{-1}\,\cdot\, \Delta\,)\ .\\
\end{eqnarray*}

As
$$
\Delta \ = \ O\,\left( \, \xi_{1_{\,|_{\,n}}} \, \right) \ + \ O\,\left( \, \rho^2_{\,\nu} \, \right)\ \ \ \ \ \ \ \ \ \ {\mbox{when}} \ \ \ \  2\,\nu \ \ge\ 1 \ + \ \kappa \ ,
$$
we are led to \\[0.1in]
(\,A.9.28\,)
\begin{eqnarray*}
& \ &  \bigg\vert \
\int_{B_{\,\xi_{\,1}}\,(\,\rho_{\,\nu}\,) \ \cap \ \{ \,|\,y_n\,| \ \le \ C \,\Delta\, \}  }\,[\,n\,(\,n\,-\,2) \,-\, (\,{ c}_n\!\cdot K\,)  \ ] \,\cdot\,\widetilde{[\ \cdot\,\cdot\,\cdot\,\ ]}_{\,(\,A.9.9\,)} \ dy \ \bigg\vert \\[0.2in]
& \le &
      C_4  \,\cdot\,(\,{\bar\lambda}_{\ \flat}^{\ 1 \, + \, \kappa}\ )^{\,\ell} \ \cdot\,\   {\bar\lambda}_{\ \flat}^{\ \kappa}  \ + \ C_5\,\cdot\,(\,{\bar\lambda}_{\ \flat}^{\ 2\,\nu}\ )^{\,\ell} \cdot\,{\bar\lambda}_{\ \flat}^{\ 2\,\nu \ - \ 1} \\[0.2in]
& \le &
      C_6 \cdot {\bar\lambda}_{\ \flat}^{\ \ell } \cdot \left[\  {\bar\lambda}_{\ \flat}^{\  (\ \ell\,+\,1\,) \,\cdot\,\kappa }  \ + \  {\bar\lambda}_{\ \flat}^{\  (\ \ell\,+\,1\,) \,\cdot \,(\  2\,\nu \ - \ 1\,)} \ \right] \ \le \  C_6 \,\cdot\, {\bar\lambda}_{\ \flat}^{\ \ell } \cdot \left[\  {\bar\lambda}_{\ \flat}^{\  3\,\cdot\,\kappa }  \ + \  {\bar\lambda}_{\ \flat}^{\  (\ \ell\,+\,1\,) \,\cdot \,(\  2\,\nu \ - \ 1\,)} \ \right]\,
\end{eqnarray*}



%
as $\,\ell \ \ge \ 2\,.\,$ Recall also that $\,\nu \ > \ {1\over 2}\,.\,$ Note that
$$
\widetilde{[\ \cdot\,\cdot\,\cdot\,\ ]}_{\,(\,A.9.9\,)}  \ = \ \left(\ {\lambda_{\,1}\over {\lambda_{\,1}^2 \ + \ \Vert\ y\,-\,\xi_{\,1}\,\Vert^2}}\  \right)^{\!\!n }  \,\cdot\, \left(\ {{\lambda_{\,1}^2 \ - \ \Vert\ y\,-\,\xi_{\,1}\,\Vert^2}\over {\lambda_{\,1}^2 \ + \ \Vert\ y\,-\,\xi_{\,1}\,\Vert^2}}\  \right)\ .
$$
In a similar fashion, we have
$$
 |\ {\bf E}_{\, (\,A.9.17\,) }\  | \ \le \ C_6 \cdot \left( {\bar\lambda}_{\ \flat}^{\,(\,1 \ + \ \kappa\,)} \right)^{\,\ell \ + \ 1\,  }\  \cdot {\bar\lambda}_{\ \flat}^{\, \kappa}   \ + \ C_7\,\cdot\,(\,{\bar\lambda}_{\ \flat}^{\ 2\,\nu}\ )^{\,\ell \, +\,1} \ \cdot\,\
{\bar\lambda}_{\ \flat}^{\ 2\,\nu  \ - \ 1}   \ . \leqno (\,A.9.29\,)
$$
[ \ Recall via (\,1.30\,) of the main text  that  $\,2\,\nu \ \ > \ 1 \,.$\ ]

\vspace*{0.35in}

 {\bf \S\,A\,9.\,f\,.}\ \  {\bf First order terms in}\, (\,4.44\,) {\bf and}\, (\,4.45\,) {\bf of the main text}\,. \ \ We first take note that (\,here \,$L\,$,\, $M\,$ and $\,N\,$ are given positive numbers\,)\\[0.1in]
 (\,A.9.30\,)
\begin{eqnarray*}
& \ & \int_{B_{\,\xi_{\,1}}\,(\,{\rho_{\,\nu}})\ \cap \ \{ \ y_n \ >  \ \xi_{1_{\,|_{\,n}}}\}}  \ \ (\,y_n\,-\,\xi_{1_{\,|_{\,n}}}\,)^{\,L}  \,\cdot\,\, \left\{    \,\left(  {{\lambda_{\,1}}\over {\lambda_{\,1}^2 \ + \ \Vert \, y \,-\, \xi_{\,1} \, \Vert^2}}\  \right)^{\!\!M} \,\cdot\,  \Vert \, y \,-\, \xi_{\,1} \, \Vert^N  \,\right\}\ dy \ \ (\ >  \ 0\ )\\[0.2in]
 & = &   \int_{B_{o}\,(\,{\rho_{\,\nu}})\ \cap \ \{ \ {\bar y}_n \ >  \ 0\,\}}  \ \ {\bar y}_n^{\,L}  \,\cdot\,\, \left\{    \,\left(  {{\lambda_{\,1}}\over {\lambda_{\,1}^2 \ + \ \Vert\  {\bar y}\ \Vert^2}}\  \right)^{\!\!M} \,\cdot\,  \Vert\ {\bar y}\,\Vert^N   \,\right\}\ d\,{\bar y} \ \ \ (\  >  \ 0\,) \ \ \ \  \ \ (\ {\bar y} \ = \ y \ - \ \xi_{\,1} \ )\ .
\end{eqnarray*}
Whereas the ``\ lower half\ " can be expressed as\\[0.1in]
 (\,A.9.31\,)
\begin{eqnarray*}
& \ & \int_{B_{\,\xi_{\,1}}\,(\,{\rho_{\,\nu}})\ \cap \ \{ \ y_n \ <  \ 0\,\}}    (\,-\,y_n\,+\,\xi_{1_{\,|_{\,n}}}\,)^{\,L}  \,\cdot\,\, \left\{    \,\left(  {{\lambda_{\,1}}\over {\lambda_{\,1}^2 \ + \ \Vert\  y \,-\, \xi_{\,1}\,\Vert^2}}\  \right)^{\!\!M} \!\!\cdot\,     \Vert \, y \,-\, \xi_{\,1} \, \Vert^N   \,\right\}\   \ \ ( \ \leftarrow \ {\mbox{+\,ve}}\ )\\[0.2in]
& = & \int_{B_{\,\xi_{\,1}}\,(\,{\rho_{\,\nu}})\ \cap \ \{ \ y_n \ <  \ 0\,\}}  \ \   (\,-\,y_n\,+\,\xi_{1_{\,|_{\,n}}}\,)^{\,L}  \,\cdot\,\, \left\{    \,\left(  {{\lambda_{\,1}}\over {\lambda_{\,1}^2 \ + \ \Vert\ -\,(-\,y \,+\, \xi_{\,1}\,)\,\Vert^2}}\  \right)^{\!\!M} \,\cdot\,     \Vert \, -\,(-\,y \,+\, \xi_{\,1}\,) \, \Vert^N   \,\right\} \\[0.2in]
& = & \int_{B_{\,\xi_{\,1}}\,(\,{\rho_{\,\nu}})\ \cap \ \{ \ y_n \ <  \ 0\,\}}  \ \   (\,-\,y_n\,+\,\xi_{1_{\,|_{\,n}}}\,)^{\,L}  \,\cdot\,\, \left\{    \,\left(  {{\lambda_{\,1}}\over {\lambda_{\,1}^2 \ + \ \Vert\ (-\,y \,+\, \xi_{\,1}\,)\,\Vert^2}}\  \right)^{\!\!M} \,\cdot\,     \Vert \, (-\,y \,+\, \xi_{\,1}\,) \, \Vert^N   \,\right\} \\[0.2in]
& = & \int_{B_{-\,\xi_{\,1}}\,(\,{\rho_{\,\nu}})\ \cap \ \{ \ z_n \ >  \ 0\ \}}  \ \   (\,z_n\,+\ \xi_{1_{\,|_{\,n}}}\,)^{\,L}  \,\cdot\,\, \left\{    \,\left(  {{\lambda_{\,1}}\over {\lambda_{\,1}^2 \ + \ \Vert\ (\,z \,+\, \xi_{\,1}\,)\,\Vert^2}}\  \right)^{\!\!M} \,\cdot\,     \Vert \, (\,z \,+\, \xi_{\,1}\,) \, \Vert^N   \,\right\}\ (\,+\,1\,)\  dz\\[0.2in]
& \ & \  \ \ \ \ \ \ \ \ \ \ \ \  \ \ \ \ (\ z_1\ = \ y_1\,, \ \cdot\,\cdot \cdot\,, \ \ z_{n\,-\,1}\ = \ y_{n\,-\,1}\,, \ \ z_n \ = \ -\,y_n\  \ \   \  \ \ \ \ \ \ \ \  \   |\,{\mbox{Det}}\,J\,| \ = \ +\,1 \ \ \uparrow\ \, ) \ \ \ \\[0.3in]
& = & \int_{B_{\,o}\,(\,{\rho_{\,\nu}})\ \cap \ \{ \ z_n \ >  \ \xi_{1_n}\}}  \ \   {\bar z}_n^{\,L}  \,\cdot\,\, \left\{    \,\left(  {{\lambda_{\,1}}\over {\lambda_{\,1}^2 \ + \ \Vert\ {\bar z}\,\Vert^2}}\  \right)^{\!\!M} \,\cdot\,     \Vert \, {\bar z} \, \Vert^N   \,\right\}\   d\,{\bar z} \\[0.1in]
& \ & \hspace*{4.3in}  \uparrow  \ \ \ \ \ (\,{\bar z}\ =\ z \ + \ \xi_{\,1}\ )\\[0.2in]
& = & \int_{B_{\,o}\,(\,{\rho_{\,\nu}})\ \cap \ \{ \ {\bar z}_n \ >  \ 0\,\}}  \ \   {\bar z}_n^{\,L}  \,\cdot\,\, \left\{    \,\left(  {{\lambda_{\,1}}\over {\lambda_{\,1}^2 \ + \ \Vert\ {\bar z}\,\Vert^2}}\  \right)^{\!\!M} \,\cdot\,     \Vert \, {\bar z} \  \Vert^N   \,\right\}\   d\,{\bar z}\\[0.2in]
&  & \ \ \ \ \ \ \ - \  \int_{B_{\,o}\,(\,{\rho_{\,\nu}})\ \cap \ \{ \ 0  \ < \ {\bar z}_n \ \le  \ \xi_{1_n}\}}  \ \   {\bar z}_n^{\,L}  \,\cdot\,\, \left\{    \,\left(  {{\lambda_{\,1}}\over {\lambda_{\,1}^2 \ + \ \Vert\ {\bar z}\,\,\Vert^2}}\  \right)^{\!\!M} \,\cdot\,     \Vert \, {\bar z}\, \, \Vert^N   \,\right\}\   d\,{\bar z} \ .\\[0.1in]
& \ & \hspace*{0.2in} \ \ \ \ \  [ \ \ \longleftarrow \hspace*{1.75in} \ *_{\,(\,A.9.31\,)} \hspace*{1.75in} \longrightarrow  \ ]
\end{eqnarray*}
In performing the change of variables in (\,A.9.31\,)\,,\, we note that the boundary $\,\partial B_{\,\xi_{\,1}}\,(\,{\rho_{\,\nu}}\,)\,$ is governed by the equation
\begin{eqnarray*}
& \ & y_1^2 \ + \ \cdot\,\cdot\,\cdot\,\ + \ y_{n\,-\,1}^2 \ + \ (\,y_n \ - \ \xi_{1_{\,|_{\,n}}}\,)^{\,2} \ = \ \rho_{\,\nu}^2\\[0.2in]
& \Longleftrightarrow &  y_1^2 \ + \ \cdot\,\cdot\,\cdot\,\ + \ y_{n\,-\,1}^2 \ + \ (\,-\,[\ -\,y_n \ + \ \xi_{1_{\,|_{\,n}}}\,]\ )^{\,2} \ = \ \rho_{\,\nu}^2\\[0.2in]
& \Longleftrightarrow &  z_1^2 \ + \ \cdot\,\cdot\,\cdot\,\ + \ z_{n\,-\,1}^2 \ + \ (\,z_n \ + \ \xi_{1_{\,|_{\,n}}}\,)^{\,2} \ = \ \rho_{\,\nu}^2\ .
\end{eqnarray*}
It is convenient to record it here\,:
$$
\widetilde{[\ \cdot\,\cdot\,\cdot\,\ ]}_{\,(\,A.9.9\,)}  \ = \ \left(\ {\lambda_{\,1}\over {\lambda_{\,1}^2 \ + \ \Vert\ y\,-\,\xi_{\,1}\,\Vert^2}}\  \right)^{\!\!n }  \,\cdot\, \left(\ {{\lambda_{\,1}^2 \ - \ \Vert\ y\,-\,\xi_{\,1}\,\Vert^2}\over {\lambda_{\,1}^2 \ + \ \Vert\ y\,-\,\xi_{\,1}\,\Vert^2}}\  \right)\ .
$$
From the consideration in (\,A.9.30\,) and (\,A.9.31\,)\,,\, the ``\,upper half\,"
$$
\int_{B_{\,\xi_{\,1}}\,(\,\rho_{\,\nu}\,) \ \cap \ \{ \, \,y_n\ > \ \xi_{1_{\,|_{\,n}}}\, \}  } (\,y_n\,-\,\xi_{1_{\,|_{\,n}}}\,)^{\,\ell\,-\ 1}\cdot \xi_{1_{{\,n}}}  \cdot\, \widetilde{[\ \cdot\,\cdot\,\cdot\,\ ]}_{\,(\,A.9.9\,)}  \ dy \ ,
$$
and the ``\,lower half\,"
\begin{eqnarray*}
& \ & \int_{B_{\,\xi_{\,1}}\,(\,\rho_{\,\nu}\,) \ \cap \ \{ \, \,y_n\ <  \ 0\, \}  } \ (\,-\,y_n\,+\,\xi_{1_{\,|_{\,n}}}\,)^{\,\ell\,-\ 1}\cdot (\,-\,\xi_{1_{{\,n}}}\,)  \cdot\,\widetilde{[\ \cdot\,\cdot\,\cdot\,\ ]}_{\,(\,A.9.9\,)} \ dy \\[0.2in]
& = & -\,\left[ \ \int_{B_{\,\xi_{\,1}}\,(\,\rho_{\,\nu}\,) \ \cap \ \{ \, \,y_n\ <  \ 0\, \}  } \ (\,-\,y_n\,+\,\xi_{1_{\,|_{\,n}}}\,)^{\,\ell\,-\ 1}\cdot \xi_{1_{{\,n}}} \cdot\,\widetilde{[\ \cdot\,\cdot\,\cdot\,\ ]}_{\,(\,A.9.9\,)} \ dy\ \right]
\end{eqnarray*}
cancel out, leaving behind the  "\,thin layer\," [\,marked as  $*_{\,(\,A.9.31\,)}$\, in (\,A.9.31\,)\,,\, with $\,L \ = \ \ell \ - \ 1$\ )\,,\, which is estimated by observing that
$$
\Bigg\vert \  {{\lambda_{\,1}^2 \ - \ \Vert\ y\,-\,\xi_{\,1}\,\Vert^2}\over {\lambda_{\,1}^2 \ + \ \Vert\ y\,-\,\xi_{\,1}\,\Vert^2}}\   \Bigg\vert \ \le  \ 1\ \ \ \ \ \mfor \ \ \ \ y \, \in\, \R^n\ ,  \leqno (\,A.9.32\,)
$$
and
\\[0.1in]
(\,A.9.33\,)
\begin{eqnarray*}
& \ &
\xi_{1_n} \cdot \int_{B_{\,o}\,(\,{\rho_{\,\nu}})\ \cap \ \{ \ 0  \ < \ {\tilde z}_n \ \le  \ \xi_{1_n}\}}  \ \   {\bar z}_n^{\,\ell \ - \ 1}  \,\cdot\,\,    \,\left(  {{\lambda_{\,1}}\over {\lambda_{\,1}^2 \ + \ \Vert\ {\bar z}\,\Vert^2}}\  \right)^{\!\!n} \  d\,{\bar z} \\[0.2in]
& \le & \xi_{1_n}^\ell \cdot \int_{B_{\,o}\,(\,{\rho_{\,\nu}})\ \cap \ \{ \ 0  \ < \ {\tilde z}_n \ \le  \ \xi_{1_n}\}}  \ \   \,\left(  {{\lambda_{\,1}}\over {\lambda_{\,1}^2 \ + \ \Vert\ {\bar z}\,\Vert^2}}\  \right)^{\!\!n} \  d\,{\bar z} \\[0.2in]
& \le & C_1 \cdot \xi_{1_n}^\ell  \cdot (\,\lambda_{\,1}^{-1} \,\cdot \,\xi_{1_n} \,)\\[0.2in]
& \le & C_2 \cdot {\bar\lambda}_{\ \flat}^\ell \cdot {\bar\lambda}_{\ \flat}^{\,(\,\ell \ + \ 1\,)\, \cdot \, \kappa} \\[0.2in]
\Longrightarrow & \ & \bigg\vert \
\int_{B_{\,\xi_{\,1}}\,(\,\rho_{\,\nu}\,)  } (\,y_n\,-\,\xi_{1_{\,|_{\,n}}}\,)^{\,\ell\,-\ 1}\cdot \xi_{1_{{\,n}}}  \cdot\, \widetilde{[\ \cdot\,\cdot\,\cdot\,\ ]}_{\,(\,A.9.9\,)}  \ dy \  \bigg\vert \ \le \ C_2 \cdot {\bar\lambda}_{\ \flat}^\ell \cdot {\bar\lambda}_{\ \flat}^{\,(\,\ell \ + \ 1\,) \,\cdot\, \kappa}\ \,.
\end{eqnarray*}

\vspace*{0.3in}

 {\bf \S\,A\,9.\,g\,.}\ \  {\bf Estimate the errors coming from various lower order terms in}\, (\,4.44\,) {\bf and}\, (\,4.45\,) {\bf of the main text}\,. \ \ Via (\,A.9.32\,)\,,\,  let us consider
\begin{eqnarray*}
(\,A.9.34\,) \ \ \ \ & \ & \Bigg\vert \ \int_{B_{\,\xi_{\,1}}\,(\,\rho_{\,\nu}\,)  } \ (\,-\,y_n\,+\,\xi_{1_{\,|_{\,n}}}\,)^{\,\ell\,-\ 1}\cdot (\,\xi_{1_{{\,n}}}\,)^{\,2}  \cdot\,| \ \widetilde{[\ \cdot\,\cdot\,\cdot\,\ ]}_{\,(\,A.9.9\,)} \ | \ dy\ \Bigg\vert\ \ \ \  \ \ \ \  \ \ \ \  \ \ \ \   \\[0.2in]
& \le & C_3 \cdot (\,\xi_{1_{{\,n}}}\,)^{\,2}  \cdot\, \lambda_{\,1}^{\,\ell \ - \ 1} \cdot \int_0^\infty {{ R^{\,\,\ell \ - \ 1} \ \cdot \ R^{\,\,n\ - \ 1}  }\over {(\,1 \ + \ R^2\,)^n  }} \ dR\\[0.2in]
& \le & C_4 \cdot {\bar\lambda}_{\ \flat}^{\,\ell \ - \ 1} \cdot {\bar\lambda}_{\ \flat}^{\,2\  + \ 2\,\cdot\,\kappa}\\[0.2in]
& = & C_4 \cdot {\bar\lambda}_{\ \flat}^{\,\ell \ + \ 1} \cdot {\bar\lambda}_{\ \flat}^{\,2\,\cdot \,\kappa} \ .
 \end{eqnarray*}
 Similarly\,,\,
\begin{eqnarray*}
 (\,A.9.35\,) \ \ \ \  \ \ \ & \ & \bigg\vert \ \int_{B_{\,\xi_{\,1}}\,(\,\rho_{\,\nu}\,)  } \ (\,-\,y_n\,+\,\xi_{1_{\,|_{\,n}}}\,)^{\,\ell\,-\ 1}\cdot \,\Vert\,y \ - \ \xi_{\,1}\,\Vert^2  \cdot\,| \ \widetilde{[\ \cdot\,\cdot\,\cdot\,\ ]}_{\,(\,A.9.9\,)}  \ | \ dy \ \bigg\vert \ \ \ \ \ \ \ \ \  \\[0.2in]
& \le & \bigg\vert \ \int_{B_{\,\xi_{\,1}}\,(\,\rho_{\,\nu}\,)  }  \,\Vert\,y \ - \ \xi_{\,1}\,\Vert^{\,\ell \ + \ 1}  \cdot\,| \ \widetilde{[\ \cdot\,\cdot\,\cdot\,\ ]}_{\,(\,A.9.9\,)}  \ | \ dy \ \bigg\vert \\[0.2in]
& \le & C_5 \cdot {\bar\lambda}_{\ \flat}^{\,\ell \ + \ 1} \ .
 \end{eqnarray*}

Likewise,
\begin{eqnarray*}
 (\,A.9.36\,) \ \ \ \  \ \ \  & \ & \bigg\vert \ \int_{B_{\,\xi_{\,1}}\,(\,\rho_{\,\nu}\,)  } \ (\,-\,y_n\,+\,\xi_{1_{\,|_{\,n}}}\,)^{\,\ell\,-\ 2}\cdot (\,\xi_{1_{{\,n}}}\,)^3  \cdot\,| \ \widetilde{[\ \cdot\,\cdot\,\cdot\,\ ]}_{\,(\,A.9.9\,)}  \ | \ dy \ \bigg\vert \ \ \ \ \ \ \ \ \  \ \ \ \  \ \ \  \\[0.2in]
& \le & C_6 \cdot (\,\xi_{1_{{\,n}}}\,)^{\,3}  \cdot\, \lambda_{\,1}^{\,\ell \ - \ 2} \cdot \int_0^\infty {{
 R^{\,\,\ell \ - \ 1} \ \cdot \ R^{\,\,n\ - \ 1}  }\over {(\,1 \ + \ R^2\,)^n  }} \ dR\\[0.2in]
& \le & C_7 \cdot {\bar\lambda}_{\ \flat}^{\,\ell  \ + \ 1} \cdot {\bar\lambda}_{\ \flat}^{\,3\,\cdot\, \kappa} \ ,
 \end{eqnarray*}
 and
\begin{eqnarray*}
 (\,A.9.37\,) \ \ \ \  \ \ \ & \ & \bigg\vert \ \int_{B_{\,\xi_{\,1}}\,(\,\rho_{\,\nu}\,)  } \ (\,-\,y_n\,+\,\xi_{1_{\,|_{\,n}}}\,)^{\,\ell\,-\ 2}\cdot \,\Vert\,y \ - \ \xi_{\,1}\,\Vert^3  \cdot\,| \ \widetilde{[\ \cdot\,\cdot\,\cdot\,\ ]}_{\,(\,A.9.9\,)}  \ | \ dy \ \bigg\vert \ \ \ \ \ \ \ \ \  \ \ \ \  \ \ \  \\[0.2in]
& \le & \bigg\vert \ \int_{B_{\,\xi_{\,1}}\,(\,\rho_{\,\nu}\,)  }  \,\Vert\,y \ - \ \xi_{\,1}\,\Vert^{\,\ell \ + \ 1}  \cdot\,| \ \widetilde{[\ \cdot\,\cdot\,\cdot\,\ ]}_{\,(\,A.9.9\,)}    \ dy \ \bigg\vert \\[0.2in]
& \le & C_8 \cdot {\bar\lambda}_{\ \flat}^{\,\ell \ + \ 1} \ .
 \end{eqnarray*}
 In this familiar fashion, we estimate the remaining higher order terms\,:
\begin{eqnarray*}
 (\,A.9.38\,) \ \ \ \  \ \ \  \  & \ & \int_{B_{\,\xi_{\,1}}\,(\,\rho_{\,\nu}\,)  }  \,\Vert\,y \ - \ \xi_{\,1}\,\Vert^{\,\ell \ - \ 3}\cdot (\,\xi_{1_{{\,n}}}\,)^3  \cdot\,| \ \widetilde{[\ \cdot\,\cdot\,\cdot\,\ ]}_{\,(\,A.9.9\,)}  \ | \ dy \ \ \ \ \ \ \ \ \  \ \ \ \  \ \ \ \ \ \ \ \ \ \ \ \ \ \  \\[0.2in]
& \le & C_5 \cdot (\,\xi_{1_{{\,n}}}\,)^3  \cdot\, \lambda_{\,1}^{\,\ell \ - \ 3} \cdot \int_0^\infty {{
 R^{\,\,\ell \ - \ 1} \ \cdot \ R^{\,\,n\ - \ 1}  }\over {(\,1 \ + \ R^2\,)^n  }} \ dR\\[0.2in]
& \le & C_6 \cdot {\bar\lambda}_{\ \flat}^{\,\ell  \emph{}} \cdot {\bar\lambda}_{\ \flat}^{\,3\,\cdot\, \kappa} \
 \end{eqnarray*}
and
\begin{eqnarray*}
(\,A.9.39\,) \ \ \ \  \ \ \  \   & \ & \bigg\vert \ \int_{B_{\,\xi_{\,1}}\,(\,\rho_{\,\nu}\,)  }  \,\Vert\,y \ - \ \xi_{\,1}\,\Vert^{\,\ell \ - \ 3}\cdot \,\Vert\,y \ - \ \xi_{\,1}\,\Vert^6  \cdot\,| \ \widetilde{[\ \cdot\,\cdot\,\cdot\,\ ]}_{\,(\,A.9.9\,)}  \ | \ dy \ \bigg\vert\ \ \ \ \ \ \ \ \  \ \ \ \  \ \ \ \ \ \ \ \ \ \ \ \ \ \  \\[0.2in]
& \le & \bigg\vert \ \int_{B_{\,\xi_{\,1}}\,(\,\rho_{\,\nu}\,)  }  \,\Vert\,y \ - \ \xi_{\,1}\,\Vert^{\,\ell \ + \ 3}  \cdot\,| \ \widetilde{[\ \cdot\,\cdot\,\cdot\,\ ]}_{\,(\,A.9.9\,)}  \ | \ dy \ \bigg\vert \\[0.2in]
& \le & C_4 \cdot {\bar\lambda}_{\ \flat}^{\,\ell \ + \ 3} \ .
 \end{eqnarray*}


\newpage

{\bf \S\,A\,9\,.h\,.}\ \
 To obtain the other estimate of the lower order teams in (\,4.44\,) and (\,4.45\,) of the main text\,,\, we may replace
 $$
 \ell \ \to \ (\ \ell \ - \ 2\,)
 $$
in the argument found in  {\bf \S\,A\,9.\,g\,} \,.

\vspace*{0.3in}

 {\bf \S\,A\,9\,.\,i\,.}\ \   {\bf Estimate of} $\,|\ {\bf E}_{\, (\,A.9.11\,) }^{\mbox{Lower}}\  | \,$. \ \ Combining the discussion in  {\bf \S\,A\,9.\,f\,}\,,\, {\bf \S\,A\,9.\,g\,} and {\bf \S\,A\,9.\,h\,},\, we see that\\[0.1in]
 (\,A.9.40\,)
\begin{eqnarray*}
\bigg|\ {\bf E}_{\, (\,A.9.11\,) }^{\,\mbox{Lower}}\  \bigg| &\le &  C_1 \cdot {\bar\lambda}_{\ \flat}^{\,\ell  } \cdot {\bar\lambda}_{\ \flat}^{\,(\,\ell \ + \ 1\,)\, \cdot\, \kappa}\ + \ C_2 \cdot (\,{\bar\lambda}_{\ \flat}^{\ 2\,\nu}\ )^{\,\ell} \ \cdot\,\
{\bar\lambda}_{\ \flat}^{\ 2\,\nu  \ - \ 1}  \ + \ C_3 \cdot {\bar\lambda}_{\ \flat}^{\,\ell \ +\ 1  }  \ + \  \\[0.2in]
 & \ & \hspace*{3.7in} + \ C_4 \cdot {\bar\lambda}_{\ \flat}^{\,\ell \ +\  3\,\kappa  } \ .
\end{eqnarray*}

\vspace*{0.3in}

{\bf \S\,A\,9\,.\,j\,.}\ \ We combine the discussion in this section so far [ \ specifically\,,\, (\,A.9.4\,)\,,\,  (\,A.9.15\,)\,,\, (\,A.9.17\,)\,,\, (\,A.9.18\,)\,,\, (\,A.9.21\,)\, \&  \,(\,A.9.29\,)\ ]\,,\, and come up with \\[0.1in]
(\,A.9.41\,)
\begin{eqnarray*}
& \ & - \   \left[\ \int_{\R^n}\, [\ (\,{\tilde c}_n\!\cdot K\,) \ - \ n\,(\,n\,-\,2)\ ]  \cdot\,(\,W_{\,\,\flat} \,)^{{\,n\,+\,2\,}\over {n\,-\,2}}\,\cdot\,\left(\,\lambda_{\,1}\,\cdot\,{ { \partial \,V_1 } \over {\partial\, \lambda_{\,1}}}\ \right)  \ \right]\\[0.2in]
&  = &   \Bigg[\,+ \,  {\hat C}_{\,1\,, \ 1} \cdot C\,(\,{\bf p}_{\,\xi_{\,1}}\,)\,\cdot\, \lambda^{\,\ell}_1 \ + \     {\hat C}_{\,1\,, \ 3} \cdot C\,(\,{\bf p}_{\,\xi_{\,1}}\,)\,\cdot\,\lambda^{\,\ell\,-\,2}_1\,\cdot\, \eta_1^2    \ +  \\[0.2in]
     & \ & \!\!\!\!\!\!\!+ \     O \,
 \left( \ {\bar \lambda}_{\,\,\flat}^{\,\ell}\ \right)  \,\cdot\,O \,
 \left( \ {\bar \lambda}_{\,\,\flat}^{(\,n\,-\ \ell\,)\, \cdot\, (\,1\ - \ \nu) }\, \right) \ + \     O \,
 \left( \ {\bar \lambda}_{\,\,\flat}^{\,\ell}\ \right)  \,\cdot\,O \,
 \left( \ {\bar \lambda}_{\,\,\flat}^{\,(\ \ell\,+\,1\,)\, \cdot\, (\  2\,\nu\,-\,1) }\, \right)
  \ + \ O\,\left(\,{\bar \lambda}_{\,\,\flat}^{\,\ell \ + \ 1}\, \right) \ + \\[0.2in] & \ &  \ \ \ \ \ \ \ \ \ \ \ \ \ \ \ \ \ \ \ \ \ \ \  \ \  \ + \ O\,\left(\,{\bar \lambda}_{\,\,\flat}^{\,\ell \ + \ 3 \cdot \kappa}\, \right) \ \Bigg]\ \cdot \ \left[\,1 \ + \   O\,\left( \ {\bar\lambda}_{\ \flat}^{\ (\,n\ -\,2\,)\,\cdot\, [\ (\,\gamma \ + \ \nu \,)\ - \  1\,] }\ \,\right) \ \right]\\[0.2in]
& \ & \hspace*{1.3in}  + \  \    O\left(\   {\bar\lambda}_{\ \flat}^{ \  n\,\cdot\,(\,1\ - \ \nu\,)  \ - \ o_{\,+}\,(\,1\,)}\ \right) \ \ \ \ \ \ \ \ [ \ {\mbox{note \ \ that \ \ }} \ell \ \ge \ 2\ ].
\end{eqnarray*}
In the above\,,\,
$$
{\hat C}_{\,1\,,\ 1}\ = \ \left(\ {{\,n\,-\ 2\,}\over {2}}\,\right)\,\cdot\,\int_{\R^n} |\,Y_n\,|^{\,\ell} \cdot  \left(\ {1\over { 1 \ + \ \Vert\ Y\,\Vert^2 }}\ \right)^{\!\!n} \!\cdot \left(\, {{\ \Vert\ Y\,\Vert^2 \ - \ 1 }\over { \ \  \Vert\ Y\,\Vert^2 \  + \ 1  }}\, \right)\, d\,Y\ ,
$$
and
$$
  {\hat C}_{\,1\,, \ 3} \ = \   \left({{n\,-\,2}\over {2}}\,\right)   \cdot {{\ell \cdot (\,\ell\, - \, 1\,)}\over 2} \cdot \int_{\R^n} |\,Y_n\,|^{\,\,\ell\, - \, 2} \cdot \left(\, {1\over { 1 \, + \, \Vert\,Y\,\Vert^2 }}\right)^{\!\!n}  \!\cdot \left(\, {{\Vert\ Y\,\Vert^2 \, - \, 1 }\over { \ \  \Vert\ Y\,\Vert^2 \, + \, 1 \ }}\, \right) \, d\,{\,Y}\ .
  $$

Recall from (1.28\,) of  the main text  that
\begin{eqnarray*}
\eta_1 & = & {\mbox{Dist.}}\ (\,\xi_{\,1}\,, \ {\cal H}\,)  \ \le \   \bar\lambda_{\,\flat}^{\ 1\,+\,\kappa} \ .
 \end{eqnarray*}
Moreover,
$$
o_{\,+}\,(\,1\,) \ = \ 2\,\varepsilon \,[\ (\,\gamma \ + \ \nu\,) \ - \ 1\,]\ .
$$

\vspace*{0.3in}

{\bf{\S\,A\,9\,.k}}  \ \
 As a conclusion, we obtain\\[0.1in]
 (\,A.9.42\,)
\begin{eqnarray*}
& \ & \left(\ \lambda_{\,1}\cdot {\partial\over {\partial\, \lambda_{\,1}}}\ \right)  \, {\bf I}_{\,\cal R} \,\left( \ (\,\lambda_{\,1}\,, \ \xi_{\,1}\,)\,, \ \cdot \cdot \cdot\,, \ (\,\lambda_{\,\flat}\,, \ \xi_{\ \flat} \,) \ \right) \\[0.2in]  & = &    -\, {\hat C}_{\,1\,, \ 2} \cdot \left\{ \ \  \sum_{l\,=\ 2}^\flat \left(\ {{\lambda_{\,1}^{{\,n\,-\ 2\,}\over 2}\,\cdot\,\lambda_{\ \!l}^{{\,n\,-\ 2\,}\over 2}}\over {\Vert\ \xi_{\,1}\ - \ \xi_{\,l}\,\Vert^{\,n\,-\,2} }}\  \right)
 \cdot\left[\,1 \ + \
  O\,\left(\,{\bar\lambda}_{\ \flat}^{\ 2\,(\,1\ -\ \nu\,)}\,\right)\ \right]\,\,\right\} \ +\  \\[0.2in]
 & \ &   + \ {\hat C}_{\,1\,,\ 1}\cdot C\,(\,{\bf p}_{\ \xi_{\,1}}\,)\,\cdot\, \lambda^{\,\ell}_1 \ + \   {\hat C}_{\,1\,,\ 3}\cdot C\,(\ {\bf p}_{\,\xi_{\,1}}\,)\,\cdot\,\lambda^{\,\ell\ -\,2}_1\,\cdot\, \eta_1^2  \ +  \ {\cal E}_{\,(\,A.{\bf 9}.42\,)}\ + \  {\cal E}_{\,(\,A.{\bf 7}.55\,)}\ ,
  \end{eqnarray*}
  where\\[0.1in]
  (\,A.9.43\,)
 \begin{eqnarray*}
  {\cal E}_{\,(\,A.{\bf 9}.42\,)}   & = &  O \,
 \left( \ {\bar \lambda}_{\,\,\flat}^{\,\ell \ + \ \mu_{\,\vert_{\cal K}}}\ \right) \ + \     O\left(\   {\bar\lambda}_{\ \flat}^{ \  n\cdot\,(\,1\ - \ \nu\,) \ - \ o_{\,+}\,(\,1\,)} \ \right)  \ + \ O\left( \  {\bar\lambda}_{\ \flat}^{ \,n \,\gamma \ - \ \sigma}\ \right)\ .
   \end{eqnarray*}

   Here\\[0.1in]
   (\,A.9.44\,)
   \begin{eqnarray*}
   \mu_{\,\vert_{\cal K}} & = &  \mbox{Min} \ \ \bigg\{ \ (\,n\,-\ \ell\,)\, \cdot\, (\,1 \ - \ \nu)\,, \ \ 1 \ -\ o_{\,{\bar\lambda}_{\ \flat}}\,(\,1\,)\,, \ \ (\,n\ -\ 2\,)\,\cdot\, [\ (\,\gamma \ + \ \nu \,)\ - \  1\,]\,, \\[0.2in]
   & \ &  \ \ \ \ \ \ \ \ \ \ \ \ \ \ \ \ \ \ \ \ \ \ \ \ \ \ \ \  \ \ \ \ \ \ \ \ \ \ \  \ \ \  \ \ \ \ \ \ \ \ \ \ \ \  \ \ \ \ \ \ \  (\,\ell \ + \ 1\,) \cdot (\,2\,\nu \ - \ 1\,) \ ,   \ \ \ \ 3\,\kappa\ \bigg\}\ ,
   \end{eqnarray*}

   $$    o_{\,+}\,(\,1\,) \ = \  2\,\varepsilon \,[\ (\,\gamma \ + \ \nu\,) \ - \ 1\ ]\ , \ \ \ \  {\mbox{and}} \ \ \ \ \
  o_{\,{\bar\lambda}_{\ \flat}}\,(\,1\,) \ \to \ 0^+ \ \ \ \ \ \ {\mbox{as}} \ \ \ \ {\bar\lambda}_{\ \flat} \ \to \ 0^+\ .
   $$
In (\,A.9.42\,)\,,\, the bubble parameters
 $\,\left( \ \lambda_{\,1}\,, \ \cdot \cdot \cdot\,, \ \lambda_{\,\flat}\,, \ \xi_{\,1}\,, \ \cdot \cdot \cdot\,, \ \xi_{\,\flat}   \ \right)\,$ satisfy the conditions in Theorem 1.33 of the main text\,.\, Recall that
\begin{eqnarray*}
{\hat C}_{\,1\,, \ 2} & = & {{ (\,n\,-\,2)^{\,2} }\over 2}\,\cdot\,\omega_n \ .
 \end{eqnarray*}
For other indices, one may   change the expression in (\,A.9.41\,)\,: $\,1 \ \to \ j \ (\,= \ 2\,,\, \cdot \cdot \cdot\,, \ \flat\,)$\,.\bk
For  the spherical symmetric case in which the   bubbles are arranged ``\,evenly\,"  close to a great circle of $\,S^{\,n\,-\,1}$\,,\, we refer to Proposition 3.2 in the work by Wei and Yan  \cite{Wei-Yan}\,.

\newpage

 {\it Note}\,.\, A.9.45\,.\ \  We observe that
 \begin{eqnarray*}
  \gamma \ + \ \nu \ >  \ 1 & \Longrightarrow & \gamma \ > \ 1 \ - \ \nu\\[0.2in]
 & \Longrightarrow &  \ \   {{n\,+\,2}\over 2}\,\,\cdot\,\,\gamma \ + \ {{\,n\,-\,2\,}\over 2}\,\cdot\,(\,1\ - \ \nu\,) \ > \  n\cdot\,(\,1\ - \ \nu\,) \\[0.2in]
 \Longrightarrow  \ \   O\left(\   {\bar\lambda}_{\ \flat}^{ \  n\cdot\,(\,1\ - \ \nu\,) \ - \ o_{\,+}\,(\,1\,)} \ \right)  & > &  O\left(\   {\bar\lambda}_{\ \flat}^{ \   {{n\,+\,2}\over 2}\,\,\cdot\,\,\gamma \ + \ {{\,n\,-\,2\,}\over 2}\,\cdot\,(\,1\ - \ \nu\,)  } \,\cdot\, {\bar\lambda}_{\ \flat}^{  \   - \ \varepsilon \,[\ (\,\gamma \ + \ \nu\,) \ - \ 1\ ] }  \ \right)
  \end{eqnarray*}
   when $\,\epsilon \ > \ 0\,$ is made small enough.

\vspace*{0.5in}

{\large{\bf \S\,A\,\,10\,.  \ \ Derivative with respect to  $\,\xi_{\,1}\,.\,$}}\\[0.1in]
(\,A.10.1\,)
 \begin{eqnarray*}
& \ & \left(\,\lambda_{\,1}\cdot {\partial\over {\partial\, \xi_{1_{\,|_j}}}}\  \right) {\bf I}_{\,\cal R}\left( \ (\,\lambda_{\,1}\,, \ \xi_{\,1}\,)\,, \ \cdot \cdot \cdot\,, \ (\,\lambda_{\,\flat}\,, \ \xi_{\ \flat} \,) \ \right) \\[0.2in]  & = &  -\,\int_{\R^n} (\,\Delta\,W_{\,\,\flat})\,\cdot\,\left[ \,  \left(\,\lambda_{\,1}\cdot {\partial\over {\partial\, \xi_{1_{\,|_j}}}}\  \right)W_{\,\,\flat}     \right]  \ -  \   \int_{\R^n}\,(\,{\tilde c}_n\!\cdot K\,)\,(W_{\,\,\flat} )^{{\,n\,+\,2\,}\over {n\,-\,2}}\,\cdot\, \left[ \,  \left(\,\lambda_{\,1}\cdot {\partial\over {\partial\, \xi_{1_{\,|_j}}}}\  \right)W_{\,\,\flat}    \right] \\[0.2in]
& \ & \hspace*{5in} \ \ \ \ \ \ \ \ \ \ \ + \ \,{\cal E}_{\,(\,A.10.1\,)}    \\[0.2in]
  & = & \left\{ \ -\int_{\R^n} (\Delta\,W_{\,\,\flat})\,\cdot\,\left[ \,  \left(\,\lambda_{\,1}\cdot {\partial\over {\partial\, \xi_{1_{\,|_j}}}}\  \right)W_{\,\,\flat}     \right]  \right. \ \ \ \ \ \ \ \ \ \ \  \ \ \ \ \ \ \ \ \ \ \ \  \ \ \ (  \ {\mbox{refer \ \ to \ \ }} {\bf \S\,A\,8\,.\,g} \ \,\uparrow \ )
     \\[0.2in]
     & \  & \hspace*{2in}\left. \ -  \    \int_{\R^n}\, \,n\,(\,n\,-\,2)\cdot\,(W_{\,\,\flat} )^{{\,n\,+\,2\,}\over {n\,-\,2}}\,\cdot\, \left[ \,  \left(\,\lambda_{\,1}\cdot {\partial\over {\partial\, \xi_{1_{\,|_j}}}}\  \right)W_{\,\,\flat}     \right]  \ \right\} \ \\[0.2in]
  & \ & \hspace*{4in} \ \ \ \ \ \ \  \uparrow \ \ \cdot\,\cdot\,\cdot\,\,\cdot\, \ \ {\bf (\,C\,)_{\,(\,A.10.1\,)}}\\[0.2in]
 & \ &  \ \ \ \ \ \ \ \ \ \ \   - \   \left\{ \ \int_{\R^n}\, [\,(\,{\tilde c}_n\!\cdot K\,) \ - \ n\,(\,n\,-\,2)\,]  \cdot\,(W_{\,\,\flat} )^{{\,n\,+\,2\,}\over {n\,-\,2}}\,\cdot\, \left[ \,  \left(\,\lambda_{\,1}\cdot {\partial\over {\partial\, \xi_{1_{\,|_j}}}}\  \right)W_{\,\,\flat}     \right]   \ \right\} \ \\[0.2in]
  & \ & \hspace*{4in} \ \ \ \ \ \ \  \uparrow \ \ \cdot\,\cdot\,\cdot\,\,\cdot\, \ \ {\bf (\,D\,)}_{\,(\,A.10.1\,)}\\[0.2in]
& \ & \hspace*{5in} \ \ \ \ \ \ \ \ \ \ \ + \ \,  {\cal E}_{\,(\,A.10.1\,)}
 \end{eqnarray*}
 for $\, j \ = \ 1\,, \ \cdot \cdot \cdot\,,\ j\,.\,$
Recall that \\[0.1in]
(\,A.10.2\,)
 \begin{eqnarray*}
& \ &  \left(\,\lambda_{\,1}\,\cdot\,{\partial\over {\partial\, \xi_{1_{\,|_j}}}}\  \right)\,\left[\ \left({\lambda_{\,1}\over {\lambda^2_1\ + \ \Vert\  y \ -\  \xi_{\,1}\,\Vert^{\,2}}}\  \right)^{{n \,-\, 2}\over 2} \right] \ = \ -\,\left( \, {{n - 2}\over 2} \, \right)\,\cdot\,\lambda^{{n  }\over 2}_1\,\cdot\,{{ 2\,(\,\xi_{1_{\,|_j}}\, -\, y_j\,)}\over {(\,\lambda^2_1 + \Vert\  y \,-\, \xi_{\,1}\,\Vert^{\,2}\,)^{{n\over 2 } }}}\\[0.2in]
\Longrightarrow & \ & \Bigg\vert \ \left(\ \lambda_{\,1}\cdot { { \partial  } \over  {\partial\, \xi_{1_{\,|_j}}}}\ \right)V_1\,(\,y\,) \ \Bigg\vert \ \le \ C\,\cdot\,V_1\,(\,y\,) \ \ \ \  \ \ \ \ \ \ \  \ \ \ \  \ \ \ \ \ \ \ \ \ \ \ \ \ \mfor \ \ y \,\in\,\R^n\ .\\
\end{eqnarray*}

\vspace*{0.3in}

{\bf \S\,A\,10\,.a\,.} \ \ Let us consider $\,(\,{\bf C}\,)_{\,(\,A.10.1\,)}\,$ (\,this part is similar to  {\bf \S\,A 8}\,.\,)\\[0.1in]
(\,A.10.3\,)
 \begin{eqnarray*}
& \ & -\int_{\R^n} (\Delta\,W_{\,\,\flat})\cdot\left[ \,  \left(\,\lambda_{\,1}\cdot {\partial\over {\partial\, \xi_{1_{\,|_j}}}}\  \right)W_{\,\,\flat}     \right]    \, -     \int_{\R^n}\, \!\!n\,(\,n\,-\,2)\cdot\,(\,W_{\,\,\flat} )^{{\,n\,+\,2\,}\over {n\,-\,2}}\,\cdot\,\left[ \,  \left(\,\lambda_{\,1}\cdot {\partial\over {\partial\, \xi_{1_{\,|_j}}}}\  \right)W_{\,\,\flat}     \right]   \\[0.2in]
& = &  n\,(\,n\,-\,2)\,\cdot\,\int_{\R^n} \left\{  \,\left[\  V_1^{{\,n\,+\,2\,}\over {n\,-\,2}} \ + \  V_2^{{\,n\,+\,2\,}\over {n\,-\,2}} \ +\,\cdot\,\cdot\,\cdot\, + \ V_{\,\flat}^{{\,n\,+\,2\,}\over {n\,-\,2}} \ \right] \right. \\[0.2in]
    & \ & \left. \ \ \ \ \ \  \ \ \ \ \ \  \ \ \ \ \ \  \ \ \ \ \ \  \ \ \ \ \ \  \ \ \ \ \  \  \ \ \ \ \ \
    - \ \left[\  V_1\ + \  V_2  \ + \,\cdot\,\cdot\,\cdot\, + \ V_{\,\flat}  \ \right]^{{\,n\,+\,2\,}\over {n\,-\,2}} \  \right\}\,\cdot\,\left[ \,  \left(\,\lambda_{\,1}\cdot {\partial\over {\partial\, \xi_{1_{\,|_j}}}}\  \right)V_1    \right]\\[0.2in]
& = & -\,n\,(\,n\,-\,2)\,\cdot\,{{n\,+\,2}\over {\,n\,-\ 2\,}} \,\cdot\, \int_{B_{\,\xi_{\,1}} (\,{\rho_{\,\nu}}\,)} \left\{ \  V_1^{4\over {n\,-\,2}}\,\cdot\,\left[\   V_2  \ + \ \cdot\,\cdot\,\cdot\,\ + \ V_{\,\flat}  \ \right]   \right\}\,\cdot\, \left[ \,  \left(\,\lambda_{\,1}\cdot {\partial\over {\partial\, \xi_{1_{\,|_j}}}}\  \right) V_1 \right] \\[0.2in]
& \ & \hspace*{4.5in} \ \ \ \ \ \ \ \ \ \ \ \ \   + \ \,{\bf E}_{\,(\,A.10.3\,)}\,,
\end{eqnarray*}
where $\,{\bf E}_{\,(\,A.10.3\,)}\,$ takes on similar form as in {\bf \S\,A\,8.\,b}\,.\,
%
%
%
We proceed to the interaction team\,:
(\,A.10.4\,)
 \begin{eqnarray*}
  & \  &  -\,n\,(\,n\,-\,2) \cdot{{n\,+\,2}\over {\,n\,-\ 2\,}} \,\cdot\, \int_{B_{\,\xi_{\,1}} (\,{\rho_{\,\nu}}\,)} \left\{ \ \,\cdot\,V_1^{4\over {n\,-\,2}}\,\cdot\,\left[\   V_2     \ \right]\,  \right\}\,\cdot\,\left[ \,  \left(\,\lambda_{\,1}\cdot {\partial\over {\partial\, \xi_{1_{\,|_j}}}}\  \right) V_1    \right]\\[0.2in]
 & = & - \,n\,(\,n\,-\,2)\,\cdot\,{{n\,+\,2}\over {\,n\,-\ 2\,}}\cdot\left(- \,{{n\, - \, 2}\over 2}\,\right)\cdot\! \int_{B_{\,\xi_{\,1}} (\,{\rho_{\,\nu}}\,)} \!\!V_1^{4\over {n\,-\,2}} \cdot\left\{ \, \lambda_{\,1}^{{n  }\over 2}\,\cdot\,{{ 2\,(\,\xi_{1_{\,|_{\,1}}}\, -\, y_1\,)}\over {(\,\lambda_{\,1}^2 + \Vert\  y \,-\, \xi_{\,1}\,\Vert^{\,2}\,)^{{n\over 2 } }}} \, \right\}\cdot \! \left[\   V_2     \ \right] \\[0.3in]
 & = & - \,n\,(\,n\,-\,2)\,\cdot\,(\,n\,+\,2\,)\,\cdot\, \int_{B_{\,\xi_{\,1}} (\,{\rho_{\,\nu}})\,} V_1^{4\over {n\,-\,2}} \,\cdot\,\left\{ \ \lambda_{\,1}^{{n  }\over 2}\,\cdot\,{{  (\,y_1 \,-\,\xi_{1_{\,|_{\,1}}}  \,)}\over {(\,\lambda_{\,1}^2 + \Vert\  y \,-\, \xi_{\,1}\,\Vert^{\,2}\,)^{{n\over 2 } }}} \ \right\}\cdot  \left[\   V_2     \ \right] \ . \ \\[0.1in]
  & \ & \hspace*{4in} \ \ \ \ \ \ \  \uparrow \ \ \cdot\,\cdot\,\cdot\,\,\cdot\, \ \ {\bf (\,C\,)}_{\,(\,A.10.4\,)}
 \end{eqnarray*}
 With
 $$
 {\bf d}_{\,1\,, \ 2} \ =   {{\Vert\  \,\xi_{\,1}\ - \ \xi_{\,2}\,\Vert }\over
{\sqrt{\,\lambda_{\,1}\,\cdot\,\,\lambda_{\,2}\,}}} \ ,
 $$
 we have

 \begin{eqnarray*}
(\,A.10.5\,) \ \ \ \ \ \ \ & \ &  V_2 \,(y ) \ =\ V_{\lambda_{\,2}\,,\ \xi_{\,2}} (\,y)\\[0.2in]   & = & \left(\ {{\lambda_{\,2}}\over { \lambda^2_2 \ + \ \Vert\ y\ - \ \xi_{\,2}\,\Vert^2  }}\ \right)^{\!\!{{\,n\,-\,2\,}\over 2}  } \\[0.2in]   & = & \left(\ {{{1\over\lambda_{\,1}}}\over { \left(\ {{\lambda_{\,2}}\over {\lambda_{\,1}}}  \right)\ + \ {{\Vert\ y\ - \ \xi_{\,2}\,\Vert^2}\over
{\lambda_{\,1}\,\cdot\,\lambda_{\,2}}}  }}\ \right)^{\!\!{{\,n\,-\,2\,}\over 2}  }  \\[0.15in]
& = & \left(\ {1\over {\lambda_{\,1}^{{\,n\,-\,2\,}\over 2} }}\  \right)\,\cdot\,\left(\,{{{1}}
\over { \left(\ {{\lambda_{\,2}}\over {\lambda_{\,1}}}  \,\right)\ + \ {{\Vert\ (\,y\ - \ \xi_{\,1})\ + \ (\,\xi_{\,1}\ - \ \xi_{\,2})\,\Vert^2}\over
{\lambda_{\,1}\,\cdot\,\,\lambda_{\,2}}}  }}\ \right)^{\!\!{{{\,n\,-\,2\,}\over 2}  }} \\[0.15in]
& = & \left(\ {1\over {\lambda_{\,1}^{{\,n\,-\,2\,}\over 2} }}\  \right)\cdot \left[\ {{{1}}
\over { \left(\ {{\lambda_{\,2}}\over {\lambda_{\,1}}}  \right) \ + \ {{\ \Vert\  \,\xi_{\,1}\ - \ \xi_{\,2}\,\Vert^2\ }\over
{\lambda_{\,1}\,\cdot\,\,\lambda_{\,2}}}   \ + \ {{\ \Vert\ \,y\ - \ \xi_{\,1}\,\Vert^2\ }\over
{\lambda_{\,1}\,\cdot\,\,\lambda_{\,2}}}   \ + \ {{\,2\,(\,y\ - \ \xi_{\,1})\,\,\times\,(\,\xi_{\,1}\ - \ \xi_{\,2})\, }\over
{\lambda_{\,1}\,\cdot\,\,\lambda_{\,2}}}}}\ \right]^{{{\,n\,-\,2\,}\over 2}  }\\[0.1in]
& \ & \hspace*{1.5in} \ \ \ \ \{\  \uparrow \ = \ {\bf d}_{\,1\,, \ 2}^2 \ [\ \gg \ 1\,; \ \ {\mbox{cf.}} \ \ {\bf{\S\,A 8}}\ ]\,,\ {\mbox{dominating \ \ term}}\,\}\\[0.01in]
&  \ & \hspace*{-1.25in} = \ \left(\ {1\over {\lambda_{\,1}^{{\,n\,-\,2\,}\over 2} }}\  \right)\,\cdot\,{1\over { {\bf d}_{\,1\,, \ 2}^{n\,-\,2} }}\,*\\[0.2in]
&   & \hspace*{-1.25in} \ \ \ \ \ \ \ \ \ \ \ \ \ \ \ \  * \left[\ {{{1}}
\over { 1 \ + \ {1\over { {\bf d}_{\,1\,, \ 2}^{2} }}\,\cdot\,\left(\ {{\lambda_{\,2}}\over {\lambda_{\,1}}}  \right) \ + \ {1\over { {\bf d}_{\,1\,, \ 2}^{2} }}\,\cdot\,{{\Vert\ \,y\ - \ \xi_{\,1}\,\Vert^2}\over
{\lambda_{\,1}\,\cdot\,\,\lambda_{\,2}}}   \ + \ {1\over { {\bf d}_{\,1\,, \ 2}^{2} }}\,\cdot\,{{\,2\,(\,y\ - \ \xi_{\,1})\,\,\times\,(\,\xi_{\,1}\ - \ \xi_{\,2})\, }\over
{\lambda_{\,1}\,\cdot\,\,\lambda_{\,2}}}}}\ \right]^{{{\,n\,-\,2\,}\over 2}  }  {\bf .} \\
\end{eqnarray*}
As in {\bf \S\,A\,8}\,,\, for simplicity sake, we denote
$$\ \
{\bf T} \ = \  {1\over { {\bf d}_{\,1\,, \ 2}^{2} }}\,\cdot\,\left(\ {{\lambda_{\,2}}\over {\lambda_{\,1}}}  \right) \ + \ {1\over { {\bf d}_{\,1\,, \ 2}^{2} }}\,\cdot\,{{\Vert\ \,y\ - \ \xi_{\,1}\,\Vert^2}\over
{\lambda_{\,1}\,\cdot\,\,\lambda_{\,2}}}   \ + \ {1\over { {\bf d}_{\,1\,, \ 2}^{2} }}\,\cdot\,{{\,2\,(\,y\ - \ \xi_{\,1})\cdot(\,\xi_{\,1}\ - \ \xi_{\,2})\, }\over
{\lambda_{\,1}\,\cdot\,\,\lambda_{\,2}}}\ , \leqno (\,A.10.6\,)
$$
although in fact the expression depends on $\,y\,,\,$ $\lambda_{\,1}\,,\,$ $\lambda_{\,2}\,,\,$ $\xi_{\,1}\,$ and $\,\xi_{\,2}\,$.\,
Recall that (\,see  {\bf \S\,A\,8\,.\,c}\,)\,,\,
$$
{\bf T} \  = \  O\,\left(\,{\bar\lambda}_{\ \flat}^{\,2\,\gamma}\,\right) \ + \ O\,\left( \, {\bar\lambda}_{\ \flat}^{\,(\,\gamma\,+\,\nu\,) \ - \ 1} \ \right) \ . \leqno (\,A.10.7\,)
$$

Similar to {\bf \S\,A\,8}\,,\, we make use of the Taylor expansion
$$
\left( {1\over {1\,+\,{\bf T}}}\  \right)^{{\,n\,-\,2\,}\over 2} \ = \ 1 \ - \  \left(\,{{\,n\,-\,2\,}\over 2} \right)\,\cdot\,{\bf T}   \ + \    {1\over {2\,!}}\,\cdot\,\left(\,{{\,n\,-\,2\,}\over 2} \right)\,\cdot\, \left(\,{{n\,+\,2}\over 2} \right)\,\cdot\,{\bf T}^2 \ + \  O\,(\,|\, {\bf T} \,|^3\,)\ $$
for $\,|\,{\bf T}\,|\,$  being small. Contrasting the case in {\bf \S\,A\,8}\,,\, we pay special attention on the term involving $\displaystyle{\,(\,y_1\ - \ \xi_{1_{\,|_{\,1}}}\,)}\,.\,$ As such, we separate it from the expression for {\bf T} in (\,A.10.6\,)\,:\ \\[0.1in]
(\,A.10.8\,)
\begin{eqnarray*}
{\bf T} & = &  {1\over { {\bf d}_{\,1\,, \ 2}^{2} }}\,\cdot\,\left(\ {{\lambda_{\,2}}\over {\lambda_{\,1}}}  \right) \ + \ {1\over { {\bf d}_{\,1\,, \ 2}^{2} }}\,\cdot\,{{\Vert\ \,y\ - \ \xi_{\,1}\,\Vert^2}\over
{\lambda_{\,1}\,\cdot\,\,\lambda_{\,2}}}   \ + \ {1\over { {\bf d}_{\,1\,, \ 2}^{2} }} \cdot\! \left[\ \sum_{j\,=\,2}^n {{\,2\,(\,y_j\ - \ \xi_{1_{\,|_{\,j}}})\cdot(\,\xi_{1_{\,|_{\,j}}}\ - \ \xi_{2_{\,|_{\,j}}})\, }\over
{\lambda_{\,1}\,\cdot\,\,\lambda_{\,2}}} \ \right]\\[0.2in]
& \ &\hspace*{2.3in} +  \  \,{2\over { {\bf d}_{\,1\,, \ 2}^{2} }}\,\cdot\,\left(\ {{ y_1\ - \ \xi_{1_{\,|_{\,1}}} }\over
{\lambda_{\,1}  }}\  \right)\,\cdot\,\left(\ {{\, \xi_{1_{\,|_{\,1}}}\ - \ \xi_{2_{\,|_{\,1}}} \, }\over
{ \lambda_{\,2}}}\  \right) \ .
\end{eqnarray*}

We continue from (\,A.10.4\,)\\[0.1in]
(\,A.10.9\,)
 \begin{eqnarray*}
(\,{\bf C}'\,)_{\,(\,A.10.4\,)}
 & = & - \,n\,(\,n\,-\,2)\,\cdot\, (\,n\,+\,2)\,\cdot\,\lambda_{\,1}^2\,\cdot\,  {1\over { {\bf d}_{\,1\,, \ 2}^{n\,-\,2} }} \, \int_{B_{\,\xi_{\,1}} (\,{\rho_{\,\nu}}\,)}[\ \bullet\, \bullet\, \bullet \ ]_{\,(\,A.10.9\,)}  \,\cdot\,\left(\, {{y_1 \,-\,\xi_{1_{\,|_{\,1}}}}\over {\lambda_{\,1}}}  \right)\,*\\[0.2in]
 & * &\!\!\!   \left\{ \  1 \ - \  \left(\,{{\,n\,-\,2\,}\over 2} \right)\,\cdot\,{\bf T}     \ + \ {1\over {2\,!}}\,\cdot\,\left(\,{{\,n\,-\,2\,}\over 2} \right)\,\cdot\, \left(\,{{n\,+\,2}\over 2} \right)\,\cdot\, {\bf T}^2 \ + \  O\,(\,|\, {\bf T} \,|^{\,3}\,) \ \right\} \ .
\end{eqnarray*}
Here
$$
[\ \bullet\, \bullet\, \bullet \ ]_{\,(\,A.10.9\,)} \ = \ \left(\ {\lambda_{\,1}\over { \,\lambda_{\,1}^2 + \Vert\  y \,-\, \xi_{\,1}\,\Vert^{\,2}\, }}  \right)^{\!\!2}\,\cdot\, \left( \  {{1\,}\over {\,\lambda_{\,1}^2 + \Vert\  y \,-\, \xi_{\,1}\,\Vert^{\,2}\, }}\ \right)^{\! {n\over 2} }\ . \leqno (\,A.10.10\,)
$$
We observe the following symmetric cancellation,
$$
\int_{B_o\,(\,R\,)} f\,({\bar r})\,\cdot\,{\bar y}_1^{\,\alpha_1}\ \,\cdot\,\cdot\,\cdot\,\ {\bar y}_n^{\,\alpha_n}\ d{\bar y} \ = \ 0 \ , \ \ \ \ \ \ \ \ \  \ \ \,\bar r \ = \ \Vert\ \bar y\,\Vert\ \ , \leqno (\,A.10.11\,)
$$
where one (\,or more\,) of the indices (\,$\alpha_1\,,\ \cdot \cdot \cdot\,, \ \alpha_n\,$,\, non\,-\,negative integers\,) is {\it odd}\ . Here $\,f\,$ is a continuous function defined on $\,[\,0\,, \ \infty\,)\,.\,$    \bk
It follows that\\[0.1in]
(\,A.10.12\,)
\begin{eqnarray*}
& \ &   \int_{B_{\,\xi_{\,1}} (\,{\rho_{\,\nu}}\,)} \left(\ {\lambda_{\,1}\over { \,\lambda_{\,1}^2 + \Vert\  y \,-\, \xi_{\,1}\,\Vert^{\,2}\, }}  \right)^{\!\!2}\,\cdot\, \left( \  {{1\,}\over {\,\lambda_{\,1}^2 + \Vert\  y \,-\, \xi_{\,1}\,\Vert^{\,2}\, }}\ \right)^{\! {n\over 2} }   \,\cdot\,(\,y_1 \,-\,\xi_{1_{\,|_{\,1}}}  \,) \times 1\ dy\\[0.2in]
& = & \int_{B_{o}\, (\,{\rho_{\,\nu}})} \left(\ {\lambda_{\,1}\over { \,\lambda_{\,1}^2 + \Vert\  {\bar y} \, \Vert^{\,2}\, }}  \right)^{\!\!2}\,\cdot\,  {{1\,}\over {(\,\lambda_{\,1}^2\  + \ \Vert\  {\bar y} \, \Vert^{\,2}\,)^{{n\over 2 } }}}   \,\cdot\, {\bar y}_1  \times 1 \ d {\bar y}  \ = \ 0 \ \ \ \ \ \ \ \ (\,{\bar y}\ = \ y \ - \ \xi\,)\ ,
 \end{eqnarray*}

 and\\[0.1in]
(\,A.10.13\,)
 \begin{eqnarray*}
& \ &   \int_{B_{\,\xi_{\,1}} (\,{\rho_{\,\nu}}\,)} \left(\ {\lambda_{\,1}\over { \,\lambda_{\,1}^2 + \Vert\  y \,-\, \xi_{\,1}\,\Vert^{\,2}\, }}  \right)^{\!\!2}\,\cdot\, \left( \  {{1\,}\over {\,\lambda_{\,1}^2 + \Vert\  y \,-\, \xi_{\,1}\,\Vert^{\,2}\, }}\ \right)^{\! {n\over 2} }   \,\cdot\,(\,y_1 \,-\,\xi_{1_{\,|_{\,1}}}  \,) \times \Vert\ y\ - \ \xi_{\,1}\Vert^2 \ d\,y\\[0.2in]
& = & \int_{B_{o}\, (\,{\rho_{\,\nu}})} \left(\ {\lambda_{\,1}\over { \,\lambda_{\,1}^2 + \Vert\  {\bar y} \, \Vert^{\,2}\, }}  \right)^{\!\!2}\,\cdot\,  {{1\,}\over {(\,\lambda_{\,1}^2\  + \ \Vert\  {\bar y} \, \Vert^{\,2}\,)^{{n\over 2 } }}}   \,\cdot\, {\bar y}_1 \times \Vert\ {\bar y}\ \Vert^2 \ d {\bar y}  \ = \ 0\ .\\
 \end{eqnarray*}
Consider the zero and the first order terms in (\,A.10.9\,)\,:\\[0.1in]
(\,A.10.14\,)
 \begin{eqnarray*}
 &   & - \,n\,(\,n\,-\,2)\,\cdot\, (\,n\,+\,2)\,\cdot\,  {{\lambda_{\,1}^2}\over { {\bf d}_{\,1\,, \ 2}^{n\,-\,2} }}\,\cdot\, \int_{B_{\,\xi_{\,1}} (\,{\rho_{\,\nu}}\,)} [\ \bullet\, \bullet\, \bullet \ ]_{\,(\,A.10.9\,)}   \,*\\[0.2in]
  & \ & \hspace*{3in} *\ \left(\, {{y_1 \,-\,\xi_{1_{\,|_{\,1}}}}\over {\lambda_{\,1}}}  \right)\,\cdot\,\left\{ \ 1 \  - \  \left(\,{{\,n\,-\,2\,}\over 2}\  \right)\,\cdot\,{\bf T}  \ \right\} \\[0.2in]
 & = &0 \ + \ 0 \ + \ 0\\[0.1in]
  & \ & \ \ \ \ \ \   [\ \uparrow \ \ \ \ \ \ \uparrow \ \  {\mbox{cf.}} \ \ (\,A.10.8\,)\ ]\    \\[0.1in]
 & \ &  - \,n\,(\,n\,-\,2)\,\cdot\, (\,n\,+\,2)\,\cdot\,\lambda_{\,1}^2\,\cdot\,  {1\over { {\bf d}_{\,1\,, \ 2}^{n\,-\,2} }} \, \int_{B_{\,\xi_{\,1}} (\,{\rho_{\,\nu}}\,)}[\ \bullet\, \bullet\, \bullet \ ]_{\,(\,A.10.9\,)}  \,\cdot\,\left(\, {{y_1 \,-\,\xi_{1_{\,|_{\,1}}}}\over {\lambda_{\,1}}}  \right)\ \,* \\[0.2in]
 & \ & \ \ \ \ \ \  \ \ \ \ \ \ \  \ \ \ \ \ \ \  \ \ \ \ \ \ \  \ \,* \left[\  - \  \left(\,{{\,n\,-\,2\,}\over 2} \right)\,\cdot\,  {2\over { {\bf d}_{\,1\,, \ 2}^{2} }}\,\cdot\,\left(\ {{ y_1\ - \ \xi_{1_{\,|_{\,1}}} }\over
{\lambda_{\,1}  }}\  \right)\,\cdot\,\left(\ {{\, \xi_{1_{\,|_{\,1}}}\ - \ \xi_{2_{\,|_{\,1}}} \, }\over
{ \lambda_{\,2}}}\  \right) \ \right] \\ [0.2in]
   & = &  + \  n\,(\,n\,-\,2)^{\,2}\,\cdot\, (\,n\,+\,2\,)\!\cdot   {{\lambda_{\,1}^2}\over { {\bf d}_{\,1\,, \ 2}^{n\,-\,2\,+\,2} }} \,\left\{ \, \int_{B_{o}\, (\,{\rho_{\,\nu}})} \! [\ \cdot\,\cdot\,\cdot\,\ ]_{\,(\,A.10.14\,)}\cdot   \left(\  {{ {\bar y}_1 }\over { \lambda_{\,1}  }}\  \right)^{\!\!2} \ d\,{\bar y}\, \right\} \cdot\left(\, {{\, \xi_{1_{\,|_{\,1}}} - \, \xi_{2_{\,|_{\,1}}} \, }\over
{ \lambda_{\,2}}}\  \right)   \\[0.2in]
& = &   + \  n\,(\,n\,-\,2)^{\,2}\cdot\, (\,n\,+\,2)\,\cdot\!   {1\over { {\bf d}_{\,1\,, \ 2}^{n} }}   \,\cdot \left( {{\, \xi_{1_{\,|_{\,1}}} \! - \, \xi_{2_{\,|_{\,1}}}  }\over
{ \lambda_{\,2}}}\  \right) \cdot  \left\{   \lambda_{\,1}^2 \!\cdot \int_{B_{o}\, (\,{\rho_{\,\nu}})} \!\! [\ \cdot\,\cdot\,\cdot\,\ ]_{\,(\,A.10.14\,)}\cdot  \left(\  {{ {\bar y}_1 }\over { \lambda_{\,1}  }}\  \right)^{\!\!2} \ d\,{\bar y}\, \right\} \\[0.1in]
& \ &\hspace*{-0.5in}\  \left\{ \  \leftarrow \,   \left[\,    = \ O \left(\   {1\over { {\bf D}^{\,n\,-\,1} }}\  \right)\ \right] \ \,    {\mbox{``\,unit\," \ \, for \ \, this \ \,calculation}} \ \rightarrow\   \right\} \ \ \ \ \ \ \ \ \ \ \  [\ \uparrow \ \ {\mbox{see}} \ \ (\,A.10.16\,) \ ]\\[0.15in]
& \  & \hspace*{-0.5in}\ \left[\ = \  {\tilde C}_3\, (\,n) \,\cdot\,   {1\over { {\bf d}_{\,1\,, \ 2}^{n} }}   \,\cdot \left(\, {{\, \xi_{1_{\,|_{\,1}}}  - \ \xi_{2_{\,|_{\,1}}} \, }\over
{ \lambda_{\,2}}}\  \right)  \,\cdot\,\left\{ \ 1 \ + \  O \left(\ {\bar \lambda}_{\,\,\flat}^{\ 2\,   (\,1 \,-\ \nu\,)}\ \right)\ \right\} \ \right]\  \ \ \ \ \ \ \ \ [\ \leftarrow \ \ \ {\mbox{see \ \ (\,A.10.20\,)\, }} \,]\ .
\end{eqnarray*}

In (\,A.10.14\,)\,,
$$\ \ \ \   [\ \cdot\,\cdot\,\cdot\,\ ]_{\,(\,A.10.14\,)}
\ = \ \left(\ {\lambda_{\,1}\over { \,\lambda_{\,1}^2 \, +\  {\bar r}^{\,2}\, }}  \right)^{\!\!2}\,\cdot\,  {{1\,}\over {(\,\lambda_{\,1}^2 \, +\  {\bar r}^{\,2}\, )^{{n\over 2 } }}}  \ \ \ \ \ \  \ \ \ [\ \ {\bar r}^2 \ = \ {\bar y}_1^2 \ + \ \cdot \cdot \cdot \ + \ {\bar y}_n^2\ \  ] \ . \leqno (\,A.10.15\,)
$$
Note that the term
$$
  \lambda_{\,1}^2 \cdot \int_{B_{o}\, (\,{\rho_{\,\nu}})}  [\ \cdot\,\cdot\,\cdot\,\ ]_{\,(\,A.10.14\,)}\,\cdot\,   \left(\  {{ {\bar y}_1 }\over { \lambda_{\,1}  }}\  \right)^{\!\!2} \ d\,{\bar y}\ , \leqno (\,A.10.16\,)
$$
which appears in (\,A.10.14\,)\,,\,
involves the first bubble parameters $\,\lambda_{\,1}\,$ and $\,\xi_{\,1}\,$ only [\ independent on other bubbles $\, \ l \ = \ 2\,, \ \cdot \cdot \cdot\,, \ \flat\,$ \ ]\,.\bk
Summing up the last expression in (\,A.10.14\,) with respect to other bubbles (\,they are $\, l \ = \ 2\,, \ \cdot \cdot \cdot\,, \ \flat$\,)\,,\, we obtain the following
 \begin{eqnarray*}
(\,A.10.17\,)\ \ \   & \  &   + \  n\,(\,n\,-\,2)^{\,2}\,\cdot\, (\,n\,+\,2)\,\cdot\, \sum_{l \ = \ 2}^\flat \left[\    {1\over { {\bf d}_{\,1\,, \ l}^{n} }}   \,\cdot \left(\, {{\, \xi_{1_{\,|_{\,1}}}  - \ \xi_{l_{\,|_{\,1}}} \, }\over
{ \lambda_{\,2}}}\  \right) \ \right] \ * \\[0.2in]
& \ & \  \  \ \ \ \ \ \ \ \ \ \ \ \ \ \ \ \ \ \ \ \ \ \ \ \ \ \ \ \ \ * \,\left\{ \ \lambda_{\,1}^2 \cdot \int_{B_{o}\, (\,{\rho_{\,\nu}})}  [\ \cdot\,\cdot\,\cdot\,\ ]_{\,(\,A.10.14\,)}\,\cdot\,   \left(\  {{ {\bar y}_1 }\over { \lambda_{\,1}  }}\  \right)^{\!\!2} \ d\,{\bar y}\, \right\} \\[0.2in]
& \ & \hspace*{3in}[\ \uparrow\ \ {\mbox{independent \ \ on \  \ }} \ell \ ]\ .\ \ \ \ \ \ \ \ \ \ \ \ \
\end{eqnarray*}

\vspace*{-0.3in}

We continue with\\[0.1in]
(\,A.10.18\,)
$$
  \lambda_{\,1}^2 \cdot \int_{B_{o}\, (\,{\rho_{\,\nu}})}  [\ \cdot\,\cdot\,\cdot\,\ ]_{\,(\,A.10.14\,)}\,\cdot\,   \left(\  {{ {\bar y}_1 }\over { \lambda_{\,1}  }}\  \right)^{\!\!2} \ d\,{\bar y}\,  \ = \  \lambda_{\,1}^2 \cdot \left[\   \int_{\R^n} \ -  \int_{\R^n \ \setminus \ \emph{}B_{o}\, (\,{\rho_{\,\nu}})} \ \right] [\ \cdot\,\cdot\,\cdot\,\ ]_{\,(\,A.10.14\,)}\,\cdot\,   \left(\  {{ {\bar y}_1 }\over { \lambda_{\,1}  }}\  \right)^{\!\!2} \ d\,{\bar y}\ .
$$

Let us compute the integrations in (\,A.10.18\,) separately\,.\, \\[0.1in]
(\,A.10.19\,)
 \begin{eqnarray*}
 & \ &   \lambda_{\,1}^2 \cdot \int_{B_{o}\, (\,{\rho_{\,\nu}})}  [\ \cdot\,\cdot\,\cdot\,\ ]_{\,(\,A.10.14\,)}\,\cdot\,   \left(\  {{ {\bar y}_1 }\over { \lambda_{\,1}  }}\  \right)^{\!\!2} \ d\,{\bar y}\,\\[0.2in]
& = & \omega_n  \cdot \lambda_{\,1}^2 \cdot \int_0^\infty \left[\   \left({\lambda_{\,1}\over { \,\lambda_{\,1}^2 \, +\  {\bar r}^{\,2}\, }}  \right)^{\!\!2}\,\cdot\,  {{1\,}\over {(\,\lambda_{\,1}^2 \, +\  {\bar r}^{\,2}\, )^{{n\over 2 } }}} \ \right]\cdot  \left(\  {{ {\bar y}_1 }\over { \lambda_{\,1}  }}\  \right)^{\!\!2} \,\cdot\, {\bar r}^{\,n\,-\,1}\ d\, {\bar r} \\[0.2in]
& = & \omega_n  \cdot \lambda_{\,1}^2 \cdot   {1\over n}\,\cdot\,\int_0^\infty \left[\   \left({\lambda_{\,1}\over { \,\lambda_{\,1}^2 \, +\  {\bar r}^{\,2}\, }}  \right)^{\!\!2}\,\cdot\,  {{1\,}\over {(\,\lambda_{\,1}^2 \, +\  {\bar r}^{\,2}\, )^{{n\over 2 } }}} \ \right]\cdot  \left(\  {{ {\bar r} }\over { \lambda_{\,1}  }}\  \right)^{\!\!2} \,\cdot\, {\bar r}^{\,n\,-\,1}\ d\, {\bar r}\\[0.1in]
& \ & \hspace*{3in}\ \ \ \ \ \ \ \ \ \ \ \ \ \  [\ \uparrow \ \ \ {\mbox{via \ \ symmetry}}\ ]\\[0.1in]
& = & \omega_n  \cdot \lambda_{\,1}^2 \cdot  {1\over n}\,\cdot\,{1\over { \lambda_{\,1}^2 }}\,\cdot\,\int_0^{\pi\over 2} [\ \sin\,\theta\ ]^{\,n\,+\,1} \cos\,\theta \ d\,\theta \ \ \ \ \  \ \ \ \ \ ( \ {\bar r}^2 \ = \ y_1^2 \ + \ \cdot \cdot \cdot \ + \ y_n^2\ ) \\[0.2in]
&  = & \omega_n \,\cdot\,{1\over n} \,\cdot\,{1\over {n\,+\,2}}  \ \ \ \ \ \ \ \ \ \ \ \ \ \ \ \ \ \  \ \ \ \  \ \ \ \  \ \ \ \   [\ {\bar r} \ = \ \lambda_{\,1} \cdot \tan \,\theta\ \ \ \ \ \  {\mbox{cf. \ \ (\,A.10.18\,)}}\ ]\ .\\
\end{eqnarray*}
Here $\,\omega_n\,$ is the volume of the unit \,$(\,n\,-\,1\,)$\,-\,sphere in $\,\R^n\,.$\bk
As for the ``\,outside term\,"\, indicated in (\,A.10.18\,), we calculate as in {\bf \S\,A\,8\,.\,\,h}\,
and (\,A.9.15\,)\,:\\
(\,A.10.20\,)
   \begin{eqnarray*}
  & \ & \lambda_{\,1}^2\,
\cdot\,\int_{\R^n\, \setminus B_{o}\, (\,{\rho_{\,\nu}})} \left(\ {\lambda_{\,1}\over { \,\lambda_{\,1}^2 \, +\  {\bar r}^{\,2}\, }}  \right)^{\!\!2}\,\cdot\,  {{1\,}\over {(\,\lambda_{\,1}^2 \, +\  {\bar r}^{\,2}\, )^{{n\over 2 } }}}    \,\cdot\,    \left(\  {{ {\bar y}_1 }\over { \lambda_{\,1}  }}\  \right)^{\!\!2} \ d\,{\bar y} \\[0.2in]
  & \le & C\,\cdot\,\int_{\rho_{\,\nu}}^\infty \left(\ {\lambda_{\,1}\over { \,\lambda_{\,1}^2 \, +\  {\bar r}^{\,2}\, }}  \right)^{\!\!2}\,\cdot\,  {{1\,}\over {(\,\lambda_{\,1}^2 \, +\  {\bar r}^{\,2}\, )^{{n\over 2 } }}}    \,\cdot\,    {\bar r}^{\,2}\,\cdot\,{\bar r}^{\,n\,-\,1}\ d\,{\bar r}  \\[0.2in]
  & = & C\,\cdot\,{{ \lambda_{\,1}^{n\,+\,2} }\over { \lambda_{\,1}^2\,\cdot\,\lambda_{\,1}^n }} \ \,\cdot\,\int_{  {{ \rho_{\,\nu} }\over {\lambda_{\,1}}} }^\infty \left(\ {1\over { \,1 \, +\   \left(\ {{ \bar r }\over {\lambda_{\,1}}}  \right)^{\!2}\, }}  \right)^{\!\!2}\,\cdot\,\left(\ {1\over { \,1 \, +\   \left(\ {{ \bar r }\over {\lambda_{\,1}}}  \right)^{\!2}\, }}  \right)^{\!\!{n\over 2}}    \,\cdot\, \left(\ {{ \bar r }\over {\lambda_{\,1}}}  \right)^{\!n\,+\,1} \ d \left(\, {{ \bar r }\over {\,\lambda_{\,1}\,}}  \right) \\[0.2in]
  & \le & C_1\,\cdot\, \int_{  {{ \rho_{\,\nu} }\over {\lambda_{\,1}}} }^\infty {{R^{\,n\,+\,1} \ dR }\over {R^{\,n\,+\,4}  }} \ \le \ C_2\,\cdot\,{1\over {R^2}} \ \bigg\vert_{\,{  {{ \rho_{\,\nu} }\over {\lambda_{\,1}}} }}^{\,\infty} \ \ \ \ \ \ \ \ \ \ \ \ \ \  \ \ \ \ \ \ \  \ \ \ \ \ \ \  \ \ \ \ \ \ \  \left[\ {\bar R} \ = \ {{ \bar r }\over {\lambda_{\,1}}} \ \right] \\[0.2in]
    & = & O \left(\ {\bar \lambda}_{\,\,\flat}^{\ 2\,   (\,1 \,-\ \nu\,)}\ \right)\ .
  \end{eqnarray*}

  \newpage

It follows that the expression in (\,A.10.14\,) can be written as\\[0.1in]
(\,A.10.21\,)
   $$
   + \  \omega_n \cdot (\,n\,-\,2)^{\,2}\,\cdot\,  \sum_{l \ = \ 2}^\flat \left[\    {1\over { {\bf d}_{\,1\,, \ l}^{n} }}   \,\cdot \left(\, {{\, \xi_{1_{\,|_{\,1}}} \!\! - \ \xi_{l_{\,|_{\,1}}} \, }\over
{ \lambda_{\,2}}}\  \right) \ \right] \,\times \left[\ 1\ + \ O \left(\ {\bar \lambda}_{\ \flat}^{\ 2\,   (\,1 \,-\ \nu\,)}\ \right)\ \right] \ .
$$

\vspace*{0.3in}

{\bf \S\,A\,10\,.b\,.} \ \ {\bf The second order term} in (\,A.10.9\,)\,. \ \ As in {\bf \S\,A\,8\,.\,\,e}\,,\, \\[0.1in]
(\,A.10.22\,)
\begin{eqnarray*}
{\bf T }^2 & =  &   \left[\ {1\over { {\bf d}_{\,1\,, \ 2}^{2} }}\,\cdot\,\left(\ {{\lambda_{\,2}}\over {\lambda_{\,1}}}  \right) \ \right]^2 \ + \ \left[\ {1\over { {\bf d}_{\,1\,, \ 2}^{2} }}\,\cdot\,{{\Vert\ \,y\ - \ \xi_{\,1}\,\Vert^2}\over
{\lambda_{\,1}\,\cdot\,\,\lambda_{\,2}}}    \ \right]^2  \ \ + \\[0.2in]
& \ & \ \ \ + \   \left[\  {1\over { {\bf d}_{\,1\,, \ 2}^{2} }}\,\cdot\,{{\,2\,(\,y\ - \ \xi_{\,1})\cdot(\,\xi_{\,1}\ - \ \xi_{\,2})\, }\over
{\lambda_{\,1}\,\cdot\,\,\lambda_{\,2}}} \ \right]^2  \\[0.2in]
& \ & \ \ \  \ \ \ + \ 2\,\cdot\, \left[\ {1\over { {\bf d}_{\,1\,, \ 2}^{2} }}\,\cdot\,\left(\ {{\lambda_{\,2}}\over {\lambda_{\,1}}}  \right) \ \right]\,\cdot\,\left[\ {1\over { {\bf d}_{\,1\,, \ 2}^{2} }}\,\cdot\,{{\Vert\ \,y\ - \ \xi_{\,1}\,\Vert^2}\over
{\lambda_{\,1}\,\cdot\,\,\lambda_{\,2}}}    \ \right] \\[0.2in]
& \ & \ \ \ \ \ \  \ \ \ + \ 2\,\cdot\, \left[\ {1\over { {\bf d}_{\,1\,, \ 2}^{2} }}\,\cdot\,\left(\ {{\lambda_{\,2}}\over {\lambda_{\,1}}}  \right) \ \right]\,\cdot\,\left[\  {1\over { {\bf d}_{\,1\,, \ 2}^{2} }}\,\cdot\,{{\,2\,(\,y\ - \ \xi_{\,1})\cdot(\,\xi_{\,1}\ - \ \xi_{\,2})\, }\over
{\lambda_{\,1}\,\cdot\,\,\lambda_{\,2}}} \ \right] \\[0.2in]
& \ & \ \ \ \ \ \ \ \ \ \   \ \ \ + \ 2\,\cdot\,\left[\ {1\over { {\bf d}_{\,1\,, \ 2}^{2} }}\,\cdot\,{{\Vert\ \,y\ - \ \xi_{\,1}\,\Vert^2}\over
{\lambda_{\,1}\,\cdot\,\,\lambda_{\,2}}}    \ \right]\,\cdot\,\left[\  {1\over { {\bf d}_{\,1\,, \ 2}^{2} }}\,\cdot\,{{\,2\,(\,y\ - \ \xi_{\,1})\cdot(\,\xi_{\,1}\ - \ \xi_{\,2})\, }\over
{\lambda_{\,1}\,\cdot\,\,\lambda_{\,2}}} \ \right] \ .
 \end{eqnarray*}
Via (\,A.10.11\,) [ \,symmetric cancellation\ ]\,,\, we have\\[0.1in]
(\,A.10.23\,)
 \begin{eqnarray*}
 & \ &  \lambda_{\,1}^2\,\cdot\,  {1\over { {\bf d}_{\,1\,, \ 2}^{n\,-\,2} }} \, \int_{B_{\,\xi_{\,1}} (\,{\rho_{\,\nu}}\,)}[\ \bullet\, \bullet\, \bullet \ ]_{\,(\,A.10.9\,)}  \,\cdot\,\left(\, {{y_1 \,-\,\xi_{1_{\,|_{\,1}}}}\over {\lambda_{\,1}}}  \right)\times{\bf T}^2 \\[0.2in]
 &  =  &   \lambda_{\,1}\,\cdot\,  {1\over { {\bf d}_{\,1\,, \ 2}^{n\,-\,2} }} \int_{B_{o} (\,{\rho_{\,\nu}})} \left(\ {\lambda_{\,1}\over { \,\lambda_{\,1}^2 \, +\  {\bar r}^{\,2}\, }}  \right)^{\!\!2}\,\cdot\,  {{1\,}\over {(\,\lambda_{\,1}^2 \, +\  {\bar r}^{\,2}\, )^{{n\over 2 } }}}  \,\cdot\,{\tilde{\bf T}}^2  \,\cdot\,[\ {\bar y}_1 \ ]   \\[0.2in]
 & = & 0 \ + \  0 \ + \ 0\ + \ 0 \ + \  {2\over { {\bf d}_{\,1\,, \ 2}^{n\,-\,2\,+\,2} }} \,\cdot\,\left\{ \ \int_{B_{o} (\,{\rho_{\,\nu}})}[\ \cdot\,\cdot\,\cdot\, \ ]_{\,(\,A.10.10\,)}  \,\cdot\,    {\bar y}_1^2\ \right\} \,\cdot\,\left(\, {{ \,\xi_{1_{\,|_1}}  - \ \xi_{2_{\,|_1}} \, }\over {\lambda_{\,2}  }} \,\right)   \,*\\[0.2in]
& \ & \hspace*{3in}   *\, \left[\ {1\over { {\bf d}_{\,1\,, \ 2}^{2} }}\,\cdot\,\left(\ {{\lambda_{\,2}}\over {\lambda_{\,1}}}  \right)  \ + \  {1\over { {\bf d}_{\,1\,, \ 2}^{2} }}\,\cdot\,{{ {\bar r}^2 }\over
{\lambda_{\,1}\,\cdot\,\,\lambda_{\,2}}}   \ \right] \\[0.2in]
& = &  {2\over { {\bf d}_{\,1\,, \ 2}^{n } }} \,\cdot\,\left\{ \ \int_{B_{o} (\,{\rho_{\,\nu}})} \left(\ {\lambda_{\,1}\over { \,\lambda_{\,1}^2 \, +\  {\bar r}^{\,2}\, }}  \right)^{\!\!2}\,\cdot\,  {{1\,}\over {(\,\lambda_{\,1}^2 \, +\  {\bar r}^{\,2}\, )^{{n\over 2 } }}}    \,\cdot\,    {\bar y}_1^2\ \right\} \,\cdot\,\left(\, {{ \,\xi_{1_{\,|_1}}  - \ \xi_{2_{\,|_1}} \, }\over {\lambda_{\,2}  }} \,\right) \,*\\[0.2in]
& \ & \hspace*{3in} \,*\,\left[\ O\,\bigg(\ {\bar \lambda}_{\,\,\flat}^{\,2\,\gamma}\ \bigg) \ + \  O\,\bigg(\ {\bar \lambda}_{\,\,\flat}^{\,2\,(\,\gamma \ + \ \nu\,)\ - \ 1}\ \bigg) \ \right] \  .
\end{eqnarray*}
In the above, $\,\tilde{\bf T}\,$ is related to $\,{\bf T}\,$ via the change of variables $${\bar{y}} \ = \ y \ - \ \xi_{\,1}\,.\,$$
Calculate as in (\,A.10.20\,)\,,\, and  note that
 \begin{eqnarray*}
(\,A.10.24\,) \ \ \ \ \ \ \ \ \ \ \ \ & \ & \int_{B_{o} (\,{\rho_{\,\nu}})} \left(\ {\lambda_{\,1}\over { \,\lambda_{\,1}^2 \, +\  {\bar r}^{\,2}\, }}  \right)^{\!\!2}\,\cdot\,  {{1\,}\over {(\,\lambda_{\,1}^2 \, +\  {\bar r}^{\,2}\, )^{{n\over 2 } }}}    \,\cdot\,    {\bar y}_1^2\ \ \ \ \ \ \ \ \ \ \ \ \ \ \ \ \ \ \ \ \ \ \ \ \ \ \ \ \ \ \ \ \ \ \ \ \\[0.2in]
 & = &  C_1 \cdot \int_0^{{\rho_\nu}\over { \lambda_{\,1}}} {{R^2 \cdot R^{\,n \,- \ 1}\ dR}\over { (\,1\ + \ R^2\,)^{ {n\over 2} \ + \ 2} }}\\[0.2in]
& \le &  C_2 \cdot \int_0^\infty {{R^{\,n \,+ \ 1}\ dR}\over { (\,1\ + \ R^2\,)^{ {n\over 2} \, + \,2} }} \ < \ \infty\ .
\end{eqnarray*}
Summing up and making use of the inequality
$$
\bigg\vert \  {{ \,\xi_{1_{\,|_1}}  - \ \xi_{2_{\,|_1}} \, }\over {\lambda_{\,2}  }} \ \bigg\vert \ \le \   C \cdot {{\ \Vert\  \,\xi_{\,1}\ - \ \xi_{\,2}\,\Vert \ }\over
{\sqrt{\,\lambda_{\,1}\,\cdot\,\,\lambda_{\,2}\,}}} \ = \ C \cdot {\bf d}_{\,1\,, \ 2}\ ,
$$
we obtain
 \begin{eqnarray*}
(\,A.10.25\,) \ \ \ \ \  \ \ \ \ \ & \ & \sum_{l \ = \ 2}^\flat   \lambda_{\,1}^2\,\cdot\,  {1\over { {\bf d}_{\,1\,, \ l}^{n\,-\,2} }} \, \int_{B_{\,\xi_{\,1}} (\,{\rho_{\,\nu}}\,)}[\ \bullet\, \bullet\, \bullet \ ]_{\,(\,A.10.9\,)}  \,\cdot\,\left(\, {{y_1 \,-\,\xi_{1_{\,|_{\,1}}}}\over {\lambda_{\,1}}}  \right)\cdot {\bf T}^2 \ \ \ \ \  \ \ \ \ \  \ \ \ \ \  \\[0.2in]
&  \le &  C \cdot  \sum_{l \ = \ 2}^\flat {1\over { {\bf d}_{\,1\,, \ l}^{n \ + \ 1} }} \ = \ O \left(\  {\bar\lambda}_{\ \flat}^{\,(\,n \ + \ 1)\, \cdot\, \gamma}\ \right)\ .
\end{eqnarray*}


\newpage

{\bf \S\,A\,10\,.c\,.} \ \  {\bf The third order terms} $\,O\,(\,|\, {\bf T} \,|^3\,)\,.$\, We look at the term\\[0.1in]
(\,A.10.26\,)
 \begin{eqnarray*}
  & \  &  \bigg\vert \  \lambda_{\,1}^2\,\cdot\,  {1\over { {\bf d}_{\,1\,, \ 2}^{n\,-\,2} }} \int_{B_{\,\xi_{\,1}} (\,{\rho_{\,\nu}}\,)} \left(\ {\lambda_{\,1}\over { \,\lambda_{\,1}^2 + \Vert\  y \,-\, \xi_{\,1}\,\Vert^{\,2}\, }}  \right)^{\!\!2}\,\cdot\, \left( \  {{1\,}\over {\,\lambda_{\,1}^2 + \Vert\  y \,-\, \xi_{\,1}\,\Vert^{\,2}\, }}\emph{}\ \right)^{\! {n\over 2} }   \,\cdot\,\left(\, {{y_1 \,-\,\xi_{1_{\,|_{\,1}}} }\over {\lambda_{\,1}}} \,\right)\,  \bigg\vert \ *\\[0.2in]
  & \ & \ \ \ \ \  \ \ \ \ \ \ \ \ \ \  \ \ \ \ \  \ \ \ \ \  \ \ \ \ \ \ \ \ \ \  \ \ \ \ \     \ \ \ \ \  \ *\  \bigg\vert \,  {2\over { {\bf d}_{\,1\,, \ 2}^{2} }}\,\cdot\,\left(\ {{ y_1\, - \ \xi_{1_{\,|_{\,1}}} }\over
{\lambda_{\,1}  }}\  \right)\,\cdot\,\left(\ {{\, \xi_{1_{\,|_{\,1}}}  - \ \xi_{2_{\,|_{\,1}}} \, }\over
{ \lambda_{\,2}}}\  \right)   \ \bigg\vert^{\,3} \ dy\\[0.2in]
& \le & C_1\,\,\cdot\,\lambda_{\,1}^2\,\cdot\,  {1\over {\ {\bf d}_{\,1\,, \ 2}^{\,n\,-\,2 \,+\,2\,+\,2\,+\,2} \ }} \,*\, {{\,|\, \xi_{1_{\,|_{\,1}}}  - \ \xi_{2_{\,|_{\,1}}} \,|^{\,3} \, }\over
{ \lambda_{\,2}^{\,3}}} \ *  \\[0.2in]
& \ & \ \ \ \ \ *\int_{B_{\,\xi_{\,1}} (\,{\rho_{\,\nu}}\,)} \left(\ {\lambda_{\,1}\over { \,\lambda_{\,1}^2 \ + \ \Vert\  y \,-\, \xi_{\,1}\,\Vert^{\,2}\, }}  \right)^{\!\!2}\,\cdot\,  {{1\,}\over {(\,\lambda_{\,1}^2 \ + \ \Vert\  y \,-\  \xi_{\,1}\,\Vert^{\,2}\,)^{{n\over 2 } }}}\,\cdot\,\left(\ {{ y_1\ - \ \xi_{1_{\,|_{\,1}}} }\over
{\lambda_{\,1}  }}\  \right)^{\!\!4} \ dy \\[0.2in]
& \le & C_2\,\,\cdot\,  {{ {\bf d}_{\,1\,, \ 2}^2}\over {\ {\bf d}_{\,1\,, \ 2}^{n\,+\,4} \ }}\,\cdot\,\bigg\vert \  {{\,  \xi_{1_{\,|_{\,1}}} - \ \xi_{2_{\,|_{\,1}}}  \, }\over
{ \lambda_{\,2}}} \ \bigg\vert\,\cdot\,\int_{B_{\,\xi_{\,1}} (\,{\rho_{\,\nu}}\,)}\! \left(\ {1\over { \,1 + \,{{ \Vert\  y \ -\ \xi_{\,1}\,\Vert^{\,2} }\over {\lambda_{\,1}^2}} }}  \right)^{\!\!2} \!\cdot\!   \left(\ {1\over { \,1 + \,{{ \Vert\  y \ -\  \xi_{\,1}\,\Vert^{\,2} }\over {\lambda_{\,1}^2}} }}  \right)^{\!\!{n\over 2}} \ * \\[0.2in]
  & \ & \ \ \ \ \  \ \ \ \ \ \ \ \ \ \  \ \ \ \ \  \ \ \ \ \ \ \ \ \ \ \ \ \ \ \ \ \ \ \ \ \  \ \ \ \ \ \ \ \ \ \  \ \ \ \ \     \ \ \ \ \  \ \ \ \ \  \ \ \ \ \  \ \ \ \  *\, \left(\ {{ y_1  - \ \xi_{1_{\,|_{\,1}}} }\over
{\lambda_{\,1}  }}\  \right)^{\!\!4} \ {{  \lambda_{\,1}^2 }\over {  \lambda_{\,1}^{n\,+\,2} }} \  \,dy \\[0.2in]
& \le & C_3\,\,\cdot\,  {{1}\over {\ {\bf d}_{\,1\,, \ 2}^{n \ + \ 1} \ }}\, \cdot\,\int_0^{\,{{\rho_\nu}\over {\lambda_{\,1}}}  }  \left(\ {1\over { \,1 + \,R^{\,\,2}  }}  \right)^{\!\!2} \!\cdot\!   \left(\ {1\over { \,1 + \,R^{\,\,2}  }}  \right)^{\!\!{n\over 2}}\,\cdot\,R^4\,\cdot\,R^{\,n\,-\,1} \ d\,R\\[0.2in]
& \le & C_4\,\,\cdot\,  {{1}\over {\ {\bf d}_{\,1\,, \ 2}^{n \ + \ 1} \ }}\, \cdot\,\int_0^{\,{{\rho_\nu}\over {\lambda_{\,1}}}  }  \  {{d\,R}\over { R}}\\[0.2in]
& \le & C_5\,\,\cdot\,  {{1}\over {\ {\bf d}_{\,1\,, \ 2}^{n \ + \ 1} \ }}\, \cdot\,\ln\,\left( {{\rho_\nu}\over {\lambda_{\,1}}}\ \right)   \\[0.2in]
& \le & C_6\,\,\cdot\,  {{1}\over {\ {\bf d}_{\,1\,, \ 2}^{n \  + \ 1 } \ }} \,\cdot\,(\gamma \ - \ \nu\,)\,\cdot\, |\,\ln\,{\bar\lambda}_{\ \flat}\,|  \\[0.2in]
& \le & C_7\,\,\cdot\,  {{1}\over {\ {\bf d}_{\,1\,, \ 2}^{n \  + \ 1 \ - \ o_{\,{\bar{\lambda}}_{\,\flat}   }\,(\,1\,) } \ }} \,\cdot\,(\gamma \ - \ \nu\,)  \ \ \ \ \
\ \ [\ \gamma\  >  \ \nu\,; \ \ \    o_{\,{\bar{\lambda}}_{\,\flat}   }\,(\,1\,) \ \to \ 0^{\,+} \ \ \   {\mbox{as}} \ \ \ {\bar\lambda}_{\ \flat} \ \to \ 0^{\,+} \ ]\ .
 \end{eqnarray*}
\newpage

Using the symmetric cancellation in (\,A.10.1\,) [\,or calculate in a similar fashion as in (\,A.10.23)\,]\,,\, we can estimate other terms. Summing up\,,\, we reach the following the conclusion

\begin{eqnarray*}
(\,A.10.27\,) \ \ \ \ \ \ \ & \ & \sum_{l \ = \ 2}^\flat   \lambda_{\,1}^2\,\cdot\,  {1\over {\  {\bf d}_{\,1\,, \ l}^{n\,-\,2} \ }} \, \int_{B_{\,\xi_{\,1}} (\,{\rho_{\,\nu}}\,)}[\ \bullet\, \bullet\, \bullet \ ]_{\,(\,A.10.9\,)}  \,\cdot\,\left(\, {{y_1 \,-\,\xi_{1_{\,|_{\,1}}}}\over {\lambda_{\,1}}}  \right)\times{\bf T}^3\ \ \ \ \ \ \ \ \ \ \ \  \ \ \ \ \ \ \ \ \ \ \ \  \ \ \ \ \ \ \ \ \ \ \ \ \ \ \ \ \ \ \ \\[0.2in]
& \le &  C \cdot  \sum_{l \ = \ 2}^\flat\  {1\over {\  {\bf d}_{\,1\,, \ l}^{n \ + \ 1 \ - \  o_{\,{\bar{\lambda}}_{\ \flat}   }\,(\,1\,) } \ }}\\[0.2in]
&  = &  O \left(\  {\bar\lambda}_{\ \flat}^{\,(\,n \ + \ 1\,) \,\cdot\, \gamma \ - \  \, o_{\,{\bar{\lambda}}_{\ \flat}   }\,(\,1\,)}\ \right)\ .\\
\end{eqnarray*}
We arrive at\\[0.1in]
(\,A.10.28\,)
 \begin{eqnarray*} {\bf (\,C\,)}_{\,(\,A.10.1\,)}
 & = &
\left\{  -\int_{\R^n} (\,\Delta\,W_{\,\,\flat})\,\cdot\,\left(\,\lambda_{\,1}\cdot \btd_{\xi_{\,1}}  W_{\,\,\flat}  \,\right)    \, -     \int_{\R^n}\, \,n\,(\,n\,-\,2)\cdot\,(\,W_{\,\,\flat} \, )^{{\,n\,+\,2\,}\over {n\,-\,2}}\,\cdot\, \left(\,\lambda_{\,1}\cdot \btd_{\xi_{\,1}}  W_{\,\,\flat}  \,\right)  \, \right\}  \\[0.3in]
 &   \  &  \hspace*{-0.75in} = \  {\hat C}_{\,2\,, \ 2} \,\cdot\,\sum_{l\ = \ 2}^\flat \left\{  \, \left(\ \ {{ \sqrt{\lambda_{\,1}\,\cdot\,\lambda_{\,\,l}\,} }\over {\Vert\ \xi_{\,1} \ - \ \xi_{\, l}\ \Vert   }}\  \right)^{\!\!n}\,\cdot\,\left(\ {{\,\xi_{\,1}\ - \ \xi_{\,l}\,}\over {\lambda_{\ \!l}}}\  \right)  \ \right\}  \,\times \left[\ 1 \  +  \  O\,\left( \,{\bar\lambda}_{\ \flat}^{\,2\,(\,1 \ - \ \nu\,) }\, \right) \ \right]\   \\[0.3in]
 & \ & \ \ \ \ \ \  \  \ \ \ \ \ \  \  + \ \ O \left(\  {\bar\lambda}_{\ \flat}^{(\,n \ + \ 1\,) \,\cdot\, \gamma \ - \   o_{\,{\bar{\lambda}}_{\ \flat}  }(\,1\,)}\ \right)\ \ +  \ \,{\bf{E}}_{\,(\,A.10.28\,)} \ .\\
  \end{eqnarray*}
  Here
  $$ {\hat C}_{\,2\,, \ 2} \ = \    \omega_n\,\cdot\,(\,n\,-\,2)^{\,2} \ .
$$
See (\,A.8.45\,) for  the expression of  $\,\,{\bf{E}}_{\,(\,A.10.28\,)}\,$.\bk
In (\,A.10.28\,)\,,\, the bubble parameters
 $\,\left( \ \lambda_{\,1}\,, \ \cdot \cdot \cdot\,, \ \lambda_{\,\flat}\,, \ \xi_{\,1}\,, \ \cdot \cdot \cdot\,, \ \xi_{\,\flat}   \ \right)\,$ satisfy the conditions in Theorem 1.33 of the main text\,.\,


\newpage

{\large{\bf{\S\,A\,11\,. \ \ Estimate of}}} $\,(\,{\bf D}\,)_{\,(\,A.10.1\,)}\ .$\\[0.2in]
We come to consider the term in $\,(\,{\bf D}\,)_{\,(\,A.10.1\,)}\,,$\, expressed in vector form\,:\\[0.1in]
 (\,A.11.1\,)
 \begin{eqnarray*}
 {\bf{[\,{\vec{\,D}}\,]}}_{\,(\,A.11.1\,)} & = &  \left\{ \ \int_{\R^n}\, [\,(\,{\tilde c}_n\!\cdot K\,) \ - \ n\,(\,n\,-\,2)\,]  \cdot\,(W_{\,\,\flat} )^{{\,n\,+\,2\,}\over {n\,-\,2}}\,\cdot\,\left[\ (\,\lambda_{\,1}\cdot \btd_{\xi_{\,1}}\,) \, W_{\,\,\flat}\ \right]  \ \right\}\\[0.2in]
 & = &  \left\{ \ \int_{\R^n}\, [\,(\,{\tilde c}_n\!\cdot K\,) \ - \ n\,(\,n\,-\,2)\,]  \cdot\,(W_{\,\,\flat} )^{{\,n\,+\,2\,}\over {n\,-\,2}}\,\cdot\,\left[\ (\,\lambda_{\,1}\cdot \btd_{\xi_{\,1}}\,) \, V_1\ \right]  \ \right\}\\[0.2in]
 & = &  \left\{ \ \int_{B_o\,(\,\xi_{\,1}\,)}\ \ [\,(\,{\tilde c}_n\!\cdot K\,) \ - \ n\,(\,n\,-\,2)\,]  \cdot\,V_1^{{\,n\,+\,2\,}\over {n\,-\,2}}\,\cdot\,\left[\ (\,\lambda_{\,1}\cdot \btd_{\xi_{\,1}}\,) \, V_1\ \right]  \ \right\}\ * \\[0.2in]
  & \ & \hspace*{2in} \ \ \ \ \ \ \ \ *\ \left[\ 1 \ + \ O\,\left( \ {\bar\lambda}_{\ \flat}^{\ (\,n\ -\,2\,)\,\cdot\, [\ (\,\gamma \ + \ \nu \,)\ - \  1\,] }\ \right) \  \right] \\[0.2in]
  & \ & \ \ \ \ \ \ \ \ \ + \ \,O \left(\  {\bar\lambda}_{\ \flat}^{\ n\,(\,1\,-\ \nu)} \,\cdot\, {\bar\lambda}_{\ \flat}^{  \   - \ 2\,\varepsilon \,[\ (\,\gamma \ + \ \nu\,) \ - \ 1\ ] } \, \right)\\[0.1in]
  & \ & \ \ \ \ \ \ \ \ \ \  \ \  [\ \uparrow \ \ {\mbox{ estimate \ \  as \ \ in}} \ \   {\bf  \S\,A\,8\,.\,g\,}, \ \ {\mbox{ see \ \ also}} \ \  (\,A.9.3\,)\,.\,  \ ]\ .
\end{eqnarray*}
As the integral
 $$
 \int_{B_{\,\xi_{\,1}} (\,\rho_\nu)}\ \ [\,(\,{\tilde c}_n\!\cdot K\,) \ - \ n\,(\,n\,-\,2)\,]  \cdot\,V_1^{{\,n\,+\,2\,}\over {n\,-\,2}}\,\cdot\,\left[\ (\,\lambda_{\,1}\cdot \btd_{\xi_{\,1}}\,) \, V_1\ \right]   \leqno (\,A.11.2\,)
 $$
is invariant under translation and rotation, we follow the arrangement in {\bf \S\,A\,9\,.a}\,,\, and present it in the  setting
$$
\xi_{\,1} \ = \ (\,0\,,\ \cdot \, \cdot \, \cdot\,, \ \xi_{1_{\,|_{\,\,n}}}
\ )\,,  \ \ \ \ \ {\mbox{with}} \ \ \ \xi_{1_{\,|_{\,\,n}}} \ > \ 0\,.\leqno (\,A.11.3\,)
$$
Afterward, we bring the expression back to the original coordinates.
With  (\,A.11.3\,)\,,\, we have\\[0.1in]
(\,A.11.4\,)
 \begin{eqnarray*}
& \ &  \left(\,\lambda_{\,1}\,\cdot\,{\partial\over {\partial\, \xi_{1_{\,|_j}}}}\  \right)\,\left[\ \left({\lambda_{\,1}\over {\lambda^2_1\ + \ \Vert\  y \ -\  \xi_{\,1}\,\Vert^{\,2}}}\  \right)^{{n \,-\, 2}\over 2} \right]\\[0.2in] & = & -\,\left( \, {{n - 2}\over 2} \, \right)\,\cdot\,\lambda^{{n  }\over 2}_1\,\cdot\,{{ 2\,(\,\xi_{1_{\,|_j}}\, -\, y_j\,)}\over {(\,\lambda^2_1 + \Vert\  y \,-\, \xi_{\,1}\,\Vert^{\,2}\,)^{{n\over 2 } }}}\ = \  -\,\left( \, {{n - 2}\over 2} \, \right)\,\cdot\,\lambda^{{n  }\over 2}_1\,\cdot\,{{ 2\,(\,0\, -\, y_j\,)}\over {(\,\lambda^2_1 + \Vert\  y \,-\, \xi_{\,1}\,\Vert^{\,2}\,)^{{n\over 2 } }}}
\end{eqnarray*}
for $\,j \ = \ 1\,,  \ 2\,, \ \cdot \cdot \cdot\,, \ n  \ - \ 1\,.$\, As for $\,j\,=\,n\,$, we get\\[0.1in]
(\,A.11.5\,)
 \begin{eqnarray*}
& \ &  \left(\,\lambda_{\,1}\,\cdot\,{\partial\over {\partial\, \xi_{1_{\,|_n}}}}\  \right)\,\left[\left({\lambda_{\,1}\over {\lambda^2_1\ + \ \Vert\  y \ -\  \xi_{\,1}\,\Vert^{\,2}}}\  \right)^{{n \,-\, 2}\over 2} \right] \ = \ -\,\left( \, {{n - 2}\over 2} \, \right)\,\cdot\,\lambda^{{n  }\over 2}_1\,\cdot\,{{ 2\,(\,\xi_{1_{\,|_n}}\, -\, y_n\,)}\over {(\,\lambda^2_1 + \Vert\  y \,-\, \xi_{\,1}\,\Vert^{\,2}\,)^{{n\over 2 } }}} \ .
\end{eqnarray*}
Note that, as in (\,A.10.2\,)\,,\,
$$
\bigg\vert \ \left(\,\lambda_{\,1}\,\cdot\,{\partial\over {\partial\, \xi_{1_{\,|_j}}}}\  \right)\,V_1\,(\,y\,) \ \bigg\vert \ \le \  C \, V_1\,(\,y\,) \ \ \  \ \mbox{for} \ \ y \ \in \ \R^n \ \ \ {\mbox{and}} \ \ \ j \ = \ 1\,, \ 2\,, \ \cdot \cdot \cdot\,, \ n\ . \leqno (\,A.11.6\,)
$$
For convenience, we reproduce (\,4.44\,) and (\,4.45\,) of the main text   here.  For $\,y\,\in\,B_o\,(\,\rho_\nu\,)\,$ with $\,y_n \ >  \  \xi_{1_{\,|_{\,n}}}\,$ and $\,\ell \ \ge \ 3\,$ [\ cf. (\,4.41\,) of the original text for  $\,\ell \ = \ 2\,$]\,,\, we have\\[0.1in]
(\,A.11.7\,)
\begin{eqnarray*}
 \Vert \,y  \ - \ {\bf p}_{\,y}\,\Vert^\ell
  & = &  (\,y_n \ - \ \xi_{1_{\,|_{\,n}}}\,)^{\,\ell} \ \ \ \ \left(\ \leftarrow \ \ = \ |\, y_n \ - \ \xi_{1_{\,|_{\,n}}}\,|^{\,\,\ell} \ \right) \\[0.2in]
  & \ &  + \ \ell\,\cdot\,(\,y_n\,-\,\xi_{1_{\,|_{\,n}}}\,)^{\,\ell\,-\ 1}\cdot \xi_{1_{{\,n}}}   \ + \ \\[0.2in]
   &  & \ \ \ \ \  + \  {{ \ell\,\cdot\,(\,\ell\,-\,1\,) }\over {2}}\,\cdot\,(\,y_n \ - \ \xi_{1_{\,|_{\,n}}}\,)^{\,\ell\,-\ 2}\,\cdot\, \xi_{1_{\,|_{\,n}}}^2 \ \leftarrow \  \cdot \cdot \cdot \cdot \ (\,*\,)_{\,(\,A.11.7\,)} \\[0.2in]
   & \ &\hspace*{-1in} + \ \left[ \  O\,(\ \Vert \, y \ - \ \xi_1\,\Vert^{\,\ell\,-\ 3}\ ) \,\cdot\,  O\, (\,\xi_{1_{\,|_{\,n}}}^3 \,) \ + \   O\,(\, \Vert \, y \ - \ \xi_1\,\Vert^{\,\ell\,-\ 4}\ ) \,\cdot\,  O\, (\,\xi_{1_{\,|_{\,n}}}^4 \,) \ + \ \cdot \cdot \cdot \ +  \ \xi_{1_{\,|_{\,n}}}^\ell \ \right] \  + \\[0.2in]
    & \ & \ \ \ \ \ \ + \ \, O\,(\,|\,y_n\,|^{\,\ell\,-\,1} \,)\,\cdot\, \left[\  O\, (\,\xi_{1_{\,|_{\,n}}}^2 \,) \ + \ O\,(\, \Vert \, y \ - \ \xi_1\,\Vert^2\,)\ \right]\ .\\
\end{eqnarray*}
Likewise for $\,y\,\in\,B_o\,(\,\rho_\nu\,)\,$ with $\,y_n \ <  \  0\,$ and $\,\ell \ \ge \ 3\,$ [\ cf. (\,4.41\,) of the original text for  $\,\ell \ = \ 2\,$]\,,\, we have\\[0.1in]
(\,A.11.8\,)
\begin{eqnarray*}
 \Vert \,y  \ - \ {\bf p}_{\,y}\,\Vert^\ell & = &  (\,-\,y_n \ + \ \xi_{1_{\,|_{\,n}}}\,)^{\,\ell}  \ \ \ \ \left(\ \leftarrow \ \ = \ |\ y_n \ - \ \xi_{1_{\,|_{\,n}}}\,|^{\,\ell} \ \right) \\[0.2in]
  & \ &  + \ \ell\,\cdot\,(\,-\,y_n \ + \ \xi_{1_{\,|_{\,n}}}\,)^{\ell\,-\,1}\cdot (\,-\,\xi_{1_{{\,n}}}  \,)  \ + \ \\[0.2in]
  & \ & \ \ \ \ \  + \  {{ \ell\,\cdot\,(\,\ell\,-\,1\,) }\over {2}}\,\cdot\,(\,-\,y_n \ + \ \xi_{1_{\,|_{\,n}}}\,)^{\ell\,-\ 2}\,\cdot\, \xi_{1_{\,|_{\,n}}}^2   \\[0.2in]
    & \ &\hspace*{-1in} + \ \left[ \  O\,(\ \Vert \, y \ - \ \xi_1\,\Vert^{\,\ell\,-\ 3}\ ) \,\cdot\,  O\, (\,\xi_{1_{\,|_{\,n}}}^3 \,) \ + \   O\,(\, \Vert \, y \ - \ \xi_1\,\Vert^{\,\ell\,-\ 4}\ ) \,\cdot\,  O\, (\,\xi_{1_{\,|_{\,n}}}^4 \,) \ + \ \cdot \cdot \cdot \ +  \ \xi_{1_{\,|_{\,n}}}^\ell \ \right] \  + \\[0.2in]
    & \ & \ \ \ \ \ \ + \ \, O\,(\,|\,y_n\,|^{\,\ell\,-\,1} \,)\,\cdot\, \left[\  O\, (\,\xi_{1_{\,|_{\,n}}}^2 \,) \ + \ O\,(\, \Vert \, y \ - \ \xi_1\,\Vert^2\,)\ \right]\ .
\end{eqnarray*}
Due to the presence of $\,(\,y_n \ - \ \xi_{1_{\,|_n}}\,)$\, in (\,A.11.5\,)\,,\, we find that the argument here contrasts the one in {\bf \S\,A\,9}\,,\, in which the key term is the first term in (\,A.11.7\,) [\,or (\,A.11.8\,)\ ]\,,\, and to a lesser degree, the third term\,; where here the key contribution is the second term [\,indicated by  $\,(\,*\,)_{\,(\,A.11.7\,)}\,$ in \,(\,A.11.7\,)\ ]\,.\,

\vspace*{0.3in}

{\bf{\S\,A\,11\,.a}} \ \ {\bf Derivative in the $\,\xi_{1_{\,|_n}}$-\,direction\,.\,}  \ \ Let
$$
{\widetilde{\widetilde{[\ \cdot \cdot \cdot \ ]}}}_{\,(\,A.11.9\,)} \ = \ \left( {{\lambda_{\,1}}\over {\  \lambda_{\,1}^2 \ + \ \Vert\,y \ - \ \xi_{\,1}\,\Vert^{\,2}\ }}\  \right)^{\!\!n \ +\ 1} \ \cdot (\,y_n \ - \  \xi_{1_{\,|_n}} \ )\ . \leqno (\,A.11.9\,)
$$
Similar to (\,A.9.11\,)\,,\, We have\\[0.1in]
 (\,A.11.10\,)
\begin{eqnarray*}
&  \ &     \int_{B_{\,\xi_{\,1}} (\,{\rho_{\,\nu}}\,)}
\,[\,n\,(\,n\,-\,2\,) \,-\, (\,{ c}_n\!\cdot K\,)\,(\,y\,)  \ ] \,\cdot\,
{\widetilde{\widetilde{[\ \cdot \cdot \cdot \ ]}}}_{\,(\,A.11.9\,)}  \  dy \\[0.2in]
& = &     \int_{B_{\,\xi_{\,1}} (\,{\rho_{\,\nu}}\,)}
\,[\ C\,( \, {\bf p}_{\,y}\ )\cdot  \,
 \Vert\ y \ - \ {\bf p}_{\,y}\,\Vert^{\,\ell} \ - \ {\bf R}\,(\,y\,) \ ]  \,\cdot\, {\widetilde{\widetilde{[\ \cdot \cdot \cdot \ ]}}}_{\,(\,A.11.9\,)}   \  dy \\[0.2in]
& = &  \int_{B_{\,\xi_{\,1}}\,(\,\rho_{\,\nu}\,)   }  C\,(\,{\bf p}_{\,y}\,) \,\cdot \Vert \,y  \ - \ {\bf p}_{\,y}\,\Vert^{\,\ell} \,\cdot\, {\widetilde{\widetilde{[\ \cdot \cdot \cdot \ ]}}}_{\,(\,A.11.9\,)}   \  dy  \ + \ {\bf E}^{\,\mbox{Remainder}}_{\, (\,A.11.10\,) }
 \\[0.2in]
 & = & C\,(\,{\bf p}_{\,1}\,) \,\cdot    \int_{B_{\,\xi_{\,1}}\,(\,\rho_{\,\nu}\,)   }    \Vert \,y  \ - \ {\bf p}_{\,y}\,\Vert^{\,\ell} \,\cdot\, {\widetilde{\widetilde{[\ \cdot \cdot \cdot \ ]}}}_{\,(\,A.11.9\,)}   \  dy \ + \ \, {\bf E}^{\,\mbox{Remainder}}_{\, (\,A.11.10\,) } \ + \ {\bf E}^{\,\mbox{Shift}}_{\, (\,A.11.10\,) }
 \\[0.2in]
  & = & C\,(\,{\bf p}_{\,1}\,) \,\cdot  \ell\,\cdot\,   \int_{B_{\,\xi_{\,1}}\,(\,{\rho_{\,\nu}})\ \cap \ \{ \ y_n \ >  \ \xi_{1_{\,|_{\,n}}}\ \} } |\ y_n \ - \ \xi_{1_{\,|_{\,n}}}\,|^{\,\,\ell \, - \ 1} \ \cdot \,(\,+\, \xi_{1_{\,|_{\,n}}}\,)\,\cdot\, {\widetilde{\widetilde{[\ \cdot \cdot \cdot \ ]}}}_{\,(\,A.11.9\,)}   \  dy \ + \  \\[0.2in]
   & \ & \ \   + \ C\,(\,{\bf p}_{\,1}\,) \,\cdot  \ell\,\cdot\,   \int_{B_{\,\xi_{\,1}}\,(\,{\rho_{\,\nu}})\ \cap \ \{ \ y_n \ <  \ 0\ \} } |\ y_n \ - \ \xi_{1_{\,|_{\,n}}}\,|^{\,\,\ell \, - \ 1} \ \cdot \,(\,-\, \xi_{1_{\,|_{\,n}}}\,)\,\cdot\, {\widetilde{\widetilde{[\ \cdot \cdot \cdot \ ]}}}_{\,(\,A.11.9\,)}   \  dy \ + \  \\[0.2in]
  & \ &    + \ \, {\bf E}^{\,\mbox{Remainder}}_{\, (\,A.11.10\,) } \ + \ \, {\bf E}^{\,\mbox{Shift}}_{\, (\,A.11.10\,) } \ + \ {\bf E}^{\,\mbox{1st}}_{\, (\,A.11.10\,) } \ + \ {\bf E}^{\,\mbox{3rd}}_{\, (\,A.11.10\,) } \ + \ {\bf E}^{\,\mbox{Lower}}_{\, (\,A.11.10\,) }  \ + \ {\bf E}^{\,\mbox{Thin}}_{\, (\,A.11.10\,) } \  .\\
\end{eqnarray*}

Via (\,A.11.6\,)\,,\,
  $\,{\bf E}^{\,\mbox{Remainder}}_{\, (\,A.11.10\,) }\,$ ${\bf E}^{\,\mbox{Shift}}_{\, (\,A.11.10\,) } $\,,\, ${\bf E}^{\,\mbox{Thin}}_{\, (\,A.11.10\,) } $\, and ${\bf E}^{\,\mbox{Lower}}_{\, (\,A.11.10\,) } $\,,\, are estimated as in {\bf \S\,A\,9\,.b}\,,\, {\bf \S\,A\,9\,.c}\,,\, {\bf \S\,A\,9\,.e}\, and \,{\bf \S\,A\,9\,.g}\,,\, respectively\,.\, We note that the argument in  {\bf \S\,A\,9\,.e}\, works here once we replace the bound in (\,A.9.22\,) by \\[0.1in]
   (\,A.11.11\,)
\begin{eqnarray*}
& \ & \Bigg\vert\  {{\lambda_1 \cdot (\ y_n \ - \  \xi_{1_{\,|_n}}\,)}\over {\lambda^2_1 \ + \ \Vert\,y\,-\,\xi_{\,1}\,\Vert^2}}\ \Bigg\vert  \ \le \ \Bigg\vert\  {{\lambda_1^2 \ +  (\ y_n \ - \  \xi_{1_{\,|_n}} \ )^{\,2}}\over {\lambda^2_1 \ + \ \Vert\,y\,-\,\xi_{\,1}\,\Vert^2}}\ \Bigg\vert  \ \le  \ 1 \ \ \ \ \ \ \mfor \  \ y \ \in  \ \R^n\\[0.2in]
\Longrightarrow  & \ & \Bigg\vert\   \left(  {{\lambda_{\,1}}\over {\lambda_{\,1}^2 \ + \ \Vert\, y \,-\, \xi_{\,1}\Vert^2}}\  \right)^{\!\!n\,+\,1} \!\!\cdot\,  (\, y_n \,-\,\xi_{1_{\,|_{\,n}}}\, )\  \Bigg\vert \ \le \ \left(  {{\lambda_{\,1}}\over {\lambda_{\,1}^2 \ + \ \Vert\, y \,-\, \xi_{\,1}\Vert^2}}\  \right)^{\!\!n }  \\
\end{eqnarray*}
for \,$\ y \ \in  \ \R^n$\,.\,
  Whereas $\,{\bf E}^{\,\mbox{1st}}_{\, (\,A.9.10\,) }\,$  and  $\,{\bf E}^{\,\mbox{3rd}}_{\, (\,A.9.10\,) }\,$ are related to the first and third teams in (\,A.11.7\,)\,,\, respectively\,,\, and are estimated in the following subsections.


\vspace*{0.5in}

{\bf{\S\,A\,11\,.\,b}} \ \ {\bf First term\,} $\ {\bf E}^{\,\mbox{1st}}_{\, (\,A.11.10\,) } \,.$ \ \ Thus compare
\begin{eqnarray*}
& \ & \\
& \  & \int_{B_{\,\xi_{\,1}}\,(\,{\rho_{\,\nu}})\ \cap \ \{ \ y_n \ >  \ \xi_{1_{\,|_{\,n}}} \}}   \!\! (\,y_n\,-\,\xi_{1_{\,|_{\,n}}}\,)^{\,\ell} \cdot \left\{   \left(  {{\lambda_{\,1}}\over {\lambda_{\,1}^2 \ + \ \Vert\, y \,-\, \xi_{\,1}\Vert^2}}\  \right)^{\!\!n\,+\,1} \!\!\cdot\,  (\, y_n \,-\,\xi_{1_{\,|_{\,n}}}\, ) \right\}\ dy  \\[0.2in]
 & = &  +\,\int_{B_{o}\,(\,{\rho_{\,\nu}})\ \cap \ \{ \ {\bar y}_n \ >  \ 0\ \}}  \ \ {\bar y}_n^{\,\ell\,+\,1}  \,\cdot\,\, \left\{    \,\left(  {{\lambda_{\,1}}\over {\lambda_{\,1}^2 \ + \ \Vert\, {\bar y}\,\Vert^2}}\  \right)^{\!\!n\,+\,1}  \ \right\}\ d\,{\bar y} \ \ \ (\  >  \ 0\,) \\[0.1in]
 & \ & \hspace*{4.3in} \ \ \ \  \ \ (\ {\bar y} \ = \ y \ - \ \xi_{\,1} \ )\ ,\\
\end{eqnarray*}

with
\begin{eqnarray*}
& \ & \int_{B_{\,\xi_{\,1}}\,(\,{\rho_{\,\nu}})\ \cap \ \{ \ y_n \ <  \ 0 \ \}}
 \ \   (\,-\,y_n\,+\,\xi_{1_{\,|_{\,n}}}\,)^{\,\ell}  \,\cdot\,\, \left\{    \,\left(  {{\lambda_{\,1}}\over {\lambda_{\,1}^2 \ + \ \Vert\, y \,-\, \xi_{\,1}\Vert^2}}\  \right)^{\!\!n\,+\,1} \,\cdot\,  (\, y_n \ -\ \xi_{1_{\,|_{\,n}}}\, )  \ \right\}\ dy \\[0.1in]
& \ & \hspace*{4.9in}\leftarrow  ( \, \uparrow \, {\mbox{-\,ve}}\, )\rightarrow  \\[0.1in]
& = & \int_{B_{\,\xi_{\,1}}\,(\,{\rho_{\,\nu}})\ \cap \ \{ \ y_n \ <  \ 0 \ \}} \ \   (\,-\,y_n\,+\,\xi_{1_{\,|_{\,n}}}\,)^{\,\ell}    \ *\\[0.2in]
& \ & \hspace*{2.3in}   \ \ \ \ \ \ *\  \left\{    \,\left(  {{\lambda_{\,1}}\over {\lambda_{\,1}^2 \ + \ \Vert\, y \,-\, \xi_{\,1}\Vert^2}}\  \right)^{\!\!n\,+\,1} \,\cdot\,[\  -\, (\, -\, y_n \,+\,\xi_{1_{\,|_{\,n}}}\, )  \ ] \,\right\}\ dy\\[0.2in]
& = & -\,\int_{B_{\,\xi_{\,1}}\,(\,{\rho_{\,\nu}})\ \cap \ \{ \ y_n \ <  \ 0 \ \}} \ \   (\,-\,y_n\,+\,\xi_{1_{\,|_{\,n}}}\,)^{\,\ell \ + \ 1} \,\cdot \,\left\{    \,\left(  {{\lambda_{\,1}}\over {\lambda_{\,1}^2 \ + \ \Vert\, y \,-\, \xi_{\,1}\Vert^2}}\  \right)^{\!\!n\,+\,1}     \,\right\}\ dy  \\[0.2in]
& = & -\, \int_{B_{-\,\xi_{\,1}}\,(\,{\rho_{\,\nu}})\ \cap \ \{ \ z_n \ >  \  0 \ \}} \ \ (\,z_n\,+\,\xi_{1_{\,|_{\,n}}}\,)^{\,\ell \ + \ 1}  \,\cdot\,\, \left\{    \,\left(  {{\lambda_{\,1}}\over {\lambda_{\,1}^2 \ + \ \Vert\, z \,+\, \xi_{\,1}\Vert^2}}\  \right)^{\!\!n\,+\,1}     \,\right\}\ (\,+\,1\,)\,\cdot\,d\,z\emph{}\\[0.1in]
& \ & \ \ \  (\,z_1\ = \ y_1\,, \ \cdot\,\cdot \cdot\,, \ \ z_{n\,-\,1}\ = \ y_{n\,-\,1}\,, \ \ z_n \ = \ -\,y_n\ ) \hspace*{0.8in}   |\,{\mbox{Det}}\,J\,| \ = \ +\,1 \ \ \uparrow\\[0.2in]
& = &-\,\int_{B_{o}\,(\,{\rho_{\,\nu}})\ \cap \ \{ \ {\bar z}_n \ >  \ \xi_{1_{\,|_{\,n}}}  \ \}}  \ \ (\,{\bar z}_n\,)^{\,\ell\,+\,1}  \,\cdot\,\, \left\{    \,\left(  {{\lambda_{\,1}}\over {\lambda_{\,1}^2 \ + \ \Vert\  {\bar z}\ \Vert^2}}\  \right)^{\!\!n\,+\,1}     \,\right\}\,\cdot\,d\,{\bar z}_n \ .\\[0.1in]
& \ & \uparrow \ \ \hspace*{2.9in}   (\,{\bar z}_n\ =\ z_n \ + \ \xi_{1_{\,|_{\,n}}}\,) \ \  \uparrow  \ \ \ \ \
\end{eqnarray*}

Thus we observe the self\,-\,cancellation, except the thin part in between, which is estimated as in {\bf \S\,A\,9\,.g} [ \ see also (\,A.11.11\,)\ ]\,.

\vspace*{0.3in}

{\bf{\S\,A\,11\,.c}} \ \ {\bf Third term\ } ${\bf E}^{\,\mbox{3rd}}_{\, (\,A.11.10\,) }\ . $\ \
As in (\,A.9.31\,)\,,\,  compare
\begin{eqnarray*}
& \ & \xi_{1_{\,|_{\,n}}}^2\,\cdot\,\int_{B_{\,\xi_{\,1}}\,(\,{\rho_{\,\nu}})\ \cap \ \{ \ y_n \ >  \ \xi_{1_{\,|_{\,n}}} \,\}}  \ \ (\,y_n\,-\,\xi_{1_{\,|_{\,n}}}\,)^{\,\ell\ - \ 2}    \ *\\[0.1in]
& \ & \hspace*{2.3in}   \ \ \ \ \ \ *\  \left\{    \,\left(  {{\lambda_{\,1}}\over {\lambda_{\,1}^2 \ + \ \Vert\, y \,-\, \xi_{\,1}\Vert^2}}\  \right)^{\!\!n\,+\,1} \,\cdot\,  (\, y_n \,-\,\xi_{1_{\,|_{\,n}}}\, )  \,\right\}\ dy \\[0.1in]
 & = & +\    \xi_{1_{\,|_{\,n}}}^2\,\cdot\,\int_{B_{o}\,(\,{\rho_{\,\nu}})\ \cap \ \{ \ {\bar y}_n \ >  \ 0\ \} }  \ \ {\bar y}_n^{\,\ell\ -\,1}  \,\cdot\,\, \left\{    \,\left(  {{\lambda_{\,1}}\over {\lambda_{\,1}^2 \ + \ \Vert\, {\bar y}\,\Vert^2}}\  \right)^{\!\!n\,+\,1} \,\cdot\,  {\bar y}_n  \,\right\}\ d\,{\bar y} \\[0.1in]
 & \ & (\ \uparrow \ \ {\mbox{\,+\,ve}}\,) \hspace*{3in} \ \ \ \ \ \  \ \ (\ {\bar y} \ = \ y \ - \ \xi_{\,1} \ )\ ,
\end{eqnarray*}

with
\begin{eqnarray*}
& \ & \xi_{1_{\,|_{\,n}}}^2\,\cdot\,\int_{B_{\,\xi_{\,1}}\,(\,{\rho_{\,\nu}})\ \cap \ \{ \ y_n \ <  \ 0 \ \}}  \ \   (\,-\,y_n\,+\,\xi_{1_{\,|_{\,n}}}\,)^{\,\ell\ - \ 2}    \ *\\[0.1in]
& \ & \hspace*{2.3in}   \ \ \ \ \ \ *\ \left\{    \,\left(  {{\lambda_{\,1}}\over {\lambda_{\,1}^2 \ + \ \Vert\, y \,-\, \xi_{\,1}\Vert^2}}\  \right)^{\!\!n\,+\,1} \,\cdot\,  (\, y_n \,-\,\xi_{1_{\,|_{\,n}}}\, )  \,\right\}\ dy \\[0.1in]
& \ & \hspace*{4in} \ \ \ \ \ \ \ \ \ \ \ \ \ \ \ \ \ \ \ \  ( \ \uparrow \ {\mbox{-\,ve}}\ )\\[0.1in]
& = & \xi_{1_{\,|_{\,n}}}^2 \cdot\int_{B_{\,\xi_{\,1}}\,(\,{\rho_{\,\nu}})\ \cap \ \{ \ y_n \ <  \  0 \ \}} \ \   (\,-\,y_n\,+\,\xi_{1_{\,|_{\,n}}}\,)^{\,\ell \ - \ 2}  *\\[0.2in]
& \ & \hspace*{2in} *\, \left\{    \,\left(  {{\lambda_{\,1}}\over {\lambda_{\,1}^2 \ + \ \Vert\, y \,-\, \xi_{\,1}\Vert^2}}\  \right)^{\!\!n\,+\,1} \,\cdot\,[\  -\, (\, -\, y_n \,+\,\xi_{1_{\,|_{\,n}}}\, )  \ ] \,\right\}\ dy\\[0.2in]
& = & -\,\xi_{1_{\,|_{\,n}}}^2 \cdot\int_{B_{\,\xi_{\,1}}\,(\,{\rho_{\,\nu}})\ \cap \ \{ \ y_n \ <  \ 0 \ \}}   \ \   (\,-\,y_n\,+\,\xi_{1_{\,|_{\,n}}}\,)^{\,\ell \ - \ 1}  \,\cdot\,\, \left\{    \,\left(  {{\lambda_{\,1}}\over {\lambda_{\,1}^2 \ + \ \Vert\, y \,-\, \xi_{\,1}\Vert^2}}\  \right)^{\!\!n\,+\,1}     \,\right\}\ dy  \\[0.2in]
& = & -\, \xi_{1_{\,|_{\,n}}}^2 \cdot\int_{B_{-\,\xi_{\,1}}\,(\,{\rho_{\,\nu}})\ \cap \ \{ \ z_n \ >  \ 0\ \}}  \ \ (\,z_n\,+\,\xi_{1_{\,|_{\,n}}}\,)^{\,\ell \ - \ 1}    \ *\\[0.1in]
& \ & \hspace*{2.3in}   \ \ \ \ \ \ *\  \left\{    \,\left(  {{\lambda_{\,1}}\over {\lambda_{\,1}^2 \ + \ \Vert\, z \,+\, \xi_{\,1}\Vert^2}}\  \right)^{\!\!n\,+\,1}     \,\right\}\ (\,+\,1\,)\,\cdot\,dz  \hspace*{4in}    |\,{\mbox{Det}}\,J\,| \ = \ +\,1 \ \  \uparrow \ \\[0.2in]
& = &-\ \xi_{1_{\,|_{\,n}}}^2 \cdot\int_{B_{o}\,(\,{\rho_{\,\nu}})\ \cap \ \{ \ {\bar z}_n \ >  \ \xi_{1_{\,|_{\,n}}}  \ \}} \emph{} \ \ (\,{\bar z}_n\,)^{\,\ell\ -\ 1}  \,\cdot\,\, \left\{    \,\left(  {{\lambda_{\,1}}\over {\lambda_{\,1}^2 \ + \ \Vert\, {\bar z}\Vert^2}}\  \right)^{\!\!n\,+\,1}     \,\right\}\,\cdot\,d\,{\bar z}_n \\[0.1in]
& \ & \  ( \ \uparrow \ {\mbox{-\,ve}}\ ) \ \ \  [ \ z_1\ = \ y_1\,, \ \cdot\,\cdot \cdot\,, \ \ z_{n\,-\,1}\ = \ y_{n\,-\,1}\,, \ \ z_n \ = \ -\,y_n\ , \ \ \  {\bar z}_n\ =\ z_n \ + \ \xi_{1_{\,|_{\,n}}}\ ]\ .\\
\end{eqnarray*}

Again, we have the self\,-\,cancellation, leaving the thin part in between, which is estimated as  in {\bf \S\,A\,9\,.g}  [ \ cf. also (\,A.11.11\,)\ ]\,.

\vspace*{0.3in}

{\bf{\S\,A\,11\,.\,d}} \ \ {\bf The second term\ } ${\bf E}^{\,\mbox{2nd}}_{\, (\,A.11.10\,) }\ . $\ \
Here we have the following coupling.
\begin{eqnarray*}
& \ & \ell\,\int_{B_{\,\xi_{\,1}}\,(\,{\rho_{\,\nu}})\ \cap \ \{ \ y_n \ >  \ \xi_{1_{\,|_{\,n}}}\}}  \ \ (\,y_n\,-\,\xi_{1_{\,|_{\,n}}}\,)^{\,\ell\ - \ 1}  \ \cdot\,\ \xi_{1_{\,|_{\,n}}}  * \\[0.2in]
& \ & \hspace*{2.3in} *\  \left\{    \,\left(  {{\lambda_{\,1}}\over {\lambda_{\,1}^2 \ + \ \Vert\, y \,-\, \xi_{\,1}\Vert^2}}\  \right)^{\!\!n\,+\,1} \,\cdot\,  (\, y_n \,-\,\xi_{1_{\,|_{\,n}}}\, )  \,\right\}\ dy  \\[0.3in]
 & = & +\, \ell\,\cdot\,  \xi_{1_{\,|_{\,n}}}\,\cdot\, \int_{B_{o}\,(\,{\rho_{\,\nu}})\ \cap \ \{ \ {\bar y}_n \ >  \ 0\,\}}  \ \ \ |\, {\bar y}_n\,|^{\,\ell}  \,\cdot\,\, \left\{    \,\left(  {{\lambda_{\,1}}\over {\lambda_{\,1}^2 \ + \ \Vert\, {\bar y}\,\Vert^2}}\  \right)^{\!\!n\,+\,1} \, \,\right\}\  d\,{\bar y} \\[0.2in]
 & \ & \hspace*{4in} \ \   \ \ \ \ (\ {\bar y} \ = \ y \ - \ \xi_{\,1} \ > \ 0 \ )\ ,
\end{eqnarray*}


\newpage

with
\begin{eqnarray*}
& \ & \ell \cdot\int_{B_{\,\xi_{\,1}}\,(\,{\rho_{\,\nu}})\ \cap \ \{ \ y_n \ <  \ 0\,\}}  \ \   (\,-\,y_n\,+\,\xi_{1_{\,|_{\,n}}}\,)^{\,\ell\ - \ 1} \,\cdot\,(\,-\, \xi_{1_{\,|_{\,n}}} \,)\ *\\[0.1in]
& \ & \hspace*{2.3in}*\, \left\{    \,\left(  {{\lambda_{\,1}}\over {\lambda_{\,1}^2 \ + \ \Vert\, y \,-\, \xi_{\,1}\Vert^2}}\  \right)^{\!\!n\,+\,1} \,\cdot\,  (\, y_n \,-\,\xi_{1_{\,|_{\,n}}}\, )  \,\right\}\ dy  \\[0.3in]
& = & \ell \cdot(\,-\, \xi_{1_{\,|_{\,n}}} \,)\,\cdot\, \int_{B_{\,\xi_{\,1}}\,(\,{\rho_{\,\nu}})\ \cap \ \{ \ y_n \ < \ 0 \ \} }  (\,-\,y_n\,+\,\xi_{1_{\,|_{\,n}}}\,)^{\,\ell \ - \ 1}  * \\[0.1in]
& \ & \hspace*{2in} *\, \left\{    \,\left(  {{\lambda_{\,1}}\over {\lambda_{\,1}^2 \ + \ \Vert\, y \,-\, \xi_{\,1}\Vert^2}}\  \right)^{\!\!n\,+\,1} \,\cdot\,[\  -\, (\, -\, y_n \,+\,\xi_{1_{\,|_{\,n}}}\, )  \ ] \,\right\}\ dy\\[0.3in]
& = & \ell \cdot(\,+\, \xi_{1_{\,|_{\,n}}} \,)\,\cdot\,\int_{B_{\,\xi_{\,1}}\,(\,{\rho_{\,\nu}})\ \cap \ \{ \ y_n \ <  \ 0\,\}}  \ \   (\,-\,y_n\,+\,\xi_{1_{\,|_{\,n}}}\,)^{\,\ell }  \,\cdot\,\, \left\{    \,\left(  {{\lambda_{\,1}}\over {\lambda_{\,1}^2 \ + \ \Vert\, y \,-\, \xi_{\,1}\Vert^2}}\  \right)^{\!\!n\,+\,1}     \,\right\}\ dy  \\[0.3in]
& = & \ell \cdot\xi_{1_{\,|_{\,n}}}\,\cdot\,\int_{B_{-\,\xi_{\,1}}\,(\,{\rho_{\,\nu}})\ \cap \ \{ \ z_n \ >  \ 0\,\}}  \ \ (\,z_n\,+\,\xi_{1_{\,|_{\,n}}}\,)^{\,\ell }  \,\cdot\,\, \left\{    \,\left(  {{\lambda_{\,1}}\over {\lambda_{\,1}^2 \ + \ \Vert\, z \,+\, \xi_{\,1}\Vert^2}}\  \right)^{\!\!n\,+\,1}     \,\right\}\ (\,+\,1\,)\,\cdot\,dz  \\[0.1in]
& \ & \ [ \ \uparrow \ {\mbox{+\,ve}}\,, \ \ \ \ z_1\ = \ y_1\,, \ \cdot\,\cdot \cdot\,, \ \ z_{n\,-\,1}\ = \ y_{n\,-\,1}\,, \ \ z_n \ = \ -\,y_n\ , \ \  |\,{\mbox{Det}}\,J\,| \ = \ +\,1  \ \  \uparrow \ ]\\[0.3in]
& = & + \ \ell \cdot\xi_{1_{\,|_{\,n}}}\,\cdot\,\int_{B_{o}\,(\,{\rho_{\,\nu}})\ \cap \ \{ \ {\bar z}_n \ >  \ \xi_{1_{\,|_{\,n}}}\,\}}  \ \ (\,{\bar z}_n\,)^{\,\ell }  \,\cdot\,\, \left\{    \,\left(  {{\lambda_{\,1}}\over {\lambda_{\,1}^2 \ + \ \Vert\, {\bar z}\Vert^2}}\  \right)^{\!\!n\,+\,1}     \,\right\}\,\cdot\,d\,{\bar z} \ .\\[0.1in]
& \ & \ \ \ \hspace*{0.8in}\leftarrow \ \ {\mbox{unit \ \ for \ \ this \ \ calculation}}\hspace*{0.8in} \ \ \rightarrow \ \  \ \ \ \ \ \\[0.1in]& \ &  \hspace*{4in}\ \ \ \ \ \ \ \ \ \ (\,{\bar z}_n\ =\ z_n \ + \ \xi_{1_{\,|_{\,n}}}\,)\ .
\end{eqnarray*}
Similarly we calculate the lower half. Note that
$$
\int_{\R^n}   (\,{\bar z}_n\,)^{\,\ell }  \,\cdot\,\, \left\{    \,\left(  {{\lambda_{\,1}}\over {\lambda_{\,1}^2 \ + \ \Vert\, {\bar z}\Vert^2}}\  \right)^{\!\!n\,+\,1}     \,\right\}\,\cdot\,d\,{\bar z} \ = \ {{\ \lambda_{\,1}^{\,\ell} \ }\over \lambda_{\,1}} \cdot \int_{\R^n}   (\,{\bar Z}_n\,)^{\,\ell }  \,\cdot\,\, \left\{    \,\left(  {{\lambda_{\,1}}\over {\lambda_{\,1}^2 \ + \ \Vert\, {\bar Z}\Vert^2}}\  \right)^{\!\!n\,+\,1}     \,\right\}\,\cdot\,d\,{\bar Z} \ ,
$$
where ${\bar Z}  \ = \ {\bar z}  \cdot  \lambda_{\,1}^{\,-\,1}\ .$\, Recall that $\,\ell \  \le \ n \ - \ 2\,$ so that the integral is finite.

\newpage

{\bf{\S\,A\,11\,.\,e}} \ \ {\bf Main term in the derivative in the $\,\xi_{1_{\,|_n}}$-\,direction\,.\,} Combining the discussion in {\bf{\S\,A\,11\,.a}}\,,\, {\bf{\S\,A\,11\,.b}}\,,\, {\bf{\S\,A\,11\,.c}} \,and\, {\bf{\S\,A\,11\,.d}}\,,\, we obtain \\[0.1in]
(\,A.11.12\,)
\begin{eqnarray*}
&  &  \int_{\R^n}\, [\,(\,{\tilde c}_n\!\cdot K\,) \ - \ n\,(\,n\,-\,2)\,]  \cdot\,(W_{\,\,\flat} )^{{\,n\,+\,2\,}\over {n\,-\,2}}\,\cdot\, \left(\,\lambda_{\,1}\,\cdot\,{\partial\over {\partial\, \xi_{1_{\,|_n}}}}\ \right)\\[0.2in]
 & = &  \left\{ \ \int_{B_o\,(\,\xi_{\,1}\,)}\ \ -\  [\,n\,(\,n\,-\,2) \ - \ (\,{\tilde c}_n\!\cdot K\,) \ ]  \cdot\,V_1^{{\,n\,+\,2\,}\over {n\,-\,2}}\,\cdot\,\left[\ \left(\,\lambda_{\,1}\,\cdot\,{\partial\over {\partial\, \xi_{1_{\,|_n}}}}\  \right)\,\ \right]  \ \right\}\ *\\[0.2in]
  & \ & \ \ \ \ \ \ \  \ \ \ \ \ \ \ \ \ \ \  \ \hspace*{2in} *\ \left[\ 1 \ + \ O\,\left( \ {\bar\lambda}_{\ \flat}^{\ (\,n\ -\,2\,)\,\cdot\, [\ (\,\gamma \ + \ \nu \,)\ - \  1\,] }\ \right) \  \right] \\[0.2in]
  & \ & \ \ \ \ \ \ \ \  \ \ \  \ \ \ \ \ \ \ \  \ \ \  \ \ \ \ \ \ \ \  \ \ \ +\ \ O \left(\  {\bar\lambda}_{\ \flat}^{\ n\,(\,1\,-\ \nu\,)} \,\cdot\, {\bar\lambda}_{\ \flat}^{  \   - \ 2\,\varepsilon \,[\ (\,\gamma \ + \ \nu\,) \ - \ 1\ ] } \  \right) \\[0.2in]
  & = &   \Bigg\{   -\, \left[\  C_+  \ + \ O \,\left( \ {\bar \lambda}_{\,\,\flat}^{(\ n\,-\ \ell\,)\, \cdot\, (\,1 \ - \ \nu) }\ \right) \ \right]\,\cdot\,\lambda^\ell_1\,\cdot\,{\eta\over \lambda_{\,1}}   \\[0.2in] & \ & \ \ \ \   + \  O \,\left( \ {\bar \lambda}_{\,\,\flat}^{\, \ell \ + \ 1 \ - \ o_{\,+}\,(\,1\,)\,}   \ \right) \ + \ O \,\left( \ {\bar \lambda}_{\,\,\flat}^{\, \ell \ + \ 3\,\kappa\,}   \ \right) \    + \    O \left(\
   {\bar\lambda}_{\ \flat}^{\ \ell \ + \  (\,\ell\,+\,1\,) \,\cdot \,(\, 2\,\nu \ - \ 1\,)} \ \right)\ \Bigg\} \ *  \ \\[0.2in]
  & \ & \ \ \ \ \ \ \  \ \ \ \  \ \hspace*{2in} \ \ \ \ \ \ \ \ \ \ \  * \ \left[\ 1 \ + \ O\,\left( \ {\bar\lambda}_{\ \flat}^{\ (\,n\ -\,2\,)\,\cdot\, [\ (\,\gamma \ + \ \nu \,)\ - \  1\,] }\ \right) \  \right] \\[0.2in]
  & \ & \ \ \ \ \ \ \ \  \ \ \ \  \ \ \ \ \ \ \ \  \ \ \  \ \ \ \ \ \ \ \  \ \ \  \ \ \ \ \ \ \ \  \ \ \  +  \ \,O \left(\  {\bar\lambda}_{\ \flat}^{\ n\,(\,1\,-\ \nu\,)} \,\cdot\, {\bar\lambda}_{\ \flat}^{  \   - \ 2\,\varepsilon \,[\ (\,\gamma \ + \ \nu\,) \ - \ 1\ ] } \  \right)  \\[0.2in]
  & \ & \ \ \ \ \ \ \ \  \ \ \ \  \ \ \ \ \ \ \ \  \ \ \  \ \ \ \ \ \ \ \  \ \ \  \ \ \ \ \ \ \ \  \ \ \  +  \ \,O \left(\
  {\bar\lambda}_{\ \flat}^{\ \ell } \cdot  {\bar\lambda}_{\ \flat}^{\  (\,\ell\,+\,1\,) \,\cdot \,(\, 2\,\nu \ - \ 1\,)} \ \right)\\[0.2in]
  & = & -\,C_+\,\cdot\,\lambda_{\,1}^\ell \cdot {{\ \eta_1\ }\over \lambda_{\,1}}  \ + \      {\cal E}_{\,(\,A.{\bf 9}.42\,)} \ .
   \end{eqnarray*}

Here
$$
\eta_{\,1} \ := \   {\mbox{dist}}\,(\,\xi_{\,1}\,,\ {\cal H}\,) \ ( \ = \ \xi_{1_{\,|_n}}  ) \ \le  \   {\bar\lambda}_{\,\,\flat}^{\,1\ + \ \kappa} \ ,
$$
and
$$
C_+ \ = \ C\,(\,{\bf p}_1\,)\cdot (\,n \ - \ 2\,)\cdot \ell \cdot
  \int_{\R^n}   (\,{\bar Z}_n\,)^{\,\ell }  \,\cdot\,\, \left\{    \,\left(  {{\lambda_{\,1}}\over {\lambda_{\,1}^2 \ + \ \Vert\, {\bar Z}\Vert^2}}\  \right)^{\!\!n\,+\,1}     \,\right\}\,\cdot\,d\,{\bar Z} \ . \leqno (\,A.11.13\,)
$$
\hspace*{2in} \ \ \ \ \ \ \ \ \  $\uparrow$ \ \ see (\,A.11.5\,)

\newpage

{\bf{\S\,A\,11\,.\,f} } \ \ {\bf Derivative in the other directions\,.\,} \ \
$$  \int_{\R^n}\, [\,(\,{\tilde c}_n\!\cdot K\,) \ - \ n\,(\,n\,-\,2)\,]  \cdot\,(W_{\,\,\flat} )^{{\,n\,+\,2\,}\over {n\,-\,2}}\,\cdot\ \left(\,\lambda_{\,1}\,\cdot\,{\partial\over {\partial\, \xi_{1_{\,|_j}}}}\  \right)\  \ \ \ \  \ \ \ (\ j \ \not= \ n\ )\,. \leqno (\,A.11.14\,)
$$
Proceeding as in (\,A.11.10) and using (\,A.11.4\,), here we find ourselves dealing with terms like (\ cf. \S\,A\,11\,.a\,--\,c \ )
\begin{eqnarray*}
(\,A.11.15\,)_{\,+}  \ \ \ \ \ \ \ \ \ \ \ \ & \  & \int_{B_{\,\xi_{\,1}}\,(\,{\rho_{\,\nu}})\ \cap \ \{ \ y_n \ >  \ \xi_{1_{\,|_{\,n}}} \}}  \ \ (\,y_n\,-\,\xi_{1_{\,|_{\,n}}}\,)^{\,\ell}  \,\cdot\,\,{\widetilde{\widetilde{[\ \cdot \cdot \cdot \ ]}}}_{\,j} \ dy \ , \ \ \ \ \ \ \ \ \ \ \ \ \ \ \ \ \ \ \ \ \ \ \ \  \ \ \ \ \ \ \ \ \ \ \ \ \\[0.2in]
(\,A.11.15\,)_{\,-}  \ \ \ \ \ \ \ \ \ \ \ \ & \  &  \int_{B_{\,\xi_{\,1}}\,(\,{\rho_{\,\nu}})\ \cap \ \{ \ y_n \ <   \ 0 \}}  \ \   (\,-\,y_n\,+\,\xi_{1_{\,|_{\,n}}}\,)^{\,\ell}  \,\cdot\,\, {\widetilde{\widetilde{[\ \cdot \cdot \cdot \ ]}}}_{\,j} \ dy\ ;
 \end{eqnarray*}
 \begin{eqnarray*}
(\,A.11.16\,)_{\,+}  \ \ \ \ \ \ \ \ \ \ \ \ & \ & \int_{B_{\,\xi_{\,1}}\,(\,{\rho_{\,\nu}})\ \cap \ \{ \ y_n \ >  \ \xi_{1_{\,|_{\,n}}}\}}  \ \ (\,y_n\,-\,\xi_{1_{\,|_{\,n}}}\,)^{\,\ell\ - \ 1}  \ \cdot\,\ \xi_{1_{\,|_{\,n}}}  \,\cdot\,\, {\widetilde{\widetilde{[\ \cdot \cdot \cdot \ ]}}}_{\,j}\ dy \ \ \ \ \ \ \ \ \ \ \ \ \ \ \ \ \ \ \ \ \ \ \ \  \ \ \ \ \ \ \ \ \ \ \ \ \\[0.2in]
(\,A.11.16\,)_{\,-}    \ \ \ \ \ \ \ \ \ \ \ \ & \ & \int_{B_{\,\xi_{\,1}}\,(\,{\rho_{\,\nu}})\ \cap \ \{ \ y_n \ <  \ 0 \ \}}
\   (\,-\,y_n\,+\,\xi_{1_{\,|_{\,n}}}\,)^{\,\ell\ - \ 1} \,\cdot\,(\,-\, \xi_{1_{\,|_{\,n}}} \,)\ {\widetilde{\widetilde{[\ \cdot \cdot \cdot \ ]}}}_{\,j}\ \ dy  \ ;\\[0.2in]
(\,A.11.17\,)_{\,+}  \ \ \ \ \ \ \ \ \ \ \ \, & \ & \int_{B_{\,\xi_{\,1}}\,(\,{\rho_{\,\nu}})\ \cap \ \{ \ y_n \ >  \ \xi_{1_{\,|_{\,n}}}\}}  \ \ (\,y_n\,-\,\xi_{1_{\,|_{\,n}}}\,)^{\,\ell\ - \ 2}  \ \cdot\,\ \xi_{1_{\,|_{\,n}}}^2  \,\cdot\,\,{\widetilde{\widetilde{[\ \cdot \cdot \cdot \ ]}}}_{\,j}\ \ dy    \ ,\\[0.2in]
(\,A.11.17\,)_{\,-}    \ \ \ \ \ \  \ \ \ \ \, & \ & \int_{B_{\,\xi_{\,1}}\,(\,{\rho_{\,\nu}})\ \cap \ \{ \ y_n \ <  \ 0\,\}}  \ \   (\,-\,y_n\,+\,\xi_{1_{\,|_{\,n}}}\,)^{\,\ell\ - \ 2} \,\cdot\,(\,-\, \xi_{1_{\,|_{\,n}}}^2 \,)\ \cdot \  {\widetilde{\widetilde{[\ \cdot \cdot \cdot \ ]}}}_{\,j}\ \ dy  \ .
\end{eqnarray*}
In the above\,,\,
$$
{\widetilde{\widetilde{[\ \cdot \cdot \cdot \ ]}}}_{\,j} \ = \ \left( {{\lambda_{\,1}}\over {\  \lambda_{\,1}^2 \ + \ \Vert\,y \ - \ \xi_{\,1}\,\Vert^{\,2}\ }}\  \right)^{\!\!n \ +\ 1} \ \cdot (\,y_j \ - \  \xi_{1_{\,|_j}} \ )\ .\leqno  (\,A.11.18\,)
$$
Note that the presence of the singular team
 $$
 (\,y_j \ - \ \xi_{1_{\,|_j}}\,)\ , \leqno  (\,A.11.19\,)
 $$
causing internal cancellation (\,``\,left\,-\,right\," cancellation versa the ``\,up\,-\,down\," cancellation in {\bf \S\,A\,11\,.\,b\,} and {\bf \S\,A\,11\,.\,c\,})\,.
\, For example,   in the case of  $(\,A.11.15\,)_{\,+}$\,.\,
 Recall that
 $$
 \xi_{\,1} \ = \ (\ 0\,, \ \cdot \cdot \cdot\, \ 0\,, \ \xi_{1_{\,|_{\,n}}} \ )\,. \leqno  (\,A.11.20\,)
 $$
 Thus\\[0.1in]
 (\,A.11.21\,)
 \begin{eqnarray*}
& \  & \int_{B_{\,\xi_{\,1}}\,(\,{\rho_{\,\nu}})\ \cap \ \{ \ y_n \ >  \ \xi_{1_{\,|_{\,n}}}\ \}}  \ \ (\,y_n\,-\,\xi_{1_{\,|_{\,n}}}\,)^{\,\ell}  \,\cdot\,\, \left\{    \,\left(  {{\lambda_{\,1}}\over {\lambda_{\,1}^2 \ + \ \Vert\, y \,-\, \xi_{\,1}\Vert^2}}\  \right)^{\!\!n\,+\,1} \,\cdot\,   (\,y_1 \ - \ \xi_{1_{\,|_1}}\ ) \,\right\}\ dy  \\[0.3in]
& = & \int_{B_{\,\xi_{\,1}}\,(\,{\rho_{\,\nu}})\ \cap \ \{ \ y_n \ >  \ \xi_{1_{\,|_{\,n}}} \}}  \ \ (\,y_n\,-\,\xi_{1_{\,|_{\,n}}}\,)^{\,\ell}  \,\cdot\,\, \left\{    \,\left(  {{\lambda_{\,1}}\over {\lambda_{\,1}^2 \ + \ \Vert\, y \,-\, \xi_{\,1}\Vert^2}}\  \right)^{\!\!n\,+\,1} \,\cdot\,   (\,y_1\,) \,\right\}\ dy\\[0.3in]
& = & \int_{B_{\,\xi_{\,1}}\,(\,{\rho_{\,\nu}})\ \cap \ \{ \ y_n \ >  \ \xi_{1_{\,|_{\,n}}} \ \& \ y_1 \  >  \ 0 \  \}}  \ \ (\,y_n\,-\,\xi_{1_{\,|_{\,n}}}\,)^{\,\ell}  \,\cdot\,\, \left\{    \,\left(  {{\lambda_{\,1}}\over {\lambda_{\,1}^2 \ + \ \Vert\, y \,-\, \xi_{\,1}\Vert^2}}\  \right)^{\!\!n\,+\,1} \,\cdot\,   (\,y_1\,) \,\right\}\ dy \ \ + \ \ \\[0.3in]
& \ &  \ \ \  \ + \ \int_{B_{\,\xi_{\,1}}\,(\,{\rho_{\,\nu}})\ \cap \ \{ \ y_n \ >  \ \xi_{1_{\,|_{\,n}}} \ \& \  y_1 \  < \emph{}  \ 0 \  \}}  \ \ (\,y_n\,-\,\xi_{1_{\,|_{\,n}}}\,)^{\,\ell}  \,\cdot\,\, \left\{    \,\left(  {{\lambda_{\,1}}\over {\lambda_{\,1}^2 \ + \ \Vert\, y \,-\, \xi_{\,1}\Vert^2}}\  \right)^{\!\!n\,+\,1} \,\cdot\,   (\,y_1\,) \,\right\}\ dy \ \ + \ \ \\[0.3in]& = & \int_{B_{\,\xi_{\,1}}\,(\,{\rho_{\,\nu}})\ \cap \ \{ \ y_n \ >  \ \xi_{1_{\,|_{\,n}}} \  \& \  y_1 \  >  \ 0 \  \}}  \ \ (\,y_n\,-\,\xi_{1_{\,|_{\,n}}}\,)^{\,\ell}  \,\cdot\,\, \left\{    \,\left(  {{\lambda_{\,1}}\over {\lambda_{\,1}^2 \ + \ \Vert\, y \,-\, \xi_{\,1}\Vert^2}}\  \right)^{\!\!n\,+\,1} \,\cdot\,   (\,+\,y_1\,) \,\right\}\ dy \ \  \ \ \\[0.3in]
& \ &  \ \ \  \ - \ \int_{B_{\,\xi_{\,1}}\,(\,{\rho_{\,\nu}})\ \cap \ \{ \ y_n \ >  \ \xi_{1_{\,|_{\,n}}} \  \& \  y_1 \ > \ 0 \  \}}  \ \ (\,y_n\,-\,\xi_{1_{\,|_{\,n}}}\,)^{\,\ell}  \,\cdot\,\, \left\{    \,\left(  {{\lambda_{\,1}}\over {\lambda_{\,1}^2 \ + \ \Vert\, y \,-\, \xi_{\,1}\Vert^2}}\  \right)^{\!\!n\,+\,1} \,\cdot\,   (\,-\,y_1\,) \,\right\}\ dy \\[0.3in]
& = & 0 \ \ \ \ \ ( \ {\mbox{change \ \ of \ \ dummy \ \ variables \ \ }} y_1 \ \to \ \,-\,y_1 \ )\,.
 \end{eqnarray*}
 Likewise, the other five integrals in  $(\,A.11.15\,)_{\,-}$\ ,\,   $(\,A.11.16\,)_{\,+}$\ ,\,  $(\,A.11.16\,)_{\,-}$\ ,\,  $(\,A.11.17\,)_{\,+}$\ and  $\,(\,A.11.17\,)_{\,-}$\,,\, are each equal to zero.


\newpage

{\bf{\S\,A\,11\,.\,g}} \ \ {\bf Expression in vector form\,.\,} \ \  Returning to (\,A.11.1\,)\,,\, with the setting that
$$
\xi_{\,1} \ = \ (\,0\,,\ \cdot \, \cdot \, \cdot\,, \ \xi_{1_{\,|_{\,\,n}}}
\ )\,,  \ \ \ \ \ {\mbox{with}} \ \ \ \xi_{1_{\,|_{\,\,n}}} \ > \ 0
$$
[\ cf. (\,A.11.3\,)\ ]\,,\,
the above discussion  in  {\bf \S\,A\,11\,.\,b\,} to  {\bf \S\,A\,11\,.\,f\,} shows that
$$
[\  \vec{\,\bf D}\ ]_{\,(\,A.11.1\,)} \ = \ \left( \ 0\,, \ \ \cdot \cdot \cdot\,, \ 0\,, \ \ {\hat C}_+ \cdot \lambda_{\,1}^\ell \cdot\left( {{  \eta  }\over {  \lambda_{\,1} }}\  \right) \ \right) \ + \ \vec{\,{\bf E}}\ , \leqno (\,A.11.22\,)
$$
where [\ cf. (\,A.9.43\,)\ ]\\[0.1in]
(\,A.11.23\,)
\begin{eqnarray*}
\vec{\,{\bf E}}  & = & ( \ E_1\,, \ E_2\,, \ \cdot \cdot \cdot\, \ E_n\ )\,,\\[0.2in]
{\mbox{with}} \ \ \ \ \ \ \ \ E_j & = &     O \,
 \left( \ {\bar \lambda}_{\,\,\flat}^{\,\ell \ + \ \mu_{\,\vert_{\cal K}}}\ \right) \ + \     O\left(\   {\bar\lambda}_{\ \flat}^{ \  n\cdot\,(\,1\ - \ \nu\,) \ - \ o_{\,+}\,(\,1\,)} \ \right)  \ + \ O\left( \  {\bar\lambda}_{\ \flat}^{ \,n \,\gamma \ - \ \sigma}\ \right)\ , \ \ \ \ \ \ \ \ \ \\[0.1in]
 &\ & \hspace*{3.9in} j \ = \ 1\,, \ \cdot \cdot \cdot\,,\ n \,. \ \ \ \ \ \
   \end{eqnarray*}
   In the above, the number $\,\mu_{\,\vert_{\cal K}}\,$ is given by (\,A.9.44\,)\,,\, and $\,\sigma\,$ is presented in (\,1.32\,) of the main text\,,\, which is
   $$
2 \ \le \ \flat \ \le\  {1\over { { {\bar\lambda}_{\,\,\flat}   }^{\,\sigma} }} \ .
$$

As the positive $\,\xi_{1_{\,|_n}}$-\,direction is in the direction of the normal to the hypersurface $\,{\cal H}\,$ at the point $\,{\bf p}_1\,$  given by
$$
{\vec{\,{\bf n}}}_1 \ = \ {{ \ \xi_{\,1} \ - \ {\bf p_1} \  }\over { \Vert \ \xi_{\,1} \ - \ {\bf p_1} \ \Vert  }}\ . \leqno (\,A.11.24\,)
$$
 we come to the conclusion that\,,\, in the original coordinate system\,,\,
$$
[\  \vec{\,\bf D}\ ]_{\,(\,A.11.1\,)} \ = \ \  {\hat C}_+ \cdot \lambda_{\,1}^\ell \cdot\left( {{  \eta   }\over {  \lambda_{\,1} }}\  \right) \ {\vec{\,{\bf n}}}_1   \ + \ {\vec{\,{\bf E}}}_{\,(\,A.11.25\,)}\ . \leqno (\,A.11.25\,)
$$
Here $\,{\vec{\,{\bf E}}}_{\,(\,A.11.25\,)}\,$ takes on the similar form as $\,{\vec{\,{\bf E}}}\,$\, [\ in (\,A.11.23\,)\ ]\,.

\newpage

{\bf{\S\,A\,11\,.\,g\,.}} \ \ {\bf Derivative of the reduced functional with respect to $\,\xi_{\,1}$ - conclusion\,.\,} \\[0.1in]
(\,A.11.26\,)
\begin{eqnarray*}
& \ & \left(\ \lambda_{\,1}\cdot \btd_{\,\xi_{\,1}} \, \right)\,  {\bf I}_R\,\left( \ (\,\lambda_{\,1}\,, \ \xi_{\,1}\,)\,, \ \cdot \cdot \cdot\,, \ (\,\lambda_{\,\flat}\,, \ \xi_{\ \flat} \,) \ \right) \\[0.2in] & = &  {\hat C}_{\,2\,, \ 2} \,\cdot\,\sum_{l\ = \ 2}^\flat \left\{  \, \left(\ \ {{ \sqrt{\lambda_{\,1}\,\cdot\,\lambda_{\,\,l}\,} }\over {\Vert\ \xi_{\,1} \ - \ \xi_{\,\,l}\ \Vert   }}\  \right)^{\!\!n}\,\cdot\,\left(\ {{\,\xi_{\,1}\ - \ \xi_{\,l}\,}\over {\lambda_{\ \!l}}}\  \right)  \ \right\}  \,\times \left[\ 1 \  +  \  O\,\left( \,{\bar\lambda}_{\ \flat}^{\,2\,(\,1 \ - \ \nu\,) }\, \right) \ \right]\   \\[0.2in]
& \ & \ \ \ \ \ + \  \ {\hat C}_{\,2\,, \ 3}\cdot \lambda_{\,1}^\ell \cdot\left( {{  \eta   }\over {  \lambda_{\,1} }}\  \right) \ {\bf n}_{\,1}  \ + \ {\vec{\,{\bf E}}}_{\,(\,A.11.26\,)}^\tau\ ,
   \end{eqnarray*}
\begin{eqnarray*}
(\,A.11.27\,) \ \ \ \ \ \ \ \ \ \ \ \     {\vec{\,{\bf E}}}_{\,(\,A.11.26\,)}^\tau& = & ( \ E_1^{\,\tau}\,, \ E_2^{\,\tau}\,, \ \cdot \cdot \cdot\, \ E_n^{\,\tau} \ )\\[0.2in]
& \ & \hspace*{-2in}E_j^{\,\tau} \ = \       O \,
 \left( \ {\bar \lambda}_{\,\,\flat}^{\,\ell \ + \ \mu_{\,\vert_{\cal K}}}\ \right) \ + \     O\left(\   {\bar\lambda}_{\ \flat}^{ \  n\cdot\,(\,1\ - \ \nu\,) \ - \ o_{\,+}\,(\,1\,)} \ \right)  \ + \   O\left( \  {\bar\lambda}_{\ \flat}^{ \,n \,\gamma \ - \ \sigma}\ \right)
 \\[0.in]
 & \ &\hspace*{-1in} \ \ \   [\ \uparrow   \ \ {\mbox{cf. (\,A.9.44\,)}} \ ] \\[0.1in]
 &  &   + \   O\left( \  {\bar\lambda}_{\ \flat}^{ \,(\,n \,+\,1\,)\, \cdot \,\gamma \ - \   o_{\,{\bar{\lambda}}_{\ \flat}} }\ \right)  \ + \   {\cal E}_{\,(\,A.{\bf 7}.55\,)} \\[0.in]
 & \ & \hspace*{0.7in}  [\ \uparrow   \ \ {\mbox{cf. (\,A.10.28\,)}} \ ] \\[0.1in]
& \ & \hspace*{-2in}\ \ \ \ \ \ \ = \       O \,
 \left( \ {\bar \lambda}_{\,\,\flat}^{\,\ell \ + \ \mu_{\,\vert_{\cal K}}}\ \right) \ + \     O\left(\   {\bar\lambda}_{\ \flat}^{ \  n\cdot\,(\,1\ - \ \nu\,) \ - \ o_{\,+}\,(\,1\,)} \ \right)  \ + \ O\left( \  {\bar\lambda}_{\ \flat}^{ \,n \,\gamma \ - \ \sigma}\ \right)\ \  + \  \   {\cal E}_{\,(\,A.{\bf 7}.55\,)}  \ .
 \\[0.1in]
 & \ & \hspace*{0.7in}  \left[ \ \uparrow   \ \  O \left(\  {\bar\lambda}_{\ \flat}^{(\,n \ + \ 1\,)\, \cdot\, \gamma \ - \   o_{\,{\bar{\lambda}}_{\ \flat}  }}\ \right)  \ \ {\mbox{being \ \ absorbed}} \ \right]\ .
   \end{eqnarray*}
Here $\,\mu_{\,\vert_{\cal K}}\,$ is defined in (\,A.9.44\,)\,,\,
$$
{\hat C}_{\,2\,, \ 3}\ = \ C\,(\,{\bf p}_1\,)\cdot (\,n \ - \ 2\,)\cdot \ell \cdot
\int_{\R^n}  \,|\,Y_n\,|^\ell \,\left(  {{1}\over {1 \ + \ \Vert\, Y\,\Vert^2}}\  \right)^{\!\!n\,+\,1} \  dV_{\,Y} \ ,
$$

\vspace*{0.1in}

 $$
   o_{\,+}\,(\,1\,) \ = \  2\,\varepsilon \,[\ (\,\gamma \ + \ \nu\,) \ - \ 1\ ]\ ,
   $$

  $$
  {\hat C}_{\,2\,, \ 2} \ = \   \omega_n\,\cdot\,(\,n\,-\,2)^{\,2} \ . \leqno {\mbox{and}}
$$
In (\,A.11.26\,)\,,\, the bubble parameters
 $\,\left( \ \lambda_{\,1}\,, \ \cdot \cdot \cdot\,, \ \lambda_{\,\flat}\,, \ \xi_{\,1}\,, \ \cdot \cdot \cdot\,, \ \xi_{\,\flat}   \ \right)\,$ satisfy the conditions in Theorem 1.33 of the main text\,.\, Similar expressions for derivative of the reduced functional with respect to $\,\xi_{\,l}$\, for \,$\, l \ = \ 2\,, \ \cdot \cdot \cdot\,, \ \flat\,.$

      \newpage



{\large{\bf{\S\,A\,12\,.\, \ \ Estimate of the reduced functional $\,{\bf I}_{\,\cal R} \,.$}}}

\vspace*{0.2in}

Similar argument as in  {\bf \S A\,7}\, leads to
$$
{\bf I}_{\,\cal R} \ (\, W_{\ \flat} \ + \ \phi\,) \ = \  {\bf I}\,(\,W_{\,\,\flat}\,)\ + \ {\cal E}_{\,(\,A.7.55\,)}\ . \leqno (\,A.12.1\,)
$$
Here
$$
 {\bf I}\,(\,W_{\,\,\flat}\,)\ = \   {1\over 2}\,\int_{\R^n}\,\langle\,\btd\,W_{\,\,\flat}\,,\,\btd\,W_{\,\,\flat}\,\rangle\ -\ \left(\,{{n\,-\,2}\over {2n}}\,\right)\,\cdot\,\int_{\R^n}\,(\,{\tilde c}_n\!\cdot K\,)\,W_{\,\,\flat}^{{2n}\over {\,n\,-\ 2\,}} \ .\leqno (\,A.12.2\,)
$$
See (\,A.7.57\,) for the expression on $\,{\cal E}_{\,(\,A.7.55\,)}\,.\,$

\vspace*{0.3in}

{\bf{\S\,A\,12\,.a\,.}} \ \ {\bf The interaction part\,.\,}    \ \ Consider the first term in (\,A.12.2\,)\,,\, that is\,,\,\\[0.1in]
(\,A.12.3\,)
\begin{eqnarray*}
  & \ & {1\over 2}\cdot \int_{\R^n}\,\langle\,\btd\,W_{\,\,\flat}\,,\,\btd\,W_{\,\,\flat}\,\rangle\ \  \\[0.2in]
  &  = &   -\,   {1\over 2}\cdot \int_{\R^n}\,(\,\Delta\,W_{\,\,\flat}\,)\,\cdot\,W_{\,\,\flat} \ \ \\[0.2in]
  & = &  {1\over 2}\cdot n\,(\,n\,-\,2\,)\cdot \int_{\R^n}\,\left(\,V_1^{\,{{n\,+\,2}\over {\,n\,-\ 2\,}} }\ + \ \cdot\,\,\cdot\,\cdot \ + \ V_{\,\flat}^{{\,n\,+\,2\,}\over {n\,-\,2}}\,\right)\,\cdot\,(\,V_1 \ + \ \cdot\,\cdot\,\cdot\,\ + \ V_{\,\flat} \,)\ \ \\[0.2in]
  & = &  {1\over 2}\cdot n\,(\,n\,-\,2\,)\cdot \int_{\R^n}\,\left(\,V_1^{{2n}\over {\,n\,-\ 2\,}} \ + \ \cdot\,\cdot\,\cdot\,\ + \ V_{\,\flat}^{{2n}\over {\,n\,-\ 2\,}}\,\right)\\[0.1in]
    & \ &  \left[\ \uparrow  \ = \   {1\over 2}\cdot n\,(\,n\,-\,2\,)\cdot V\,(\,n)\,\cdot\,\flat\ ;  \ \ {\mbox{to \   be \   combined \   with  \   similar  \   term \   in \ \ }}  (\,A.12.12\,)\ \right]\\[0.1in]
  & \ & \ \ + \  {1\over 2}\cdot n\,(\,n\,-\,2\,)\cdot \int_{\R^n}\,V_1\,\cdot\,\left(\,V_2^{\,{{n\,+\,2}\over {\,n\,-\ 2\,}} }\ + \ \cdot\,\cdot\,\cdot\,\ + \ V_{\,\flat}^{{\,n\,+\,2\,}\over {n\,-\,2}}\,\right)  \ \  \ + \ \\[0.2in]
  & \ & \ \ \ \ + \  {1\over 2}\cdot n\,(\,n\,-\,2\,)\cdot \int_{\R^n}\,V_2\,\cdot\,\left(\,V_1^{\,{{n\,+\,2}\over {\,n\,-\ 2\,}} }\ + \ V_3^{\,{{n\,+\,2}\over {\,n\,-\ 2\,}} } \ + \ \cdot\,\cdot\,\cdot\,\ + \ V_{\,\flat}^{{\,n\,+\,2\,}\over {n\,-\,2}}\,\right)  \ \  \ + \ \\[0.1in]
  & \ &  \ \ \ \ \hspace*{1in}: \\[-0.1in]
   & \ &  \ \ \ \ \hspace*{1in}:  \\[0.1in]
  & \ & \ \   \ \ \ \ \ \ \ + \  \,{1\over 2}\cdot n\,(\,n\,-\,2\,)\cdot \int_{\R^n}\,V_{\,\flat}\,\cdot\, \left(\,V_1^{\,{{n\,+\,2}\over {\,n\,-\ 2\,}} }\ + \ \cdot\,\cdot\,\cdot\,\ + \ V_{\flat\,-\,1}^{{\,n\,+\,2\,}\over {n\,-\,2}}\,\right) \ . \ \ \\[0.2in]
\end{eqnarray*}
In the above,
$$
V\,(\,n\,) \ = \  \int_{\R^n}
{\bf V}_o^{{2n}\over {\,n\,-\ 2\,}}\  d\,Y\ \ \ \ \ {\mbox{where}} \ \ \   {\bf V}_o \,(\,Y\,) \ = \ \left(\  {1\over {\,1\ + \ \Vert\,Y\Vert^2\,}}\  \right)^{\!\! {{\,n\,-\,2\,}\over 2}} \ .
$$

\vspace*{0.35in}

{\bf{\S\,A\,12\,.b\,.} }\ \ Similar to  {\bf \S\,4\,a} in the main text\,,\,  and {\bf \S\,A\,8}\,,\, we have\\[0.1in]
(\,A.12.4\,)
\begin{eqnarray*}
  & \ & {1\over 2}\cdot n\,(\,n\,-\,2\,)\cdot\,\left[ \  \int_{B_2} V_2^{{\,n\,+\,2\,}\over {n\,-\,2}}\,\cdot\,V_1  \ + \ \cdot \cdot \cdot \ + \ \int_{B_n} V_n^{{\,n\,+\,2\,}\over {n\,-\,2}}\,\cdot\,V_1\ \right]
\\[0.2in]
& = & {1\over 2}\cdot n\,(\,n\,-\,2\,)\cdot C\,(\,n) \left\{ \ \sum_{j\ = \ 2}^n {{\lambda_{\,1}^{{n\,-\,2 }\over 2}\!\cdot\lambda_{\,j}^{{\,n\,-\ 2\,}\over 2}}\over {\ \Vert\ \xi_{\,1}\,-\,\xi_j\,\Vert^{\,n\,-\,2}\  }}    \ \right\}  \ *   \left[\ 1 \ +  \  O\,\left(\ {\bar\lambda}_{\ \flat}^{\, 2 \,(\,1 \ - \ \nu\,) }\ \right)\ \right]  \\[0.2in]
 & \ &     \ \ \ \  + \  \   {\cal E}_{\,(\,A.12.4\,) }  \ ,
\end{eqnarray*}
where\\[0.1in]
 (\,A.12.5\,)
$$
{\cal E}_{\,(\,A.12.4\,) } \ = \  O\left(\   {\bar\lambda}_{\ \flat}^{ \   {{n\,+\,2}\over 2}\,\,\cdot\,\,\gamma \ + \ {{\,n\,-\,2\,}\over 2}\,\cdot\,(\,1\ - \ \nu\,)  } \,\cdot\, {\bar\lambda}_{\ \flat}^{  \   - \ \varepsilon \,[\ (\,\gamma \ + \ \nu\,) \ - \ 1\ ] }  \ \right) \ + \ O\left( \  {\bar\lambda}_{\ \flat}^{ \,n \,\gamma \ - \ \sigma}\ \right)\ ,
$$
$$
B_j \ := \ B_{\xi_j}\,(\,\rho_{\,\nu})\ \ \ \ \ \ \  \ \mfor \  \ j \ = \ 1\,,\ \cdot \cdot \cdot\,,\ n, \leqno (\,A.12.6\,)
$$
and
$$
C\,(\,n) \ = \ \omega_n\,\cdot\,\int_0^\infty \left({1\over {1 \ + \ R^2}}\  \right)^{\!\!{{n\,+\,2}\over 2} }\,\cdot\,R^{\,n\,-\,1}\ dR\ = \ {1\over n}\,\cdot\,\omega_n \  . \leqno (\,A.12.7\,)
$$
Note that\,,\, due to symmetry\,,\, the term
$$
{{\lambda_{\,1}^{{n\,-\,2 }\over 2}\!\cdot\lambda_{\,2}^{{\,n\,-\ 2\,}\over 2}}\over {\ \Vert\ \xi_{\,1}\,-\,\xi_{\,2}\,\Vert^{\,n\,-\,2}\  }}
$$
appears again when we estimate
$$
\int_{B_1} V_1^{{\,n\,+\,2\,}\over {n\,-\,2}}\,\cdot\,V_2 \ .
$$
 Cf. (\,A.12.16\,)\,.

\newpage

{\bf \S\,A\,12\,.\,c\,. } \ \ {\bf  Estimate of}\\[0.1in]
(\,A.12.8\,)
\begin{eqnarray*}
\int_{B_1} V_2^{{\,n\,+\,2\,}\over {n\,-\,2}}\,\cdot\,V_1\ d\,y & = & \int_{B_1} \left(\ {{\lambda_{\,2}}\over { \lambda_{\,2}^2 \ + \ \Vert\ y \ - \ \xi_{\,2}\,\Vert^2 }}\ \right)^{\!{{n\,+\,2}\over 2} }\,\cdot\, \left(\ {{\lambda_{\,1}}\over { \lambda_{\,1}^2 \ + \ \Vert\ y \ - \ \xi_{\,1}\,\Vert^2 }}\ \right)^{\!{{\,n\,-\,2\,}\over 2} }  d\,y\\[0.2in]
& \le & C_1\,\cdot\,{1\over { \Vert\ \Xi_{\,\,1} \ - \ \Xi_{\,2}\,\Vert^{n\,+\,2} }}\,\cdot\,\int_0^{{ \rho_\nu}\over { \lambda_{\,1} }} {1\over { 1 \ + \ R^{\,n\,-\,2}}}\,\cdot\, R^{\,n\,-\,1}\ d\,R\\[0.2in]
& \ &\hspace*{-1in} \le \ C_1\,\cdot\,{1\over { \Vert\ \Xi_{\,\,1} \ - \ \Xi_{\,2}\,\Vert^{n\,+\,2} }}\ \cdot\  O\left(\  {1\over {\lambda_{\,1}^{2\,(\,1\ - \ \nu\,)} }} \, \right)  \ \ \ \ \ (\,{\mbox{similar \ \ to  \ \ }}  {\bf \S\ A\,8\,.\,\,b}\,)\ .
\end{eqnarray*}

 \vspace*{0.35in}

{\bf \S\,A\,12\,.\,d\,. }
\ \ {\bf Outside} \ $\displaystyle{\int_{\R^n \,\setminus\,(\,B_1\,\cup\,B_2\,)} \!\!\! V_2^{{\,n\,+\,2\,}\over {n\,-\,2}}\,\cdot\,V_1\ d\,y  }$ \ . \ Similar to   {\bf  \S\,A\,8\,.\,\,h}\,,\, in particular, (\,A.3.36\,)\,,\,  we have
$$
\int_{\R^n \,\setminus\,(\,B_1\,\cup\,B_2\,)} \!\! V_2^{{\,n\,+\,2\,}\over {n\,-\,2}}\,\cdot\,V_1\ d\,y  \ = \   {1\over {\Vert\  \Xi_{\,\,1} \ - \ \Xi_j\,\Vert^{\, n\,-\,2 } }}  \,\cdot\, O\,\left( \ {\bar\lambda}_{\ \flat}^{ \ 2 \, \cdot\,(\,1\ - \ \nu\,) } \ \right) \ . \leqno (\,A.12.9\,)
$$

\vspace*{0.2in}

With the result in {\bf \S\,A\,12\,.\,a}\,,\,  {\bf \S\,A\,12\,.\,b}\, and  {\bf \S\,A\,12\,.\,c}\,,\, we arrive at
\begin{eqnarray*}
& \ & \\
(\,A.12.10\,)  & \ &  \int_{\R^n}\,V_1^{{\,n\,+\,2\,}\over {n\,-\,2}}\,\cdot\,\left(\,V_2 \ + \ \cdot\,\cdot\,\cdot\,\ + \ V_{\,\flat}\,\right)  \ \ \\[0.2in]
& & \hspace*{-1in}= \  C\,(\,n)   \left\{ \ \sum_{k \,\not=\,1} \ {{\lambda_{\,1}^{{n\,-\,2 }\over 2}\!\cdot\lambda_k^{{\,n\,-\ 2\,}\over 2}}\over {\ \Vert\ \xi_{\,1}\,-\,\xi_k\,\Vert^{\,n\,-\,2}\  }}    \ \right\}\! \cdot\! \left[\, 1 \ + \ \emph{}O\,\left(\ {\bar\lambda}_{\ \flat}^{\, 2 \,(\,1 \ - \ \nu\,) }\ \right)  \,
 \, \right] \   + \     {\cal E}_{\,(\,A.12.4\,) } \  + \   {\cal E}_{\,(\,A.12.10\,) }\ .
\end{eqnarray*}
Here
$$
{\cal E}_{\,(\,A.12.10\,) } \ = \   O\left(\   {\bar\lambda}_{\,\,\flat}^{ \ (\,n\ - \ 2\,) \,\cdot\,\,\gamma \ + \ 2\,(\,1\ - \ \nu\,)  }   \ \right)\ ,
$$
and
$$
C\,(\,n) \ = \ {1\over n}\,\cdot\,\omega_n \  .
$$

\newpage

Note that
\begin{eqnarray*}
& \ &  (\,n\ + \ 2\,) \,\cdot\,\,\gamma \ - \ 2\,(\,1\ - \ \nu\,) \ \ge\  (\,n\ - \ 2\,) \,\cdot\,\,\gamma \ + \ 2\,(\,1\ - \ \nu\,) \\[0.2in]
\Longleftrightarrow & \ & 4\,\gamma \ \ge \ 4 \ - \ 4\,\nu \\[0.2in]
\Longleftrightarrow & \ & \gamma \ + \ \nu \ \ge \ 1\ .
\end{eqnarray*}
Similarly we estimate the other combination of indices, leading to \\[0.1in]
(\,A.12.11\,)
\begin{eqnarray*}
& \ &  {1\over 2}\cdot \int_{\R^n}\,\langle\,\btd\,W_{\,\,\flat}\,,\,\btd\,W_{\,\,\flat}\,\rangle\ \  \\[0.3in]
&  =  &   {1\over 2}\cdot n\,(\,n\,-\,2\,)\cdot V\,(\,n)\,\cdot\,\flat \\[0.2in]
  & \ & \!\!\!\!\!\!\!\! + \  {1\over 2}\cdot n\,(\,n\,-\,2\,)\cdot C\,(\,n) \sum_{j \,=\,1}^\flat \left\{ \ \sum_{k \,\not=\,j} \ {{\lambda_{\,j}^{{n\,-\,2 }\over 2}\!\cdot\lambda_k^{{\,n\,-\ 2\,}\over 2}}\over {\ \Vert\ \xi_j\,-\,\xi_k\,\Vert^{\,n\,-\,2}\  }}    \ \right\} \,\times \,\left[\ 1  \ + \ O\,\left(\ {\bar\lambda}_{\ \flat}^{\, 2 \,(\,1 \ - \ \nu\,) }\ \right)  \ \right]\\[0.15in]
  & \ & \ \ \ \ \ \ \ \ \ \ \ \ \ \ \  \ \ \ \  \ \ \ \ \  \ \ \ \   \ \ \ \  \     [\ \ \uparrow  {\mbox{to \  be \  combined \ with \ similar \  team \ in \   (\,A.12.16\,)\,}}\  ]\\[0.2in]
 & \ &      \hspace*{0.5in}\ \ \ \ \ \ \ \ \ + \  \  {\cal E}_{\,(\,A.12.11\,) }\ \ \ \ \ \ \ \  \ \ \ \  \ \ \ \ \ \ \ \ \  \ \ \ \  \ \ \ \ \ \ \ \ \  \ \ \ \  \  \left[ \ C\,(\,n) \ = \ {1\over n}\,\cdot\,\omega_n  \ \right]  \ .\\
\end{eqnarray*}
Here
\begin{eqnarray*}
| \  {\cal E}_{\,(\,A.12.11\,) } \ | & \le & \flat \cdot \left[ \  {\cal E}_{\,(\,A.12.4\,) } \  + \   {\cal E}_{\,(\,A.12.10\,) }\ \right] \\[0.3in]
 & = &  O\left(\   {\bar\lambda}_{\ \flat}^{ \   {{n\,+\,2}\over 2}\,\,\cdot\,\,\gamma \ + \ {{\,n\,-\,2\,}\over 2}\,\cdot\,(\,1\ - \ \nu\,)  } \,\cdot\, {\bar\lambda}_{\ \flat}^{  \   - \ \varepsilon \,[\ (\,\gamma \ + \ \nu\,) \ - \ 1\ ] \ - \  \sigma }  \ \right) \ + \ O\left( \  {\bar\lambda}_{\ \flat}^{ \,n \,\gamma \ - \ 2\,\sigma}\ \right) \ \\[0.3in]
  & \ &  \ \ \ \ \ \ \ \ \ \   \ \ \ \ \ \ \ \ \ \    \ \ \ \ \ \ \ \ \ \    +   \ \  O\left(\   {\bar\lambda}_{\,\,\flat}^{ \ (\,n\ - \ 2\,) \,\cdot\,\,\gamma \ + \ 2\,(\,1\ - \ \nu\,) \ - \  \sigma  }   \ \right) \ .
\end{eqnarray*}


\newpage

{\bf \S\,A\,12\,.\,e\,.}  \ \   As for the second term in\, (\,A.12.2\,)\,,\, namely\,,\,\\[0.1in]
(\,A.12.12\,)
\begin{eqnarray*}
& \ & \\
  & \ & -\ \left(\ {{n\,-\,2}\over {2n}}\,\right)\,\cdot\,\int_{\R^n}\,(\,{\tilde c}_n\!\cdot K\,)\,W_{\,\,\flat}^{{2n}\over {\,n\,-\ 2\,}} \ \ \\[0.2in]
  & = & -\ \left(\ {{n\,-\,2}\over {2n}}\,\right)\,\cdot\,\Bigg\{ \ \int_{\R^n}\,[\ (\,{\tilde c}_n\!\cdot K\,) \ - \ n\,(\,n\,-\,2\,) \ ]\,\cdot\,W_{\,\,\flat}^{{2n}\over {\,n\,-\ 2\,}} \ \  \   \\[0.2in]& \ &  \hspace*{1.5in}    + \ \int_{\R^n}\,  n\,(\,n\,-\,2\,) \,\cdot\,\left[\ W_{\,\,\flat}^{{2n}\over {\,n\,-\ 2\,}} \ - \ \left(\,V_1^{{2n}\over {\,n\,-\ 2\,}} \ + \ \cdot\,\cdot\,\cdot\,\ + \ V_{\,\flat}^{{2n}\over {\,n\,-\ 2\,}}\,\right)  \ \right] \ \ \\[0.2in]
& \ &   \hspace*{1.71in}  +   \ \int_{\R^n}\,  n\,(\,n\,-\,2\,) \,\cdot\,\left(\,V_1^{{2n}\over {\,n\,-\ 2\,}} \ + \ \cdot\,\cdot\,\cdot\,\ + \ V_{\,\flat}^{{2n}\over {\,n\,-\ 2\,}}\,\right) \ \  \  \Bigg\} \\[0.2in]
  & = & -\ \left(\,{{n\,-\,2}\over {2n}}\,\right)\,\cdot\, \left\{ \ \int_{\R^n}\,[\ (\,{\tilde c}_n\!\cdot K\,) \ - \ n\,(\,n\,-\,2\,) \ ]\,\cdot\,\left(\,V_1^{{2n}\over {\,n\,-\ 2\,}} \ + \ \cdot\,\cdot\,\cdot\,\ + \ V_{\,\flat}^{{2n}\over {\,n\,-\ 2\,}}\,\right)  \right.  \\[0.1in]
  & \ & \hspace*{5in}\uparrow \ \cdot \cdot \cdot \ (\,{\bf A}\,)_{\,(\,A.12.12\,)}\\[0.2in]
  & \ &  \ \ \ \ \ \ \ \ \ \ \ \ + \     \int_{\R^n}\,  [\ (\,{\tilde c}_n\!\cdot K\,) \ - \ n\,(\,n\,-\,2\,) \ ]\,\cdot\,\left[\ W_{\,\,\flat}^{{2n}\over {\,n\,-\ 2\,}} \ - \ \left(\,V_1^{{2n}\over {\,n\,-\ 2\,}} \ + \ \cdot\,\cdot\,\cdot\,\ + \ V_{\,\flat}^{{2n}\over {\,n\,-\ 2\,}}\,\right)  \ \right] \ \   \\[0.1in]
  & \ & \hspace*{5in}\uparrow \ \cdot \cdot \cdot \ (\,{\bf B}\,)_{\,(\,A.12.12\,)} \  \\[0.1in] & \ &  \hspace*{1.5in}    + \ \int_{\R^n}\,  n\,(\,n\,-\,2\,) \,\cdot\,\left[\ W_{\,\,\flat}^{{2n}\over {\,n\,-\ 2\,}} \ - \ \left(\,V_1^{{2n}\over {\,n\,-\ 2\,}} \ + \ \cdot\,\cdot\,\cdot\,\ + \ V_{\,\flat}^{{2n}\over {\,n\,-\ 2\,}}\,\right)  \ \right] \ \ \\[0.1in]
  & \ & \hspace*{5in}\uparrow \ \cdot \cdot \cdot \ (\,{\bf C}\,)_{\,(\,A.12.12\,)} \  \\[0.1in]
& \ &   \hspace*{1.71in}  +   \ \int_{\R^n}\,  n\,(\,n\,-\,2\,) \,\cdot\,\left(\,V_1^{{2n}\over {\,n\,-\ 2\,}} \ + \ \cdot\,\cdot\,\cdot\,\ + \ V_{\,\flat}^{{2n}\over {\,n\,-\ 2\,}}\,\right) \  \Bigg\}
\\[0.1in]
    & \ & \hspace*{2.2in}\Bigg[\ \uparrow  \ = \   -\ {1\over 2}\,\cdot\, \left(\,{{n\,-\,2}\over { n}}\,\right)\,\cdot\,n\,(\,n\,-\,2\,)\cdot V\,(\,n)\,\cdot\,\flat\ ;  \\[0.1in]
    & \ & \ \ \ \ \ \ \ \ \ \ \ \ \ \ \  \ \ \ \ \ \ \ \ \ \  \ \ \ \ \ \  \ \  \ \ \ \ \ \ \ \  \ \ \ \ \   {\mbox{to \ \ be \ \ combined \ \ with  \ \ similar  \ \ term \ \ in \ \ }}  (\,A.12.3\,)\ \Bigg]\ .
\end{eqnarray*}


\newpage

{\bf \S\,A\,12\,.\,f\,.} \ \  {\bf  Estimate of (\,${\bf C}\,)_{\,A.12.12}$}\ \,.\, \ \ As in {\bf A\,8\,.\,a}\,,\,\\[0.1in]
(\,A.12.13\,)
$$
 {{\ V_2\,(\,y\,) \ + \ \cdot\,\cdot\,\cdot\,\ + \ V_{\,\flat}\,(\,y\,)\   }\over {V_1\,(\,y\,)  }} \ = \ O\,\left( \ {\bar\lambda}_{\ \flat}^{\ (\,n\,-\,2\,)\,\cdot\, (\,\gamma\,+\,\nu\,) \ - \ 1}\ \right)
\ \ \  \ \ \mfor \ \ \  y\,\in\,B_{\,\xi_{\,1}}(\,\rho_{\,\nu})\ .
$$

With this\,,\, we proceed as in {\bf A\,8\,.\,a}\,:\\[0.1in]
(\,A.12.14\,)
\begin{eqnarray*}
& \ &   \left[\  V_1\ + \  V_2  \ + \ \cdot\,\cdot\,\cdot\,\ + \ V_{\,\flat}  \ \right]^{{2\,n}\over {n\,-\,2}}  \ = \ V_1^{{2\,n}\over {n\,-\,2}}\,\cdot\,\bigg[\  1 \ + \ \bigg( \ {{  V_2  \ + \ \cdot\,\cdot\,\cdot\,\ + \ V_{\,\flat}}\over {V_1}} \ \bigg)  \, \bigg]^{\,{{2\,n}\over {n\,-\,2}}}\\[0.15in]
& = & V_1^{{2\,n}\over {n\,-\,2}}\,\cdot\,\left[\  1 \ + \ {{2\,n}\over {n\,-\,2}}\,\cdot\,\left(\ \ {{  V_2  \ + \ \cdot\,\cdot\,\cdot\,\ + \ V_{\,\flat}}\over {V_1}}  \ \right)\ + \ O \left(\ \  {{  V_2  \ + \ \cdot\,\,\cdot\,\cdot \ + \ V_{\,\flat}}\over {V_1}}\  \right)^{{\,n\,+\,2\,}\over {n\,-\,2}}\  \right]\\[0.15in]
& = & V_1^{{2\,n}\over {n\,-\,2}} \ + \  {{2\,n}\over {n\,-\,2}} \,\cdot\,V_1^{{\,n\,+\,2\,}\over {n\,-\,2}}\,\cdot\, \left[\   V_2  \ + \ \cdot\,\cdot\,\cdot\,\ + \ V_{\,\flat}  \ \right] \ + \ O \left(\,  \left[ \, V_2  \ + \ \cdot\,\cdot\,\cdot\,\ + \ V_{\,\flat} \, \right]^{  {{2\,n}\over {n\,-\,2}}}\  \right) \ . \\
\end{eqnarray*}

It follows from (\,A.12.14\,)  that\\[0.1in]
(\,A.12.15\,)
\begin{eqnarray*}
& \ &  n\,(\,n\,-\,2)\,\cdot\,\int_{B_{\,\xi_{\,1}} (\,{\rho_{\,\nu}}\,)}    \left[\ W_{\,\,\flat}^{{2n}\over {\,n\,-\ 2\,}} \ - \ \left(\,V_1^{{2n}\over {\,n\,-\ 2\,}} \ + \ \cdot\,\cdot\,\cdot\,\ + \ V_{\,\flat}^{{2n}\over {\,n\,-\ 2\,}}\,\right)  \ \right] \ \ \\\\[0.2in]
& = & n\,(\,n\,-\,2)\,\cdot\,\int_{B_{\,\xi_{\,1}} (\,{\rho_{\,\nu}}\,)} \left\{ \  {{2\,n}\over {n\,-\,2}}  \,\cdot\, \left[\   V_2  \ + \ \cdot\,\cdot\,\cdot\,\ + \ V_{\,\flat}  \ \right] \   \right\}\,\cdot\,V_1^{{n\,-\,2}\over {n\,-\,2}}  \ + \ {\cal E}_{\,(\,A.12.15\,)}\,,
\end{eqnarray*}
where (\ similar to {\bf \S\,A\,.8\,.\,d}\ )
$$
| \ {\cal E}_{\,(\,A.12.15\,)} \ | \ \le \  C\,\cdot\, \int_{B_{\,\xi_{\,1}} (\,{\rho_{\,\nu}}\,)}\left(  \, V_2  \ + \ \cdot\,\cdot\,\cdot\,\ + \ V_{\,\flat} \, \right)^{  {{2\,n}\over {n\,-\,2}}}   \ = \    O\left(\   {\bar\lambda}_{\ \flat}^{ \  n\,\cdot\,\,\gamma \ + \ n\,\cdot\,(\,1\ - \ \nu \,)  }\ \right)   \ .
$$
Argue as in {\bf \S\,A\,8.\,g}\,,\, we arrive at
\\[0.1in]
(\,A.12.16\,)
\begin{eqnarray*}
& \ & -\ \left(\,{{n\,-\,2}\over {2n}}\,\right)\,\cdot\, n\,(\,n \ - \ 2\,) \cdot \int_{\R^n}\,  \left[\ W_{\,\,\flat}^{{2n}\over {\,n\,-\ 2\,}} \ - \ \left(\,V_1^{{2n}\over {\,n\,-\ 2\,}} \ + \ \cdot\,\cdot\,\cdot\,\ + \ V_{\,\flat}^{{2n}\over {\,n\,-\ 2\,}}\,\right)  \ \right] \ \   \\[0.2in]
& = &  -\,n\,(\,n \ - \ 2\,) \cdot  \sum_{j \,=\,1}^\flat\  \left\{ \ \sum_{k \,\not=\,j} \  \int_{B_k} V_k^{{n\,+\,2}\over {n \,-\,2}} \cdot V_j\ \  \ \right\} \ + \\[0.2in]
& \ & \ \ \ \ \ \ \ \ \ + \  \int_{\R^n\,\setminus\ \cup \,B_j}  \left[\ W_{\,\,\flat}^{{2n}\over {\,n\,-\ 2\,}} \ - \ \left(\,V_1^{{2n}\over {\,n\,-\ 2\,}} \ + \ \cdot\,\cdot\,\cdot\,\ + \ V_{\,\flat}^{{2n}\over {\,n\,-\ 2\,}}\,\right)  \ \right] \ + \  \\[0.2in]
& \ & \hspace*{3in}\ \ \ \ \ \ \ \ \ + \ \  O\left(\   {\bar\lambda}_{\ \flat}^{ \  n\,\cdot\,\,\gamma \ + \ n\,\cdot\,(\,1\ - \ \nu \,) \ - \ \sigma }\ \right) \\[0.2in]
& = & -\, n\,(\,n\,-\,2\,)\cdot C\,(\,n) \,\sum_{j \,=\,1}^\flat \left\{ \ \sum_{k \,\not=\,j} \ {{\lambda_{\,j}^{{n\,-\,2 }\over 2}\!\cdot\lambda_k^{{\,n\,-\ 2\,}\over 2}}\over {\ \Vert\ \xi_j\,-\,\xi_k\,\Vert^{\,n\,-\,2}\  }}    \ \right\}   \,\times \, \left[\ 1  \ + \ O\,\left(\ {\bar\lambda}_{\ \flat}^{\, 2 \,(\,1 \ - \ \nu\,) }\ \right) \ \right] \\[0.1in]
  & \ & \ \ \ \ \ \ \ \ \ \ \ \ \ \ \  \ \ \ \ \ \ \ \ \ \ \ \ \ \ \ \ \ \ \    [\ \  \uparrow \    {\mbox{to \   be \  combined \ with \   similar \  team  \   in (\,A.12.11\,)}}\ ]\\[0.2in]
 & \ & \ \ \ \ \ \ \ \ \ \ \ + \  \ {\cal E}_{\,(\,A.12.11\,)}\ .
\end{eqnarray*}
Note that
$$
O\left(\   {\bar\lambda}_{\ \flat}^{ \  n\,\cdot\,\,\gamma \ + \ n\,\cdot\,(\,1\ - \ \nu \,) \ - \ \sigma }\ \right) \ \le \  O\left(\   {\bar\lambda}_{\ \flat}^{ \   {{n\,+\,2}\over 2}\,\,\cdot\,\,\gamma \ + \ {{\,n\,-\,2\,}\over 2}\,\cdot\,(\,1\ - \ \nu\,)  } \,\cdot\, {\bar\lambda}_{\ \flat}^{  \   - \ \varepsilon \,[\ (\,\gamma \ + \ \nu\,) \ - \ 1\ ] \ - \  \sigma }  \ \right)\ ,
$$
$$
\gamma \ + \ \nu \ > \ 1 \ \ \Longrightarrow \ \ \gamma \ >  \ 1 \ - \ \nu\ . \leqno \mbox{as}
$$


\vspace*{0.3in}

{\bf \S\,A\,12\,.\,g\,.} \ \  {\bf Estimate of }
$$
 \int_{\R^n\,\setminus\ \cup \,B_j}  \left[\ W_{\,\,\flat}^{{2n}\over {\,n\,-\ 2\,}} \ - \ \left(\,V_1^{{2n}\over {\,n\,-\ 2\,}} \ + \ \cdot\,\cdot\,\cdot\,\ + \ V_{\,\flat}^{{2n}\over {\,n\,-\ 2\,}}\,\right)  \ \right] \ \ \ . \leqno (\,A.12.17\,)
 $$
 We modify the argument  in {\bf \S\,A\,8\,\,h}\,.\, Observe as in (\,A.8.35\,)  that
 \begin{eqnarray*}
(\,A.12.18\,) \ \ \ \ W_{\,\,\flat}^{{2n}\over {\,n\,-\ 2\,}}  & = & W_{\,\,\flat} \cdot W_{\,\,\flat}^{\,{{n\,+\,2}\over {\,n\,-\ 2\,}} }\ = \ W_{\,\,\flat} \cdot W_{\,\,\flat}^{\,{1\over 2} \cdot {\,{{n\,+\,2}\over {\,n\,-\ 2\,}} }} \cdot W_{\,\,\flat}^{\,{1\over 2} \cdot {\,{{n\,+\,2}\over {\,n\,-\ 2\,}} }}\\[0.2in]
& \le & \bigg\{ \ V_1 \ + \ \cdot \cdot \cdot \ + \ V_{\flat}\ \bigg\} \,* \\[0.2in]
& \ & \ \ \ \ *\,\Bigg\{ \,
 \left[\  V_1^{\,{{n\,+\,2}\over {\,n\,-\ 2\,}} }\ + \  V_2^{\,{{n\,+\,2}\over {\,n\,-\ 2\,}} }\ +\,\cdot\,\cdot\,\cdot\, + \ V_{\,\flat}^{\,{{n\,+\,2}\over {\,n\,-\ 2\,}} }\ \right] \\[0.2in]
  & \ & \ \ \ \ \ + \ V_1^{\,{1\over 2}\,\cdot\,{{n\,+\,2}\over {\,n\,-\ 2\,}}}\,\cdot\,\left[\    V_2^{\,{1\over 2}\,\cdot\,{{n\,+\,2}\over {\,n\,-\ 2\,}}}  \ + \,\cdot\,\cdot\,\cdot\, + \ V_{\,\flat}^{\,{1\over 2}\,\cdot\,{{n\,+\,2}\over {\,n\,-\ 2\,}}}  \ \right]\\[0.2in]
& \  &   \ \ \ \ \ \ \ \ \ \ \ \  + \  V_2^{\,{1\over 2}\,\cdot\,{{n\,+\,2}\over {\,n\,-\ 2\,}}}\,\cdot\,\left[\    V_1^{\,{1\over 2}\,\cdot\,{{n\,+\,2}\over {\,n\,-\ 2\,}}} \ + \,\cdot\,\cdot\,\cdot\, + \ V_{\,\flat}^{\,{1\over 2}\,\cdot\, {{n\,+\,2}\over {\,n\,-\ 2\,}}}  \ \right]\ + \ \cdot \cdot \cdot\ + \\[0.2in]
& \  &  \ \ \ \ \ \ \ \  \ \ \ \ \ \ \ \ \ \ \ \ + \  V_{\,\flat}^{\,{1\over 2}\,\cdot\,{{n\,+\,2}\over {\,n\,-\ 2\,}}}\,\cdot\,\left[\    V_1^{\,{1\over 2}\,\cdot\,{{n\,+\,2}\over {\,n\,-\ 2\,}}}  \ + \,\cdot\,\cdot\,\cdot\, + \ V_{\flat\,-\,1}^{\,{1\over 2}\,\cdot\,{{n\,+\,2}\over {\,n\,-\ 2\,}}} \ \right]\ \Bigg\}\ . \ \ \ \ \
 \end{eqnarray*}

\newpage

 Thus\\[0.1in]
(\,A.12.19\,)
  \begin{eqnarray*}
& \ &  W_{\,\,\flat}^{{2n}\over {\,n\,-\ 2\,}}  \ - \ \left(\,V_1^{{2n}\over {\,n\,-\ 2\,}} \ + \ \cdot\,\cdot\,\cdot\,\ + \ V_{\,\flat}^{{2n}\over {\,n\,-\ 2\,}}\,\right) \\[0.2in]
& \le & V_1 \cdot
 \left[\  V_2^{\,{{n\,+\,2}\over {\,n\,-\ 2\,}} }\ + \  V_3^{\,{{n\,+\,2}\over {\,n\,-\ 2\,}} }\ +\,\cdot\,\cdot\,\cdot\, + \ V_{\,\flat}^{\,{{n\,+\,2}\over {\,n\,-\ 2\,}} }\ \right] \ + \ \\[0.2in]
  & \ & \ \ \ \ \ + \ V_1 \cdot V_1^{\,{1\over 2}\,\cdot\,{{n\,+\,2}\over {\,n\,-\ 2\,}}}\,\cdot\,\left[\    V_2^{\,{1\over 2}\,\cdot\,{{n\,+\,2}\over {\,n\,-\ 2\,}}}  \ + \,\cdot\,\cdot\,\cdot\, + \ V_{\,\flat}^{\,{1\over 2}\,\cdot\,{{n\,+\,2}\over {\,n\,-\ 2\,}}}  \ \right]\ + \ \\[0.2in]
& \  &   \ \ \ \ \ \ \ \ \ \  + \ V_1 \cdot V_2^{\,{1\over 2}\,\cdot\,{{n\,+\,2}\over {\,n\,-\ 2\,}}}\,\cdot\,\left[\    V_1^{\,{1\over 2}\,\cdot\,{{n\,+\,2}\over {\,n\,-\ 2\,}}} \ + \,\cdot\,\cdot\,\cdot\, + \ V_{\,\flat}^{\,{1\over 2}\,\cdot\, {{n\,+\,2}\over {\,n\,-\ 2\,}}}  \ \right]\ + \\[0.1in]
& \  &  \hspace*{2.6in} \ \ \ \ \ \  : \  \\[-0.01in]
& \  &  \hspace*{2.6in} \ \ \ \ \ \  : \  \\[0.01in]
& \  &  \ \ \ \ \ \ \ \  \ \ \ \ \ \ \ \  + \ V_1 \cdot V_{\,\flat}^{\,{1\over 2}\,\cdot\,{{n\,+\,2}\over {\,n\,-\ 2\,}}}\,\cdot\,\left[\    V_1^{\,{1\over 2}\,\cdot\,{{n\,+\,2}\over {\,n\,-\ 2\,}}}  \ + \,\cdot\,\cdot\,\cdot\, + \ V_{\flat\,-\,1}^{\,{1\over 2}\,\cdot\,{{n\,+\,2}\over {\,n\,-\ 2\,}}} \ \right]\ + \\[0.01in]
&\ & \ \ \ \  \ \ \ \  \ \ \ :\\[0.01in]
&\ & \ \ \ \  \ \ \ \  \ \ \ :\\[0.1in]
& \ & \ \ \ \ + \ V_{\,\flat} \cdot
 \left[\  V_1^{\,{{n\,+\,2}\over {\,n\,-\ 2\,}} }\ + \  V_2^{\,{{n\,+\,2}\over {\,n\,-\ 2\,}} }\ +\,\cdot\,\cdot\,\cdot\, + \ V_{\flat\,-\,1}^{\,{{n\,+\,2}\over {\,n\,-\ 2\,}} }\ \right]\ + \  \\[0.2in]
  & \ & \ \ \ \ \ + \ V_{\,\flat} \cdot V_1^{\,{1\over 2}\,\cdot\,{{n\,+\,2}\over {\,n\,-\ 2\,}}}\,\cdot\,\left[\    V_1^{\,{1\over 2}\,\cdot\,{{n\,+\,2}\over {\,n\,-\ 2\,}}}  \ + \,\cdot\,\cdot\,\cdot\, + \ V_{\flat\,-\,1}^{\,{1\over 2}\,\cdot\,{{n\,+\,2}\over {\,n\,-\ 2\,}}}  \ \right]\ + \ \\[0.2in]
& \  &   \ \ \ \ \ \ \ \ \ \ \ \ \ + \ V_{\,\flat} \cdot V_2^{\,{1\over 2}\,\cdot\,{{n\,+\,2}\over {\,n\,-\ 2\,}}}\,\cdot\,\left[\    V_1^{\,{1\over 2}\,\cdot\,{{n\,+\,2}\over {\,n\,-\ 2\,}}} \ + \,\cdot\,\cdot\,\cdot\, + \ V_{\flat\,-\,1}^{\,{1\over 2}\,\cdot\, {{n\,+\,2}\over {\,n\,-\ 2\,}}}  \ \right]\ + \ \\[0.1in]
& \  &  \hspace*{2.6in} \ \ \ \ \ \  : \ \\[-0.01in]
& \  &  \hspace*{2.6in} \ \ \ \ \ \  : \ \\[0.01in]
& \  &  \ \ \ \ \ \ \ \  \ \ \ \ \ \ \ \ \ \ \  \ \ + \ V_{\,\flat} \cdot V_{\,\flat}^{\,{1\over 2}\,\cdot\,{{n\,+\,2}\over {\,n\,-\ 2\,}}}\,\cdot\,\left[\    V_1^{\,{1\over 2}\,\cdot\,{{n\,+\,2}\over {\,n\,-\ 2\,}}}  \ + \,\cdot\,\cdot\,\cdot\, + \ V_{\flat\,-\,1}^{\,{1\over 2}\,\cdot\,{{n\,+\,2}\over {\,n\,-\ 2\,}}} \ \right]\ .
 \end{eqnarray*}

 \newpage

 Proceeding as in {\bf \S\,A\,8\,.\,h}\, and {\bf \S\,A\,8\,.\,i}\,,\, we obtain\\[0.1in]
(\,A.12.20\,)
 \begin{eqnarray*}
& \ &  \int_{\R^n\,\setminus\ \cup \,B_j}  \left[\ W_{\,\,\flat}^{{2n}\over {\,n\,-\ 2\,}} \ - \ \left(\,V_1^{{2n}\over {\,n\,-\ 2\,}} \ + \ \cdot\,\cdot\,\cdot\,\ + \ V_{\,\flat}^{{2n}\over {\,n\,-\ 2\,}}\,\right)  \ \right] \ \ \\[0.2in]
 & = & \flat \cdot
      O\left(\   {\bar\lambda}_{\ \flat}^{ \   {{n\,+\,2}\over 2}\,\,\cdot\,\,\gamma \ + \ {{\,n\,-\,2\,}\over 2}\,\cdot\,(\,1\ - \ \nu\,)  }  \ \right)   \\[0.2in]
  & = &         O\left(\   {\bar\lambda}_{\ \flat}^{ \   {{n\,+\,2}\over 2}\,\,\cdot\,\,\gamma \ + \ {{\,n\,-\,2\,}\over 2}\,\cdot\,(\,1\ - \ \nu\,) \ - \ \sigma } \ \right)   \ .
\end{eqnarray*}

\vspace*{0.35in}

{\bf \S\,A\,12\,.\,h\,.} \ \   {\bf Estimate of\,} $(\,{\bf B}\,)_{\,A.12.12}$\,. \ \ This should be of at most equal order to $(\,{\bf C}\,)_{\,A.12.12}$\,,\, because of the presence of
$$
  [\ (\,{\tilde c}_n\!\cdot K\,) \ - \ n\,(\,n\,-\,2\,) \ ]\ .
$$

\vspace*{0.5in}

{\bf \S\,A\,12\,.\,i\,.} \ \   {\bf Estimate of\,} $(\,{\bf A}\,)_{\,A.12.12}$\,.\,  Namely\,,\,
\begin{eqnarray*}
  & \  &    -\ \left(\,{{n\,-\,2}\over {2n}}\,\right)\,\cdot  \int_{\R^n}\,[\ (\,{\tilde c}_n\!\cdot K\,) \ - \ n\,(\,n\,-\,2\,) \ ]\,\cdot\,\left(\,V_1^{{2n}\over {\,n\,-\ 2\,}} \ + \ \cdot\,\cdot\,\cdot\,\ + \ V_{\,\flat}^{{2n}\over {\,n\,-\ 2\,}}\,\right) \ . \\
\end{eqnarray*}
Consider\\[0.1in]
(\,A.12.21\,)
\begin{eqnarray*}
 & \ & -\ \left(\,{{n\,-\,2}\over {2n}}\,\right)\,\cdot  \int_{\R^n}\,[\ (\,{\tilde c}_n\!\cdot K\,) \ - \ n\,(\,n\,-\,2\,) \ ]\,\cdot\, V_1^{{2n}\over {\,n\,-\ 2\,}}\  \  \\[0.2in]
 & = &  -\ \left(\,{{n\,-\,2}\over {2n}}\,\right)\,\cdot  \int_{\R^n}\,[\ (\,{\tilde c}_n\!\cdot K\,)\,(\,y\,) \ - \ n\,(\,n\,-\,2\,) \ ]\,\cdot\,\left(\ {{\lambda_{\,1}}\over  { \lambda_{\,1}^2 \ + \ \Vert\  y \ - \ \xi_{\,1}\,\Vert^2}}\  \right)^n\ d\,y\\[0.2in]
 & = & +\ \left(\,{{n\,-\,2}\over {2n}}\,\right)\,\cdot\int_{B_{\,\xi_{\,1}} (\,{\rho_{\,\nu}}\,)}\,[\,n\,(\,n\,-\,2) \ - \ (\,{ c}_n\!\cdot K\,) \,(\,y\,) \,]\cdot   \left(\ {{\lambda_{\,1}}\over  { \lambda_{\,1}^2 \ + \ \Vert\  y \ - \ \xi_{\,1}\,\Vert^2}}\  \right)^n\ d\,y \ + \\[0.1in]
 & \ & \hspace*{2in} (\ \uparrow \ \ {\mbox{see \ \ {\bf \S\,\,} \ \ of  \ \ the \ \ main \ \ text\,,\, and  \ \  {\bf \S\,A\,9 }}}\ \,)\\[0.1in]
 & \  &  +\ \left(\,{{n\,-\,2}\over {2n}}\,\right)\,\cdot\int_{\R^n \,\setminus B_{\,\xi_{\,1}} (\,{\rho_{\,\nu}}\,)}\,[\,n\,(\,n\,-\,2) \ - \ (\,{ c}_n\!\cdot K\,) \,(\,y\,) \,]\cdot   \left(\ {{\lambda_{\,1}}\over  { \lambda_{\,1}^2 \ + \ \Vert\  y \ - \ \xi_{\,1}\,\Vert^2}}\  \right)^n\ d\,y \ .
\end{eqnarray*}

Argue as in {\bf \S\,A\,9}\,,\, and obtain \\[0.1in]
(\,A.12.22\,)
\begin{eqnarray*}
 & \ &  +\ \left(\,{{n\,-\,2}\over {2n}}\,\right)\,\cdot\int_{  B_{\,\xi_{\,1}} (\,{\rho_{\,\nu}}\,) }   \,[\,n\,(\,n\,-\,2) \ - \ (\,{ c}_n\!\cdot K\,) \,(\,y\,) \,]\cdot   \left(\ {{\lambda_{\,1}}\over  { \lambda_{\,1}^2 \ + \ \Vert\  y \ - \ \xi_{\,1}\,\Vert^2}}\  \right)^n\ d\,y \\[0.2in]
  & = &    {\hat C}_{\,0\,, \ 1}  \cdot C\,(\,{\bf p}_{\,1}\,)\,\cdot\,\lambda^{\,\ell}_1   \ + \      C_{2_{\,+}} \,\cdot\,C\,(\,{\bf p}_{\,1}\,)\,\cdot\, \lambda^{\,\ell\,-\,2}_1\,\cdot\, \eta^2 \ \  + \ \ {\cal E}_{\,(\,A.12.22\,)}   \\[0.1in]
 & \ & \hspace*{2in} (\ \uparrow \ \ {\mbox{see \ \ the \ \ notes \ \ in  \  \ {\bf \S\,A\,12\,.\,j}}}\ \,) \ .\\
\end{eqnarray*}
Here
$$
| \ {\cal E}_{\,(\,A.12.22\,)} | \ \le \ \flat\,\cdot\,  {\cal E}_{(\,A.9.42\,)}  \ , \leqno (\,A.12.23\,)
$$
\begin{eqnarray*}
 \eta & = & {\mbox{Dist.}}\ (\,\xi_j\,, \ {\cal H}\,)  \ = \  O\,\left( \,{\bar{\lambda}}_{\ \flat}^{\,1 \ + \ \kappa}\,\right)\ ,\\[0.3in]
C_2 & = & {1\over n} \cdot \omega_n\,\\[0.3in]
{\hat C}_{\,0\,, \ 1} & = & \left(\ {{n\,-\,2}\over {2n}}\,\right)\,\cdot\,\int_{\R^n} \left(\ {1\over { 1 \ + \ \Vert\ Y\,\Vert^2 }}\ \right)^{\!\!n} \,d\,Y\ , \\[0.3in]
{\mbox{and}} \ \ \ \ \ \ \ \  C_{2_{\,+}} & = & \left(\ {{n\,-\,2}\over {2n}}\,\right)\,\cdot\, {{\ell\,\cdot\, (\,\ell\ - \ 1\,)}\over 2}\,\cdot\,\int_{\R^n} |\,Y_n\,|^{\,\ell\ - \ 2}\,\cdot\, \left(\ {1\over { 1 \ + \ \Vert\ Y\,\Vert^2 }}\ \right)^{\!\!n} \,d\,Y\ . \ \ \ \ \ \ \ \ \ \ \ \ \ \ \ \ \ \ \ \ \ \ \ \ \ \
\end{eqnarray*}


\newpage

{\bf \S\,A\,12\,.\,j\,.} \ \   {\bf Conclusion\,.} \ \ It follows that (\,A.12.1\,)\,,\,
(\,A.12.11\,)\, and \, (\,A.12.22\,) that\\[0.1in]
(\,A.12.24\,)
\begin{eqnarray*}
 & \ & {\bf I}_{\,{\cal R}} \left( \ (\,\lambda_{\,1}\,, \ \xi_{\,1}\,)\,, \ \cdot \cdot \cdot\,, \ (\,\lambda_{\,\flat}\,, \ \xi_{\ \flat} \,) \ \right) \\[0.2in]
 & = &   {1\over 2}\,\int_{\R^n}\,\langle\,\btd\,W_{\,\,\flat}\,,\,\btd\,W_{\,\,\flat}\,\rangle\ -\ \left(\,{{n\,-\,2}\over {2n}}\,\right)\,\cdot\,\int_{\R^n}\,(\,c_n\cdot K\,)\,W_{\,\,\flat}^{{2n}\over {\,n\,-\ 2\,}}\\[0.2in]
 & = & (\,n\,-\,2\,)\cdot V\,(\,n)\,\cdot\,\flat \ + \ \\[0.2in]
 & \ &  \ +\,\sum_{j \,=\,1}^\flat\  \Bigg(\  {\hat C}_{\,0\,, \ 1}  \cdot C\,(\,{\bf p}_1\,)\cdot \lambda_{\,j}^\ell  \  - \ {\hat C}_{\,0\,, \ 2} \cdot\, \left\{ \ \sum_{k \,\not=\,j} \ {{\lambda_{\,j}^{{n\,-\,2 }\over 2}\!\cdot\lambda_k^{{\,n\,-\ 2\,}\over 2}}\over {\ \Vert\ \xi_j\,-\,\xi_k\,\Vert^{\,n\,-\,2}\  }}    \ \right\} \,*\\[0.2in]
  & \ & \hspace*{3.8in} * \ \left[\ 1  \ + \ O\,\left(\ {\bar\lambda}_{\ \flat}^{\, 2 \,(\,1 \ - \ \nu\,) }\ \right) \ \right] \ \Bigg)\\[0.2in]
  & \ & \ \  \ \ \ \ \  + \  \,{\hat C}_{\,0\,, \ 3} \cdot \sum_{j\ = \ 1}^{\flat}  C\,(\,{\bf p}_j\,)\cdot \lambda_{\,j}^{\,\ell\ -\,2}\,\cdot\, [\ {\mbox{Dist.}}\ (\,\xi_j\,, \ {\cal H}\,) \ ]^{\,2} \ +  \\[0.2in]
   & \ &  \ \ \ \   \ \ \ \   \ \ \ \   \ \ \ \ + \ \ {\cal E}_{\,(\,A.7.55\,)}     \ \   + \ \ {\cal E}_{\,(\,A.12.11\,)}     \ \     + \ {\cal E}_{\,(\,A.12.22\,)}\ .
\end{eqnarray*}
Here $\,\left( \ (\,\lambda_{\,1}\,, \ \xi_{\,1}\,)\,, \ \cdot \cdot \cdot\,, \ (\,\lambda_{\,\flat}\,, \ \xi_{\ \flat} \,) \ \right)\,$ satisfy the conditions in Proposition A.6.8\,.
Note that
$$
 {1\over 2}\cdot n\,(\,n\,-\,2\,)\cdot C_2 \ = \  {1\over 2}\cdot n\,(\,n\,-\,2\,)\cdot {1\over n} \cdot \omega_n \ = \  \left( \, {{\,n \ - \ 2\,}\over 2} \, \right) \cdot \omega_n \ .
$$
In (\,A.12.24\,)\,,\, the bubble parameters
 $\,\left( \ \lambda_{\,1}\,, \ \cdot \cdot \cdot\,, \ \lambda_{\,\flat}\,, \ \xi_{\,1}\,, \ \cdot \cdot \cdot\,, \ \xi_{\,\flat}   \ \right)\,$ satisfy the conditions in Theorem 1.33 of the main text\,.\,

\newpage

{\large{\bf \S\,A\,13\,.} } {\large
{\bf Formal derivatives and relations on the leading constants}}\,. \\[0.2in]
Consider the main information in (\,1.34\,)\,,\,  (\,1.35\,)\,,\, and  (\,1.36\,)\, of the main text\,,\,  and denote them by [ \,$l \ = \ 1 \ , \cdot \cdot \cdot\,, \ \flat$\  ]\\[0.1in]
(\,A.13.1\,)
\begin{eqnarray*}
& \ &  {\bf F} \,(\,\lambda_{\,1}\,, \ \cdot \cdot \cdot\,, \ \lambda_{\ \flat}\,; \ \xi_{\,1}\,, \ \cdot \cdot \cdot\,, \ \xi_{\ \flat}\,)\\[0.2in] & = &  \!\!(\,n\,-\,2\,)\cdot V\,(\,n\,) \cdot \flat  \, +   \,   {\hat C}_{\,0\,,\ 1} \cdot \sum_{j \ =\,1}^\flat C\,(\,{\bf p}_j\,) \cdot  \lambda_{\,j}^\ell  \, - \  {\hat C}_{\,0\,,\ 2}   \cdot  \sum_{j\,=\,1}^\flat \ \sum_{k \,\not=\,j}   {{\lambda_{\,j}^{{n\,-\,2 }\over 2}\!\cdot\lambda_k^{{\,n\,-\ 2\,}\over 2}}\over {\ \ \Vert\ \xi_j\,-\,\xi_k\ \Vert^{\,n\,-\,2}\ \  }}   + \ \\[0.2in]
  & \ & \ \ \ \ \ \   \ \ \ \ \ \ \  \ \ \ \  \ \ \  + \ \,{\hat C}_{\,0\,,\ 3}  \cdot   \sum_{j\,=\,1}^\flat  \ \left[\  C\,(\,{\bf p}_j\,) \cdot \lambda_{\,j}^{\,\ell\ -\,2} \cdot \, \eta_{\,j}^{\,2} \ \right]   \ \ \ \ \ \ \ \ \ [ \ \eta_{\,j} \ = \  {\mbox{Dist}}\,(\,\xi_{\,j}\,,\ {\cal H}\,)\,]\ ,
\end{eqnarray*}
(\,A.13.2\,)
\begin{eqnarray*}
& \ & {\bf G}_{\,l} \,\,(\,\lambda_{\,1}\,, \ \cdot \cdot \cdot\,, \ \lambda_{\ \flat}\,; \ \xi_{\,1}\,, \ \cdot \cdot \cdot\,, \ \xi_{\ \flat}\,) \\[0.2in]
& = & \!\!   {\hat C}_{\,1\,,\ 1}  \cdot C\,(\,{\bf p}_{\,l}\,) \cdot \lambda^{\,\ell}_{\,l}\,   -\,  {\hat C}_{\,1\,,\ 2}  \left\{ \  \sum_{k\ \not= \ l}   \left( {{\lambda_{\,\,l}^{{\,n\,-\ 2\,}\over 2} \cdot \lambda_{\,k}^{{\,n\,-\ 2\,}\over 2}}\over {\Vert\ \xi_l\ - \ \xi_{\,k}\ \Vert^{\,n\,-\,2} }}  \right)
\right\}   + \, {\hat C}_{\,1\,,\ 3}  \cdot C\,(\,{\bf p}_{\,l}\,) \cdot \lambda^{\,\ell\,-\ 2}_{\,l} \cdot  \eta^2_{\,l}  \ ,
\end{eqnarray*}

\vspace*{-0.2in}

and\\[0.1in]
(\,A.13.3\,)
\begin{eqnarray*}
& \ & {{\vec{\,{\bf H}}}}_{\,\,l}\,(\,\lambda_{\,1}\,, \ \cdot \cdot \cdot\,, \ \lambda_{\ \flat}\,; \ \xi_{\,1}\,, \ \cdot \cdot \cdot\,, \ \xi_{\ \flat}\,)\\[0.2in]
& = &     {\hat C}_{\,2\,,\ 2}  \sum_{k\ \not= \ l}  \left\{  \  \left( \ {{ \sqrt{\lambda_{\,\,l}  \cdot \lambda_{\,k}\,} }\over {\,\Vert \ \xi_{\,l}\ - \ \xi_{\,k}\ \Vert \,   }}\  \right)^{\!\!n} \cdot \left( {{\,\xi_{{\,l} }\ - \ \xi_{k}\,}\over {\lambda_k}}\  \right)   \  \right\}   \ + \  \left[\ {\hat C}_{\,2\,,\ 3}\, \cdot   C\, ({\bf p}) \cdot \lambda_{\,\,l}^\ell \cdot\left( {{  \eta_{\,l}  }\over {  \lambda_{\,\,l} }}\  \right) \ {\bf n}_{\,l} \  \right] \\[0.2in]
& \ & \hspace*{4in} \ \ \ \ \ \ \    \left[ \   {\bf n}_{\,l} \ = \  {{ \xi_{\,l} \,-\, {\bf p}_{\,l} }\over {\, \Vert \  \xi_{\,l} \,-\, {\bf p}_{\,l} \,\Vert\, }}\  \ \right]\ ,
 \end{eqnarray*}

where $\,l\ = \ 1\,,\ \cdot \cdot \cdot\,, \ \flat\,,\,$
$$
 {{\vec{\,{\bf H}}}}_{\,\,l} \ = \ ( \   {\bf H}_{\,l_{\, 1}} \ , \cdot \cdot \cdot\,, \   {\bf H}_{\,l_{\, n}}  \ )\ ,   \ \ \ \ {\mbox{and}} \ \ \ \
  V(\,n\,) \ = \  \int_{\R^n} \left(\ {1\over {1\ + \ \Vert\,y\,\Vert^2 }}\ \right)^{\!n}\, d\,y \ .
$$

\newpage

\noindent In (\,A.13.1\,)\,,\, (\,A.13.21\,) and (\,A.13.3\,) \,,\,\\[0.1in]
(\,A.13.4\,)
\begin{eqnarray*}
  {\hat C}_{\,0\,, \ 1} & = & \left(\ {{\,n\,-\ 2\,}\over {2\,n}}\,\right)\,\cdot\,\int_{\R^n} \left(\ {1\over { 1 \ + \ \Vert\ Y\,\Vert^2 }}\ \right)^{\!\!n} \,d\,Y\ ,\\[0.2in]
  {\hat C}_{\,0\,, \ 2} & = &  \left( \, {{\,n \ - \ 2\,}\over 2} \, \right) \cdot \omega_n\ , \\[0.2in]
  {\hat C}_{\,0\,, \ 3} & = & \left(\ {{n\,-\,2}\over {2n}}\,\right)\,\cdot\, {{\ell\,\cdot\, (\,\ell\ - \ 1\,)}\over 2}\,\cdot\,\int_{\R^n} |\,Y_n\,|^{\,\ell\ - \ 2}\,\cdot\, \left(\ {1\over { 1 \ + \ \Vert\ Y\,\Vert^2 }}\ \right)^{\!\!n} \, d\,Y\ ,\\[0.2in]
{\hat C}_{\,1\,,\ 1}  & = &  \left(\ {{\,n\,-\ 2\,}\over {2}}\,\right)\,\cdot\,\int_{\R^n} |\,Y_n\,|^{\,\ell} \cdot  \left(\ {1\over { 1 \ + \ \Vert\ Y\,\Vert^2 }}\ \right)^{\!\!n} \!\cdot \left(\, {{\Vert\ Y\,\Vert^2 \ - \ 1 }\over { \ \  \Vert\ Y\,\Vert^2 \ + \ 1 \ }}\, \right)\, d\,Y\ , \\[0.2in]
  {\hat C}_{\,1\,, \ 2} & = &    {{\,(\,n \ - \ 2\,)^{\,2}\,}\over 2}  \cdot \omega_n \,,  \\[0.2in]
  {\hat C}_{\,1\,, \ 3} & = &   \left(\!{{\ n \ - \ 2\ }\over 2}\! \right)   \cdot {{\ell \cdot (\,\ell\, - \, 1\,)}\over 2} \cdot \int_{\R^n}\! |\,Y_n\,|^{\,\,\ell\, - \, 2} \cdot \left( {1\over { 1 \, + \, \Vert\,Y\,\Vert^2 }}\ \right)^{\!\!n}  \!\cdot \left(\! {{\Vert\ Y\,\Vert^2 \, - \, 1 }\over { \,\Vert\ Y\,\Vert^2 \, + \, 1 \, }}  \right) \, d\,Y, \\[0.2in]
  {\hat C}_{\,2\,,\ 2}  & = & (\,n\,-\,2)^{\,2} \cdot  \omega_n\ ,\\[0.2in]
 {\hat C}_{\,2\,, \ 3}  & = &   \left(\,\emph{}n\,-\,2\,\right)   \cdot \ell \cdot \int_{\R^n} |\,Y_n\,|^{\,\,\ell } \cdot \left( {1\over { 1 \ + \ \Vert\,Y\,\Vert^2 }}\ \right)^{\!\!n\,+\,1} \ d\,Y.
\end{eqnarray*}

 Note that
$$
 \ \ \ \ \ \ \ \ \ \ \ \ \ \   {\hat C}_{\,1\,, \ 2} \ = \ (\,n \ - \ 2\,)\cdot {\hat C}_{\,0\,, \ 2}\ , \leqno (\,A.13.5\,)
$$
$$
 \ \ \ \ \ \ \ \ \ \ \ \ \  \ \ \ \  {\hat C}_{\,2\,, \ 2} \ = \ 2 \cdot\left( \, n \ - \ 2  \, \right) \cdot {\hat C}_{\,0\,, \ 2}\ , \leqno (\,A.13.6\,)
$$
$$
  {\hat C}_{\,1\,, \ 1} \ = \ \ell  \cdot {\hat C}_{\,0\,, \ 1}\ \ \ \ \ \ \ \ \  [\  {\mbox{established \ \ in \ \ {\bf \S\,A\,13\,.\,c} }}\ ] \ , \leqno (\,A.13.7\,)
$$
$$
{\hat C}_{\,1\,, \ 3} \ = \ 2 \cdot\left( \, \ell \ - \ 2  \, \right) \cdot {\hat C}_{\,0\,, \ 3}\ \ \ \ \ \ \  [\  {\mbox{established \ \ in \ \ {\bf \S\,A\,13\,.\,c} }}\ ] \ ,  \leqno (\,A.13.8\,)
$$
 $$
 \  {\hat C}_{\,2\,, \ 3} \ = \ 2 \cdot  {\hat C}_{\,0\,, \ 3}\ \ \ \ \ \   [\  {\mbox{established \ \ in \ \ {\bf \S\,A\,13\,.\,d} }}\ ] \ .  \leqno (\,A.13.9\,)
$$
Beside formal justifications (\,see {\bf \S\,A\,13\,.\,c}\, and\, {\bf \S\,A\,13\,.\,d} \ )\,,\, we seek to discern the relations in (\,A.13.5\,)\,--\,(\,A.13.9\,) via another (\,simpler\,) route\,,\, namely\,,\, by conducting the partial differentiations
$$
\lambda_{\ \!l} \cdot {{\partial \,{\bf F} }\over {\,\partial \,\lambda_{\ \!l}\,} } \ \ \ \ \ {\mbox{and}} \ \ \ \ \lambda_{\ \!l} \cdot {{\partial \,{\bf F} }\over {\,\partial \,\xi_{\,l_{\,|\,j}\,}} } \ \ \ \ \ \ \ \ \ \ \ \ \ \ ( \ l \ = \ 1\,, \ \cdot \cdot \cdot\,, \ \flat\,, \ \ j \  = \ 1\,, \ \cdot \cdot \cdot\,, \ n \ )\ ,
$$
we match the corresponding constants in (\,A.13.2\,) and  (\,A.13.3\,)\,.

\vspace*{0.5in}

{\bf \S\,A\,13\,.\,a\,.} \ \
{\bf Relation in}\,  (\,A.13.5\,) {\bf concerning\,} $\,{\hat C}_{\,0\,, \ 2}\,$ {\bf and\,} ${\hat C}_{\,1\,, \ 2}\,$. \ \ In (\,A.13.1\,)\,,\, the expression for $\,{\bf F}\,$,\, consider the sumation term with the leading coefficient\,:
$$
{\hat C}_{\,0\,,\ 2}   \cdot  \sum_{j\,=\,1}^\flat \ \sum_{k \,\not=\,j} \ {{\lambda_{\,j}^{{n\,-\,2 }\over 2}\!\cdot\lambda_k^{{\,n\,-\ 2\,}\over 2}}\over {\ \ \Vert\ \xi_j\,-\,\xi_k\ \Vert^{\,n\,-\,2}\ \  }} \ ,
$$
in which we can find two expressions of the form
$$
{{\lambda_{\,1}^{{n\,-\,2 }\over 2}\!\cdot\lambda_{\,2}^{{\,n\,-\ 2\,}\over 2}}\over {\ \Vert\ \xi_{\,1}\,-\,\xi_{\,2}\,\Vert^{\,n\,-\,2}\  }}
$$
namely\,,\,
$$
{\hat C}_{\,0\,,\ 2}   \cdot  \left\{ \ {{\lambda_{\,1}^{{n\,-\,2 }\over 2}\!\cdot\lambda_{\,2}^{{\,n\,-\ 2\,}\over 2}}\over {\ \Vert\ \xi_{\,1}\,-\,\xi_{\,2}\,\Vert^{\,n\,-\,2}\  }}    \ \right\}
 \ \ \ \  {\mbox{and}} \ \ \ \
{\hat C}_{\,0\,,\ 2}   \cdot  \left\{ \ {{\lambda_{\,2}^{{n\,-\,2 }\over 2}\!\cdot\lambda_{\,1}^{{\,n\,-\ 2\,}\over 2}}\over {\ \Vert\ \xi_{\,2}\,-\,\xi_{\,1}\,\Vert^{\,n\,-\,2}\  }}    \ \right\} \ .
$$
\ \ \ \ \ \ $j \ = \ 1\,$ and $\,k\ = \ 2 \ \ \uparrow $ \ \ \ \ \ \ \ \ \ \ \ \ \ \ \ \ \ \ \ \ \ \ \ \  $j \ = \ 2\,$ and $\,k\ = \ 1 \ \ \uparrow $ \\[0.1in]
\noindent So the {\it total}\, is
$$
2 \times {\hat C}_{\,0\,,\ 2}   \cdot \left\{ \ {{\lambda_{\,1}^{{n\,-\,2 }\over 2}\!\cdot\lambda_{\,2}^{{\,n\,-\ 2\,}\over 2}}\over {\ \Vert\ \xi_{\,1}\,-\,\xi_{\,2}\,\Vert^{\,n\,-\,2}\  }}    \ \right\} \ .\leqno (\,A.13.10\,)
$$

Focusing on the term
$$
 {{\lambda_{\,1}^{{n\,-\,2 }\over 2}\!\cdot\lambda_{\,2}^{{\,n\,-\ 2\,}\over 2}}\over {\ \Vert\ \xi_{\,1}\,-\,\xi_{\,2}\,\Vert^{\,n\,-\,2}\  }}  \ ,
$$

\newpage

we carry out  the partial differentiation \\[0.1in]
(\,A.13.11\,)
\begin{eqnarray*}
& \ &  \left( \lambda_{\,\,1}\cdot {\partial\over {\partial\, \lambda_{\,\,1}}}\ \right) {\bf F}\\[0.2in] & = &  \left( \lambda_{\,\,l}\cdot {\partial\over {\partial\, \lambda_{\,\,l}}}\ \right) \ \left[ \ \cdot \cdot \cdot   -\, 2\cdot {\hat C}_{\,0\,, \ 2}  \cdot \left\{ \ {{\lambda_{\,1}^{{n\,-\,2 }\over 2}\!\cdot\lambda_{\,2}^{{\,n\,-\ 2\,}\over 2}}\over {\ \Vert\ \xi_{\,1}\,-\,\xi_{\,2}\,\Vert^{\,n\,-\,2}\  }}    \ \right\}\ + \cdot \cdot \cdot \right]  \\[0.2in] & = & \cdot \cdot \cdot \  -\, 2\cdot {\hat C}_{\,0\,, \ 2}  \cdot \left( \  \lambda_{\,1} \cdot {{\partial }\over {\partial\,\lambda_{\,1}}}\ \right)\ \left\{ \ {{\lambda_{\,1}^{{n\,-\,2 }\over 2}\!\cdot\lambda_{\,2}^{{\,n\,-\ 2\,}\over 2}}\over {\ \Vert\ \xi_{\,1}\,-\,\xi_{\,2}\,\Vert^{\,n\,-\,2}\  }}    \ \right\} \  \ + \ \cdot \cdot \cdot\\[0.2in] & = & \cdot \cdot \cdot \  -\  2\cdot {\hat C}_{\,0\,, \ 2}   \cdot  \left( \   {{n\,-\,2 }\over 2} \ \right) \cdot  \left\{ \ {{\lambda_{\,1}^{{n\,-\,2 }\over 2}\!\cdot\lambda_{\,2}^{{\,n\,-\ 2\,}\over 2}}\over {\ \Vert\ \xi_{\,1}\,-\,\xi_{\,2}\,\Vert^{\,n\,-\,2}\  }}    \ \right\} \  \ + \ \cdot \cdot \cdot \\[0.2in] & = & \cdot \cdot \cdot \  -\  (\,n \ - \ 2\,) \cdot  {\hat C}_{\,0\,, \ 2}    \cdot   \left\{ \ {{\lambda_{\,1}^{{n\,-\,2 }\over 2}\!\cdot\lambda_{\,2}^{{\,n\,-\ 2\,}\over 2}}\over {\ \Vert\ \xi_{\,1}\,-\,\xi_{\,2}\,\Vert^{\,n\,-\,2}\  }}    \ \right\} \  \ + \ \cdot \cdot \cdot\ .\\
\end{eqnarray*}
By comparing with the similar term\,:
$$
 {{\lambda_{\,1}^{{n\,-\,2 }\over 2}\!\cdot\lambda_{\,2}^{{\,n\,-\ 2\,}\over 2}}\over {\ \Vert\ \xi_{\,1}\,-\,\xi_{\,2}\,\Vert^{\,n\,-\,2}\  }}  \ ,
$$
found in $\,{\bf G}_{\,1}\,$ [\,refer to (\,A.13.2\,)\ ]\,,\,  we see that
\begin{eqnarray*}
& \ & -\,  {\hat C}_{\,1\,, \ 2}   \left\{ \ {{\lambda_{\,1}^{{n\,-\,2 }\over 2}\!\cdot\lambda_{\,2}^{{\,n\,-\ 2\,}\over 2}}\over {\ \Vert\ \xi_{\,1}\,-\,\xi_{\,2}\,\Vert^{\,n\,-\,2}\  }}    \ \right\} \ = \  -\   (\,n\ - \ 2\,) \cdot {\hat C}_{\,0\,, \ 2}   \left\{ \ {{\lambda_{\,1}^{{n\,-\,2 }\over 2}\!\cdot\lambda_{\,2}^{{\,n\,-\ 2\,}\over 2}}\over {\ \Vert\ \xi_{\,1}\,-\,\xi_{\,2}\,\Vert^{\,n\,-\,2}\  }}    \ \right\}\\[0.2in]
\Longrightarrow \ \  & \ &   {\hat C}_{\,1\,,\ 2} \ =\    (\,n\ - \ 2\,) \cdot {\hat C}_{\,0\,, \ 2} \ .\\
\end{eqnarray*}
Compare with (\,A.13.5\,)\,.

\newpage

{\bf \S\,A\,13\,.\,b\,. } \ \ {\bf Relation in}\,  (\,A.13.6\,) {\bf concerning\,} $\,{\hat C}_{\,0\,, \ 2}\,$ {\bf and\,} ${\hat C}_{\,2\,, \ 2}\,$. \ \
For the same expression found in (\,A.13.10\,)\,,\, let us consider the partial derivative \\[0.1in]
(\,A.13.12\,)

\vspace*{-0.25in}

\begin{eqnarray*}
& \ &  \left( \lambda_{\,\,1}\cdot{{\partial}\over {\partial \,\xi_{1_{\,|_{\,1}}} }}\ \right) {\bf F}\\[0.2in] & = &  \left( \lambda_{\,\,l}\cdot {{\partial}\over {\partial \,\xi_{1_{\,|_{\,1}}} }}\ \right) \ \left[ \ \cdot \cdot \cdot   -\, 2\cdot {\hat C}_{\,0\,, \ 2}  \cdot \left\{ \ {{\lambda_{\,1}^{{n\,-\,2 }\over 2}\!\cdot\lambda_{\,2}^{{\,n\,-\ 2\,}\over 2}}\over {\ \Vert\ \xi_{\,1}\,-\,\xi_{\,2}\,\Vert^{\,n\,-\,2}\  }}    \ \right\}\ + \cdot \cdot \cdot \right]  \\[0.2in] & = & \cdot \cdot \cdot \  -\, 2\cdot {\hat C}_{\,0\,, \ 2}  \cdot \left( \  \lambda_{\,1} \cdot {{\partial}\over {\partial \,\xi_{1_{\,|_{\,1}}} }}\ \right)\ \left\{ \ {{\lambda_{\,1}^{{n\,-\,2 }\over 2}\!\cdot\lambda_{\,2}^{{\,n\,-\ 2\,}\over 2}}\over {\ \Vert\ \xi_{\,1}\,-\,\xi_{\,2}\,\Vert^{\,n\,-\,2}\  }}    \ \right\} \  \ + \ \cdot \cdot \cdot \\[0.2in] & = & \cdot \cdot \cdot \  -\, 2\cdot {\hat C}_{\,0\,, \ 2}  \cdot \left( \  \lambda_{\,1} \cdot {{\partial}\over {\partial \,\xi_{1_{\,|_{\,1}}} }}\ \right)\ \left\{ \ {{\lambda_{\,1}^{{n\,-\,2 }\over 2}\!\cdot\lambda_{\,2}^{{\,n\,-\ 2\,}\over 2}}\over {\   \left[\ (\,\xi_{1_{\,|_1}} \, - \, \xi_{2_{\,|_1}}\,)^{\,2} \ + \ \cdot\,\cdot\,\cdot\,\ \right]^{\,{{\,n\,-\,2\,}\over 2}} \  }}    \ \right\} \  \ + \ \cdot \cdot \cdot \\[0.2in]
& = & \cdot \cdot \cdot \  -\, 2\,\left(\,-\,{{\,n \ - \ 2\,}\over 2} * 2\,\right) \cdot {\hat C}_{\,0\,, \ 2}  \cdot \left\{ \ {{\lambda_{\,1}^{{n }\over 2}\!\cdot\lambda_{\,2}^{{\,n\,-\ 2\,}\over 2} \cdot  (\,\xi_{1_{\,|_{\,1}}} \, - \, \xi_{2_{\,|_{\,1}}}\,)  }\over {\   \left[\ (\,\xi_{1_{\,|_1}} \, - \, \xi_{2_{\,|_1}}\,)^{\,2} \ + \ \cdot\,\cdot\,\cdot\,\ \right]^{\,{{\,n\,}\over 2}} \  }}    \ \right\}  \  \ + \ \cdot \cdot \cdot \\[0.2in] & = & \cdot \cdot \cdot \  +\, 2\,\left(\, n \ - \ 2\, \right) \cdot {\hat C}_{\,0\,, \ 2}  \cdot \left\{ \ {\ {\lambda_{\,1}^{{n }\over 2}\!\cdot\lambda_{\,2}^{{\,n\,-\ 2\,}\over 2} \cdot  (\,\xi_{1_{\,|_{\,1}}} \, - \, \xi_{2_{\,|_{\,1}}}\,) \  }\over {\ \Vert\ \xi_{\,1}\,-\,\xi_{\,2}\,\Vert^{\,n} \  }}    \ \right\}  \  \ + \ \cdot \cdot \cdot
\end{eqnarray*}
By comparing with the similar term\,:
$$
  \  \left( \ {{ \sqrt{\lambda_{\,\,1}  \cdot \lambda_{\,2}\,} }\over {\Vert \ \xi_{\,1}\ - \ \xi_{\,2}\ \Vert    }}\  \right)^{\!\!n} \cdot \left( {{\,\xi_{{\,1_{\,|_{\,1}}} }\ - \ \xi_{{\,2_{\,|_{\,1}}} }\,}\over {\lambda_{\,2}}}\  \right)   \  \ \left( \ = \ \ {{\ \lambda_{\,1}^{{n }\over 2}\!\cdot\lambda_{\,2}^{{\,n\,-\ 2\,}\over 2} \cdot  (\,\xi_{1_{\,|_{\,1}}} \, - \, \xi_{2_{\,|_{\,1}}}\,) \  }\over {\ \Vert\ \xi_{\,1}\,-\,\xi_{\,2}\,\Vert^{\,n} \  }}   \ \right)
$$
found in $\,{\bf H}_{\,1_{\,|_1}}\,$ [\ refer to (\,A.13.3\,)\ ]\,,\,  we see that
\begin{eqnarray*}
& \ &   {\hat C}_{\,2\,, \ 2}  \cdot \left( \ {{ \sqrt{\lambda_{\,\,1}  \cdot \lambda_{\,2}\,} }\over {\Vert \ \xi_{\,1}\ - \ \xi_{\,2}\ \Vert    }}\  \right)^{\!\!n} \cdot \left( {{\,\xi_{{\,1_{\,|_1}} }\ - \ \xi_{{\,2_{\,|_1}} }\,}\over {\lambda_{\,2}}}\  \right)  \\[0.2in]
& = &   2\cdot  (\,n\, - \, 2\,) \cdot {\hat C}_{\,0\,, \ 2} \cdot    {{\ \lambda_{\,1}^{{n }\over 2}\!\cdot\lambda_{\,2}^{{\,n\,-\ 2\,}\over 2} \cdot  (\,\xi_{1_{\,|_{\,1}}} \, - \, \xi_{2_{\,|_{\,1}}}\,) \  }\over {\ \Vert\ \xi_{\,1}\,-\,\xi_{\,2}\,\Vert^{\,n} \  }} \\[0.2in]
\Longrightarrow & \ &   {\hat C}_{\,2\,,\ 2}  \, = \,   2\cdot (\,n\, - \, 2\,) \cdot {\hat C}_{\,0\,, \ 2} \\[0.1in]
& \ &  \ \ \ \ \ \    [ \ {\mbox{can \ \ be \ \ seen \ \ from \ \ the \ \ corresponding \ \ numbers \ \ in  }} \ \  (\,A.13.4\,)\ ]\,.
\end{eqnarray*}

\vspace*{0.35in}

{\bf \S\,A\,13\,.\,c\,. }   {\bf Relation in}\,  (\,A.13.7\,)\, {\bf concerning\,}   $\,{\hat C}_{0\,, \ 1}\, $ and $\,{\hat C}_{\,1\,, \ 1}\,,\,$  {\bf and \ in }\,  (\,A.13.8\,) {\bf concerning\,} $\,{\hat C}_{\,0\,, \ 3}\,$ {\bf and\,} ${\hat C}_{\,1\,, \ 3}\,$.\ \ Let us consider the coefficient in front of $\,\lambda_{\,1}^\ell\,$ in the expression for $\,{\bf F}\,$,\, cf. (\,A.13.1\,)\,.\,
Here we apply the partial derivative \\[0.1in]
(\,A.13.13\,)
\begin{eqnarray*}
& \ &  \left( \lambda_{\,\,1}\cdot{{\partial}\over {\partial \,\lambda_{\,1}}}\ \right) {\bf F}\\[0.2in] & = &   \left( \lambda_{\,\,1}\cdot{{\partial}\over {\partial \,\lambda_{\,1}}}\ \right) \ \left[ \ \cdot \cdot \cdot  \ + \  {\hat C}_{0\,, \ 1} \cdot C\,(\,{\bf p}_{\,1}\,) \cdot \lambda_{\,1}^\ell   \ +  \cdot \cdot \cdot \right]  \\[0.2in] & = & \cdot \cdot \cdot \  + \  \ell \cdot {\hat C}_{0\,, \ 1} \cdot C\,(\,{\bf p}_{\,1}\,) \cdot \lambda_{\,1}^\ell \ + \ \cdot \cdot \cdot
\end{eqnarray*}
By comparing with the similar term
found in $\,{\bf G}_{\,1}\,$ [\,refer to (\,A.13.2\,)\ ]\,,\,  we see that
$$ {\hat C}_{1\,, \ 1} \ = \ \ell \cdot {\hat C}_{0\,, \ 1}\ .$$
Here we provide a formal justification. From (\,A.13.4\,)\,,\, we have\\[0.1in]
(\,A.13.14\,)
\begin{eqnarray*}
{\hat C}_{0\,, \ 1}   & = &   \left(\ {{\,n\,-\ 2\,}\over {2\,n}}\,\right)\,\cdot\,\int_{\R^n} |\,Y_n\,|^\ell \cdot  \left(\ {1\over { 1 \ + \ \Vert\ Y\,\Vert^2 }}\ \right)^{\!\!n} \,d\,Y\\[0.2in]
\Longrightarrow \ \ {\hat C}_{0\,, \ 1} \cdot \lambda_{\,1}^\ell & = &  \lambda_{\,1}^\ell \cdot  \left(\ {{\,n\,-\ 2\,}\over {2\,n}}\,\right)\,\cdot\,\int_{\R^n} |\,Y_n\,|^\ell \cdot  \left(\ {1\over { 1 \ + \ \Vert\ Y\,\Vert^2 }}\ \right)^{\!\!n} \,d\,Y\\[0.2in]
& = & \left(\,{{n\,-\,2}\over {2n}}\,\right)\,\cdot\,\int_{\R^n} \,|\,y_n\,-\,\xi_{1_{\,|_{\,n}}}\,|^\ell\,\cdot\,   \left(\ {{\lambda_{\,1}}\over  { \lambda_{\,1}^2 \ + \ \Vert\  y \ - \ \xi_{\,1}\,\Vert^2}}\  \right)^n\ d\,y\\[0.2in]
& \ & \hspace*{0.8in} \left( \ {\mbox{change \ \ of \ \ the \ \ dummy \ \ variables}} \ \ \ Y \ \to \ {{ y \ - \ \xi_{\,1}}\over {\lambda_{\,1}}} \ \right)\\[0.1in]
& = & \left(\,{{n\,-\,2}\over {2n}}\,\right)\,\cdot\,\int_{\R^n} \,|\,{\tilde y}_n\,|^\ell\,\cdot\,   \left(\ {{\lambda_{\,1}}\over  { \lambda_{\,1}^2 \ + \ \Vert\  {\tilde y}\,\Vert^2}}\  \right)^n\ dV_{\tilde y}\\[0.2in]
& \ & \hspace*{0.8in} \left( \ {\mbox{change \ \ of \ \ the \ \ dummy \ \ variables \ \ }} \ \tilde y \ \to \  y \ - \ \xi_{\,1} \ \right)
\ . \\
\end{eqnarray*}

\newpage

{\it Two ways to find the derivative.}\ \ Here we record down the  partial differentiation\,:
$$  \left( \lambda_{\,\,1}\cdot{{\partial}\over {\partial \,\lambda_{\,1}}}\ \right) [ \  {\hat C}_{0\,, \ 1} \cdot \lambda_{\,1}^\ell \ ]  \ = \ \ell \cdot {\hat C}_{0\,, \ 1} \cdot \lambda_{\,1}^\ell \ . \leqno (\,A.13.15\,)
$$

Recall from (\,A.13.4\,) that
$$
{\hat C}_{\,1\,, \ 1} \ = \ \left(\ {{\,n\,-\ 2\,}\over {2}}\,\right)\,\cdot\,\int_{\R^n} |\,Y_n\,|^{\,\ell} \cdot  \left(\ {1\over { 1 \ + \ \Vert\ Y\,\Vert^2 }}\ \right)^{\!\!n} \!\cdot \left(\, {{\Vert\ Y\,\Vert^2 \ - \ 1 }\over { \ \  \Vert\ Y\,\Vert^2 \ + \ 1 \ }}\, \right)\, d\,Y \ . \leqno (\,A.3.16\,)
$$
We find the same derivative via (\,A.13.14\,)\,:\\[0.1in]
(\,A.13.17\,)
\begin{eqnarray*}
&\   & \left(  \lambda_{\,1} \cdot {\partial\over {\partial \,\lambda_{\,1}}} \right) \left[\  {\hat C}_{0\,, \ 1} \cdot \lambda_{\,1}^\ell \ \right] \ \ \left( \ = \ \ell \cdot {\hat C}_{0\,, \ 1} \cdot \lambda_{\,1}^\ell \right) \\[0.2in]
&  = &\left(  \lambda_{\,1} \cdot {\partial\over {\partial \,\lambda_{\,1}}} \right) \left[\  \left(\ {{\,n\,-\ 2\,}\over {2\,n}}\,\right)\,\cdot\, \int_{\R^n}   |\,{\tilde y}_n  \,|^{\,\ell}    \times   \left(\, {\lambda_{\,1}\over {\lambda^2_1 \ + \ \Vert\ {\tilde y}\,\Vert^2}}\, \right)^{\!\!n }  \ d\,{\tilde y} \ \right]  \\[0.2in]
& = &    \left(\ {{\,n\,-\ 2\,}\over {2\,n}}\,\right)\,\cdot\, \lambda_{\,1} \cdot   \int_{\R^n}   |\,{\tilde y}_n  \,|^{\,\ell}    \times   \left[\   {\partial\over {\partial \,\lambda_{\,1}}} \,  \left(\, {\lambda_{\,1}\over {\lambda^2_1 \ + \ \Vert\ {\tilde y}\,\Vert^2}}\, \right)^{\!\!n } \ \right] \ d\,{\tilde y} \\[0.2in]
& = & \lambda_{\,1} \cdot  \left(\ {{\,n\,-\ 2\,}\over {2\,n}}\,\right)\,\cdot\,   n\,\cdot\,{1\over {\lambda_{\,1}}}\,*\\[0.2in]
& \ & \hspace*{1in} *\,   \int_{\R^n}   |\,{\tilde y}_n  \,|^{\,\ell}    \times   \left[\   \left(\, {\lambda_{\,1}\over {\lambda^2_1 \ + \ \Vert\ {\tilde y}\,\Vert^2}}\, \right)^{\!\!n }\,\cdot\, \left(\  {{\Vert\ {\tilde y}\,\Vert^2 \ - \ \lambda_{\,1}^2}\over {\ \lambda^2_1 \ + \ \Vert\ {\tilde y}\,\Vert^2\ }}\  \right)\  \right] \ d\,{\tilde y}\\[0.2in]
& = & \lambda_{\,1}^\ell \cdot  \left(\ {{\,n\,-\ 2\,}\over {2}}\,\right)\,\cdot\,\int_{\R^n} |\,Y_n\,|^{\,\ell} \cdot  \left(\ {1\over { 1 \ + \ \Vert\ Y\,\Vert^2 }}\ \right)^{\!\!n} \!\cdot \left(\, {{\Vert\ Y\,\Vert^2 \ - \ 1 }\over { \ \  \Vert\ Y\,\Vert^2 \ + \ 1 \ }}\, \right)\, d\,Y\\[0.1in]
& \ & \hspace*{2in} \left( \ {\mbox{change \ \ of \ \ the \ \ dummy \ \ variables \ \ }} \ Y \ \to \ {{ \tilde y }\over {\lambda_{\,1}}} \ \right)\\[0.1in]
& = & {\hat C}_{\,1\,, \ 1} \cdot \lambda_{\,1}^\ell  \ \ \ \ \ \ \ \ \ \  \ \ \ \ \  \ \ \ \ \  \ \ \ \ \   [ \ {\mbox{recall}} \ \ (\,A.13.16\,) \ ]\\[0.2in]
\Longrightarrow & \ & \ \ {\hat C}_{\,1\,, \ 1}  \ = \ \ell \cdot {\hat C}_{0\,, \ 1} \ .
\end{eqnarray*}
Similarly procedure can be applied to the relation
$${\hat C}_{\,1\,, \ 3} \ = \ (\,\ell\,-\,2\,) \cdot {\hat C}_{\,0\,, \ 3}\,.\,$$

\newpage

{\bf \S\,A\,13\,.\,d\,. }   {\bf Relation in}\,  (\,A.13.9\,) {\bf concerning\,}\,   ${\hat C}_{0\,, \ 3}\, $ and $\,{\hat C}_{\,2\,, \ 3}\,.\,$ \ \ Recall that the gradient of a function points toward the direction in which the function increases most rapidly\,.\, Via formal differentiation,
$$
\btd_{\xi_{\,1}}\, [\ {\mbox{Dist.}}\ (\,\xi_j\,, \ {\cal H}\,) \ ]^{\,2} \ = \ 2\cdot [\ {\mbox{Dist.}}\ (\,\xi_j\,, \ {\cal H}\,) \ ]\, {\bf n}_{\,1}\,,
$$
where $\,{\bf n}_1\,$ is the unit vector pointing from $\,{\bf p}_{\,\xi_{\,1}}\,$ to $\ \xi_{\,1}\,.\,$
Consider then  \\[0.1in]
(\,A.13.18\,)
\begin{eqnarray*}
& \ &  \left(\ \lambda_{\,\,1}\cdot \btd_{\xi_{\,1}} \ \right) {\bf F}\\[0.2in] & = &    \left(\ \lambda_{\,\,1}\cdot \btd_{\xi_{\,1}} \ \right)  \ \left[ \ \cdot \cdot \cdot  \ + \  {\hat C}_{0\,, \ 3} \cdot C\,(\,{\bf p}_{\,1}\,) \cdot \lambda_{\,1}^{\,\ell \ - \ 2} \cdot  [\ {\mbox{Dist.}}\ (\,\xi_1\,, \ {\cal H}\,) \ ]^{\,2}    \ +  \cdot \cdot \cdot \right]  \\[0.2in] & = & \ \cdot \cdot \cdot  \ + \ 2  \cdot   {\hat C}_{0\,, \ 3} \cdot C\,(\,{\bf p}_{\,1}\,) \cdot \lambda_{\,1}^{\,\ell \ - \ 1} \cdot  [\ {\mbox{Dist.}}\ (\,\xi_1\,, \ {\cal H}\,) \ ]\,  {\bf n}_{\,1}   \ +  \cdot \cdot \cdot
\end{eqnarray*}
By comparing with the similar term\,:
$$
 {\hat C}_{2\,, \ 3} \cdot C\,(\,{\bf p}_{\,1}\,) \cdot \lambda_{\,1}^{\,\ell } \cdot  \left(  {\eta_{\,1}\over \lambda_{\,1} } \right)\,  {\bf n}_{\,1}  \ = \  {\hat C}_{2\,, \ 3} \cdot C\,(\,{\bf p}_{\,1}\,) \cdot \lambda_{\,1}^{\,\ell \ - \ 1} \cdot  [\ {\mbox{Dist.}}\ (\,\xi_1\,, \ {\cal H}\,) \ ]\,  {\bf n}_{\,1}
$$
found in $\,{\bf H}_{\,1_{\,|_1}}\,$ [\,refer to (\,A.13.3\,)\ ]\,,\,  we see that
$$
{\hat C}_{\,2\,, \ 3} \  = \  2 \cdot {\hat C}_{\,0\,, \ 3}\ ,
$$
Let us provide a formal justification based on direct calculation.
Recall the constant [\ see (\,A.13.4\,)\ ]
$$
{\hat C}_{\,0\,, \ 3} \ = \   \left(\,{{n\,-\,2}\over {2n}}\,\right)\,\cdot\,{{\,\ell\,\cdot\,(\,\ell\ - \ 1\,)}\over 2}\,\cdot\,\int_{\R^n} |\,Y_n\,|^{\,\ell\ - \ 2}\,\cdot\, \left(\ {1\over { 1 \ + \ \Vert\ Y\,\Vert^2 }}\ \right)^{\!\!n} \,d\,Y \ , \leqno (\,A.13.19\,)
$$
and the constant
 $$
{\hat C}_{\,2\,, \ 3} \ = \   \left(\,n\,-\,2\,\right)\,\cdot\, \ell\,\cdot\,\int_{\R^n} |\,Y_n\,|^{\,\ell }\,\cdot\,\left(\ {1\over { 1 \ + \ \Vert\ Y\,\Vert^2 }}\ \right)^{\!\!n\,+\,1} \,d\,Y \ . \leqno (\,A.13.20\,)
$$

Hence
\begin{eqnarray*}
{\hat C}_{\,2\,, \ 3}&=& 2 \cdot {\hat C}_{\,0\,, \ 3}\\[0.2in]
\Longleftrightarrow & \ &    2 \times  \left(\,{{n\,-\,2}\over {2n}}\,\right)\,\cdot\,{{\ell\,\cdot\,(\,\ell\ - \ 1\,)}\over 2}\,\cdot\,\int_{\R^n} |\,Y_n\,|^{\,\ell\ - \ 2}\,\cdot\, \left(\ {1\over { 1 \ + \ \Vert\ Y\,\Vert^2 }}\ \right)^{\!\!n} \,d\,Y \\[0.2in]
& \ & \!\!\!\!\!\!\!=  \  \left(\,n\,-\,2\,\right)\,\cdot\, \ell\,\cdot\,\int_{\R^n} |\,Y_n\,|^{\,\ell }\,\cdot\,\left(\ {1\over { 1 \ + \ \Vert\ Y\,\Vert^2 }}\ \right)^{\!\!n\,+\,1} \,d\,Y \ .
\end{eqnarray*}
That is\,,\, we seek to show that \\[0.1in]
(\,A.13.21\,)
$$
\int_{\R^n} |\,Y_n\,|^{\,\ell }\,\cdot\,\left(\ {1\over { \ 1 \ + \ \Vert\ Y\,\Vert^2 \ }}\ \right)^{\!\!n\,+\,1} \,d\,Y \ = \ \left( \ {{\ell \ - \ 1 }\over {2\,n}}\ \right)\,\cdot\,\int_{\R^n} |\,Y_n\,|^{\,\ell\ - \ 2}\,\cdot\,\left(\ {1\over { 1 \ + \ \Vert\ Y\,\Vert^2 }}\ \right)^{\!\!n} \,d\,Y \
$$
for $\,n \ - \ 2 \ > \ \ell  \ \ge \ 2\,.$\bk
First consider the case when $\,\ell \ = \ 2\,.\,$ Let us demonstrate that\\[0.1in]
(\,A.13.22\,)
$$
\int_{\R^n} |\,Y_n\,|^{\,2}\,\cdot\,\left(\ {1\over { 1 \ + \ \Vert\ Y\,\Vert^2 }}\ \right)^{\!\!n\,+\,1} \,d\,Y \ = \ {{1 }\over {2\,n}}\,\cdot\,\int_{\R^n} \left(\ {1\over { \ 1 \ + \ \Vert\ Y\,\Vert^2 \  }}\ \right)^{\!\!n} \,d\,Y \ .
$$

Via symmetry
\begin{eqnarray*}
\int_{\R^n} |\,Y_n\,|^{\,2}\,\cdot\,\left(\ {1\over { 1 \ + \ \Vert\ Y\,\Vert^2 }}\ \right)^{\!\!n\,+\,1} \,d\,Y & = & {1\over n}\,\cdot\,\int_{\R^n} \Vert\,Y\,\Vert^{\,2}\,\cdot\,\left(\ {1\over { 1 \ + \ \Vert\ Y\,\Vert^2 }}\ \right)^{\!\!n\,+\,1} \,d\,Y\\[0.2in]
& = &  {{\omega_n}\over n} \int_0^\infty {{R^2}\over {(\,1 \ + \ R^2\,)^{\,n\,+\,1} }} \ R^{\,n\,-\,1}\,\cdot\, d\,R\\[0.2in]
& = &   {{\omega_n}\over n}  \int_0^{\pi\over 2} {{\ (\,\tan\,\theta\,)^{n\,+\,1}\,\cdot\,\sec^2 \,\theta\ d\theta\  }\over { (\,\sec\,\theta\,)^{2n}\,\cdot\, \sec^2 \,\theta  }}\\[0.2in]
& = &   {{\omega_n}\over n}  \int_0^{\pi\over 2}  (\,\sin\,\theta\,)^{ n\,+\,1}\,\cdot\,(\,\cos\,\theta\,)^{ n\,-\,1} \ d\,\theta \ .
\end{eqnarray*}

Similarly\,,\,
$$ {{1 }\over {2\,n}}\,\cdot\,
\int_{\R^n} \left(\ {1\over { 1 \ + \ \Vert\ Y\,\Vert^2 }}\ \right)^{\!\!n} \,d\,Y  \ = \  {{\, \omega_n\, }\over {2\,n}}\,\cdot \int_0^{\pi\over 2}  (\,\sin\,\theta\,)^{ n\,-\,1}\,\cdot\,(\,\cos\,\theta\,)^{ n\,-\,1} \ d\,\theta \ .
$$
Therefore, in order to establish (\,A.13.22\,), we need to show that
$$
 \int_0^{\pi\over 2}  (\,\sin\,\theta\,)^{ n\,+\,1}\,\cdot\,(\,\cos\,\theta\,)^{ n\,-\,1} \ d\,\theta  \ =\ {1\over 2} \cdot  \int_0^{\pi\over 2}  (\,\sin\,\theta\,)^{ n\,-\,1}\,\cdot\,(\,\cos\,\theta\,)^{ n\,-\,1} \ d\,\theta\ . \leqno (\,A.13.23\,)
$$

We continue with
\begin{eqnarray*}
& \   &    \int_0^{\pi\over 2}  (\,\sin\,\theta\,)^{ n\,+\,1}\,\cdot\,(\,\cos\,\theta\,)^{ n\,-\,1} \ d\,\theta \\[0.2in]& \ = &    \int_0^{\pi\over 2}  (\,\sin\,\theta\,)^{ n\,-\,1}\,\cdot\,(\, 1  - \ \cos^2 \,\theta\,)\,\cdot\, (\,\cos\,\theta\,)^{ n\,-\,1} \ d\,\theta \\[0.2in]& \ = &    \int_0^{\pi\over 2}  (\,\sin\,\theta\,)^{ n\,-\,1}  \,\cdot\, (\,\cos\,\theta\,)^{ n\,-\,1} \ d\,\theta \ - \  \int_0^{\pi\over 2}  (\,\sin\,\theta\,)^{ n\,-\,1}  \,\cdot\, (\,\cos\,\theta\,)^{ n\,+\,1} \ d\,\theta \\[0.2in]
\Longrightarrow \ \  & \ & \int_0^{\pi\over 2}  (\,\sin\,\theta\,)^{ n\,+\,1}\,\cdot\,(\,\cos\,\theta\,)^{ n\,-\,1} \ d\,\theta \ + \  \int_0^{\pi\over 2}  (\,\sin\,\theta\,)^{ n\,-\,1}\,\cdot\,(\,\cos\,\theta\,)^{ n\,+\,1} \ d\,\theta\\[0.2in]
& = &   \int_0^{\pi\over 2}  (\,\sin\,\theta\,)^{ n\,-\,1}  \,\cdot\, (\,\cos\,\theta\,)^{ n\,-\,1} \ d\,\theta \ .\\
\end{eqnarray*}

Thus it suffices to show that

\begin{eqnarray*}
\int_0^{\pi\over 2}  (\,\sin\,\theta\,)^{ n\,+\,1}\,\cdot\,(\,\cos\,\theta\,)^{ n\,-\,1} \ d\,\theta & = &  \int_0^{\pi\over 2}  (\,\sin\,\theta\,)^{ n\,-\,1}  \,\cdot\, (\,\cos\,\theta\,)^{ n\,+\,1} \ d\,\theta \\[0.2in]
|| \hspace*{1.2in}& \ & \hspace*{1.2in} ||   \\[0.2in]
-\,\int_0^{\pi\over 2}  (\,\sin\,\theta\,)^{ n }\,\cdot\,(\,\cos\,\theta\,)^{ n\,-\,1} \ d\,[\,\cos\,\theta\,]\ \ \ \ \  & \  & \ \ \ \ \   \int_0^{\pi\over 2}  (\,\sin\,\theta\,)^{ n\,-\,1}  \,\cdot\, (\,\cos\,\theta\,)^{ n } \ d\,[\,\sin\,\theta\,] \\[0.2in]
|| \hspace*{1.2in}& \ & \hspace*{1.2in} ||   \\[0.2in]
-\,\int_1^{\,0}  (\,1 \ - \ u^2\,)^{ n\over 2 }\,\cdot\,u^{ n\,-\,1} \ d\,u\ \ \ \ \  & \  & \ \ \ \ \   \int_0^1  v^{ n\,-\,1}  \,\cdot\, (\,1 \ - \ v^2\,)^{ n\over 2 } \ d\,v\\[0.2in]
|| \hspace*{1.2in}& \ & \hspace*{1.2in}  ||   \\[0.2in]
 \int_0^1  (\,1 \ - \ u^2\,)^{ n\over 2 }\,\cdot\,u^{ n\,-\,1} \ d\,u\ \ \ \ \  & =   & \ \ \ \ \ \int_0^1  (\,1 \ - \ v^2\,)^{ n\over 2 }\,\cdot\,v^{ n\,-\,1}   \ d\,v\ .
\end{eqnarray*}
With this, we obtain  (\,A.13.23\,)\,,\, and hence (\,A.13.21\,) when $\,\ell \ = \ 2\,.$

\newpage

{\it General case for} $\ n \ - \ 2 \ > \  \ell \ \ge \ 3\,.\,$\\[0.1in]
(\,A.13.24\,)
\begin{eqnarray*}
& \ & \int_{\R^n} |\,Y_n\,|^{\,\ell }\,\cdot\,\left(\ {1\over { 1 \ + \ \Vert\ Y\,\Vert^2 }}\ \right)^{\!\!n\,+\,1} \,d\,Y \\[0.2in]
& =  & \int_{\R^n} |\,Y_n\,|^{\,\ell }\,\cdot\,\left(\ {1\over { [\ 1 \ + \ |\,Y_n\,|^{\,2 }\ ] \ + \ \Vert\ {\tilde Y}\,\Vert^2 }}\ \right)^{\!\!n\,+\,1} \,d\,Y \\[0.2in]
& = &  \int_{-\,\infty}^\infty |\,Y_n\,|^{\,\ell }\,\cdot\,\left[\ \int_{\R^{\,n \, - \, 1} }  \left(\ \ {1\over { \left[\ \sqrt{ 1 \ + \ |\,Y_n\,|^{\,2 }\ }\ \right]^{\,2} \ + \ \Vert\ {\tilde Y}\,\Vert^2 }}\ \right)^{\!\!n\,+\,1} \,d\,{\tilde Y} \ \right]\ d\,Y_n \\[0.1in]
& \ &  \hspace*{3.5in} ( \ d\,{\tilde Y} \ = \ dY_1\ dY_2 \,\cdot \cdot \cdot\,dY_{n\,-\,1} \ ) \\[0.1in]
& = &  \int_{-\,\infty}^\infty |\,Y_n\,|^{\,\ell }\,\cdot\,\left[\  \omega_{n\,-\,1}\,\cdot\,\int_0^\infty  \left(\ \ {1\over { \Lambda^{\,2} \ + \ {\tilde R}^2 }}\ \right)^{\!\!n\,+\,1} \ {\tilde R}^{n \, - \, 2}\ d\,{\tilde R} \ \right]\ d\,Y_n\\[0.1in]
& \ &  \hspace*{3.5in} ( \  \Lambda \ = \ \sqrt{ 1 \ + \ |\,Y_n\,|^{\,2 }\,}\  ) \\[0.1in]
& = &  \int_{-\,\infty}^\infty |\,Y_n\,|^{\,\ell }\,\cdot\,\left[\  \omega_{n\,-\,1}\,\cdot\,\int_0^{\pi\over 2}   {{ \ \Lambda^{n\,-\,1}\,\cdot\,[\ \tan\,\theta\ ]^{\, n\,-\,2} \,\cdot\,[\ \sec\,\theta\ ]^{\, 2} \ d\,\theta \ }\over { \Lambda^{2\,n\,+\,2}\,\cdot\, [\ \sec\,\theta\ ]^{\,2\,n\,+\,2}  }} \ \right]\ d\,Y_n \\[0.1in]
& \ &\hspace*{3in} \left(\  \hat R \ = \ \Lambda \cdot \tan\,\theta\ \right)  \\[0.1in]
& = & \left[\ \omega_{n\,-\,1}\,\cdot\,\int_0^{\pi\over 2} \ \ [\ \sin\,\theta\ ]^{\,n\,-\,2}\,\cdot\,[\ \cos\,\theta\ ]^{\,n\,+\,2} \ d\,\theta\ \right]\,\cdot\,\left[\ 2\,\cdot\,\int_0^\infty {{t^\ell \ d\,t }\over {  ( \ \sqrt{\,1 \ + \ t^2\ } \ )^{\,n\,+\,3} \  }} \ \right] \ .\\[0.1in]
& \ & \hspace*{1in}\ \  ( \ \uparrow \ \ {\mbox{coming \ \ from }} \ \ {\tilde R}^{\,n \, - \, 2}  \ \ \ ) \hspace*{0.2in}\  (\ \uparrow \ \ {\mbox{via \ \ symmetry}} \ )\\
\end{eqnarray*}
In the above,
$$
  t \ = \ |\,Y_n\,| \ \ \ \ \ \  {\mbox{and}} \ \ \ \ \ \ \ {\tilde R}^2 \ = \ \Vert\ {\tilde Y}\,\Vert^2 \ \,= \ Y_1^2 \ + \ \cdot\,\cdot\,\cdot\,\ + \ Y_{n\,-\,1}^2\ .
$$

\newpage

Likewise,\\[0.1in]
(\,A.13.25\,)
\begin{eqnarray*}
& \ & \int_{\R^n} |\,Y_n\,|^{\,\ell\,-\,2 }\,\cdot\,\left(\ {1\over { 1 \ + \ \Vert\ Y\,\Vert^2 }}\ \right)^{\!\!n } \,d\,Y \\[0.2in]
& = & \left[\ \omega_{n\,-\,1}\,\cdot\,\int_0^{\pi\over 2} \ \ [\ \sin\,\theta\ ]^{\,n\,-\,2}\,\cdot\,[\ \cos\,\theta\ ]^{\,n } \ d\,\theta\ \right]\,\cdot\,\left[\ 2\,\cdot\,\int_0^\infty {{t^{\,\ell\,-\,2} \ d\,t }\over {   ( \ \sqrt{\,1 \ + \ t^2\ } \ )^{\,n\,+\,1}  }} \ \right] \ .
\end{eqnarray*}

\vspace*{0.2in}

{\it Integrals involving $\,\theta\,.\,$} \ \ We want to express both integrals
$$
\int_0^{\pi\over 2} \ \ [\ \sin\,\theta\ ]^{\,n\,-\,2}\,\cdot\,[\ \cos\,\theta\ ]^{\,n } \ d\,\theta \leqno (\,A.13.26\,)
$$
and
$$
\int_0^{\pi\over 2} \ \ [\ \sin\,\theta\ ]^{\,n\,-\,2}\,\cdot\,[\ \cos\,\theta\ ]^{\,n\,+\,2 } \ d\,\theta\ , \leqno (\,A.13.27\,)
$$
in terms of the ``\,unit\," in this part of the calculation\,:\\[0.1in]
  (\,A.13.28\,)
$$
\int_0^{\pi\over 2} \ \ [\ \sin\,\theta\ ]^{\,n\,-\,2}\,\cdot\,[\ \cos\,\theta\ ]^{\,n\,-\,2 } \ d\,\theta \ = \  {1\over 2^{n\,-\,1}} \,\cdot\,\int_0^{\pi } \ \ [\ \sin\,\psi\ ]^{\,n\,-\,2}   \ d\,\psi  \ \  \ \ \ \ \ \ \ \ \ \ \ \ \ \ \ \ \ \  (\ \psi \ = \ 2\,\theta\ )\ .
$$
We start with
\begin{eqnarray*}
& \ &   \int_0^{\pi\over 2} \ \ [\ \sin\,\theta\ ]^{\,n\,-\,2}\,\cdot\,[\ \cos\,\theta\ ]^{\,n } \ d\,\theta \ = \ \int_0^{\pi\over 2} \ \ [\ \sin\,\theta\ ]^{\,n\,-\,2}\,\cdot\,[\ \cos\,\theta\ ]^{\,n\,-\,2 }\,\cdot\, \,(\,1 \ - \ \sin^2\,\theta\,)\ d\,\theta\\[0.2in]
& = &  \int_0^{\pi\over 2} \ \ [\ \sin\,\theta\ ]^{\,n\,-\,2}\,\cdot\,[\ \cos\,\theta\ ]^{\,n\,-\,2 } \ d\,\theta \ - \   \int_0^{\pi\over 2} \ \ [\ \sin\,\theta\ ]^{\,n }\,\cdot\,[\ \cos\,\theta\ ]^{\,n\,-\,2 } \ d\,\theta\\[0.2in]
\Longrightarrow & \ & \int_0^{\pi\over 2} \ \ [\ \sin\,\theta\ ]^{\,n\,-\,2}\,\cdot\,[\ \cos\,\theta\ ]^{\,n } \ d\,\theta \ + \  \int_0^{\pi\over 2} \ \ [\ \sin\,\theta\ ]^{\,n }\,\cdot\,[\ \cos\,\theta\ ]^{\,n\,-\,2 } \ d\,\theta\\[0.2in]
& = & \int_0^{\pi\over 2} \ \ [\ \sin\,\theta\ ]^{\,n\,-\,2}\,\cdot\,[\ \cos\,\theta\ ]^{\,n\,-\,2 }  \ d\,\theta \ .
\end{eqnarray*}

\newpage

With the change of variable
$$
\theta \ \to \  {\pi\over 2} \ - \ \phi\,,\,
$$
we have (\,the ``\,mirror reflection\,"\ )
\begin{eqnarray*}
 & \ & \int_0^{\pi\over 2} \ \ [\ \sin\,\theta\ ]^{\,n\,-\,2}\,\cdot\,[\ \cos\,\theta\ ]^{\,n } \ d\,\theta \\[0.2in]
  & = & \int^0_{\pi\over 2} \ \ \left[\ \sin\,\left(\, {\pi\over 2} \ - \ \phi \right)\ \right]^{\,n\,-\,2}\cdot\,\left[\ \cos\,\left(\, {\pi\over 2} \ - \ \phi \right)\ \right]^{\,n } \ ( \, -\,d\,\phi\,) \\[0.2in]
  & = & - \, \int^0_{\pi\over 2} \ \ \left[\ \cos\,  \phi \ \right]^{\,n\,-\,2}\cdot\,\left[\ \sin\, \phi \ \right]^{\,n } \  d\,\phi\\[0.2in]
  & = & \int_0^{\pi\over 2} \ \ \left[\ \cos\,  \phi \ \right]^{\,n\,-\,2}\cdot\,\left[\ \sin\, \phi \ \right]^{\,n } \  d\,\phi\\[0.2in]
  & = & \int_0^{\pi\over 2} \ \ \left[\ \cos\,  \theta \ \right]^{\,n\,-\,2}\cdot\,\left[\ \sin\,\theta \ \right]^{\,n } \  d\theta\ .
\end{eqnarray*}
It follows from (\,A.13.28\,) that \\[0.1in]
(\,A.13.29\,)
\begin{eqnarray*}
& \ & 2 \times  \int_0^{\pi\over 2} \ \ [\ \sin\,\theta\ ]^{\,n\,-\,2}\,\cdot\,[\ \cos\,\theta\ ]^{\,n } \ d\,\theta \ = \ \int_0^{\pi\over 2} \ \ [\ \sin\,\theta\ ]^{\,n\,-\,2}\,\cdot\,[\ \cos\,\theta\ ]^{\,n\,-\,2 }  \ d\,\theta\\[0.2in]
\Longrightarrow & \ &  \int_0^{\pi\over 2} \ \ [\ \sin\,\theta\ ]^{\,n\,-\,2}\,\cdot\,[\ \cos\,\theta\ ]^{\,n } \ d\,\theta\\[0.2in]
& = &  {1\over 2} \,\cdot\,  \int_0^{\pi\over 2} \ \ [\ \sin\,\theta\ ]^{\,n\,-\,2}\,\cdot\,[\ \cos\,\theta\ ]^{\,n\,-\,2 }  \ d\,\theta\\[0.2in]
& = &  {1\over 2} \,\cdot\,   {1\over 2^{\,n \ - \ 2} } \cdot \int_0^{\pi\over 2} \ \ [\ \sin\,(\,2\,\theta\,)\ ]^{\,n\,-\,2}  \ d\,\theta\\[0.2in]
&  = &  {1\over 2^{n}} \,\cdot\,\int_0^{\pi } \ \ [\ \sin\,\psi\ ]^{\,n\,-\,2}   \ d\,\psi  \ \ \ \ \ \ \ \ \ \ \ \ \ \ \ \ \ \ \ \ \  \ (\,\psi \ = \ 2\,\theta\,)\, .\\
\end{eqnarray*}

\newpage

Next, let us consider the other integral in (\,A.13.27\,)\\[0.1in]
(\,A.13.30\,)
\begin{eqnarray*}
& \ &   \int_0^{\pi\over 2} \ \ [\ \sin\,\theta\ ]^{\,n\,-\,2}\,\cdot\,[\ \cos\,\theta\ ]^{\,n\,+\,2 } \ d\,\theta \ = \ \int_0^{\pi\over 2} \ \ [\ \sin\,\theta\ ]^{\,n\,-\,2}\,\cdot\,[\ \cos\,\theta\ ]^{\,n  }\,\cdot\,\,(\,1 \ - \ \sin^2\,\theta\,)\ d\,\theta\\[0.2in]
& = &  \int_0^{\pi\over 2} \ \ [\ \sin\,\theta\ ]^{\,n\,-\,2}\,\cdot\,[\ \cos\,\theta\ ]^{\,n } \ d\,\theta \ - \   \int_0^{\pi\over 2} \ \ [\ \sin\,\theta\ ]^{\,n }\,\cdot\,[\ \cos\,\theta\ ]^{\,n } \ d\,\theta\\[0.2in]
& = &  {1\over 2}\,\cdot\,\int_0^{\pi\over 2} \ \ [\ \sin\,\theta\ ]^{\,n\,-\,2}\,\cdot\,[\ \cos\,\theta\ ]^{\,n\,-\,2 } \ d\,\theta \ - \   \int_0^{\pi\over 2} \ \ [\ \sin\,\theta\ ]^{\,n }\,\cdot\,[\ \cos\,\theta\ ]^{\,n } \ d\,\theta\\[0.2in]
& = & {1\over 2}\,\cdot\,{1\over {2^{n\,-\,2} }}\,\cdot\,{1\over 2}\,\cdot\,\int_0^{\pi\over 2} \ \ [\ \sin\,(\,2\,\theta\,)\ ]^{\,n\,-\,2}   \ d\,(\,2\,\theta\,) \ - \ {1\over {2^{n } }}\,\cdot\,{1\over 2}\,\cdot\, \int_0^{\pi\over 2} \ \ [\ \sin\,(\,2\,\theta\,)\ ]^{\,n }   \ d\,(\,2\,\theta\,)\\[0.2in]
& = & {1\over 2^n} \,\cdot\,\left[\ \int_0^{\pi } \ \ [\ \sin\,\psi\ ]^{\,n\,-\,2}   \ d\,\psi  \ - \ {1\over 2}\,\cdot\, \int_0^{\pi } \ \ [\ \sin\, \psi \ ]^{\,n } \ d\,\psi \ \right]\ \ \ \ \ \ \ \ \ \ \ \ \ \ \ \ \  (\ \psi \ = \ 2\,\theta\ )\ .
\end{eqnarray*}
\hspace*{1in}$\uparrow \ \ {\bf A}_{\,(\,A.13.16\,)}  \ \ \ \ \ \ \ \ \ \ \ \ \ \ \ \ \ \ \ \ \ \ \uparrow \ \ {\bf B}_{\,(\,A.13.16\,)} $\\[0.2in]
As  $\,{\bf A}_{\,(\,A.13.16\,)}\,$ is already in terms of the ``\,unit\,"\,,\, we move on to $\,{\bf B}_{\,(\,A.13.16\,)}\,$ first.  \\[0.1in]
(\,A.13.31\,)
\begin{eqnarray*}
& \ &     \int_0^{\pi } \ \ [\ \sin\, \psi \ ]^{\,n } \ d\,\psi  \\[0.2in]
& = &  -\,     \int_0^{\pi } \ \ [\ \sin\,\psi\ ]^{\,n\,-\,1}   \ d\,[\ \cos\,\psi\ ] \\[0.2in]
& = &  (\,-\,1\,)\,\cdot\,  (\,-\,1\,)\,\cdot\,(\,n\,-\,1\,)\,\cdot\,   \int_0^{\pi } \ \ [\ \sin\, \psi \ ]^{\,n\,-\,2 }\cdot [\ \cos\,\psi\ ]^2 \ d\,\psi\\[0.2in]
& = &  (\,n\,-\,1\,)\,\cdot\,    \int_0^{\pi } \ \ [\ \sin\, \psi \ ]^{\,n\,-\,2 }\cdot [\ 1 \ - \  \sin^2\,\psi\ ]  \ d\,\psi\\[0.2in]
& = &(\,n\,-\,1\,)\,\cdot\,  \int_0^{\pi } \ \ [\ \sin\, \psi \ ]^{\,n\,-\,2 }   \ d\,\psi \ - \ (\,n\,-\,1\,)\,\cdot\,\int_0^{\pi } \ \ [\ \sin\, \psi \ ]^{\,n } \ d\,\psi \\[0.2in]
\Longrightarrow \ \ & \ &     \int_0^{\pi } \ \ [\ \sin\, \psi \ ]^{\,n } \ d\,\psi  \ = \ \left(\,{{n\ - \ 1}\over n} \, \right) \cdot\int_0^{\pi } \ \ [\ \sin\, \psi \ ]^{\,n\,-\,2 }   \ d\,\psi \ .
\end{eqnarray*}

\newpage

By combining (\,A.13.30\,) and (\,A.13.31\,)\,,\, we obtain
\\[0.1in]
(\,A.13.32\,)
\begin{eqnarray*}
 \int_0^{\pi\over 2} \ \ [\ \sin\,\theta\ ]^{\,n\,-\,2}\,\cdot\,[\ \cos\,\theta\ ]^{\,n\,+\,2 } \ d\,\theta
&  = &    {1\over 2^n}\,\cdot\, \int_0^{\pi } \ \ [\ \sin\,\psi\ ]^{\,n\,-\,2}   \ d\,\psi   \left[\ 1 \ - \ {1\over 2} \left(\ \ {{n\ - \ 1}\over n} \ \right) \ \right]\\[0.2in]
& = & \left(\ \ {{n\ + \ 1}\over {2\,n}} \ \right)\,\cdot\, {1\over {2^n}} \cdot \int_0^{\pi } \ \ [\ \sin\,\psi\ ]^{\,n\,-\,2}   \ d\,\psi\ .
\end{eqnarray*}

\vspace*{0.35in}

{\it Integrals involving $\,t\,.\,$} \ \ We seek to show that
$$
  \int_0^\infty \ {{t^{\,\ell } \ dt }\over {( \, \sqrt{\,1 \ + \ t^2\ } \ )^{\,n\,+\,3}  }} \ = \  \left(\ \ {{\ell  \ - \ 1  }\over {n \ + \ 1}} \ \right)\cdot  \int_0^\infty \ {{t^{\,\ell \ - \ 2} \ dt }\over {( \, \sqrt{\,1 \ + \ t^2\ } \ )^{\,n\,+\,1}  }} \ . \leqno (\,A.13.33\,)
$$
We discern this ``\,shift of order\," via the integration by parts formula\,:
\begin{eqnarray*}
& \ &   \int_0^\infty \ {{t^{\,\ell \ - \ 2} \ dt }\over {( \, \sqrt{\,1 \ + \ t^2\ } \ )^{\,n\,+\,1}  }} \\[0.2in]
& = & {1\over 2}\,\cdot\,  \int_0^\infty \ {{s^{\,{{\ell \ - \ 2}\over 2} \ - \ {1\over 2} } \ ds }\over {( \, \sqrt{\,1 \ + \ s \ } \ )^{\,n\,+\,1}  }} \ \ \  \ \left(\ \ s \ = \ t^2\ ,\ \ \ \ \  dt \ = \ {1\over 2}\,\cdot\, {{ds}\over {\sqrt{\,s\,}} }\ \right)\\[0.2in]
& = & {1\over 2}\,\cdot\,\left(\ \ {2\over {\ell \ - \ 1}} \ \right)\,\cdot\, \int_0^\infty \ {{ {d\,\left( \ s^{\,{{\ell \ - \ 1}\over 2} }\ \right)  }   } \over { \  (\,1 \ + \ s \  \, )^{\,{{n\,+\,1}\over 2}} \  }}\\[0.2in]
& = & (\,-\,1\,)\cdot \left(\ \ {1\over {\ell \ - \ 1}} \ \right)\,\cdot\,\left(\ \ -\,{{n\ + \ 1}\over 2} \ \right)\,\cdot\, \int_0^\infty \ {{ { s^{\,{{\ell \ - \ 1}\over 2} }\  ds }   } \over { \  (\,1 \ + \ s \  \, )^{\,{{n\,+\,1}\over 2} \ + \ 1} \  }}\\[0.1in]
& \ & \hspace*{2in} (\ \uparrow \ \ {\mbox{integration \ \ by \ \ parts}}\ )\\[0.1in]
& = & \left(\ \ {{n \ + \ 1 }\over {\ell  \ - \ 1}} \ \right)\,\cdot\,\int_0^\infty \
 {{ {\  s^{\,{{\ell  }\over 2}} \cdot \left( \ {1\over 2}\,\cdot\, {{ds}\over {\sqrt{\,s\,}} } \ \right) \  } \over { \  (\,1 \ + \ s \  \, )^{\,{{n\,+\,3}\over 2}  } \  }}}\\
&  = &  \left(\ \ {{n \ + \ 1 }\over {\ell  \ - \ 1}} \ \right) \,\cdot\,\int_0^\infty \ {{t^{\,\ell } \ dt }\over {\ ( \, \sqrt{\,1 \ + \ t^2\ } \ )^{\,n\,+\,3} \  }}  \ .\\
\end{eqnarray*}
Hence we obtain (\,A.13.33\,)\,.\,

\newpage

\hspace*{0.5in}Combining (\,A.13.24,)\,,\, (\,A.13.25\,)\,,\, (\,A.13.29\,)\,,\, (\,A.13.32\,)\,,\,  together with\,  (\,A.13.33\,)\,,\, we arrive at
\begin{eqnarray*}
& \ & \int_{\R^n} |\,Y_n\,|^{\,\ell }\,\cdot\,\left(\ {1\over { 1 \ + \ \Vert\ Y\,\Vert^2 }}\ \right)^{\!\!n\,+\,1} \,d\,Y\\[0.2in] & = &\left(\ \ {{n  \ + \ 1  }\over { 2\,n  }} \ \right) \cdot \left( \ {{\ell  \ - \ 1  }\over {n \ + \ 1}} \ \right)\,\cdot \int_{\R^n} |\,Y_n\,|^{\,\ell\ - \ 2}\,\cdot\,\left(\ {1\over { 1 \ + \ \Vert\ Y\,\Vert^2 }}\ \right)^{\!\!n} \,d\,Y  \\[0.2in]
& = & \left(\ \ {{\ell  \ - \ 1  }\over {2\,n}} \ \right)\cdot  \int_{\R^n} |\,Y_n\,|^{\,\ell\ - \ 2}\,\cdot\, \left(\ {1\over { 1 \ + \ \Vert\ Y\,\Vert^2 }}\ \right)^{\!\!n} \,d\,Y \ ,
\end{eqnarray*}
which is what we want in (\,A.13.21\,)\,.

\vspace*{-6.8in}

\hspace*{1.5in}{\Large{\bf [\,for the e\,-\,Appendix\,]}}


\end{document}